\documentclass{memo-l-single}

\usepackage[T1]{fontenc}
\usepackage{lmodern}
\usepackage[utf8]{inputenc}
\usepackage[english]{babel}
\usepackage[protrusion=true,expansion=false]{microtype}
\usepackage{mathtools,amssymb,amsthm}
\usepackage{array,booktabs,tabularx,longtable}
\usepackage{enumitem}
\usepackage{tikz-cd}
\usepackage{imakeidx}
\usepackage{graphicx}
\usepackage{xurl}

\usepackage[
  hidelinks,
  pdfborder={0 0 0},
  pdfdisplaydoctitle=true,
  bookmarksnumbered=true,
  bookmarksopen=true,
  bookmarksopenlevel=1
]{hyperref}
\usepackage{bookmark}

\hypersetup{
  pdftitle={Homology Nilpotency in Rational Homotopy Theory: Cell Attachments, Retractive Towers, Divergence, and Stabilization},
  pdfauthor={Paul-Eugène Parent},
  pdfsubject={Rational homotopy theory, homology nilpotency, rational LS category, and Andre-Quillen methods},
  pdfkeywords={homology nilpotency, homotopical nil-length, rational LS category, Sullivan model, Andre-Quillen homology}
}
\pdfobj stream attr{
/Type /XObject
/Subtype /Image
/Width 1024
/Height 1536
/ColorSpace /DeviceGray
/BitsPerComponent 8
/Filter [/ASCIIHexDecode /FlateDecode]
}{
78DAECDD7554945BDF07FCB5869CEEEE640266868E6180A1BBA45B1AE90E012514045450ECEE46B13D767777EBB18FAD3433C3ABE7DCF7733FCFFBE7
7B8E27EE777FD682A5CB3F645DC377FF7E7B5FD7DED7D818000000000000000000000000000000000000000000000000000000000000000000000000
000000000000000000000000000000000000000000000000000000000000000000000000000000000000000000000000000000000000000000000000
0000C0BFE8FFD777E0475EE8EF7EBBD2FA5FBF818B0E00FF3F88BD4EA7D36AB5DFBE7FFFD3F7BFEA7F03060000F8EF4EBFEE7BE4BF0F02DA9181AF5F
FA07FAFB47F563FA7F031708F8CBFB7FF05BF8E3D2AFFD56F3875E1E5BD1B268DEFAC3BBAFECDD75E4D2CB4F43BA7FB500E00A0120FFFFBD9756ABD7
7FB9D1BB38596A270C4E4B6FF775109BA9FDF20FBFD1E9477F1B01C0350280FFDACE7FAC7F7DA64A46A0842432894161181C012B4C8E0D5444AE7FAF
FBD619E840FE81BF41FD077EC055FD7E5975D7720878918C641DCCC248F8503C03611B27C1092D63BCCBAEEA74DF5B0070A100E0BFB1F6EB5FDE7DD4
2CA15098663216158EC4B3E0682AC93D8ECE9609D3A6585BE45D1B18D3810E0000FE0BD3AF1F1B391E57E6052773E864970409164BC723E0587AB827
CCB34E8CB2C1E0D98AC5FB6FF5EB41FE01E0BF2DFDDF9AFFABA5DEE65608DB22295E33C3070FC5D00848ACD2562108A956529A6AE9743999EE1CB4EB
CB18C83F00FC974DFCFB8F763B8A1D48305AAD17029B93668AF0E021101432DD29AF30575CB0607190190B0DA5DB73D236BE071D0000FC7715FF2B71
9161728100C50E7085E1589164D8E4062C8C668527D64F62098B5728E1686E7C6DDDDA240CBBE8ECA016E41F00FE4BD2FFEDABAF931558AB642B48B2
1C4B28D59545C65677924CE9420A293C02E9511ECE54865817AFAF756321B93E1EB1BDBA6F23060000FF0DF11F7DBA3A06C79B8433B1B1762D742293
6D1224245930D19417C040474E40F994BB6313279BB9343A21602C71883F97E575690CDC0404807FF2845FFFEF89FFC7CE14111EEEC6838A1DE041A1
26F676AE5608511C1661EEC0E7557B0A539B4C108D0DB4E8855104A97FA4C20A6519C3ECF802260000F08FCEBF5EAFD3EA75BAA7B91C1983E4E18F50
27F2F831E6442B615E00A6A889C07310864D9F685F1B44B5696856A4D44A49D6A91AB902CDDFD44E0DBD00F20F00FFE0F4FFBAB757AF1B3C944A9209
4D53A69351212E464ED144730C42C5F25D2CC1B151E3D60632DA3CA099732268724B96C83152CAF2A68B366E150B6D5A0646C1080000FFC8BAFFEBE6
7EEDF0D0DBC72FB7FB108539D6B822BE89298D888F8EA38BB138147169039AC5762D91C13C9C5573E732887CE5A4DD935D3C824AD2F179ED38BB0227
B7C303E0116C00F867C6FF7BFAFB3F9C5A777A5B86088BD710788E1023BA14C18F5610894C8A91DF060A9A8C897042C82C020E5723C82521BDAF16DB
994565C6C843A2E8768E4C1E356DED08680000E01F197FAD7678F0C38ECA0DC713C82424824D0E0E308029CCD1DEE9542285CF66CE1D876689EDDD24
F2C4C06257A8F5B2FCE0FD110845636DE4B8C20617099E51EA489E3DEB0CB8050000FFD0FC6B1F2F8A6A3D15C314E105BE5CA85FB089218E040D0B37
C6B3A9C659F3A12648BA138B91E0656B65E1BDD213DF1EC12C5ED9E2631E11A9B0C499AFAF870A27DAB70C83331800E01F95FFB1DFE2AFD31E4C7388
D8BF4CC061DB77A4E1891E85663C2B2CB9D0054A6453ACCF7A1B2B828462654095D472F9C27473FECA1E4DC1CA609EA0B82B70FCBABCE65C13042F8E
5F313C0A4E6101807F58FAF57AADBE7F9997543EF5B833CD026D974F2651F851A69A00036C0AD71483466C9D6D283AD88C2685867258AB776A98963F
6D51723552BA57F1C2889C5BBD7E66483A2E76BB3AF8E9E0AF678182650000F807557F9DEEEBD4702759F0AE6601574321E0610AAEB38B81A387A955
882986088FBC696EE4546C8656A7F270D37AADE88AAA7A0E5A1D621DDCD3C0763BBB43849633B9F25D3D3CD71BC303A00300807FD2ECFFFB99DEF733
89D64249ED3A9905DB9A6BC5A70AD1C15E54371A26D70F8232221E49354012317807773CB76B119B57E0EE48904EA89788CA0A04617B3709D1EE5242
4E745B2A3EE6C4F99B7DC3E0405000F8E70C005A6DFFEE481A47CA8B5993CEF6B4E7F1139C983C7E888FDC9D800B32476221ED7788A898682EDDDF52
5CD5C6A4649530384E93DB5D6DE7B7F835FCDCCA25F9059977EF8BD8D1EAE9EF679DFC7C503F06F20F00FF90FC7FABFE9BD49E85CE1861418790E46C
4DE6B8326042175B16CD0C650247A00D54AF5C2016EEF60A9F5CF3307BB475BE2B95ABAC4A43594ECB09FDE9E57C89227B9A5DC4D29443E7E2DD612A
47CE7AED28A8FF00F08FC9BFFE96869495C0C45ACD08C18B82994C71820A6919C4C112CDB8780318D4F8C07C03ACC22A7B4391848086B9253ACBE2D5
D69676138ED6459F7E53420D5958E3C3CB9D73FABC8A2D8DBCE1436C1D19D682FC03C03FA3FBD78DBD8BC771E41222C6B7940813C4E048A4092EA662
078E258FCCA61B9A42DCCFE00C4DC85D3DF96434C5CE360667D7DBCAE439CEDF1A13FDEC4280B2B33B91EAD7BD65E14209CD51736E038518BDFF8B7E
140C0000F00FA8FD3ABD6EB895EC988461D3788DFE484240AA19071540835A3A58AAD90485106A64D0B3D108619631538396BA874E90D1532FCC27D2
27ADCCC0563F38EA1BB566ABADFFA40D7392497C8B26EFFAA3E604D1B8A9A7467420FF00F00F68FDF5A3C3DB04FC6A0716D144388B65822970766051
ACF1282E5F262472B86804C4F1A5B789A6299C848B4C748813D0AAAF374039CB9629151D37D21CD79C89F5AA5B332D8AE7139C3F235D30B7052DB4D4
4C99F5CB1058010480BF7BFE75DA8F0F5E9F4F44DA6752145463692E1ACA9DE1E2C6B75363997C9E375F10AD809B181D384694A5C850B8906A336739
B5EBDE7803D1D62662CCED332AD4A4A5F6FEAB9A1D04A9135764264FACDA7D20D127273AACB3FAF51058010480BF7DF3FFF1A78D67D2A8B6096C7B8A
0994CB373656E44983C8167204092799246365CB8C21516FED0CB14C24262C8727B2C0B6DD8A3472DE908BAE7AD3CB14B56799C41FC9A5C5B46FA971
0E59FBF0CB69AFC6E02A8BD08959573F8F81D38001E06F1D7FBDFEDDEEBAB5AD3C6443356B1CCF14E780809904E4621C700802920EB79BC93316E20C
702F561991A908926B2A1FAF11FAAE0A83FBAFF4A5CC7D530F375F13436D3818C59D78A6CEDB73CEE89DDD572CC3DB7CD74A9776969D7D3F066E0100
7FE35F7EF084EAB7EAFFAA277BC21A5F0261861F4D4D352279C160866E753C2789319AC03716142121C6A69069FD42280E46488921D35C4343A215AA
E26A0BDAD2E311C69A4E8D7CD30E6BF1ACC5F10E1DCFFAE7665E98254D1076BB955C1937FEC087212D98010020FD7FE3E65FDFBF25CDAFB3518C41DB
B2782C1C016B614240F022A2C2EDE94C9A182FCFE1C23110DE40972106E5571CC5A47B4F0C34F7CC2A49C73ACC9A6606F39CA4A9393C8928689AE551
7A7C646F4CCDC0334B61F0D419158F22DDF20F0E0E0C825B80C0DF36FFE012E8B463D7F3EC5297B8E0E1248634920585115438131C9A13EF89E79913
314A958BB7890164D6579E212E6DAA078A1616E2220C689A2865792F8923335DF24BD6EF0AC2B96E680CDED437B82AB857F7D08F91B02D34FF55353F
A5F6C4A7A1970F46C0000080FCFF3DFB1FBD56FF79A247C8FC8966584AB49DDF54111A4A53E1A166380B0F169C64C64286274BBE4D00D09F961B20B2
A7314CF961133C84E336A6B324894D7E249AA663F38E895C42E6E9F599BBFAFB7A1AEFE8CE3832A31ECD8E7BD9C50CAA5E76FB97D1FBC75FF5E9C169
4000F0376CFEF55AED50578443C65A5F06D4AD59E09B8E2642F15438BBD9DED6150BE5F3F0685F377E2C0492D1CF8390A25930B30C3F67AA736F039B
3DBED18D24A85CB7B5CD1EE35077AE2164D7D74F6736FC32F4939219747AB9DFB3CFF6AAE8F6938F5F7C3AB5FAD0C5511D180000E0EF58FE75876CBC
C366773BF271F61E9CAA34281A8FA322A4532D38DE6292404EE266383824410C9FF642506814CACA27404A735CDD4C60C465D85123965D9C1C245086
2DDD9415B6BBFFF3CD9F3FF52D16F20BE62D4D7B399C661757B8FBD12FAF5E9C5ABD6BEB3DEDF7FC83210000FE56E9FFD6FDDFF466DA54ACCAE631C8
381CBBC0C688888572713CA129DD9C6B8115E2711909C17C88FF171EC4084AE4FB9299F6B2FA4E16DE3FC59C5CF772B3BFA5DCA2687A635CC3998FF7
EF3DEA7F52C674AE9A57D0F058DB68E613BFE8E1F3176F6FF7CC5ED8BCE2D78D40E0340000F85BB5FFDAB117490C46F4B29DE328742E91289FCC35E5
93515C969C02139829DD58523EC3ADB8886A70EB9611DA40186B45C4B9102D9C89B8D430267DF9BB4EA950E997DB9955B0B9EFDD9CD9AFFB8EA8D829
4B0AD27A46DE254647C44D39F7EEF59B57E7CADC5B82131F7E7FC058FFDB8960E0BA03C0DF824EAB6F17D89A272C9E2181C1B13838A79C0563B01004
9ABD1541C092A85CB8B61294BBBBA1B3361102752BF7C418D952E18E507671191CBFF0592A59E01338B92479DACD81379D53AE7D9C26F05F3D3564C2
A5D11731B53D31B56B1EBD78FDF1F242A7F402BE64D157AD4EFFEF11000080BFC30440AB7FE524B356E7CE4EC22259164426DD01851709F06C84B986
2A5490DD33444A29C1D11176F28481B1F3743B08DC8A6FE56F4B296E407B4FAA5763158113537D337FEA1BE999BAFCEDC914F5F49E52F7796F462EC7
D4BDA98C99BDF3D58B57B737A7869D54B169D54FFF27FF0000FC2DBA7F9D6E3857E640B4CDAC0D30A157449008780C9A28A4F03910A12D54EA8AB0CF
139388347B94DFA00BC4502581E015E60C5B3C756219A660832394E55E96E1E9DFF16AF0D6A4392F9FCD8DEF38B825ACFAFCC0E71ECF82DE9AC8CA19
3BEE3EB8BF3BD3F7FC71269E16B6FDF3BF67000000FC0DE2AFD38D8D5D12D8D8F182B293F988905C069341A123197296A399A1B53B46E14557572A55
16745BD285314BA85B2C9EA646D80570FD77CC342B5A2134E215568429F32F8EDE6B2EFD69F4FEF8B853EF274FD8FD61E04E9155414F7642494BD7DE
B3B7B7D508165EF12551A9FE2D17F5A0010080BF4FF1D7EBFB7E294113501641FE4AB4918F3D4C2A41D0904426C39A03E17A4299AE02BB09FCA8380E
8CF26999B1B188892423241EA64E8F7E12FAD6118C5C8A2698F96D1EBCBF2C62D1ADE73BFD0A1F9E4999FAA27F6097BB79F38A82F4F4D4EACE35DB77
1770B35F46E1F951569A86EE9723A0010080BF4DEF3FFAFC66B7B9CC02EF9892A61423B342A876961426162175B061C1C5617889A7D0371F27891543
A7BCA79A404D8DE0505B2FA664D71633BE1D1C19D2AD26E63EEB3B3ABEF3E2CBD6F894E56F36DA757D1A7ED6C0759ADF3361F2BAC2FC8EF9AB565659
B93F3B4CA54424B11591B6D96F74E0302000F85BC47F74F4EBC123171330814512DBFC28B29149723C4B4661918DD0016931248826036DAB2279A642
61721CE2C5625304068EE45ABA602075070548172131B652489CFBF1DEC4CC9DA3F7436DD71F785AA1D8383AB0CF931EB163716AEBCE0CBFC6F6BC86
0247C5D5C76AB2C2C38CCAB7B4B53E3D02360200C0DF62EE3FDAB7A762FF660D9424B34A4B25C08C8CA26A041C3C958D4206D6C7F30CA2E3204C0766
58360E47361DAF15C25070B44DAE27869A35478C1059797B07A3DC4EBF6D73D05C183DE31873FECDFAD8904B630F32988E4B777415EDFCB421A17E66
71CADA184ED783208AC29121B590B0ECE8A55FF4E02C1000F81BE47F74F85165E3C14E27281C464F0933313435D09410C514041301530408E810474B
08DA19E9960DC5604DAE1F83A24D4D02729C718CF825115059A885C60C55D0774A23AE6E3B33433DFDEB928CC68EC7C3CBEC6533CF1E685AF8706853
F2B22E7F8BB69D2CCF5D4A8C8548C0120898581137F4E488166C040280BFBEFC6B3FADCADC7C249D8444C225B5A53C13AC61540137D01CCBC7E0246C
A203C64501E178A2830A9178A4F590CA144E08C9303336F22BB7975885A55349AA456F67A05D97AF7657861EBCD5E1B1FEEC9CE5290E93EEBF59D973
7F70A0B7E2C2426F5CEAA375012B82B19A009ED29E1B1644A6F0CC167C1DD5EBC0D50780BFBAFC6B8F85A6F44E139A98D269FC40051D8A87AA13942A
4B01194D5733398570974008C5D63432CD0061B8E1B489A9B02A0665646DCBE5F8FAA965247EFE9C5DB1F29C9E248FE4B4355B22DA1FEEF0185F5378
6CE4C2E27343237D2D4DB7BB54E890A74F429B4F0606854BE5CE14F7F3F104A673F6E457233A907F00F88BE9742F326C0B36A7E2C8FC719E6C94A1A9
B18981438DAB77800A07637909D5AB1449F510FB48E3F25208DCF68B8F095A610185DB260871D6394E508BF4B6191E53DBEB5B54E11B8B038B57DE1D
3DEA9E7B65EEB9CFF35A1F0C8F3DCB9FFDB5D49913FEFA5D9E656D90D8436CE7EEAC3C3C17671D1E9FEB7C5A3F029EFF0580BFB8FCEB47179BDB55B4
3B600451FE3C67A611D6D8C0405661EBA3B680E1E9364CCB5972DB7A88D778487515C468CF0B381C0687E27D1C70C29C208EC4E7C07EA5D3EEFDA976
E6A1D313DC37BC197E365BEEB528A7F67A59E693B1B153BE3B3E56F947853C3CEA81127251967449A810517247C64CE12044F2F5FA1130FF0780BF38
FF63B78325AE2D795C1AC5CC411AAE44A28C89248BBC60A5888122539DBCA3678B25D30CA3CA4DE2BC2196436B8C90283233D01FA9189F486585AC4A
93A49E98269744B6944775DCEF3BD7A666DBCF98DB79A43CFFC5D8E772AB036F13CA4F8C2BCC219BB035348DB7D4D3176DBFC60F979143A1D83834F5
E9C01D0000F8ABCB7F8719533DC51F4F438AB2BD443C8429AA321616126CA361E10826EE45996D727491A15D060A0783940F6920DC10818D0C1D57CE
878AA2A8488F637B43E9B1F98B16361F1EB85F28613393969E9A5192B5EAED9D0D6555CF4EFA2F7B532152A051D62151E1210485052C69978A101443
A788F1048BCBBA51907F00F86BCBFFC31015DFBD41456661E40DF10A0196605A158B080908493787914C030A05B579EC8924D5641AD20872E325C2C8
5A214B5386C523F05149B6B8E6AF57AC68F9EB264F09A97B785C8D33F36DBBBF3BD3BBE1C5FBC6DA13A39FEAD3AFBE28A2F178ACE0FAE8E0EE68854A
50F2240BE798CF44D0A818A2D57E907F00F86BE3AFD377AB73491E45623C17C109F3F09122C8063E29C6AA78D7681B53B4897F88894726B6086F3511
6504B1E96F86408D98B164741086569607531FF8309DCF8EAB9E5D2DE3B66E73E5A5E41D7ED8E2947E6A604FCA52EDD8A5E0257D977DE8521F595497
BF5943B4953BABF25315CA2654C6A010B8D373C6F580470001E02F2EFF77ECAC72D001796C82906C29575A93987C93C062A85FB687AF051E4BC84D84
D3DD19E544C64C1788C1E13E862105E12A320E56E3C72758365C5DE88BE1159F981DC6F0DD7B21DABEB877C3461FBF8D9F2F55649FD16BD7D45D7CD7
2CA1AB53EA7D54D1252DF3D3122DAC77B4E0782EA136EE1C78580F073F7F68142C0000C05F59FEBFD42119A130EFF1280299ED8527118CE9165061BC
A995C641C6151823F2A2380E76D218A6E1A4348872703E044141C291B99568BB40CBA5BDC168BAF5A205753C61CBCFA71215455B3765D856BCBFD910
D7F97EE4627DF7F3256ABA755C53A29BC2B3F940B17986C2A0F527962CCE52E122CCB40F56A0C9730647753A70071000FEAAE2AF1B3EE4839644C223
53A1548A244984C661700CB8249328B32059DAAB88564D353661B1CE396A132F1B48DB3B22C41022F7B18EA2432DF1BC643AD12C6BC91A4759CE85CB
BB429CDA7EDE64A93AF47655D6BCDB034F8BFD9794BBE245412DD581E6E9D3B7E487547405601B763BCB53C6DB85A8E5FB5368689A70ED20A8FF00F0
D7E55F3BFA62020EEA1A62E437CE90481285324C11382886284B86E2C450736B47555A533C5F91C50A7587B110903BAB208604872909781303BCD0C3
91C28D0E6D2A1184EE7BDC632F89E979B34A9A74E744CEF42B83EFD7B914CE0AC4D1AD8B5BB28419738EF694CBAACE4DA4CD3EE880B5F1B7704A4447
6F222279E3693DFD5A3DA8FFFF198C01E0CFCEFF520A0AD61865ECAA36E209DC0B452412066BE3A19964CF509AD938861617CFF397B8942B7C22A91C
E3B8C10608BD63A5066162C4C84D72A52BA6B72B1C3D8AAF7F6EC720D46B1E770B6A3FF766EDFE3AF4A8D866C6653728CF3BBB55ADECFCB035D12E6F
4D972567CD2E33184535CEDEB78E59966A829C5D203A39A803F907803F3BF6FF36AABB632FE63116DB1BB104185F6640A590C625509DAD30B15D12A1
BBB934A5C5AE50CD732DE50A120968C37D0372085D638634C1AA2C0B38466E3BBBC57265FDFB4F5D7491A6FA4C8574E9BB4569E7B5237B82C42DBD4E
5056F0E40881D791172D727EE5834D4CD1AAC54C2CDB223190E01267E98634295D8476BA36A805FD3F00FC45F1D76BFB9A50A5A116D32D9004242B8F
EBD4E0C91233304424C2726DA9C48E42F5E8300F0E6325AFE45B34C810F2819F20086308DE1437A9884645A5AE5E2AB54E735DF7A61E2959326769B6
CBAEAF35F9B7C79E34DBF88E2F343392D714D8B3C6DFBAEA26554DFF72CCDAF3FC56164A93AA560B85D509E6745CFD742AA9F2E50038030800FE9AFC
EBF4A3A3B755B02A9120020D2518C23C69C2C93E141A868047C9A04553E164164E3C4BE69C078BECA5B05AA546953A6F08828BC5F09A9AF0285CD56E
3B92EBCE893EC191A4CC9DFB73AC622EDD4EA9FF30B8CDD725AB2B868A8D9F9F45F75CF1618599D7B482B73B587E3FEFE032BC33268BA829E5FE5618
E39A663453B2F163BF0EACFF03C05F30007C2FBCA343AB58D83C5A4E2B1643801AB318B673C751E8240A49D6454DEF3016C8482EBDAAA012D3F45D78
4E1B0776FB31128A6593AD27A49BD295597325866E49EDB64A43CD96F7F379D8AA1BDB833ABE9E4EB24BDB322D80E6D4BD4E2D6F7FF97E1236E5E1BC
298558C9866E335ADCC238855565998FABD238A3D848E6E470626044A703A78002C09F9EFFEFC1D369DF669A126DC9EB7713E128A43197219F1E0DE3
8AD030AB2D0CE71203BEB5B1FD2EF3C05CA3C46526CC22A8ABBED0008EC7F20BA28DF91D9EB1F65CB776194A5451F6E0748A216FF6CD12CDE26773D4
01EB7696095C8A56B798859F1ABCE147EF7CBFC83E401EB830966F595D1827CA5E3A99EC9566A05C41150658265ED08FFEF61A10F08900C09F9BFF6F
238056D76B47A21089AB0F99A1B148A4BFB9D9825C94D8864CF2DDE3A5AAC69A3BA2DCCF5869EAA045F34D24F1D09DA32C5318C7C5DF8EC468984463
BB2ECAF5823B9D19BA564E30723FB6499D74E46ABCC5E4230B1D1CA61D9A68C36E78FB751147F4D3D32C6AF8EC5939AE4E75C5157C61E7CD78B48D2F
53B95C4588124A5A5FEB467F1D86C6DE3E06930000F8F3F2FF2DFBDF66FF43C5F62E441C3C6E310D4936258F671AA66722047CAA24647238AEC982CB
87090FA8250D4699D546422947B7DF1841B361DBD81333F3B0CADCAA2A36C2F7F9C87E336351F79919F2A6C71B1C8296FE146B5E7CF35C2886B15177
350E9F726FAF85B8F5C2527B5668E1146BB2C3CE03BE38B98646983DC1C8C61AE5B47B7444FB2DFF3AEDD8F136701B0000FEBC01606CE817BD56FF3C
D0D51D4EA33BA77209142CD7968096A772ECD868555051066596DA4C45423726AA56E263D3612AE294B10863B40317971968512821572D0E26D15BBE
F635220C13F61C488CEFB9D7645F746AA19DEFEEA1DD0AA2C596CF1D0CD98A1B25D6C9DBCF4EA10BB267240AED5B8FAE1211EC3CDD6CEC7C603ECD48
7CC1F37F957F9D6E66F420C83F00FC59F1D78DBD7CACEFD71FE5CB6CE196D9C9F124340E6EA5446089690421CF44C4AA9C838D6510D57443AF5CAF35
229B18281BF1F03D0E8A27902A12713E168C9C222FA4CB31DD7E8D91B2F56C8566D2FDE3D156CBF6640ACBDE7EADA233023BF7B971F3AE6F0ECDBDB0
3C31769C6D51AE3D25EFC8E9589A79B493BD5265C188BD140FE56E1C1CD469B5DF1A91912B7E69FDBADF862500007E7CFE75879F8E0DEBE692243E84
88D2EC520A9E8888AB4193CD8B956E6C840B7DD2219C9D35D3CA17EF96E8F0937F7C351C6A35B6C49480C64FA8C57083D851AD026AF69BDB19506CDC
85459E7E3BAFB6F39C66AEB1B4EC193CE28B144E5991C6CE79B8A73CE7E0C011B6FB8E8E9204A2F5BCFBF52C5A5E51B2CC35C8227BD18C68383AE1C9
C8B0F65BFE47F577325415A78741F801E0CFA1D30FEDFC383A3A908965A989D1AEE3E26138122C73391A434E52A70B8C3504EF6E03B7953684744EA8
BD494FF4F45208A461CC1546B02A0E43288A50EEF904F399CB5AA948BB9EED29E4F223C78358DE1BAB28818F9ED47051562B965979EEBF3F2B61C7F0
5087B4687F9D6BB091E7F16DEE88D0D91D564EFE6A49EAD975684353E192AF03DA11ED886EE44985997D5CEED73170143000FC39E55F7F6EA77640F7
3280247425A5A6258C43E128B0A06A0C971AEBAA126034560A7F63E77B531041B46877D8328F087763E8D3A77004BFD1D3D4BAD1859BC9121C984546
D31A7A5659D92C7D5DC7E2A7AEC88725BDDD6D8124C76E68964CB97BB074DD9B916DF1DEAB7A34D6CD56EA950570CBEA1531A2BC6A1BCD9C0FFB5908
3825E9D6F0E0E8C88876F8680657691191F612E41F00FE9CF2AFEB6B3BA8FDA83B674524ABA03E69762238920653F8A384283B07A5BB29C2D78E4543
77CD33B424BB061A2C72301640BC471B8C104C96A9620E469A42569C9C45200777F5FAA96B6FBC6C24F36BBABD8C335E4DC111C32A96478DDB79A6BA
EBE9D899A4D075BBE7994DB83C3562893D3D6B738585C5C48E008FDD174E4F841BA3EC977CEC1F191CD6FDB2DC8F4F5331C645EE01270103C09F34FB
EF49793CD4A7DB614E26AA30D6622E078566B1923389761CB19F269B85ED4AE3A8618A485C8CB95D1A74931B4A0CDD3E626582869B4A2A5D3DDA2C2D
574D2408E37B93F93673D7CDCA14F1232779E36B5F1463C7EFBC59ADE9FED093B37FEC63A37DEBB323C9056B6FD7B84EB5F7DCBCCA1DE55FD764E5BB
731DCFBD5C84B39B7867B07F50DB7FB7C682401107C942C3570E813B0000F0C3C33FA6D78EBDF0C9F8DAF775A4838345D0881E1C8ED8C894A0985482
B2E013C22C233D043DA184895CDE386E82D42CD160560502CEFE70D5C8040A27FB59E42D91392F8A46C52F0C7347AA0FAC6793DC4202EB45D8799F26
D2E67E7C9267B3FBE9EC0977B4A7BDFC77BF684C58FDE194AB6D8767F6A1161A3D7F913725FAF82C3A322A83CE4C3BF0EDBF1FFC65471CDB9D800D4A
1033EDCBDE8C819D0000F027947F6DB770FAD097E1AFE522368DCA4832F78C348692155D53797E96C40817AF427AED38566E415C86C5380DBD9A5079
846A5C3E5C6980C3FABBE16A5A79F4B4304CD8C164ACA0F4F1F14056FA9E9DB1E6AC957D6DFC4DFD973C5CCE5CC899F5FAF36C65CAADBD89937F7E3F
13A7595439FD642ACA7A66939562E2896ABC3435C623ACE4C08BAF7DEF8ECD742384D793C8B1F6049438EE39580000801F17FB7FEFFAD3EAAFFB33E6
8EBCD53E0CB14AF06033EC88E1510686546AFC1C529017C1CA4795691A5081E6362C8EC0C747C12A5DB8CB85A6A7FB451054E0787C5A3381182B85FB
EEF547C9166EA9B754C69DB85787911CE9AFA72DEB5B63117B6D49FCFEFECB05018BCE76BACD18BA1F451ABF2979EDC9505CF0823C7CF8FA8DE3D092
796942F7EA3DCFDEBC39BF32D315CEAB0F33212A690434D3E90E78172800FCB8B6FF5F7B7E75DAE1291AD932ED07DD091F811919AEB1C3268C333066
715D36493C6328AE812E33B9219BD8DC80AE445A720EAEA21EE6827719380231362BC605B709891A1FB4DD2A2F047BD14C228AAC5246F90922EEBDC9
84A5EDCD92B79D6FAC79F86C5E5CD3E19369F1CB7F5E626F33BFABF8F83E0D2E644520B1F4760913E1BE74829947FEDE7B4FAFEF28752570304C4BAC
2943212508AD5CEF83FC03C08FAEFFBA11FD4D1F5BFA5AED6BED026F3E16C50E333752981918B079F0764F762AD66E3C727ABC7C16C391199DEBE693
868CEC86A10CDA869321501BA97F231E1F6403E3CCF4810B5B73F03C3B7B4F099ED4D97F5783AC7EDF1EB7FA7CCBE137AB221BB73D9EE5577C746786
5BED9A45735E6CB46057CE5221263D4837C1A6EE8C65A98B165E3A777C7618D6DCC7968080C249C16670CE244DD02B70120000FCD8F8EB74BAD1A12A
AC84BF4DFB71A8DADF02871957AE44F1C806860C33D3B22479B258309718D3C89B4A97895C278425A590C267A251A8471F891053AB804806CA270CC7
991901A7554F9207394BB3DAE32CAB7E5E41E16F7E5A1DB1F1E4F5D77BD5899B3E9F9911DE787B6D72E7F10FE71F9C28A4BACC5A6026587C29001DBC
6AB33FCD22BBF7F8C1DE0639366F4F0A014511B242C753F1599594898320FF00F063FBFF6FF1D79E1072F9C2DEC1BEF7515236121159C4C139399A1A
980A21BED1EA125B46BC85BC91B0C4122B4BD4C8AB7C897E8570E3D8E185102388A33FCC384C03E52D2D4250623C48E16916CD313CC78E5D25DC94DB
3762CA6F7E1A389BE8D0F2E6EA9AE4B4F63DB99557742F8E2ECB2132A27B26D13D4E6C700A3FF864BE90261EDFBDB9BB339F8859BA586662C44DB4B0
AB641A896B88C41983600B2000FCD0FCEB74DA919721141FA178F797C1471AA18DC6353599C70A486118A054C6911348B33441054C9703D29DA978F9
CC707E7E0027211F8A3A3EAC8640310214DC2B0C259C534766A6760A6C36CE6B4C41781FFDFAA86EE597F9E3E60C0C5EAFB50C3BFF7167D6F46D5716
85B47FF9B2D0C646427098B5349C917D719262FE9727893896207EE5B4EC6A1F63797503DA182E4FF346FAB8C0A54D2A1C7F46FFF7F93F180100E0C7
E57F74F4CB062E5528956CFFD4774AC977F1899A9026329387F320C860439B4C44062FEEA635F59055E70C047D6E27C633151E966B6235741B666C40
C2C22D139882EA2434D53B4FE5BC6B61A59971E187BE2B7121CB1323CF7D3A56E460BFFCEDD59CECBD1F3FADB3691C7E1202E32748023AA7881DDB0F
852A8E6AB7489024696CED04B7F0608FECAD7E3028465EAC30E64B1CD29735A1C8DCB6BE7F9D040E860000F831F37F9D76749B1B99A94EB3E8FD38B0
4DAA3083E2636B05D6F230B501AD1487B7A78DF7609F9A82DD165DBC974C5F7989A72E4786E7984ED54D83E0912462C2CC18B750299629B2C7D13775
5947522ABEDCC9E088EA1D42AFF68698DB4E7D70BF21AEE35ADF2F2DF68DFD171CE06E699661EB9684E71EE856BA5D186C236268F2C0AA64DFA404DE
BC8BE3A0505E690A0F2915258C977A48CDD58265A33A3DC83F00FCC0FAAFD53E4B92E1E86565B48DAFBFCCA5D94A4D4D6561282A421366802E47C25D
30D6D3A07E9BF04D95F687C598C9C7A5D65964574F930BC30A08D6C424A22810C7112319610E627C60BE62736CDBFB497C4ED7A5CED84D9576C19B4E
3D9E1951BABFEFEB914CFF6D2F560A317EFEA2A6B38B9C345D8D61136F0E9722B11CB62AC0515C3C53CD1A6F8523A74D2BA69B3ADA0B73E4C60EDECE
5E9CE5DAD1B15F4F0204F907801FD4000C2CB58F910716157337BDF8D044F6728791E96E7802D42EDE883E9D80ADB709DA81C02E5496AF33DB1741CA
AEC0463608D43CC9C875388C809114B020089CC833C2959B9D9B6477B427FF5A1032E7DA973A5E7D62FAC90B979B6B6AF77FECBB34256CC1ED132150
72614CDABEC38966898736CD3E38FC3C174FE659F285E2C8D95D7C076BABF29CCA03F170929F23352994179F40813B29F7E8B5FFBA43093E2A00F821
EDBFFE762039901B5FE144DFFCE0753146114536960442491886B721720EC7A8C033753F06DD65A3394CA9C9C2140518C4D613042625639321261076
0815C9A2CAF81A776BC9782FA76B57ED1798B3560FFC5246CB6A58FCE5457D4BE7F1B7230FE684B69C7952C314355D3AD27879AA79ECA913AD138E0F
5DF5C1D0B82C332A25E0F81657CEF4AC6DD7D23B43D0ECC064B27FBB9826C392DDFC3887F4DA31907F00F87103804EBBC682E0E0589422E1F4DE7996
8693E484701D83C81C8EC009059BEB0919EF879D4527AE4C529EE307B6B2A6541A8C5B4780136E8FC59A224CA32261669E5C1BACA95B53A94DC29D2B
66215CCD35EDFD045CCAB53397CEAF5CFBFAC3D0F335FE93F7DC6977158D3FFC6C5A69A793CDFAC72B7C2A1EF52D1220995CB640655BB0BE9168BF66
71C6511F9A0B4F36AE3C2AB4568A56A624BA36BB08AE8FE940FD07801F187FFD9B7154323A3897CC636CB87527148BF7CCE4B0951CBEA9C00A625012
09C9CA85548418D6771236086DE75193A641CCCE490CE3F5D7D1A6E289C97826958A8561A47139D2989B4FEC6DB86E7DDA5E4B4AC5995375ED475E8E
3DBF793C3E65C5DB95EE8A29ADCB2BE5F1D536454FCE47866C7B7DBD0883E34808B4949DAB174F92786DDAEBE8E06659E1631D53E869563A81A11EBF
D68F202585BCD3FD3BFF600000801F907FADFE189F2CA5A584E145A225772F7A934996712C9E8544405484E20C2A4B21E357A326D418265DE24EF7B0
ED760C9D6702EFB1323A38960231749D23C739AB2CE864CBA23873F76B472CCD7C52CF9FADE6E20B1EB73A56DEFC3C7637B939B4FAE5E38982C8C5A7
135511CBB7C4ADFEB0503DEDE7AB751624968C84B75FFAF27987AB67D1A90D96BCACE0854D2297EA62556D6F5154944043C652849BFF67FA0FF20F00
3F20FF3AED8752AB9C1044A41943446FBF73D09544C286384345B66892FB249541D834087D264EDE6CC839E8929DC3CB1BAF9A448604E1487DCF1186
1082A3A9223990CB8D09B357CF39D5604157B53C5D44478B5A9714BBB4BED6E90E5AA57B757F7C3A9E9579A0499DB0BDA32A79D5FDCE98C3CF962870
2C5B095E50F560684F914F70E7C99942F59C19997904DB863879DBE51245B4D4D06FAE3D3FF88EFE3FE51F0C0000F047577FBD56B7DF6EFCEC54EC04
F37185ACB63B3B9CA9647C5234CAD65B2833F78B94B8AC1553C3641EBB6DF06D95290D82E02E97F14168222453BBC4006648A6F15D54342BEF0C57CD
96BEE95864ECF2DE06022278D6E9E2ACDECFBAC1198A5985E59F0EB9F3CA569486CC39322375F3EBE18B332F1EF4A40A9C855C66E691D18BB9B34A02
B79CABA3F9EDEBB1B2727059D96ED578299FC44F446397CC67314BDFFECF0E45907F00F801B37FED9766B4536B003450103553D47873BD358584F649
44C9ACD95CA4A94448E988E344ABCC77CC30746E0F2E52582D53D957391B437E1A098618925D385453B8CD0417B4ECD09B5222D1BF73732835ACF248
55E0CC3B7DBA9361CE7BAA025FAD32376F9C67977B727D48DBC7B1E1D9BED3AB056CB58B14E3B3A2BF6F91DFF64EF7AD17831839D7374A65AD85AB73
49E5CF6B7014472A2C699E3942B96D58A707FD3F00FCB8FCEB2F694CA55E04589AAF5385A8F6D63A2B1A1D17514EB114A1B158638E5A58D2EA387B1C
67FF4D01AF2934365A363F92B3BA1CC2FEF496610A97A95186CEE1B6E16264C5F314083A734E838D65C7F3ED7EBC395FB54F2632C3EE35A8AFCE80DA
6EABB1AF3B5A9072686CEC8A93C5D69944ACAF5BC4D40D57EF9FAB8E7FD82DDF7FD6C96AD1CFDD14B7BDDDC979A2ACABAD5C86444DC898402493E35E
E941FE01E047CEFE87560AF05C0616591687B265565D5E25157148B175384B5714864B94F8E3E42B54055146B31F3B509664280B09C569C895059090
912D06260419CA2478928D230312B1CCCD8090BE3AD35DBAFAD946A64BCFE08D4E6B4EC9CB34F5B35A74D98164D9D265AAD281B1CF4D58D5D94A6CC2
B6C6BA9FB5A7AB36AD6BBBDBE27C783937EACCA34C52F6E56AF9F4B4FA27F368FC20A691EF247CBC27A16344FBBFF20F060000F8A3F3AFBB9F65C656
FAB94A6727C28892A24B2BA5E612AABC4A1E364788178B2CFC04ACF5E9AA3C52E52B7FE6F5F9DC499CC82A7A673864E148A8119A4C42B9D589689E34
D729F6107AD7CC1069C1F2EBE9A2FC87CF72290E612BCE94479F28E62EBF1114B9AFD461C5C8C00E1784EFEE82D493172636BD79D69A7DEB71F59C86
8CF3E5CCAEAF471DF19DA702C8D5DE84C933E576253148EBC9E8F107D4B697FED7EA3FC83F00FCE1F91F19DD17A044133332D95354508A45F2E9A532
0E1D830DB40F5EEE8D2131DD5C45A8DA46760B56F5201279E8082F83E950CC0CE618DE788E31461BB3CB23A154B5897AA9058457B5D0555A70798D9D
E396377BCC48A5B3EB67B7C64FAF8D5ABF7B7C49BB79CC93270B73C5A4E4A3BDFBAEB6B9AE78B62B6FFBD75DDE33CE2CDAE2C9EA19E9E148D61CE0D0
4AE2656DCB85F279C948456F90B8D31F37B54FA703F907801F3600E846BF2E503209E62DD9AC162B539E3CEAE84A2B068340B0740F585141A38B353E
21DCF8AD96136DC4678B91EB2EB062D4964572014C337A168941F2BB4A49303B3B724E2044B0628ADC6AC1EE7266C2CD372550B7738B4236AE0A4F5E
3365F29E5569CE56AB1FD7D8D7279BB7DCBBB328D1BEE2C62F33EB5EEA66C61CFEB4245FA13ED3378B283F345B9575ACCDFFC4DBD9CA825CACF9DE0E
5E181BE77AE7FFC61FE41F00FED8FC6B753F67523048BB766B58B1931159E9B27B9D944AC6A33962FE8CA9141A4B2C4F638B7B9D3DE2106D738C5236
A02255D4F654BCC912DD22B8A9607B3992AA1141DD3C21E8AC454EEE4B8E6928ED1F4EA80CB3CF4D531F5925F1E92D2C393855CC2CBDB946107365A5
EBEAE75BBCB959E7868EA6767EFE5855FAF2CD64A55DCCC5679578E7B54D4E3B9E2DB6DC7235C1BB28006573623A5E224311170F8DFEDFFC83110000
FEE0FC5F0B1DC746C93339E4AE1463A6BDE59E5E858885267B7398F59B6C186C1F4D913D7353A5B89E98B087EFDEC84853A337AD40733EE94361CCE6
6E0A4A99C81744A9A8515D02CF9D73CC551B6ED420F1150BAD926EB5E123B737749CF4A5A7EE3D1287ADFDBCC066CD852245C9C1816B05E3B60E3F2E
9AFFE9412A3772C9D98391F2BCD3D3E36EBDCC6436AC14297ACB881687675128816C52F273DDFF3BFE20FF00F087E67F789EDF02470C452E6057B60B
7966665B769A73C5D4C8990C94F356370CB92BC3CD0351D7C3CEA1586FF7E0E559FA3B437DE34CFDF4E7E146D448265A2A942994C2C88C0877D98299
44EF239B982681EBB270D5A72620336F2DA92B9778EC3B3D91425D3B3845B46C6F44E2B6BED74BA2173E1D7BDEBC5F7FD353DE78E472AD5BF0D18FBD
B9CF4EBB2042A6B0C85B97A1A43B661030F9C1A64E4787FFCFE2FFB7A90A180000E08F9DFF0F248774D85279E6CEFC8A5E7F1147B8F2800DDDCE22B5
58C0524CF52651E7C47293A87517315ECEC23539ECB2F0104F380A075D34966542255159CE56588D9A917D2C16376E5D2125F96C9A31BF7DBE3D67E3
4111ACE3D10CB5A5F7865F965920A447866A79CBF7A6B53E1DB855D3717944FF64E659FD556FD7B97717167B4FBF31BCCA76DA622ECCAF51405EBD83
4F9E504246E71E94F067FEA2FFBF8B7F9F9F0E83470001E00FA4D3BD7208AA1510BCADC4B8091B9DACD09CF6BDB65829836EAFE1F2FD3C71F8309563
1C3262B7A945B8694D33A132D22691488023AFE945C644183ED50766554620AE9C83B1EE9E4C4C39EC8E749F5B66E4717903D67AEB99209657CFEB23
9E48B4EAC240BE60FDEDBC9657BF6C5DB0E7CBC8E8A180D5DABD81E34F9EAD4F5BBCE1DDDD3A41DD2C193D318F0D5F70588EB055A019536E05A2124F
0D69FFF700A01B7B7409BC0A1400FED0FC8F9DE6C534D871C3AD64A4F8A9EA5C0CB7FE8023532916BAE546B93B868905AEB6E16BDC15B39982684CEC
02764DB2731E8F8172D13FA29269ECB876735C6420CFDD972E5A5023F3A90B09084FC5938AEE16C04B1FCF9779CE7D7135120F37733E727D9CD7FEFB
0D8B5E3F5ABFF5C9D7E19F2749279D9B1534FBC1FDA6F43D77EE9DF2F3D9F428CCB67BB69030FB9A9BB1D85D8A8F3E1C467259FDF6FFDEFBD7E98F2F
1F06F107FE916DF6DFF387FA9EAAF998C019C13C6BBE035E5EE0528A24C6F6CA718E4CACC7F81895D039D33ED24779B10DB928036643552CA0C7267A
06491006356375464879539702EB2646A5E611ECA617DA2B2AA754A549F9F4C55783500B9F4C35CBBFF5663203E6E463E93B27A1F1E74F7357BD3BBF
F2FA97E1AF7BBC030E0CAD293CD8BFBF21A9FBF6DBBDF6498F4FAA5C0E6E16A2A7DFF635B4706ECC6609041859CB8BFF5BFEF55AEDAAE62FE0550000
F047E5FF7BAA86B3B041AD0E2C3F3B1BB6ACD2B72E80E236554A0B92091B8378AE66760B348E13282D3B695B0F622D1516DD32F7C20C67190C794DA7
34360D5E1B8072F5C15897595B2C5D43674F5A9CEB12D8B2B9BD98EADF7B33CDAFB6B9D205AED4F891D5F3D7F60E3C5B76E4CBD1233F8FF4DF9D605E
F57AB0A7E16AFFCDF68E4D575F1C084ABAB547967D6B218FD37CD2DB20B2B2638155121F2B987A7744A7D3FF67FD5FA7D3FD52967F1B2C0002C01F13
FDEFC57F4C77DBCE5450E4484E0AE0B005F10EAD1D04BABB18A7B4736E17B3B218B49966C43AA6FB2EC2C4F3CE6C5746A32BB1A48C4883980DDD3386
19CABDD0AEEE30DF628945406B04B7A8C71E5BF3ECE434362EE8DE0697C89D79189AAB6F1C4D34F1E18DE5BF7CDCB8E6D9B58BC3DAC10521913D6FBE
AE2A3FD7B7A3B572C68DE7AB82532ED488DBAE57905D576F76C23477C5B6C3B23A0884C287C3DA5FDF49FADDAFEF26D58DAE0E9DF513C83F00FC41E9
D7EA47BF7E59C741B372826DCA8BDCC5B2E282297970928508AF894FAE11DBD6CACDA68AA939B176875CA26E3673E2F88DD1A8E21226093249BBD000
85F6B335CFE47033D4046B85A579C59C545AC9C0150514DBBAAD4C52BD23919F1960AF4A2C3CFF6979E1B5DB0D8B9E7E7DD53F7ABFD067D3E5BBB736
36DD183A5E3577E78D63B3C2E65CD86C3BFBA00DAE66FF0C99EDD2F5160DAED1DD3454D2E3D1EFA55FA7FB9FFCEB464EC7B4AC9D3F02F20F00BF3FFD
DFDFF739FA61DFB9C7F9508AC039D2356E7232175B3145E04E44522968F32C55408C3CCA031FA7448B977B372752366EA3F8301B538C12FD08C6900B
5A6F88A1B4CE82ABC488E9247ADA44F5F88D2EA4898387645075CD74867741AA9B6F5A0CD663D9EBF77BBC63EE3E8A8A79A43D7F79F0526AFEC1330D
FBBE9EBAD2BFA26CF3F31B3F8526AFFFB0C97EF52691DDCCCBF9C2F493B5E285EDF2192C64E803BD56F7AFEAAFD569B5DFBEEB1F15FA1D48CEFF02F2
0F00BFBBF6FFF6BEDF3D932F5EF7C44D88E36406A634252BA82D130D27E623A85434BD35DDBECDDBBDC07C721A967D7E4E50A372FE518E97647219D6
2F9C04B11CFC4833316672A0FE4A130A8950766BAA5B759AD277FD7A9E89665D142A619B3D213B8F4E6F7C31703F923BF9D5314DF2A59F6F3CD00EED
5DFC70BDFFFC51DDC5C39B165FF97CB726A8E0C2D38D5E8BE791FDCEEE8992776C4D57742EB7F3C1A3DCAFFF6BE1EFDFF9FFF6F5656B644CB338E58B
0EE41F00FE3FA7FF3FAFFBD66A3F34D4DC3965854E8D46C4C932D27D54C6B1330DA6DEE6E34868C4946E626330A3D5A7B01B47E9FA495A1F9BBA4C24
56C44F43E0E779422A4777418CA146107396914208717BD94177F5567616CA30B48A292E460907659A063525ECF8C8EBC996BE3BBF9ED0145D7979AB
7F78F0C8FCEDF362660D0E6C5FBDFBEE40FFEE28F5BCC767C258D30FBB161C99A1085FB95A12787806CEC21CAA38A6D3FE3BFFDA5F0D6BBF1ECBF489
5225547E06F90780DF5BFD7FCDBFEE6AECB447876C3078093DD32929DF37996CBB37B0FEAD2B828C45F81452DA7298F5694117DDF18E2723C695CB52
6325C1CEDD3CE38D519093A3F90670439A6528DEDA036B73B20A238BD3E4969AC3ED7AB709A139877DA312D9EA0D7D7DDD1EDE6B1F0FECF19AF9F4C1
D3370337164C58565D71A0FFE49C2D2FBE8CDEAA0FCBDDF7A6832199BBD83C726B9055E3DE4AF3B2C3B5E2201E92B7F93F677EFC9AFFD191A1BE63A9
9681D12B0EA47FD582FC03C0EFEFFEB5A3A3335C173C5C2FC619C922348EF153FC3BC4462507B30F0799E2D128929C5DD086F59B80DD9B6142D89C87
9D66E5D44909E0B4B94002298CFE01BE214C96E92B8C8F8250F735C25889FEDE13E8B8E0A5B5045CCD113F5F0F7CD5A3B16B01968BAF0DBD5B376ED5
E7B797BF3C6E9B7CE0C6F9F9576E6C9A717570EC55BBC27BC7A397A5088F43472D0445DE419B970509BB2E461BA516C360F3BE6FFAD3E97F5BFDFF96
7EEDD0C8979DA16A6B69F3D709F1DFF20F060000F8DDF55FABFF39C0B7E749378F8A34CFCD9930AED2B9C1D5D8620F7B62279688C71378D693B7133D
DB88130270E2F603AE8B3335AB559E9CC66C181AD23C7ADA048A8E773715A6A1188B728DC91AD7142F2C63D1C30823F4B49D2A19D3A6573BBAD5C2F3
F0970FAF9BC6FDF4E5E1EB2FC7A62EB8397C79EFDDC3B3F7BD1CD31E716339EEB9B66392BAEEFD5D474E72EBE42585B6A55B8F859AA25DF0F089FDFF
79ECE7FB30353232F271838F63A43AE7E93169F330C83F00FCCEFC7FBFB5363CB6831E74E04125010F37716EA99245C9828B10880E9CE309128288A4
5812AAD623147596856A3AD9E552927F807379AA35352FCFD808F1405B6604834BA1A62A2ABF2DDD90244F5259A3BD0F3FCF84D12B9A145242E8D3B1
5FAA69E36E3E7A772929EAFEC0D7A7B78EED3BF9E5C5C1956797763F19D50D4C57271C3B3F3F6ADBC17DCFCEC4D0C7F75647FA95ADBAB2C4CA54E64B
35CEFDFA9FA77EBFDFF6D78EE87E596C1B142FCE7DD193EAB8FED77F041F2300FCAEFE5FAF1BD54DC2C65DB89F8D25A2D1E4940855946FF6065BFA1A
35EF82279280086FE7B62CA45A67B8ED2A259A39373AD1CCCC33BBD5BCE21ABCA9D7E8B0D2044FA2B364916C77770839705C8A885770E8901B8CD3B8
C5874EC8EB1B5AA82296DC7AF7F3EE98F6BBC38FEE5CB87FF1C9BBA7978F1E3DB7FBF5F0C0A98ACCA3C3436D9A6DEF57E6E746BB37AE2AE345FC74BF
2B8642F1ED4CA4567CFA5F5BFE75DF1FFA1DBDDFEE33AEDE26F6D336CBAEC40B7ADD18C83F00FCCEFAAFD3EA2EDA6063AFDE0D3772F1353596065744
4B926655909AB34CE7D799A08D931E0786D7D0C551B48E7D5062420DC3CA07E5D12963C52DC341266B2F1B1BC331084A14874EC5885DD3DDD1CCFD5B
94566441D3540D5C386BE852105ADEDDF7E1F0B2F9DBDF6DAB9B7EE9DA8BA7EF3F5EB9F8FACAADE74FCECFCB3EF0457F222EE1C2CFE502ABF4F6C973
35F88C676723F058C79A420F5CD907FD7FB6FCFF3AFF1F39556816B856E9F7F1027B5676D44B70040800FCFEF9FFB798B5D0F0E36FDDF585664EE3A1
6D13EA4A2DBD579E089BB81059BB0B8983E1E7CEB24A53F0132D5CB628A98587DD35994C4DBB0BB7F012D7F888AED31085842A93FC8C691E4A4D783A
8DF8D3761692659FD5CA34CBDEFF722E8D3BE1D22FB73A1BF60C3C6B51E61EFBFCFC44EB86D3673EF79F9E3D79D1CEB3CF472ECCF46E79712D5C983C
B3BAC82BC4B6FED52A2E9CEA9C15496034BC19D5FEA7F6EBB4A3BAB7EBFC5C6C2AFDDD9E3E52CE5ECFEC1A00670001C0EFCFBF4ED7E787C74EB879C3
D5C07D8AC2D823C8AAD3095D7BCDCB7603D5FC18DED8D890DA294A4984B9F919E44E32915D4FC587DBF00B839089FB499897A33EC63034213A0345F7
F53515A79AC04EAD8713E5ECA03084DBFEEBD37D65D9679FDFAF882FB9D4FFBE2967F9AB972FCE2D6DDBBFEDDCA7735D15B3EEBC79D7D39ED179F0D5
7A77FBF28E7896ED8449E5CB82A15867911D17279DFA48FBDB797FBFAEFEEB46B5838F6BF911A1299AA0E78F9D665E930AB76B41FE01E0F7977F9DFE
0C9F4DA9BE7DD5C39096EF87D184F256A5A02DD6F84A77BB116E6A0CE0A6C695A1F9E54C677F8C3C8D265A502E94091C67C431CCFDA001BA674428A1
28DB856E51680F77AA1761D7F6D2E17C4F6B6B6AC2A37B81E8A07D6F5F16BB7BF73ED55D2F6FBAF3E4EDA19F369C7931F2F0C5E5ED171FF72CDFD532
7EFBE5F77B33A511950DDE9E05EBA7D99038CE8B6EB43070C4A0D9B78647B5FFCAFEB7F4EB46DEAFF5662D5F2C096AFEF82678E17D676AF083DFEE0C
80CF11007ED7F47F74B8CC56C39B78E3848311AA6E320A13249D30854B68988C6EAA33DCBBCAC008669412363E5EC6F49049DB9CE13EF94A575BBBE6
18261165326F6C3DC4C8AC8487E3458461CC67AA39A95504348B694BE4CFFEB4DB4F1977F6F3AE00FBEE0BA35F4E164D79F061FFCE6DC7B6DF7D747B
D7E1F31F3F1C89AD3C73E9F8B10D4592D892B94D15C92955130884EA87CFAFD430519C98A5D7074646B5BA7FE55F3FFA76670A9E5CB4DA72DCB691D7
A95BDF7962848D5FFE7567107C9000F0FFB9FC7F2FAEF7BC942A56E3ED9FCCE188840D6AEAC478CF434996A9FBD9391791F503E68670CE74EF8A4C27
7E88C4F351B349F624EED4C996B3E6CBA464F6F3B13C4334D20429292AE54A2B7C6DBB1550B2C496C6715B757B85327BD2DE770B1CA28FF4F5BD9CEA
31FFC38BDB278EDEFDF8CBE0BDAD47DEBFFBB83EB3F7E52F6F5A2CCDE5E91B274F1B2F902637FBFB1C1F7B178241C019A5DBEFF58F8CE8B4DAEF4B93
3AFDF0DB5D9944B3A29561E48A1B834FABF67D4920093D8F0E83438001E077E5FFD73D35A37DCD5409875C7F6D2313664AAACDE315C4D257E66AAC77
5AF3AE84F2DF1D373270396B963A31C4C85D4C3C728919340933B11B5BB6C99F061FAFFB6A610087E3ACE3AA2DE15EE3A0A129A65085335E59B77D89
576CCBB1FA73D58A8A275FDFADF151B6DD797D60D7EEB3CF065EED5E73EB8BF67549C2E6475FCEE5F269B10D530222636921DB9FCE4FB93A76580A83
9A48DACFFC3C3434A4FD3E027CFB19FB9FEFC9E1E3FD2E7E9D65D1F972E05AD9AEC15A3C5D54766B480FDE040E00BF37FF3AED5D07270F02BDEEE64A
2A1E890AEF32B70EE6F7CEF44F2C69A06DEE416D180C3240AD77F39E9266ECE36532ED5D84DB5617AF2894C7EE6E16E3C4D87E049EE1E21B57CC3755
245A9162CCB17C2F6EE8E1E351A2EA9BA7DAA2E31DD63FFD79759665C6A617B70EF41EB9F975F8F1D6432FC7F4A7E2D50BBF9CDF924C919634643A55
2C4B2E3BB7285239E7E70E3212810F5E7DFBE3C0E0B0767874647878F0D585A5D9522631F5E5E7D298A3436F8F361DFA3C8F6B2609DC780FE41F007E
4FFC75BF2DADBD9C2D4AE51288552767926068A2799D9DA89157BC8735CD76352EF699B9F9C75E0383E63A42692A5C180D31BBE72B58E8E31C4E474E
EDC5A806C7F20D713E05EEF63E68663CDB2472360FAD648AE7AC7555CCFFE56CB8854BC7BE9BD3D3432B16DDBB757EEBD6E3AF87876FEFBAA81DFBD2
EE1171F1FAAE521B5E584CB1C061EBB5CCE0F9A1AAC5A7EFA6119874B3FA03CF3E0FF4F70F0E0DF4FFF2F2DAFA222F228D452BEB7F9253FEFFB0F75E
41516D5DA37615084DD33439E79C738E92244A928C8A4431A262423181604E185154CCA80898504154040405442527C9199AA673AFD8EB6FF67BBEEF
9CFFE67DB7FBAB73736A3D0D377D41CDAAC533C71873CD39C704ABEDC1F54FEC37C6AAF6B61955DF59FF55FFE3FAE3E0FCA3E88F2C9DFC79EE149124
27AB7EF4FD0115B28AADC3D678935396BBE60DAE7BEC8E349F3C267897662798D96C90B8C556619B04F16E8CCB0B9F95171DC5569F96C8C198262429
373BF1E53A061E46922BAED8109D931D13F2B6EEBBD874C7D3F448C5C4A588E0A38FFAEA9A6A6A3FF53390D153DB5B30F68D35E91DB36D99563A7EE1
7B1D64E26ADEB9686E5A75BA7FFCBEAF8295B6F5C9AFE31C2693CD65B366BAAA0B579B90E4C3B73B1EA04CE55C62B32BDEFC9EA0D4469938BB27BC7E
53CDFEEF9301F8D3C4C1F967FE4360EF4E9D2DD65232AABB5F1F56975074B6F13DA81EEA14C8F44CCAD039255EDA2CE8C0DC2468389A629192410EB5
1258B54ABD7C8D518197B49AB9580D564B10931413F571140EF19532DA604FD05BAB115B7B27F2F297F5CA118FE9AFA30CF63F1BFAF0A4B8BEB27990
CBAD5A7F6218F91CB5AE1DA51FB7D7F64D381DADA918539BABE5F6E8D9D7E993862A86F6A61E771AA7990C16C0614DB695E606A88890D6D414C5DD64
759DAD407B5FFEE68D8EF7466B7B389A9E6ABB53C0FAAFB30178028083F3CFAA7F189C3D6161B85B5B524D35E9E5491D594D5DEBE0C3CB035778401B
0C9EAAE6D8EC5D70177F6C2D285C5B691ABB3F6879BCD0B61D320DE7D44F271A2893D4A9D871612979113D0752429A9692A1B2B0C9C6C8C827F1AEC5
5FB66AADE9FE99A1EC76B5B3F6F5858A91CEAE05EEEC891DD520E358E053043A662469197B3DCFD3387ACF115BD3E3FD8FB26F381255E2628CFC4B7A
E798741693D6F7F678B2A381B8E2BA57375785562FFCBCF791D5F60DC0BE3DFC1067E41062E7FE7D3CFF0EF85FF11FF71F07E79FF90F725AB6D9CA06
A89165C9510FB3D525CC95B542D202635D14864E0A17391A6EF2FE72473C8C2028183F19EA1F7E784D985CD201C2EDE3E4CD7BE5D4447C30CC4B584C
54C39E1096212D6B2625699565A4979FA87FB33D4EE5E664998D52FAB777F70E5FFB3ED1D831C7EDD876729EDBBA66EB2FACD343D279C393FB49AA2E
05EFD728EEEAFA1E27B232DDC132C05B2FF0D5C0DCECDC507B594EB0956980BDFC86B1565F97A7A33335773E2C54FFC2D0B2B3A3396405074FE51D13
9F0F7C42FFFB70003E01E0E0FCA3FC1FA21C5FE72EEAA8236124EB5F74D050D1D6D2D4E3D476EFF50ACFBE8B6FD82A14B23EA7DB7DAF2841486BF09C
6F8C7F81A59A4391B2679266FA256345E2598CAA202EADEC62E2B3538964EEAAE893E261736EF3869B3B8DBD3E23A7650D0BDA6A524F7F1E99E8FB3E
C77DBBAE6261E192C7192AF25ECFF8ECEFE79BB53C4F7E2AF40E289E7A60AC77F8D309171509E39D55C3FDED958F77F8CACB2AAFD96A207D86576F70
94C269CE3FFD9B3D3A81B11BEFFCB863A06A61E9B0BA9FF9E8C2F0FFD90F047F9E38387FEA3F8FEFFFEFD596CAA2E67696C6A2F6B93BAD15F403DCF5
E3F2BCD3C90133CE2A7B098A27632AF664CB091205EFD47BA5AA5C0C509448F250F2760B396145227EC13E1025656C34DC8E6B132D4C35DCC2356C5E
C4D9DC8C93F61D41B2C56DCB670FEE286BA3CD55B431D80549AD8CEE74DB32805B6A12DBF02446CE24A9AAC07FEDABDE57E1B1679BCA03A52415575D
FA56579EBD5A5F914C925859F22DD6F43ADA1EF284FBF3C5E5CA5F6CDA046FB66D6E723689ACA925BFBC0A184ABD39FFFF6B078C3F501C9C3FF71F7C
E561A4A616EDE2A92A6AB76F8F97BA59A88B8A79EAB6F5165A3D0F954F5A10779CDBF9E67288A0BC48D4AF987D9619D71C55BD83F4E3423CAEFA1334
18D84EA2B49EA34D96B5A8818782F52A1B87BBDBFC2E59C99F5B985D4FF6A9AB3F5738F6EDC6AFAE612EFC725517B56EDFD12656EB2E8590AC0009BD
3D553757863D9BEE4D0D7CB6F8799FB6A4D3FAF47347B6BB9A8991C434E4F573BA3F1EB931C8E87F55DFFDF6C5AF292E75BA67B4778CF3A3F0A99B91
86D6CA6B546E4E56151DE6E1FEE3E0FC8FFC5F8897097375D965A92725A917B9DB47DBD8C956CBD33E7C9B8BE0154ACCDA28419DAF09D70B0A342489
9AAF369E72B4FE106B99932419ED217F2E665908069B8ACA196BA6078AEA58A818BA8B2A173C34750F5EF98E556623B6B9EDC3BA22C6F75DA5D54330
AF757D37E3EEA527ACE13D9AF231F96E7A714F4A124D4F32A60F845C612E9C52221A669CDA136CA728252C6D79F24371E28DB60F05DFC707EB3F55FE
689DE0808BB3036DDD3323ECF78E966E4ACA723BB2BF400FA3AF542E803C7C02C0C1F91FF9DF62AF16ADEE15AFEDA423AEBA32C5C1D8CAD2D7C8CFCE
FCCC2691339CA6941431C9C643BEE55F4F4B68686E8CBB156F539EA679FDBABEE32AF5833B050BB12E92A4AEF68A2459650B1D1D270BA5839774944C
8F8F56C589C89E9BBEBAF9DCD8CBED3D5D0F66B0F1AD5F994517BA581F4D96E9EE3E60A175EACD11E3D8F6D95C8BC354E48B0D41544E935F3DC8A86B
38865E6919B9BAB372BCA1BBEDC797FA96BEB6790C65D0C7E699A3C3ECFE6D86EA0A8A5256DB4F6DA031D3CE9CAFA67271FF7170FE47EB7F8B3BE46C
22442C527582A28DA40D0D14AC3DB543CCB5820C927245962F8C876FD510D8932FFFF8797BB08A9ABA6356A6756AA4946582B5A1BB8CC37AE521ECBA
908CA47EAAAE949996A2B58CCA992253F180B2736BA59759D57EDBBCF9FBE3A49323F5DBC6B0DEF535DCBCF3C3DC6B1284D07DF1F26BB3CEBB98DFE6
9699BB34638CADA242622449790925B5E033B73B7B3AEF9CB835421F1E98FFF4B8F9F70C82A04CEAF0C8DC7CD7E4E2F315D25246AA5229ED35A64FD0
CAFB77722BA739FF674370FC15000ECE1FFB5FAF430EDB20EB7CD4392851575C414CD62AC2D5CB5927C456797F88725BBB7E762CD1F5957DCED7B144
B2869A79CA637F1B1B051919AB581B7925B578085D2BAC24171E26AF6F6560AC6D5C725B476ACBD7635202C4F575E51EA75A9EADCAE96FDED4C67B11
D7317BBC98D29E41545E9F6BEC73F95DBCC9AEEEDF5B5C9FA1BC3BFA22A2622482085161DD97F1B9CEDEF24BB50C567FC302A3EC71FB3C8A618CB28A
A171CE48672FF55582BCBCB6AAF84626762C84C32DB977FAFC8B19DC7F1C9CFF81FF28029E2199E6EEB1774D355AE12A42961515578E0CF77070B035
978848163830EF61972DA45E1161D636BC936814A0E47ED76955A28AB4986DAAA9A882C03998AEB34CC8365E59D356D1504BFE559D05296222478C64
9BFB252AE871534644C5EC29EFB2C55B870786B2AA7E3FB056F7C98DD75FF3E4A4FBFAC6D14BAB6E03D8771F228148241044CD923FCDD35FBDADBB56
41999BF8D4B4B8D058D1C6E062D8ECDB773D8CD19F3D933DC792D71948494A2630C096A85A46C3AD3DA99985C30B285E00E0E0FCE309004527C288A6
7B0FE8EAA55AFBAC101513279164D7C76E586363EAAA1F795E6727A758EF81BE79E54DD30FF4323DB783C61EE551B1799E2EEA1E999BE564494D501B
594431D9996C61AD67A072B2C55ADCE367B684684CD755F78CEEDB7127DA1A0FAF79C7ADBC325FFBB0B238CC30E0E85ECF90F38702D735509F399FA0
A1DDA962041111024139B56C86C1A8BDDD325C5531CD78FD6E8ADDFBA39F4A0351E6AFAFC3D85C4F47F7AFE6DD0AFEF95B544C8F4F433359373873A5
9F8BB71DBD353203E0FAE3E0FC53FF51147AAE4A261AEF94D5CCF6D1F02293C488A2A4AD36BB8F6AA807AB6A3E596E39356A981343DED3EE953A31E8
A59F161FFE30C97287E9160B855D87A4058C69E0150199687782BC96A19A886F7B88B84ED5053199FD9F8E591CEB3B1FFAEEFB23AF6363F44785B467
B91DA70CE2F61D5AAD9550BC552B8B351EEED7CDB9EA40261048C2C2AABB0730CE42CBD30F939FCE7FA0763C68E78C3735CDB00060ACABB913A54D8F
CEF5B45FF0B3DA92B1DA7D471D481FBB7895CA297D011D3AB4E3E1D8D0040FDF018083F38FFD9FD9242246D0DD6FE8F830452B458B44269154AE3AC7
97B9AB853819D75CD4EAE5EC385620B189B2C5EA279CA11E9EEB999F6D7468CDDA388BDD794A0279003750C0224C5CCACAD2524D756F8484DAC54219
B9E3ADA12AF9FDE9415F66762795CDB727DBDE3F927225C824EDDE4E83E813D73C23DFCD3F71390B753A0B098B880B0A2B640C61C0686565F3F8424B
C6959196273D9CDFAFBB39009BDAFCBD630E624ED32803F7939D36ECC9770B2F67CC50A66E5D5B00DFE58E4CBE7B9B5E3ADFDAC04370FF7170FEA1FF
50998994FC32A5ADF616D70F2AA69A89F2E3BFD23E27B9E77BA4F5A3925E5488ED625F0E6DB53669382A7780FD44D232DD70E35585959996E9B6211B
05095DC0005180244314D2329390DCB486A4B8A7505DEED46B5BAD973F4313DA9AD3D37F8C5D704C3B6217DB9AAF107C29DAEBF4B73C8F5C7ACB6AAF
16F48AB8901899A0E67ABD07667694BC6F1899A1D49DAC987B5731C7ADAFA00000E757ED107D9032B3C0415A822C93B69DF50A2AFE39D0D9DE73FD31
85579D5B3259FD33EBE8AB89C2872884AFFFE1E0FCC3EA9FBA4F396D1359F76C88C78B6BAA5BBCF8F15F425867D5B223CFF5C5571CD9F37EA5C17C9D
E574817AC9671B9F598A8D61C2F29817E6565156B13EEB02049CD9C063013109A2B88ABDB6B8F76D73F1C8AB361269AF6D2CBF4E45780DF56CCAFDDD
BDD2E15865F4AEB98558E7ADE61EEF1E797BBC62DC32CFE5422904A284B080F7270AFB6D595155E7C214ADF149E977FAC077D6DCEB374C08EC7BD5C9
9DFFDCCB6222DCFB76B66B2FC4BB1CEBABBF5AD1F6EDE15B1A5076AAFDE7C8CFBCB08A4F4D6B5EF1FDC7F0F88F83F3E7FAF3506C204CC1CB4D523FC7
DBE9C34552B0235154842C645EB84CBF5857506983D3F36BC47BA3B2AF0664235B9CE45E628724C2432D8BED651CC235ECD3D505D2B85C4F01929080
46928BA0F70D4FD2FA2C03F5DC431A56E3CC78EFD1CFC197817EAFD0E60BE1456CCA510B2BAB7527E295325BEAA2ADDE63C3EE22245161C2DA396E75
DC9AC6F6F9D9AECAB2E7DD93F36F9E8E0F977D0610E0D3BD29E6CF1A0A7BAAB16075C8A687EBCD76748DDEDE5CF1E3D9CB8139A432F317AFA569F1CE
B6FB8DD7836B50043F028483F3CFFCE755DB89090B11AD4E04ABDFBE464E4C1493248B0BAF682093762A2C134D8BCAF9EEB68AEB99BCB8DCFA7BB272
3AD2AB1B916AFBD855D731C52A245D5AE029F89B2448D0F68B5CA1AE9AB6811454E8A470ABC150FEDDF45AC7BE49FFA3D0FC2ABF8132979D4F9A32A5
74E212CE5B39968CECB3CD9BC19A0D4544887A59CF877E1C8BB93749FDD9FEA2A4AE7761EEC5BD8A9EF1A159149CAE6963CED78DD0FB9F6418ABF97D
FC10A07619E61E3A3CF0E3C13013024AD3FBC0D9113AFA25A1F0FBE62DCD3CDC7F1C9C7FA63F0A9F5753921697524E8F20A43F155B775551CA587599
468DA5608887207947866DE359C5C15312F527C48FE6CA1976C3DB3D6375EE192A78BB9867AE1210EA034B0404C8DB932549E272C66A6BF7E92A15D6
38AA3F7B6661D7D3E9B61301D27C06AE1BDCD920AD2369169513AB99D25F9316FE05C16E4A11C51CAEF5231D0777DDEDA4D515DC7CFB7166AEE3DAAA
CC761003306CA6F1E730343E4E1FB89CE8EC1C73B570A3B267033671F4E9FCE010079C660E6CA84247C62627B9D7BCEFBDF1DFFD0DFBDFFEE3F53F0E
CEDF959FEF0BCA636C51F6F1D056D6DCBC4769CF4F95C86CB2848181884C45CA32CFC344EDE86B5E8F7B4C2BDA25B77638061DB3D5BE8F7E708C953E
E9AABCD15EF5A09F8037001D141095521310955749F5D7D9ED247DE283A7E28D5706DA4D134EA17424C3AFEBA1EAA60245B2A1E9C62807CD7D1345AB
1A685D397BC509BA676630F6F5E83B7D8CAE0B3B8A3B66673B4F78B8DF6560C8C4F070F9EB29880352876E477A7B27BDFFB25DCFE73E847C2B6C64D2
27406081C538FF149C1D40E738B4F56B9E9D5975F63B86E0EB7F38387FAEFF52F93F16A8E0A1A3A42AED952BE73D187CFC961889AC445E96795A507C
CB72B2F4FD6DEB46FDC2E69C94EBD68B85A6CB05D0879D1305A39C8522CC457CAC05CEC20C5701828418D1D271759681AE81F2A1125BADE2B786D6B5
5F4C6D67B9716EED9714FC6E192BAECE58E1AAE2F4687677E638346A49241113FA30AC3326BD91D9FB3235A799355F9BE5A59BD88A21EF0EDFA9785D
45C37870CF057F9755EBCA7FDF36F77A0E63C3B72B46D199619039C5059FE70D8F0F334629D855BB2D37A3726ED4E1FEE3E0FC13FF5114C53AEC4D62
D574CCB5DD8E5959BFB6DB39E32E6EACA3B82CA2545ACCEBB0A254FE7DA7DF578DA7AF293EDA2DEA74C1D6E41367F506EDF86D120E9EEA5AAA222D68
A732515A425F582FDD575CCAD33EFD8491FEA97237B5FBAD96363D8C58FB978765ED8B6F39F8A4FB29DAEEE99D4838C605AACD848564F2208CB6CF2F
9F3DFF31F75CC5D45049B0534052211B2C0F88A95A9C82A0810F35E703ADA277BE6B6ECB087F02A0B3374EFD0016E9F320C2E270EB77F54314009906
DE79869FCA4C7BB7BF92875F008683F387FAA3FFEBDAFF6A35835029494F33F9D410A55DFA4ADF1F92750DF4088AF7CC85C857D70BDA566B9F6E92BE
D9A5E290A3221D1F2CB18593E9171F5A2022957E536F99EC387A5C505C7A9994899DBB1431D2DD2CC754EDC1551D833337CC2D3AE004878F6BA4826E
1C37D4F550D0DDDC38F7D1FE2C405D2F2EB4CCA50A463F44847FA5B5E615D40F35E4AFF54A2A2F793FF333636D31178331EE8B7D9772E2CE5EAF2EBE
79F579138BDE5974BE6C82C99A63B0787418AACC68989B0790C1A98900F3ED87039FDC71ADFAEF0600B8FE38387FCF7F0C86FFF21FBAAFE27F4A4BFF
5CAAE68502E3834E6297A657C83BBA6912F6A711C83BDB25355A8283A8BE6BB81B95F7EFD2B35A6313B6F0C03143799791A4E76F5F010F108922C891
D5FD027D9424921E99185858DF7AE8201F73D95CF733726A5D4584ACCFADA3C6A13B9C022B667BD3425E336EEA0B0B086C9E017EEC4B291EEEBB76EB
2B77E8DEE1BB8515DD2FCE55B797D78118C2C6BA8E9E6D1B6CAA6E2ABB935934050E5F3BF2A0618ACA5DE0203C0043C0CEA23A0E93010EBEA49CB20C
C85EB3B762B96B07DE000407E70FFD47E7C1A57BFFC1C52312098F14F49F6F352D6BB5DC11246AD39E2B27EF61256C962F4154E94A5EB6F79A42F325
C9A6A72493439E76C7424C3ED71906113CD2C48885BE02D908C58A20AC1CEDABACB32CE183AF8CB1EEE50A53DD33EB14E54AC182B4EA95A4808295AA
BB5A0E24BD6AD81E7078BEDA5A88B88C94C7993A1A7EA367F64CCA1B046C3C5E39F1F956E3D9FC3EF624CC036918F074EBCDA189E72DEFEF5FA89E82
C62E245E180419DC293A8872418437FDEDE31C63860A0C5C1AF9616B199FBCFCC37955AFA9BFFCC7A33F0ECEDF0FFFDC7110596AFA37B941ED6883AD
E3CB209B9A2FAA8979AA2A459B249474ACB46577AF14112BFA2515CBF6DFDFA79AD0E126E76469762149BD64D4C35CD3EC882EC9465FB69BD7AF2E29
EF61276AADE7F12A5950CB37FAAABF76CE156595B2A9E32125E192769BE2DD3377AE483C67AF1C35C83C4A16240828954CDC8BCDAA1A2DCE2CF835FA
FEC29B91996B07EBAB1F4D406C0A4283783D47326FD7DE78D1F9BEA0650146678BCAA641109A986480100AA0D8CF12FE44B0C8A0D33A7AA6B73AF944
5B66B04A2C77B3F0F7FE38387F048A8D7622FC0A004206A265AD4A371AED3054385D241171DE4CFAFA6D290D457307A2DE311921DB491FB9B1EBCAB5
FE522D35AA6A563AB109B2D17DB1DA01529B3D09CA443B16F29C20E9612CEC6F20773D4E4829687B8CBEF5C50DE24E8D3F7C8C6E868898AFF20A0D0F
5B5970D55DE10A6DC45D8424289ADADA90B2F1615F75DEF1BAC9C6DC0723F0CF2DD9A3DCBA1160680646504EF9CE938D65D9391525357DFC312EF453
61089A1B1E99E5C25C00C1FAAB66785CFAD4C208853D9AEAE8E7A99D585171D0F02984E2FAE3E0FCD1EA1FEFDB6F0441F8E57F93AF0A39AECCD9CEDC
E0C01583B8DBC94A3747B4C5AD8C9C8DA5E25364543AEE8AE5D528ECDB47BE8965CA3B9B1F38AB6CF2FE904184DC86443593656910942EA6AA2A6EB5
5A362658582133C24E433FEB987CC4C0574BD2F1635AF2BE190976A1A52D3B55B636CDA6AB9008A4A86FF31702F3BF553FC82F9F1E2A7A358571EFAF
38BF08F5775187671098FDF1E4D90FB5653BB3EA86A730D6740F75A41F86688B9D8D6300040330B6F8998D21DD6F46D8B3F4E164DD95DEA661433BBC
021D7FA0B8FF38387F16FEA7DE3321184160E099B19CB8CC6503E32C877517DCED8F670AAF9E4C2358AA18B969489F765D76A681A8D310A67953CA6E
BCC75C5FE560BB91D4BE425D57ADB4686D2D811BE09C0E515158D3462C7297986962A22AD126274862C342ADBE54EC1E2D4953AF4015FF5FE321B2D7
D11E1792B498492958EBB5EA4D4DC9E1DBFD7315B7FB305E47DAF2C3438CC91108A4020BF777E755B57C39B731BB03E621532DCFCAE7FA8BEA6629FD
5D933488C38178D46723185C77F12B6D727176836268ACAA63DF0BD34F896EBF517CED0F07E7CF56FF9EBE0320148511F096B6A2BED41177D7EB317B
6F6F30BFBC49C86CB24622C84CCDDB5926E438396DDE4CAEE98BD6937089625A80BE71C6FC3AF29A0736F626391B4DB4E4BBA01E093555AD1443B524
6DB32C370575B7D3D9AA6BC7073C05BDF65B8AD86F75950FA8EF0FD27D335866471227FB2F2CEC504AFBD45AF0A0A5FBC3B5D245DED409B7554F061A
DA4018A17C7A91B5ED695DE38392372DF3088C0C74431C00A17CFFCD60D25910C0E142D07C7123BA5853D00CD1E638AF0C03B2961B7EEE367B31A4B7
7E165FFCC3C1F9A3F08F76EC68406008414056AE9C8336C12F4679CBEE8DC74F2A9D5929285CDDBB2C7AB584DD720389630ECAFD3B458F30BC3D2E88
BBFE5C2567E6D5719FEC722A504B252F5783E04A878A05C82A6BFC54137DD5F4ADAD6D121E6D25C58E4E878BACDCAC2FEB1BE927E5F2B5D3DFB3EE83
29495C5B7727BBDB42F9ECC4FDFCBA9EE7179E0F01CC1B5E9EB93F19952FE778E8D8F6BD4FABABAB6F3F6F9D654074F6E0AB6E08E101143A17450018
02B81008F5DD7ECF9B7D79A787C9A1B0DAFC0C82F4575C694F28A6AD5138C7C497FF7070FEC87F7671D63CCA2FFF5188B14B3BDB8C1850A8117763FB
95C726A7C2A5C4B23B97B9DD96B70A0A538D2AD4AFFC22B5827D51EE4EACC289BDAAEE562F5AB58C731354A58E3D51143B0E825B052575EC24961FD2
970C3913E6F6F0AAAC7B0F75F732C7C240996D8FFDC5023F15DB6E9CD923262627450A63BCD7F2FCD87CB6A8A9F2F6FD71CE40618CF3F6928EA6DA06
160F6D4B3ED659FFF666D5100B8098ECD9BA0B5F6104600E0FB11108066108E2C0AC6F65BFD943D595DD18678EC9BA9494B13CF8C18D3B9FA0932A76
15F05FD93F3E01E0E0FCDDEC7F6C4D1600C33C1406E7D6999DD410F6FCE97AE85CD4D5B73E29ABA4F41C2E8B2ABC72500BDA682DF3D0CD793C9CF4B9
851C72453CF212D9C466CB0F378D6D5BCD09F1CF89B2DF208A9A808AAC68408A16C9A361B7ACD5160BE397A3BB44569DF714DEFC3159DAF7D579D115
BD29C2E2326461DD77DF576EEE7D12FDE5EB8D073DE058D6C69C8CCCD2B29685F1692EBB61EFA9AF1577BF4CF2E7238031363E3134CE0FFADC89712E
04B061985F087087DE7C9F46C73EBE9B00D91410B87969784F60E591FC09F083AA914B2B1EFD7170FE28FCC325FE0FF8E11F81516432CE31DB5E58F3
B4526641DCB99F59A109E44DDEBA066A65BBE4576C8B55D97E55E5F64752F2A4BDF27D4F8F37864AAB56966F335A7DCC9B68738660CF863B96116DD5
1D6EFB895A946610653C65E4AF5D32534DACDF21E357B84ED6EBF67AE1D8AF8E021262E4F0AC825D5AE97DC7770E565DF909C36501D73F2605DF1FE2
F2A62ABB7EBFBF7CF57D69CB02DF611EBDB9E637950131594C2E83C185610484607EEE3FF7F22B44EFABFF398D21B36CE47566DF6631F7CB9769C094
B9A273EC386E3F0ECE9F847F84BA2FBE1A86501805C13E17B1D8AD24E940A26563F8FAEEE716AB979D2A12B5933B50266FB52A4943F37398E6CF6DE4
FA1481F0D3D2796182C11E3B0AB51CCEC588A97A096F80917A1145924E76B8BA75DC6AA264C25AE115D7B484333A3648C5ED72D5DE732040744F7B98
008140383257E2A47D73F2DC8E9A9BC5348CB9CDF343A9776A171B82DACB3B267EDDBB54D2BA88A2188F565FFA6D9AC16574F5F1433E1D00412E0441
0C1A973D3BC361F5D67F6723AC05047BBBABFF0CC9FEDEDB510AB051445F2B878EBFFDC3C1F9A3F4BF3134A50D4178280F865ACC898125CE12D7D564
DAF7ADA54C3884904FF66ABA58C6D4DB9925A7FB997E7B23BBFE97FDC92A219D47011B7289FE9B57173BD91C3F2463A04C7C8DF11E491B1AAFDEA698
7C405B4422384ACE619FBE7A6EED6943AF4C6BA76D17F534AE0D860A9049C27B0772A5967F1B48DA597AAD05C39AADAC3F1CB3B8460738DD679F4C2C
BCBD52F0760EE12723C08F17DFA9308F034CBC19870080CBFF403C68B17F8A0EB081F9C69A611A3F3FE0600DF9D54519977F4DBE9E026A08628EC6EF
0114D71F07E7EFEB8F3177EBAF19583AFE83C0D05B0D11CF6B2EC22F12091B4EC54D2E78794AC48D3BCA46D8176F54763EECAFB87DD492FCE2996B93
A140E49DB02C05D57D3E7B3619475E57912549CE60808710C1355A6CD5293D09F1D01D229A0FDC448A862395B6A728B95FD92AAFFEF18DF1320291EC
39BC572CA2EB81CFB55F0F1730285FC1FF5B8AEE5B80B5F0DCFD18ADFDE6E5DB6F17F8693E34DED4CEC4307EC807A7860080CDE5B201900373E7A6D9
1C9035F4BAA60F8438B35CF879CEC8894D5DDCA187D7E728F604630DE741080FFF38387F12FE3FFA04EC9EF8EBF42F8294AB8A471FF417297B45948B
3B34C58A725B6535BB573CDCEFC9157B8BC848B5751361526BE7936E1D24183CB4B3095558EF1178C1CBE6B6ABB4A8130F9B141711313037DA692F2C
65B5DF51668FA7E28DF295F23667BCB5776C91D27D785961998898FDDD921485A48F999E6F4BD7FCC2A61364F70DAED5790FD0398CCC1B8C4F6B124A
3B28C822421D1A5844308887C05C900F00811C2E07E1B2163930CA197F57DE4465814C1A85597AA6F660E6E8CC6455F32CA3DDCC3942EA240DC1D7FF
7070FEBEFE3CF64103A703B3D0D2AD9920B34056286C6D9CF086765191CD070B5947544E8BDEAD201A98BF2AB177D6396CAE71274E55EAC33DFBE76A
84F509E453FAB67A16FB234959B144919D3CAC544858495932C25348C2325E59E8D225A1BD35BAA48DF1063AE136C2062F5308221202A2B73BBC350E
36AD096FC80BE9E7B5D81A963252EDEA0016BC78F0EA60C18123FC347E8E429B1BA5B291A56A9F1FE4F9E673F9F53F1B64B1A8542A0F43BEBF6C9867
821CCA0C17F99D79EEDCF909B0B1770C9E19DF26AA2AAAFF958BBFFDC7C1F993F47F28DED6F5341DFAEBF43F3D5B8EE0BA394B6C2FE0A8F123DFBCBF
45F3A1A9CD47137393F323AB82E2F7256ADDCB33D74A1FB6BFBD4EC0A7D87AEBAA2047EB733B48BBF388A2EF51DE7602D9CBDA66B321C929DC5B68CB
3753A3343D89C0423F47579B65B64FA30924A2BCA8FB436BA3C71F53639AF76F5D403EC9B976D336FBF66218B334FCF0F54DAF2629DCE9FEEE89F969
064867B238B4292ACA433830BF0A8099B3C3138B0088511A7FD20110A00E4C30809AB8B395EFC1DFED834C5ECFE12A1582A4D8D10914F71F07E7EF83
F22A6C63ECF3D9D0D2EE1F88914412534A3FACBE8FB7DE9C522E749C1D9E714870B7A5828A2BF5A2F1AE75476C931FEA6BAA75A607140949DFCAB4F0
5ABD51F9C841C9E557C575693CD89620636094EC2864BB5E4B78C7902BD14F7999DB4175F535520281EF2D089204ED9D2B6212C24B5F3BA67E3C7884
09FDF65EBD381C9A3C83F11A322EF57FD958C305D8E38F2B67B974F60293394F87B9531404E5B0594C94C318FEDE4565C3D3D52F6651940BD2EBDB38
D49DD27E3F16298DBFE6E7B091948B17A27D24F5BFD271FF7170FEA8FCCF3339B2FC3677C97F1462C54BAA8B8715BA1A7E3940BAD6A4EDDEF92EB65A
FE40829C957AC94F2F97AC6711498DFE4EDAE5DD4E6F02971DB9B82ECDF6A866C81E55693FF1340C034C45B51C52761B2A6F0D26ACA624886868900D
B7279B854892B794EA0B1185AD4BCEE8EC7D5A5610BA7FE854E430307FEC38BDDC3B8F83CD5C3CF3137D97FA15E4B26975D5B35C90C2664EF7522074
694302C86621204C1DEC195B600DD4DF2F9BC3500ECC696DE28235AA6B071717DE7F66C1C0D0F1E6DB8A360612B93330EE3F0ECE9F94FFD440C51DCE
C52C108550009C8B54D3254AE7BB0BC5BF97F1F8EA64B06536F37D68F8495103D5AC89436A3ED73D0D1F17D89B7ACF3AC6DC23FB27449C1649D6333B
A62A2C4948C7B04175598394AB9672A921C2F6536B85553D098A6B3D64D7EACB6F3F4912161459D3755E34B9E9F08E8D07470E87B4B207D3B25BF7DB
3F4581B7E73F61EC1BE1355C2695313A07420C804B1BAD6D61F3CB7F1ECCE1F098DCF9BEF6291A489DECED9B5F6A4F0AB3677A58DCF98B272781A917
6FA6C6E1EFA5338D8A2415B2531307E1E113000ECE1FF8FFCB52FBB0D70B7EFC8711009A8BF63390246E5E2512425FAF57BBDACEAEAABCE8965BA7B9
B2E129CE573FE7F22CC55D5DFEAE1AF73693EEAFB4D9E575472CC2D3205F475C56E80486DD92305C9FE3AFE4B25145EE6D1A4131D4573B30313E2CD5
34C98C2820A8F8F8D75E11976F39C1D7AF97C5AEEB9EDC6CEF9AE7E0F411FE7EAE700AEBDCB2A916A203C0D824C08180F9EE96AAE6791042200441C0
B9960F5F7ECFC30087C301F95AA300854BA731181FCF3CA150E63E3F988399E30DE09CB98CBE89CE0D0A8AE25D7F7070FE20FB076F286A6F72780942
100A429CA9E040670992768698664FBD78E67A53DD4D63D9954A1FF693E4DDA699179C4F94C8E93F0973515CF9402878A77AC1B6EDC2216B65F7EB0A
8A115A30249028BD75ADA84EB026217E9F88CA467B9D10C5F8A7DE9AA1FACB0404C2BE745808AFF8BCCBEEEECFB3FEB6C38D8E0E2527EC370C4F3E3C
F6898B3E8B3BF38C0AA3486F2B1302C1DE37CDF3F35C10E00230C2EA1FA8AFE8A16330777E96018008880093E35C163C987DF8F33083D270AD7991DD
5B3B8A6C242B7B49850C8028DEF50B07E74FFC9F4C5672DA67FB0E5CBAFC0BE24E058747492A93B65B49D553571D2D34F774E8BFDE18B9FB21D9CCE0
06BD7F7F784784726C52EC4AB3CF966AA7CD8F3C0B25056FD5C8B61016B285B1294329715F5759AFE512169992723BB7C9A6AC712FBE686A61401494
3DF235499EBCB138DF3977BCD03CBD22DF76FF8F0BBA87A6DAF39E4E4C763FDB7BA68ECEC3E0F65626082DBE2C1E851118E42CBDF5A73C78D83BC5E2
8F12A68E51201E04223CEE14C083BB77EEFC014E4FB67D6EE3CCB47C5B808BC49432566B960128DEF50F07E7EFE7FE3C94FBCCD878558A4D05B8946D
73C07EEFC08D921A247D47818D0B396B7E18FBEA9EAFBE70CBAA564FD9C0E1537D955DFE3131D39D1BBC48E7B3450E859AE7F9A99A249A66190B08E4
C3D87BA2A892848C8F8F98D43A7789B404C5C828DDE2BB8A6A5684653E2DD7C44C02C2EF447B1EE9BF66726C70BD6D715D94E6F1DA8B577AB9356F5B
775DEC46318CD1DC4687995DD75F72500486612EC0E04E7F7A3B0EF0C7B834323682F2672884334347D1E93D87A76933A35DDFA8B4D1E70F7F301AB4
1CEFED930E1B4770FF7170FE24FAA3C3EB357DB5D5DDABF8FAC31000767A071F52909311B65FE64EFB943696E2E2163D7CA86555EB41E9E59E4FF26A
EF4616683A665D4A970BFEAE9EB457C231D44F27C42ADF4F40618203171165CC448D0F99CB79996B1ECE960ABDA21570DD47DC404FC865F296907EC3
598715CBCFD49DD4B94C5B67F0B5C8FDD0EB27FBDF72A9376F0E15DEA0F38732F8718807CDBD7DD9B7740919DFF789BE5ECAC4101BE12E598F2C6DEA
45F81F943E0371078A8AE79199E1E1690A77B6B7A699311E97F07983A27639938777FDC5C1F903FF11F4A3B9928D9CB2DB7B04E24F005C5683D3AA03
B2E21A02FE56A4E61F316F5EEBDB38F4EE7F7B30BD5E41DDE36561E262CCEE08759FD707E464AB5C2C1FA91944583BE8A95C5F2DE0C4A2B3D7126549
16894E44F7B562969B8C4DF6052946984B19AB09AB76DE13522A6BB1D68C3D5A192D7797BDC5FCE921BFE69E3DD1DF819F7B2F3577F5A24C90F5A36A
0EC37A9BC7211886101499681D68EF9D5F04D800C08FFBFCDA8407A3FC1C81B748E97AF6A1BD1581077F8D4F01ACE1DAEFEFBADB77650C444BC9AFFE
C54597EEFDC0FDC7C1F95BF6F3E33F7CDB485BD321D2E423BA641EC0AEB54BC8D391B297B13B24F183965A34EA6D63D37025A7D79F924652DBDCE77F
EF78749E917BFF031BD98C5DEA159B2D32ACBD8C08C55B05F67298837A324EC19B5D486A016E467206F1AB6DE58D9D75D41D64D4EF954A12773787C8
24DD7FECAC7FBB7B8BE1991DC739DD3E699D8C3BF1B728F3E0F470E9D7C6DE456CBC656C29F0837CE1E7BA862727A82C1600F0537F0846962A0214A1
B4B6CFFDDA903B81F2B09EBA3184034F7FF83DDFDF5790D3B65D52CBB6A8E35FEFFE31FCEE3F1C9CBF05CA9B4ED55712DF7C53A9020510188068EF2D
35D3FD840D8D0576889E63E7067E396EAF70F44D0C25E4D8337149A3EF973DF21D931DF59FBC7293B73E2F92F8CE352748D54FE05892C04380F34680
BCF3AC01C13A8ABC5C453E2BC022D248CCC445936CF5F19684FEF90BEBF412EE5FD6B3FF54B9DC71FDF51AA4C1F3187DF16CFCCB19EEF0A797BFBE52
E780F686271D3C0C022108849973CC85613AC465F1C70481208280000C539F9D7B4BEB5F7792C6B77BA4998AA0E8EFBA8EE9C1B18B1195A512322ADE
E55D18F2BFEEFDC4FDC7C1F95BFEBFD27337165F7357BE140161FE0FE3BD8978E04E315D4B9100C720764DE2CB0BDB0D823EA5D1AF2BE45B29E91776
276445B8B818B9176DF5502A70D1B9B17CFD61FB20A2AF85400704DD1514775017D2DF1B68EAAB949E41F43AA2266B6D2A61FDF689BCD1ED9A508DF8
F24C4DFFEE0273F3431F188B4F56DC436917F21AD853AF1FB7CD4210AFF36E61DD285F5B90B9548680133FC759207F2E40FEAA4A40805F00CC1717F4
22BFB7E47350189BFFC1E20FBEEFF118C85878B6EA4C85B18CB1EEFEB26E0CC5337F1C9CBF9FFEC3488E7AA4AEACDD0EF275801F77F9F57FB39598F6
1E12599F4C4E23D476EF193994686EF032E8CB94515C86B2B4DFE829E30823BB08C5E347AD84B7A50B3A7BA8A4C7DB28C80A28D1412059485E8860B0
D3C8C58B18FE548DEC602EADAB4EB2FD7E44D2B5F2E20A8B98C20C29EFB60B2A86D78669DFA21D1B796DC78AA9DCF61315B3201B613EC92C5FE0D7F7
E0C22C1B627379DCBB6534085ADA91B474DB27087139C078E5C705B4F5DC3BFE77BCF96E367FF863E5A3F022F031F62C751B414946E1CADE16DC7F1C
9C3F990178434181AB352D4D1C44CF2CC210CC8FB455FA521E970D24EC6D8947ED8E77695F290D8CD0B81F55819DB18E925132FA5CA29FE4E6B8CB2A
EB61886ECC6E926198AACFFA1833258138081EB327EA086B040768EC763138E323AEA920A365226E5E9A2DBEB6E0948343F4E56D5A81DF8ECB1A3EA5
216FB52C3AB937375703706DF2732E0401BD47B2DA61FE70B8FD231C7EB50F238CCFD3FC827F69410282411084863A3BDBFBA1C5B74F47F909013033
B214FD995D0B1897D5B4A59875495456879CF861EB2082AFFCE1E0FC7D105E85618ABBACB36D92C8E621145C5A00BC6D28679ABB52C42694B032D3A1
D7CAA23DF180566C42246F78ED0D53A24CE2375FD360753F078FE63DC64A2122AA89B2449F1C9765028F20B8689988B0B4AFAE54EA0AB57D0744557D
E455B5480E2FA3952ED4FBE979ED4A592FEFD376515CEB210D6AB3F5FCD61C9139C4E3BED9FE9A0D81EC471BCB1808C644E7DBA61094D532C3AFF861
1E0C2CB9CF9F1C3830D05BF9730800FB6E55B25018644FFD06311E873D8162D8D4C77D4FB94564613925C3DA03BED31086F7FDC1C1F9FBE53F7CDA74
A3975AACDB35B188361EB4F40AA0C8D55036F884B662B8B9E6C9D87B65DEA5B74FC4FAE6E8B6617915270972F6BD057E67CDACA2D6302EA84BDB589A
67588ADB3F4A11509983A0FD02221ACB554596AF36D0DA1729651360A827E6FA2450B66A748DB2C3A69381CAFE5FF214EC5EB281418FD0E99BFE9508
D2B5EB6407C49C2EDD76BC1F46995C3AF5F7228FC77A96F690CECFFA91A52D40FC0C000001A0EBED100870ABCF37F0BF9CEC9FA33379C034870761D8
AFB2E2C7DC3B443169C32B55670D52993086C77F1C9CBFBFFAB7102919B6566DA5CE6569B7BAA5F80FC08FF59425C41274051C3C25432EDE6E8B72FB
9A93782457C21FFEE85CA44696CB9D5E7BF8A27F90FBC0557173E50DBA6611A2D2292B04C2B8742857503AD45A6045A86890A29A818EBE91BA885763
805815B243CC7E5F8AB74660FD21955383A323156E7BC62EAD9E42C776BBDCA473B9B31FCFBF016116974BE1201C0C1DFFF8E6C7ECD24A0408F1D0BF
D62360F6A7720A84CC9CCEED80E8EFB3CF0E7261903ABD542D2CBC78D1F588D62C499410B66BFF94B92A8B09E3F93F0ECE1F94FF3DAE62BEDB2D2315
9FD89BBE47001801782F351414649CDC85CCCCDD9C56760CDEF62E7A1A7FE1ACAA5C3DA06314A86B6CD872CBFF6BD272D5BC97A69B6C0B1395769A12
E524850A6106B047DCDC4A40759FB16B9482849A858125D9AA294DFC057440517BC78570AD0D852B239A78ED617EEED90B0733D80BC757865580BCE9
498801C05C3607E27097B6F9767F9900FFF5BE7F29F747107EE94FAF7F3C81A2BF765E1F05C1F99A976300C89E5F60F1D873AD0F2BE0175DB09388B4
B84A4ED9E7B9AF7BA9F8C95F1C9C3FF01FBE6E2A69981A9A44BE1EA75404023008F0DEA8A990098AE90431438B04A9F4893ECB80E6839E57DD244280
F3A249E69A9A914DC1392BBDF8727B44FAEF4913DC705E9A2422F21986BE2A10D484C8B1F2DA09D2CA9A62FA6A64A5EA0BA427BC83921A1BF7B86B9E
B8A11F378F755AC9BB3CA01C2AA05C0DDBDB3A002E2C42C0047D8E320F33BA072114A57E1EE4DBCFE58B0F825C14E6C11083CB1EA91944D0DFB73E41
008DBBB41E08B0A874081AA9AE7DDC0B7FAEC40E10C51575CEBCB930CB3DBE918A6078E35F1C9CBF1DFEE9497A5A26DB1D62544FE52A9D6682100AF2
1A9C5C4964A94C4FF266BB043BB9079450EBE71FA2AF1ED2729EA07884241A849966DD4EB78A37B07D696BE7AF11460E6B4D9726994C22C8191151B2
6490BBD91E0B929191948196C8B9069DBBC371B2A62167EC358ABFE96EE16223B62A7B4A178BF23EAECEFA31CB1EE4B0697D03CD2363009BD95ED8C1
4367BE4EA1FC9C1FE27F4004E641101701B8F37314764F430D03E5D2D810C2E5FC752C18ECB8DA5435CA6B7C823C11170FD99D3EF2A619FAA1BB7D01
C2F0F08F83F377CB7FB4C5C57FB9D156C350C9E087D2DBC74110E6A23F6D122D84851CB334CEBA6A6C56D56C2E914D19DD175DED2D7E857D59629BBD
538AC5FDFB16EE46666F366BC449052A124EB69A0A843061A6BB9018C1CD51E8609218594959414E684395FAD6B7B6649DF4B5FA3A15C3DE4751ACD7
55F604857222B6B8F023B3A76F1EC1167A5A862701181AED1C681AA051A716794BBBFF601800408431F06B164091C92E2AEDFDD96606C69A98FCEB1A
50EEC22C9BD55EDE53DD0F5415019FE5C5E2AF457CEC3D3E01ECD0499BFBD7E51F78028083F377C23FF0D039CD5923C977ADA2CB4BCDC4017EC90DA2
C3AEA9C90482CA117F3F2BF5BD2B54F72F7839B47D887CFFCC2566F1B765F83AC5ED2B3734AFB1F7362BAEF7D96CB0C64A2078D09B70158067CD2414
F4D6A87B6E264AE9CAAB6A13E27EBA7A24AB4B186DD9A4E15AFDC9EA14867DD257BCC3FE9516F1813150F3BA134491AEFA312E44E7327AEB06D8E8D4
381B8141108521800322D04277673FFF8BB1F6F989074F2650049E9BE32C7527610FF53081E1EAB9BE01F6CBDBD47613725CB153115A52017CF3D87B
6801BFFA0B07E7EFDA8FA23D2936C92AAADB563AABC9DFF2766980D808004CBA5BAD242C5B16B12938D7CA30405FA5F70E692F7563F4AF3CC572EA56
D9D552D6E9BA778FCAD8C906766F8B93F37417D07EAD2FD6C2027EAB104457B811636D0415B4E46484535B967BBE5291D448BA6618367275792E63F8
B484D6736EA1EB968EC50FD9AF293017EA7DD28E8C01C050455D3B40A7B391A5ADBE4BAFFCC72600884DA3A3280675BFE9A756D5B0501E6B7211E1A1
00CC19EB9B654E368FB54E724A1F71A73D44233E3AED471BCED1E7C3F6A56CA0C1F8ED3F38387F4F7FBEFF75AEF6294A5659C9417AAA25BB2CDE825C
10E4CE876959899089F647138A768767AF5139D8A3A7F7FEF3EA9E5F0E2B676B7502BD8CB6BA6E2AD2B0B032799CE6EE691E29A4932B69BD48E3B6AB
892BAC34F7582E24A963A02CE2D7B749E6CB3E497D9320DFFDEDEB227E02974D245634BD8875BF3F3B7CE5D4377E9DC198ADED65D2A766DEDFA96333
38D3B33006C34B353FCA199947203A04F0A37D5B6D7B73F91082806CCAD27E5F1E0232690C90393ACFA02C3E7FCB5D8C10F16D4D5B350DDE9BE4DE8D
AC33DCCFC5FDC7C1F9DB1300ED8A4FD02A75D3A413A1F2CA37CF281C5F5CEAB5331F6EB9514E8C40CADE149BB7FBD5BB30AD8A1D423E1DBB1E0F6C24
9431D72846CB7A86AA5EF790F6908A49753AA6E0405074100860D2D9D904714343D50833511D09A965869FF62EBBF5444CC5538FB4A26985D72C5622
23BBB3F1A8D9891FCCF7C75FF1E37BEFA7A98909A47FA6E3F8C10E0E02FDB832C883601ECC2F001080037161104180DECFB523AF2E0E2110C2A102FC
01C3300AD3E759CC59EA0267B2BC9A4B8913726A3A68DF40FF32C6AD4ACE7BAC749503F3965A98E0FEE3E0FC47FF516C34392E3AC6CC627B9E8FA1D1
D15A9BB831880D711692B58F5911244556EE892FDB14F0A56E53CC6347B5EBBB031BEFA8AE6157E9461ABB6E31BA7054CECBD97AAFCB0D375B717D55
816C2E7BCE5448D9D0C3C5DBC5D64B559C70F2A1B077B3B7A86D94A6C95537AF45F889AEF286D2A43D75B4F15FEF9B16A1F9DADCD71C68820E7C497F
3A8BD039734FEB211E8C623C085D7AEF0FC008BBB7F253DD48F3F99BE3FC8980B608FE355C7E0230BF08B3695C664B792B77264AC0BCF29653316770
14AB7E9A3D7356E7D3520353DC7F1C9CBFE33F8275F9DBE825B8CA6F5DAD686AE83E7848BF9ECBE1B029DB84D29C04A54454D75B3EC8965C4F1DCBC9
BA4B74D9AD53D21025F5693EC059572EC3D4E785AADE469DAD86C97B930DC5C5855E00DC6E5182886192815DA0A789A25CF265B2D6113751230F258D
F3A1E673D03909A95527CF144FCF4C777D78F77B86527FEF2D83320732DE9D6B8060081A691E43317E84E73296623D0C2008E563D50FFA62F3859205
1406195480872CDDFD83B1FA86B9340E9BF5E1623BD0ED2568F9FA99D569CAF06F6CBCB8BC8AED6ED681207F35304571FF7170FE73FCAFB55194D9BE
5123325869B9976ADB33E5071C0E97B3902DB23A42444156725DD4FA4213FD0A7832B0C040FA56E2A1FAA77207D81BEC23B4B6EC75280D56DF1BB333
34F4484696B2A85C1F04FE2089134D9CEDE2DDB40C14BC6E7948166409EB24E96B5EDBE7F27B2C49462AFCF1FB9F738CD9EFD9777A983F5F3D1863B1
A739C397EFCCC1307361AC8B8672593C6CBA6E1A03517E92CF9B9E9884696FF2CEF367071066D14184C7971AE32D74B4CC212C90F3F8443FABCE4ED0
F3E70BD75D838BD330F8AEBB1F6D538A9FFF97FFF0D23660FCF9E2E0FCFBF8CFBBA5AF2596BA5D3120515FD3432AFD91E8FE7926C8A29C1536592B22
AD286A9D6C782D58C19FCA3108DC27B0F976D4C5CE60A38678D95435F36293B8DDE20E4743938E449ABD4810941986C08B028262CA520131F2524EAA
76D626F9791A8A6E3A5AF9A9D173CD1624B9D493E75FF74ED05F5C7BCEEA7EFEB9A3973D0601B5DBF29928C8AEFD35C1046106429FEF9D581A16BDA3
FAC72C067C4ADC585E37C9E582CC39CED26DA0180F9B6FFDF49BC141B9CFCF75B1CF6B8B7837FC4E4DF84E1B02B1C62A4E3FB0573C9B5F37F0F54799
BD20EE3F0ECE7F8AFF508E86BCDC91632EE957C394776818DE558BFBCD0658B49B0A721B8C6555156C77DB5DCC34B33E319CA853A7E2F260D7A7CF65
165BC265236C155375CD0B6CD5B79BEA1C5C2D75E529C17801EED21026892B1B445A281A3B68E8AA1F7FA32E16E56B793A32937A5145CC29EDCADE9A
51D6FCA37BD3CC96838F29E02C85CDBC1F7D9D0EA3C305AD5C360270E9B4413A8C41544AF3BDAB8FC759BD87824E0FCE76F7316098CE02F9FEF330F6
446BF708C0C3E8A5450BD44CB25A58E7C8BEEDAD2C061D1DAE0329BF7F582B97A1FC790245782D1F70FD7170FEA3FF63A1F68A6E8FD3EC126E24489E
B455F8B545B98CCBE250EF6A8A84A5EAA8A9C864B8B8E4A8869F7F5229BE6D9341A67B4F598DFD8697662A4E641347910D8972B66A64DFA3EA2B0ECA
EDE342CF0596890B90C2DC45550DA4ADB786ED5E21E1B8C221D4E2D4CC3E1919DFFB47B79E9E00460B3F3290B7C1772766E717A1A9E38915D3203AF5
AC1341E8207DB8636E114221101BC8BF37CA9CBF1EBDAD768239BFC086000613E402080A218BF3935C0443FA4B3E03EC74B9D0F58F46DFE5D4020B13
9CB9BA19E8D5DCAC9161136FA95B18025EAC4671FF7170FE83FFD80F979C40FF1396AA61D9170C4EA5C874978BEC5A6471E88FCD8866311A3AE6CA31
B9DE7906065F6BC6E355135337895C1C7C1FBD9D72582550CD6ABFF2DA6C75571359D75B1ED2AAAA5FB9C07B09B28AB87D9A91A4A681DAFEA62475E9
CD4EC65E2E5739BB092AB165A5B91F3ABA6FA6D7A18B977CAFD18049106A3E70F737D237C7ECA4F2B80C04189F64B01930344C9F29FFC9423F6F5E53
3E3BFC6371E9FC3F07844000E2FB0F2CADFF6194CAE6598CB1D9A370E3F6DEDEE26F341A85C15A686F3AD0883549C68CA34BD704A30BFB7F21B8FF38
38FF1E94F7D42C3F4676958B9ABDF1091BC70384CC26A2F73073915DED44965692543297573CB73A6105217EE46DAD8262F63909CFF6373BF4DF0C7A
EBD99B1CD6084897305B29A67CD8864454F8CD64BFD3561294DBE04C9455155FDD7848C7C43B59CE82748E932D2EBBBE3E67CBA396EF27422BA0D95D
564FE1691AC82E3F5BCB037FFEA24D804B37FBD3C6FBE91C16C468EAECB8D183F5E7F91DE9C5061B98E85F87FD20E61075695720C6F71F6DAA1CC3B0
EF51DB9FEE09AA9979F1749C313BC906EBAF14958058B0781E8D5F27C0285C15F703F71F07E73F847F1E7846FB52A27854869589DC8E64E38B6ACA1F
EC544B5974A02354425E4653CFC64CCE3F50639396DEB7324A9490F339FFD535856F42778DE768AF745DA3BEA122C83D4ADD64BF2E99E4CFE5B0F70A
498807C5CA12350DD4EF562A681FD2765FAB2DF1EE2E4172E3D393791D1F7F969D6DE1340607D7B0A81464EAECA53904FDFE721E61A2543644E9E330
400630563704BDEF076E06A5BE065036852FFCD2394088F3E3FDC2D20DA0FC927EEEC35704E35DB4CF7B1C7762843D50DF4BA32C3280C9BC7A04836F
881B952C2DFFC1C8EFF59B6BF0FAFFFF9D2415E7FF56F9FF335A2CEEA0DEB648FFE50A2B4FC8640409E59C175E4F5D648EA7CA28CAF839CBF81B99E6
68448593E22FB4BF22C978D89C9BB8FAE1A6DCEE5C5339FB7051FF1F5B226DF45DCFB911898710EEE20A0101ED5419016555D1D03B4E96A956869132
CBF22B654557BFD873A5B9F7D7DDC4078BED6BA33BD05180F7F55CF91C46AF7C3C83D266211E6F611EA1B3C091E64A2AF8F2F0B7DD29151086F10B7D
0C81977A8042ECC979106443004CEF681CC7B01FA94177CFAEFEC861D35E574C5366D9ECC5D2CF18F6F189ADE4CA3AE8AF0BC3CE6D3F7817FFB7F97F
C57EFC41FE5F0BFFC0651B099F7B0147E3D7EE324DAA324F2D908DEB905DDE41A7332E69AB2806FA4AED4A0F3D9E68EE61E017564AF55576540F991E
BD5860A8E46BAF6CB1434F3BDDDEDC5D5DBBE89CA44A2B86FE365BA6B0D15840D24F5D3125CC7AE772D3AD8EB2D7E71385573CBF7FF367E7D743F9F5
132D670ED52D0C80D4E7B76A79B4D1D66F0B18651E42D18156041C60B47FEC59643D0FCF5D736116C3A0A5C96969232F7F060041E8AF0381D48E0F7D
6C0C3C177CB57A775A3B0F62B4DEED63D368ACC5AE068871F2468D9271FE52CF6098DBE87825BA02FFB7C167009CFF10FE87C355C55715381E3DE296
A5E3FFCEC6F0B6BE72D36A99A70B146E95859CBA5E10396C6DB6CDFA58198B2391D194FD045B0DD1BBE8E6D4C34A99E924A9602B92A98B94BBB654FA
31197F260F7B43108C0A1110D071970C7270CA3356DF6E2E5102642F737F77637B4947C981E343932D370B3E7067E6A66EE42FB2B860E7AB6184C640
11466B5B675BEDD4C4EFE9BEEEEBFBDE1EBB4643D8109706A02082F27FC1A57B00012EC81AAF691C81B0FE78BFCAB68DDBC7301632995BB9304DA7F5
B4B75138C72E62751A1BBFB04010841B03D3EE19B4E0FDBF7070FEEDB48A226F2D3C1513EE44E6EDB43D1FE6FE345AE36288DC9B9752C72669ECDE35
B2B67A718ABE661B56ADC8D792DEFFD066F68BD94E0BB1005A616047CC817A03E51526721E3B3C1D9C6C74F4881BF87FEFA2B08C9B8800C95BD679B7
76CE265268B8CC15F038D1ECCDB37DCF664BB63D1F1F2D38563E3D3AF0EEDB95FA19844BED6D9D44203682FEAE18851A3E0C23EDACF745DF9FE79D7C
47A532B9DCC9092E0FFB57C74F008041900B0CD47732B89C97DE69BFDEEC38368E01BCC5DB0FD9DC797A5BF9208B7ABB0C46F72C2F9D0101903D9C12
D399EC3B8DEFFFC5C1F9B7D93FC43CEAB05DC6E77EF4CA04B3BCED923BB289F1C984AD7DAE5EBF6834CA01929DC52A11175FF31BC631EB489AB5CEE7
07ED9EE592649EFE36797057397FAF94A69D9B6ED43AE9E8ADEA64C2317EC61EB64C5C4C4078B996D41A99A8BD24CF309178F892A8C68B7CF3A38D67
7754B13A0EA43E7939457B9CFFBA0E8258CD252594A5CE03B4860AC6CCE71E8CDDD85E5D397A73C7952100A02E2CDDFDC5FBCB7C14E4C2000780D93D
3FE628ECC9C3CEF9ADCF33EFF0A33F80F555CF32E758036FDA288C92023A34ED77B08DC664837DBB0DAE5C944967E1FAE3E0FC7BFF27E2937394BD1F
846A65989FBDA9BEBA5C2B3847C57CE0A2F65B3A83734FC5DC7A85AAC106F3AC238EA986C4FDBB0ED12CF6771BCBC47336E77E1537A9F451705EBD62
55AE51FCE14029856A0C9BB3232A2E5B1614495C1F697CC34D3DD7CCE8D32949F9FBCD0E51EF8E6D6B62D72566562C0E0F96DC793FC08099853BDF2F
60C0426BE7978EF1BE1680477DF2857AEF44DAC6FF8FBDF70A8A6A6BDB766BB39602927306C92022414001892A0A2AA88811411141498260C6843961
8E18312744104582029241720E4D6EE81C679EB3E7DFBCFB3DFA3ED73AD9FBE8AF79515054C30155CDFD3CF73DC61CE3690261BE500C63384962B854
FE12129F9EFB4960DCF671B1B0257265DED0CD6345A2E933FEFC7140C2678EB67588D9793B5B00E0EAAA92499E903F74487B67C66CC3FB38A5FFFFBB
CC2AC5FFEFE91FAD9F1793ACED726FADD931ADB8A2050E85CB0C32DDE41FF4CEBFCCE141E5CE9AB6F67EB2CB161ADC3A11E3A462191CD061AB551BAD
685A56E87EDE52E15CB4ACAEEF19B7ED9A5667CF291B76814481BA9A9CECC63D33DD62550FEE96DBB2D0E2499292EAD9A10DBECF225287D1D2E09491
C15FF93967AB3B61947620A50E1E67E65F6B054790660E46949E1F119F0974BDCD81C5008F2D02F9D3F77E833C48EA2A0889842FC23126432CAA0CDB
DEDB7CF5FA6F3E21260970420460E2B1DA21E650CBBDBC29B8D8E978DFE4147FE0BCF6F212BB394E2F28FD5350FCBBFEE13CEDFD19169ED776BB9C36
0A6ADC655B7C56E7F25A75EFE16D5B0744206D8DA29BF33A45AF8DC6EE292FE2D59764BADEDC38EB448691D1153420E184E6C6084593F96766CD3757
4BDE3B6BB94808272B1B29AE3FA3A39ABC34EC989A75D492944433A3A889D71E4FAF1E6591955E7B865937579E28A9E180506EE4791EDCD5FEF9C130
C9EDE043DC4F37B23943A717EC6E4060800FF2456231882038018B90E9F3496465098AA12298F931FD7A5DE98DF75D080F25493E072631B4A77352C8
EA296D12303B362EFD3A31C9A43D3276EA3A6310E3F883BA019482E2DF0C9584E0DE50389039C7E1F029CF44679DFC53EA8FBEAB6C8BD1D3FA78DF3E
472860EF9737B6B631345967A86E9FFEDBD3969EE5996D607AD84BD58B17B76CD0D734D15CD164BDACB79BCC6C8FBFA240016FC92CF9804BDAB3562D
3339E0A1E11D15EEE1ADEF5D53E475ED5C149D6CF789618FA5B92435F54F4113C7C26B20A2AD72B00214D58F4990E22DA165A33FB6EFFBCA93206204
E4B0C4BCE9118418C013A22086E11D99341226F0C627A5EFDFE4BF2CE771590449024C505A04067E8F02E2013A7D629479DE26BD7B98DEFFDCCCBA78
C0D467AD7D1FF5F83F05C5BFEA5F424F547C9CBBC8E1E4B53507576AEDBC6D9D5AEFB4F98AB57526CD3E6D9C2F7A66A2B758DFD3D24153C7C1FCD30B
8B1F4CEFD2D3EAC7934C2D1E3E5FC5CAD4BCBC4BC12D4E7143C80C4D1DB96738C2F757B4CFB09773D8653227CCD12B234169C95AD367232BD38AB654
12EDBEB11CFC91FF9951A610AA587B7C0C03AB0A81E1FED1EF6209FBFCCA2BBF697B367CE54848118C020C1E028308280200110822308101FD4C69A5
1AB997DD95BFE7697DBF106089095234268D0624BD8A4772057C368339D9B23CF85DFF68D7F539862F80333A37F537F2A8FE4FA57E8A7F83900C86C9
9F2F5862117769C7F69D1A0EAF837C3EC6D83C59A7BD81BEC3F91B97576A31332978A19591A9B29DE192E7062BC5E7B6166BFA6D3057B6B9BB30EBB9
FCBA4C4D9B04F905FE0AB27FABB649C87E63E59D3EF26601F63A814B67AFDE34376895EA99DEB0956F37BC40AA56A60AA1CFC107DAD8422077E56510
11BF7FC69D2812F58DE0F5919B5E22BF36EEFC2521F9388E313A581080C0ACCE6E1102E1423E8600382A95FBA3C807ACAB3EB9A334983F2E7507C0A4
40EA0178F50C924050CE248DC9396C935837DC7259DF280B6BB7483A257F01FABF5EFF943628FEBFFD0349DAFCD522BEEEB6D979F5484C9AADDDDBBD
61C75F3A1D7EEDB268B4D472DF30BF33446EDF3153133B476D333F87F536332F77658C39799FB4B0D40DF5F57E64E7732F68CE7613A3A00015192F90
24BFE99BBB2B9AAE08D5763DEBA6B8C86A69AC6678CD01C7C70FEED21E6D7E00E065B119FD74B1E861EC0F18A93FF78005F40A44151D9FF79FA173DE
6C3C3F2EE14B0880236A69073000144D8D8BA6078008411C172328397CF4D8504FE4B2DC31164738C580A4C9404C4A2BC304878451883DCE9BE27E75
5F9BDBD5B847CFF91BD495E0F5DDC1AC5C42F97F0A8A7F6B1F98A4D256D7E8CC31BDA5F117D6C5CF354D3F1917F57CB7EEA354F55BAC70FB2F9383C7
64BD37993A2B3BCF9EB57483D76EB9D93F57569D50BCBE7896F6EEF38AFB5314B747C907DAFD659BE523B31D25B08D0A7ACAE60BBD4DD577AE56D672
0F0930B02DBDA21FFF3CBBF2E4A9121C7C1373F897001DBFF8741083F3123E09C0DA5EB0E8747175FB647F72E247882FCDFB537D23029E500C419034
07600884880588104224C497C422A22C24B992D13B013318A01810A3080E9342168948F0EEFA21BAA07DD5FC8B75DF23740EF74BD26E1CD97F467931
8DD23F05C5BF767F5CF2D578F6FCADD77D1625DE4B495A661B9775E442CA03DD558FBD57334ACC927AE8EF8DE784FA441938CFD5F24970BAE5A458BA
7A459FDBE15B0A4A27BBE6447CD25F163D6BB9B7BCD6DBB5324F7062CA5A41CB74C37C73D5F06445F5A57EDEA6C1376EB8075D6E7B1FF15E0020AF96
9FEF60A25DA78B310C7993500E8FFD6898EAFE54C3FAFE316DE7D52914E7E39CC6A2718904164C08206C7A0E30004E5FFA27C6C8B1FDE1ED4856D20F
883EC68685420806611883618968627A0960AAAA738C458B743A94FF345CFF214A3E09EEDE722442F9AC907AFA8F82E2DFE48F62F8033D93AB7E2712
3C6EE41D5ABB4DDFFA43D2834307FD5453D3B4B2D158CBB7E3D52EDA0B234257E9CCB5D04BB50B5D37336EDB8CB759262F2D64FC5A0EBA9E74725B3F
CB7AB9AED27AF3BF6A61B445436E5E84A7A19E7D82B5A24DB0AD5C70C755F7B9A9951F173F0151F10DCFB49FE3D0EF0B4D84047B98D0C6ED7CFD932F
EE1DEC3810B277FD6B2E09E258DB8FAE462E296136D320B1088411004084000FC0F1D73EDBC7CB4FE68118C24131988B48E53F7D148818BF314C2224
7F92C9E6746C364F787B63B3EA75186FF01F494DCD5BA2F11AA3AEFFA5A0F817FD13304A64EACCEF8FDABACB6CC3C5ED4BB6D9E95F397FBCEEF00A1B
FF749715A24A8BB8AEFEADFA6EEF63D3FDAD57981E0A753FEAE0B34EC5B9D23E78D50C83FC37ABF705046ED1330EB3D3D2F8CB8029827F9BDA6E0DD2
D4F30AB152B45AE8ABEF5053E36115D7FAD3FE1C20125CB54AEEE817F6DFE99660EC6B09CD5859763397C7BA79E168C68EE89F382621FA4A1AFB0876
DBB0786412C111A9FBC72008E10901B423D2F37BE5EDC71C02E762240989310CC651084490C6E84F248982228910688C30D974EB5288C6011ED4185C
773BAA264AC3BE5242E99F82E25FE42F951296A1E2399C3027C84673D385E88B166A3EB723AAABCEAED25CB553FBE3D4429BD7C327D48C2F3B856499
FA5A2D3DE094E860B4C56EC6B59B4AF196323B8EB964456A472E50F334B75794F114B2C1A786AE8E6A5AAEBBA4DD7F95AB86DEABA1D5FA6B1AEB1DB7
F1E1F118C33DA59D68EFB50E82A0C5C40FB1F21F8DF1F9DF4FBC7CB7EFE8B12E0929E9795ED389B1715E710F8C93125C0CC3DC4136C80110DE71F753
8D151F26309227C209090180040A6390341D746597A1388093C2413E74543FF4EA497FF94362B8757DD9AF651F52D5957632A7F54FBDCB1414FF64FF
A5491A4D57F11BCFD03BB8DB6ECD93839F971BAC39B6A583F1FCD2069F37C19B68DBB40F767F9A67ECE36271E598B7FBE2E7EBE20366C706690474FA
F8EC94778F554C3CAD1DBD62B69FD102BDBFF6827CE1DA198ADA3A9E89F3344D96ACB252DC357E447B6955D342771A54B74E2BBAB797D772B793C48A
430F0E335E7DA48B1B2F6F7CF32BF64893182704F70E754A68BF7A09808E4BFF301C2124BC8A5F200603AD9B57BD11F6F66184800F4FCFF3993E0D44
602884B19A3B4004142318C6E27233AD1C4F5DB1B3B82DE64F1D29180F7AF6C34ED5E0A3B45A50FB631414FF9CFE4738222859D1B1E3BCFCB1930ECB
128FEF3B6F7FEB907AAD30744BFFEE3B3FED9FBD54F7F9FC7B89C6EC040BE7DC9D76568FD66EBF61332FCC4EEDF51385F085B247CC358E5AFBAD37B3
75F593FFFB29221E3295B5D0304D5CAA68BAD4DB5C7565DE495DBF9F53C1FA65E817378DE8B6F1C6AA174324FE21F422AFF96905C4FC109791752DF0
1980A078E7D1ABCC515A63EB146F542C75F8042EC120661F1B1CA53FF6F3CA138F03382086A4BA47494C8448A66F0117D4150CF2403E878932C5E8F0
8779DA41692EDEBF00FAE8D507DC9D27BF266BABF9B551F69F82E25FF44FF2FAF90098A26553F9DC34EEE9EAC06389999FE6ED3FA1704EE86C447BB3
612C7D43F772BDF3ED29D626495B6CB20E2DB65A13B975305D276695D6C66627CB63069107676D5FE3B673AEE9E63D4AAA1D18D260AC6C60B5758D9A
9E57A09D7648D16E8D05F5CC04A5EBE8574BCDB57D83AF3E944D92AC2B1BF3273FBF1A861AD213CFB4A77B1620100C7F38D3D55F56C54571215748A0
F8F4D03EE11407191999187AB0EB6E0F8B0BA3620427100C8745B084C47182F3F6F520C412F2F922482C198E75D4777035DE4683DE64D5BD9AB87BB4
B9748589FD6501257E0A8A7F8690482DB9801F6BAA9559E2E490B3C5F95DC0990633ADD3FA9653E18ACF6B2D8B1A1DEE949B06163D70D05970D27B73
90AF9FF54EA7A60E4BFF9D463A5FCEC86C58A773457B4E8A4198877678A29C1D07C39EC8AA1AACDEAEA6E937DF5A777EC9553DC732E0B47238F46DB6
4A444375E6B776046F4FDB53D279F51B847E48DC975D91B6BC081143F0E5C32CF4F36F820031A1D4D08B09121630C7E1DE4E3A54FBF865AD08E4205C
3E8C01001F9E9EFE8B612439F9B480C76681EC51D11443420B56D7B0D6538D13C0CF7731BEFE285AFBA1F5BABE8E63014EA57F0A8A7FD13F51DA2B84
C7D6BB5B0456869B649FF07B1D7D6E6C93FA797FCDEACFDA7103F6C71971F387EED8DECC0FB130CF4EB15BE816E51669FC89E9ED76D2D320ADCDCEE6
B961D07C9D749BC5DB2D9D9C65379238764C463764B3BEFA8275F366CFB997676D9523BCA2E7D55D305B3BADEAC5D96A2E2C7C187196F6EB629B847D
25EAE8D7D103E103B810E3249D8684F57C0981F3410C475118120A18F844DDEFD1DF8F56EE1F802050000A01A188C5856001864D5F0AC0286E14036C
C1049D8F7138F0E7456616DA066705F89B5D82C9F3BF775FE3DE77353075EF26092AFD5350FC39FA4B48092ECAECE0818321CBE6A8BEBD289BF266D1
E943DB46DF696D38A2B677CCD1F847D0FCCE8F8A97F8C1CB9FEF75303C7BD224798DB163A8D3FAD1258A89EBCC4CAB2FFD7D71F30C17A5707BEF1DB6
866AB229242E0C92515C6E3B6BAEAFBB816EE2571FDD2CE88AC182923796C6B78B526F0E08E09ED4A89F9DD5999D92E19B870F3736ED3ACC9740446F
448600E748CD039B2900C42C041BEE42259CBEDAEEC62F914E997C1CE6000286186433F8080203084848C8A1AA519225E2D4B5B191F1617238665568
4AC21B127F769C3F759253B8A2E6D9610B45E54D7C6AF59F82E29FE43FCDF0C12E2E448B58642D77E683D2AAF2F529172287E92B6C2BD706700F6A5E
48317ED3ED3C5F503AF7D0C940B34DB981377217E8853CB63BBBCF202C2650E9EA3B99B077CA8E0EEBAC0332ED8D9466BD97E0F52A336769EA38F96D
B4D3F74C5BA976867DDED0F3F363DB39CF0A56BFE443D0974D673B459F3E8C49EA8F3F3DF6A330EC2944E0C4CF884C805DDD000D0F8B51112C953832
3884B3BBFABA9A0FEF3E59229A9E010E0B8528CCE522A858046388B8ABA60B8730F158511307455B04DC8D8A6A71BFAB39EC77D7B85066CD44C8E39D
979BECE4D5EF11D4EA1F05C53FEB9F9094C5774E21FD9BEC6CE442BF595A7D4CF47D1D923B7653E5DA238BA6665DA70CDDD5DD4766BE1106B81E0A35
357A1F6BFFE1A2AE71CE1AA3EBD67ADB77CAF9EEFD4BAB34F0EFE5ABEDCCEF5B29CFD2EE03911C1919D9199AAEABDD66CD8DF553DF23BCA3ED599C63
EFF0FAFDE25752179FB5FCB5B0E3E16710BBB7BBE276D6C5E07C094EC20F7695D33AAB98623197100B309C0BA1FCAE89EEFEDEB6ABFB4E0F8A081000
70842946210084C542310261A2E20F2C9C84B85DBF7B21142EED447628385EA63D2E9EAC2EC2B899ED53277332D6321FCE5234A890E094FE2928FE49
FF242E7997DCC74687362F4D555CF96395E98B9C0DC56752595FB45C8B434AC42186B7032DBF94EA0508B30C776CF2D03E7641CB2C25442D28DDE279
94EAC6F4F926DED69AD98F55FC7786EA5D7057950F41F860B1A1B299F5DC0DEBB4B5576DD0DF2AF8696E5BF873F98267796EF7A5C9FDC4E2525EEBBD
CF62E1B1E5F5552742432A0892A839787180F6BE1F9B1EF749F2D8382640C43F2B455D05BB537E01984424553F810A38102E3501D30B851848EB1C16
9338C41818C25014AB7F3A715C27A1B73FBF5632CE96E0AF8BC50F7F95BA373568A968ACE82029FD5350FC49FDFFEDFFD8F5E80E0134181CF2C9D4F9
F30EB553A3C9C7BFBA3435DAAADD3AFD58FC5425365921A123796601CB69DE066FCD05F79C55838F28AB1CB48B786D641B1D3AD32054717EC102EDF0
7D72EBE7CACBED41D9C065358D657EF3637514FC02547D8727BCCD5FFD5CB6E67DBEE7111132B53DA0839B7BE89D98B5C7BBB226CA67470741B28FC4
7CE6943F1DC6311C2791C921903F0A30BF14F35A4F39A7D349121AE7C3620C05C69930024ABD3F8488A7CAEBC708E9B7623E0A4118D19A9DE7A3BA6F
9076BF43221C23918FD5D0E36F5DBE7780A53395D593A9DD7F0A8A7F0D0013E7A2DA85E8F096A5B9AB2C5F65CC4E184DF61E5AFE8A17A3BD33FF0287
BDD5FB51A0D3974AF3B5C075A3D54156731FED71093AB950F370E8BCA65895C86D0656DB6DD4AFC45B05EC5174335596CBC2C53C0F19252BB3CDA19A
D6F1BE0E555892DAB1C1CD8165A5011B26F0DE353EB4EEB35B5E7027F638BD2D0DF53C35219194071F1CA77F783786E338CC9A1A024443654DF4A6FC
DA8741AE8F0852828D31710243B85C0E88A0100AC1426ED3AF010445A5398044212E46B6E68E47989CECFCFC7E00072092F72617AD78CC3F94023C99
A9A4A11FD4FF9FFE4F15000A8A3F1700F2D7B688361E321A617D6BFD8C0337CD3D2B6FBA0E25244F66EB3B3E5F5908D5D9EF7D6B7BA8235DF9C7A483
C50A27ADD51797DA2C3DA8E218A31C5BE369B3D84ED33D42DF396391DA92F9211AB272DF097064CEDFAA4A867176CA4B17DAE62057945655EFF6FBF5
615D6025DEE0193854179B9ECFEA8DF47EFD7EB9FF7B1E01667A3F857AEED6F2301C6DADAEE9144AC425DF07066BEEAF9B77B449AAFE291A4F82E204
C815C0B8581AFEB9BC490613C67114E442042475FB487ED6E0C9233F07DAAB38E2812912BAF70A6EBAC97A9DC01ED29D394BD9C1A3777AFB9FD23F05
C53F25806FAB36D508C4BC0356F19BD502337DE73EE98C2AABCA6C69F45C7CC033890BDE5E9477697961837D98F884A1BF97B97DBAEBEA2DFB8DECAE
DAA9BD2D72F65AAE67FDE0B0E99EB566865E978DFE36669058B3B19CB1D982DD2ED6A1665EC38546F34A6FCCCBAF5E7E32935BE318CCF818973F2EAC
59135A9CE31D5A03E2A3DB16576115AFB9040F99A8EF26D82C09FBC987017169A475C44F8244BA7F0D4E0FF9060639B034E74308AFAF714488601832
7D211006308430EBFC9ED617D77AC0A6D649B493498239F731D65D5A677CEF94CF4C5979E788CD4C584215000A8A7F2A0084E4ADBB7F8E88079E525B
7F718EF6DE58D32D538F4E09B34F0EEE50BBB4DA60141CF05AD29D9C307845F16BB18A81AFBEFADE859AEB32E38D77FBCF72693D64B4D04029F2E8D2
C541C1329A571D655612B8244F7196B2F16EAF197E167FC5149B2A5F79A77F6E60F99DC4CF7D0BD632DE277608818E15DBDADF2ED8D4076395DEAB27
A0C2EF84B85B34D48C22A34304ED459E809FE5B0E0314092A2A11EC1F4EC4F4E6E1588A1000283CCAE01218AC252AF0061A8709C014FDD3F577BEFBE
78B0BB1FA4354D9164E92D21FB71D7587CA578D94C59AFF9A13B36B150CAFF5350FCB301C0AE6E097E29E22159DA9155DB35D61F755F90D77607C8F2
6C2AD0DBF1C6F6353C7550FE6A6374468F67D8E80A7D573BC380C4393E7B9EE8DBF96A6A6EDC6AE2E26CA86BE2B87BD92E15D94833994C1227AF2AA9
EB875CB1733861E8D313FFB7EFCDB9CBBAA2131F26B3372CEE7B9E4883842DEB43CAEFB82DEF40C0E756BBA7046F2B08D1AF096404C67B46F19EDC5E
9893A4B7B38B204522064C48081C1FFB3180E1308A20202016C22884A228826328679205F2EA2F7FAA6C1188869922611F8BE0B67EE89794FCC6B3DE
A08767CA5BBEF2DDF3F4B0345848A8024041F14FFE1F38BA6BCD7D1E037DA9B0B4E694B2C3A14DEE115FBFF77FF0BDD79AA4FE35EF2CD876D06EDE68
897FF93B5F5AB1899EA3A6C6C9E84506F1F3F40F59A92CD93F5BC565AEB99CC2AE954B1D652C94646E4B48EE6A1515B3079B14F6EF902BAA939BF7DC
D9AC283EECC7D28EEB1ED5B75677E250C3F2E51F0FCFF3AE82E09B4619E29677ED04BB56800FF0A16E00ADAC63A0DF027D3E80382092004202954C8D
8827D8240103280A02308CC31008213882210C9A48447F76EB4A11844FF483C2BE31113AD2F4A68DF5BB52F0FE3CE7A792C28CC80CDD1B97E3A6086A
018082E21F56FF4842428FF70FC9E2F0B0777A4E5F9F98399D4DF30ADD59DDD57A2EF16A83FDAAC9929E7AE7ED731F837BD278875E0369DA73F5D54E
65EACDF2DCA4773C74D6BACE756A6E9E1B94955CF7FBAE94539E29F30AC5AA6DB43437BD996D76D5366562A972D61383B3579C2B6F1CB8B3B0E4D6CA
4652F863F9CA4BC7EC9CBE42E06D83EBECCA3C013E522C207A2A509C0135BFAEBC11B73C8B8181220CE733601269A8114C1FF79D9EFB0D82180E2110
0AC0D24F46DF145FD478E641274D8488270088C184458585EDD51C823B467BDBCB74D3D6B6CAB25EDF1376962FA1FA3F05C51FD52F99BEF7B3CCD13D
E012938B3E555258F3DC51FB64BCD939EFACA1CE4FE98B7FDC5638CEEE98F0591819C42B5B585D13D9D5B540DB54DEF79ABEF2FC7D7AB30F6B6966E7
D819CE4B305350D9B1C65751E16F9536102D9CADAAB66FA59CEF3CDBC63D3357BF30DF757E7EF153EF5B9B9E673997E3A3179606DF4F7670FE00405F
CDAF718B4A40A2BB502069F8C14069E31505D7036D2F0DE0B0588C40028684E4BD2F460818975A7F0212084018C7A73300CCE96BEC1B87C672535F08
5088CB628B21DA140296DD6B28659064D5F8D33CE494AABAEA859DFADF6ACCEE0BFEDFF64FE99F82E28FFA7FEC7C69CD11260FCED751D67B95B8F0F2
75B75727E398639FEFEF8E69F3D1CBE709DEA9A57B97E191E948F23646B699A996EAAE68474B4F5FCD3D2B94C3BFAE73B15EE6A9ABB9F1949FB6DA5F
4BC4205666A361B2DAC46089DAE312559BB2632E67BCEF76796D2B7A75C12F172908F1D87F3E6DA957AE18E8F43F3E52D2898ACA3F01704303202C7C
DE921D3D7F6D1582222204114F0200D8F5A9538263D2AC8FC142A150041104C897A67F64B26508E435DF4DFF21FD016D0241855C26282ABB3BDE4A27
C9F6873FF2916653539DE42F736E092FDA578094FDA7A0F827FD1308BACF2DDB3F7A840D7F3154907DFEC9F768815F529357D1546DE673EF47CFCD4F
BC9D0476AD88D90EE52F1C2CB47CCA592BAB2E33EFCDA2993B2FAB5BC4A8999D5CE36CA8B3C940D6F1AA9FEC2C99A3101BFD64A86669646C3A2BA9D5
5123A67883EBE24F155EEB5B3EEE5A9C8B3FF1BB10197B74ADD77B1138BEE918ED7A09C27CF06A08CA6B808917BB861AF7191F1EC120018C201C061B
EE7AFA0B22090499DEEC137011445AA168E54308820A2659A0B0FCCCDEEFD260305ECF4620BA08E6D77EA8F9D6449213D77E674F12E11A461BCAFDA3
26DBDCFD4A01EAF93F0A8A7F6AFF04C68D72BB131831C083DF9ACF56BFDD70E2C8F7F31BBAD6C571B8CF2E1E08295BFD34F93632B5F5FAD26278ED2D
38DE67B440575D59EF71BA6C7A9B9BEA513FCB47FB8CBCF4B7F8A9691F3BA23A53260F16C1D7F4F5EC743D6DAD1B5395031BEF5B3BDC6E58BAB4AB23
707B39F7644AD133BF2D29612F85A2C9E864F6CDFBACE2B36FC0BEEFBF88A90F71EFFAB72C7C04A23C1804C5433498FEEA422F49A2188EA3E8F47EBF
04E78282D2523E31FDF82F32FA26E3690D866300538C88209170AC3167ECF7478424EBAB4563A23373EDF77CDDB7A097B3D369CB4F946AFF1414FFA4
7F14EBF5364858BCB889057F5F60ADB6E0C715A74BBF6CDEBC9EDB20FCEC777BF3EDF577BFB8F5E2BFAF1D8E071F2D9B2A55392608FF5BF52FF717DA
9E459EF261A1CA871FEA79E8DBC5E8C8389579CBFC5D0E439321CA9AAA568E4A679F2B19DCEA70D4CFCA7659D44D0B496530C3E2A1DF4E3BCE2CBE24
82FBC253197733F9F7D35AE88DEF7AB0863B9722362E5AFE1943380214E136FC620EBD28E491848424E0E9B53F0497A07C00C7593C1CC685C3C3752F
DFD6611242CC1363D2A800F32B2A065ADB9F8CE0928E6AA64472C8DA7EEBF875CD8B409E7380572B8151F19F82E21FF48F41B91E86F1FBDCBE70C5BF
9798E81A7CBBA39BD2EA93CC0BCE1055DA9C7B7EE850E860F80E00AB2E08783D19F992E5398F596FA43E4BE144A44E88DB1CBF089BA5792116739C8E
CE9797CDDF2523FF1B163FD455D7D18D725D58BD547541F351A5D4CAB96695632B43F9F425DBF9B4D09559B14798705D6832EBE7A989FC533D5051FE
10DE1AFEE2828552120D06201811D0AA8A47AA4AC77169E72709A9DE710CC170708227AD55042284899EEF3FBED6C1B8048144088E0310D05050D027
F876BD0E4547EBA47E21CBD067C18B3E37E7417E54BCBF6FD77F97FFA9024041F1BFF48F6123BBADD5D38AACAFB239DD8BB4CCE50E553A6FCC8F9833
FCD27F62D46179C1FE8F2B2E96D93D2678EDBB22FB3EAD19C99C7115889295FDCBE38A825AA8AB829FB742DA7E25A388B71EA6322E73656CC701E601
0D4DF5B507350E5C52B38AC834D8501CA4F562C0D36DA267D17A51CB8AAD6F166FA44B4ABC12263EA68DE59CE530BF9608D0F2F8DF9F4C9C3F00381F
E581ECD61F9D5DB5FDB0D4F3E3129823C67118C588E1269E04C361588089DAF2AB06400407A5B24649546A0F6A2F7FEF6CABCAADE78DD31A4770F299
89AF55D06886D643FE8B6507B596F4E194FE2928FE49FF4091E31C9D94CE8587A678831B16D82945B11276D7DF302E980CCA11ED722E3DF7F9C9AAEA
F3CB9AA1A2A38733FA621F4F59040A4FFDA52A6F71C5DEE0E252F915DB0D0277E9EB2F6BD9BD5446CA5918441E1B2B3B5D9F6FFE29D0F2EA2D8BF925
7BD4CF0E046A7FEFF68DE8BDB9EDEAF9B08D7DF847AF1D1D8571BFBB2F0C0B72BEB290A2ED0D0D7EDBEA504404B37FF5750D8D422393D3CFFBE1044A
1F1449709810D5D5F00989340570195079C1088E1200886220424EEF0EF4DE2F614F0ED772E1A9728618227E5ABB463B154381914CE4E099588D9011
4AFF1414FFA47F647487E21C75EF96255B86B823EB9DD6A89BD61D59D091A7B64D98B159F84CE5D9BBA3AD1B37F46E089D1A8CFA7AF643F9BA9E731A
AD550A0A0AF24101DA5BE2D42CB7DB18381BE969DD488C929B2123F30A8139314AF28191B3623E59AEF9B2DBF6E925C5D8E1355AEF738D0F318F5817
E49B05B0F06CD7DD9D5F763EAD395EDB975D0E72AFADF879C3E6D0182C12118CB79513747CB2933DEDF809049F6C194561021F296E010950FA12FB6D
F6875F80440242084C703A4089D4FF4FBD2E9C9A6819A78981B2660024EA1CDC76BA66A2673DEB26EB63EE9A6B050FE292FFEEFF53158082E27FE81F
ED7454DDE2645EBDD9BF97C7DCBDE294B1F2E19CADADDFBCBC479B16B7F65BA7B665B4DF722E7DEDFF1EBD7FE0C705C1B5CCC1C08BE021650D35B325
36694FF5AC3237E818189A19066CCBD091FB4BA616463E28CC50986FB4F091876F5B96D1856A6BBFFADD0A077A16450DEC997FE3984758E5C4FD453B
7EF79CBF3D54F4839D57050A531DEEC7CDBEC11589712EADB89A25FDDA2624088CC0A51F832D2292C069151384B4DD633CE0F5E2937DD3B715E31846
484493D25FC1B9BF7BBB9B8B7B18202C1E1640E8D4D2458F76ED839B3DBFD0BA9E3D39AF6DBA858613D4F33F14147FD6BFF0EB5CA3E36BE5EFEF36FD
C9E326079ED499B5A4EBF2AD73B1BAAFC43B5204DB1D2B6EDE2CB78BEED91DD4493F7C626B136BD7EF8B2E238C402D394557CD559775B476A6CE529D
ABA06B17966424FB97F994187EA73653CF60F69D2885C4668FB91FFD1754A629EE29F0DFC18ED7B8956172BAEFE97ACFD4B6928CE3D54D6FBBDF944B
18DBACCF45ACCB8745204CAFA8E8C4C6C7695C8890367F0206095157178049E86D188E23D85849EDC4E57B80D4F24F5F130441840447618CCD1BEBAE
1F98404410441BE6C182E0656F0EC53027D63DE61CEE49F8E4A8AEB9850E1304B5034841F147FDF71EB5330A4F57BD744EEF1D5F7CD62FDA447D3EFD
A7DFEAFB9EBB817CFFCE42EFF28EF3AD91CB463F05BCE8F8F831F5269218DB147605DAABA0A933CFC67C9981A6E14E0F5D172DDDA54107D6CE944900
8560999DAABDE3852C15B5F404DD13278D32EFA9ED60A6BAF477B96E8E0A2DC76E597BDC2CFE7AF1721752F6B2BD1B1ED9A1BF7BCFC9518CCE412ADE
B70A85F47ABAB4F34FAB9FC5C561D1249B40D94C8400307AE9ADCCCA31B68498BEF59F40090CC549708C85F2798DE52C69218010815800F2FDFD2BE2
D7B662FBF74B4E95FD3A9BA7E0AE73504CE99F82E2CF058090346D37B3F2BBA317F354ED0243F0C8D8739EBADA838910AFDFBB0DAB195BD2B95107A1
2F1FAF5A5EA01DBBDDB8F367CBE6B16FD61FEF6B577C5252D1D05FA1AFE56E26336F9F81AE8DBCAF73E82A19994C902BBAA5AEA6E3981D281F724E3B
E6A3DEA27DC6E115679DDF0F6D8F88983D20B96FE678B1FFE0E557FD48E5EDB60166A98FF5DA95B94278009DFCFA99C106BB2A445281E31271D7AF31
B14448E3C002E1248AE278C7DDDC9EBA3A9E54F9D3D680C42084C430FAA0186694557050144450802186861D026BC33734237723842F6F11712DAB56
67195F8791FFFA7FAA005050FC0F03803C5B601699F639765B89C9BE716191BD9BBFB98E53E1CBD0C28279A78539EBCA9F2C6FA6BFA85D16C2F9B9E9
DB9560EEC9B8C91D2B7F2D8B61C56B6958A77BCA6D70D7D7DDB3C5C24D75FEB2C5DE327F178342FA7A19432DBB0D4A466FF7F8B65D315B6A35BF2EC9
F426E7D1C26B8BE7BDDAA3BBE05D55746E5D0B5A78B385F162FD8AA0F5076B279AC540DB876A16D6DB57C7210884C0AA5E75F02418481F17823C0C21
04E0F8F70E21CAE44888696F40A2BCEE9671CEC0C824048CD7F5F0A60F08A0E38D0298BB60CD7844148BFD756D67D731E2DEA3F615557E16DF10ECBF
F2A7F44F41F13FF42F3A6C6C10B3EF52E68A52E7903E7693BBABB79D995A48FD896D53497E2D23A7D77C3D9DC82829B933FFF1C89A05E51EA913AB9E
E61B9FFD68FF811DA2669C12AD6D3377D5DFF32F7999CE9AB3C64C5366660B246A73FE7FFEFE7B839B6CFC61ADEC3A138355FBF706F89D38FB7A9567
C41C4743B39027DD516F07DEA1DFCF4C0C9C73F27EB2757DD9582786D47EEEC3269AF89D341CC589A9BB3F2052824193038034D5E3D850C9E7631F79
420E84E1D8F43C608C3330C86409B8A058C0A1D3A4514000A2AC321ACCF00B65EEDB3D456B3E56D01FC56BDAC18A78B14EDEE53748C99F82E21FF23F
2BCDD43AECD08DEF9BBEAFF42B9FA2AF0FDE606A6F6173B9E0684BC5FE8FFD0DC1C747CE7FEAC9EB8F5AF1E18ADA81B7B31F5727E4A5A44FA5EE98BA
A5A5E79AEAEB65196B37376DAE92927392859C8C3F0F866A6CE5D4DCCF5ACD7FECBA662A45D1E5C1E345730ABF397A2C4B74D5755B79E47BD7DE1743
5950D585918673A72276C6ED1FC505B8F0530E0B1276F1DA6B858484F89DD923212508832B10C2624CC4AC6D1F78F7950D2108F69F8D3C12E14D417C
102525A470820F620424CD0B386318ECF65A5815B79A33D93156541DF38D16D6F8EAC4175595A061CAFE5350FC43FCC7AA62D55D4C83320A17DF89B3
B835C0DDECB77DA18DC30AB7CC0777B82FD21E0F3EB67E013C2EFCD85EE0B4BEFEA4635B9AF5CFBBD74E2DFADCECBDFF95E52C9D701BBFA5F3237516
3BCAA9291D5CF1B7CC715808BED3955139B4426E7DC8AC173FB40DAEA7CF5BF5A535D4DEE360A4EAD6C5EBCED16E1C6DFA28693921EC4BB9F9295A6B
3903AD689BCA7F3D887777B0699D5251C31F2FF1A7D7EB063B614C04B1EBCF5D2A647433701C86510C4630129EA44F0ED119A804931042364A485F25
4962922EFE3E7771E3FD554D438D08FDD6A5A2DFEBBF7E5A4BDBA5AABD912DF9EFEA1F55002828FE47FF67DFDE6ABCCED63F3A67E38D87A6C7DB84E9
6E11BB7DCD23976C7C7B8A5E107AB88CBEC6A18A5BF8E1FD7892C34BC1BEE35321FEE76F5F75D9C7DDEBF03D49D134CA532B64F13E4BDB952A9A7207
5367A8D62202344B4DD663AFBC66FC3CE7CF21F25B9E59AD2EA385691C75DB1C1E70D665594F754C4D0DD279B2BDE3C1C733CB4C563CA93872A963A4
86CDE86A07587CA9FCB9654D1881A22439C496CA5D5475FB5973791B13871074FA10202E4DFE5DE35C2126213182944008F19FE1452436C82AF1DFDE
DEF67E9C573102E7A47516C53EEEF02B6BD236B7BE0448A8DD7F0A8A3FA91F875A3203F4D66C3D90FAF9CACB16DB95858C87664B130375E79FD87279
CD1B5ACAA5275F3E6BAD1183BFAFD517DA2CEE1CD9F1B3D4CEFF7C618851799EDBED4F06B33CB628EA7BED31D08B99AD28B7206C86C1182814ED95F9
7BA5FDCC0037D95D67643DCE2DF46EADDA6898B2D2F6A0B7AF7F4441ABE7D396C99FD143FDDED177DD5DA3C7BB4A1A10C18FEF3DE54D02F194348C7C
1E4324384E4848900EC3485F790FBD3DBF4584494B018A496B00817575892162FAD08254F7300C6338419250DBEFC9BC65691D9F2FF7F25E5443BD09
95DD5B9E75AEBE0B6FB3F6D47CF19FD3BF54F3A7A0F8DFFAE7BDB81CB83C22F4C8DA9C175E1531F639F49F9E8B0E273A995F5CBDEFCC2E61C591AE87
0589BACF49B4E95A6B825606909B36707EE1B9AE3B8E2F99097B7B36AACC09D55234DFE8AD1B19A4A2AE6228EB060270AB8B9C96BB86E37E5D975FEB
5D1F6E322D1F5D3133EE9E72669241F0D1ECDB4BBF3496FC88AA6B4B9DBBD0C93D678C56220DFC0D7513DD9DA034C84B8A0FB74930A9FC5192E07248
38EF603DBBEC613B864FAB7FFA0E1084DDDE2F26108CC0A63701A79F0992E02459F1A88C5E9C9A59F9FAFD30D85FCC85B3AF0167AFD1A3CF20B99E17
16989513D4EC6F0A8A3FEA1FEA397C6361FC5A952D06D70B669DFA6E98D6DFB4DA70CBE38D8ADB0E3AD4AD7DC8DF54223AF075EB9C7112BB75A5CA72
610F2BEA6AB387C3B74AAB607E965DD943659D3906CACAFEDB949D924C8D559464034016F84C79A6B292719C83CAE98B361929336EF3E3E46DBF7A78
A4589CCC3BB43BADA82E7620AD9A1FB8334445EF25547FE617CA6DEF47055D629023C1DF1D1F2109181737F7E32C01247A137B756274908D21D2E88F
A2080AF23A7ED573A7CFFD488D8064BA24083169F4AFC9A5F36F247E282A7BDA03099E76C0CD17FB3E3C151C8FE74F789CBCADE4DE4BC99F82E2CFFA
17E5445F734C59267FDAF36CAF75382F7A55EF609CBDCBAE5DBAD61F7D3FBE4B406E1D20CB523ABCF6924457DCF069B312207B6D4BB8E2B9017FBBC6
DEE08B830BB42C8D343417A5189B1CF5709CAF39EB12C685EEEA1819294526A96DBC323BBEDE692B3753DEE248D0925D8A47C097F6677ACBD78F95DF
E66DDFBF4D795E0EB3FC613520EE1D46263E3582082EE92EC048A9B827BFD433409128EFE683C1CE1E0043D1E9EE2F2D0108BDB67A9C0742D3A33F30
E9CB12E95792441BFA49F1DE2DBF854D398300EBC243F1D0D996D11C76CDE64E2069CEABB9AA295C82923F05C51FF53F76686DA251C409F9F501DB7B
C3CDFADECCFE3A74DAD9DE79DB16D99823CBE9FB731A7CDAC8987D2DEE3F49785776BFED5E418FFB856BF23679210AFB272FF80C1E9861652267303B
C253659BBFC92A4DA5328C27CC30B6B50E4C51767960A67537DEADEDAE86C6F6233EBE73B6BEF8E07FB4B36C73C168289638FF82896715F4FE6E3DDA
5F45C77B2F354E27FBE62A92C42578C3F76104C5BBF71C98120D560C40A834E3A31002837C3E7D8C8B8A45A0D40F2028024B7F20F5FEBD4D5C9215BF
A35FCCEC6C1021154F5993BBCFD2F2C75887DE132F54376CD7D1CD02A9F64F41F167FD576FDA10A1B7B96575E881A8F61746F79FE9DDA67D707258E5
792FD2FE8AEBA3BC747EDC39F29751E1AFC8013227467CCE6F42B87775F562FDE3296601B46ABF9C2A2DEB390AC6161B37A92EDA64EDA36DCC20C011
3F250D9F172BE4CFA4AA2C396D7EEC86A15AF49723D6EA9FAA1D8C529905A959C2F0BA9706C1264942E685EB832C4EC31401BDFB89E298A8A61222A5
FDBFB9560C6303FBBCB310B2BC0142006C7AE11F12B3384C212081C5202A2D0708828162692620B1F25A9CECDF7984870B68B50032DE0A31F71FA87E
D40C753DE00C38199C72D8E0568711D41B4D41F1C7E5BF97AE8B57ABBBE5AE493ABF20F7F7BAE8C79A914DBFE65B8699AEA56DDFBC76DE8F8BCF9EFA
D1C840537A698680B9B1B6C72357F8D121EB8882FD6957E56B1311CB4777CBCDD152B4D910A96A156A6929E728C4E06E4B7959DFD3F2413735EDF7DB
6DB9ABA9BBFEC3312BB94CD10EB5C0FAE2CD0FB9D9FBF375830D3713E8FED866E077C508417436E2B8A8A1AC8BC024E464D50480C26FFCA22B494E7E
3B46007C9E78FA2EE0D1DE091047A5BE1F026069EB87786C0824D0AEFBBF50B226E2368C4806BF4DA07DCD02E84BC648CB372E3FBB4E18A512E8BFED
CDA26182F2FF14147FD4FFD46107EF40FD80C36A092F5D9F773E39DAE6B7A4B229D43C6485CADEE31B33169E6ADB7FF3DC194E9E6612F4F509F4E406
702753C88CCFF8E1669AB2DD6AD5C803CB3775C63616EAFA0B36D9B96F327550384EE2789FADAAB2AFB5D96B1FD59D67030EBB2E08BB7F4057ED1CEF
8252D8E3BAC32F06B8DB5E46AC74D88D4EC5C5B6F634FE66A350E3DB3E6ECDF7561142E0FCF23A3E8E0C2545D748243F5E0FE0125822EC622100BDAF
6F48844B305C820CD3A57E00138E7125185A7FEC918828DEF18A20F9BDB57CA2B715813EDE9F1A681189AAF2A72EEBDA452C6DBDB4934F2DFF5150FC
49FD84A426D8C8D65969C56393253767A7D617EE6A4ED7CB6E4ED0718A51D0760AFA10635BFE715179DBD7B1C3C61F892B1F5AE347682BABD02F610D
491AF6BBFD35EEE6B839D59C51355355D00F7333DD61A6A1F88D84C84133D9591E9A49476778DC71B60B4F4D0E3AE534730B2F6356E0AB8CCD8786F8
875C539D8DA3306658029DF7BA6A844E2B2C148D9535802826C227AB2A204CF422E014888F9C7F3B4102288030C62194DDDF3E06FF67CB0F9F6C9990
DA0104608B618857DFC483EA53DF4B8B4E6F03404E8808ACF1F9F0C0872968E066E3570BF3348BA81AD33BD3A37FA9379B82E20FF6FF8DBFCF6C2793
A08EF52BB317AE6E6CDEFEADC4EC4CF3158BF92973940D821F662F4C9B083F2D19AAECDCB155D87C01B8FC80E71AC6468E7EAF7331BB98E990FA2B4E
F758E5ECD93A8ABADBFDB5AE062B584E1020F95A739666F286E506BA0712B57C1FEFF0B8BC5D61736F82B247F915C780FAA91CA7C576DA21E0F89A8B
A2BE8266C6ABFBD5E320D84127081682D29A9820D61EEEF1848732DE154BC41084E3D09408C3C43C0EF09F2B3F24E8F8A80885210C4130542014E082
4F19951262BCA89D24713149B6BF168D36B38089AADAB2C5FA491956EF520DCA494AFF14147F80207969CB6E1B5A1A2DF9B1CE2A67B1F1E7C1C0E3FD
01CBCAB3ADB463C3F557ED5A5296AA95FBC3BE0DED2DFAEA76073EFFEB83FB5092F215223F61F4B0D28E1C6BD7AE74AFB96F2D64CC95647DD7C8EDD8
A8E84A87F1111F1559E7CF8B6768274569CFBDB0DCE2FAF95981C3894A964F5FDA84E4893F5A6FDCFA97FBE4B8D759889B5B06E49F2EC3495A0F8F40
F92C567DBB00C77F2C89AB82C6D9B00811F70F0128CA1B63C1203C3DFF6BFAA9A0E9EB7E300886601811417C4844BB73A78140874B5B09090B25F1AE
D291E15A04EACAEF9CDAA3EC71D9E15287D9BC216AF40F05C59FBABF84EC0FF2FBE46A6715D9703BF8FB8EB93747938399E75DF27316596D4AF33D58
9EFEE58D4DC0687834847D7F786ECD48CFC59AF917ABE7395773933B7E5BB97C4C72293ABC6ECDE12459034DD5793B940DCCD5F6C002B4D0444DF55A
BA9C7250B2B1C1812415EFF31ACE2D0714CCCE7E981F56379EB3C2277CBE6707CDE73A2C6056747E2C18C38109068E81FD633D1D7D2C02FB1972071A
2CA3093162727C5400E3A880258241080520707A0430313D08009A9E058A82100416A43C11C18C29904F8CD048126DFA8670DA60A0B2B16EAA788EF9
E9E46DAC7766C930F5E82F05C59F0A002E2972B0CF5AE9E5E476B922FD418A6B78E73B9B9646C7732F42741C33FC9DEABFC414C6CFFC52ABF312129C
7A1E750848CCCE5ACEDBA6B51DCABDD3ED3673FF2DB9B82C956D81851BE535672904CD95D39C751E9B40EAE6CC30BBA2A9BA75B5AEF2F163DA2E5BD5
67DF0851D7DFFADD67711BFCCA2220D93A813BBEEE298A8F5C2C1E2C1261F4463E49327ABA011EAF9D87666DCA9FA87FD12082787C0E2A6DE8089F2F
4484108C88797C588C483010E04B9B3F84887842A03723A144883121801433B9243955D4CEE38F93E2079FB9A25ABFB9713716BE6D739DFD95A0DA3F
05C59FF44FA0B7EC5CDE6FF072365F969D12F83466553D6DD153E1E6D067E1DAB64776B96ECC70C9AE304A86A2034548DBFDF298A9C721CCB80FB9F6
0BDBC85B39E98661B92E7EDF5C820E245D9B2DAB26EB17AEAEABF81167226516B2B676B336241AAA1ECA71D7D9E06D72C66F9659F8A724E72AA8CDD3
6873E43D317DF3271CA2A51D9E1AE1C3EDE328868C76231C3E42BFFCE672423EF7479B34E0637C0E268D2788D4E3632884226C0620FD1E23C45C1044
5002000021D075F8482B2264831221838E1064EFD726422224B1D737A780810D8B0EE65C5978F99EF5CAA9E9016294FE2928FEB7FF6727185A7F586A
E06BE9F830CDAA39C3EAED48623070DD72DFE6D99AC1591B03E237467525EAD12BD53F8B814B3FAFDD14AC3F93B7B132D838852C0EAF58A8F172AFC2
89AD7F3D3DE76227AB2B6390A6364BEE0BC6101D99A5A12BBB2E4547F54099B7B69F85CA93730AC6C1D9A7751F4926FC4D3DAC1E4293EB9FA0586BF4
C99E91AEB1092E01A3A35D62BAD4089C880ECBA0093852778F8A85221881241291408442308A725AB8840483008EB414A0121410428C9EEAD4759F04
4C3EC61F1C628849B825B78524A556E2C7E571A826C8E5DDE3BC396BBEBDD4390052077F2928FEBCFCD7BD6ABED367BFD91BEC973D7CEEDE56E77C71
F89BC358B7E7AAD5AE3A0B7EA5ACBF96169DDAE0FD0EF689E08B4A937F9F62D6045F7DDB7B749E732F1A517157F762D1DC2B79B35206437D558CFF32
DCAF3E4BAF0D17F142E5AC0CE79D9B3BD3E55384BCD71A8DFD253A1A4B2E3D774A4338F13A5BD62D1B9D5A7908868A224E750A7E350A086264825B52
2B44B81DD723775D1B946084040760911812A3124CCC47700CC6517864089B9EE009F025084A822D6322FEED6D89B75AC50208105455F03004A9CDEE
27499C2407DE76A34DDEDAF74AB223163CEA7E62F2957AF88782E28FE0C43BDBE4803B4B8C826D83127E99BFE22C0FFE4DF7B98B25D8AD76D1D6C93D
6E9EBFEECEDCBCDA14F60D8D6F0CF6EE4B9FF64E1525DF1E7DE1A69F44DEDA3530D7FE7BE2AA6C2D8DA7C50FFDFE9E25BF4949C67618160E792BAAA8
EF5C2AEB1418A66EEBAFEC5B6F3FD322FD4950E0B02859DD7BB9C54BF6CA1D7CF8CD8A6BCCDEFC3A018177B7619D6F864138FFE8FE9B27072542A104
C5F82C8690212030D684787AC51FE20B44F8F438008C27003012F8D285B7EFDF787510436036C42CAB66A1385855304292104132BFF522B5EBC26FBF
BD7BC1EACDA9D474972E6AF38F82E2CFFD1F4D553BEA717CE79CC54E5B63DB5D7E2089A1D9EC5B91926F7357B819EB9C7F6296B5A62065A5F8E33D56
E05E0158973CB42EBAF7D7D642FA0A2777EE440C7048FFEB03D9B0A572015FEF9D54519BB146EFAFCD1884B296696BDB1DD13039B04565F66A6BCD5B
31B2AA9BB2D70757B1F65AEE3BBBFCC6AFE56BB9E23B2B1F8E8EBE6D8330A2F73BC6FED60777DD3BFE23FFFB380A100439D546E70AC4420403A7F8B8
B41A882739FF39EA4F10E8F4E55FFC2F23447BECAEDF188E09614438D02DC2006E4D0357826312126AECC26ABDE3075BA2EE6F3F516EBD7FE95616A5
7F0A8A3FE77F71B0EEE5459B0ECCB3758A0F2B337DC63F151BD24D5BD22A5E69B374B6AAE71793E0A4BD25AACFC198AA6F4BDBA6C48FAFE76AA489F7
2D17DD36B5B94A9CA8CD534E7BA7B9E4A5A265F2DC7B0B64FE5A672E73442226DA5D55B576AE50F0488E9C1DED21E39A3E53CE2D71E5D21A619CDE0D
30259CBE6B5E2F7C27EC3550FFBA171512632562714E0F5E7AF259EEF7E7ED30974F4ABABE0F88313E1F12F3F902094E4093C31C58EA0224048A8228
0EBDCF16498ACEBE9CC2710044201E2844609455D18A60048891706333DEB9685E3DEBF2FE6B097D114E39EA47206AF58F82E2CFFA6707F8FF3E1E78
6CA1A5D3FE902AF79B604486CF6DECE27EF281D58A39DA4685E1D66F23F23DEC193FCE8A4EE68804A39B73934C4A7B7CB347169A06E09F62E8419BBA
B7F854F9A89FF3083FA3F477A8A3CC0309021CD3D29E73C4D4E15C92FBBA6366F23BDCFF323E161B5228786C1036FC6577D1A145796D6777163206F3
69080FE75472899F8D48CEEDC2DEDBAF011C919AFDB6FC7112160900F63013815036BD6F944B60D3C77C71088571C6BD0C70F8660E4F8088201499DE
19801178B86410C64010C1495AB3683C5CE701A320A9FA44E96FF5C737B5DE518B7F14147FD63FD1B550FFEBF3B9D1CB8C6C62DC0BE34327B6BDDCEB
383C1CD84377775EA0AF723A4FE7CEAA03EF544EF11617572D6E07C1AC4D4D01DB84371CAB0E6A1BBEEAF6E83F6C52F7516DEF911969674C73F565AC
4D0DAB50A1204347C5D15FD9F6AC9FE1DE45329EDB6554D6AFB57CC8CF364AAE48B2B99335FB6473D2819F407DFE103885317F4EE25D3F8437EEF1AB
2F972138464779BF4BF9D234CFE5831383428CD5DA3A3A864A9110E34C1CC2A18A02A6F0E1BE3631D83B89E03001320108C2189FBBA41E402C22C8C1
CA095E9ACE1966C7E99CFBE9F42DDE1D0E763D54F3A7A0F887FEFFCBD5F0DDAB3987F61AD89F5E5B55153990F9B4CAF2B5E4D401F298B58799E626EE
DA0DA7D6D2E316D2AF468AC3024621666AC1EB4D83935E21CFFD2C8245BB73EBCC2F4C2E72FA6CB8E059C8030F194DB5650254406619A858B8EB7B24
B92C3960A2BDCFFDEF75517AF19C5E9F951377DD361DF63A5878E158C350F5FB2E0445873FD151C1C39C970F2AF3CE15090941E740FFAF2111890BA7
4410C467D119ED9D93D2D68F48937E7BC11881D65EAC41E0470F7BB955ED43220CC7D1B13E0043473E75202824067164A8B287FFD42E758AFDF64569
74578D5DF135E51811A57F0A8A7F58FEBB6FA57629C7ECE21335932B812F790985755BDA77EC14B72E1A6F369C6DA136BBFD82EEE3A5D7BB16ED1B5C
74AF4A79A700BCB5AE7A550A70CFE6C24E03836F0F978F795A343D334D7395DFFA60AF858CC2DFDE2C113AB6465579DEEC3921B60EB10E729B56CAB8
27CC9C37D81E1E5F97E2FAEEA26F44C9A9A027EDC5F7DB1010EBBD518D42D7F7BD3FF7FCC22B1636FCE5FDAFE1A6310986331F95211C11305AD73E04
8018849312C668E728D49591528357DCEDEABC92DBCAE06010848A790214A6D7F7A33000036264A4B2152EF7D939C12DBC57B8BB92BECAFBC63CEDD7
D4D15F0A8A3FF4FEE94F6E9CADF9831AE7FC06758B0B8BBF8AC2AF8A2FBEF8143380455D4082B5ED676B9D2B30BBB2771DF82A70E4B55FD735838748
8EE1FB7BAE4DDC6D710F161926D2BDCA8F68DE66ACF5DA626873F7A6A1A28EDC261840BF9869CF7676DDBD52DB33C1C4355CC9788FB9C6C3329F79F5
1775337EB9873FCF4ABCDD56FD406AF627BE1FAD100BBFC53FCB4A7ECE042564CBB30E6C1C443082D176AA92804146593D43DAD3A5C91FE2324428FA
7AF56976EFA5D39DDDDF0AC7852880C322004370CE041DC620181620C050D3A0B06D675497A8E961CE87B281CB161B8EF8040C92D4EA1F05C51FDBBF
A46BC15CAD136D76B9CDF25AE9DE4F059E9BCAB36F0E1C7802BEF2103D5330B551985BE4697B65411EEFC4094E4412B0DFAA7DCC25B0D671353FDFED
FA2E6593B6CD6BCA34AC5AC3D512E7C8FB94AC53D098B90FE2B3E355D5ADB542D3CD745CACE784A9CAAEB6507CF046C7ED61ECDCAC9B8EA967E3B7C5
D6341D79CBA78F3D7B5B2A107F49CEBE1B71476AF90956DD183A388CE278EF20C4E4629070B093FD9F353D42380A48F3FFDD35C58C27E14F7E0FD251
1407C42828C21042C89814C3388C022C1E9FD6DE098CA48535C35D2FF38B5AEAAFB9F89615E81D1552677F2828FE41FFDF9D63ADCF4CF87FEB53D74A
3F749F1FBCAD71C7DEF1FC4CDED88A4AC85BC7565FFB5EA2D6931D91E281B0E2A694417ED8F6E1175E9519DACF859BE29E39AA7D295BD81AAE909162
9816A8AB77BF4257452317134D84285A39CFDD17A46DE9A2BB7A9592DF56A53D9F2C1DEADE6B9CE8F539321CBFE3644E51DA0B16DC36D00FB31915E7
6F7C4C7F3D800BD1C96E1E31C0985EFEEF220812E032F91C1806411823983D004EA29777B6B5C6C4FF648EC0A404024114816142441BE743048CA220
9D231A6F67F3FBAEC457211DF9BF26BA1B2F442F2A1DF632FD8149484AFF1414FFDBFF4B08F49ADF2DC32D8C903BA306AA01D7A2FA4EF90DA539960F
9F2885CF47E06FE48D8C15D6DDD1DF7CDAF00967995B7FD903B4C2F6143F36B6CD72E1E4399D4BC9335677863C7DA8E6F9DADD39CA5E656195ABAC49
0B2A64259838A807271B68CD33325FACBE30DD32F6AC9ED1AFE1F92B3E2DBF069E0ABEFEB93FF6B88856C7908831D6A39437770E94E340FE131A5BC2
6F65123852528E9010CC981C17C1A81800309CD60A926453F449D6FD65E7FAC6783801C31084A022046DFBD88F112042C013C35C513F8DC16E4CDFD6
0074178EF2BBFF0F7BEF1915C5B635EC8E81DB40CE20493228884451141130111491A0A28801CC888A8A624609E6808A22A062024150C951A2E49C73
6C72E7EECA555D17F7F9CE79B7EE7DDE31BE1FF7C7BDA31EC6A0FB6717CDB3E69CAB56CD39F2F89E5706EFFA823DE3B38B0955015050FC43FD8F73AF
AE3FA36042BB12CC369253B7557AD778B03AD3208499F796DDBFBA986F25A92563F3554FF6AEA76D7B84D8696E7A3B7A7D455AF7AEAF0F14BF761BAC
7EAA21F9C1EBDA88EBD2F238ABC035468AC78DE6AEE6A1D0CCC1652A1AE75CE4F4B597AE58B3227883D50919E168B68F69EDA13DDCB8A58F07DA2E5C
A3F53603203EDDFAF8E4B3A4AB3F103C2FBA4D00F0DA66DD1E4FADC05101DC3FC20770369F8FA0F4AC2A82C46E3BA6B69CF6CE1D9F62113806A31802
201852F96A0247201803BA3BD8ACA9A9E9FE675EE195DCA66A0639CD6B7B76F313BF414B3605A1623F05C53FAF00F8D0E1353757887F7BB6113C22A2
AABDF058CBC5C8BE7D26D9DC5B39E05D27FE37517979C51B060ADEA1FA97BB0D64BEB213194D9B3635647B7E58730A39BFC06FFB9CE3271777854B6C
2DF3F036DB68B054F58FED30172B35529A67122CA7642269B358F9BCABB489EEA61525612AEF6F18777D56BD3856F16C674BF7470099992E3877AAF8
C9995A8C9FF20581C75A073804DEFAB51B2188C9B62158C0627066CBFAE6D4325C50B3656F47E1E1243AB7630A434104450114C311D6108821204440
3D0D7C627466A223FDE89D615A630B8EF27066E893D8715EC01F9BA7A8F04F41F15FCA7FA2C5C722E5A3DAFD0C7B6E9ACA328F33FB7292A3692986A7
E9392F59E31B3FC17EA2AA9267B76905ED397B8A7955D2BEBB338316E776BCFB9AC7493F7E8F855594A95DD1E2C07C8B8547F63A6D5CE9BD4251F409
CAC7A26594A58F6D90325CBA64B5966F8BC7429D17CFCD5729BE4AD0FF5C6FBCA13BEDC2CB02B8B409AC2D29BF72212FE6440D027D4EC4C8D6310CC0
D0DA6F6C1C820BBF8E09B8335C00022A0B9B5818E7EED6C78D57AF8D123395E3C8ACF70007C1499C3503CC8A0E72597CD638071E65F6BCCD6BAEEEEE
AA6162240310C45DC92E425A95240BA8A3FF1414FF1CFD67CBFFEF2E5A973A2C36D7AE2A1F309DB7E78DD3AD82D8DA4637ABD723419FE0A455DC1E55
99799B13A4B7F87A7B24D42F928C2032CA3EBFDF163C7CD0DA2C060A14093C29FCF2AC4A5AC07CB3632B4FC8EE721196C88727A047CA22F67E22CA12
CBACA5B51B3F2BEB9C4C5B2FBBE8E0479D27BD6B569484EFCCA301DFDF4C7F2ECFF3BB5A1576A00D03DE65E0FC0E8264C1486E0617C5E198FB633836
CEC405436FD2193052BBEB74D7179787083A5E3400A308CC98D55B80D246410C8441EE3487CF61207CFE54FC95F4A19AB2090063B23032CBB7E4692F
1A3867E7CFCD7F2AF85350FC53FD4FA0A9B64B0F0DF95B7778BEE5FB88EE7F13FAA4333F9FF6D0E5C2E087EB1CFEE67BF86961E1658DBBB52F5CF5DF
3C71567A650350179F98B03AB17BBBB6E34C99DEAAA7863629BA3E6F35B4CE6D0C5CB271CF7CB56E7814FEA82FEFA4242BAEE2AA23175C6F22772C7E
85844980DF92E33D97B766BCDD558491C33BDE8F94E56F3DD314B6B71D47539231621C27DB5A58A55F411C1D7DF2924532FBA748B82DA51A44E177B6
9F7AC29C3E237059096DD67E369D8109484EFF0C8661B3AB01C0FA3915109FEE2B0E8F6DCA6F9E4198ED3338D1B6BF32E7ABA05C41B5928AFE1414FF
18FC678B621C78626EB1B2EEFEA2A2E341D317E7187EAC0B6D6C7E3DF8F6F093B7CD7E5970AEE648BDF05C915B79BEE11B8FF95F6BD395392FA0551E
290F5BD191A6A715DDB54264A38DE8ADD34A717BE65BEF73DABCD4EB0FAD612E7DEAB2BCB892B8A8E2C1757F58146F15758CB05E68EFEB21B9E7C79E
D5B9CFF6B4A11870FD4E7769D2862BADC1A7FA0934FB2B820D20646116DAFE7D0A06EA1E266330A36F06A39794F7E358CF3E9BEFDFDDCED69363B90D
7404634CD1015480D386405C80202802C108CC42C8DED8B2DCF4C9BE1600A25574F3F02EBF94B147FC89D5F3CEC254E30F0A8AFFE6BF60D2DF6CEFFA
8A1493BA67BECC1445BD2B3FC2EE76250FB7F8BE7FDB9812CB0737F82176F3846D0B42B3F69B7FBBDE1FB1D0F6078F13727F4877EFC821D558F082EC
061B59F3BB66E131626AFEB67BF43DE7DAB238609E81B0948C94DCDE8372664F83C52C2FDA4A6EB8B14CF150CD59858B9FED73716C34F9235418EC14
5E773E7C8CC04B725062BC8BA8C800D83400A7F7D40DA2D8C01048B47EAD0305CC280B8FAA48974404EACCA807510673928B10C454DD344EA2184EFC
1C0208C330417F95C4E9EAAB1F01215A490F088D5D794E64B713FE9286BD54F8A7A0F82FFE0B08A26DABFE29D797B546C565AE3DB56AF20AE1A9BE1D
1949335EFB8AE36ADF77A10952FD1173FE907EBEEE5CF5528FA2174D0E863B66808A2D4DEEA2314512DB59F96AC78E4A691CD968136F32C7D5DE6BF5
EA3F0E70A7C19722A20BE7CAEE775251BAF04E41F7DA3A09A76B2BC5760E47281F7CACE9C31FCA8AEEE3371F587FA7F0F0073E8115E562042D935B9A
850CA48D73000CE662033D9338921DDF89088A9C940FA65FBFDAC5E9482D1A4027992C1E8413C4D0C71912C3711CC3501C4179303193FA7DB0A27E6C
020147CB3BD9202FF231D4914ABC53920CA336FF2828FE6BFCC7E19C0D96F7829EB76B1F1F395E4A339394DEDBF2A2B42D7E38D2213F3DAFA694CDB6
8FAB15169DBFDB55BAEEAE5CD2D7AAE3168BA278E845DF3009E39A1DDAE5DC1D8E39CB15EFDE507CE92564B2DE2D70E9FC57300B7A20A1AC22BC2E40
4ACC2D6A8D98CF4EA90DD1EB64036A9FAABBA73BAE297DE3F87688A8DD6A1D5113F80526D02F29B8602C89C12D6615A734C038D0D934401BE2E333CF
DF4F11C8F325562F07BF65B1475E7C1A60F533A65004C7097CE81B4DF0E714200C24081806B1C18AC2BCEC111402C181D23E2E879719C50092B8BDC6
32F623D4CE1F05C57F89FEB3F19F71DF7289EFE127954BD794DDBCC30EFC4361714D4D282F27BD7FEBE9B66B05ADC5E88BED4CB33973D41E4AED68D1
34AA4CB91BB5D16000FCBEFC9D8BB04F98E4A6B114B5976132EE99EA3B4FFCA1E1A61DB254E2074C63DC9497986315286D28E27A7E81958FB0CD756B
CDF09A938BF6C66EB26DFDEE9ECEE6B6FB2EBF5E79310525C0D84FB860E0F30C5E56F735610643A71BD346D071006BBE9306125301CAFB0746E2DE8D
954597217DB5035C02C5088C3BD4324160C8EC7B9C4547D1D997B6F4D277F53C920581ED65FDD31CA0E0F514F4AA99DC23A95E443DF84741F1DFFD17
D0FDDC37BB3ABF6079BAE6DCB83055A5A6AE5BCDD9553193CC8C3FD5FEE4E54C5D6F7FD0F85B616989507FF3EAEB5251F9C76ABF1BDCE5704E1E8A51
3589DEA21831B239304347F3F371CD27AA9A27E5B6192D1B4347A6A396484A781AEBB9699C73523960A67BD65E293E7783B447E696650DBCFBD9083A
74D8FC5EDAED4C5CC07AF01E276A33A7F1A6F0DC2636CA628D4C72A16E16EF4B580D8A7EB7D7881C88DDF38E599D38804CA53562E46CEA8FD13AE828
86FDEC034CF24798B3D900372B272BBE03277190DF534A03784043F6149C9A44BC9795BF0152FA5350FC6FFE3B07FABB6F38C3F25B3170E4DEDB9E35
F36422A0DD1E404BC6F8D5ECBA0329DCB892EF994CCB797357BC5F72B2D9CCA272ABFFCC2DA554F0BB65F466D17DC7357532A3CC5356897AA5AADCDE
20EAADA524BA89C186DA0F2A2AAACB291F90743F34CFD04CFFA8A548748DB9D496340FC57CFED9AB9C9EC90FA627BE5C4D139003A19F497676254AD4
9DFD3E89C0D0001D12B0D2BBBA239E8CA03D8132E6F9EDFB7C9A90E6AF4CACE7ED0F0E268048B0630099757F36FC6304CCA6132454F936B3F4610B89
73317CA899CE03A0A2AF6CA4E1219FB752781D8DEAFA4D41F1BFF93FB3D1DAD8F3C2F1B18F963D612FBEA7FACFD7F6423F18152099CD1FA2C1DBDEC3
0587E307392F543434DFDF726D7A6B9817A2F0A544D2B20B387628C9CCEEC12EFD5D2D3E978E2A187C71DBEE2DBA799784D029741A295D2765A42664
E7A97DDA42CE66A9FB1671FF4E67919539DE9277D18BE67DDCFECC75FEA54F5E6164FBB95272EC45152118984D0A8091818949041C7B9B51FF2413E2
3F5617F19FC8D8F6808756E441707D62090F07D83C78720A13CC16FE38461038000BB875699993E5913D28800A04CC7E108358F1EFE9302D638CA893
572EFCB9F9472D001414FFC57E42405BB5C5CFE2AA65F1A4EDC3EF3E938FA29455352BC6BC2FE58E6756DC2FECF489E6C7AFCB19E9725693D859EE7C
AC6DBBDB779DF5C5878583D8994B9EEE973C78698D7272B2DD556349DFC38607E788BEDA24741B634EDED7969490563197F1B199676CB9D861EEA616
2F91C5D101E23E9CBB764D35CF9B7D7796DCDFDF45369F6F269B425B0504E773391F6E8CCA6D9C845B1E15D35E3720349F052B8A80771E85385E5DC4
038B9F16A0047F6C0C80E13FD54731F4E7E06FA4BAA0698253993D06C33C0C9B6C65400894F69105732B265096FB829B304150FA5350FC2FFE8F98EF
7B64926C130A7CBA3E16D09B9A7F58654934713FC6779ADB96D70C7DB88CF082B78C434F9799AC2C8C732BCBB7FF76DB3CBDCA6AE90FA69F67AAAD5D
D82E83C393C78EEC53B639BFFA8884D0AD0B73B23178E49C9AB6F29A2DDACB4E2B4A6F586BA16098714264C9F3406917F6E3657DDDA70A53777CCB70
793CDC70B29F6C0D6B27C9DEC82288D3189230D53E355958D61ED7860F3A899F64324F1EEB13903D25ECE1D791151C8CCBE7F17001F6E70010884050
84647E29E5D7B78C7670081E02321ACA660438DADE06239331FD5889B984CBD47F9AFE528B0005C53FFADF6BBAE1C6DA4C57FD2EDE9BE6DB950D919F
3574FC8027D5010FC999EC5668D43303EDB37C08172DB1543E5CEBB26BE8BA6BC11A9FD23D22C7B9F906378FAB6CF53655CFCE5DBB4FD33044DDD144
48DF527D00600CDED09176BDBCEC8F1DE6F3F5D69ACAC87F3C25ACFFF29AC4F2A1B7162D4D9719F55BC373BDEE4D76F8D490B5CF8649E2DDE94C3E32
FA2666686618A2B77FB9D1820F392E8C02CA1D03A731B2FB15ADEC6E5C1D8632F804301BF1510CC10100E4715192F675881E99C6033178362DC07834
EECF83017D2832EDFB11CD5FF8C7DA2EE2CF9EFF02CA7F0A8AFFE67FD70ACBC4CB8921924FC1E2B0D725D093AA0B861B473FB427AD68C76B937944AC
3B84BC36AF993CBCDAD2F2DD6DFD689A53D4FDF5D94F4CEDDBD901E6F64B57B8EC553DC3F2DA7670F54D7B2BEF794242960032CD8D31D1F0DE20ACE4
BB50D1C164B1B0F77D298963914AEA9F134CF2C79FB3EAFDAED55EBE384ABC7F4C76DF6290D0D5C3F520B726BF64B4B71F60673C8B6DC2BB1C967CEE
8DB0BCC6C0B1CAE733C5CF1B1082984DE5F91C14E4E33897C705B95C8898CC9B1E0AFA2C204114C2669785D96B81A1AE5614651E0E474616CD5D538D
0BA8EA9F82E27FF59F68B3524B7AE1F55E65713B23E075C4CCE7874DAB8C52AA62DBF4F6319951B5D884E10B747AEB7EE08DAAB3A665AC954AC543E3
AF068E99FBE57DB9990EDB7597EBEF5DACD67A7DD11D139F630A1E0A4242663C1E6BE2BA868C95AC90D32A212317D90516B715A4DD4F2A1B475E37CD
EA7A0357AD9E2D21F6F620A31E0DBD17FB48867F24DED8FFE5552F3431C8AD8848FADC8E55D85A7E2D7659133E8EC07151D0780163B6DE4770820F10
388420009BC3C7B83308FD07BBF5669E0083381002033FA78210F4C22C168A56BD0110AF056645B32B86E0DF337FA9D9BF1414FFE43FDCB45EE5F217
CB8CE3524FA184C77155CC636527F52ECEBC1B08B7FFC26A2A448873AB397081EF78BFE7B6B53A672FC91DEE5A1F1CB1E18CD322A5A49990676B8D3C
B638A905446ADF76DB7D45D463CD3C210790C7AF705432B496533C2427BC659D88FC6D5FADD3C1AB349F3F72CDEEF98A94D9995D8A3A9807E041C7C0
8BA524FD6022595CDA15DB0220954D91376347721AA3ACB79745DA6D8EA7A163179FC1E3A52081E31881F341FCE7041018015180D933C465D30B1207
719CC7C7502E886308820285B93482C79A44D0570B34BFF188FF73897FD60094FD14147FD31F27A6F36CE4ECEBB69DCFD2DB3931BE2B31160B7A90A3
B1A6BDE4EDF4D1A0324EC210F15D2E15A48744D2DF2F7734520D7712FDFC5CE7A5FBE975DA0BF5FAD3221E6AD9AD75B2587C64E9263F8793F3B738CF
113AC59F66BFD3115F2126E1A52F6EE9A820E975C8A6F485AEED8323BB86C6C239D35B165F7C6213C646BF2EABCD7C2218399A4D7EF802BEAC06A1FC
D72FFCAA9B6E3CDDA171ABE9AE91FB573E3A7DE3C114BB8D47FC3CE9838D8FC2F0041DC76755C7464BABD97C564D199BE4F32101CC6702180413C0F7
1F7C90D93FC840CB641462D9C47FFCFFB3CC995D3AA8358082E22FFACF5AC1E9FDB252C5A0E3B6CFC499C076CE8DC8244EA5FFA0876E32103D5E782E
8DFEA34100D8EF44F1CA504EBDB9ADADF2E567EA2174DF4BFB932FC938E94476BCADB4545A1EE0617A6BDB627FAB20258F63F3658A6106EFED523D6B
D90DFB45B55DCC447C2EC947A52F31880C73692623EF4F5DB209C87D1CD00E77983C0263471A0372898F5FF0847C14483FFDFC663BFBDA8E8D76C9D3
27558E754048CBD54488DB3F1BFD6745065ADA41707A8485CDBEC7269A5BC6D80D4945B3853F820B4818FA7304205A5D0D118CAE2100ADD095B94DFF
ABFD029220C69994FF14147FF19FC005C060C387E57AD2EF8AD6D4DC38D4DA5E7136B38C0C4C7C2AEBC5A87835F1E47E25A37482BC235F074C452474
84992D56357864BA6426F1F6CAED5D070D4C17575E4FB929ACBCF3A6DEA66386FB94EDD798EE13D2998427C194B57A327E97E5944C75EC9C0ED9D8DF
325238136C90447EF71F7DBFF17EF503A70FECB1CDBB78E537F222EA5B9E94A31FDF03FCA7DB23DFF471CF3A6EDED331BA53F1DA24C67E1D98C30468
33088E0B04F4DA8E093E0F8270982071E6682FABF5EDDD1206F6537F144167633B8413EDB5384E6F1807B16F2A0B2E4FE3BFF82F10408D0CCA7F0A8A
FFF19F988DAC9C9AE6D7C6163217997B623FD97433072E54F491F9175ACC16A7F2636BCBC3BF4DD78C933D3AE750B8E4D2E0C055F5657A59D1324F8E
8787CA67954A98AA46BADBD6AC9532FC7641FD4D9CEFE24D9BCDF6CEB1E2A2FCF13B46D25A45FB458DD494F73BEADE0858A4E01DBEEA3467DCAF9ABC
BCBBF4FB9E081AFBD8BA41E6FE9BF570F5FE214142F424F0C835AF9C3D73446F6BC84CF58AC5EF51B8DCEF6227C69EE1E1288103FD6D1DAC693E86FD
BCF787728618BDA9C1CF5A2014C6481C47678B83D924809CEAC00450FD3046C429881DA561BFE98FB5D42104F5B55350FCC57F8C6024FFC8B2D65FE0
CCBCE43CB2FC24C44E08650C62178B0F483B77673FEC490EA9AD2C61935716F5323BCFDF87869C7414AFD016DB84AE1DB73549125750F58A978E0E99
2F73205DC66FDC4F64C7BA65D642261310F86EF76221C318056D6361BB55CA27EF998BBA3C596A33049C7E8F7FB5DA71D1EA40F7CC15AB1AD615E72E
F4C7D6462CE626137ABCA9AA3D65D46D61F8838638758B36C1F065AFE720084C0E4EE338D69ADBC9981EE562288CE0C0E0D0F478796474290B41F828
FEF310A000C7719264D10402BCA10511C44ACF77ED44F15FEC2748DEDB6E82F29F82E257FF675E9516AF58A2B162BAD87EE8EACA16C66444131B4F48
2CD6537BD01659D67D2FA1FF5EA760D22A189DCCDDFF9593EDB66A51C445DDAAC3B1457A6B74644D97D7DA19452BABC9DF3F26977845DC3768F9A279
FB996C66889BF17CDDD5626BF5B41CD4F77E3294317B71543F0DB8B3071A5F6DE3AFB7BF653ACCF40BF866453ADC68FF198B393B8AC63A7CE76E3BEF
B1F8C38FEBDBE50346FBAEECB932042193233D4C141BCCCE1D00B01916F6F3C42F505BD13D9670BF9607C15C0805B93F2782CD42A2A3CCD98B1AAE00
C94459C9CDF588E0D7F04F08AAA3E9D449000A8A5FFDC7C6EFA57E56B7BBAE5635B62DA9DDEE7636F2ECAC00EA3AD7B35DCAA9272E8456173CFA3D8C
4E162DAEE1D6DCB5EC98715DA5A872F1F2F9987DD3A1AB7648A92884864B7B5B5A2C567B6FE818AE6073D6484AEC3306B7DDDD2E27A2A8B6447DA995
FCD60C4751ED6BA79522F851BA35C02E7DAFF5278658B7973DE266DB3CE754EF48E23EF519C0B3CCE2417FDD15BAC5559BB40CE3A0976B2E77B28999
D181091A0FA1955433493E870B20088A718B0B69B5E79ED15064860D831C1E8062188E1124B7767AF69A187D7CB244497CCF6C21F09BFE047036994F
95FF1414BFC5FF81BBF15F57EC6A5A1A09474441E7439E8E767AB609B04F45C92A762DCD77EBF809A5F8FDAF02FC92C374639E730812385F597179FC
B9F3FEDD9DC75EAE1297DB5BB06AC5AAC54785F75D5309DFA8B86E91A8EA10060F361D149757373159E4646416E5ABB0D0FD96C5AE9952ED20F8ADC1
099F2BC3DC8FCBAEB0CBDC1FF3EB9D32919BDB7AB18215A1E063052BB3F4C6B5921EEDFCCB2EF990803DD436CDC2A0C1A1293601F379300423284AAF
68E066393F06710E97CFE7F2400CC7919FFB7CCCE69FD19F394990E55A3237867FD75F8013255E3504D503FCFFDBFFB1D49FE0FF05FF9B032FE6D858
776D59CBCEF5871E5D4B7A80ECF181A1E1F0B66D5267D9C9CF47BA9290328F0978CCF8527F4AC6F2AE6A39F9250BEDAE271F8D603C3CF17AA1F8E298
E54A3EE25BD7CFF55DE2EB2B62AA2A6C320371C7BB360A292F96375D6FA078F5BCA292DD690B93A626A375631F2CAEB81C99E27F32F41ACB75BEC369
DB93C1BEEED401A698044CA4E86BEBBC695921EF3734B12770860090EEA4761486F33ED01000E5417C2E1B46E1B6AA61E8E5FA640203013E87CD84FF
6C02420000B79B3E7B416334824C57947FC3F97961BF867F9CB3FB5A2DF5242005C55FFD1710185EB9D7F3CB16D3B6273A1D9CE0BE96A0869BD3EF4C
6A51E4F9E76FAAE6E55DCF5BE1EC04C6E5E710F2D228F3D9E4E9209E8BC812235D8F77F13BFBEBDCF3B748AA061E50F2D192D92DA7657C284ACF5853
643BC687AB72ED4434B474F6D8AABA5E325336B9B05B239AB34B3E25C7EA4EC0EE364EF262DB9ECCADF7E89DBBE3E16BDEADC31775CE0C7DD6D25973
B17283CEE5E9DCEDF778301FE767952382E9C731E328CEE640080702F1899C3A68E2847BBD80806657063E08421C369D0B30E9008B470A700E8D249F
CA6AA4C182BF83E4BA163551FE5350FCEA3F0E7CB030BA69299351A9178FDE7F455E7B9BFB70C8CC8B0DD51F2ADD27E539FD25016CBB491B3A59CB9B
3C7EFD745ABB737195BAF8B22D9EEB724F1E65ECF04AD6FCC338DE40D26D8181B5E252AD5B4765172ED885B1E027072C24C585BD764BE95F5B2D67BC
FF98D46EE63DC9E0F7DA77EEE956F3228D369454BB46B22A77BDEFF70F60F5B92D8F2B3BB7D4EDA84FC6FAE5EFDA5E7AA7433088D10BDB31AC23E41B
42421C069B39CD4288FEE41EA4C0E13A9DF879CE970F231802C1DCB131161FC2511041B800895F973429F87DE7EFCF346770F3EB0FAD02827A168882
E2AFFE0B98412B2D82D7CB44D1B68463658FC81F5778AF06AFAF9E86807B51EFCD75520773B8405611DE1EDF0BD6C65E729EC809E1C5C868DA5E5F13
F0CEE147B25DD95E61B507D7A59CEC166AAEB3110E8C9356110EC2D9D07D5B7B5931557F7529CF3D7A4BB7872FD5AB1933752E37DFFC71691CEBC9E2
0D0D79FB5F018D7631F40327F84D8EDB0AE8BBC43724B8079E5C93D173F1662701E36043CE202E28B9DB42E2180EF331E60C2218F93ACDBFB2E62B82
227CCE340B132018861338C883511443594C00243B1CA43695C37FB37F567FE6B533DD67C7A919601414BFF88F09FA5D766D71B35BE890E8EB38CABA
D142DC2EAF4818D07F0500C5873FEF5EEE4A2FADC5FA6EF6717F340025976A3DC3908BCFD0F352AA7B4F59DF3FEAD8E37EE6A1C202F3D4931A9E0642
6B77CD317D2029291C83D3E1CCED2B84C47C6CE6981E5CA4BCE5A4ADE45BFE3EB3AF1BB563B50379F7B43DF3B276BC817BDD9F746CF3E3D7D99F1CEE
DF6D74739BD9F1A48D9F2BCFBF676302A822BA721267FE489F2405003C1BDD510E2EE8CD026AD76FAC8301067D9ACE801114C75062F617CC83111062
3170F2B39AB87BFBDFA33F21C09042CB9607FE3C2AFC5350FCEA3FF1C3EA4EE846EBE5869BF75A6790EF72C896B0BE67BCBB87202EFD4C66988F5EEA
50250CA77EE2C0C555E5FED519F645CDDEDDE3F60A868157B63CD899F26D45E05AE9854762566A9BC9E81D52567AB07C81EA20C1C28A8E18CFB73F2D
23EBBD4DDEC2D551FC08FB81CCD3C762078E6C1989523A5C11BEEE0130BCE32EDBCB9B53607A6CA0C9D9E04D884650FCFA3BAD97D391D915E968D810
06014D2D18C94349012C6076CC20ED655D31E6A707100E9BCD81501487301CC5F0D9C41F826008E27351E4BA98E4F101FC1F823F81F73A9C6D338AA5
A6805050FCE63FFA46EDC96DA3A5166BCF476DD8456B7ECC209ED5BE7EC1F6A8E60311C109BEF6DBA19606BC21A81F2E89E7951C6E79BC9EFEDC0BCA
5454B1BC1BECEE7DAD4067D50D3305B15DA1B2D26BE76FB49C17E83577AD408020CF774B195CD096707791B3F05A26B266245EE660A682C99DF545D9
26DB13C3D577D39ABD6ED38EFA01B5068E5D492636695EC6DF0BD41CAA42B3119CAC72783F9BE3CF344FE1C8ACDFB00063D6D773E1A6C4A3C10510CC
0367F3014C800118323E8363080623089F0DB0D1C95DF3E4C299F83F24FF1831755EBD21795D13F5FC3F05C5AFFEA3E0B5D52DD78CB72C76FBD116B4
2289FFAD43D0903DEA37F63970B600F02C3C1ABC650CEDE6B3E20A20DA3B16FBE43BC6CED76CAB17F86145C38D8F2D54B71F76713EE3A2A1A8EFE727
B75ACD72BFE29A95736D5198A085EE90B5B05AE070D46899A78394CEF76C0DD3420FD9B090CD5F1EBCF9BADD64E3A318DBD3CC8BCEDC7613FB1F4F96
5A26DE74AA1A58E7187BAB1442894ACF5201064EB44EE30290CFC7017C666A14192E8A3C51C845707036EF4766AD47719C37CA9A5D1E5080C9E6B340
A471F55CFD54E0EF1B7F3FEFFCB3EE6B5FFB7CF50EF5F01F05C5AFFE13E8989BF6F763469E8A3665B5B9CE9E13C5B9F8CC5BE8FD11C6A14FD08CEF83
A79BD75D11F0FA58756143505E1A2F2F69F0F1CA99E845830DD21AAA9E7B17E9FB3A7B6CB4D03736568BDC266C64E0A3A6A8F887151B63973F549356
5C606F67B9D3CE4D41F6F577BD4531D7E6DBBC56D957CFAAF4B6DD111F667A68F08AD548EF72939228BD15519EFA030C07D3A48FA50080371D6C2041
B82F7F0C054098E0E3A3AD7D30DAE9637B360B2140049DCDF7D19F45FF6CCE3F9B07C018D8DF38CA0679B1F20B5CEAFFCCEFFF96FC63488ACAD3CACB
97526001D5038482E217FF897A5BC9FB219B83179BE7F21B9FAE7DDF564887F2CA191BF2DBB674605F0E545DDEBC7A009FEE1D484C0169698C8AA2FE
0AE76F90F36E8E8BAC81CD0347E5E88B1617F61A1AECB1F00C57D7B63E6CA022FDC74558C02A8F959390D57333BE60BBC65CF25489B558509AA1EEDB
B34E159CF675CE216F5EAFF11A885B5D33B47AE9B73BFA1BAF9DD378CD3FB4EC7DF5088B03575D6815E013F5F923080EF16192DD5ED28940E58E4BCE
7C63B10162B6E44710F4E7B97F8C0F62388CC3EC9EBE29A0D5E10FA9B3E3FF9AF2F537FFA1C69541C0C9D2E0514CF0EF4EA01414147FA6FFD8673DC5
0FD9810735657CD0F2D2ED9E9DB9B9704728E79327FBD2019875E851F451BB23D344C3F7D1B01CE8C748B36F17C3D30B2E52F9F6554255F6EC7B4587
B78BD63D5C25BA79EBDCC3DE73F5BD55C416086510446789BF88A8A1AA71B8B786EDBC6D39F6C2CE099B16453C5D5208B7591867DC8E32741EFAE156
C370507F1BAAEE131CA468D71AB424BD35639003A79EEA2327DA2A07A6309033C3A54FF734D7F2D8B11B0EC54F8D7300C1ACFF2882B26BDA80D9171C
463182D53FCC417A172F508BE7FEABCBD76FF6E338D8E97A8415181B7B1BC604542F500A8ABFFA0F8391CBACAA3B2E9BA9A8190F8F64BFB3CFEB291D
01529B00FFD8AEADD964DAA5A890BBA609F04405A32AA26F74080E8E11E4AEAFC66EF802AE5272071B35E4130E193DDD29296527B734405BDF6FF902
61B1760C9BEAF415565DAC71E0A29CCB32A30F7B450C9E5DD23D58EF6C99D1B645E14381AFF5D6AADA83E5D8D185A191EA178E789C77BCBA694D465E
CA109F75EBC0283256D5C1C4B8289FC6C1683F6A39B4E493E7BE72DABBD8B39F98C7C788D98CBF8681A07F8EF4C21923630052B3496CF16784F84F8F
CF5F4EFD234DA7DCD9590F86DCDAFED50E80F29F82E23FFE4333FBD6AFC91C0E545AA22E9583DCEFDC70726CB89EC7CD01AE9B0EA7EDA433C2E2F785
5E3FD43A6B22FF4304BB9D3570A977CA79336DC2B132475C42EB9DC69CFDCF16EC3AA330CFC26B81B783D0CE834242D22330385276424C4E7467A4AE
ED72C9C7A785B50E9EB7357D715121E8E94AA927F50796AD6ECA32B8C4BBA2E8FB5237E090D9E3C3877C77D6453DE9E20367AD7AF1FE263AC1E741CC
F69EC9EEF2210CB9677E9733535CCB250404BD61924011908BE23FB7006645E60E8D62DCAB1262EB32E17F99FDDBDEDF6CF4AFDBB5ACAF2105BC7290
8B0928FF2928FEEA3F01F7BBF8ADFF40BFB972C3B2C559E483EE2B1BC6F9C3ED53799D9DDB0AC1E858415ADE0B9FA2C8B00A7EDBD45070E54C07919B
86A6EBDC00D21E4CAD16933CE72164FF4E63F3ABC5C2ABC3967AFA4BBA5E1512D2E34058C5BB13D20A2A67EDB5F74BEFBCA5A6B82168BBE6FEC26DEB
9BD2C51D69B70D2DBF7ED35C5911A761FF4C6FCF49DDD8678647925BEEC7CE40EC6BAE4DE0401F0BC799C04C7F1FB7B5BD8F84AF18BD0099492D2886
41ADB563108143188621D89FFD3E78430C12F21113756E2704FFAEED7FC9FE31A866AB7ACABB0740A565224A50FE53FC5F1AF2FFFBED3FA8D2EAEC9A
7393A92B6DFDCCAF919FDF57E81630C18C276359F0FD33D0D4E5EED15BD55E27AAB6DAF6B31B8058B7C9F65EF0533BDB5FF93E67DFD797C273B45E49
093F3FACF47CA7F062AF2D1BAFDA68EF16125A05B1D0C4C7AE0B242CD6C9B8CA1B3FD297B03F6D2D679BEEAB1C3FEAA214FF7191D1C75C43F36F2F94
56DFB6B33F2D77314D33A0AFEC58020883B1BBDBB09EC6418C3B0DF37B86999D25748C7E78F957B233630A01610C1C63902482623886C1300F119063
7D2039B47581B86FDFBF7BFCFF6DEBAF6997524C86733FD779F71041A5FF1414BFF90F14AFBCEEB0875EB56DABAF8903DA7B0FB27CC499693B5D34C4
EE0B61104DDFB0B7B96F2D23AFCBDC04C6B8FD563726DBE1A14A56B5D5CAECAF97FAB62A4A9D392977F8A3A543C022CD2DBA2611A79454E609ED4378
FC37C5DB4564D597EFD8A412B25BD174979F9AD6D3285987E118D9FDD5762AD179AB96E56419ADBEE1B5E9B0CA81788D9D9C22B78F303A70764F3F67
B0858E617C3E9F3E355957C7043E1FB8CA64D7577150040231822067D5277F9EF8072010A595779282387509A5B0E97FABFDEBFE1F81E3AD3EA6C9DF
4CD2E1CC85F719A4809AFF4141F1EB02C04ED03BE36DD3D47FC468859166317A9913787ABA8A5371F3471FFC2E1AE7BC1E6AB9DC7EC9252E68491BD0
033D544EE9A840722B4623155DBA8F5FCBF45496BCB15DEADC79E1B57B25EC564B39075ACE11167A8480C38FEAAD1728883BED13F1B925B1F4A8A3E6
A2AB8FE535DE3C54DF9EB241E25AC93A8D777D2EBA57F69B1D94734B367218CEDC982D603707F8D68095B5C3B3B93D8B3BC266659EFBCE3FB6BE89F7
237508C5F93CF8E7D03F02984DFA277E940CC3085A9F364C72BCC5C417C6F3FF72E8E72F15C0ECA2D0B3C7A8A3CD229839715CB7102228FD29287EAB
FF5977342FDF342F86EE1A2DB35F76878C4C2FBD0A74F5F1936ED5802307BB045FEAC88FE96D61571BF646F287D8B42DDB6B8AE983C535356717DFAF
0DFC7A7CA5D0862011BD2B6B74F7193A9D52510A38242526568EC0231F5F2E9258286FA66DF2C84ACF67BFA172D85BBDF957F3558C138E897A54EC90
7BCCF05B1810E6BB43F948B6ADC5C00F871472E4716536B324A7843D8EE3F0248B315979FD49DBBE752DA3A1D91809226C1E08F1414C8092647F49E3
2882C12D5500D963252129739DF7EB7DBFFFAC6C3831E8A95E3761B5B36F94176CFF0DA5FCA7A0F83DFF9FDAAD72E18D5EF474BA99CA7A1B27A0EAEC
E49D5676EDF0C8CBE821F0ED39ACEA286B26BCB8C9EF52DDA152A817C95E91589B38D354F8A2E7BA71ED879B662E160B1E1D105F7B65878DA1FE51F3
79EB4E480BA9F4439CBACF1EF31545F425C50377893BF82D960A2E583ED72661F3BAE7FB65CDB382A57D182FE4ED825D9CE57C9ADC9412CB367D205B
C29A60A0F0763B7D7A066FA99D610CA69DC9A4EDB4E92D8E1D2351DE0C8707CC704002830962A478901410CCEA1682FCA821A22CB56CF0674F8FBF33
9BFCF77B2C4CC7FD94D3785DC9675C13312AFDA7A0F8DDFF116795B0EF76E13363E74CAC3DAD0B60FF8AC216DE6029BBF5512D30115C4A84C693B971
1577AC4BBE3E04590CFEB30BBD6F1A4647EEDD99F0DE71E68CB35BA8FC9D7C1BF52DDE219E0BB7ED58A8B55F55C8920D81795FBCC4959457285BC658
EAEF7750DADFBA61AE6E6C9CE3B7EBCA4A77E264ECBA0AB58CCEEC58A3E159BF57F4C2D0EE6841D3B93612FD78A20F6903D19E9C0EEEC4B3E86CFE87
632D999524894E32611400B83044407C62A070822471FEC03881078A8AEB6ACABE4305827FD6BF7BAF460E39BDE7584DF7F9A541BE8594FF1414BFA7
FF78ED3A950B559BF70F629F4DB5B66D3D85BD3A5F9DCEE25637B3BEBC6570328BF032F7BE89F0B03863C7F1D71FF1B1A1D198EC81679D8C6AF5A40A
AD15779E2ABBAF347EE1AF66237F3842DEC04F5D6AB79C903D9FCF4A7EBC72AEB0A99AF8A1C362EB3C256C8BB7092BECBDB129FDE3C2C5EEA9C66AC5
35864A8E3E5692D6D587445CBE1C08054BAF8E62D51F0FD6008D0052D6C120FA1F3E1E62DDF5A86B6B2649C61444A010078011F664D700429FE49300
9D0F0978DBE78BAB998AEF6213C43FE94FA02DBE0615C4787D5E5B6BC80AF76B36F9022AFFA7A0F8DDFF948D6B0F771DBF508166EC36B40E751A1EF7
FEF261089C2A19A6254E612DF12876399C9F19FCEED4E2E2B1339D48FFBB967B632F930738FB543B2296BCCD5AB9FCACFAAA0B9B0C97C86EF534DCB5
422D4049E804C41F29CCB59A2BAB276E9760A37278834AC8FE054A1EFB5D634A5C14B63EF4907A37E428BE6EF7DEE586A921226BB3824EB33F3C8227
9FA4277C05DB46E945BD24D4723B89413B7EB483DF8513D38CD9B80D8030024D15D6F601A46088C767F30464FFEA0512521A4AD603FF6C3F01FCF0B3
AEE1339B0BCA9A1B028EC71EF1EDFAD37FEA4BA7A0F81FFF71ECAEF339BBF1BB9FAE7717DE715278E1740F8AB85C93CC84BA7338D5E93010104B9638
34D15EDD4A36DA3871C593899D7DF4F1D5686809BD4539A07F8D4FA487C50D0B31C7BB56AB74C57C365899C8BAA909BD8280B6A2A2A57F4849AB5FBD
A8BC72A3AC83BB98965D80E48A2A6FC9F5877DC4C227F649993BEDD597FC962BA173F95978ED8B38A82AAC77E46EDD6003ABAC13C36AB65DE60DEC0A
9A020627712E8B1070610485E0E1AC0E082389A649C10C9F243F2F5A20AEE6B542B9F49FA23F8113E36FAC3DDA7B8B2612E3AACB5F87647C5CE2334D
457F0A8ADFCA7F1CF4F7C8B1E47DE9708E9FFC7845E6E239676ECB9191D21604F85ECB7DDF4494AECF06036279ADB125574C32B3756E42A5EE3F9EF5
14A70E0361BADDE196B941D2F6AEBA9A8E112BDD94CCB7BB9B8A2F145BD006436D51411AD2F3D646FA291B1D36D10A31D0396C6F67109EBA6CF9F943
C62758B7E6991FDE6C25F594B652764F654D726C2DF8E5058771A51419A74D0DA3CC371181A9797B636690D60102E593DC69808011B03D7F08126093
1DFD240090FCE3E2A2D24BAE5DD28F477EBDE3FFAF2D0D01D4B0536D5563412C3BED547E6FD3C3E2FE54C35B1075F48782E277FFA79C5CC68D2AB3EA
0F444291F71557BE36CAE6FADCE6958CA153EF46C6DF4D93C95B1B13EDDB986931C5B6EBEB2F5B550BEE8527C62455E6F50DEAF8A71AC7EED6B45D6F
A4A7BCCB67BDB9F0AAEBEBE72D1092674240C32D177121F957414A22B6368A7B572FBD6C686B66FECC4A33E0B8DCAAB10C651D073B8BF93BA01BF337
56B63C8C1885DE7F86F9C1C968FFF7F17E363D767F5271F8B9CF2073848DF2F9D0C830FCB3F2AFAAEB8771B8BA731C63C082166BE185520AF7836582
B9BFC4FD7FCF33652519EA86641D7C36C08CB89AD8FEF14529F3E4C24C8C6AFD4741F19BFF4493A527BAEA4E7DCEA73BDCF71F6C940BDC4FF1E37670
69DF58785D1ADE5D2920DF0694BBDDE58CBC2C7EB92AAAF3CA4370FA6E70437AF940297045FEE9F177B6DAA7366BAE3456DCB9D3516AD1F9108DF942
5630CC6D49DA2726647E668198919D8CC32D952027CD88CD07FC55B6DD5E231BDB68BAEEF2E6CDF246B41FD26AE9A30F5EB5B1C25F23BC07B1E058D1
1003EA3ABBAFA2E7C68B5600A6B3000C20C6C76101ABB4B0B59F358AD0731B49FE148EBD5E24A22C277F2464B16BF73F24FF183C744B7B4F2EE3D8EE
A9E6FCDBAF3ACBC33B38A3EB77F650E53F05C5DFFCCFD45E031F0EE8BBD7EF57D6FCE58068E475D3FA4E9B2F6879398EBFAB21EB7A08D6B598B7D6F9
60C9937C87259F5A0F86C283FBF3A75F760E76E62E3A7428E1DE32EFE5864A36EA7A46DE167356C46E1512F24580F1C67A7F21A1C5F2B21ACEC22A4F
8D3CB6C99D59AE765271E5555DD1CBD5D6D6E71DAE588B950C6B493E2FBE15CF609EB8CE65C424229D2FCAC6C6C61FAC4F045E5FA5E1A3C840CE108C
F40F0A04F4E2A23E60AA0560958F10389F18F715115DA820E1F6CCD22003FE7BA75F9851BAD7F0197BEAB2FF6861C697D80AFAB9343A1EA01C364652
F19F82E237FFD1379A86BDF1E75B1E338F27F16B4EC8D947DB4532EE1C86F1FC6E72F00E02E52164F6A5F18BFBFBE9DFBEDFD5585D9062FA1CF871B6
37F95CD7E0989F6FE083C96B47D729AA5859A82CDA7C4866E1959BC242CF31909593E92A24B564D13ACB1552A167D6DC96737B28E77CCBF4E00EA9FD
593E5677979AEC927F306AA99E99722365020AF1E4A2B9F9F0D0BBE231565BF0AE82D1677706710E772ABB044006E8303991D30972279BFB7813A000
16C09F4C44E43434B464D66DD4BF44FFFBCE1FA7E7A1B556146F38E4C2102DB8B8AD9577F51688DE57B68DE7518FFE5050FCA6BF807F5B5DFE4DE5F3
94F3ADD10FB9DDC714153FB8EDF8311C9045F29AC604092FC9C92E41CDAE8FD5BB8F4EB1BE7DB314DFD7EBA094CE0F3AC8DCB67F7C2675F5F90D23D9
17EE4B2B29592B1AC91DB19B6B1C2A2B948FF159F74E180809E9EEDE24AF78E082F2737BF993763BB2DCB7D84805E52F77B8B172D16649CF344D93DE
628F2C04BC6AD50C7C2D44BA9ED4733B1AEFF924425F42BBB1511C48AF46A63B2648B2A7781080672A7286272154400E1E545DA82A6EAA2061ECAE6F
DB4110BFD5FEE8649ABBA6FC3E4647E0335E7B40D2703B9CE83E81D6A818D9C62054FA4F41F1FB02C0B8A02B1D59FBBA3638B7F562EF58BC89C6ED33
5B321A3B4E641182612EFBD614993D0AC65DAB2DB53BC3182C7BA866D9FD4EC6B279C4FBEE13CDD7EC899397B5AFD31EC76D50947171B09676BC2023
E1A926D143A0037E178CE7C87AAC9591B57193F6B928EBB0523B37C366D3E2CB1F2C34DF1E15B5725977426169F797753910FBA6712E519BC36BBF51
421FFA1014DC36FCFE423186B3472A3BB9602F2820ABA37AD1B1FED2B221908708F02F0E0A06D28A0E6B55EC43BD0DB291BFE5FE7D7797A9AB9F67D7
1C4B660D1EBC3A3ECC6E58570902DE0B6D2CF2509CF29F82E237FF478EBA498696874EDEF0ED7FFD9E4EF354DC7A7FC3B787A32D561F488801B74F92
89AFC9F63D713311D2E16073FD6E912B1D0672B970CDD692837A0D60C20D5D95F674D7FB168A8BDD562D59EEBF74BEE6C2459304D6EB1BE9F98759A0
C682F59BFED87A4D62E95903CBF3FB5DE583327464EEDD94D0F5DA7A58CD61E486551AB3E7A079195CF36E90F6208BDB7EDE33F8CB68F28B42EE24D2
5DD7814FF53348B238AE1F446945853C88C327E8B7561A6B08EB9F75115F17ECA9F58CF797F04FFEBCEB0716B89869AE7D4EF484E7C175EEF7069BC7
D87BDF40FC4C052BDDE58D2075FA8F82E2D70580100CEC3BAB73B5F96943CAB6315EFC30F052D3E1CBAE875971E4832D5324413454931D0100167981
5EA5A658000E6559AE1A3D3DF73C0F48391DBFF10958187D40F63ADD6DAFA7AA88ACF932A7FD4E728A92F6B000EBFEF67AB1E4364B29D31D862BE32D
85DD4EAD34503BE8EC9E6826E75BE5A8E6786CA7FAB1DE0BDBDAC74F6DF7A8409B1286F16F1F267A43AD8F5777252433F9645B65079FCB67100276CE
871992D15A3200406C106B3BEEEDA02FA27F67ABF0BAF53AEA37C67FDDFB2305D3A1ABCE059E6A2406A3BBE01EAFE0998901EEE36B209DEDB4F880D2
3E1A9FDAFFA7A0F835FE13829EED47F403EBE3623BFCCBA19A37DCDEE506CF8EEF99892D859CAF9302F2C78E49322A46307022A6D6566C038FD61325
FBA84646A1161CBB57786CF570C3ED5C35E59667AACE4AF384CDE5179DB2975098E38A0B786D45FBE6CF5BA724B34D635EF801515DD7E51AEB83BD14
427CC497573F13591660B370DFF43EC769E48AA47707DC9E328AA53D1E643E73BA3BC22CA9A111504B71371F9B62A0C4D0D34F0484F6548CA0200B86
93765FBFB2649EC6C36D0BB605A82E3A3F8CFD3AE807C87674ED2E383445D6470D1063FB6EB387DB793F4E4DCCF063C4BDED16DC1CA20BA8028082E2
37FFFBDD8F6F3CD313F778E2927B4FFBF356FE590BFFE0439DA36FC70AB6672124C3EB233979B51DCDDF7F31D04FFC29BFAB63BF59E711A1C3CCE6E0
C8049B2430B4718FD8F51E8FF5CBC465CC4D179EF490D15C7092247909F7ECE6496DD058B95371E34D5105AF9D7A664937149D9FAC36DFFFD9523AF0
928C6BC7A99DD3D01DD5930DBCCEEC4EA2FAF1303BE6481A971D19320963FDC53482C360C248FEC37216C9EB6987308807CF0479A6DE3714D6BAE3AD
14745143F160CF2FF61348EB61E59B38EA5F2048BF3B48A08FEFB180AE317E4E0B97D3ABA31FA8A9963E320153050005C52FD53F26E876D875F440FB
A3279DAE9EA5B7BFE60D662EDB7CDBFB36DEFD819DE2554392511B07C8DA970CF6AB6D61B9D69A5D135F52256EE689C97C9A78B321CBD58EF9FE6992
B84651ACF65AAD7932DBB46DFCE525C462098417F5E2E47C610D83F52A1A976CE79A1D31160D4956518CF015DF767D858C47A0EC9AEEA35B38CC8B9A
A1CCAE928F4D48F58D0ED69BF335DCA99B61F5007FA081C70767F848F78BC43A72669A3ECC1300285ABC696F6D98AAB8D9AD9DA6D14FF4A4F735233F
FD27FE75DA9740C65FAA1AA492FCA769CCB8339D2496F00A64F50E4E967F9F66F28FAB9DB21036CC69FD31F573BB80BA034841F117FF7B9D363FF46A
791E3DB87BFF44DA9DA2C9FEBDEEB78F070CE0DFFBD04F0500D4BD2F14C37F5473469F1ECFCFD53D0ABCBFEB6D50BB7BDE03801D19FDD6E2EBE88522
0B99C0861D1BEDE5A56DAC96DF5FB140AB9340D81F1A4ECF5D60E665296EB2537C91E776E975853612ABCEB81D88CF36B3F9A061507FD7AB2DDF7543
EAC84471DE8F9936BF5CFCE3951630C3F73D8F3B54350E932C169FF13EAC04EA9B9902B9A400C4A1DB6621DDA7A4E5ACAF3A98BDBDA424BBB316FE33
FCFF1CED87E328FDBB83B2C7E8F8A7BC8CF23BD75A09A0308B074F0CB33855034CE09DD6F1538B54B737BC8EA5E154024041F197EC1FC5894613C3B7
5BCADF04B5841C7F3771F9E4F850A2D3BE905D211C76CA3802D147810CBB5C8C9737C4A9DC7D81F9606166EDC6A74BFC13E6D8D167BA03D3D66D6287
4606482946DF5631D6992363B968BF9DF00A2E828CA46668CFFD63F976513D5D65754F276585B7174474D6EB3834B4BAEB04AC5549083835126ABE29
1528F8FEEEDB78CDF687A38FCEF530EF6D4D07E1E1FE5101CE65239D0E1B68F844FB184942028C1CD86D95357A4046DAFCACB994C716057987EF2CEC
5FE19FC071845E15B8D8B67E32E77A63CEB3C88CF289A1FC62923731CA636754F2F9093AE687B5F72A1CA8BF153541F94F41F197F08F3171A24849B9
E4F4A7D29B3FE28A121AB33DBE42ED870E9CDDB3BD18ADAF40C8AAD7E04CE8D601B83E6566F4C19EDEA14D8EEDA78F5ED67BB55EBD9C076624C619E7
7CBA90A9AFB637485E557181B2C522752D611B9087F457C7CE9D3377A1EC22072D89ADFBF4857D3E2A496C7655FFC87AA8EDE86B7CD53D603A42CBF7
43674E635EEE30F3D1F5BE948B8D50E88676647A741822493A1F2BF73852CF19EBE493284A906491F5D109DA360929533F3BDDE05BCB64F55E4FA03F
27F9100481639CEAC0E55619BCE1D70D63C5C79346B83365F1BD084EA7B3B1E2B70C285D5BFE80B577B8EACB6FF79E52F19F82E2AFFEC38300FE5D5A
A639E2EA58EC8784E6D444DAE5033CDAF34D07E38E5DEEC0927AC8EF3B47197D0E4758ECD799CC74FBAB8C873A51F966775C373E925C3F43A7DF2FDE
EA56EF54E82FBC2BDB7DAEBCF47C0D4B1125919D1807292D089FD55F4CD4CD50D4C46695885EA4E97C53C785CAA56F642C762BF93AECE1C668BD7C7B
AFA9282F9AC1BE7DB3E5E5B15AFE9B939D00A3B36B9AA0F360A23DF20BC248A9E4C0387FF6733E59F912EBB313175F7ECDDDE6FDE8D3A58A21A3204A
E0027CD67FB0E59CA6E4791690993F09DE77A8E47541D309F524CC9E04909267D35086918ABBDBF264F7BD5FBCB6848D6394FF1414FFF19F6437CF20
3992AA7D9FFC18E5F7877ADB5BDB320EE642253B7C621EC6D48DD6D5110CFF28C6D437EDE7FCA1C43EDA15BBBA260FF7E2539ED146113B54DF31D11B
0F3E1A571EF4ADD7B16FC85BA5A729B6EA90BAA27828C985BF575D1612DDB46CF951FDA507565BA878582F303CBF5E6BEFD9A50B77EF72B9E431DCBE
3C847D2AB7E5694E7FE38DA35F4376C40C56E44FC0F5552D30007310F0CB953184D79947C341982407029C6B040DCB17483BDFDAA61B961EA0A37EB1
0FC6D1D91F1CE70F3E3096347903F5DFAD4587765BBD059853ECDA6E01CA0521012D630AA7392AFA44BBDD4FDA9AFFD6DD2C8A41C57F0A8AFFF84F08
7A4AC6F8B9C292B555E65D43C7BFD1EAB845C92F23215AC896C3516558EDA7965AA2C8B1A36D34C4A60B2A4D9A2E5816C08D597A2173E5C9FD16CF34
8CBB99F9EEE53B0EDE138F0B97B8D17C7F99A5B8ACAFBF8AE84B92C7492A741092B45DBF454BDAF9A8BED86EB7398A2E1774946F6D96DDBAD5D4CFAB
99B36D4FEF9537C5A78A469B771F2E89B1BD177BBF1DEBCDCC1A9B4147A708ECF5964CA87FB09BC98758042BE7D81D3A99A8FA87F8B6D7764601070C
550D43FA5008FE39049CD5F1C441765BC6685F55722FB770B7DD8399AE314173390F67327092F56D18C742358FBED97C2377F3A3ECBC20DB6F3C4A7F
0A8A3F23FFBF1E97292B1A83CB2544B3A7ADB2800767988D8303AF1373EBF0BCED47C263B9F4E886F65EE8D2B9CEEF333B1E4153D98D534126D9B400
C7B4F0154F75EE1D57784F03CFDC4DB5B8AD6B95B1DE3EBBD0D55E5F52FBE062995C129B7C556D272424A6A52565656BAB6F156428B276FB0A7DFF72
A715172CB7C67E673F58555112937E297BEADBE1DBA51FEFBEAABC1D5B937BBF761A9BEA62936377F7A6B7E51570F82003E6B6BE4B6D21C9FB32A292
F6B736583FF55CB4C8367618C6101C06663A1ED84AC958A7B7C53E62D0FA3FED3A515E3B4463D13A791804A224FD4533017FB23A937C2DB8EA827B5D
7DB4997D154E9DFFA3A0F84FB72C027D9B46430AE70B87724F1D063E19D43717C16F4E4D154D4DDFB2BB187A87188AEC4B1D9938559A5BD4E1D18674
C64F96E8DB4F943A9C29B772DDA3F74261230465D8A4F8BA6D143E15F1C781C6038ACB45C4D464A52A48A235E2BBAE90ACADC23C0F7711FD4D57CDE7
ACBCACB550EBED27531B339BA0D7FC62DD87CD918FFCEB86C30E3D6B8EB8DEDB78A998F525B11DE58D7470597DF73C53F905BD3C680C065B338B5B60
828C961797D7D8B829F08EFD2203F7E451048451FE54FB8335F216D6CB37BD6A189DCE7FF0E942545543460B6BA48DC3E68FE3F874CC779C2C5E7AA0
607B445DCC9A475579BA8B5C0704D4F3FF1494FDFF03782B751A2B1015F365167B0CF4AD8E1FAC418BB68CD11BD9D95E5E09115DC4D72F759F81D667
DCE7835FEFB0B939698CD33297A6123745C5685D368B3BA336C86BB63FFFF9BCE52AFD576B9655BC365C2125AB2CA3374292358F1A8D84967A499B84
698AAF74F696527C7142C26889D75E7B4B8798A436BED7F6AED4ABE732DB03D605543CBAD2CAB91A8BF36B2661707864044AF33BDADC56C122003A7B
AA30A382C1C305714A128B1CF6B83DCC7752B5F0CF9D62820834DDF6748DF8964C7EDFDD220884BFD9A455C77697148D32C7C6E8E434172299A955A4
6068AB73C19D83534FAE1DA82C77565812C4A1FCA7A0F89F23B302EEC9F753F077D92566DD5337E24037176EEBC0F48637702799F1DAE368752D3C12
3CF8F53BF4FCDBD433CEEB0FE0D4EBE6521DF914768069E349039783CF847DB099B3AB1FA4C9ADB7D81B2B7BA778F5626D5105B10DA080CC7AD0A43B
574E58E78CF61C13D7556A12779F0AEB6D9037B6D1F12E4CA8808357657C397AA5A1D5D1F241A8FFDE0E4E68186FB4B40B417B9A90A9077B62EADBD2
C7082E03EE4FC89D1A66E140849498C1F18BD617D39D54CDFDF37B7900C8E87EB246C42A091C7975B10441E0DCBD798DF143A98F26D0B1762E639C85
F3C75ABB20927BDCA122DEBEE3D9CD1CD7C6D86DEADA2F104A7F0A8ABFF8DFE7164D47EB8DCFE895D04B9F325E9974D1878893EE643D8B1FFBF0C8FD
473C382E9313353E16D05EF97DF4490EB3FD43DF5531DBB13A3D9F8E3506078F2BC8D6009F3CB69718AA1E737EEFB1BD2AD078B9B4B2D859122653DE
3768492DD4F45E33477E9593B1A84791D53C3B5BCDAD6A0E55677CB9595A4FFB2F5D6E18F7D5BF9A7FDAF1E3C4F3F019B06386807ED4B2061FEE4CE0
7D2B40490693539B5C4BE7D2A0113F1935ABAB37D69C2B7251B50AC99DE082CCCEFBD69226CFD99C54FF581A0A4F6546140DC60E947D1C8278333C2E
93033213AA7924CCBDB5BEAE6BE3BB111766C8C5CE1B3B552C7FFC2CFFA9AF9F82F2FFFF4CC8C01296054F20ED16B74C1E0E0F2595726DFC8951E2B5
6117FDFF61EF2DBCA2DCDEFFFDB5B0490B415004C542C5C02EEC6E8F85DDADD8DD8A280A088808A80848772A484B770DCC5003CC305D4FE73CBF9913
EFCFF77FF83D97B8F4CC72E9B0E65CFB7EDD7B3F7BEFDF4CE679CFCB4EC578F96DB0260C4D3D29CF2AED7D9188FD496A5C35D81BBF6B1A9464EBB477
8EE91128D56FA6D7458339A757BE99F72D6ADA34CB11835E3238E3111E676479FBEDCC7E96ABC68F1FB7D3C3D96EE2A4C5330DA7E77A0D0F2D9C76AF
F4FE9962E553E7E080FB9BBDDADFB8AB257F5A4159CA871EC5DD7DF59CA44A8AEA40D1F8D06E44D04D36AE1F6E35C7F5DD82592FAE3BEE799FD3A151
70BCE70D5FEE2F037B7ED7CA08BC33BF2C53D01AD25A5880935D5142AC574C023E914A9516725B5FD6BAF39E2CB3B2C3991BB16EA6C59A46D67F1696
FFE7965CD877D1A55E5CB0F9D499B3A2EAE612EAFB3C3182743AFB30593DA2A7071FACBE87689E35A1A93588EB435989A0E87C0396559C3361265C33
7A76F5D19147AF58D87564E59C5FFD6DA4B1D7F483076FFE5E3D71ACF9E01006A71F65FA0CB179BD7588C5E235A346ED485A3DC07989FD9A09563F0A
469FA8DBB8BFD9637FA424F4409AF4FBC5905CDF1F12A4A1498600455CDE933BD5351920A38118E5F74429209111BCD526B66B6F853A4F7CFF75D705
EFA44E3EC77F85D9AA1F6AF1EF3F7294C2852979BD88E07796B457D73EA862B361B51A457E7A833400F86ECE929E78DEC13B25E90E56EC769C35F544
B396F59F85D5FFFFE2BFF4D9E2FD3DB860A3B3D726594211B705DF7C879232EEBB988A4075E1F9B38F1C2A90CF01B03C18A8760A53C78A7D0E88D4DF
5B9387FB9387877EF01B3EC7739CA187A75F83D33D9701570E9B5E5F137CC3D4D264D02F0601CEFFB93DD87A8689F99A71B6D666EE4F8C562EB15934
6EC8FBDA29B31BAECDAAF55B1BDFFAF64C59D78F6D6732DFBE13F5FC690030E25772C5BE8DA9990D5A840B32C0D748392ED1504DCE8346EC7E7EE7DE
9187412EBB1FE734D57C3D60E110A0C1B2DF16E35A4A931DD98042F9C10502511B85F5140B71AC578AA4BCE855A829FF83A5CA0BEE3DC8D78580E457
F8C4B9F6532F573214AB3F0BEBFF7FA7E5D0CD578F1CEE224427E7441E6CAC495737D3B113B231827F5A28FAD8DE99FBEACDEC3B68D92B04AF8EC3DC
5696738AF8BB9EA0B53EB0DBD8A65CABDD05DBA7BDD93AF4C6C7D9EDDEFB52A6EC0919FBD7910B5F2DC68C18DDA1C5357778D747CE1865B66AA3C91C
C7759F6DED1F8CB79F6776A9758B7160B3E3B3840DC1CD2F76C6A802E65F88F38A6817D5348A302CD8ABF9FCE24FD5226D5F11C234DD8EC061194CA5
4F1C38FAE29D4349857177972DFB9C5713B8C57CE2652E0986064974A5BC2DB618C0D4E12F2A30058EE33C9E98C021158AC495C06A24C3B51C7B7516
6E92BD4FC3DEA6EE9FB361E2EC0D51ACFF2CACFEFFCF1579196B1EEFAA41A44FACBD1E7D830B1B7E4B31D7836A82490DA7E2E2DA1BBE795F9EC9015C
9A0830B8BAC9F54D7B1A2F6B410AF12A556DBD02B83822FA99E5A21B6397E65BDFECDBE4B3D8316891C583055F361819AF8611B2FB5EC152F309A6AB
F7994D5E33E2B99BD172A7B1530D36F51C35989BBB7DF3C7496E6DCF9625167B1CF6C80CB95FC1CB2ED32092276E3FD76FF82901F0BE925E0AAB2C84
3015053F331F3CE9F587136E8D3EEBC66F7B911BB8DBC2F1218F20B9B1B908A9054AB2840421FE1ECC1376B5239D6D4A0243400D86377EE58244F2CD
6EE8F5E15669746905F3E364B1D3D2B9CB16CEF422D8FCCFC2FAFF3FA888834927731130D6F672840F939A915B457685FF166991240D3FBBBD322ABE
7AE30534EC3A46F233EAAA6B8A62FE68DCD774F21F283C4CDC6B675D4EDCE6F47A9D4DE6D3798DBB365FB0FFBC6FF82DE70BEE86834F124ABCE3B5DF
D83176ABEE2D1DBA7AF582EB73179C1F3F7D9843E3F38166E73F4D3EB3EBDEAF8B6703726FDC0E17BDBD5CC3CF6B95A0DC0327F2761E6F8485584789
4AABC6692DA8A2443B0D8DCDDD7A7E44BDDB3B6DE67EDF57BB6D367D13318C323B4D4850202FB39ED092CAD2720DD1549299D00E9114A4ABFF38525E
8A68429EB5E25EEB8AA8025E03293A5314B4E2D2E5349F694118EB3F0BEBFFFFF4C73D4F669E8D81E0B2096BD3DCA157516A4FB1124DF7D5321DD968
46B630FF8DF4FBB078E468128D37248681D25FE1458213EFA96FA1E022F3FCD74E8197669EDC3DF844EB8AAB8F16079A9DBF3874C7B9E9EE6606BB11
1952F1DC7FBCA9C3A325C61336DB382DEB77306849FF51BF43071BCE383B77F1A5D0A88D173975079E74A47BDF882B4FF903813EF35F043A3FEF162B
D0D6020188C0B40602B58D0B06191A5DEFE5255C1FEBB8E5D2A935165B1370DDBBE7FFA88408849F55D2ABFB0F594D3B0E77D7178414633488E2B842
09D766C39C8F7E4AF4C5E63FF28A021528BC1E163BEFE6A51FBCDF53BFA214BBFCCFC2FAFF5FFE479E6D49791D08439D6BE7FEFC1115F49D082D9542
0DBB1A18BA885B15D7A34ECCE99835B98F73AC8FA98B3C938E082AF38182F772D985DE2CE375393BDE5F379EFE60E6AC863B569FCFF92D5B1B307DC9
CB398B871B1CC40134FD73C9BCE16BB71B0F5BBE65DCA895D6F38F8E1C161D6762347FC3EEA52FABBE2F0D902BCE1FECE53C7B94212868C654A74D4F
7E5CED2E57F091724F3E51DF48525284FE336D88D1D08325351F66DBAF38B6618199530AC930B43A3BA98B4495D9BF3A489A21FAAA1518D1919E5258
87D2248643DD2D6A718FB43D9D8BE00967CA716EA59496DE7EF567CDA3C0D3F9CAC71322F4E7FFB3FF03B0B0FEFFB3FCDFB66BA2DF27579552E23226
8A7BB6EA778624465E08DC3D8D31B551505414CABBCAF7EAF7924CBBCF48E3C2F777C9538B78B2FB514C882B72D9C8E38AC3CB5163AFEF1E7C2564E4
F1D00BAED65E2B87B9DE3036337884AB90D0B0F431E6F3865BCDDF646AB967ECA8B926FD5C73ACFB39DD9E38DCB52B614118AEB9B4A2B6E0DA474E4D
8600C53F3AB8BCBF1CD9D92AC4736E640B6539D920405029B6830D873F8C38B5C169FE9ED5338D1D0254144D12F1A1950849762796410443007210A5
A08E9AF63600D3A21845CA3820D6A3E1460B312CEDA918AAFF2DC3E5DE6E3587AE94EF2D9078DA587C67F33F0BABFFFFFC2763B6AD7E9970480828EE
DA86004F7E7102D4059C6849DDD628AD265EC10DE886DEB935CF5A24A69E643025BCBBC16095C71349F85A3EBEEF937CEEBCD00907F78FDEE53E7E74
E48589DE2F3F8FBE7977D48AC8D9FD0DF308004ECEBA66623B67E48C95D3A62F786A316E53FFD1F76DFAAD3EB7C57C7343E5C2D790E0864D02FFC2FD
3AA147909A6838F0A22237BEB25DA32E7D988CF19076AE86C0832D061B5ADFF258337EE636E71946633DC438896991E07039C920FCF23688423005A0
65A8E652A1AE5B00308AD2E23804A26AA0F1DB1F58E976AA07E555F64148F85DD1EB83AD4FE388381BAB31BF2892CDFF2CACFEFFEA0FBD727D7135FF
AF0A8DD4DFF27E7B8AFB9FAEF0F6EED204F1EFEB2D744480E6C73798B737FFABD50F9A7B4BD3F7BB7075099C713E09BB7A8DF935A73879D8D3E3C36F
580FB975CEE040D8B0657E016BB6BE996AF2E87DFF69625D610E8A5830D0DC64DEB4B136CF1FCF1A3B6EA4E59559FDCE5E993D715C8078DB4EA07ED3
887B75276E55E5DEF2EFEEFDB4C4B5BAC037A70504939E94F02BA4C242250E3D341A30C4F6D41907EB050BE68F9CFFBA07C3554D65BDF9050843A23D
1C310E02380EEB54EE0CEBA6C47D5C09CDD0148AA3B844DE171BAA109C9CFC136AC9E1C3AA029FA2F0DDD91EB1E8AF8796A3C6E6B2F3FF2CACFFFF5B
FD53DCBFFFE54EDDCD58A02F6AE6BAD4AC583F22A692D49414F3CA73C8C2BD1CFEB70E34F251E7B933001D9982C52BEE9F42D15C1FAAFB543EB16F87
FAE2CCD0A5BB369ACDFBBA7CE8C7FD631E265E3D17B8DA64C62B53170CC779B763979B59CF383963ACC38D4BD6BBA70FFCCBB3DFBEA8C9E60E77B897
D6D437AD1DBCB2E8D291C218BF2702E6C9988369C921C52A4DA3A7571E511D27148831DEAEFEFD874C5EBBCCD671FA5C9BC55F54040188CAEBCAD212
7B35B49AD3AE0251882429DD37D39924C7843D2A00552B309A22011885EAFFB44BA37725C0A2562E80B6F8D5349CFB58E94E16353EB7B4BEC4C519D6
7F16D6FF7FF4A7FACE5D0C3F55F1F97A5F6FE6AA499E5FFA0ED7415C6D5B1637106F2CEDD8EA09737E69BA8E67966EF981CA3F75F37A2B6CDC493028
904E3DAD281813285ABEF7F692BDE6834EC58F5999683125F1D4C6900BA6C31D061DC701A2EC7CD6D481236E1E349CBDDC7ADAC6A366F697C76C2E59
65386F5E7AC09A32FEAAA1CBFDAF6C6C885AFFA903FFBEF0607E5940BE1C2B7DF6AAB5AFF96B612F44FCB2EDD7DFC46EB9BD9D9DD3B8894FA4342549
88EE904B22B3A534AD68D4DF042C453194D232BD79126D9B403F2D289591348AA2385217CF6B0CBDFF1B96D6F42980A6ECCC9CCF57FFDC6FE88E29B5
1B6E79BC93F59F85F5FF3FFF69D10597ECA375F91BAA45E5C7A77FFA823EF84470F0AEA7D2A432EAA7D0FFB6026C9121F11E9DD1FBD388A274A518B9
62D54477BD2CD4BC0984CF4F69CE5B76E4C8B369C3C686BBDB7FBB6CE9E76DFF3666AAE9D0FE67B43855EA55603E7893BBBDC9CECD43A6EF9E38F3CC
44FBA26726D30E9EDC73BB4670C0F478CEDD79F96D1BFF4A41B2A75CAEFE12584F00C9CFF21AE066DF3288A0324CFB0D315FEE6069397984DDA97286
817E5D0C1481C5E95D24CD902A158128241081E35A6D75234DD488685AA734456B498CC0C0DF49F5FCF0D834A9462C4151B05ED299F92CA1D01DACAF
386D6935614D0F3BFDCFC2FAFFCFDC3F4D525D8727FE3816D7362F54D2E46AF53429AB715B5F4B2613F2ADCF1F6989E8B81A86C03CA4E3786AF7E143
62594561035165FC174C655D68EA3E9CD8627518B867BDF8C50143ABA92F8EAD0D1ABFE9E78C35694EA6A6FD2E3028DD98E1DDCFF8E9CA7E93663B5A
2C32377FB2AD5F50E89019EF4ED94DE788B68D712B74B38D111F38959796B4EE4645BC3707C6E2DC8BD4AA8A4F5510833E311FD8CF6896C308AB1123
F6E6918CB0252AAC1115E7556B708A201886C0D472058A123456D0A0E90A4B00B514A9A5684C435274CD9B2465875FAA94FB9BA38041794703C8755B
597A2306A87CB76194CB7D1780A2D9F93F16D6FF7FFCC71B4F8CF2F77F2C5E7743D0F971C2FEBE24E5A158348A563C2BADAE4123CBB3AF348275F59A
A0A71D112BC3DBBBFB0A21FCF0C8261CF6F4A2A2F7CADCCCD3F21C67DF4A186BB1D871D5A5E70BE7255E9D1CBADB70FC407786605A8A2EF7B73A32CA
66B9D9FC0D0BADFE721DF7306AD6C8CB112B6D623AF6193E153C1BFF46ECBAAD21EFCEFDCF9D5FDDB88CDAEF95B00FAC7BDE4CD0D0BE410347EEBFB3
C46A94E1D42F28C34F28E44B28A0FD4FB72EDBE3BA524FC32042A2A80255569696FBE694F6E1FA357D9A81512D9C7CB71205E3A3E55863911845643C
19D21674D2A7F4408F366BE7EE25D16F5CF53D03EB3F0BEBBF567F6D0612B8CCD2AB657BCBD3F9953DD18E933B934A0BE2E9C23C26FCA998FBFBA7BB
48FF0C7D624BD566DFB6FD67BB53DA643554F1B06B84AACCB9123B7AB177C942FEA541730B9EF73BEB6EB6C86D8CE98DF831073D4D8799FEA021E667
C2DE0136E38CA79AAC3F646FB9CF65F89B1A07D38377A6D94572570FBBF6FB9AD5A3DEF3DB7811BB02B952CFA7AD281172AD01420AEFD662B474E320
93C1BB23F75A0D1FB8A79DA173637B1906E49737887014C374A6E3084AA03808A284A42E24BC93C6289CD08D640CDA07A1C59FBA6820F03B28EBE0CB
51A8432E02355F7CCE8A02AEA0F4A3150B8F459EF0C2D8FCCFC2FAAFFDF7DA1CF8C2CCABB9C4ADECD6393FC4B95BC726FC0C41EB0138992EBF1D2AFE
1378A51C090A45055590970B276E7F775D06C827E823D31448DFF1EB78C1A2D404FBF8FCC956A11C2BFB927B2607E60E5D95B2D3FCC50E43DB0A06C5
3E676D1D3A75B8B5DD96A025A3E6DDDD7CBDDBC57089EFC9E91F65A78D4FD69E1875B9C86571D9EF2D1F21E8CD955612F4BBC593A37F2E54E344B3F3
A0913367AF19633E70943B4E3586D7E8DE7047EA9F0E5DDEC7084297FF5184A0301425290A2C13502489E1144DEACA3F0DA34A114589DE1D2942C47C
040355324055F5F45B72AA36A54CDBEA747A7BE0EFC351ECF21F0BEBFFFFEABF62CF2CBF602AD61D3D79BCBDE6F8A813B94FE5FC248693C30BF7FE44
A4DDBDDC0D3CCC209B9B389B1FAA5CDFE3E585CA5EB271F2336143FEF458FCDA06C1E1792567FAAF29DED2EFD41F876507071BDE783664ADABD18C0E
92EE7D5F3863C46873B3912FF70C1AB776D1B6967726B6FB0E6C8EEBBB6CB23BE788D5C592B5737E7F99EF83299FEE2F8491B0574D08927EBE0A2732
ED075A4FDD34C37AC0988D15445B680EC030C8EF1F8D508FB40F26499C402104C7718C409B6333C5F50D000A533445E9AF0723E46A82A45BEEECCC10
0AC500000AF962ACEDD895CA977C2CAA0F3BEC707B4ECC53FB28F6F46F16D6FFFFFCA73B57AFBABA44D2FE04495DFE87FB74CAB2FC6F1C24A1154A6D
2CA97CD05A9EFDEA95A4FE5119D9228B3FD0D2724F8494720A4B88AF1BEFBBCB3C7676D63885144F3F173469C4E58386B6C1C71CDDAD8CD6072D18BB
66F0660CA33B62AA6C478F1D6D3D77EDB0D14E4E4E31859623B679AE3F551A31C22EE8FED8DD9F2FEDF9F9CDE1A11ABCB5AA1893BEBDD3A9C2E3EFF0
08D4C77CD0A81953C70CB1FEEB37A00AFB2160189C535EA106C51D22254EE00802E92A3F0A13D0D3830922582CD2C7005DF5D73234C897C18CF0A7C7
AD7C04870112D54850B0C53D24C1B79468CC2592CC6F3CBBE635CD3E8DBDFC8385F5FF7FCBFF4D738FA54F0BD55CAE122CFFC20F9E6CF7253D9E68FE
0C4B8AC3447FEA62B8EA087FB4E969135E213C751C8AFC46C82A6A53AA25492B66A708971E073D9C4B5D4C1E3C1DB1F0E594819BFD2D76AE30B47E73
68A0B5C121544DB527C40F1F696334E2EA2CD305738C46079F369CE1B1C3DC25D879F4C92BB687DE8F9BDB983CF58E5AF568759EA6FBD6430123F2F3
1290F8F38103CDC78C1F357C66A8064972AF2028A0B898CFD09A0E25491124ACEB01743F10940483BD255A4A4132B47E4580246906EEEC64B40551BF
DDAA194C83A388024548C43D80772F036C4D0405331DE3B6BD5C32627C1EC5EACFC2FAFFEFE65FAA76DD5DF8C609F87D0C78E1666FEED24917AA8AC4
9AE472429958892A3392082CBB585BEBC90338B98E5F355F6544B34A96D087784F5D531F61172EBE1896BC7855F291299F6E0C9DFA7DDBC27BD3ACDF
3E301DD9CF1D5391BDBF3E0C1D6D3D62EFFD31539CE78C5AF7C66AE4EB20EB25B5C196CB3F1D5EF0FC6CFFE5FECB4EF7014FFF6AC01BCE78A909F0E6
0D9886AF0C186866616E32E35D355A7D3912C0C986C04684617AB83843D3384652FAE80FE30456984733926E94D6EAABBFEE2745AA412D121DCCFF5D
8EE9BA048A2421906C0E8A89F1AAD3D497F710F746B945EF0D7B61E7C465FD6761FDFFCF7FA260F9D19EB05582D487D298850535BB6D16E63456805D
71AD5A5E9414E59F2B23898C1626EF21565CFEC0AEA0228454D411BD7F1A148F2CEF761F9AD45A7922E2B1D9F9B809CE7E3306EFF9B4F0CBCA91671F
0E373588C0E4445DE43543CB210B3C6DCD66CF34717AED607CDA7FD2A2CC68BBA9AE27A6589F709CE63A7131B7F3F2F11679EB691F98446E5F5033ED
1B0D0C8C4C870CBF5587CA7EDFCAC670283B4A4033587B1BAA6574559E26709CC440485E1C524D03AD25629A26755F0449692914617A7C12E99E2A8A
46518A02210AFB7E2E3FCABD16E264F7A1D156EB3D4F5C6EE44D5AD44769D9E7FF58D801E06F482C69CECA86C6FDD5DD1EAD5D9B3C5A9ECD9FE556DD
DAD3DE9528A08A73102C746F0D2DFF4E301E3EF2B0F21D9B8A3FB531B53D8C20A287FBD75A7EFDA433CA67CBBF398FF6BD3DFEE56DE331C1DB4FAE1D
76C86BD250D3024C4324849C301A6DEDB5D7CC76A193CDE34BA6AB221D477F4F9D6575F0F5EEF167AB9CEFBD5F97DD757041BDBAE2CAC33E94F47B80
301FADFA0D1C3460C09A6C4693F427A90747395115A81604383C8AD14F55D0344110148E285EDFAA87D5DD9D72F29F5729DDA0D44E4199EFABB5848A
D1C50106C7682D9EBEAF8AEFDB8B35B7F542B5E3563E3F7CAE42FD79F261959666FD676107807FFC07BD1C6DA3FACEC7687EC4493F9FA9FEEAB87045
726377A5A8260C57FFEEE46B82CFE753DC12867B82CBFD1437CF35FB850A8D42C8E8206186ED3DD06D22AF6EDC6977B385B10B66F83A0F5C7BC6E154
7FFB4FCB0D46D6A22AFCB5EFEE41C6E7EF1A5A2F98647CF89DC9944733ED9FDE729CB271978BF9B8AED82977E705AA2E2CFF5395E879B716865E3EC1
996F4306180E3618721B64F80149953081C57C1191A0462494E98CD5477C14A7188CCBD5881BD5A802D4BDA8335F97FE49126DA853B67F0915331A05
A62575A3024E13DAD62315FC6005CEAB94006DCEB33FBAAEBE23291F37FA9C9A61FD6761FDFF37FF830F9C67DC937B7F55FF79A2EA3E9D50B471E5EA
E7690A598AE66B2E2EED48E7C0696772984A3113EB43E5043CD9911C134E977AF70ABC7E291FD826498E14130FA6466FB77CE831EEEC5DE35177169E
1F6B1178D0C0BE1383809B113B074C7832C174F2CA89937DD70C3BBC75DA3B4FDBD19BAE9C9C65EA2999B6FAE59B9EC8ED19D2CCBC2725A0E6FC5E84
C91A6E6834B0DF849F14DDF3DE0FD090686C2C8263BD228864189DE65A82C019821FF4A55D21C7945254F7CADFB55FDFFD833055175A4D615205A67F
B6972610AA9D5390DDD088D09DBD5289FAAC63F2B5E7F7E2A0B346361F6186CDFF2CACFFFF8C00A4F2CCE14D8B38B957B89D7BB334A7F656DCB7BDFA
F40587FC9CA77CDC86290ABD1550FAF562AAB9BBFE621D9EF476DFCE2ADF6CE66D645FA9776999955D6BF26DA273EAC68FA3ED3F1E9DFA665ABF333B
67AD1F726AAF81930845F8A79266F5779A613CDA69E2E873878D962E1B1FF5D576C9D1D9871618FF257D31E2E9E7EECA6DEF4B12632E4501C09DA98D
4CC1E821A683066EE1D278FDA75C991C44CA53615C285011244390A4CE7E4457FCA35D53001C922985528CD2AFE3535A5C57FE953D40494A2FACE4F4
E23445915A949096A77701AD792423ED1474416193C23EECF1DBD9D637D9714E3ECEAEFFB1B0FC9B0168F9859B77EC12052722D4DE3B5471B3427F2E
71CD8F4A53D77EC47EFBC288F8690EA28ABB55DB5728CE78CBC363F76C0EEEFC82F7FA4557D67F6D7A367043C9B64AADE798F7BB46B9F82F38F1C87C
A9C7CCFDD6931DFB1D85619CFFFA97DD107B0BDB49EB27393E9A6CBB606E60D402A7C3C7CF1CB65CD129B259D99E52B8EF72634BC2B93018739F5B47
7D3537321D6CFC02A7251969A594A202EAAB83A5BC6EFD7E1F92D0B9AECBFE34E7D5BD2AADCE6EB95805E95FA3B5BA2F82C21BF30A4B6578475527A2
CB043AB71151CD9F0E4614DCCED4A7937265DBCA3B859B8A5C43E16B33225C5AB514C3FACFC2A287D64A8EEE7E3FE6BEE4FE23056F7916B47EFEAFF3
334A940165784C26E45180A145DE7D22717A6073516157D8371CF0DAB42CA9381FFB11EB2DFF13F467D3C0AFAF16CB258E73CF4F73BA7B73AAFB0AC3
D73B17DB0E33EFF70857E16D9FDC860C727018B5CAC662DB227B87610FB3E74F3A647FFDAEDD8878FCC49C129F772F36E6B7079C0D01342F8E72115F
634393C136690C93F5395E23A94CE402504B499900C0618C20B5348A62B4B6F989BB88A131442D532224A60FFF34A5A5411CCC8E6B45A5422E1FA1F4
A18024D49C7A54AB09CD235BDFF14031F4E458D5FE98CFC755B9D609DF8ECAB46CFD6761F92F01288F4C0D5A7B8C1FFB86AFF67F48858F791B36F52C
58968049E239D28876142ACAD248B88D6D35FE25AA6FF14CEF259727A541B2EFA55FA3E0B2E874AB1982BD8F88776687F6AF5CE4FED7E1C72397FFE5
BCD86AB4412009E06D8F57F51F766CEECC0DA3E6BB2E9F3461EDA545D3B7ED7DF260DAE4996F4F3AB5F0AFB43DFF248F3DF849AC3AEE222C3A3862B8
E5E0A5F50C99ECDD2C107DFFDC4B83ED31891214C6218CD2523846926849A59466284A0B49D4847E1DE01FFF09581718849A3E4EB712C32952170828
B04708D3647116CA34D43280ACE4724EDAC98A5D952AE77D8AED5FFEBEFC9BD59F8545EF3F293B3FE69BEBECCAAAD7855D1DFB39AA654BA2765824C0
BEF9B02850A04912628A84369C534508D35FB789AE3731495BCF847BE473DE886E94A39FFCBCFBF935CD4AD32C9C79609BED6ADFB977E71A6E5E7665
7AFFFEB13880FF3C3FDFC0FEA2E5CC111B1EAFB0365B77D672B4CD8A90D849F657B7DA8C8A57BB75051CE1C7AEF396F46DD80A850D1C6269627A4CC5
6812E2BAC9DA0FE90A926E0DCF1322981AC2489CD0D94F757ADDD5E8DF308E016A8840097DFED70D0114AFFA674A97BAAE518220BA57290A27416137
4A4295D9200396E24CD79F1799A25D7ED7BF43C1568DE90B2A289A61F33F0BCB3F50024FA76F6FE664B405FD6A87FDBF6BBFCDF47A34E6085AF84A8D
467B60AD0528262804F1A6121910E8AA697AD0ABDAE390197711FC9410E383356CF9B5DF51E9EDDC9D3676E5E9ADC3CF1D773E3F7ED64AAFF9038635
931011FDC0D9C0D9C97CF29A1717AC47399D98693BDDF688FF923187DD9C1726F4AE7C58B2A2B16DF32BB56085B38A37CAC0A8BF439296A98E2EA654
B9DE7984565D1856A52235004C1024A6CBFAEA3CAF3B397FAB8BC3304960FA677EF5CFF15182B2ACCCBE8EBA0E3546E8F706A304A6118A009A682E21
183AA484E184C5170271DBE33C4502871BD0CE937DA4F6EFB94FF6936761D13264F5E991EF2216FBD5642557F6763DE0A9D76FBD3B675CB1EA6E1AD8
B9B29AA8ACC3885F39145E548D0A0EBD63B21E95C7EC5E5776EA031C59F1AA04BF7AACD4ECA568F931E0E2B0536F468FF798767CD7C8F9E7670CB457
130A22EFE9C8A10EFD86ACF1596BE1306DBEBDF574CB23A7ADC65E3F37654E96E4C4DA56B738CEEAC700779E935A39DE6080C1FA36068BFAD446E009
CF4A281CE2FEA9865490467F961F8651B838F2D1772EA9A5113549E827FC08FD310034CD907249F5EF9AE6560544EA5A0402474048C1E723B4EA7781
54CB64E730EDAD41B164F3FC97618DD84B9BEE3C9B374AE2DF5B4F595858FF193C6FD990B74D07FDF34A2BF278E4CF102662D6DEDDE31E11C9F7B8E8
AB9338D2ACA2BA02D594A6B981CA76AED3065ECF4C59F1A5787776BB28F915C45D537E694E6FA24D40B3E3EC8FEB47BCF5597F72EEBC89B6430ED028
88165CEC3FC27E9083EBAA5163C7DA6E5C6CBDF0E4E1E923B7BE996F1FA3B8B689D37C27CFC55551397BAD843FD5C0C0E0AA86ECF60B97331D9F9E56
905A5CD73E6098528D13108AD338DC12F9B9A8974210B54ABF784FE22485639456A36488FAD4F04A8E0020315D4EC0091CD348DADB95A4D4FB839A61
84F90CD3D2FABE133AB0BDA74C2D767C8EBA4E7CD747FFB3F0C97EF62C2C5A062B983924A8ED695CF40F41461C22BFD6801D98B06FBC431BF749A410
5C1541839D085E9481928A4C1172E0AC56EDFFD8F3FA2DD5833D791AA97B1CF9E242B9F5E38E29E68DEF071CBA6A3EBD68D7FC8D9B868F1C7C835610
AA97BB070D31187AFB2F43CBB1864B1F4D197EE4D618AB53872759860ACFAEFAC3DF7DEBD30375EA8C17789AA94E7F37142F70CFA4B1FCFB1E0D8CAE
C0631A1C1229F45B7DF5E71317A7963541302EE4EB4FEED1D22449102889023D324EC2B32FF5204CA03841EAFA040241A59C4639D2FD2954CC30DD59
88B6A0F1579BF4C2AA8AAA66ECD9CC9696450E4FFA68B6FCB3B0FCE7BF22D5C6AA92FF2036E79DA2398A437E7DC3FC9AB27AC998D71D9F73B2F09FCE
DDB4AC95C6E31B48B4B11C8D59CA2139FE7ECF36A57BDFBCF98DAAD8D9CE5DDFEA3123F9F8C82BD58EB3DC571ABABE9EE57C64AEB9913FADC63BEF1F
363232D8E539D162BCF9D2474BACF79CB21B79F0C514F3D7A22B8E590AB7D10FBE0BD316253279463AFD3F6150E8D31A4C1E742B11601894AA4E1511
801CA55082D0CA0B4B5BF9A8440DF4346BB40CC5576B494A7F00505B4379D46DCF120D4511248251B47E0F300E7437B5A97BAB39FA4783EBDA196131
08714F1D2FE6C643FC050F95A1F3ED0E0B2936FFB3B0FCE7BF20CEC628AEE7E83B89571692F103E51DE821AF8C5E39C629CBBB36B38E7A7A99248A9B
E8465F3E0E85B5421BCF13E27791DF475F8DBC1BB9B992D87D937A75B173EC8650BB615977066EDA6364173467F4DAA7B38704532A3CC3E7A08181C5
2B079389C6764F5C4C57DE9866EA72DA665C60F341E79896A86D6E01B545F6494C9E697F8381E158DF9D871D38EFF57B0EC1D002489A5205499BA404
06D364B56F2C9F908B74D55FAAA1702D50A9A6090C256035E7E7B3C04A80A2758D3FA26BFD49DDA84060B2E6C6D2CA1E84D612780707E35762DA3ED7
F50DC4E762F0E1B69F75576C2CB70A683600B0B0FCE77F6F98FDA0D579BED1BD75D912696227F33282699F65ED6873F555BE2240833E6CA7F19C4EA4
E45B2D961D4DA538B52371BB33766F6DDF5BF0D355F27369ADF470E563E77A1743DFA811561B1CADDEEDB39AF2F6B451320D923FE30E1818AC3B6460
BBC87CC21CF319D7D69A6EF59B6415DFBE6FF26F79E8B25782E4FCF52FE978B301FD8C73F19CBFDE48C0D2D7D11A06078A9B750DBFB2ABAD4FD7E253
6855469108D488D5B046431014490372FD29001859E5FF25B240A1FBBD7EC50F2128FD6A2041E2B2FAA45C314C61A456DBAB225B6568CEC6BB6DCA8A
0245C3D9C0DF7517165B5D556AFFF59F1D0158FE7F0FCDF09E8D1A60919C9E5E582FAC81EAE3A8A66B201D6B3C71DC5F179F20D9BF988A1486119562
DCBCC01AE44D0DB2F5362638EA116BEF1B7448F3EE11E6BA074974CF99FC2176F894E06D83761C31DAF674A4CDAA5D66F9B412AFFEB6D4A0DF0AB371
F36CCC461A0DDFBBD168A9D712CBF4A66553D25B1FBCB8FAADE8C1B2CFF4F32183075837137EBBE255D2389F7C0CD724FDE8D2A9C9D135F530859370
4F258A6A502524AE68C3298CC02035A12511022979F8B9488C6A715C57F66914C4F48700501423A9ACAD5390188E5394A24BD2A44282A6B9F4A29519
00F5EEA67FD5AFD56F677C84E8FF763EB39F3E0BEB7FD5CD0946A32B6A5EB6E527CBDA80E42E26F017855D33B6D813B2B71A0A1732CF6B184603D5A4
7063A58D1F889F2B3BE1B4EB7D9B57F7EC7ADC74F157EB7C5F8D67C2868D1D5BFABFF4315D7366CC1CDF7113C7D94E14306AAA366C96C1E859E3FF9A
DECF74F28819772659BFBB6CE2CF593ECC5F71F2ECEF8FE549BBBDF034A34183E6CA882B67EA34F5DFBF374178C6136F9EEE5D55858A180AA3317EB3
1C520A418AC638BF0404296C5111BAA44FE10437BC0645690CD35B4FA0108C53B4FE04C0F417F53049EA1A024C565022270548D769976E4DA250A215
9C090B515C71AD5F944A50ACFF2C2CFFE67F6DD9A37D66433EB69D69EEFB582B90D7A760BDCFA458DF2A43ABC40BF7F186448677BA8741D48A5B859A
5C282C073FE50AC92F45466C11BE1EF5E3D7B6A6A8799535F7569BFE491DECF475C5D80D2B469F5B3062A6D10254ABA6937DEC0C0C972F9F6D613671
62FF45AB06AD3D3FE868D5DCA17E797BAF763C7E1EBFFE0BEA6B3C68C839AD78F76B8DA0E2E3972E4C1CEC5ED8C9E08C348DA3A5719AEEC890916D79
521A2094657C5C509CDFABEBF0294497F285B0AE3540F47BFDF4E780033041EBB426E25E7668295D77801082209F4AB08BDF9118256352CB309879F6
24B8B8C431E5DDE86AF2DFFACFE67F1616864C7AF669BEC955F87D0C5CE32143FAE2BAE9B0200CFD6969783B6B4D8D36AD99C97B28A60034ED8EA88F
DFF519AD9B9F84A65C2A3F199634ED08F7CE33CD95AB9284EDA697C1E34353DCC6ACBA60B3E4FCF83923F7E8820599EC6535D0CC65C5B0610E4E5663
F79BD85D9D3F357B57FF9B3DAB9CDACBD61485EE0C7F3DA2FFA0C74CD7CEF79022DB27BC07EEFBE6D74E23786F77959C21A5B0166CE243BC9F8DBACE
1EE27171C5AF6A88D2CFF2E138ACD265039CC0FF7EFA1F25485417FF498AACAED67D3FFA910053E67E2991403C594FAA8C6908D3424CEDFE8C20C5C5
C3E90E4BA5ACFF2C2CFF2BFFD4F733D1C7072F12673CED033E4613DD1589B0E27C0D069FED37AEC9F510DE9688303F9E092165DFA150A41B898CC46E
4FAA01AF05041FCA7F35E749A14B44FDCA00D93DA3618D2DC3AF47CE9C73CE62C48325238D2FEBFA8A9ED03743FA4F586731D2DCD9CC60CBBA518726
1845F80F3CCC39601C223FED9BB5CD7DE7E001FDCF31355BBEE3ADC95E3FC470DBC750B54622C2A18E3686E117AB48580616E671480222E4CD6A495E
2F45513AD131FDC2204A1004A49FFA234098A47194C029B295A36F67308CA20589DFF930D626EE7E5FC468EE36E390FAE4B777F5A56B534E9B5D92FF
7BF7176B3F0B8BCEFF1097E4832673057DAF6BC09A375CBE38BA012F7D0D600DE3067B364FFD41563611EA98300045235D8A2065EFCB82AE392E40C9
BDFA276F6B6FCE2DCDDD9F1BB4AD34D5D1F01EBA7771D9493B97B9C3AE5F3432F26308E6D7B7D743064CB0ED3F73B6D3A0092776ACDA387867AD9D73
CB6BD3BB22CF5B1D3B77C46EE8DF7F93BA727BB85A54F93D4BA64E732F52A2D26E9022411C570AC5BA3ADE9650A88110A150A151913D5DBACAAF2BF2
04064100A9AFFEFAD50052ADA4745183C4491AECC4F5DF11A5257B52AB4518C86D45F2A331E6E31706A53EBAFE0C07EEBCAB5B3AE78982DDFBC7C2F2
9FFE0CFC6EE58BC343C794E2B15F7B5501A97DCD75C17D58441E8ABD1A30A9FBE552A9B854A896E676F6C0FCDDD75422B8F48230779C37E91FFE7B75
44ED8610C4636BABDFADBA8B46A3F23F0CCD7333DDEE3260C9A3418312180AF78DBC3C60A0A1D1B8A9B65613ED9C0E590D7449D8E6501060B44BECBD
9A933E78FCCA9156D7C194D539744396CF77B132E0462EA20075D51CEFEC2424BDBAB70673D22AD56467AEEE1F47501CA249FD413F2482831A042428
447FC61F8A03807E473085E21880313481110CD3F4BA88405A130B35955F14DA440F580DB7DECAFE262D74697D1C7AF4988AADFC2C2CFFF35FE37BE0
83E78A49DF89566F2E52FC51592D4E2F4455994A4C3867F0FD9EA59F541C51AB58C11548549F4FFA026A22EA2BE23EA7A2DBB3F5D082C6804D5CF90E
3FF05164A6B5F1DADC69815993379E193AE3B6E1A86686813FA7EF1C3064F0EC1D632C96CCB3BDBA7BC8AEAC2B53BFE78D9D5F13B334B57B83B1C900
EB9F60F2C64CB033F465B442F0FE7D8B062D8E139078595C3B0E930CD2915EDC4B8085DF0B5B94322957A2535F57EFF5CFF8811215A1A5682D4DA118
A18BFEFA3B4C7042836B69404330B4C4F383A03725AE0395E70AA95A4F1985B526B564D42157439B3CF19DF761D67F16967FA1B592E7670AE3E64CD8
2EC14293E48A97B99D5C61601DD69729A77E9A8CADFFB0A4A5958B35F21045767B83FFEA1CA24F76BB44B3611F9A9E946C78EED730173075BEA0F654
EE6E23CBC2436B2A77CDDE35DD68ABA19D0A47E1E89869FDFBDB3F73183273AED1B92F862BCA979B7C2E9C32E64F9EC307FCB2A1B9B1C1253C697A12
56FAC1AB008093BDFB7024DB3347A3EE4B2D86608656A6A6D4E0982C3AAC9120284D5DA99CA6308224F45FAA3E8D56BF874F23C508844075CEEB6F04
05000253E20423FD1E24E21634C19832B503AF75EFD2606555DD2DB954D429C1EB06E9947092D59F85E5DFFAAFF3FFF69EA8E8C576333298FA600D19
EE27AB90D58670F40FFB53AF4C2E7216BDEF288295B56AB0B452D9E7B9852BC3538FAB7F4DCB5787F61EB4FFB565588472D10BF28D7BF68491A9A9B6
B96EF607D7F79F63B89C40D4EAB86FE3FB0D7DE13868DCDCC9B33F2CB6F4FACB6049F316637FFEF64BA5415696660397F6362E7153F4FA7FE4928AE7
C79AD49D4DE10DA40CEFEA044898288FC80430BCEA478A4A3FCFA7E06BB47FAFF4EB320082A835FA137EC96E2EACEB06681CC0F52700C31801A02403
D46716C951B106C7A57F04589D4B15A9EECA924A3868F75F29B5BE54E4942AB6FB6761F95FFED7769ED97C30EEB5C5E2AB20F88187355EE334D6209C
6C359E5B45A3272657BBAE2C892FC1851C445C02C31D4E271049EFBED7D0A95DE8CFAC4A539700D339C29756CD6DFBEBAE0E39DB39F3D98FF1AB370C
B2ECB70A56134D7E1F2CFB9FDB32C07CDBF8C1EBF6192CBA63601AFE61D0B9D61517F86BCD2D8619AE52B52FF156663F0E1112F20B97796455C92FA9
BA4F242C92A2544774588352DB9E93DB0BE318A212F581FAF3BF75C91FD6DFFC8D3334C320BD7C44DFEE231801118CAE0B00D4108A6B4BC3AAB5F246
19812A2A7970DFFE7018EA4A12CA2A61ED9BB322DF6278DD02011BFF5958FE57FF9996A3177744974EDBBCBB589B544C82CFB9A272106DE113B2A436
A67BEDD3B0EDA105EF7AA16E14FBF513877CECD25450CA86969A69BF80F7E243269F768F8DE33B3EA29EFA578C5ADC757C69F2B629FB9D6D8C2E620A
BCD2E3A3A5F5E6C186F3B6DA4D5A35C1EEEB61831D95F367B73C9AC74B18606434CA5BA9387545C1BBE9A9A62ACFF8E3686AAA4C4972BAF0562949D4
B9C77730CA3F691D38599FD107AB452A8AA6182D4E202886E004A56B5B547D6A5A970674EDBF4AA5C4080A2515BD280D1684F3704D9D1027A45D7D20
FAE915A8E6A594533D2AA6745B6AAFA7863BEC26C8EACFC2F2BFFE9F2ED91170E941F9BE49279F409C7014BF150E36D72150AD8ACA8D556B1F8E7FEC
76A9C13F1091B711A57B1AD4D21D3BE442F1D123921BCEE0F7F86C93E5AFAD960BAFD8B5E4AFE31C3009F71D137B63F092DBD30CEE630A20FFF1ABB1
438D06CF5D387AC631BB015BBD0D4DAFAC9E9D7ECBB4A8D772C0E081C18C64FBC9F6D46BA10A49D8F6541A7FF1484D890532B4BC16A378377FE8F2FB
CFCF422D597A39450D6A2008D35FF28351984CD809E84FFCD5C84092D66FF94711618F148568A4BB410468123CBA19692E1747C4223588559D9620D2
A4544828633A4FDF6E8BFC08948F0CC369F6436761F92FFF33D97FA5FA5F6AF01B77DA251F8D10E19FEE228A8A5E54DC4BB767F0887A07A7F7EE71ED
7E5D2857A1F8F8A9589DBFA94809FC5CF03ADFC9A7F1A5E0BA53A9A7736981AD9B667B7484C9A9DC792F03C6D8B92C32788F2148C21B2FCB41832DCF
CC1A7569FF30BBED36838F9E1EEE5538F6357CCD60F0C073247FCFE12E4E600AD0BD6F422A53BEEE448F96A815C94ABA71B2F35D0A4601315F3AB454
F6B56C0886215DDDA7B524AD2BEF255DCD4A12C36118A57192A2089C42FAF8008D43FCEA76024A0AEC26847F5A084AA346318C7F201FD5B4FC966272
52FEEC6C6EA36703193EB19062EB3F0BCBFFF5FFD1EB9242D615A7DA38DF7A02E5E7629C53EDCA8E5211D50E365734F5927E66DBA39E56E56703C20A
54D51C1BAFF4BCD0276C39B72128D035FD4251A9CD8DCEB33BAAB7D857FAFE55643BB368F3CA8899C6CB8D0D4251B532F8E58341FD074C5EE7B07197
85F5D1F9060BDEDA9D8DB0BF24BB6330D0603DD6346D5319EF6709CCD9691B4D962FBFA166DAFEF476FA171344AF5F1641887343450C5C18D980A21A
1006D5284A223852142DC4351A3582A35AFD01A0E4DFF95FA9413154DAA942A1FA4C1125286DC16099920011D4F7032215B5080125D1EEF321B9203E
02554C76E4D3EC73BF2C2CFF3700C46F4A2B58F3B366DFEC7B570AC5716ACC330B816B1A2800A9AC934A54F21BB60F12D39B782A314FA1A19A3FB776
5D0E5148BF3FFCF43EA3FDFB85E2794669F9D363A2AC6334AB42AE4D88BC3DC7FFD8D8D52386F3310D1AFBFE747FC3514ED3F71D19E3B87C8FCDE2D7
574F645FDE5F7AC0608081AD983B7D5952D69544B076AD7D369175FB378374E40894C9DF5448D3C74A8CAEA8EA2188969206800234208A4008AE2149
94DF826811508312144911FA0701299C2250FD7E4000452955663A88C18086A05018233025AF58A91689E598126EF30EEACACF2DE862DE0F3C84D1AC
FD2C2CFF3701F863AB7FF5CA074D118EEBDEBF1214556345EF0054CA1528A9E25F3298A714EE9EFD332EB71D559689BBDBA9462F71E9BE6CB46E7FE1
FB0B609BD3898B837636CD71AE98791C39B23DC6EEEAD7E90783A75819CE0251084B78B1B95FFF19C3C75FB61B7E62F690459F364C0C5E7495B7C4C0
C86028B7C96A6D43C4CE13D9C93397578B133E0A08EDAF4CA4E76B0BD8151C514852BF7E600CFFB9573B09A10004422081A238D8D92162085CA6C6F4
AB8004A19F0BD47502180291984AAE1B0EC455620C80491AC3640258C32BE7616A8904540068E7C3B8A6CAC8D27C40E56816A75FFD6357005958FE2D
FFD4B763576BAF1DAE6FF96BC673CF14798612F0E752A842D8A614140861295F9DB3F46E60522248B637E15CADE6463CE875A489BA76B165A51F7CC8
A264D3C8EAA0F1E9771DCBBF38A6ECDC18BB7AE5F7C316838E90308EC77AAD1E683ACF7AEFA9518BB7988D5EB271C2A5BBFB725E1B0D36185AD06E3D
B331DFF5625BC4E825A5C293AF95746D5325A66AAE277B3F7F6B9253A5E7354CAF77921C2375CD3F0E23A87EC77F5260551FA95200B8FEDC7F9AA611
00A3485437CEA024AC0048554B1F8681088EA104DADA20ADFBA3526A644A185543AA676F7A732373549D54B9C92699F69FB3FFD911808545EF3FE97D
765756F4DCE8C6A0854E91D125D585688D3F4A0955080F6F2855693A39528FB9C7928A3308A21EAA17C8FD5ED6F6ED3B0654DB6686D89515995FFC3C
FE569EED228F690793EC5D5F397C3D347BA7C7847EE7299CE4BB05CC1E30DC64EE05738B0396E6BBCCFB9F4DB8D576CA785AFF098D22ABB9F581674E
17F49C39D49579C94BDDCDAB5140489D80AAF02C29ABC5A1142E58C7D5503408A000ACF31C43082AF79384262100C5291CD7EFF5A55100D56AF5E780
A1CA5E84688EA880715D56C02198D4748B3A7A1498482445F57305DFAF37E6849789E43CED9D1161F8DFFB7ED9A33F5858FEF19F703B7A35B867C3D3
CED6CBB32F17A47725F4E05F73B59A2E0A568AAB440A4D734DE3CED90F35116D5A41FE9F0CBCD4E3AD3CCBE619EE7112FC6BA7F8B643C27DDB920726
0FDD66949C9913B3DCF3D42ABB9BF307BDD2E274DA8398394616764F0E18CDDD307AC12123FBA2F0EB27CD870DB76EC0368C2F4ED9BF3758E671A343
F022550DC696027479702355F4224F538068B23A7ABE3413140842288C9228AC5FE69335110C030A71AD7E0790CE7F2DA3EBFD758317891272990AE3
C4D6600481E10882122A1C168909B853A082490024ABCFE4377C2C956321BDAA49EB05FFECFA67070016967FFC476E2D7DB85FB9EE7877F9F7E7333E
714A4AD340D1730E236D674871774D2FC44F69CE5E38A950F34DC6549604760BD32203C4372DE3C883312DB63ED5C357273B7ECA1AB82A7176C15BA3
E00367D76C1BBD60EA402F8626E29F7CB13734BAE86E3466CDB059CEA3C6BFB9F47C8BC114DB413FA82B138A42D65C7A0EC71E2CE00714E2255F0A50
38CFAD589EFD3CB5A31203DC4344059D24A454AA410424081C9754F6C100CDD02D95284D12B4DE7F5DDBA2A54882A018A45D882B72BFF2682D8AE318
A8845004ED6E5513C26229828200D9733EA4CF375984551F549DB4F043B4F4BFF99FF59F85452782EAC1ACAB7BE2EE1FE7734B4A0EEDE1D770531A99
AC8F0C23E8D5129AC65691A6365110627E846888A549694D32274B111C5AE5BCA2B7F92A1CB0A0E6D6D0F8EF8E1F1DC77BECCCFF69EE7EF7F4E2D373
6C268D4AD6FDCDB57E3EB6266617560E59BF6AD881D5C60FDF1CAC5C62BCC76474579A7D78F2A6E3D1ADF15B3308490C2CF62D5669429F74CAC32F87
C8C554EFA3345DCF8E4B3BBB2528A02BF730014706762024A12EAEC7F51B7D294AD7071014C5FC3D0C00256DBDF12F225A743D817E514021D60586AE
16845416B5A0380AC1F0B3FB9AB440A91ABB12D3326C7113FDF7EA3F7BF60F0BCB3FFE335D3726ADF539133A23561D1ED5B02E525AD95C20C35ED631
24AF97C6143D956A756AA6FC8A591A9398CE68D0E29C0031F1F99BDFC88B4CAC3B74F4B9D87665EF9495772D976F381AB961CBB52B63F6CD19356C44
39033325EF9E0D35369F3872F60A73FB6D23D6040C2FCCB1583162906FEDF49018C77D7E959F5D7228F5FD9FED6F8BF0B677A1BC1ABF93010A4C5B7F
BF90A255822E719F442C43684CE7784796085521FC1621A3BF168822485AD30151384E503002A41621318FB31404A96BEC290C918A615CC35792640F
17266139067FB920AC7EC757A34597701FB3EB4A92F8CF7F1616168666B8B7562D6A7812BEC01779BB05F0BAA469AEE86D616282182DD0DE8E6110B7
512D88AB6F5EE50051817C959216876433E01BF72323FCE8A7E9A51BF37FDB783CDAFBC27AC53DFBBFBECCB9E26371F9808DC5AC2E2DCC647C786239
728EBDCDDE7926BB5DEC032EADFA633ECE65D04E3CC9BF69DFAA973947DE344384F7732A3D1DAD3DFF4692F1E7AE571FC524DC69D3A2AACA46050E80
6A98A27186692D9148C5AD9D725A7FB28FBEC8EBC60739A5F39FD2CA783DDD1D59F510AD1B27480225715D62406199FE3E003541631A00CD3FCD813E
D443A0E20D175B639B8AE1E43F277FB1FEB3B0FCED7FFDB9630B4BDF7ADDFAABAB6A666DB3E36F7561B398A95E26D722544D1354AFAEAFD170CAAAF3
D6BF609465F2749AA9BB2567BACE7D9A67D354BDA3E8856DCC87C5675F7F9D36CC6DB1D56F87F91F8D4F86380E5CD78752E05BCF0B26E633C74F5834
66E2F149FBBE0CF1D8DD6F91F1B452FF7792134B7D5E3D3CDE8162897715116FA08613CF4A22EBF71D913194D703094DF4A41742BACE1D512A609A11
86A789115EAD00A3B424CDE8F4C775291FD03FFAAFC5202500F57EC82628942029542444701243847248F787FA00AD5629C65B1FD4E30105881CC928
C563CC9C1B3110A7FF59FF6707001616FDF11FF5FBEFED4C78EB5D7EE41BECF65ABDF8A8AAB3424CD2873E33088556D715268A793DADEAA6C48E331E
6D281CFC8BC19E65334CC1A5A7E617A8A063790E16FEEF176D8A7E69F1CE7F56DCFE091F26AD0F9EDEDF158368C5934F67468E1C317FD314B3F5271D
3EDC1EF761F414AB1111B9E611BF57BECABDF1A8532DF878B9B7FE5C5BCBD5FB3FEB72379F936A813B6F31BAA734A15E43406A39A8C24828F359B616
EF6BE823285DE5A669ADAECC93140EE857003158A321B3FC6B89BFEFFAA051B5CE7B4ACCD1457C9250F6E9D2020AAA130AC8642F4C85359660D2D9C3
6FF4697A702D3BFDCFC2F23FFFE9C6EDAECBD27F1CEEBBBFAA5B7156FE7D6CB6BA960730B98721404381A5A2B0623550D3ADCDF729B55A2763DA0F2B
98BCEB284C7FFA70C1B61CD87A2B7EF8BCCA13E6C7022D56A52DDDFB74DC15179BFDE6FDBD083526FEE0B3DCD8D46AFBBAB1E66BE69CFD3A7CFCD629
238625754DDF59BAEF54F0219F4648B9C34128DD96C93D7AB9ACC17FFB171525381AA095CA9A2A9A484A23051144066AF2BEE99A7E49B5824070FDED
9E14A5F71FC1088A20758303D11AF1B99DA6715C7FEFAF7E3210EF6C05F4AB037D3C9021181A2EAEC30B2F883458573940BA19DB7D92F27E43FFD47F
567F1696BFFDE7ED3DBB24ADEEB1B474E9033431045B764A2CED4069E27339C84348B950535A09AA790479D6FFF984470CE3F195219EE66AB55CEFD6
BD27D519CE658F86BE8976B23EBE6376C8FB192FCF6D7C643F6AF0805842832105DF67185B5D3A6E35CD71D7DA47176D8F2E1A621FA07C322EA37895
E7F787754AF59D315FEA6FC5B5B9F95628DFEDAD85D0DE0D5E04D0DEA7446846A5047190271057B4428C3C3E132475FEEBE280AEFB27747A63B8AE39
D0D57C69EDE704A5FE79000A275082D6529446846B5112157728F4019FE6B461D24735841CA853501CDB319B8B44BF13202DEB3F0BCB7F68E98E9DA7
96DC90FC286EDE3D288EF481E2A6E5E35D3D3853F703CA4D45A9760D5920A2BB3AE8FA8BD541B35A19FE3E39537B834F91A9BF2287BEC62E6E4C739C
947ADBD8F1AAC9829F8E7B3F589E3F6134A07F3CAE263AE343C61A1DF09E3462D53487CDA7672EDA3C74584AE39A0557520FAD7C77E96B1B123DD4AF
7C451AFF4A2A26713DD10A82C081DD80AA438176EBA2830823B0B01715BC265C55FADC5D452304DA51C2D335FE5A0A81704CDAD80EE0DD0901695C4A
4BEA0F0187215C3F2F48125A1A2791B666DDA88031DA9A325553751BA654717B70C2C568FAE5465E600AF077F3CFF6FF2C2C7FFBAF156CBB716E9330
ED67C55DEB190A65A17C87A71C680369B0A899E3FD87823B28059F54B589995CB7BA37AF49E6DD0F86F8F6AE8B941D895C6D9DDAE874FFF90CAFA885
8ED71D867DBEE2FC65EEBA0F56FD8D727190087C1C6F6571DAD978E6B1E9632F646E1C6F6EE62BDC617ABC20664364E6FE0C38D3664FFB9998DA8BB9
247C734E312A057D1E57A7D7406A592F048B704613EA562D56CA1A3E7FCF03B434A9A9FAC351D25A9A861500A16A6E83A08E978763217D1ED0AF0860
88FE89601D5A5D8F0008540C43D1704D92FC7B491F012BC06E800C369EBCE963D3AFA751F03F0FFFB0018085E5EF0640BCF6CA9959850181EA886B4B
EF923CD0EB5ABB4ADDA922B981FC966F6D142C823AD45A6E8D0CF77A25B8F491915D93D18D89C14998D7A508B345AD0FECDC762CFD7A70D8922316EB
3CCDEE1EB27A32DAC0B41403704F8FAF66C36DEC26CDDB683E2BF98DE588412E8A9323977824ACB8A17969112B5EB62EDFE5B3F47000A9BE30D2478C
494F7CC332638568732F4D35778915E559721204F23CB351990C2188BE2E5D5F8FA030A4C06851B31A261B3CDF55E05AFD3300A47EF61FC7F497FFE9
7FC109FDE94038C1C8CB84E9A12A14EE90B78989529B118EBB52AB9EADBF8A6859FF5958FECF7FC546B70793EA2B9FAB9B1F7CB9FE554254BB71D4EA
2E9044C313604E24C868A4848C015B5B20BE7B6CF1DA2F4CEE2FAABB3377552170ABF88C996FE38AC3EEBB7C3F2EB5BEB7D336CCFAA6D7A81516FD26
B7611A2C23FBBDD1B099D31D5CB65A6FD9376C94D9B8C652EBFF8FBDB7FEAA6A5DFFBFC7C16D20A288288AA0A22826D88589DDDD851DD8DD622B1222
82A24888282821084A487777E7622D56F7EC7EE6649FEF397BEF33C6F30F7CE64B07E062FDB01C73BCEF2BEE2B16F926DEDFD7986F3D2DEAEC8CE26F
17F561DEFA76A7A177A540FB9A4B605BBC0215B641CAB2E28C64A19A840824FA7DA5162C17D234C80A1DC150504F51A8A8434D223F9EA4000C85B311
0145B0BF8360942451899EC05156FD9C1F40328C263D0BA33050A16A8384F38CA62E7B9AF6F0D1CA8B18FD9FDDBFFC11C0C3436BD7BD491EF345FCE4
1771FF4CC3A30052F02B06C45B55089EB6AF11C9CAA06825AEEDA4C4097570EE6991FBA82AF2631B92A97DE684145EF1E933A9FAB6C59597371F1E37
DC70C5F2BACD6ED7DEC34C0CE62A103D9C107FAC9BF93CA3C977EC6C66F4303537F055CE355BFBFDF1D42CE9E2A1AF1F0C7893BDA5F3CD27B464E6C4
6349847CF10EADE64389B25A0B34C6247ECBC1B8485FF426B41A13550A489A75EF59FD430041D2884A4191AD8F1EB7728500DC4D1F8E22A84200E130
AEABEB80706E1F30B724A4AD5D2B2E841912C0F4428D766B8FF10E733EBA1EFE3ADEEBFFEDFEE4F5CFC3C3E91FDBEAD9BEE0B5C627812C59FCB3F166
6A45757C27030A204473DD0BD6C58819AA1D6C00C8E26C25F1D1B36EFE76ACE12756F8557B321A7975CFB17F60C68485CF0396EE98D4CF69EDB9958E
4E4653871AECD52330F8DA7B87B1A5C5DC8DCBFB4DDB3E7240B7C3198B2D66FA86CD7BA6DC6574B568E7EADA779515219ADAA953BC7C6089EBDE6655
4C21266957E486C717E5EBB8F1BE11CFC3E52D9DC24E9C8DE54902C7008CB5F14A11C3E0A107826192F50708CE01605F45DA245C1B30DCCAEA9FA2B8
8B4034CB3F54C27AFA044DE0881EF5E96EB6C87E7BDC26B735E63114C5F0FAE7E1F98FFEC9BD4E0D2BF63547786988FB8BC5612FF33E97A769288102
4232D755A0B51934A3946A6B09754289AEF971C175F31438B943E3D25E765A26BDB3DD685CEE0AD3A569A76C4F19CFDCB961FD98ADE683FA1BDC86F5
88E6A6FB5A83C1ABD64F9BD8DDF1F04083F54DF6DDC71EF170D8A5B8DD677BFD36ABAC8AC381A150FD4CDB9F1EE53AEFC7605544A1A654A14D0C6DCE
CD8358F9C73C8ED34B8A5AD43A1243488684A54A36CA275A15BAAAB707CE9611148A70B7FEDCBD3F8AE0001B0160280600288EE3248590E2979E7534
8333384E215AACC5BAC7849133620EEE081A322895FC8FFDE773003C3CAC35BC31FDE7AB8D09EDB74458DD4ACFA60F053F6AB3AB10508240E033770C
CBAE074156514A445C56AEC82EF839E7B8429DA54F790B7C8C01532ED89BBDBD3E7CECFDE4798736D82EB87060E2D58D66FDFA06637A44722D689581
95D390A9A3475F1C6B30AB71474F2B47EF1DCB2BFC4CA6D50759FB4377679F831AB64D7C5FFEABEC799A18AEABEC68286BFA155525AC826926FFD577
56F825953005E12483537AA912C514D9154A227CC98E50198DB2B267B58EB1413F8DA1AC738063188A1318C6D501D12825F8D4CA7A0DACFE090C41D0
BDDDAD768C0BCA9D125F38717C03CDEB9F87E72FF2279F9B3D6FDCEDACBAED4F12218BD2A3DD0B635B625A71B904404BEF161596D55456EB51BC45AA
92B4D6E7A467AC991C58D79CADF3F8ADF5CC127C9CD873C9ED31634707BF1E73DEBEEFF5EB46A79F5819F5FD890370E3397F7B83C1A6E6330D56EC30
E897FFD8C07CF1E6F5E352532D47A7449A5F4002A67BD568F7587B16A5FDBE1A05C873C05CD71241DC97F6C22C846A7EE2232568595D33C83AFD3481
0AC402804423DD626AAA7C1E6568190246510CE3D680B0D13F4462208EE0DC16508CCB05E030008B2AE534011134ACC33400E3DBA3F7DA2597DBD72D
A97FD37F93EAAFF69F8F0078FECFEB9F4103065E51BA6FAECCDC8F609AB34F8B3D425FB71797EA11855EA3880E7AE5A3D097172809B9BEA2105015F8
377AEE3873529ED2DEE8A2CF0D5237EEFF63E31D8715E397C79F39B86CE8BEF889EBCE991A8D2A2100B4EEEEBB11067F74B3B5B23D6F6670376260DF
1937CE4DF06E5D3038386DE4A6DAD0B5C7540DDB66BC682CBDFAF48306FE56AE7FE6595BF52BBFBD09C673CF7E411999A44581102081C18D95756A29
A249FF056291EB7DB40C43A26CB40F62181700B0A69F243012E70E04F63CC050841475AAAAC5144190142601281DD96269B8C47D67798C558870978D
2B4ED20C6FFE7978FEA37F2CD86679EEEFF9E1D8BA2086CE3C94FFE5C1C538B4AA15D16BDB9571B1B2E46A0628A95089841D4F3EE17054E2CFE37767
A588DE82DE6140544962C1987E57ADE72DB3BEEE75EDB0E9902F27AD17F4E931BC8ED0A2B9379E7433E839B487E1F965068B63AC7BCFDE7174C1C596
ED7D9C93A66E54040F3D5E53B663415273E9AEA3DF5BD152BFB6F898D60CB704B1528F3604FFA2F58D2D522D2E0108BC3E264FA2A65A33F23A30E4B3
53084012388523000CA308EBF1E3048369489AC46008FD3326009ADA1A322434CD8D0807141840503B0CC7F8CF7DD0BE685351EC04CBEB08C1F0E69F
87E7BFFE3F15EEB0C4ADFD6000ED375B4F41DED9950F5C7D406D9B0AD60B20C10B11239030707915D409645CF455A97FA7FAF99C3E2CBB59A57AA869
4ADB1C9F3ACEC5B1CF4EBBF97BD67A8D33DBEAEF603EB0F7783109606937AF76EFD66FBCC99A133DAD22F6FD31E196F392B93F5EF45C94BB7085443C
775256DEBCDDE558D4812731FE5A22ABB8214795F5245EA98705A93932541A938192A00EC0742591D91A0A6FFA5105C8FD8FB9D42204C2F50121308A
633049901048A13A0CE5827CAEF4074589964A49B996262992442904465026C464A89BF7AA02FF2929CD67460C7FA527195EFE3C3CFF8162E2A69D7C
52FAF0A25036E13EC3143CC63FDF738D459442046AEB207FDF20A9462DA9AD126002A0F1C15B6147EE27B72CA78CFC70C4FB93AEE0E1CC6A370FCFDE
4B66DBDE1B75E24ABF7E572EF41BD0738A0AD7A291678F18F4EC6B72F6A261B71D9E832CE6ED9C627829C078F8FBDD939A905353D2A327EE9248FD76
DEA9760B50547F4D6A51155EFA8860AA8043E15ABC33305E8D206A146BCD28AE4229A828A193495F641F0A9208CCDDF8E3A81E611D7E9C22111DFB03
0C200881A2088290586389344F4B9104C35A7F02D5C174AAF5B07769B641D58E2FE2120EF7B18D1513BCFE7978FE0B49973ADCAC4F0C395D4D7859B6
30B84F69CBD30F2ED5B0060000118A5E0863904E9456E58BB432BD36E3B9405B713C363003BEDF9A73449A56B579F7B31309A32DD7D9DE3933F7ED04
3387A039C686EB1014465E5FBF6A6038F1D6CB117F585EB6B558E8BEBFEFBEAA3DC6776F0FCFC7DE8F8ECAB33D0081176D9E943F749309EBAB1A5BE2
3E46EB15B9170E8683AAB26FF920A50251A830AD0A8615B2B24C98F019B9B19CC2B9989FE4E67FB0DE3E4272B57F38B70B80CBFFE33837FD53D0589B
0B91304D6320EBFF8338113FC1C6A7E0B447D5C6D31D5773360FDEE5DDF497FB7F1E1E1E8214ADD8228B8ABD965BAC5D7400A52A6E35790707782860
9D0481E44CF146052DE8A0D0FAEC3C5CD64447F8C1BA5B2BDD6F41EFDDA087A915D129EFE60EBABEA1E7ACE9933F8FDFEDDCDFF4E341C35E7B711D04
3D76396D60B0FDF758835176F307F53F16346761B5ABD12E9F119FB0B071C1C9A3B680C40DABABF9B10F1A4A6395AAFA2F4F6B751577FCBFAA41596A
4C11AED6E3A83839560CE2AD85CD2DBAA6C32BBE100CC15DF9C138C96DFD46119CEB04C629BC6B160086215C1D9048A84C951010C6CD07A460444F26
4D1A132E5F78ACC36341D895C8C05E63B61C2AA1283EFBC7C3F35FFD13D0F531DE5F9B7F5627E8B2E69693F8AB672941152E096A1894B0B13773C39F
211B6564476B41222E6B25C2E3C19C8D6767FF905C03CA43F409EECA8001232F8C1FB968E4CD23D35F4CB1B9E5D8D7F02AA98301EFC83D06DD9CD719
98D8CFDB3469F2BED94B7E8758CEF59FF7188E18F5257B84A39872353BD112F0B85490DA0296BF0C2961E2B6A54130961F9AAB80302D0A557CCB93A2
AAD444ADD87DCFAA0F08EBA7742DFB64AD3FC5FED1AAB95120DCC05FBCEB751C253044A8D36488710027740219CD203059BA68DA076DECEED8908319
3EE7C50ED623CF7B95FEC7FEF34F9E878761708CF035591AF95B57569E003FBAC386FB67927C7E24BBE4827067094133A5DB558CB408C66A80A828AA
45887CFE2EBBE77F7431E41D0E78D58A2E7F826EF6DDB4A0D7B4A913EF4C74729CBC7DB0712F17428974DEFDB1C0C064AB61EF713DD71DB13A6A67F8
246498D593C52BE5BE2382E2872D6AC39F585C8E7BE45314F65D8C36BBBF90C0E10B134802FD7D27578E6B75B83429BE08D6E53FF9A8171D36DED5C6
D0DCF42F9C800194EBF32711851AEDDA01DA75094812ACFF8FA81B4BDA5B0428AC23081864180C276B374CF5E84C3FD3FAE57CE0FDF09C4F4636C681
B7ABBBEEFFF9F55F3C3C5D701DF33F464CFD7CB941A68911B75C90D2C4AFF7BFBDAB83BE34CAE0F236F61D0F8269AAB183565680A1BFB16A54E39B97
F3A2654E52FB6DB020BD356F799972FDB47D236C1C479CDFBB6ACD8615FD061A06927AB464E3AB05DD064EE83F64CC78A7F97613AC76DD9E30E0C0E7
99970ECDF2FB6A79A0B6E4CAE01BB2333BF22B129B94D90FDF3742EF2786B39FE4E79E34488503F2EA8C349DAAD42B384290B2D8CA9B64BACA7A29D6
53D1E1DCE62F1464F50E73F33E58F9E35D833F10955695550E91DC2C408222190C63243B173FA850A6A628BE95BFBBDE9E3C639AD9F6A2536DC4BFF7
7FF14F9E8787D53F1B54A78C3779F5CE5DAF4BFDA27DEE4611F25B3129DF9AE393AB516DAE8C6172AE28506DB1942C2846BF94028568874FDDCBB4C7
4B81A791E06F2FFDCDB5FA0C8B45B37AD8CF19736DD998158B8D8C8C5248355AB4F5899D4177E3813D6D1DC7AE9BD87BEB87E1C633BDB7CD39BEEAF5
B5713EEAC085A36FB77B6FBCF6BC54D11C71C14DAD70B1F36155F9F370820AD66405FDFA29D1AAFD1E26E9730E59EC2AE2FA7C599DB3BA4641AEE88F
933D457219009CBBF3C309D613C0240A529329A17194E2B6032224C9B41D773C9629F91629F749CDDFEEF3F5C4BC514393DF6E10920CC5EB9F87E7FF
E91FC4E9DA55036F2A1F15AB3B3F35146E2EA4A9B48709014D82E22C3DA6289731F88B4A021356625058135AADA86AC32A420BA3D2A60566DC52B7BB
44B5CF7B069C9BB577D428C791E70F8ED8736A50AF311D94066F76FD60F3473F1373DB834317390E38F9CE7EE0E8DB0F463FCA7419B7325FFF6A9CF5
79C5EB8DAF2E667546B81F7E27D63CB47F4B32E4EF475998A42EF07112EBDD57B87CD2683C6C57844034C12DFB453002E1A600B05AA7B8A89FBB0764
CF046EE20717FFEB3A41ADB613C5698A1B07C41E1874E9C5E3CF53049E4B0B9A02EBFC0E7F98BA6D57DF778D0EBBE4346FFF7978FE03D589328A5D26
3BE4B13FAAE0245F7DF46B146FBA7FE3B887469EDDD2CE6812944CB23BAAD3D6B761751F8504246E6801CBBC121A76CE2D7AFA1DC83D2B8B9C9A9D33
6CD9967EF65BE65CEDE770D1B4FB5C3D011235C17E7D0CCCBA19EEDB6E39DDE4C4F799C643575F1E7DADE1DA80F142FDF521D38E24BF5DF4FA6A9DFC
C4E5A3EF20E5ED998138837CFF2055ABF3525235180A24BAA4E90B97DAFA68095C85B0C69EBBE553D4E958F1D35DF5BEEC3FB1AEEB40AEEF8744F1E6
5C9D86FD916BF9A5280CA7B0F0332919D158DE9982B67765B18BFCEFCD7F60741A78D877BB9AE1B7FFF2F0FC0732574313AF063A3C3F96246C16BE6C
27AA857A6DCE2B8715298A4E48A6676A7F60D89D1A02042A54AAEABC1C1805EA22F5BFFDC0670BBC7FC4A4E19EDEC0B913EAAB63579B0E5E39F3E1B2
C9470776DF4D630095F5FA9589F1803E731E5A4FB298153AB3F7986927D73B891FF699F4BB6ACFC0719B829FCF0F282A2ABFE0F7E1AD467C6B610489
9606257790D5093900480873BE7FAE575C1F7EA00E27F4FA0E25C16DFDA128580691A44CD655F5DF35F983DBF78B735B7F082AFA3D8276A502499A1B
FE8507AEFE5CFB4B9BF4A041FF231FF53DE6393D68D744917CEBA0D340D7EE1F5EFE3C3C7FDAFF84129AFA6D3EE1F3F57360A1E2D37B8850B4812AC1
85C5DB72C5EDA41CC27C12898C336A4C2395693284DFDD7538FE2118F3A8CA5975EA4EF3EBF68E2319D9B3124B2CFB5A18F49EBBE292DDFCEE064728
1C403E3DB86168FC87FDA30503962E3E606332D87EE3E4A9355F075BBEA858623AFBCAA32B73DEEBA4213762536F893ACE6C48D240E1CF520849FEEF
6231A12D0B28C95124AC991189932088A33049515DB97E6E0C085AD3821144D78D3F97FCE33A81600286F0D43612430802466986628F876F8BEFD57C
D1D43FAB849F06E82ACEBC987BF8F5E020B5D7B191C7A17FDFFDF127000F0F57FE9F1E49D3CD0EA3AAA4E743DB652DF72B719544DF89651F5B195BF3
B18981982C1F2DE97E4109AB9BA1E400E8E5538011DF1665BEC76E5DD958931244FADF53B85CD7B898DB9B181D3F3465F05A538367AC3CB54FDFBA18
9A58DCBF3DD071EDBCA57D87382C5861F9B66CF4E047A14BED361D0C3C36FB13A47F74A929F2465BB5F395124C9EFFA55AADC88D6D43752D7199719F
5D1E6EF76CA23018E5CA7B281202B8A15EDC582F8C4BF6131489735E0037145001129006ED5A0EC26D0A27D8B73360F4A6DB2DADAD85BEAD8ACFDFDA
B4E12F5FAD769B7F55F9656A82CB238CF7FD7978FEAB7F2AC50522553B8D4380A4EB4DE5BAB8BBAD84BAA4A3411D7EEA9A38EA13852A55D72211FDB1
7D52A8A042EDDE403C892419CF879A979A98AB9BDC51F7ACB6CDC1D1A352EBEC462CECB3C16BBA9943EF7F05126A527CD6FB540F13E75BFD663B9A3A
4D1F306CFD64A37DCD8EC67B43B72C3BB8F1E9011B7F4C75E98E24EE5E53DDE9C71D2AA8CAA310C98B096B520055A1197567ADAC9E0A587BCEF5F342
5CA64FA327BAAA7D4098E0B67D70153F5897DC71480592A276929B0480E3DC58008020A9ACED3704ED311579A9CAEF81706D42FCB5BD2E5BD6B62816
AF93BDFB44F1C13F0FCF5FF49FBB4F48623E03B688C4713F3452AD6B909EA86BEB1466B96E0C07E2EB190D9D793E05951FB926EFCC8084DF30F05B27
557FA0B3ACB8C6F7C50671EDCBE22F878B369F03DD2D1D06F65F7F77AD9DD198025C8D8BAE841CEE3564EB98A18B6DD7DE1935E2C0FA81532B9F182E
F15BB8EF95A3F7ED31BE78C3CE5B508CB3A0F9FAC74E50117923032A724D6FC0DA137D4A7316F79F9C8490A806A729AEC317233188ABF425294C0972
FB7FB839FF5D11001BFBE3B01A9669A9AE3B00D627801084D2A4DEF72FCE4AE8504435D7043457A60BDFDD7159393D1E281E11019DC9E1577FF0F0FC
35FF57BA289722B3879B7C96973C9337091B9CE3A1CE5C48109FEB3A3B0F2D4508184FF26D4655375E018D19507622A5AB57509EAF90B8B67779D39F
A299EEEDA73FBC1C50136130C8A69BB1A7C7E07FCD68D60170D1F9CF730DCD4D068D323FBEDF62D2FE3DE6832293FA4E889A67E3B7F0C5D3418FC89A
E9C7D1DFF3BE77DC0F52E85A5F9CCE45BFF9C5A934C91905A59D274D8EB6B2EA47546A1247D8209F40B46C948FD324817377FD5D55BFACFEB9693F18
8CC93BB9EE1F92CB0662ACFA71B8C4C3EF83F78FCA16CDC7B416CFD616D766C5817B57079E95E1A7C6757CB46B66285EFE3C3CFFCDFF35CEFB4E11E2
4D031CDB1B52D3907AD2EF8E12AB95A0E599D51B9D4A65621145912DD53ABC6EB1AF36BF1609AFA53472A6C55753A3785E72C1F627E31D99E1DABED2
A3A0A7A58349BFF9AB7AFC6B9146A7C37EDF8C18D767A8B9C5B4B367ACFBCEDB32B2FF9DC6F946AF7D875D9A73F69BC541AC66CA72B0685650D9A10F
323077FFA906E5F7E046BDF8CD4F2918B7DB2E0427D408864994DC4D3F6BE6B9C95E9CCBDF65E3BBAAFD585F80FD9160BD03990021296ED607C64D00
43B0F6E8E4A6A297D94A4CFE255EFB36BBED711EF4E671D4CCE54570ADF9A9F619AB55BCFA7978FE1A0034DBB92010FCB0BFD9AFF696A76A545EBAFA
3BAACE8791B4CAC471A7C4B282060A409166103A3DA9451954577D4344B7C9B1D84A8D30EA4D9DF5D2C6FCD34D4703AA4E166F349C66D567D8A03FFE
B55CAE6E47BE9F7C6E6E6CD277A2DB3BCB3E53270D353AD87AAEE78ACFD647372F8C1E37B5533477726BC97C9F82832F14E08FA5D71AC5DEDECD98F4
9A0FD6E63CD6FC3D85803A1491A96882A1B5F53A82EBED01D1AED99E5CA32FC11E0B044D00ACEA71580ED138D7FBC73A013082A2A59FABE88E0FED0C
0026FA367CC969B9F89B487178F370EC6B3D75A4976BE6E8B3C0FF9FFEF9A381E7FF9CFCE9F6053B642091623BE446A93ABC9AA9151DBB0622C20614
4C135CB2F7D5485340AA4524AC818AE7BD87EAC254FE2E1D402D56E88DA675BED03C36776A395F1F364519FF227ACCE429FD672EE8D7FD945EAD4033
6E5F3137B3186C77799A91F5426BA30D799F060EF97264D4F9950FB68E8C536C1A999867EF0D5E7CA9838326DD96163D7DDE8CAA2F9E55D66D1D3029
0401211CD1B31F8921E9A64FF55CE93FDED90EFD69FDBBEA7EB8CB00AD0A855B5AC1AE164004E3B201300294263588823C840C0147DE2FAB6FC8BD97
8168F63B05EFBA2EA6EACC879744AE7B01F222E7E1F9EB0120DA3A3115403AD69A2F0F6852F9138D55396BA36130A711EFC8C8765C5442FD8A64943F
9591B5C8CB250AE8E367E4FE33A8AE4D7AAAB021393A45386944F1B3FBB0FD45D17ACFE0A1138DAD8FF5FBD7714024019B936E1A1A590C1E6E6B6935
61589F19D18FAC4C0FB80E5B71C16BBEE907D99E910185739EE8DD9F8390EBD47775C94F5E5529D55736B6E7380ED85CCEC042855EAF56B34E3EA9CC
6F2649042371A8ABD907C7585F0063CD3C410032989257EAB857090C605F67AD3F505F24953DBE27A2112CD72D4FA988BE570CA10F36259D5AF34BCC
3C30BAD079F4F00BE87FAD3C7F20F0FC5FD6BFEAFED8A75ABDCEDD7C6A6820F3A2062DD079B8E961490D869755BC9E715523BE2260CAE455E9886CE3
0BA8F35683FC794A713D16FB1CCB8C7909BEB40E6ABD86DD192D0A599CB96CECD0D9F74676BB82288559EF125CCCAC06F55BBA6B98EDD40193BEED1E
6A34F18DE308BF4BABC7EFFEB9CFD2A375DF69E2C779157A7376A4CC3D26BA52AABEEE58983CCDFABE862124156D28882024D8A4D6C37457B51F49D3
5CD2EFCF493F7057E70F8C425A8CEACA09A070D72C60597E1D22CD4D5511049C1DD606E8A3CED78148F6FCC80F4702953AED4CFBF6DFEBDD7D30FE81
F3F0FC55FF6AD729CBEAB4AA448B019967C1EC4F446D62E7BD5C18CC6C22755F7FEF9E900406F990FA8A768F3630C9E6239A73575B9295EEA7D69C2D
D1C6DD8BEF70980B04A7A70C78DC30E4D2239381ABCE8DECFE1052A9134E853E1FD86F92FD3E47CBA97D0CEF799A8E1C3479CBB09DAF86D8BDB8BDD1
F662F385A5E2A4B332C5E15515E5B77E7E8F92ABAF38A4478E5F9F43326867791B0AAB29529ADEC0A99F2BEB45B98C5F97D4D97F83084150ACAF0F00
28571480E308074594F97C6B90D6572971589697DD8C81BF4F952048E7AAC86F339FD40A183FB3506CC35DEFCF247FF9CFC3F357FD032FD6CE0B02D5
A5F34DF36E7E2562C448B222EF7B3DD4908433B981B18BCFCB343E4A22BE313911D09F1B9D8B857EC848518725A0D167C1CA107F20707C75FB23DDAE
05A263334267D82E99663A240E57C24DEF726F99F55EB372A1C99C49866B3F0E1A643F69F9E07909F3A7786DB875EA707DC1A42FA2DDD9BA6DF3CB9B
4F7FAB089300B7C786BC9AF844C3D06453858CD0A238D199DF4EFFB9D91B2349EE1E807300B82DA07A48AF83343A989BFCD71514C05C4FB0FEE7A76C
496D45330C6092CC0210D4647A16A32872C757B0EF64BE1C924CDB837D9CD37C3195973F0FCFDFF27FD0F525AB368AE482A3E61FBFEF22EA22B1921C
28EB690D9699476B1F44DD191BDA529C80FF4806E22BB5B5B356A8A4E1C9C75AF42F8B744743816897D4B6257BE5D7725E1A45A5195FBB6461D6B7BB
793E21052ABC63A7F6341B3D7AEEE871FDBA7B3DEC3ED070FFE1EE1B9FCEFFE834E9A6957FE7E23565CBEE68F62CAACDDBFE2BED71077C65D4B373D3
3E9304831796EB642A8000DBEA019C42508A4008CECC231037E68B75FF111D801447B7225D9700DC9980B2F2074038330323E0262984C11D1995EC09
121529624F88D88FF0E939E94DADCCBDA9056DD33F4A9714F3D53F3C3C7FD33FEEBACA7974825C1E30E858C5D23AE447657B3200FE8A82C006982971
CBD9BBE3BC30435CFD031595A991CF36DFF17CA9770896FD1E493E256B7FF35CF562685A617CE33827EDE28D1196C67D7A8C6B2501A2E8AEE7E07FD9
F4B159673762A4F1C15583472CBA31BEBFD36A9F7B53628FDED17E191E14B65A943C2BA366CB97B4F325A09BD5C97DF3D258275FFCAD1425857A6969
8D1AA330B8EBAE9F00245A821BFACBF9FA904885137201427545FE5C15008AE0CD691A042149A91026604C95D688C2D2947400C789AA4CF4A5F5DD6A
0991393D083FBC9F4C98D5C0AB9F87E76FFA279FCD7C61B35B20291866FDF3F023A4294CF2F5830E9137A1700B85BF097F3D6E5C404B842EA8141514
6925CB57EBDA6A545E1AE2438CFE8C87AAE47E5298F9A58AC3ED7B0C934347065AF636ED315F89CA91A247E176DD7BF45EB6D868DE42E361B6FD0EF8
4F3358F17EEECA4D0BF6F82265130FBC9E99A05F1D5CB9DA2F755936F670C8A1A53BEB089C2AFD5049A92498202256472024500A532845E1B8468573
B3FF2886D46B5BA45C2D20D70CF0E7E04F129515C7C8281AA9CC97213081CBD34BB4089A9BC6060A4473061162BD2DB95CADDA7E15485BD1419D9A27
E0F5CFC3F357FD5374D8829025C37FCB1A37F63FF56EBB5C9F9ED6FCC847A2D12378BB8291A6BD996B3F3BEB7D5EBA170C9434201FEC73E82CE06404
2D792D8F5C1B93177CE2E986F1A9A7F644F5DB22B27B7CC4B0AFE10D0C06B12CAF9245867DE76FB7375B693771669F51AE7B7B9B3D7F30CCE9919955
49C5E4B5D2339BB58FEFAA2F7AC86F65C1EF466D72BCD249E0404CB8926034EA8EB26A98DBE3D71A87D15D457F5CAD2F495364638E0617B7415C5760
D7C83F4EFFC0AFBB5FF514AD894F53B24E02AACCAF861159451344AAF450BDD677C69A5F7A09E0BBB35276338B691AB34FCAEB9F87E7AF90E4FBC92F
4F0F3AD9D6FEB0EFC8A01D81FADA60BDFA4544B2825651A44EDB527275DDA8D0DAABCD1FBF22F26F8D1D4B37218D199F9634D3BF3EB73ABBBDF67D75
E0A04964DDA49FCE3D22360D8F9A67D8CF9F51C154F4B55FB6FF9AF0D8C178E130D383D3469C5CDE67C8E94B661B9A4E198F3FBDD82E347A5654F081
9AF727453EEFC8DCE18ECB1F8104260B4EC6183D2C2C2ED722388C60A4A4950BF83112C3519CA2710ACCA8A5484885A218C2AD00E17A0069FD870705
1A12C1340A9C1B12A8A893EB5140548BA13A04AF2CF75CB6E18B5C00676F8AE90C4D26C91B03DD005EFF3C3C7FD33FF1FACC4DD7E193728471D3069C
7FB653AECB2BC18A0BAA7CC4783346A972A112F78577F5016FB5AE5558F677ECDB8418ED9332EF6B88D43931342CFED417D795338F03C1A119FD377F
308FD8D5C3F003A34308BF4771C3BB2D3866347C97E58C45831EBC34EB39C77DFEF0DCF2E1962EB79CCBD2A77F6DD993DEB033B3E812A4D93061D62B
8880DB02CA2898D295E5CBF524A2D377B5F7699428D7E7C7FED572433E69124089AE997F9CF8B91B412C3D4C4DC3108C72138231A846A8C7D1F65215
46B00E842CE4C9E95BDE1A5441FABDD67EFD8CD21513A6FC26F8E7CDC3F357FF9F441EC69F7CB46CC0D39632A73EC3DD977D528BA2DBA4C5D8DBD74C
761CC354643655BDB10B93399756B84BD09074F9F6B9B2E047EA03BEF4BB93EFBC447EBB7E6CDD30EC832EA268519F483BA7C5DD7B06D25A407BE9BE
7FEF7F598C349A386D90DD88ED9E73FACD5F3CBBEF53C5DABE7BA3CE44A18F3CCAF7FF146EF82A5E5B42DF19BAD00B05052D3F2B111C1717B72A690C
856ADA70F6AB4E20EB9AF14790EDB94A6EC3875EA0F973D517D7EDA3530268C36F80824118A36092225B735A00D613482CC73084B5FE5EEF02CF5FAE
812464C695C6D207324A74CF66673B7FFDC7C3F357F933A4784FDAEB4BA7472CAF687E3BB1D7D1F317959AE2745D5149736803E455CE100DD98D0DB3
E6D7A5FB409F7F31ADEF3AF3E705A0CF5A0B7714B53EF53B56DCE0FC2AF8E08CE995B7033E0F76BB3B77668F81193480A9EE7CF13736EADF6BE8AE39
E6E6F3A30F19CC7CB460C0B4B64FC693B22E8F0BFEFDA0E9EC4BD5A5D7DAB3376117CB8B7958C78F6C41252CD308CA050C496398BE4E82139446AE67
4D3D7B2A20E591751046C2380A71F300B896609CD4473ECF914B11D6F747B86C20033465376024082A1B21420341453E11911FEF15320088DC0C4583
92487D80EF88F7BCFBCFC3F337FDD344C3E4E7197B76CF358F6E4E3FD067D89999B15A4D41823C3D14C9916ADE3563C2E62FE207035C946EBF3A6E17
30E9FED0A78D821F4FA88867D5EF93F6DD9444ADFA7167D6AAFBAB5757CC9B1439C6E28F610DB81A697B9872FE8FEEBD4C1DD78DB2EDBBDA7BB4E9F6
F5FD8DDE140E317BE66515245B76E06D54EBA1E7F8A3F9B1CE533EC2FA9AD0ACF6FA7A549C59AD42190AD70841D6BC838452419008D7FE2BCC6A6575
8EA25C2100FEE704001C52A4A488B86B40840D09309A02ABCB35348E10A806445108C9FCF0BDD2E7650AA041A8E0BBF2DF5E4A28B6F2F9446E9B110F
0FCF5F0E00BC61FEEDEC536B970D7ED056FF62D280DD47D6B5219D5128F6B69CE9A0AB0B810E716248DB0EFBCF2DEF3511B720262C9DF00C80BC4A81
BCBCDB9FBF2F4C541D3E9EBE72B247C0BC1FDEE6316B4C7B2E97C33AA2C135777B0F63C31517868DDD62B964419F2DB7C7F45E5AB9BAE7B98AB53BD1
2F23EE26E55DBE2D8C5C7860E9C23C4AF5E2478D36B352274D4ED33238292E68025192245859A310EBFB4300A047489033FC6C04D035F59724445502
1540D27F6E00E432038A0A2141D10485EA5B5458BB5A9B51F6F34D62A68AC2E8FC6BF5B2E01AB8A1015A7D50CEBBFF3C3C7FB7FF649EDD991FDEE7A7
8C5958D49138B18FDDAF55575038AC92C9BEA723EBDB6AF5BA4CF9ABEAE8698775F95F5ACFBB52F2ABC5EDF73ACB032BEB34214EB2935B5411A3A3C2
2C9C6B67AFCB1EEDEC6F6CB05DAD96A1BF0F844DEDD163EE659B7E4B470DDF31C4FCFCC61E23BE5DFC6363FD8621C582E9FB5B921CCEE89ECCBA3EDA
2C9969DF1E04B414E76A3471516AD68E177F6F85198C06753401B292D708143088C310822204D7DFCF86F604A2CEFB5AA501BBA200EE7A10A53049BD
8264308A243B5A1568672BACAFF4F99AFFB3514F30F9F71B556E257041055A3BC21BE4F5CFC3F377FD53997607E2FC3FAF5E60E52F6B3B30D82AA960
A93F5D9CCD74DECBD6EB243EC958555DDE876A97ED1FC188A89243518CFFC6ACAFAF15593EF7F5BA1BDF4A6684CB36EFEAD830AEE0FC84DCA313D267
F7BA05A93564E8DEF4A586439CE7184C5A6534F7D8F0592B7B0D3EE66632B1C8B79F377AD63AA3F394738ECF8943CB4C0219E1F2A7B036250B688A4D
1713285A95266750066FC8604D37ABFEFA1A194440885E0DE124C98DF7470914EEC8496B01099C26BAFA02481281B5420D43120CA96A13E870B41D46
8B3E0495A489288C4ABFDF48BEF7C3C06A907C393E07E79F370FCF5FB37F3445FE1A3EFD9B4BCCF139266B1BE59F06F63809A56F2F6B88828599B991
3A34E58408C96C7813FF70FFF46CEC4966D6F126918373CB9D4469C89A243CEF54DDF1458ADF16311F7A3C0F33737B3670EF84EEE7B54A19FCF6BE60
5E4FFB49DDC72E1D6C71688CE1F63906E30F8DECF7FEB3C921F1BBD1EED215CB5A1FCE0F9B31F433D3BED68714A714AAEA0BB301D6BCB767403445C1
2DB94D34C550D2AC261D8180B046DA2EC1B9661F1403317D6D9354DEB5E983AB0126290C06F5FA3F7780E905ED18A915EBA1D6B0E0B4947A0CC7128E
965185571478713E5E67B15344F2E93F1E9EBF9B7F32C36674D8AB2F770FAF1A1528ADDF643A219148F6AE6F48F8A401DFFFD2C0510F1512C1C7F0DD
FBE63B091A5FCA7E7E424237E7FE72ADAF3C7D408E7924642F8C40B7CFCF98BC2DDFFE70CC1CA301065B35804AE31A5C33AAD780DE43772E36D977A1
DFD8CB93E66CB532BF153D627D92CB8D4B4DEFC6F8D5ED775EB8B69869D9E94BE853B3140DE9121D8A08C5ED4A1AC3C8B66A8C6170303B5E46C03004
AA55AA4E15C535FF6940145357081108E5BA7FBA7A00289C3D20201245090651E8290456E8B1929474FF44804409A95701867D2B22B1BABACAE5833E
02BCFA7978FE790054DB5B79BF717F7ECDB3EFBA22D9BBBEC63352E8A8143AEF603B51F555A6D3BABD6A86B2BF876DBA30F222F1C0930849936F79A0
0F7EDE123DC7054F3E5CE2B4419F30D4DD6948E42C87904D7FF43538A85529A5FBEEB81B1A190D18EDD87FD6237BB38D76C32E3AFE6B49D2DCC5753B
D7B77C0A1A1F50B379C3A4A572BA6AEF17BA21BA0C294BEAC46034274B89A004A013EA20F6633586C4EB7056DC3A400620A04A83E2D256698B164374
4897F527BA368270537F611446601CC3F4ECFB153006172427F8A568719C5085C769D0CC6F20D05A8914F59F5685F3FAE7E1F987FF4F670E9FECF730
207ABFC7AC09CFA53553CD07EE6855144298A73F2E3C58ACD389EEBF53E99205218F77CF8D4E9CDDA8BADFE17F4ED5F620B2C27946BE64AF7F984339
EEBCEAE384A707A7F95F1D626A14006921D976B7BB3DFA993B2EB49F70E944BF595B67CEB31F38FED2E6E1F9DFE7E4C55F5DE386C42DB0D92D2632B6
46E84ABF9493F589120C05C342550889D5E428401A6368A8BA8E6020851044611D882232358C74CAA485CAAEAD9F5D793FAA6B2310816B9B3B0110A1
30508FD3B00E0572BE54F9C4E31448E83F67208A823801897668E963035D147CF68F87E77FF49F603CEEF3CDC769EBCF7A588DCF573E33B1197A9568
AFA555EF44F2192B34324D717814D2F4232930F98063FD81D38CC7F342C7AF40DCE9745FC75D328F85E1D38E9115D39F2DD9E837E084DBE06E5685B0
44233A9178A95B0FBB63E3A72C5D6D3373A4EDA1A93DC6784EE9FDFAF794A09C591F4F87842D325F25C39397462BBE84CAC9A25F628202DFBCD51010
5696DC49C00C05D477E8188C94C757C11004E975723581214AAD4200E230B7F6A36BFF27E7FBA3185697DAA647611401B46A3D85774446A605E7816C
14017E4F126BAB130484A245447F3599568CF0FAE7E1F9BBF7CFEA3FD77CFC8F474F7E5C3BF57DA5950F58376BBC9DC32F52A826DB7E6B3EF47B0E49
AB44FEADFAD2F00475FAD28BDF5797775C94DCD925163C79937FF36C5DE7BEDB27C7FE240E1C7FB6326AF4BAB0997F2C55C15AB0E365E3D13F7A8C72
303FB86FF0828B5647AEF77778B3FD8F55ED1BB7956DDD931F387FA4C982362C7A5576B5FB7B099AF95B4E309D6FA359175E55580092140A356466A8
30061597166A703D012B753A2D5C195704A0DCB66FAC4BFDACE74F113A31EBE62B1418B70704C7F57A026FF6FBD89058455010210B2B460059A91C85
E472A26DB285878EE2F5CFC3F33FFACF186AE97DD5E35388D3F3EBA6CBA5D02BB3456BB68B8846004D4FD71D1A528854945746EBF52187AA7427FA3F
79E485DF7B5939E33158E89CFAD12959EB64EF6F37B3FCF6B8E8A5F7D75A7C5864B05CABD3C0851732EC0C7AF4359ABD75ECA887CB66B9984CA8B8FE
87EDCB7D6BEA2ECCA8787ACEBA976307F1794663F3ADF7624D422EC23075EEBF31121167D602248029CBF2852A06AD4ACA2D82B8F4BF0ED6C2C2F4AF
B97256E518C27503FE7B1238A6916B516E1A00FB85444908A12ADFD564BC6EA1082D0DFDCC4660457E1D2A2917A2CA0D8316D712FCE83F1E9EFF49FF
D1A94BE678BD7F9610FAD8CBD5724C2A56EB30F5A8E353B54C880259A2D699EB94408E38BFA8A1E4C01B2CC4D4FAD11B65CE89A6FB365FE14097E2F5
3BF2DDA73E3932707F82C3C763471E5B5E9FFEC7661002D1A85DDF0776EBDBC7E2EC14E383EEC62B960C0A0B341AF4E4E5AC9F8A83CE31F71CFA9DD4
92DFA765679F091374C4E5E30CFDFD56218D6235190D388063A22AA996060435B1E94A9066558E22B84A2610EA483604E89AF3FBE716009C641D010A
462118463112671D7E4A9790D6D1102BA469822464AD3028AAEE8424118D20746AE0F4CF10CDEB9F87E79FF69FA4BFEEBEE55D7326F96558F08BDD03
F76A548F2C576E9C16AF2B1723B29FEAB4313E48DD2769699824747D8B62A5F1E6FD49E025CFE2E913735597E30327BA662C99F46AF5E4B2DD3B3D37
870D5F65F5C7721054A38187E2FB1974EFB97397D1D05D2BECA618EECE37EFE7E4B2F089F0EDA97B0B56F43D4731250B8B6A0FBFAC15C466220CE479
B39E64948989620002A1A6261186AAAA6A406E5237851208835644D6B2FE3E86229CF967ED7D9707C0C5011401A10401615D7BBF29A62CA0AAB285EB
682429499B1EEB4CEF40A1D8340DF0C574F089369CE4FD7F1E9E7F1E00241DB63B6497F4B9FF8FC0D498E051637E80856B676C997FB8419C2D418B12
01EF9D427D5C2E1C9B853D4C66A2266C5E72068E3D257C6ABA4590EE557BF25C93BBDDFB6F6BB2C2E6BE599D3ADBC6BAD75E1CD6637E2E25D63D0CA7
390F32196DE1B869C884AC7D866BFC362DCAFF3572FBD311BD1DD48CEC5841D956D7EA8AF47498427DDFA818489E1025211430D2520B006871600BD9
359988A11902AB8ECC5577E91DE7427FF4DFF69F0D0148480953DC643002C7285A9392A10120EE3F45C9046225A4CAA841E18258582D5C31705A2A8A
63FCDE6F1E9EFFB5FF6F1D5226E4A59CAE0F4808FCB475E82641E74DDB896BC73D52A57D50E89332A54E1F7171489BD8B3A9F2A24CEB34EBDAB248F9
AE80DC1966EEF08DA8B045DE39CB56A66E3E56BBC779F90BC7DEC37B1CA70008BB79E6B785C1E807630C17984FB960D5FBE517C3A9FE5B076D69BE62
E6B4DDD4389D116ECE14EFBEF5BB29FAA78E24C33EB3F65B5157D801A9717D7995146D7F77360BA5609C2248F603EA5BC5AC6B8F6004D6D50284A3FF
9EFA8F2330229391344E120842E2D0EF877124CD4028CD6865ED706B0BD8580C41A5DF3A1560BEE5A0172A8C44095EFF3C3CFF380128CAC341B2ED99
EA62626E4248C2177BFB4F9D710E93D72F9C1158EFF21DD178A67DDADD8A55A5A21901D4872F4488F5F55B07A4A17B9B1F0D5D50F7F349E3CE3595C7
2C13EFCEAD0DDAB4CD6DA3C530430F860DC76FDCFB3CA087CDAC9E430FDAEEDFDED73ED8CED4DF779843D015BBAB474D877EA21A57BCCC5AF7ACAA21
2F494350B9EE180D9796290112433415D5FA86AF276E1777F5F85314C3C8EB8510C1A05D293FD6C9EFCAFB616D75ACD70FA3380A203445B2BF40A0BC
5BBB02215A2984185C26C42145AD5421504ACA2B3B4065EEDA114E6D040168F8E65F1E9EFFF10090EBA3EA0316C902AECB7F24C5156DB05A91DE7874
98BDBDE59AE49CBB6570DEBDCAD32760F26724FE2C49732E55E8BCEDF385E7E265BED50B073DAADEE67BDBFAC3A59E37BC8C7DB2D78DFBB0AAEFA0DE
31AC0916EC767331F8E38FEE3D2C96AF5EDE6786CB983E27FDC6595C4D9EB0D87F8AD957B264EEDBC64D2E55F5BFC20518D97ABA8026D212C4848442
04851D9DA5AEEE91529C608D3B89104875B11C60B8541F2B76BCEBE29FC030A855806208CC4D07E7820012D0C319470E478944F5CD00C3C00846A95B
55A0482DAE15801AF5CBD18396150310519BC1C608FCF3E6E1F99BFE19BDEBD842C9CAC08E8BD515D9F773DDC74CBD50193371F822FB0901F2787F15
1C912539E443296EFFACB80BA53F6F4BDB75F0D3CEE4A0137AFF61679BAE5C8B99FD34D2FC70DA889DC58FC687DDE9633641CCE08C604BF0C39E8686
3D7ACF3EBC6AF4E0B7270D17C54CB77E7575F79EFD338706A279F6EFB5872EB7892B939A215CE31108B596A4A8501251E65674E0F91EE96A88B5EDDC
7C4F799B3047C170A5FE5CEFDF9FB69F20B1AEF23F0C4351929B094EA235E9EA8EE08FF51A7D470BCC9D670C2152E8618D582E020105F875A8E5FC9F
7AF6DD89B17CFE9F87E77FF42F7D3FF00310E22889F655575E7D5B306DD834EF9A5D7D2E87D8AC2E56BC0BD7169E100A2FE633A5D785611F55FE7E11
BB8D0FF91C111C0C91CF9B59F16953D69A556956E629732D23C38638C71B1BED630898116D0D3EF247AF5E06BDCF9F1CD0E3985F5F9BD04D7D8FC64F
5DF9648EC15D2471B40B7C6B6DB13239BA1AC2947E4970EEB33435ABF8E692724C97F8265B87C9F428A4C771794D47477D931CA7701AC75977A0ABDF
9FF5FE3198F5063002C138E38F11F288E4E68C5A367A40219CA6581F9F96B543A446A4D5AB2025526ADF6F4E8C1A84C8D627E5346FFF7978FEA97FC9
A731D75A15133F6AAF9529636FCB5DCDED9727C60F5921B862710FAEF3AA018E6D4665EFEA9984AFBA1F3FBF46DFDF3CC8F4C59DB0A4FD3ADF3161A2
0391DE8E799B8CBDF70FB913B3C4E1978DE16A808098AAA3E9BBBB1B191AADBE3AC764F2F759BD6E3D1872EAE2BC831EF34C5697C60C3FA20B5E99A6
2B88A9D511821F5978CB93EF22086F2D2A0741859F575A6B9B4407EB5835CB3B9A25A2D6663D41B1E69EC6BAD60073E3C0610844BAAE01B8D11F340A
A9948D7175040E6214EB0CD014434B1AD92343A6E10A82B1863926C3BCA50000C93C2F6AD85FF2F0F0FC5DFFA247B647E3B1BDCBD084C7CAA6E4CAC6
1596639DF2CE9BF8A6DBDAC601B12FDAF2AC6FD0AA3011ED572A733B5D9772F7F480C9411B726F7AD78FDBA57D71376D46F01383B9D7AD66F99E1D70
D7AEA743238ED3FE2B9396FDD1D360C28B197DFA3D79D47DA17B1FA758F3E13FAFF7DA29BA6E71421CBAEC872E2AA816269AC21BF11697300D0A1424
3629E8D27BDFE415050A1824105C2B6BEC282F5183080921DC2505812304FB8D3D06D8C09FAB02C25098A0294887685ABE15B12F72D703244E51A4AA
096148A5160701002F9C69DCDBB955AC81253E0BEE4114C93F6E1E9E7FE85F7ED7EA8EAF366C5427F846D096D1ACFB31D476B25BD2C2193E57476D6C
947F88875C07FBD1AA84BAA6026DC2B11CF5A718A701A78E1FCC395B72637C71D185826D0FD32C07BD586EEBF576D286A97D2E697082F9B4B1605B0F
E3BEE777199BDA3C1F373078D6BC5FF34C975CB13BD2E269B94F9CB82D521BEE55079155C142B2F545A41686E202DAF1F6EFD74210AC498A1328DC5A
5F2BA94E4E69813104830914A5288AC0096D07D8D5F3DBA57E0C410942DC2E6AFCE25F8971D1002B7FEE5DA80E6118584BD0084C17D9F61DB8BE5AAD
829ADDEC9D424892B7FF3C3CFFD4BFF0E5C053EE6972BB1798FF2D22E47B93E684E998C5115E83A63EBB30EC852EFB56BE789F69300986E664158B4F
EE14E5DECC5C39FCED4CBF9873A9165734373CCFAF145CEFB675A7E9A9376B4759F70AA14090C9DB9FBDBC7B8F45CEDD7A3B4C9867787DDBA0A2CD7D
0F2E32DA81BF37DDD919B2294E1B78A10245335FD6E1AD773FA8095D805F6B5DD997FDAFF592CA1A2904B7274497D654FE8A68C0511844505CDB89D2
6CE88F1625E8591F9F8D02D85301424044D7562F2D7D7AB79522615403100446922088620C0DC8511C4718C16C63D399392A15D474C6D6F76A0CC9DF
FFF3F0FC53FFB4E8EBB2CD8941A8CB2A24754A714EB6342761A8ED8C33093B473844DF704851477AB6D7AE18FE1890E467BEAAFCBCD05FF7F87DECDC
1B0FD7162E79736DA220F45ED4FCC21A879927ED16FBDE193FC9309CD2434CE4C1F2393D7A9D586460B37984D534D711AF5FF69E7CD1CCA22A73D8EC
C2941DE1E0D7CD695A38F3CC3769C60D5F2926F678DC511F94F1F375B9B4B8BC0521E4750D2AAC313E57D4B5EB07C3484406D3141BFD03EC37D6C693
2806C18046AA11D5350A1BAB94DCA040B51AE38601C03A80CB0CAA11068719F4C8000B875F2A1828701A999432B798C6F9E95F3C3CFFD47F43E401DB
3CD7EA8ACDF98847A8FC9B32A6EC6CFFB176616F470F746E3EBDB7B4E9D99B8ACA4BF35D0B554DD75E8B4ECCAECBDF1A7E6BCAFB136F1F8F4F9A7029
6D4BC6AA83CA67463BE6F7BCEE3376A86104A1C5A8F3AB5247F5E8B96CA0F10AB349B347AF39767BD09C5DE6E65FF2ACA76507ADF80A7FDB10A70692
CE277626DE09D5508D57DC15396F9A1A82631A6A9A41B542C1C5F18D690D0083A1086BEA013620E0727D48D70A50CEC5D72A348016902B33BCFD2B21
1A4741040630823D2860991A82081C62688AF5146E5A584E895600E2B825FD22B48E0BCB0892973F0FCF3FF44FE6865C9858EFF74BF52B9310856832
CA9A52F2C65B4F3EF56BDDB091994DCE6F0B92D37EA5E872DD1E244109D79B8BE73DD2BCB89FBE7CFD957BD2CB672F1D6CDD7CEF0D6BDFFBEFDB3F74
79D0FA3E1679B896A08E5C2AB1E93DC8B6DFDAEDBD4F4E727CB6A1D732771B93E7C5738644D4AE7989FC5E1DAD06C28FFDC6122EFF46E8F6BB6FD4F9
DED58A8F6F4AA4AD242D6A9291B82A3B474E776DF9216090F502BAEA7D516EE2174531884CA00620189564BB3E2B01B95D800857FF4BA024AE96EB10
983D20189A7502AE0D9D34F79D54567E6E989993FA8A8D731DC9EB9F87E79FFAC7220FB8CDFD1EFBAEBA2EA68174FFAA0DD3C6563F349D3E2BE0EE60
332720EE74757250876B04A32CBA9E87047A292E8DCE6BD9F0F9B2C59E3DDFD53B17EF2D3B362975E5F83BA3279C1D6D78DFB5FB901A5C8DE9D6BD2D
31EB39DCA4E73EDBE58EA67E27BA8FFAB6C2F066FA642BBFE495F7E0F4755190C2C32957F1D32503A19B5FC7220977EA946F3DCB8435302092E1245C
9D548E12ACD167FD00D6B083140DE9F07F0FFC255144AB017004537F0DF855A5232894B5EADCCD1F9712542B642A904031042118D18161936C9F0B3B
82160E79B02CF19BC5D2DD35047FFDCFC3F3CFF41FF8F058E53E3F81574D7B6E3591E3051596897E642CB6B43FF37D9E95559CEE895B47704BE18204
14FB794D27F66B2F5B7040F7748DDF1CE7A3133CDF8D1BFD2A60CC33BF89AF2F8F7D38ABC7FCDB43664B512DA1DF939F32D078E880C98EDD178CBE74
ABEFA0FBC70CF7D72CEABD5577698346B4F23DD27466E16F38FC4E01A9CCF6AB00131F57343C7D592EEC10E9C57A9A067333B424C555F761985E0F51
24A928D7505DE3BE6810C07098DB005AB7FB5839841308CE0DF520BA80C4ADEDD50A9C24710861DAB68E9CE1F8A4B1F48AF5C4EC9A2B9A034B776EAB
A178FDF3F0FC53FF9AED7B6A8EB9EB9F46E7B5355410AF5380CFD2AC88D0D183A77CB86D6E3A4F9EBEF157DA63F96D8746AD6CE74B302E52F1C8EC93
78D973E7C957465AA5BA0F5EFB7DB6CDB7BDE7BDCD0E9C37EE6BD96D3E086B08E1FE7A374333A31E7317F5345EFC6A54F783377A2D6CDE6638E9A99B
73AD7AE715B462CDDCC406DF3B35449A7F460D19725B28F4742D5617372835380DC1CD6D088D73BD7E18A890A2A8A843D1A92339FF1EA324A14D5C81
1F43152E3A2FC431882B03A47092C2BB3A0160442F03B9BA60906C5865BD68A567F2BB59C376D5E24F321A262C9EE950C1F0FE3F0FCF3FD1ECDBD8FC
FABE36C659500115E00D2192E24C5150F8ADC123F67E5F6939EC42E30B3FE1AD60F19A732ABD9F5D83D0B721DF6E4E6BD8995B63763B8DF3AEDF3029
F292D9DDFBA7E2672D8F996AD0DD60158AE888D65B8AABC6C34D26EE1C387ADE56BB9E33EEDA9AE77BF41CF1EAF1B418EACE4E55D98A99712D17AFD4
E35FCF54D39204BFBA1AF710A1B4AD5187907A5D934849615C613F0D03A8460E74D629299A26B94A5F0C6FFCAD211826EDFEB553815A0C46BB168153
5D5B40BAB60072BD4108D71F9C32C972D189ABCF160D5BFA5E48477B11CFC62E9E79B496E6CB7F7878FE89EECAB0A2A0CD9D2D73F27E76C80AF1D854
38A2B0D6F9DD5DB3491EFEC3668CB819B4A5F0E79E9A9FA3DEA135C38F815F1FCAAE1B5F971FF9BAC2E8CA9E8929C1234FFAD9CC7C60FFFE86B9DF35
0343833520A027CAF6951D36EE3DEEDE06C3793716751B7B618CC1CD8AC1863B63A75D0583163774AC9D9351B0FF7C07FAF5542B25F08913D6BCFBD8
0A5795776208545FA5C058E75F2BD4900842C3E24609407415FFB1D61F83104EC0D4E771235C5A701882391781C0700A87308AE46682A26CE80FC208
11DCCF70EADAAD7386587994CA81DA67C28ED90BC6CF7C56C1B7FFF1F0FC0FC0CB61D9D11B3AC9BDB784454C3620FEA129FFA17EB93AEFDCF0B5096E
53972FB97EF0A5CEEDADE2CEEC12DD699B26D1D5E49AF9D679BE17D2AC163F9E7C26FDC4DA88D5669E0B2FC65BCC4D703430B8036B41327966A4F3E0
BE73769ACF71DF6E667566478FC9CD7BBA4F78F4E05853C78210E4F8C4A4F63D273BD06F47EAE98AB759707478BA505B5DA804315956851EC3515294
DD8643340C09EB945D8E3D17E563B05C06B1BE7FC6B955677330028638B1B3AE3F46D03880D324B71080DB0ECA1E0BCFBAF71C38D1BAEFD093A5D222
351D9E467ACCB8BAF1484403C25FFFF3F0FC33FE57DC310D69DC120B7BCEC62BE44DF9584E0E1096D2B6E46CCB9A21B7EACEDA6EFEF269F9CFD25D21
D5F356ABAAADAFC2114E2DEF06AC95AD7D73C6E8DCC6D1C191D35C6FF5DE72CAE6DB74835DB116DD0210B5060F5CFF6399D180517DFB9C721D3260FD
D9FEDD3F86F734BFB77B7CBA72D369E8B94D84F8E4910E3CCAA91AFE79A18428BEFBA346D119D78820F2B47431AE27F09AB44E2EF727ED5440180141
9006420058AF6D5312A4E0CAA83DA90843A2088820DCCE1FF664C05190ABFD615F4331108070F76E3D069AF51874368FA9FDA164727DE182C93BDE2E
381A51CBEB9F87E76FDAEF6AFF7968F98D78E605B74F8F97D423E1B02A112AFAAC0C1F951B673929AE6AEDD81BFAB0DB99C1E7CA22863F821E2F9474
5EFB58BB794C6ECE8DA4FEB38FCEBB52E8BC2D64AADDA789F76E99581C1A63F001D140E063DFF6D97D2D46999A382FEAB7ECF4B8DE7BEAA61AAC8B1E
B552F1C1BEA472DA73ECD20E011DB1A6020BDC95A04C7C951AED2B282B4371303943D72C47E0E2540D8DE37A8102C728AED7470F21B04AA302E5408D
DB58FB1098E4F6FEE230CC0D0305108C0DFE099ACB00C05CE50F8A7F1E69656D687CA284A222FD558CE68D1CDA3AEBF59359210FCAF9F67F1E9E7FEA
9F6EF4B17948C43D10E0CB276448D08C3434BD5CFB2913D8B65A73C264475B92A3CD0B2AEA71C9E5F38D1E831364CBEFA191EBDB52AD57C3CF23E6F6
3ABD6F6654CC74D7B3267E876CC3C70EB136320880557ADD6EFF068B9E437B99AE5D6E346CE3D2DE96E9570C265FD93EF347828D9F68DD01D5F5D5ED
CC97E5656ACFA3F93A9FC7D5E28B97EBAAD40456972DAB2D52B527E5B241BD46D5AE6268048531148320854A8D00B5FE87AC6D5CD40C03610C7B3E80
1A6EF50F82A33842D238C2E0801EE7028267BDFB0EEAB63C1727146EC110467FFBCD3CB0BEE0B1D3EDD5F6665EFF3C3CFF3C03E88AC80387A9F2D7D5
90576F1752A7FD80B4C742D981CA3CDB88D251635F2B3E4F9DF4491A9A14B53A427E7283E4F7AA4CC9C11BE20B4343F3CEBEB2DE71D6F17CD579A7AF
53777E991C7ED462CBA04119A01A00F6BA7D3337B5E937E599FDA079CBA60F5C70F5FF63EF2DA3A34A9BAEE1B58213210E04274070770FEE1664700B
3E380CEEEE90600912820492900489901037E2EEDAEEDDA7FBB89FEF9C8699E79EF7BD797EBCFFBE95DE9315EB9E197AD1FBAA5D7555ED6AD32F784D
CF30C9E4F52A9FB1F181EBA4A67BB3734D27AECB546FFCB22551E7BF288D2C569B2CAA48033439E518AD51E310222CF823300AD1EB8D4694363D18D9
75DCD95A8643718C62399A414DFC33040370F30E409C440C184750CCB996CD9B355B031058E1834C0AE6AA82E977EE8B03AFDEAB9C7506B0F0DF020B
FE83F93F3FA7BF7B34515FE953F91D9D35113462EFABB02095F15D0E76718CFE82C3F402CDDDA13DDE28EE5DD9382B5DEA7D03F65924FEEEE1F7D1A9
67DEC58B4B5CB7CC1DFA2569EC252F97976B97BFB01DD6B357112843948BAFDC6A6DD3A2EFC119B6EE0B3DAC27DD76B37D74A5F315F191D16961C3FC
8206BF49397BB128E7AF20B8D2E796AA32F0580E6A3270E21F958662ADBE5E469A644616D1202045632489EBF53A9C25BF4CB71EF55A4E3028C6207A
9CE569CF27FD427F307F0210200CC1261D4C1004F7A4650B97C11761DC501EAD262146F4525F3C7CF2D36BDBABA39CFD01C1FEDF72045860C13FEA5F
E07FCADBE4A1A5860BE8676DC9F0EF38FEF5039E928095A441DA39576A46F438A46AF86BC0DCEF89CB4FADD96128BD94AFDD79DC706B79C8429B35B9
F7A73ACC5C31EC44E3D4DE7B7BAEBE3CF944DF21DDE6898D4A54BE2824DCA19DD3DE532E0E8B0E8F710F8EB6997D75E0B59A239322E3C66C2ABCEDB9
E65A51E39B3BE5AA843D2F8A74E1C17100401B4B0A8C2022AA6C8050506EA08C22182084861F69BE98A0C12F4BDB763A5943D206906628D84432B4D0
F2430A36803485690C386684509622B9ACB68E472B115224D74885BD810DE14A72D7DCB7C95E1F91E3EE49A8799D8085FF1658F04FF8E758E6FBF9F2
41C1CCD16CCD4B72EB5A14CAB98EEA5FEBF42905D837F7A420E7C161DA4C6FB7A5E9E7963D9AF15417BD4793373F4A3B73C5B3F6D62F3F8E6BDE65E3
C8A13F8E58FFB9AEE3C6ED6BBB0E6E3DBB0ED1196B26BF8968D17CFCB51E1D272C9C62B7A5685697ADED0F4977B87CAC9C3D21C977F1B493AAA28D57
0D8633DE7EA6942BAF8A1412948A7CABE74C5589E95A525F22D1974955852526942153C21A49FDEB196DADD71621940950A871DA041114460B4EA018
25F8011346154C5308456228A79DB7B286E3C89A480586C1045A9706300F1685556F3BAC2B9B313C0212F67FB1961D201658F00FFD79FE276D93AF7F
C205C43057CB2BBC24A0EE9108CB2CC0A405066C9F67CD999E5E25C0EBF9D38EA6FD75FAC2A010D9A1A758D861F987A9317B3B7B7D3CDCBBA7F7A651
3EF15D665C193264D7B3A11EAE4BD4A8C954F947CC3DEB561D7A75593A7FC9F8B9DFDEBBCC19323D6D4B8FDB5FD70C8DC85DE9F9BAF4D1CC8BD2F4ED
DE1115F9876E9C8F855855F0EB4602CA09CD36128A84625CA331A6172004A50808D761E173AD6D3C3F410C0823A00A664823CED1044B41905900084D
C22885A0388D7098362192E048A2325E4B0228662AAD00E9C88D89FA875B65C8E9555E8116FE5B60C1FFC57FEAFB82A26B27A9D73E54C40BFCC20B25
1C1B8129BEE8A1FA327DE3B4ABE0961E87D55517C6777BF17D98CFDE613FBE0CFB8CED3E6ADCB731A0AFC3F983BB3C869C9EE659B8BCDD991D8E4342
36D876982735289194D951D3AD5A5AB59B326ED0EEAE071207D98F9BFF61BA73486C6FB72BB99BA63FFA716272185CB0604BA9E2DEA9EF4F8E88595D
F0BB06AE363134CF402962B231A35A5AD8405268CAB53422D7BB83D3E6782D8D803806D1384212004C0B93C0808120CD2D81308199CDC20943BD4CCB
7184BAB21CC6092DA2A9A884C8A41DB970F8F868AA74E25F136E1A700BFF2DB0E0DFF29F2392D63C483D86455D46B49FD0D248A3B2F6038C65D49160
693D933221BBC673D80BC38F3D13BD2ACE1F16AD3CDA70636EF4A7EEB1155E2FBDED4EBC3EBABDFBB61DC302229C96DEEFEFB8ED944BDFA3B546004B
991733C5AA85C3CE01764BC7F7BE7AA6DDF0A5F797BB86D54F6B35EFD38EF137B2774D492095C776641B92EE27468494EA8B83E250BDA2B2BE868FFB
F51528216A90A1145DFAE1561D72CDC37A5E12CDA9741841D0248DAA8D84E006C043B8F9232802153682F3E780566A2CAFA6594E9BA52008C6A482F5
955200CB589645E68F3A08EB43CEDF1DE707200C63E1BF0516FCABFC87BF9CEA29DB5657B6B6117B938C86541872124548CD5B8854C8B4CCB565F0B7
6E533215372774BC283D929BDADF57E5D9A5E0C2901FDBBDBE8EECFA76F8917E3D8E4C1B9F36D3F986B7ADDBA8F6E39220B512CBF22E9D61D5FAC49F
565D97B79DB8C9C1FDFCE6F1AD9FA9FF6839FAD456A73D1F8E8C8B44D14FF3DE29DF9E090C790D19139E14A8A132350D5469384AAED0D64B4D38DD18
1314F1B9649D4DCF2B1A0E4C15F10CC7190A468CB030E12374FA0A3BBF84417F084250026BC893D735EA49529AA9E0381C34A05A890CC62A16BF806B
67778811A95E3F5D3ADE5F035BF86F8105FFE63F439E70F7D41EFD0EAEC9C08B3FE21F8F1A9591D518F6249E4135B508B0F5367CA2E7EABAB2D91EBD
CBAE6DD0AE5E071C6CB757E535C7ABFF9B0BEDA60C1AB6A3E7E6CD9DCFDFEFB0FE710F1B87767DB221951E7CB5B572B4558F93AECD162C1E7A6449FB
0D1113DB2C551FB01A786F53976DE9DFE79F02C06F332EE9B28FBC2A30D6E7872767C3C50A235C560AB2744D362C571B48C5F7AF25927BA346F4BB52
4C317046061FEE85FE3E3ECBA785DC1F47315CC8FD095E1580268AA028406E30F08FE9E51047D21889212A1842C427FD819A350EC7C43140DAD7C80D
3715386DBE00B0F0DF020BB85FB7FF0C72A2D3A8CC0F81DA8B3B4CD897C646EF583C39468FD7DC5733BA7A095DE2952F1ADBED82FAE468D78B7E2E5F
120FA61E5DD4332B73EF5FDD463C1864BDDB7DFFB8AE1BBB0F0A1AD6F9C19A36ED6D3D324AB341C3DE55B1F656AE4E561ED3EDE61EB2EB7564A9ED98
EF7FB61E7D7D43DB31EF0A1F0D8F86C286AC2D0B5AFBDE8889DF7DCB576616994865459E86E30A127515B50490F42ABD3E7B59B7F6CB13481AC2500D
C1D2384DC81B10C10580A2F8849F0FFC042128001C47099A424C004E51244BA21C07932C0BAAB4085E7733A056163E70634D4C2AF36DEF8DB9570CA4
25FFB7C082FF0CFF1C4B697D7A8E49CF0FAFFFB82D092F0987B3CE351A8BAB013A3501C44C0A031D76D6F0BCF390E8C4991316554F9EA3BDB96BC7BB
25DEC9D1A2530E7F5DEEB8F3C4D28D2E8B27F7B8B6D371975FC776F65E707D098E9F3DFDA56D1B1B2BBB0D635D1E5FEC73E48C83E7E7EB6D477D3FDB
7EE0E9DC5CCFADC957FA2D4F4F9CEF83E161A7CAF0BCCF2510101B5F62E28CF1A1D21A15260E0C4ACBDAD1CF6EC66B35814242830F4B9BB77A98844E
7F4698F4A17E76FE90248DF3C93F896A85D65F82A6CB75E6828649CD2705A2ABAF818A86979BD36BC38DF85FDB4ECE7887503F1B802CFCB7C0825F2C
604845A4F3E0B7970BD597E24E1B91E8322C248030E6E79493297924246D348027EFA816388E8A3BB2B273C8F3AED541BDE7271FEAF1E6528E6E4BBF
90F9CE01DB96B90C59E638FBAE6BDF37B35BB5F1C6F40A2D71F4F2072BAB165D877A365B1CD47FDEAD1EA3121E77E8ED1FDAB3C79DF4922D33A27C47
CD8A0F9E7AC404DEDD948B16879613E8D33F8B0D0455FCA6586F2053EE8748AE0FB41E7BAF91244C084E9A4C284D337CD02718F3CA2F823F01CC6B7F
493EF0F3498049AE8761FE518CA0B43F408EA65993486A24AA2FFA5797D584AEF22B0E29A622CEA80F2E2D476833F72DFCB7C0827FE23F29CF9B34CC
774E19F22EF04522A12E24E1A82C5813FF1E3465E84983548D56ADFF9C3BD675F96DDFD90BF2961DCF9D36E4C9E1B6C9AF66CAAA973EFBEAB2E4DD99
110337B90F7B33DAF5F2E9366D6E70F526843CF8E653F316CE3397751FF87249BB655336A71F19D4E741CA6C97B9C962DFC19FB2D6FF1153BC61A7AA
F1C49F557466AC9C525F3B5808419A8258B91180A21E95557ADBB9FE554B514608C701858120F552C1DE9716F67F13C2E23F8CFC39F7CB7FD0B4A151
8F11FCC9C0FFDA84300CC367001A03AEBAEFA3A9C9155DBB2C4ECCC4D467457503EF987086FD39F064F9CBB7C0C2FF9F9F59A2BE60954DD4C98F5AE3
C5806098F852CA18D3354469C437489E6022649532F0F1566DE2A04E8BAF5E773BFE6662C44BB7EEFE23E7D7763F00052F8C9E627FFFFA62F729535C
8F2FB01EF7B843CB2B5C9E8100FF087B69D5BCD368671BEF4D36BDBA0F883CDA65DAF0772FDC3AADBC7C78DABBEFF396A795AF5F2FCBDDB0BB9E2EFA
0852A6C7872A998A82C46A25D750E9FF55953AB9D9DC2C9A31E90018D22B8C380314C90881E182C32F4F7502C530AD06A3088A6668824205C7109424
1114E1D30482E2508D8191F986C8CAB2EBDF9F4D4FFEA0C1DE7E267CFA4412246389FF1658F02FFE33785DF9D9D6EF924E8A90F29B1F8B09853FCACA
8A51A02EB4924ECCA10171954973F939F8AA5FAFB94F96760BF63F94B5BFED932C87B80BAE792A8F7D773A6E8DD93F7AD09C0EE3FE70EC7C755CAB77
5C3EC8826B8A6E5AD92CD8ECECB9ADB9E3AED1134675DFBEF066C99EE17FAEB7ED7E26CB736D89E2E0FC8602AFFD9554FC96183C23F00388E6BF2EAE
A0A8F888F490CA838EAD2F42ACD1041A314C032034212E35F0BC36E7FCC255BF60F94709A3FE1461D011C2CD80A00A480A43098AE4C33F891949C0FF
2B22CF10D73C4C5486D46385F78CCA79EB25A4B0FD8BFD1FE56381054D9DFF7C40341696DFB1DB907531B301FBB2A9D4C8C47FE0E87AB1415D1D2531
061653A6AA7A20E348BCEA68570FAFCB4E338AB6DE4BEBE55EBFDABD6AF802C9E66EDF46B846DDF0741AEFD8618B63CB45335BDE6712259C6A69ED76
AB96F33D9D168FB159F4475F8F9187860C90C6751A7A7CB8F3859383B7C84A560C884B5F7A4C667CBB2DBAE1E9E13C1CFFF4A2182280674192D767A7
B6EAFE99C2353A8D1C35EA4092000BCB4C2C23ACFDE6E33FAFF249F3CD1F239C06FA423EA5E7BF27483E21C02094C0F94381040042713998D1FFD097
BDFDA6FC9004ABAEE5E2B7DD7DF40CFDABF467E1BF0516FAFF92FF9C26BDF2DD80C9DF63321210C391AB34C1F8657394425C01D696D155AF9598A6B1
541BFDE64DF2D2EEFD4EFC61BF33E3CFB3EB9B6D29750B8EEF1E943FDEEFB8D309DF25EE7DFAB8AC19D37ACAE8560FB97C3D275FAFF16A6EDDA1B5FB
4487A52BBB76DFF97E57A790AA493D36CD197AEC49EFA375F1335C6F15EF7DA8C0EE8F4A670A7C0A19FDD5C372942C38F89A7939D0CD697D35A5AB14
4B155A1380E354434ABD60EDCF502623CF71A1EA2FC47F94E0F50028D7A214CE908811E63F9B6082404896C0404876FD1186AA4D8D11696071901E79
7117AD183E3FD94433BFEEFE2CFCB7A0C9D3FF6FFE33E5E1255F27F67B764314578C89D6E5734CE53D150B29732A09711D9DF40440D469E2FA469F3B
F1935DC6EDEC6CFBF48BF31F5D9A455CEBD7B06248C99971B73A4C7831665CC7FEADFBEFEDD8AB63F3105D16C9E5EED48C6AE5D8B2A3D7D03E5B9C6C
D7171D6E7B5A3FD761D372D7CDBB86F918DF761E1CE0B73C0842433CBF810D6F7F182AAE5EAD4664CF3C8389D743ECC77FA13859891C36D0B809A74D
65596A4A50F8240921660560AEFFF1011F2270A58665089C4601238EE2186880301C238C7A93F4FA23829592D2944A3EF4976399DBCB91031E5B9361
D2C27F0B2CF825FBFFEEFEA7626FD4C6EEF0B87633060804C8ACA3245B119957C7355646CB916A2512180421E2B246ACEED4AB77FD3BCF9ED2AEE7B7
A9CE4F5B0F6A9CBAF344FB0DF707F86E1AFC7CFED2196E8E8EE7BCBA3876AEA51B192EF811D0CFA14BBB71933CD7CCB0ED76E1699F3F8AB7B41B7B65
D9C43F76BC310675770F889B781AC6E3574413753EA1BA3ABFF720F861CDBE4792DBEDDB6D9073B4BA4A4ED00CC34775556523C4E0FC0F8439EA9BB7
019B7DBE498017FB20CA9F0CB8701AD0084E213A18A610C8086BA23E601CFF5FA900093C3014531C7F8B7C5C367D800F4AB116FE5B60C14FFAFF6DFE
835D5B551B7AA3DBCA988392BA67063A2A9EFB96A7CD303594557F12417926C0271443AA1B4AA82CEFF45B8EDDC778B45CF6C9FEF9D166AF9FB69DBB
C5F5AF8D234ED86D3D357AAB9D43B3E9A7EDADBA2B8C1920B1E343A553ABF6DDDD7A6EF46CE9367766B7D189C76C879F9BB9BC7DBF7A3CADE7D0AFAF
FAAF5791712BD2A0863B2F74C0B3B70074FD5898FAD90CDB1E4F698E55E7C829146B4CAEC360B996A5CD23FE7CE64F0B5F8425A07CBE8F03F93A41F9
53A4E003C2F0E90005C3261CC3F43A882909863916951930864E7FAA87DE5C92976EFF34A2CF533E39602CBD3F1658F01F4DB02CA75F3FAF2CE6D2FC
4929F7E2A8C867249DA22D8AE65465A26FFAE2C42A6D25A67B5A8419A1862232F876E562E7FE433BD95E3938B97CF284C8717D23BC573CE87E77E99C
B0D1F37ADADBBADC1A6835B2F65B3C23DB9E14EDEAD2B163F7CD1B5D9A4DDDD5C965F6F2CEEE87EF2F5ED5EB7A55D6A431494903C657E31F97C561B1
CB7C75EAA0C732F8D67935E6D7B5F5B874216EA7576184282F34D160B6F7119AFC5010FD99F7933F7BFE20D090A9670881FD422F20451238823204A6
01201C95A66A380E9381FC6BAB7CAF278BC3AB4CDB2FE9060C78895096DE5F0B2CF845FF9F9DB034D7B06C5048E1C1B3FDAF7DF253D377E2395303F2
5ACAD54995655CE967852C1F170528B10632B150743A38C9B95DE76E2D7A5CF3F8FCB0F5A1B06E875E3B1EF01E7C7DEA2DF7010B9BD95AED9C6F3541
9C2CE58C29558FDADA3BD98C39D2AF85F39C210E3B665BBBEFB83EAB6BC7B719E3FA8E4C0DEB3B220B7C322D0B0C5FF75AAF7EFC42AE3D7556831DB6
6FE955CBFF71740965B8B44EAE07F8488FD21421E87B18C229A1E38720710A25298D0CC541E142808FFDC23C008191089F0050BA1A0D464A92A41C6E
A807198E2B0F072979BA04C974CB94761B16835296E13F0B2CF89BFFE66F18AEFCCFA14724071F7B1E288FF8412AAFC83880CA8DE73005A7A389B26F
F28A52223596D428917CB0FC7EF9E9564E9D3ABBAE9A3C3678F2D8CC53BDFD27ECCC739B3968F7929E7B3BDA5AAD3961B51C5263382D2A09EED4D1B1
D3A5E52DDB4F5AEE36C97FECA0B17BDF8EEFFAB076699B21D16F5D5C3ECAAF7927210F6627A2F157431BF581C71B64AB5BB7388A30245E99528763F5
2A9816EEFA083EB4E3244608BADF5C0314B27F0C5519854704714031662700FE5C6030D4A0D5612671B684E1147900FF32CB228C1C5C0340C896E9BA
C843DB4A30D2C27F0B2CF88FE8CFC77FF6C7FEBEC3EBCF1ED83539373DD180E7BE053923FEA584D3AA382DCD448762892558600E5386115560E209F9
ECE6ED5DED777EEE7EC0C76D75A9E7E2470BAF0D1DE531CAD775EA4C6BABFE87ADB621901EC19EBFFEE0E4D479CFF14EADBB2F74B4DE37B6DF19C771
EB1C6E68D65A7B26FBB7EB11F066D94DAD7EDFAE3CF5D76DFE125DE4F66C70B395DD539A450B3EA535D0EA063D0EA104096324A90308A1D7971616FD
09270149006A3D2F0960D27C3A08EC17B6FEA1040E436A0437967E17D38CB154CCBFCCDC37268EAD6B80C8CCAE5F1E05264F2BC4C9FFA97A58DE0516
3419BEFF2FE19FA5E96F5BA6B8E5072EBE3EE676E3670CABCF2DD06A2849B88190E11CCCA8DEC41BA2D5B5974D508A29FBA3F6F2DD861E361D5C2654
5E1C7F7B4DFB8BC183DECDB2197074BEE3ED298E2BDC9BB98E6BF90CD3C01CB2C7C7DFCE71F9F36E6DDA4F9863BDF9AAED8E430E4B9CA7482E387A7C
F9D47D7068CDC80540CD86EDA5E0319F84C4F2C4A92B55EF9A77896259C5C7B755B841566FC08C263ED8E304AD8B9308D1FDE73F3CED715251041066
DF2F827F803683627018C6F82300533596D6C324D928E75F9BE29394E38C7218374DDF55BBCC746A9088B0F0DF822619ED7F53FD13BE5044C4A691B6
4FF3965F9BB6B2B42C854B6D44658926263502416B1946C9CADF56EB231BDFBF67625F4B1E64412BC3429A3B38B63E57D377ECE3B15D5216AD386BEB
74E35C8B51DBAC068C75B1736A1E8ED5839C6CB3FFA1B6B65346B7751C39D26EE8F31EF6E7DA775BDED2FB5D971EA76E0C1AFC5DF4C7BC3AD19A9DAA
94ADA13547820ACF8E5D70A8F3DC0A062F89CBC5114371B90955CA3102E6353F55170FD2248A0BC11E27291443703D4053FCB78400735B109F020843
010C4ED388580211069508221816CC13530CA43690CCF9C1E2C3FE60B739AA5FF53FD6427F0B9A6EF4FF1FFEFFD4FF44E8AD036D0EA0C7CF6D5AF75A
9E8097A47374661C83F9A4F01900AB413871BA4C9C1F1705496EA964D148EAF8DC45CD1CDA0D8C3DEE7EE25AC7AD7F7508B8603FEACA70DBA30E1D3C
FB75756C974C6A11A6E170F6F976EDECDBD90D5AD4C3F6C2EE666B33BA8FD83672E528B7D5FE437B47988E0D2FD1EED9A72EDFFC4EF2E47DD1B55B5E
0ECD8F02B43A35870FF285B98D046D52817CB8672810536A18E26FE0F9091A1811AE0004EF6FA11A204C01D3346904705EC6300C019B7010D389F40C
8D680123C7A1008E7171EE51191BE074DBB362CC72FD6741D38BFE24F3EF33E07FC2BFF02D853DBEDBD061381CB263CFDC9BD925E57820C441EFAAB8
C46DF59C4AC38930A63CCDA0B8FE504D247D41BEC5C1C7B714756969DF7A6174CF6EF757B4BBDCAF5FDC78AB05471D1778749AE5D8C5A6650C1F7FB9
C2CD099BDB3A3BD80C993AA0C5F87B6D9A073FB7723F32ABEFACDD9B97B6DB0D3D6E1F8C9C9C5A2BD91CAEDE11A6DB357E47FB962750A2F67B0D45D7
24E7E9085A2D36523841293FE7121882091BBD5102C588C6CF35188198277D48F36D20F1F32BAC858542200BAB211CC320234A52805C4AE01C2E33E0
5C79FF8BB23F42D16D6D2FD7FED3FE6B81054D86FF0C4F02F6B7F29FE7FF834BD466977AE5CEB3B3B6658962B8BC54962CFF0260672E315C3D8AAA19
53558634E55009AE09D389E20C8DF373DEB47468E7FCF8A8EB369F9EC72EB49874A783EBA1091D4675F2EAEF6AEF96431A8D5CCE99E29976CE36230E
0DEED0F5DAB416BDAE0F6D31FDC0C0F1092F17F51F5790DF711DFC7E5A3674645B8DDFABAA53A3A7766B7703C6D3BF2A387DD8FB3A98220D32986109
44FB351A14ACFD710AC3110CC5712D6066BD701348FDEC06103C40F893C25C0724E58D3045F199004850B0488D728C2EBE9EE4946337AB6F9D32C9FA
2DBB2C62180BFB2D6862EC67F579C07FE3FF3FF93F72651B196E1D489FDFB07E7A00949BC5046A19E8473E53BA268E03CB498D02D37CF9047D7D04A0
0111505AAAFEE932DD96D6CEB6638347BA1E5EE9E4DFDB6AD7C5E6D317751AED3074AA75F345382482E9D727337BD9D80CB932AA9BC392532DDB8D1F
D0A2D7D4AE76F7D2A68DEDE2D73877766DD6D020F4FAA2E23BFBAA8F0C5AD1B65318067C083732E0A3E74A1226945526FECC521A0061E617C1111887
7102433110C7714C68F531C77FE14368002070547006A631792D28EC0A35CA5114850DFC6FC08C749453CC5F931D7220097F312A7CAD5858FE658105
4DA7EEC7032F4A23D8FF1EFFCD09000DDD59022A3DF65049ABF67A1D95EB3F0065912CA5CF95D21F0E6838B58A1499743561F5486C3E9E17014B128B
951BEECA06393AB93E38D371DB778FF57FB5717BB86EF8E8E9A36C87CD736AB558A2A823D9B797B2C7D8D91F5A653F76D8F9EE76133C7B4D393CDD65
75FAAAF653D78AFD266768373DD4DF9D940E9E7817F8C73CE7C1E93872F6094A41AF3E182906965400145DF65D4990ACF9925F98EC211148D8010E0B
D3BF94E0FA21901FE3753E411BA4269EEB04020118CDE703A0444D9290E0040A96D7A35CF50A5F4D98EF7731B52A20EE8CCE22FE2DF87FA997FDFFF3
55FCECEE634D09D5FFC74BFA77FCA7F537C7D44073074980832BFE5AE043177E67624A3993BE426FBAF784633594A1A21ECE7D011A3F95C0A1B970B6
5F69C9E29C483B17BB61A1137A85CEEEF07D94D5C22FB39A4F9FDE62E8781BEB2E61B856C1DCD95636B8758731B663BBAE5FDBD263A6ADDD93C3CD06
876F761AE1E9736FD247E0AE9FD66764584DC099DBF316D9D8E752862BFE088A47879858C4D4D8003086AC0C9130D32F047921DFC755791212C724F5
8870EBF72BFA934A0D829A94757282C470083491048111A8548E63287F585060691DC9E48E7FC949FE4AFA84E72FD2DC780DB18CE5DD6D415339C5CC
1C6758059F58FF37FEFFF32483BF67BEE4FE8804387AD18579939574088A84AB4D6C7D1D5E7E27936890C28D9526DD8D737AED87EAB4C70675D06730
F64FCD119BF64E978EF758B0DB23E2764BA7D3EB5A2ED8E43263B9ADFDB23286D23157CF9A3CDA3A58779D39667387BEB327B6DBFF639CF5D123ADDD
F7EDBE34C6CF18B5F5F5CBE567F345A7564EEED869C07744F13218C535A15F8D1869142B295A1F59C13066633F5AA8FD813A65450D429258753D6EEE
FF374B7F8A046012A9AAD3A27C160099F824002709C4A847480C254802AE6B80A8CA4187392E31AC22873D18081D8F272CE1DF82A6C57F3EE3CD8FC3
D9DFA4FFC2931856F3746C4A7CE5FECB807AFBD23D939F93BC0028FC069B88DC3ACDB77020D35F04E79621DFFA45E3E9B79441AF4CA50FB388CBF7B0
0976F6FDCEB7EF746ECC82B4490E3397B8F55C3E7DB497BDF51D4A0CE9913D41856D5B59BBB8F55AD4C375F58056131217341BBDD5B6EBA8EB4B5DF6
EBD2072E4A7DB6364074B6FD3C37879362FAE1F3780C28F38F81685C53AFC3597D4839C3D0BCBA176EFE680CAAAB03709A21701C407EF97C09957F42
A805C82A74248191C278002DB4046120466010869388B4C64848E6CCD333DF8F664713F9FB8CC90B33490BFF2D6852D9BFB0DA272A83637E93FEFFAC
0F2AAE0D48CE2FCC3E27C35F4D3E7FE0004EBED0C1F155524E55A3A9F921ABBB1DA6ADAE004CFBCF13A87F8EF659BE31EE2BA4DD22CFEEECD87D655F
F703332654EFEFB0EACF21F66BF6751CE6E6104049080C7D501BD0C2DAB687C77CEF9ED30FF6EDBC617633D7B36B3B6D3D71AAFFDA868A2943A213C6
EDD7240F18DCC1E9194D9EB92C5100FE57D33014958B40166F08CE61589EF798E0EC8FE2A232314209BD3D82DFEFDFB53F8A10CCC009632D40095DC0
A4E005641E03109602F1FF16AA561A31D5762F099D34D6CF544CDE8CE77C26D4D016F96F41D32AFE31AC3A2C9F637E13FE859F18B66C7DAFBBFA7780
6F8AA161FD8AB587E454F037AE2CBE51C2D6578B720B907C9F78715E0D91B3AD1E493A060281999A0F998CEF01CECFD666ECB4D61E5E1E376ED8B53F
D4AFD5626FEB3EBDDA5EC2D0064CB929774F335BDBB613B78D18B969B0FDB4112D5C362C6CB5F2D5669749A5659E237E7C1BB94CA6ECD3A35B9BA714
76FC088C549EBF510522E26A118EE0A509A50CCD308006C36102979514483196970206D5CF4DDF8459FA53244AE8E4C25020810A7DC24256400B43C2
FCB9C09F06A04E0E62E089D5F5047CE53052179B7B428EED1E2765184BDF9F054D8BFF5CED1BF1BFCA5EFFA2BF59FF171C58B9018BABCAF8A0A70366
2E5F9B48D5FBB350A234C350FD5DAEAE31A05F82EB1B2B50FC7302021C0E470BEFD556C6E1DAF1DFB9F56D3CF7B66BB67AD2EC5B439D270C779E7FCE
A57F3FA7B30A1D88A60C7E3FA785A383D3D64D2EE30FBB0D9CD3A3F7982DDDFB279D6A3B20AE6E71D784F2516361F88FCEDDED2F13D0C123A05212F0
4286D7E5E76B01AA36A70C61599202CBA53889E1BAFA7A90C1689234E5D40A53FE94E0F9CF03D5117846094EFCED094011C2735406824F1C589A3201
100E3FDED680A30D415AEA656DD02B9658BA1AB094FF2D685AF46769AEF0BE82A2CD1BAFD8FFF9F53FD61FC20541C6D983FD8AD49F1B9E48E8E4490F
969C50508FCAB91FD50DD1555F7365B5D5DADC488931AF9150FA899080D131D8E3BFF44991E4AD01066967EB0DEB9DFB8E1FE4BBB895C7D2B6EDEF76
E9ECDC27A2AA1002EB0E270C6B66D362D206373BF7DEED27DADACCEC61D3EAC1B70E8E2F1A573A05E74F9854A77BECD6ABED611CDA771850357CF804
A932A20AE484B4385F84730443E3EA721D05202800B1FC4F2046012248E8F91552005EE963959568430144983B7FCCA340384EE0A0096784017F9244
710289F8B31CC7CA5F88C89028E26119D2E87119B3D0DF8226C67F868D7FA137B7BDB17F5704FF1DFE79FEA7BD8F1CF09E8DA98C2D822B76BE3875BA
8A4BF8CAC98BC8941C494EAECAA0AC7C1E0F578A60F2470A2AF63EA6549D4B203E94C2A35771CF5AF4D830C575E1CC9BA7BB8DDDD4C3FEF1D04E0EA3
95DA460492F9554C6B65D3E1AA676BE79EF63327BAAC3FE8663D3D665CBB9BB2BF5A5DD02EF4C8FF3A646E37C77384F4F8559338F64BAC11897C530B
EBD33F666B599266200D8E0130A637EA099221093D800A7701E6BA9FD0F14FAA45E5350D15108EE33F6B0124896118FFA810FBF957CCEB0312F9F1A0
0645650F3288B897B42A1EE33EF548662CDD3F163435FE3341410686E2D87FD3FF9FEC9FE3D9F2294ADDEF00531205064841DF554F8F8633F2FB263C
43AB0F6C3024C991F292A80B55F5457A25703A1DAB8F493665EC5354DD07331C5E183B5AB57176F45EE6F9704CBF294BDBEEDED4C6664841AA168652
BC53965959799E726AD77E78174FBB458FFB0E597C7272B335350F5B4C2B5AE19CF9C9A17BA74EA174C5DEAFA038E041265C1A16A1C6B4E1413F209C
2558759288C201AD064061866239484CF1929EFED5F14FD284B85109EB01984230F326104117A03029AC02A24896E51500077D8F4C90E278CDB9AF58
F62D804A2A83A9CD9DD32DDD7F1634B903007F1B6514EA5EFFC67F3C8F21031F354E5C4218424CC1F958CA9AE0C07B2AFC762D23AB23D21F281B6A30
797CF1D7A48658B0018E5A5F816B447AECA93F1B9CC09DEE01DCB6B2B26A3677F75C5FAF11EE67060FF5B76DB1479A08E8F1DC9BB23FADAC46F46EE7
E43A7E70FF11FBFAB7DD32AF6BE73E9FA36CBAC69EEFF04139AC9B47F7582E7666141C7FF55EAE4112FAB9962C0EF866E2E33C85CBD28B61D2D0A8C6
100A325184B20A2068C6CC7FB31F9031A58AE6137F465803840B2D420CAC30A298791680A004FE73F8BB6D1F401C179F7C01D73D9611DA6F7ABAB4AF
473C67E1BF054D8BFF0C67B8FEC540B3ECEFF8CF3F050BD85B7BA17B051195AB79AD915D0A3FB5ED3DF6F405869418F0B74FA0E802AC3E22BE0C4B4D
8152C07DDE0080A514CA8EE5C98E3518865D867A58B5B66AB5C7F3E85C4FBB551BAD6F8D6CF5BAB80651625F7C740BAD5AB6ECDCBF6DBB81BDDC0F8C
6DBB6196CD92CE4B937A39F9DD6F75D2B074FC389B6B5CD8A82820FFC6952A436CB2445C16F7558253108D35AA24329A46F42688A6289D9E90C55632
34C54776E1FE8FE0957D4D864198F931EF0110F6005114AA33F1E93F4D5318CE0857FCC0A7DDE71B290278E08F28031B102CF2AD91BCECECF1CD12FF
2D6862FCA739E58E67269AFB4FBFAFFF23FED3C4973522CDC0CF64591C1E55064727DE9A7CD69475102454625CFF34A1F0AD082BB91B42E83FEA2B0B
8AE7FBE274941F93FFD8E81F40BFE925CF6CD1B28DCDF90533BD278CEB7EA8F798C9D6276EAB71103D7EB8AA7733DB16137BB71D30D86DC18591938F
F4F43A3AEDFB56DB9361DD877E58EC3EDAF514F77546229972EC4A4D45C20F2D5B7F398AA05182D495D402422D5FE8F8D71565A90906121B2182E5E3
3F9FE7D3A8AAD180910044D1C4CFB2BFD91150680DC049B3132843739C3820E84A1143115951A0F24D0D0A351E2F47C5C37AED945BCAFF1634B10380
E6CA365C3152BF0BFFC28D0085C64CCA20FF3885C8F3AAF263E1A4870D2B66A5A95745A2FA4A1D2E7F5A167CB4124E38F099AEFA0A0525BE191ACAC8
1F57114F3E16EE4E910FDF80DD6DE9DA62C1AACEE706CD683B777F974EADCF7E5103007EFC5CBE63332BDB01AD5A8F72723BD2DFF5B87BBBD9A3EF46
3A6D7D37B8E7CB8976DDEC6F72716B0AB18CBFFE4CAD7E9707D29A97B11843629429A68C0228CA68C420521DF7320B115A003090E015004DD0045C95
5762E47F169680903F3B01090C23CCE70545F1F93FC1B06C6958C6F7728EA3B36310EC5D2C0482D79F28B5513DD686A096F06F4153E33F5B76F4B49E
F9F7C0FF7FCA7F3EFEA309732248FF1DA8BCB2313993AAB806868DF625EF1C475063198AD71729766D2D377EF0AEC0128BEB7C128F78E67335D954FD
7BFDAB97E053D71C553747DB6E5BFA6F5E33C57D70E0785B9B9726A911216E8716B9B4B61AEAD2AC773FD775B75C36DC75E938EA6AF08CF9FE075DD6
1EEAEDDAF50A193E21970D5DBF3F2BC32FD5002111090C9FB753A2E03C94865195166231834C01F0593D2E14F369C1E49742C5A515200269F15FF700
046EB6FDC7298CC070F326009AC2158551D5853CFD89FC5808FAF61C44E1F833AA5220D6EB762D6D09FF1634B1F21F45959E3AAEFD0FD3EBFF4BFEB3
B42972F9612C75A65619A74F89561177BFE0AB47882AD664E921551D6AACC0CBA73DA7A0538795751F0D050F5F792D43B9AFC5DAE0B0ECBF12CA87CD
54AFB66AD16AB967CFFD633D9C0FF6B36F7351998323D8B1B08496560E8BAD9A4FB66B7FDFA3EBCB61B61DDD5F9EEE7BFD72B7BE23BB39F47F87878C
FD4A7DDCF922AFE25AA09E003E665004C7200D7E415ACC84AA1A101CCA4BD4900886F0919F342B7F9A24B5D50A80D48955FC2F0427704A28009A633F
8EC20841F01902837E0B2E6F48CEE45F9728CDC0D4FA94E358AD6F75E347CD9AEEE74D346BE9FEB3A089D5FFF092078F74FF87E9DDBF58C0B0BA373E
3335F219B5687CBDF16B2D15F5888EEF7A89B87BD75082988CBA64884E5F5245956C4805527FE87342AFACBE42D68BB0A25899CF13C5238FF2546BD7
7663E7BBAC5B32B4FD90CE3D5DEFE9C538811E2DFC6A65356CB655FB05D68BD7B45978C0D661D2FA1BAB370478771AD1B7F5B82C32DD3398FBB83A4A
1F79214C8B4351B13CCD99CACACA9C4A89C4A0034056F2E1931484009C6418868FFD82BB0F45A0180B03304C3094D9F19B12E84FF28F32A4E00022D8
05625FEF26D665D7201CA59750B428BC8224E0F422223937C8BE7B006261BF054D4E00A039779E2AFF37D73B8695DFCD1CFC0D9E7F05CF4D21E37D21
D5DE06727BB7EAF25595994938829E4B6798A0ED287E6093421D530757DD39BB2DA0380A46833E64EE7D5F3C7B05F247CB8E4E63874E9831B673FBBE
6DEDA259A3016F585CE06B65D5A3BD55EFCEB3B7B7B0DBDBD37AD2AC3563473D5AEFD0CDC3CE3D0D4F99F5860A59182A8F3D1B8122EA48FE7FC3822F
33000526ADAE32E928243DB4D0C00ADBFC2856D8EF45E20443503843AB3418C59A2D3F7FFA8099F77F986B83B4F9FBB48715B5654A8E46A400C54943
CB6992C8AB80A54FB3A6D9AFA9B4C87F0B9A5CF8E7B01C9F27D2FFCDF4926194BEDAB50F88C77D24C602A87E7F16FD2C8029E9B90DBC734DFFE53BCE
A41DE5F3EAA32790B4B1D1C6BA8C2AA066EF8E98FABBA9F8E3152937CF8ADF0FAAA97375ED3E6694E79CD19EFD063AB68F50D59054F60AD949AB96B6
56CEC3C75C1DD77AE1E376FD7D678E397438B06797616EDD42B0C05E57F0D085018AD85B6F60794D71394852DA3B0F0DA41C502AD518A20F7A50CF41
7A95B0DF8B615821CF177A7E692C2705E4289EED8C79E587B006FC97FD3F23188213929CBCBAAA46986071BD8963AB3F94732CD9980DE34141E5F3E7
5F955BC2BF054D8FFFC4F7DBF7EBFF97F8CFF2FC3FAEF25962920EF98055D511913741A9772D7DCAE57BF5D97AF07D01CC5E0E671971EF5BE0AA9975
E2BCCA7275E1AAC378DEAE64F1A1BDE9EBCE14CE9D8C3FB5EE38C271CE9A2EF347D8DBD80657A65174D4C2CA9556CDECDB3A8D383DBB459F33D3DA78
FDE1BC27EFD4C26E035C6C03F150B703C658AF40CDFB532F104D8A1205F81420E28D1290E91552C0808B6EDD941B3183584608557D96273A6CD6F9CA
EFD988E00AC2F0A70265DE0A26D4FF05D76FB307381EFEA4D1502382480CD62334A70A2DE609AF0815130517AA3E2C1D7510E62CC3BF163435FEB364
CC8DFB55D4FF52F96219606D5AFE9C02D2D78F826B4034B084F43BC3340CDE264FFBC1C863259C72979CA54FF690658DF49397A0BA0A38EB9C967E76
06CF5BF5E9ED8EFA8F4EFED16D1D07B61FBA69F02C4F47EBB1552605C5E61CAD9A6665D3A17DB74DDB9D9DC7CDB79D79D87D61F6F9D57DE70DB23F82
D58EDD28FDEAE50B272F083135E6AB88D4024A9C5AA4C2D57AADDE80EB138352951860C25061AF274ED0AC30E84F52486182926679B20BCB7E19EAD7
2E508A1616FAF1CA9E3125BD17896BD4C2E00F80B29C36BD9E7F6D585E21A1BF9E58E339B247106591FF163441FE7F3D7F3681F85F6CAF19DAB4E0A5
F1D419A2FAA81111E7630D858D8A25619CEFC0ACDA37F5AC2852C67DDF65A2F21C0EC0FB3C2BCB8AB9BC8F4CF45D1839FA0C7DEE55B6FD5259A72E99
C3DAB8756B377F54F7794E760F884A0461D30EE5BA5BD93B594D3BE56CBFB083FDE02B3DBA273CEEB2BE674FBB03B074E1B4C2C8196700E9AAABBACA
7888CDBE51AD48ACA0341506400AC19AC0102355554B3224CF7B82C5208AE5137D4C07692B41861701421F30214C009B978199F7FFD2FCD9664A09C8
915696EB090282400867A07C3947735C43969E7AE3A7DFDAA5CFB43A73F39FA503C8822655FE6399A80BD79ED3FF1BFF19E3E2D7F2FC13C5D4BD1F8C
B1BA5205D580418B54F0FA7BA6EC2F38FD231C667D0E1BD0B33DF38B177ED48416549EFC2E5E768D2CF4CA01BD5FC56F2ED9E550F1D5BA9BBBFDC285
3D46B9DADD2A092739EEA17766AF6676F65D0F7AB61EBDBC55CFA5DE363B03BB1E5BDCD1BA53896EFBE0A4848507D586DDBB5559F92007DC4E927CAE
A1F0FC221881C806DF7482AC8C92312C4E9324CD928830CFA7FC9025F87A08BC37BB02923481C3C4CF1D6024C3704CC5B357C5DABC1A84403104C108
4C55A2E65F1887CB31A6F89932ADD3F8F1CF89BFE79F2D6F0D0B9A0EFF6936629FDF259CF96D02C0F2F17FCD83D2D2A44B78C929D80815D5E3F5D9DA
C37E5CD8947C383C96907D8A21C8C34721F98259B9976657675C50141D2BCF9BE247BCF38222A7166C3A9DE6B20E9E6C33C07AD4763B97EEAD36E424
63147D7C679C93554BABB5EBAD9C67BAF61CEDE4B1E276A793C59D6C5B9F0297F70CC99FB3BDAEEAE2D9ACE80C9203DE7F83F20A69A23403C008AAEE
553E8E6B724BF8D49E214088A2689CE598DC0B415A1C11867D69F3F83F4DE210246C0511A6FF31FEFC4A0E8806A10639FF0C88673F08EA450A96C639
53B59A943E4FC75777F35C5EC3FDEA7EA669CB09604193E13FC3C67847DF33FE5EF7F2FC878FFB4BBEC86F17734FD228934A8341395AD13DA369F5C6
523CA41887BE179286EDB7808AD55EFE6B5F901F33B8447F227EF637ECEA33E8E0C7843DDAF35D6AE2EC867419BADBB9A35BDB9B8489C3B9872FAB7A
346B36E8A843AB61C35C57F6F3B8FC685DCF146F6BEB5DC069870B951B7616BEDF71BFFA59808AD3BD7B4FA5FE60A9EA722D058933BE372215F50A10
A730CCA8D3C042DB9FA1E16360094F745818F5A329F3BE6F4AD801069AF7FED0B0B1E673603548403A9C4251948041B1CA881014C960B929A0362607
ADF13CB7E491F9F25F38018CD0CFBD671658D024F89FE2157AB3E677936FAC39FF3F764FFFA6B0E42B527FDA881BE4285559C8E5C6D237BDF62A34BE
D928995C41D61E4B36151F5C1EB8EEBBE15115F5F82D76A757097AA93E7995ECCF47B5BDA7D5CE741CE0B0A1978D5D2B1FA38423A96B61756E566DF7
0EB67299D9BAD7AC66475216F45935C9A5CD7C3AD46E4FCABC4DE937167D4FB8FE506D2AF0CD40CADE1B99C2469ECF37B7A59045EF1B5198860C6A48
A3379228CD961E7D55CB27FD7CB417143F1FFB85516021F5275192E01580BA41525A2C223018E2150286C390A44A8D21248EC03428AA69084BD272CF
CF1F9895CFFD54400CA5305AB6FF5AD084F47FC1B247EF937E6B7C2DEC0780AE5D962526437912F6F5530CD5EA71A8D0C8641862F61C5F54F7E9603C
09A43622499F4BEA250B6EBD7AA8F9FA8A543F7A923262BEBA2140B22739C15B7FB0CDB5B7EDFAD9AF99DCC6A1639E2C99E3E83B19D59DAC1C465BB5
18BEC0798CCBB8944DAE6BBD1CDA0ED056749A59BC6BDE9B9D6B83738FBCC38A6EF827D145FE1A38B99166E12B3BEBB0B408352FF9190CD6832085E1
28557B239CE4589CCFF35114160C3E081D429B6BFF662750223BAE11E453058C3F0C285EFA938846A426481845488CC039383BD1481BD7DCE87551FF
77F7335E0558F86F4113E27FCDAE47B5AF6186F9ADFE6721DFB9098A983CB44865BC9A876B9430DDF88383D28DCB5F6E5FAF0D3CF305318514473F94
BFCB8B5F1FB6FBBAE9C277A67CF9CD804E5BD1F0FB57BDB5F31FA63BF64919EFE03A647DF316AE8D196F390E389092D4D2CAB6794B97A99D9C7B5B9D
0E75EAB878585BB75270588FF063C39F2FD95EF066DF47A4F0AFBB1230C9AF5199984F73C6C00B2A3C2C9BA0311AD71B21828FF808AECF4DD55082AB
17C167F62694140A00468C22313EFC0B778244F6D30A9282693EF0D3182ED885E2084412308623463DA44345891282DD327E46973492F9657F24FD62
327BA1590E000B9A86FE6FB875A8C14FF79B02A07941001639EA09921347375651B53106930A63C862259791FD6552DEE9773F629E8518F263AACE86
8BEF549CD9B06E7A4ADA5385E6F1CABCFD1D1EA97DCF6EACBABFDAB4CEFAC54DE7DE5DB6B669D6519499C971C60325F16D5A5BB76833DCABD538AFDE
0F67775A3CC1D13680DBDFEE71E4C4C39BB7157DDBFB421DBEF3752556F6BE9ECC2CA73979700858112C323B7BEB4AA5248DA138541F5980B1182674
FB90040A0BD3FDE64E1F41FDF3EC873120B74170FC17DA83CDFE40C206205CD81688C9CA1508A62953A0F4C50E819D8F68D85FFC67524209DA720D60
41D3E1BFE4D19CC46325CC7F3F007E2E08291AB24E56F6B011CCABC772A5B44E82734025673C56BD7753DADD2F6F7383BE8B13F435E7EAE3DE7F1A3C
63DE2A4DD4D332E9FE7D21A30797541D1AF3BC76E2A7D066D362BA74E830D4DECA253B59CBB132EFA2D0162D9AD9594FEC677FC2CB739E6DDF293D5A
AFE6021C76252C9EBC6E4BF2AB63C99217AB9E4A0CC5CF7294E9751C991F92AA967F9452144E33404D839EE0A3385EE19B422004854108C3E70238CF
738611EE0009FE318C268A42EB115E0B20089FEE0BEEA0384C087B00715E37E814558D0002E79419683FBB4B774795533FDD8F59067A164651963E00
0B9A0AFF59461B7B20F96EE46F16DF9B9F4135AC3CA9907FCAA51BF310998823D5004E1BB46CC6C5926D875FA903024D3929A579A6F46F3F42CB374E
DC39344073E62D92BF26296AF8E6C2775B4FEAFD7648FBB9A5AC71ECDADCCECAA3BE06E338C95FB2E0D6ADDAF4E836A6D5B49D1DE776EA3CC3C37980
34A7FDC88C87C3E6F9FBBF7E95A3CBF8E3B404CBDDFF1193552058D2EB7C56918F082CA68C268DA0F729B8E6633440403C9F4D18C1AB028C4031935A
62A41861E49F20734FFA2884F15F028685C91F86524B30DA7C2308030A991AC76151B98189711CFC75FA378C637EBA1F73D267C58245A085FF163411
FE93A6A4C3111F0FE968F6BF09809FFCAFDBFFAC5653FF490EA71622A51A0E94C31AC24490A1C1A137F604365C4F043362AA2A32BF8BDF66440F9C38
665443C296025D80AFF1A64768C2C6AE19C8ECE09356BECFDA74B46A6BD5B7402723B9C27DCADB2D5AD80CEF68336873EB512BDB741D61DBFA8DDAC3
E676F4CCBEF7FCF63FCB30D66FF22E44D21EE701CA3C23181DA86868084E177AF9614D898C647096D0666457A306030C82284A100C4D23008C68551A
5428FEE1A821FDC1F37C58D004208E91A8D022A8E71FE30F0A0C0341030A41A8AE584A54F7B39DB33A1466995FE6674CF2C706F497199AE52D624113
3800482CEDC6859C2D45F47F7DCBFF6445E39D07FB01435E39A6488635152CA3C3A1469020812F2F221E0D887BB5556D8ACC004A0374E501B2131DE7
0F7C851F3A636A5CF446B1C1BB3666AC37FC645A96D3DA1C0F075B6BABC972DC4873E597357FB66CED3CDAD16DD7F47687E7761BD1A9E57C78639BB1
89BB079E4FF17E542E51ED5B94A98D8BA8E5AABED61B63DF95801F6FBC046886ACF854A2A5CC5E006105185A5F2A0720981226FC21A906262861C79F
60FA95FB232F36DBA0D718CD4DC01082E0821588D00C4CD12402930804A298460E1917B88E98F310FA457EFEF8D307C5690886B5F0DF82267100F06F
799AFAF6E7F2A89B7E20F35BFDCF2AAF84F72933CA456A223B8BAEADE230150BC42908585D72EDC7E545A907CF80354F8AE137EFE9A04078E9E0C9B3
6471633F29E7752F4AEDF959BBBC551830F085779B1F7F58B56F6B35BA11876836ED5CAD670BEBAE1DBAAD9CDB6EDED4EE43DD5B38E6BF69E3782962
F0B4A81DDBCBF4BA7D93D22B7D032063755A9E2426425EF93AA811E71852925286520CC6193363441859F4B50A34E2244973B4AE5A860A9EDFC2FE2F
060A3A9E0D60180DC955A460FDA1E79F64F60913DCC16902C25018064D2291063B6EDDCFFDB89E667EEA7DFE23F34220445AF86F419311001CC3E6AE
3C10967B41FADFDFF1E6BB70F854D28E9BB8B49444B4A97AB21CE06098290F855119967BBB725340FDD6343032DC28B922814FE74BC64F1AFF1A3AB4
5E74CE6585FECC7AE9BEE6E34D9716A5B4B8F9C2C6C5BED994023D410337F7E70F6BE9E061377D9B6BFF597B1677B26F71A9A247AB09097F8C7EEFBB
24472FDF3D2EAEEEE8610D9AADCDFF91562C36BC8B20F83F865EA5C5691260909AA42288917D8A91A1286A6EF8D36871E6E7F25F9CC6D49F2E54F24C
474802416892024DA0B0F98716B2078CA429A105102364293F9448B8DBD6D913CB7FDD7C0AFC87EF9E8DA129C692FF5BD044F82FD8FB48BC7CDEC8AE
2593FFBD03C09C00DCBC96BD4A0946D66344650A4D151B3918A1D3BEE0865A2CFC8CD72079C26E95E47A369EF28411ED97BD1B376F6A4DC1F43BB973
DCC2E45E5F33865ABD96778D1DD733BA4F4B87E66B2BEB491A787EBB68604B4797218B86B85BAF7BDCDDA6E5E6B2A92DFB9FDE3DF49368FB7BD8F0E7
A06F557F6DCF16896549A13F7251D9D35090FF53D4E6EA10464B63AACC1F751C5DF22E13225110C30896A837710C4351048A50983E3DA90EC4711043
718A61688C10B6840BCD40144110348B43008C50A8AE24B501AD58EAFF6EFCF79FA38FC207C3FDD871BF82A2FEF641B7BC432C68020500C6E8ED7323
3EF509F47B0700EEF33CFDF944B2EE0B45A3F1524E934D5310864756621546CDC579ED8F03E77DA1881B2632248F4BDC12B6F384E75F0ADF1552BF21
0BA5613ED9AF3A8D8597ED7965FD62752B57BBE72448A374C5D7740F6B67A74D133D3A0D7E35DECE6E63C31EAB9E0FF7747F4F1EDB8580673D3ECAAF
EF491595C803F7A56AE4CA97A774BCF6AFCA55B200861275296508C3CA9245348E9804574FB458CE713483194C2804974416C224694271C1F6577001
31EF04A648F354308563463D4461B8AEBAD2A0AF8D497AD8C317FDA7F18965B14B9189068A612DF1DF8226C37F8EC48E6C09D82B3FD5F85BF77B862B
1998FFF9194124FEA0C8BA0823939144E30642150B9A8A0C69098B3B14268F88D7FFF992D2BE92B2B7576EBC7CA2F37DF5E6533FA6F77D217E7FDDF7
41EB800F1E996EB36FB6B275FA4C1BF8989C7A3FC2ADB5759F7903FAB7D9BCA34DF3BE15CF5A39DEBAEA7C17BFEA5144FB7779547F6DFDA37CADF6C6
1C9F2AB1FEEDE297184BE77D836814A4C8B2AC2298A3CBF3209AD2E8F8404F818955421543522C3751EAD41409429B5018359A0881F066FE53C25E50
4203F0E2C008F372005255CB24BAF2D2E293D60B54FFDC7BB20C5BB64BF909622CFCB7A069F0FFA7CAA598770B62B7AADE7FA57E3B02C8014BEF3684
5194F18584C40B7251E0411909417879352ED156E6664F3EA9D8B54C92332B99CBFF62D23F9CB5CA6760DFB89C55E90F066E7A5593362F72E680D289
019B3B86F56A3BBC9805491A7B7A35A6773BEB79DB870E749FEAE2687B33A757CBFDCFBB9D04AE596F8583C79CC8F2D977BF52033C1A7149229617EF
BA282780A45C23CB0B78B0204E4CD3C69488325C2757633CC75549557CEC97A5E72868589C5500D04610C149546B34B39FE1A53F23F09F30884C346C
4410CC20AF9768201D202E5F67D32D8F61FEB13E66C91BA1E017F4EFAB40CBFBC38226C17F864D1BFBEAE85BD136DDEF6780F06B331BDF1930EA4736
4960B56ABCCE47476811385FAA16ABDEE67D1FF6F4D388B364E0123D17FDA94476BFEF9FE7ECA78A0FACC8EE3FC7C717DDB234DCEEF28995F75B057A
B7E85FC2620425397AE7A37D6BFBD59B7ABA8DE8E8643D2E6D4D8B79CF3DD6488FDB4C8A3CB1FBE1CDE3370FC623DA07E3AE942A9439476FFC304973
13318601483CDAAF440718223EA8498D14422902AE2852300424AFAB4451B02AB98A008C080E23088E0A6ADF1CFF85BA20AED71B415E09604600D7CB
C51868248AF78D69D726F027FDCD2700CDD4CCA9FA9C41FFBC07B56C02B0A069F09F6538A9F7A9D8F3A673A1BF4D00583676706546354B6865248C23
46184FFF84195598B6DCA4432B6E34FA0CB976754E2C797E1F857C7A9E68F01EFB658BCB9B9C69A90F7AA65E79933DEADBE641FEF3DE3B1D79DA72AC
4487D064E6D9B88F3636339FCEE8D1D3C5DDCDFD4B4AFF81F7574E2AB8D17A7A63E0F1DABBF33E1F7F6952DCD819A2954A8B7DFD74E2B8123944094D
7DEAE81FB4425C52A2472A9482E3B7AC444451801E4318485E101B2BC5701306F1AAC084920C611E0530EB7F4225452804D041268CA4091AC7402A7D
9443179B43F83F8D3F1C4B1B1F5DAB3829B88059E6FF2D684AFC67119F838D7E91D53334BF170095EE61D2445CE03C0612BA06100DC8C1F50A442D57
82D08DDBFA35DDEE5C595D4AACB9C789DF1D2AAA1DBF2CC063B0F2CA26F5843F3357359C5990D065C3AEF38BBB1D6A351608939264C9A9E0972D5BEF
BEEED8BE9F93B3DD9574CF7E471738BEFCDCAE5759813F19362EFCDE46A5E2C65395CAA8FD71F6B608FF90852214A32B454CA5202A1283154A65898E
2168B4B408A6690D7FA2606055C4AB1F4AD280113AB90647210C15CC40CCB93F4E509086C01158A534200C2F0A181C61E3BB5A8F775F0AB07F877FFE
1B22C88B483AA9FFD50665E1BF054DA3FEC7B2349374B9B0765FCE8AD0DFD900B29C61CE3D3657492272008309638106ABBB61C40D624C95ABC6AB6E
8B8BFA8DBC7476B75CB33119FA71F9802EB4CFFDCB9D6F012B7C8FB925DCBA50DBC777B7C7B5F517EDDC5B8D245E355094F14ED4A3D6AD47F477E833
AC8BE38C80C56E3B6FB55B9336CEC5BFE414A35D7431EA5076EDA3CF102A05E2F7DF506B1EDF833803892566220A1D5D25ABC96E804098A4C8EA2F35
048BC98C0C41EAB33F8667210C4918157A0D880AD3BFC2D25F1AA74894FFB1BA508712C2F51F216C00C671EE730747CF01E31A59E69F7DA70C533BFC
3277208CA42DDCB7A069F19F622B36DDD6BED950B856C1FE660A98C3CFCD56CB6B0CB8A2128300A22AD3887D7C85636211589663C0E3DE12AF9DC75D
7C784A5D7420495CB5F8A67165F77B53AC23DFF47E396655EAF4C27DDD233B2C3BB870906DCBC18A400385D75C4DDCDDBC958D43C72E9DED47DF9DDA
AAFF8991C36FAC735D5977CCBBF0F0BEB89D9FF58FBEA00A03F461C3E12C38E25C0507B360660E5AA5C1F254921805C7A9210A8F7F54CB30609D96D7
01B5A9DFDE97D3948984B47A92105C7F049A0BDDC07CF4875002D20BBDC1421F83302544332F3AF55C3A656206F31FF4674D07FBE90AD6575396C2BF
054D8CFF0CAB39B646AC9B9D15EC83FE4E00FC7FECBD7554947BDB36BC3E2CC44010EC6EB1B6BAED007BDBDD819DD81DB831103B50141515A5140424
A4BB3B87A1671886E9BEFAFA5DF5CEA0FB7EEEE77DEE7FBEB5BE7FBE873904166BC0B5E65A8BE33C8FB3D9B87E9ECD520DA6E7C94918D717A831E5C5
6202AB6CC15B0C00CA2843AFF69F59F6ED6A73CCDE24F8FBAC821C8725877ACCCEDC78C9CBA9F4D98DC2C14F770F8F5BE6D8BBC336D57B8CA234CFD2
F676B0B3ED3F69D19061EF5EF5EFB7C47DE5FEA3EB1785642D775C7339E0986BB1DF5B3DA4D17E5A7EA5B1F6C7573143B7C0C5C510C1D35494D5176A
380E35D0505AB094A145157A96D2263FF52D14930449E37A13FBC9D6DBBFC0487E9AA54CB77F49CA741BDC7419C0D40D4CEB0FD9CDB8B1CE311DFC4B
FB1BBD3F953BF20338F10CA6CCDEDF8CB6C67F9AF83C3A11F5DC66D857C230FFB1F79DE51A87FEE9AED2CB30592389AA012FDD80A69E6841E5224056
AB18B99F48BBABDFC1E08BD75B220E96412777A91FD8CC1AD4FE5A88C38755AF548BC22E4CBAD323E0448F81164790BB4D2425BBF1636B87EE9676EB
E6775FECBDC47ED8894D930E076C3994B365E693C7AE67B647853827C018F27EC5D54CC193181983170879A5145A2CD5E424E5111C575A462B0BCAD5
0C282A34D2BA29E4F6AD001D30E500E56AD37A0F40FC5E02CCAA5B8C2F13AD7780A8D6A320388ED3A533AC67CDD8E1F44A6F7A54E697F6A719D16617
4DFDE348409BF96F461BE33F4BB31573CFC0CAF90559973434F39F1900EDD97AE23EAA32E06A948448D55725027B072224BF9EA8C826B8061FB86963
EF09474E7FC21F3F86CB17A78A168F5A33707AF1AA4B812755F7CE544C769B7D25B0C788F637C10F3EA0F2AF179EB6ECDE79D2D1DE93BC778F18B3EA
E69465DF6E8C7D71D4FE86E2DEB4C7917CEF60050A7D5A73AB31D42D0426A98C608D1467AB32097E108FE690B41F2AA8A8916691B8628AC3737DFCDF
86E9701D85EB74180670D32A2053B79FD1DD73CD1554EB22401AFCDA0A4CA118F831A4D79835F7F65D6AF865E9985FDE5FFF7A8A205AF4BECCDCF66B
465B63BF8904E2837DE2A8E7BB89CF5F0886F94FFE9FA51296A72D4AC5A508A406881A4F0BC760C1DD7CB839B416A48838F6C90736EF8F5107EECE7B
D9E81C049F5BA30CB05BBAC6FADED341891BC3B56B42B69DB83EFE9BAD75072F34AD9C04DE7B8B0F75EA36E8E0E02E3B1EF6EF35EFC4C2B1AFFEEEEC
FCA2F38222BF69FBA2729EBF68C6542F377B65479C0AD59344BA5723C451799918EFBB313C178697602D3F9B6856F7B39A65A4A14FBF4426884D9B3F
F57A8CF9A5FCC12F0B40B17A3DD37A15840246CB607C85C0C88FFD860CDC248C5A9280FECAFD1BBF98A68713C7FA47448AAEAB19C64C7F33DA02FBFF
1BFF39F4EB9F9B30CDF26CF8632AF84F4500E32B12A72FC19E98CE0019836C0DAE0A56A278D2A578434B39A2D1B19C606735FB75D2B4F3A737E5F83B
D5978D78A13B31E6D8E62389B38F1C8FA5A6CDF29EEA3FEBC4C0AED68998484B02CF93D2E39D3A2FBFDA7D468093FDC2A32EF3CEBD1FD2FFF572EBAF
CA3DABB2CA633EB510D27BCE018A8F27F3301CAF7857434348DA374C1ADB0C40C1073E5717DAC8B2B2A826168AB9792FF46DBCDEE8ED6585429397FF
0DD33550CAB4069C224D5782680A8641EB2DA097FD070EBFA6A89AE126FB1DE83046D044E9C24375F7E1B057C04C7F33DA9AF76F6541D4CDC1C95CD0
7152FCBE96A2FF13FF59BDAB9338B6195637C3288936A2E55F304CE6FD56C335E6E008CD71F74FD2F8E33F77E5DCDC2EBE71CA70A75F6ADEB0A91E3F
EFEEDCB4FA1EBEAD6B8CEDC5CF930775E9F253252154E0D68186259D3AFF31A6DB319701C3F71C19B1E8C7A26EAEDBDABBABEFCD09A8FF76359BD0DE
DA91A4F23FC52370A2F4B3C0800AF393E0A264030587BD9512A591328E2D8851A8E2EF1C7EFC3AA8D674F1531896A46F65BEA9EA6F5201262140D024
41FE6AFE072485E10CB8D663C4AE205CB3CCA90C6D0DFD4DDA9FA129E19E6D528F52EC520ACD98FBFECC6883FCE70A3CCEEEA371F708A82C5C4FFE87
1400CB82ECE97175E13255950C10A0A5120BCE92E76B0ACA115C0CA3952A4EB525846BD835E952C03EEFC6DD5EB2A54E195B7BCE29C85AF5F2D2B0C8
8C7EDE67C74A9DBBD8F4C96B91E01AE2E2BD8A915DBB74EBBED8AD97DDDABB5B1C5E7CB0DBF0BACF2155F8C447CAF26B9130FC664F3C14B09F0F50B2
F08388D4A7A5230655160CAA3D5F3582F404150BF9A730B17F1D7BF42E2A56C731145E51D18C50BFAEFF983A7E29D22090E2ADDD3F646BF61F374A00
8C81F6F75DEB2323C0D171319A7F77FFEAEB2BD579EF89BAE32D66F96F465BA13F60E87FF5BFD05CCBFCB73B0A19C99D42E4D10F12B0FFC15830F2A3
EB9A6FC6A3C26A63308DC628E1B8BAAB09646D3D83355202BF6AAE68B7183A62DB7FDBB55969D97BB2321DF63EE8D5EB68CB8D839E3D16EACF0FFE3A
FDE9F38E56831B4A2B3042B3EE59460F4BEB6E9DD7EEB41EBCE7E0E893D193275D1CB5529E30D24599EBE2AD05419BE3749F7796D224197A4B8C6065
799022B386927F7DFA4EDD1C9185519A7B69AC78ED9990A80F3F508E25EB4B25A6F1FE56CF6F74F73405744D8DA6F91FA3DFC7F53045B45E03601B1C
07DEAF9123C07DE05B39696AFC677EB97FE899030F766BC49FEF83FFABFDC16C07CCF8DF1CFB73AC16FBB7EE578EBCBA27E2324D365C6BC85A51FD7F
4500BF7F098F1DFEB7AFAFA145A2313AD5C65C82270C5D5E02891AF4E50686F7BA850B7FA178643778C29E0DA71A836F89BD37C41F1CFAE7D37B0EC7
E7F5FA5A3BF0D9EB49DB2C3BF4CB8BA83672716F444ACFAEB676038EAE1EBCD8658A93AFCB9083F347E4E78E9854AFBCEBAEC2634E7ED77FDF1CCB50
C49B037C521A1C8222354A4A78E34ABA4CE6938851EA4BD12C7EFA4EE6E347F124C7F1D28438C5C0A6763FA3C737C6FB9852A6C14DB3FF466340612A
D30E40E33F4E387ED9CF2AA592F41F74A38568CDFBFFAAFC21FE6343888F8F08D5BAE7A439F96F461B31003454ACE4FEC9801B79103932F264224365
DE413CFE9252CCFFE43FC30817D9F97F496DE421104E2265189CA8BAB12D5B52D4525F8971993E3895F055747F60BF191BC7BC929DF0129EF5481D37
D2F7749765E7BB0E299D392267E300BB4E23F3535B4803B1272ACCA8FE2DA71E18306CFB1AFB556787AD59DCE59560B66DA4F6C9D5165DC4E15075D0
821704A7BDEE22814541A135B53C0D837A9D49D2945C0CC480E6A837439C3E177C33B099E3A4A169C6C01F0284C17405C4A80048482AD491A63B2014
693A01869346F22324279ABD4D0115D451B983B70808E677EDDF94FB4B9D7496901DE7E9E3C6240073F46F461B09FE019D134631FF2202273C78BB64
9F01A513FDA16DBB5BFEDB35B07FB2049857979D45858559384E125A154C8BE26B8E9E48132BF106298BBDCB64A83C7FC3A3A1D3B62DD95217BEB23A
EA8FCF2E03FC52078C709FDCF9D53EDB1FDFFAF76E3757A9C560C0BAE4FFECD1A357F7550707CCDCE9B471EDE035B707EDAEDB6779C7E0B5B70A7D33
FDA92ACCE99494939D72D780B82FA9C8B7AF2A4AF7F254B6226EFB1B0368DAF794961F5C7FF76A81F13DF283780C8D021AD51A006B1AF0D56860D3BE
2FDAB4F8E3571ED0240C4842E3F746A7AED05052C7D9F9E42FEF6FFAA041F5C23FD5644210AE4EDA5DCB32E6BF0D33FEF73B7F93DFA7187DC0173503
98D60E189A85DFCD50DEF262B464545ED3CCFDFF5608FFB72E21DE30BB6C43556A134192B17C12617322057FBF93D6C258611359F65882D14FBDC8BF
FBCFDF3AFAB668D1BEFAE913C2270C2DBFDA7EF5C19E4BAFB473FA36A067BBF91A08C119744D76A875D7AEA356FE3172C9D045AEA3FB9C719A5E7CAD
E38EE6470B8AA1F7E34F17043A9E1370A2B32F2934D24F521712AC26D047DBE2959F76F84144DA5A2F50B7EDE4993315C6775490AC665B27FD61ADA9
D99FC2E58D1AF2571DD05405F87505145622B8B4167070B581A2CE0EF98CD1BF9A7E4CDE9F569CECED4F357D9209F4C10F5ACCF2DF8C3661017E2900
45F053116EA475AB01604A5679D56E6BC6D548B92E77E25105FB3F0D00835DB5DAAE45AA4330928C0EA331920EAA54C5A40B2BC9DA0608846521B470
739AF2F0E8CDCB17E6FB4CCD0B1CF1FDB1F5B1CC59C31F6E9FE03178E0B9E1FDACB7E993658015AF1046DBD80DDDBDB6CF12C7A997FEB29FB57BC10F
AFEE4E0D9F464460FE434FC8A2D63F51B095E7BEB1C4671F4CE7176B003AB76DB96ACF2DE918F67A411495BDE4E6811B128695FEAC30EA178CC2D46A
943445FF0CC008CA54EB3386020C688DFB018D898506254A7224AE24C92743762A7F37FE99821E067A3F68A1962A2BC775FC4D7E7A73F6DF8CB61301
5040E87727193306CDA60E58B6F9413FFF8033A841D6584714AFDD26FBE72020FB6F0220AB4FA7EF8422260902EACF0A56C3353E8565899585F58A5A
0C0EBD8C8B88372B9A6B170CD935F68274DB32D186A9D14E1DAF075A1FBF65BB7C57E77D4E56ED9DF50F45804B9A25FF666535DF65ECBC59F6A7F758
8F9DDCE7EFB7F60E19A17D6F62E923F7D6953A3FD1D2B5A782398DF727A2F97D010E14AECEA532D7BD15247A7B56369EB2FBFAE6C71803270634B31C
8103B8BE1A32B5FB1BC5BCE91288E902E0AFA55FA639004C8F52243F57459330A6879EF51F1643D2BFEBFEA6BEBFC869FD421961BCA699F058F71331
F3DF8CB623003800E5F9BECA402986A58DFC4713A78CCEF1F02255328D2420CC6393F0D7669C7FF7FFACF684E59FF5A4265840D3E541040B7101A134
5C91AF51E64B0C19E72274BC96633775F1A3FE3ABE2CE3E7BCF20CBB53AEF6FD5E1D18B575C1F41BA34F3A5B599E857F600C97BF890EE8D963EC9811
3B47AEDC63DD67FB1F7F3E1FD4FD79EA104759F35C4759D5A593F584F4EF4F9830FC1B2EF24EC3A9BA2BA7AAB22F1CAE22901BB3B2553E672F6EFA01
28C9872088051CAAA7511D61A4BE49CCD3E0B7F6C77E6DFCC70902815100893414499104D9EC38E0BC9A66FEE9FBA191F81D07D7505C61AD4AD77424
280B33EB7F33DA8CFB37B20028E3BD7DAB598A6E759DFC0B138E35DDFD4A6A5B08DEEA92F7EB04ECFF0D06248FB03C86A3F2FA16351917CBA949EDD9
724E9ACDC7AA630CF5D1AE4D9955B507DF23F7ED0F2FDD5E7F757BE30C3BBF953DA6C44FEC7378D0BAF58EFB3A7738AB6B60088EE7823FB5EA6667ED
3076CC81613DE6ACECBBE5B0F5B1E4D1FD72149B47E625DC39F41E373CB9A81688232A526E15D170C6C9A7095FD66ECE253467FFE2C5FDE5F968633E
4D09BE15A246272F9413C677CFB6B6FD51E4AF757F14491228C152521961FC0EC03069A435C9E99554FDD415D504FD7BE887660D91EB4B378772E2B4
7C0175FF7A412E6DF6FF66B4250BC0004DE8BEC7FA5F8BB218E8DDBEE737A5E7D328999E2DB98C7E5E52C8FC0F03A0DC6B3F299902228594527FABE2
9AB8E413629677AFD1905CA06D0A8D880C457276E7AA57F6FD6BE6F38C85B79F0FF7F41A3820DCB3C7FA15535C971DE86DEFAD967230A7F8885EB4B2
1D3E62F1EC71C3ED9C668D9F7DF1844BAADBA8C096DD43BE7D3F937F38170B72C9C1AAA36E247E4E055ACFB3EF4ACECFDA1801371DDECC4B9EBC2BF8
762D4394A6888DDCC68B2B69966C4DF6D30CA0591A81480698F6FED300E5C5A8016EEA00202896E258DA80285DFE4841FFD5F847AB7D26E77E5A4F30
C55110A3B991F4B28A36A7FFCD685B0A80D5461E0DC669539D8C647297C65FFDDC78328D31A8B9CA40D66F6208C6FE37136094CC11BD462C557332B1
96AFD04609A12AEEC55688B8E18ACBA2AA4A6569DE27F264BE6B9A527B0C1E3B3AFD94E5978B23C3FFB2585132B4D7F91EDBEE2FEC69F7ADAE94021C
EF9E646D87815D7A1F9E6BEB307952EF5D0BC77F8C3D30FEBB64E7E06F1FB70BDF5CD43CBA91AF12BC3C9CA84BAFABB9E1CAFBE4E6BC2914E5AD5927
4A98BCE5FE8D7A1AC9C883809ED25734996A7BA6F33E34696A59006A436BD59FA258ACB2524F601A9D0E23680C501C55AF46CEF6F1D653FF3C0925BB
3D20A4706615680A33463A1FE31AAFA969B3FB37A3CDA4007E7941168E3F9943986A6680D1DC3C5CB38D97EC524A1B205690CA4538DC9483FF6E0058
C586B183AFD08C8C909492D2D4DA5C1DE572106F5A1187F3324AE48DD16ECEC9861D3B75373A8F9D7C2064B253F4AC33BEDD7BA59DE8717EE184E435
B6938AF402A32FCEF2906FEB36C072C165DBD11B460D5E7867E4CEF21DBD3C6B4F0C0E0C5A5352B0BB26703F8F28BDE0168FA86EDF7AE826F0BCF6E3
6220A4DA373F3F77C7FE4FF7728130BC8024012AABD3330C30A5FD00CD50A67A1E6E94FAC0B4F997E270294C6B61855A81E0244171B2901AC273F029
25FDAF48467AA9BF9F74691850C4D55274EA0B101A8ED1E63F0B33DA9C02507D7CCF674C1B7308266B5E5AD2716DD4916C5AA7618BA399FC39DB4A4C
D5F27F2B01909133370D4CE3E475744D392D2BAFCFA4E06597A99FCBE5861FA9BCE0ACD20BCFE415A33EB44CB51EBFC1D7A5FDD64F437DF6B45FFBD1
76F2F12E4F37761C5A85343210F5FEA2E42FABCE139F8CEFBD79A4FD5A9731E333FDBA6FAEBD32E85DECBC38ED9E2FFC8DA948C3EE03C550D3E9CB09
1EEF3E3F8ED8FB40AA715B949BB3EA70E0E110B226A284001851598E5158EB910FD38C1FC39188568B9AB67B1BBFA7080CD5EB7108270D1A354E906C
6E2899376A8F00FC7A6496A6B4F7077ED26F7F4F4813CA31A2EE582A7A94479ADDBF196DCE00B04C4644A006B4CEC942B7B6111F3E50512742F528CA
E407D3DAA3C33F2828FADF070575D74F5FDD0F33CD7AAA5A49ABF08666AA61A41F79DC0BAFCED57D8D858ADF67A09FE788C27ADA1D2F7C36A0CF87AB
5B8227DA1E9ED263D3A0A52B3B8FC84BAFE418F0D61359D6BE87B747A7858E960EA71D6DBD72E72F8BBA34C2A374A51F766FBBE8AEB72EDE6553BEB6
FAD015997CFFC680A2CBAE02F5A75D29C2939B7F7C08D16724C80046C9D32B219234C0A451C91B637F8AA62104D2909CF129580A82090A6875240E28
540F93C6D84003D54D9B504CFCF31414F4759817E6F210D1FE4CC17590DF4B4DFE790D65E6BF196D8FFEACD62BF7336EEA962198C6B95FC88B0574EA
9937021CA1533F20E4A7FE1B2A31EABFB6E4336CCE881BCFBD39BC9A210B2490D02013503F069436CC6E36A4681B7CD26A1A122BE11D0730D78E83AE
5F9A6DE994B7E2EFCB03A71C1B31C871F0D8CEBDB3232A8DD1F81B0FE9748B5EA7EC07CDEF3866F654EB2359B356E55E1D7425C7E926783D2530F42C
DFCFF16C065CB8CA590DBB3ABCFC71E63ABFF1D1F11ADEAE61A77CA2549959064062E2CC4A021850C840D2BFC6FE6025CEA14A00000B14F169108D03
D30160CCF805468CBF4F2957D88521CC3FC13F1E3BFA21E576C920FDEEA3D12139DE35E8F597E6EABF196D90FF460350F421FDAB9A224CEBF1236728
4AB7F2F11A9F27E5244255BD2D26CBE78E085093F47F4D0A42177A3FF02FE2A48D92DC0C8DA8492F9411F716EADCF7A199B170D09B7271D27B79A543
20B4B19DEDE02D7F0E7A7B6FDADB0D930EDD1EB572766FAB29F5B5040733719F1A27B4B3B2EAB97A98EDDC99F68B93578E8C7A3DE458E996ADB2C469
274AB7BA7ACD3B95A38BDBEB94A07F6D77BDECF2A5B29A439B0A84A7FB2C0A6FD6C76590146C282F9552288612C6F74C53C6A08466741203899A2A01
4CC1CB449DA91A606A0520088CC44C174020679B9B06FAF75D6F9ACC9874977D715AAD4D7C5A87EB0D3F130DBAAD3F0933FDCD6893FC071FD2FCEFAA
31531C8D1D3ECC46388B94D23CBF5C424519324235F04DBB3D7908F8AF34407EAFA13FDD6A88C2F4E68A5C457E9DBA4681BA5EAE5C92A2FEBBB8FCAA
6795E8F67BC4D3A15EE4683968F94E9BC9BE13B65F9A30FCC7DDA1B3BA5ACE69917014A07D5FE4F4B7E83A78EDA4CEA3470D58ECBF7A40C8B35E7B4B
B7CE6C2971DA987E7259C09ACB7C75CACAF799F56EA3DCF26E9DC8D7DC9D1A51B0D57651A201CECB4201A5CEABC58D94864D298B56E74F230862306D
F6A739B430824F7334CD9A329AA6D13FA3ACD1532AE72EBBA5BF558CF1798BA63E6403AEC9F4C55FF8A45AA1A86A6A8E762CA1CDDD3F66B441FE9BD6
FF5FE1BBBDA9650886625B2607B121B79B6BD592A062524580C68034B474F7989B3CE25FA732A14B36BB5EF9D155E1223CBD542880F30A3195EB8BA3
C7B0E85BA2F48701B9E9CF44EA155BE1E23F86CFDC34C5EE82FB4AB795934EE53BDADA5AAED144EA695ABCFA79C5D00EFD76EDE9DDFF8FDE43FC8EF7
FB163E686EC9BDD1D1359B567E735F90E07FAC164A5FE62193BC99FCA0E9C6FE5CD983C55F1356F55E9A8EE1E9E98021D539CD803175F752C6C89F61
8C5E9E54AB6163504FB31422AC367014CDB5CEF6D0140920190A23C871EB9D22F25FA9FFB4C9AFB8887B7A54945A65341C4AADBE8EF774B7D24C7F33
DA26FF1936F942D3FEBBCDAC915474F8D002EA9BB7BE5A068515030422B154DF3224E7C0B447A27F1AE7A882DE7DBD96A4D085D14271725989BCF456
0356F8E2F2907CECD66BFDCFDD1764618FC5A923DE02CF9E76530F0E1BF774F2E6A3D38785DCEFD2D37285C15B04C8A289EF0B6C2C86388FB49BD87D
F4BA558382721CA6673E18EC235DEBF8D36BE4D384ED9978A9D30545BED7CEC8A20BEB0BA5B7267FCE74B45B594C19D262352CABCD9419637E205098
527D26F74F6B850D4A63BCCF324DD9D54AC61813700C6BDAF60F685C5327D2E98167B72DE2DF028661F1EF833DB9C867B05E939FAA8154B016A9FF26
FCFA11FFFF1BFFCDE6CA8CFF4FF86F3A00EAF391B7D95D669A9D81CE4CE6A3DFC2D42A9D36AA0431C6D8401BFA9C4FC66C9AE32FFB3D2A6878E410B2
7FA6C490942A11954525C95F3EC15A5EA6386E24AA6EA6F2EE6E7E56FE2E4A7B7781403CABFBC4534B465FDA32CF6DC58CD34923AC3B1F2213711AD7
1E09CBB06E3F6187EDE045F6F356D8DFAD741C101FE370B1E5FC786FDEE64D1597DEE2C2D55B7899570EC7484F6CCA115E98EACBDBD463331F40897E
1A8E5096AB4C7D48F5997AD3B08F6972495BC15391A451FB0BC36B602DAC474C373C68E30F290A97290994CEE9B3A189FCDDC9C4E81FDB3DE6125F40
40ADCCAB257108C2850F53E58F5218B3FF37A38DF29F61E5AE4DC5DB9FCA6902A3846B973742F1715A89AA25B9588A1010C5343F0867D4EF672C0F68
448DACC231C1EED9B7BBEF20243F139B05191FF9E27389F8A3173143B3D094876F6B5F6CCFCF3A1E51BAD5591DD2B3EFE85503673CE8BD78F7BCE15F
8F58763AAF8D03184EDEF40DEF68613FABD3C8B16386F6B89E3DD5CE3B6ACA01C5D3E1AE1927677EF039DFA23BB3B42A74DDADCCEC63DB0AEB9CC785
956CB0DE2B60C55F3F0B01A3126B4D77FC9A2BF546B1627CE31403F8A9420A60465912ADC608030AA128C91AE37FE33B255180E36C7EBF85B524FB7B
EA477AD2F63397F31CC51564ED4F01A9D691628F2043C8F21A33FDCD68A3FC371283CABF86571DF41792044C662DD8A0848A4A75A85C211394C11446
7258642ACBD5EE1AB1F44D910E1030923C64CFA2CEB7B0BAB41269737A8C26C345C8732E3CBE0ECF294DCAAE8ACE94BD78DB9CB92D507FDCD6E1F89A
8977778EDDBF6594FB5DCB4E9E86749A200C2E91D19D2DC638F45E3A717C978D8D1B3B3A17CF9D5E90306F7BD1CFF9A7AAEF87A1BE733353FE74AE2C
D8B5AFA47CFDCA20DF65F63734B8F043889A55CB24A61E7D56CAC3680237CA7F0650589D1427688014848B711D0A8CDA00333A7F8E3346FF901641B9
947E4ED5F8AF655F2C102C1F90CCC94375346C201B4A5555068A48F8A1D71D3DA833F3DF8CB6CA7FD3F9DBD76E34FF7A582302A1129FB98721B83C87
40A40620AF911B03698ECB15018E7C3DE1CF4D8F2B5114922CE9F36AAAA53F51195BD95CF7250DF67C88BA1D4D9BEA9B1584D787E5C764679D89179F
5A2D96CEE8B37AC380316F570C38356FE166CBF6D79B823002CB59101F6C6139C276E4DCD156D333AE5A4E0F5E3B2C2676D29FAFCA76CDC8483B5C58
3DDB4BE0B4A73C62FD5941D9DAFDF9670774BF8D32855FD2301C69D0B6DEE72B2B07C02824301C35C6FF72154C015254A781655A9222289A32C6041C
6028BD4405D15C88DDC26AE2F79243A278D6441E47546A199D16E81A652DF59CAAB61C33A85CBEE0E63F0833DAB20068DAEFC1943E2BA892E348F5A7
B1970C58C24F8A90EB70A20936FA528EAD317ECBA56ED8B67C8E9F044113071F8C1AD9375F959FA254E77D55CBEEA6252D8878B6ACFC433C2D4DADA9
927CF490A4EC7C0D12C73A8C1B68BDE0F6D89D2B060EECDAFE06FF2701B0F0A5B5D1969D86CD5DD0A79FAD4F7087216FEEF671D3EFEDB320E8CD0C3F
E2AC27E5B3A3C67965F9CF55379A0BD7BA666CE968ED81232509D508535FA134797F34BF1061088CA4F428066439224050A26CB9D1106849401124C5
72408FD314A657EB69E2AEDD925A9C6905AD0BEFFF978C8118890AE028A6C101C38AEB2A4905885D9863EEFD37A30D2700199AC89AE6026202AACA9A
516DC5F329A70D6A9F741A6865382636E5D339A296C0D58CF4EA8D95039C0BB4EA2B76E723462F143745FF546A42B3B1CF27450F7706BBBCCC3EC5E7
AA797145292B436A9739942353BBF4E933D172F5E2FE0BAD7BD9B5FB9B9F034178D1FEEA571616FDB6FFD9C96A40C032AB7D9EA3A656E58D1CE8E239
EE1C19B7B1055AF9F5CDB4D48085B7D4F94B4FE7CD6E37E13B4EA66623944E2BD6B02CE08C41066DBAF28BC07A0C147DCCC21828344249507205A070
DC74EF8BA5541A6340402074ED6A9BC53C9CA24CD540A2FE789F330602624BF24892D6292508CD1955438301354C9D2B32CB7F33DAAEFE37896343E0
E84D35C78BF995884E59FB6DC10DC890516D00A8C640E08CA9DA04E42440E99C8DC9D7FBCF79D922DADAEF52ECE40B1A7E70A44EF953997721BBFEE5
FDF32FB1AFD7EAF5BC6B5E95D76F6A6FD86D343CE93C68F092AE230F8E731E6F6F671BDED080A344CB83AA93161D7A2DB01A3C77D8EA8107EE3AD97B
EA36F69BF5EECC8206ED5A5F36607BD572DFCCE59765696B2EB9CE6D3F974769A24B31A3AC37AA7D6314AF17CA7116C7291AC70156F03693C22B7C53
3414A2D711142048920234CB5146978F2054C978FB55A5886915300DA0EF33FF8A83250AC04FD0D23A44AE8651B6A580E32B292EA1FB7DC45C4E33A3
8DF1FFBF0500A68178B7C1C9996B2B2AAA104CADCA5AE241E1BC6CB9C96112344A19FF03A6C6F528487AA6FEF647EF2345097D2DEF264EFD262F7C12
41144624BE8BC97A9EE2BB3210BDB7A541FE6E675CECB2E0ECA936419A895DFA0E9B35E552CF8DF32DDBDBE716E713285A76B7D0C9A293B5ADCDB695
DDFAEF7EE5DCE34FD90FDB510F5C464431E757280C6BFC6F3B17AEBCD292B0D1FB6DFF0E7FD65058580C41C30CDAA0C128966894730045690CC30191
1BD4849165E9D59401420C18451A811380310DF25100A58A87F7DA508252AD8B801BCE8D7DA337A4B580EA0C14AA53233A03443707422D15A49E5DDF
BBC0DCFC67465BE67FEB258CEAA3BBE16BAFD3F3D438421832D60432848497216A340006824CF9710CC5210288F974F6BC9E73A336B4EB99FC71674D
79F6B30C7562444089C6FF8170FFDC54E1E11B3AB1EB93F227BB2BBD1D56E8BCBB0EB1D9BDEE60F7E3CB7A76192CAA2A257082EF59B3D0A24B779B5B
5B3A8FDAB5E8D892610FD3E6DAECFA36F31419DBFF0BF9E254D3A2F02B87A5F15B13D2465AF68964B0E20403A0F56A7E1DCED0689386E570A397C749
02CE2F5231E2EC7A82C00D10D67AEC0300D41801502CD7A2E3E8987EB67B2A51DA74050C4E9D323195E4872A4153314E560848528D91860F15841060
9C61D045D8ECFECD68E3FC672843DEA137CDEAF86FD9080E4186E039A9340C4BC44551D518414A5A701AD143A891687A0DD3B2B3DBA2C396ED970A4E
3CE0A595BEAAACCA3BE36E509D7AF5DD61ADAA76AD2799B5D533F774A8F4A06D64F3241BDB5973770D9FB7C9B2FDB0FAA83CA3698976E73BB6EBDCF1
C2DD2E9DCEAE1974A6FBE50C87AE0B3D163B36354D5A2F8F5C58746FCFBB8D82C0C33FAA17749818691421E510A5176A9A9A484093957C8A41315C8D
03AAF6731140F2E31B4CD1004E40B8C900A0388A4128603579247CC3C6FE761D6EDA0A88B53C1CBCA3850C0DD011C20A5C53DC42E00A0D8607E4325A
0E52732F07569BBDBF196D3D05C032B4F053A26B3659F0B4C28042B0D475298FD4E801D1CC478C7A5A2343180A211084642098C16E759D39AA63A717
59C7C54905A2A05879F08AEF5CAE5BC1D3A95E74C4B22CF59303C2F41315B993F61B3CBA8DE8D367BDD3169FEE1633F1A4621223DEDFD32D6FD7DEF1
A3A5E5CCC393EE1DD999BEA3F328F7A343BEA3C746E4974EBBA9DC746277D68F73E561B33B2DA86029453320D0A272951AA101A22857711445400616
CB7D55422802FD9B59D234E2431835010D700C277088048C124E9D633F35504598A680D0828DE3BFD1C03D10278475FAE67A891ED269092C3D8A6248
0EE5D423AE0033FDCD68E306C0542027F33F652D4EA7CAAE26C1B05A5FBD65592D21D1223489A280E1B40A93BCC670A592A1608A8E76E862D5CE3AF1
C51D8D5F63BD4F3279674A2C1772A6FAD1FC4A62FF2EA87ACBE3EA411BE457FB94E6F7EFD9A5F7D231332F5A5BCCA7EA6A8D06E4F651D12C8BF6ABA6
771A7BC0FEC0ED9169CF3BF53DE6DCF5A2DECDDE57E8E858ED39E760B0C70569D6E48EB31B6950570433504EAE4E063058A4D2A831022564102A8B8E
69A0CA5EC71AC94E18DF1682620092EB0914D3EB081A612057DBBECB33B4A6E54040133879BD88535E0806065193B2065654648911189727EA489441
49EEB45D9ED9FD9B615600460140263E8919F58CCA5BF1450C2B34194B96D7E06A196214D71C8552900826090C87958AA6663D02AA37B6EBDE6E7AD0
DF59192970FE55A572C71FA9F4DF477E9C3D9C9BF447001ABA2F735E17CF82D15EEAADBDA7F59AB2D0D6AE5BFBD954612D20917BEE8DE32D2D6D6CC7
2EDC3EFCD034AF943156CB9F4C5B2188B1D92C393024217BDAE19F4F37572A17751A5F44C9B2F361022ACAD7A13AA02D2880294229D72110C54BE1D1
9A503F3EC92084A913880498A2B81A25304C031950963C64DFFB040F33D11FAE76197C8CE432CE66D07A41A34A0E4B338B2AF5B00A92AB71CC181A70
E93D3669CDD1BF1966FE9B56E2836F37831C2E30BCFD3B0AD52DC29885F30A08AD5A0D744A14C780BA418C61705D190C15240A5102BAD1A95DBB91C7
5E97F927EBAFB9738D0B1C788D9356DD74DA51E9BEB85679FDDAFB6EFDF3D638895FF69C63E9707840C75EEDA7113EF924825CFEBBB06FE76ED6CB67
4FDBE3D87371C3A68E4EF7C7CD2C4D1BBD3AF5D098D022C795BC4F43C3995B1D279610B5893C3DA7CD2A45A5325C9C5EAA2330698DD640D3A569285D
FF3052CFA030610C288C9F94BE42481138DE2C876156B2BAFBD877629431B2BFF9FD1F7F4672CCDBBB624651DDA094E32D594D00A774621443F438A9
A7F5936DFCCC773FCC301B805F17BE88F77B13B7C67020EEF4535113FFA7D3943412D600522CE10082C132A59A94C6E6D168E177010AE3615DDBB5B3
77CE8E086E4A9C11C5954F5DA774EDB9ECF8DCCBE5FB1F6BB28F241DE974EFFDCCFABA55F3C65B2F9FD3BD4FC7FD4C4C0BA5273EF9F38677B51DB17B
CAAD5D63C6BEB96737F0D9852E0F793B57663CEEF15CB97353317FD11124C07A4C09252E13D12494970F19EA94B5D965068646642A8A208AC2356CC6
F30206A0304292084191985EAA35E5F96B6AB514FE6550B76DF90A0A27F4822F0BED5C21AEE5730AAD130B9A15B0A224ABC5F8CC8CC60030BD49CDB0
2FBBEEAA24CC7BFFCD30F3BFB508C8207F2F2BBB91CD71E827F7C4E2CAC099C35FA140AF672BF93485901889A08430924790554125069488EE6661D5
ED4889F773EDB3A9555CF634AF8249832E5F1B7F276E4B92F6B573C2885E2133EF2B9D67ECB09BF5474FEB7697183FB591AB37BF94F7EE643DABCF9A
8061B62EFE836C562CB5734D7674ACFA3AF842CDE9D5D5C2F5F31B7D6DED5269194F025164461AA2E0AB9AA3F3691A2625329C5665C6E825B191728A
910AF5841E2650A3D3972B19025095053AAA6E83A5CDFD0683C1A0AC0A5CDEDDA99C125557C888E626BD02D6CA73621B399CD1572B01A9D311889EAA
EB3B3C32FFFFCDDD1FB35230E37FAD01306DC6A6657B96D5ECFD6AFC3B6FB8E9EA9FEABBAAE79A4C0AC7D1D4188C422014D3AA34F2B2741D2EFC1A23
4789F4811DBAF6FA907A3FA9E5F06221F7DD29DE6BD889576B77445FD9DE20750E79D923ECB963FAB191ABB7ECEDD3BD4FBB17D481069AAA5FFAADA2
AF75AFA963DF1C1E36F5CAD67E1B9C2DE78BAFF4F2AB9AE9D8F8CA315B74707EEE7B3B9B4F0C562F0238C18B51520D025C982920500A6830BC3EA65E
5711578270B0AC4C460023FD73D2953A88A669A42807023F4658F67DADD4691A2ADE2CEFDAD71F64C64B2860286B213088D22767EB0896010D751403
20A32D83F1233D5E3F798C9B377F9A61E6FFAF4D0014D9B0E62A7AFD8551138384DD272EBD3BD2DBFA31A03132D5B3029006032117E0487EBE14D726
BC10E0A06A4C47CB11F1496E3F055BF7E9C9F37F151D1FE77265F6F6E0F5D7A1C0A959AB17670C3EF8B9CFA0CDCEBD2C6D2D1EC37B0514E02D8BC9EC
DCCDBEDBC1C0C17D56ADE9ED747DF2AA608F0957CB364E48FA3ED207BE63E3FD73A8CD5B068B2FD01344E1E33CBD48A4925664D5139C5E0161454135
C8F7B7520627346A3D431A308297A4A33538832B6BAB08D4A34787C9E11295B8F8CBF6F69DF643F48DF780E3D4C56AD28022A2B44A9463A8D69BA60C
46A03044E55AAF4D9BF39334F3DF0CB301F8DD0644A3C5EB1E106F6FC98C3FD32738FFB5EFFAEC7EEB6B8C4173CDD3B0661207984E06E8A66C298266
F82B70A6616027CB5555A16E95E5CB6E81FA996EE55BD77FDAE51810BAFE9D6475BCD429E0FCB6F43F278C1AD0C3DAB6BD0F9E4C5000BD519161DBB3
EF90A0CB03A62E5FDACF79F6FA96877D8F959D19F8B96AD9A59622A715E7470EF2E798D400398DE7DCFFAE681141484D9D5CD3226AD2130D9FCB0C91
8F6B285AAF27599AD2EB217EAE8E4670C680C87540B4BD739F6BC53A594D9A876D7BA73816714D343E9FBC4245E0064C98C6FF75F6D8040A61481494
8E1E10F374713D6D1EFE33A3EDF29FFBB7237F0CCDD0843E63DD5651CCCE1898A3692CEEAF791B4ECC5F704EC9B244A15F595125A06004A2043C18C6
535E1B182EDFBA53B777A27B576A23168440BB87E6E46DDEEC36724DC3DE41BCB83360CBC22CA727AEC3265A76EEDDC9C20B91D1084D9CCB88EA6A67
732E6090C3C139E3E6CEEF7E3262C2F4A48F7DEFD72F5FD9D2BCF28F33433ADE65B89A503549D67FCB1454E4C33AA1A04121882BD3B10D7EE91231CF
18EB6BB414C991E2465D53911268210AD521049539B1D3B210AD4E96F57CAEFDBE388CA90F6A313AFABC2C1D30686014C24D1D8FADB38C1C21276982
568DEC383C63F5012165CEFF99D196FD3FF3AF5EE0D639795C9EB07E7AB5E4F9F56C8244097DF0B42907CE0D73382D2029AC39E7D5731D87CB10B23C
5E0BE1515F118E8BEDDE799A9277F051D1D33D09F786B86B12A6B9EE1CFEECCE5067AD7B49C4F8DCB33BE3A6EC1BD1C5BE9BC5276D0BA565C0B18AB8
6EB693B30EF65E77E14F87D5A3A7AD1FE810FA7AEAA1923B3332E50F266F98D2618D9A6B8AD4701C3FA389CCCDD741F9693A54515F66A00C3FA20D90
82C26194021CA9E757A8C47C14D6C21450400CF5CEDEFE7CBD565FF36ECAC003452445569A86FA6BC20B6102433002001AD0ADD6CEF8CC92168E06EC
856EFDA715AC7D2263CCFADF8C36ADFF19F6F76E5F5305800184E4C7AAB57C50E9EB91801810B2CE73CA8AFD97375CF9D18C93381A7EAF8203F26A49
C06718438303508EFBDEC3EA8AFCE3FACF71FBD63E9C31A25CE736F4EAD0BE4F5676FC5EB546BA6549D020BFDBD39C2CBA76B3F0ADCB0534D5B88117
646935E7E988C907163B4C19B4F1D61F768F52C7EF297B3436087D3B6CC9F40ECB9BB886AF46CF2D2C400D65990AB2F6732E06F8792A461F12AB6869
5660B84E0E23B424B7520D29608302028A7405A73BD875D4DB16A82578E9D0E961384DE8EAAA180E2AC95319350C8253244198DA8138C628F6195DA3
16002E7AE09B911F12277949CDC37F66B469FDCFD2FFF0BFD502D00C218BBF34CF4B0AC73A7BA31A15606B7D57CC399DCEFF722E8EA059A1670E4949
AB45915F1424161D22E1B80FD63DEE56DD70CFF9B4EDCC89015F50C9C2A53B07EE72ED375979FC73CA80C70B9DF9CBB65977B5EA10274A02002EDDA3
FCDCDDB6EF98B17B2F388C983EEEF2D6614131335715863BDE6C4E9BE9B064C0FA46BAC94BC071BA3C03539DA3214439E9423EBF508369C3C21478B3
01C0FC622D4B8B2BF8242DD2235A98E005953225337B3A7E570B63B78D98EA29A5A5A4B64C4D612A5E93512710A601411D89B568289AA14C1B83A55A
94E004FDF77B9E449FECF6D5D166FD6F469B160024F895FCA77F2DCAA219A0CDBF367F6620A1092CC7480DC47272F7EE7DDDF2DEACBD21357AD5607F
2550E6947D78DE4C51F1F78A39F66ED7DE4FB3F77EFCB967F4C9590E42F489ADDB6ACB43873AB8974CAF5D3EE3BACDA7772B275A74B48CE17F4720B6
748BE4B9650F6BBB5197B6D83AEC5A3AAEEF43E996312905337734D52E18B4B8EB521125BE5DC4B044B68CE1A541446D62A3A0F85B9E90517E0E564A
1BB4045E1DC34769594C3906C9301226E09CAF55A46FAFAE3BE2EA2276F419F9A48951B4B0EA221D406A1A218EA62800E37AA1545B51698C1A4882A6
0915AC90738605836E2D11CB575CBC0FB166FE9BD1A6F53F4D00823618FD23CBFCDE954741FC97F3A79E12C1157A0AC1B424C7A52CB69CFC23E3D4CE
709AA3F2BF66E1FAD2A2D0CBF914117D358F259DAD467FF8763FF5CEE865576DAFA3E9BD56FA380CF8B8E00F7CD3E6C7C3EE8EDB5938674E47CB5E85
FC6C12E71A4F68DCAC6C7B8DDA716B5CBFE587060F39947467DCA7A24BAB23EA5F8C5CF4C7FC0A0A7AF49326B1D4624E932C05953E0D48F5E7241DA3
FCF249810BD4884A54584950DAD462C0CAD5284128D31335B88775BF335FC3B7F7B53E984FC1554ACE90A7A4E4022986338065000610BDAEA1C9B41F
880624AA55B14A88B939ECD5B950EAA3CBDF9F0973538F196D98FDA64F02D263908AA4E9DFFC67280A12F8EDDB702A9CFFE6BD90865135C6C2DE2B26
3DAAF49F7851CE319277214A3435EBC7AD029A8DF0507182411DA6647ADEB93BB7C3C13DED9F36CE6B77E3639F553F467E4DB2F5D8E4B8DD3E64AEAD
9D85352F3F9386889AFD02174B3BAB3E6FD6775BBC72D060CFEBF3473D14BBCCFA2C7E3DC469D2D00C46F13E1ED363E1DF684D8C84E67B96CB2A1372
651CE4EBDDACE629D0FABAEC721DA90E8981B495B5720A2B492B221487BB4CD97B63612FABF521382128D6D1BA2A95AA49A65360B091FE24030C0A8D
5E4300601A062620894EAFA3026747451E46EA1E642D08A5CCE1BF196D9AFF9C69AD9E86D41BC5F12F05D0BA0E046E08B89FE49596FBEE6E358361B8
1607E2A77B76A6C40C9D514CE3FAD0EB95487C5DC1E31A964B4A05F4138B6EA744FE7F3F9E3CE2CDAC76AE5EB6A3A35CBBED7FB18ABF61D5B3D11BFA
7F38693DBCBD756556158D3275E745D7BBF4EEB1EDCBB0E1DB07DA5D7D346DF0B4A2F74B1F47E7FF3560D9E84740FF3505C7E5211F60E86B312BB99D
6A482AABD253CAEF1F0544931011E68A55284D9746AAE94601C2081372203C7682D5E1BF1DBA74DFF293044DC955285D5F4D508D5A95D64010A68D81
9CB6B25EAEC108C274299C24610DAAA5CA17C6E9B6E5B02125A933B3CCF437A36DA7FF8C1F544DA5CE6000BF3200BFB2804057753F54542900FC273F
29C2A0D3497470C995158FA2FAF54D011A2CCC3903CF68FCF944C0E1AF0A19C5540B9B1CDD9B93EE3DE7C7CE9B5276A8D39A2F36D6416EE70A261D5E
3C60C065BF4E3DDB5B65FA14511A2E73B7C2C3AA4B1FFF45D67306F6B9BA62F8F23141EFECAEAA1217D84EED7900543F2CC6B5F0CF0F4A322092D3DC
4B1547E568B59CE65B905A555BABAE2C6C4218DA20CA1343E5A532033FA186929CB4B27D766B70E71921140325A74900595B0B51B0562CC72812B080
A3B555555A588391044599760435097540B0D29B3AF8944E4CD51DDC613EFC6986D906B04463B20206A6F21FFB4F0A80D4973CFAD458544FCB3F45AB
B58D626943BD5615B5F46050EFEE3FF166A4F0666A433312E2CEA3A4B77E90A50E164B21C98C6947073D2FDAF434784C4F2707ABD5A2D55117469D73
1AE65237DCBADB10FE9B6286A0334E098F77E932F34AD72193FBAFDBD07BD3F8938AD54B92A060FB7EBDC7E64AAFC710444B74A0044F7A8B10F7C3E1
2FC918493707A58A918CEF1245B11C672885A14E4A5497CA416396988D1F6FB5C0FFA0E5187715431446579214999D0F6888D2C94CB7BF4DD74AD5F5
6294444D3B03014500A2B440852A8EDD463D5CF086147DC9C407B899FF66B4F90C20C392C9B548EBAAFC7F52803441E8738FAE7991D5426903039A15
BC4A4553AD04135DDAF0B857E728BC506988494C80F4B74F9793D92BD3099E83C565C4C762D2AE3E89612BE3EEF71D75B45FD7C2170BB2C73A9FB29C
5DBDCBA2E36C5D522D45A28147791B2DAD27F51D3CAEEBCEC51DB66D9A9DE53E3DF64BD28E21C3ADBF28CF24537245BA2FCF90F64CCD78FA20810918
C2097D0A294DDAA7F2CA7C11C1496FC7A10A5146B61CE747C959EF1E5DB67F9A6BED2A05A03032076568F18F2C8CC070CC40008AA401C5E06A15C940
3849E00007149A12DD8CA9CEDCC27E5CD013614DC4C9B109A499FF66B479F76F3400D244A192A1FE9502A44D9B73E1F21DF3DC5EE7E8F0D8009EB2BE
5AA0173493C077F7717BEB0F4D8F8554634680A2E8D88D30326C4F33CEB369F74239CF62C3EAD5C5AE77B3368CB8ECDC65056F6CE0EBF9577A0FC9F5
B6B0DC41F0E50067125F6BF677E93178D1AABE0B4E745FB466782C6FE83DCD8E7963C6DAEC0061FE80E297256549EBEE9731BEB7A0B0381210E56F73
00531D57DBF4B385A025C79F23045C64A47FB28F9876B7EE3E676BAFD5851428FE14A1A1196154402D49EA308AC08DDE9E323E0D85C0A44909A03069
7C01898AD6113A8FEB48DD7D0D5B9A41E64DD959C09AF96F465B677F6BF9BFA0A4B401637ED7008DF4A7288021617F559779BD12D2D9BEC5CAE6FA66
8DBA19A6781747B5EBE4DD7CA1823124F828E3EF7B7DA53CCFCA89DC1E166E8FDAD9DD9D7F326EC1D3F051FDDDA7760E7359281EB0DBD1C235C2C262
8FE47C1D69007151E8424B2B9B8D4346BEB09EF0A0DD26F5A6B18931C3EC27D88C6F2A7647F1466572486DCDB34832C859141949A234CF3593D15624
6ACAC2542459EF16061869AA311CF9FC584D1CED3A706A5F1B379CAD8DF22BC309FE87FB597A866C166200474DEE9FA45855A596C47002D54318A075
F1E93080DD8F28EB2FF3A9E6A81AD595D907B3CCFC37C3CC7F966338B8A22C311BC6685313306DE2BFD163920A97071C91FA220D2FFD902855885B74
2AB90118BCBB5A0C57666E2E62A89FBEF8E7D817F99A27678B88E281EDC62C1E77E8CEC99A3B4B526FF638EED17349F434EF254B8E58ACF46D6771B6
6EBEC0A8291E7BD6FED1A397DDF0E177F7743FBDB29BEF0BFBEF8173A7CFEED52B4CB7359D93D73544E5CB933FD404AFC98AFF885080FF20568F4ACB
CAC5992D80CC7B50CAA0495F731555775E62C8E62EFD465B2F4CE25A22A3AB14686D5E788E96A5195CA431397C53BCCFD0EAB01C0C000A431112A790
DA5A6304F07A5944C9991012C91168A49B8F6CE17166FE9BD1B6A3FFDF0A80C293728AF2C5AD977480C9FB1B41D0D54BDE72ACC8FB1BD6E81355A316
0B6488568992F1F616B7E98773331926E8ABE446FCA38A3AEFBF2464CD588B916767EC3D752FD6695DCE7ADB3B2B2C4F1C9FFCB1CBA6D13DD775B470
D179680001DE7A970DB7B4B7EEB4CEABD3F02536274E8C7F12B975C8F6915DEFD37B5CE01A5153408A38F6555DFCBC8F097FAB68D0FC2212AA2A9169
F15A0985A73E6BA154F75F146AD24E06512D4E5D7BF51EEC8173C9C102846AAE2A5662A69C25859B368212A801353E05529E8B3294D11C18AD198AEB
70A329F0D91495B0D3532329A96F8603E6EF5ADFC0316603608659FF1BF94F8A3C4CF9FE6603869BE267D06A00D0B4258F59D6E0735F63C888AD14AA
F9029C306841E9509B6AC9E4E9891C78965DEEEA7FDD5F707BBF8E6CE8DF6EBCFBFD95872A7CC69D0FE8FDFAB6E5E8CF8B02FE9C73D87A70370B5722
07277126365A37AF8B8D9DC3E78D9D87CF5D706ADED3C855BBD7CEE8B4810877ACD66735C57E6FAA7D9DA23975B1E87A0D0DA44FA295027E89866BCE
A4A0C06F1A203EF348260FB998030AC65AD9749A59CA617783594EDD5C26C54808238DD10B4903122320180324A45693A6422645B3A60FDC286CA28F
27491FBED2223C0524072E41070FB470E6E93F33CCFC377501E3D08FF56F7515F935B55A8430EDD5375DCE43D5D14E57081AFFF8B71C1355968954F5
99F5246AA01AA7CE313CEB3A3981935FCE8A747B7CF961F3DA031099D7BDDD1F9F1EAD704F71B47E73F6B48B65F7976E0BEE0FBCD8B38BADC541651C
4523EC376FD594CEDD3ABBBCEF6AF3D79A15639EA40F3F78645E975E65BA851120441AFDB891EF11863CDE957E269DA1642FC2F4B555EF8D123F18E2
BDCD2041CB552F59DD2BB70AF0D1CED2D272879E8BBB94C2A18D054D2850D7F3340C4B1B45BF5E8D999A7CCA4A75ADC7FF48930560591C23182EE35A
B5E15DB04151A3431AB0E23DCA5997250C63EEFF35A36D0700FF6CFF95E5DE5B18A0AA4E4F2D2F5762E42F03806924893377D6033CEC5619666816F0
A559013F8508C64A17ADCADA6C3D359DCB3F9E1276F6FDD11F55FDD628C8B7EDDA4D887F7132688BF5E2C2ADE3FAD82CFC3CD3CB66E928AB7E16C7A4
510C47B03F3E0B4776EF36276672C779E7178C385BFDD7E017D38677B8C7BCBFCA5455153D4AD4FADD91A6AF8B740FA168917B20DC58F43D0129FD82
647AD4D1A0CC3350C67FF4408A5CEFD4A1C3C08F9CDAE38E9CAD0B8ED11088BA304F69AAF7739C26554220B021EC781261AAF89324652A6A7024C950
F19E4290F21DC7440A029692DE913FFF78AF64CDFC37C3CCFF56FE2B53532EBE486F940745E7F30B15A4A95F0E10B04C98BA6EC21704CFBBEE2737E8
AA2A9B14693FFD451A0E3AB7CD6BEDC089C55CD8A6CCCFAB9FAD2FFCD06E2382EF6DD77E61D1FAE9AE53BA3CFFBBCFC19E5D1E3A2FD9F9E73A8BEE16
670CB1104B133E8F9206B7EB70C9A5E3A08D7BFAADAA3ED5F9EACA3E16CB0C7567258DC54D6E71AAC00BF5D59B3EB9791354D5A9AF704EE2AB14B6F2
1D5AFC51CB1289AF0B6BCACEBAAB945B3AB4B358DEC4165D4B02885F48410D823737CA3040D3C667E0EA0B310C21CAEEC5192301937E014633C0910A
9CC33F3F52320521085A2723101E24BF153F76F61D05679EFE33A36D0700FF540059FDAB3BFE65E96179D94F6EC466C7A56B8DEC275044DB525BE83A
6055162EF779E5273234662A30343B3155CAB1FEC72F6C9B3ABF827CB0B1E4D6BE6B2EAA1DED3ED0E251ED3BDD78D5EFF4D13E0B82E71F18DBFDE9BD
E1AF36ECB0B16EFF98C9863982FAF1317B7467AB599DFBACDEFF9743B46FD7F5276D7A74CFA5DE2693197CAF8F9ABC3349FA5BF73D3D1190B23B40CF
8BF8900FEABFA0B59F08AEC6DDBB46177BE4995AB0B87D87F67B31F0FA6F1555F4A55014942BC750C2E4E64D5466F5284D204DC9F546AF4FB75A0092
A1A1023EC7667C96B070BA82918831BA38087F54766393E35D9DD9FB9B611600BF3A0040FDB38BDF6353D4049276E76663F5D77C842030D8A01154F2
A2A6F53F2AC2AB7D5FA661FAF26C21C9E4F916B24CD2B9C59B1DE6956A8F6E2A3B7E798F5BC960BB6AF6A75507FBE82D3DF72EB53A76B4EF36DB85AF
471EF076ECDBA3FD1B32026731D8EB73895D3B1BABFE7316EFE879A960D210171B3B8B3B64B4279C50F3C65D947F2210F75BEDFE484BC5AF0BC17871
51A584F093ACE039CEC57A7C10B77CBFFE439232BA7D3B8B2DB8E8CC03024D0E90D587F855009AA149A3FB675840E10C05F482629951BB185FC54912
C728B8385CC2E407A938B6BC08570BB4B420AAE5D505C18A7AE768CC9CFE33C38CDFDB3FA9969F5F42DF7F2F212850FAF4ABF0FD779E5CA2D669654D
856585E77BCCCCC10D154F3309ACB9B91696043FAFAB9229DF3ACC1FE954C8DF723565FB9B9DA1FE5DE619D8D3969657F2876DBCDCFBCF67D34E8DEA
FBC1715AEA8E01765D8299088263D557838443BB751EB27BF41EC785093BEDB7AEE869394B27BA5027FA19793AB5F9CA4D7DC9F20B6FA454DA5FBE58
A9E78310A8E2A1A0F1A912FB78BBB45113F17742DC9A011616F6171571077E60CA4F11E2FAE03219626A5A0246F74F636A826291468906023868ED61
A0309AD40BEB159AD8500D47F37D641CA2A4C842917ABFE0C17BF5CA62F3F24F33CCECFF5D0560681AAA2C4AF677F593321C15FF2AE5DBDBB8F08CA2
EA06515D796699CFB891B118D2129402EB19ED8F0AB2499610CF72315327F7595293B8E2E1DBBDB756647B58ED6735132CEDD3DC7AAC5EDCEFC8ECA5
233B2C3E64F5E5A195AD5510EBA7A159E5418F78EB7636436CE62CB3F57A3868DC2A07AB5E95E8C560FDF31FE7DFD4DF3C2E166F73892806A9CB9F23
19577DD2C485CFF2AB1F89944FBC735B9A1F9E13FEE8DFD1C8FE74D58BBF8BC822EF74B82CAA0CA24892A54C97BE2992D0938C24AF1198067D49E36B
B4D1FBA36AB18842B2E3152C815637D3689381E115636EE1C2CDD20C473E63E6BF1966FC1305D02C569F939D70FF4C2EC571DA1FB96585F919217ECF
BFE709F8D999D12B4604629022B319CA1717C529300ECB6A065CF2D82943EE28DECC7CF7E2C1FE5B828BDDDD395FEBEE2704B3A6DF5F30ADF7B0795D
EDAFF67F9231A09F7D02178CD00CECF1F59B6D179B2163978EDFE93E77FAE4F13696FE64F02EC3C7701FE792D0DDE994DBF6B2723475E16D6DC575EF
2295F26D2EFAAE4A7BD3ABAA4C72EE9CDAE0D4AEEBF6323AFFA43FA1F5F7A8401AC30A11DC8073DC2FFE530C2CCA896B6169DAD4C10400C59098A2B4
4A8FAA728A0147A084D1CCA9C51C556128F6A2DF841171AB15E6E61F33CCE0D8FF1A03305497A467BCD87F5FCC71487502CFC82D4DFAE3ED17636A32
5352CEF57D4510B26C02CA65245F53D45C431442537EE346DA85B56C1EE5F7E1F6C68FCABD3DBE51CBBA588785CF715B33DA66E0212B8B254357544E
EE62972FF7E2084277EB4B408FFFA7E38673A3E6DD75F873B98D5D076792B7AE30E259D5B657316B0381EF8E123192BFE0BAA6FCB27BB106F6F511F8
15C18F9FA646571FBB44A31B3A77F502C4D3AD89407C786D3159F085874090069806FD5A537DFF87BDAF8EAA2AEFDE5F3F15A55B041111C542B13B50
4C8CB1C676EC1C5B47D151B01BBBB1B003314029910EE9E6D279BBEFE93EE7770ECEFB7DA7BEEBFBBE7FCF791677C9BD5C17E7ACCBB3F7B3F7678729
EB6EB892E584493FB860017092967C56108D3112FE576B018CA6D10A1DAE50569D54146ED0EBAF6F82C5F49F0811BF5F03444B2BB32B9AB2CF6D7F86
724C7154782D4E9358E979FF95CFD23E65DE1E118433BA4C126DE66429F51C151587499547DCDC8617C7B88E7EFEF1E89A6FCD7E3E95B76DED76EBA6
F5F41E31AACBDAAE6DC74DE853B7AD43A7A2CFA739926EB8575DD3E3FF99FB8FF75AB4C873F1D8CE963E4A60E113EDD9D293AB2B7F3D684899148DEB
8B266D063302830B4CCAB043F569E9FA87678B53ABB65C66F493CC9D5E9215B396D530EFC72FCE83F23F55C118021304C632609596A4749FAF65C01C
CD8B0152E8FAA3088A6C7C5381E6BE6DA119935CC10B043CAF9AD5E8899715F4D13818DB7483123F7911227E6F0038A43CEE4BAC812E3D16544870C6
DCA81715204FB1C65B7E8BCE3E4E7AE9BB42C1E94A7595951C585A4F024F3EA34DF2AB1EEDD694ADB39CF3EAC9FAED75891D771BA75AF6CA5FDFBEEB
DC85360EB33B740A68F3E17E07BBACD0131CCED4DC323477F97F76F6E683A73A2E596EEF6899C0FCB2517321F9E5D0A7A747946A266DD59A14B3E635
25EF3D1D69529E5F5751F2429F7C2AA62877ED5556EFDFC12D8AF930642BAC7B38FBB88C2EFDAAA21010C06092E50CF1D11A1C284DAC678411FF0481
F2EE1FC349AAF183846C8C977124ACD6F29200C9CE0030827E1BCB353E572BAB7DE3C5ECBF08117F14002C98FBE5E63384A3F26EBFAE2658AA29E363
198CE1B8FCE9C2157BC26F0C1B95C061D578693347352811FD87424863BAECD52B2E6544B725D16F7EB9071CF66B4AEFE211F5D9D17BF4915F46AD1E
E8B2CD66F6468B4E15599F3982A93F527FA18DB5B7BBDBFAC9AE07A75AD9FFC25C1F59FC29AC7AD6C6B8BE3788C3E38B4CDAEDA3F3F2F6DD8855E607
AE4F6F8EAAAE3AFBB42665C1790AFAC17C503E76DDEF92A6E6C8DE240A2B2900786122CCF5A00CE5292520AE954A219AA11896E1B94F10384E981ABF
B510B8D4C8F26F040445A02A36210091104799F2102D75A95FB518FE8B10F107FE3334A57C712B20388FA7863CBDB89A64286D6AAC0A453142FB7292
DFF22DE3063C614D302793214638AFB1F6468EA90E08ED3C29F1CE248F9345D7563E2CECB39A3961F103B4BEBDD3AE2F2BA78FB45AE46DDFDDD6BD42
12C7C174E1CE9A1D6D2CDB3AAE99DC61C87C072B1F59BAF7B3921345C17DEE8EF633A6F43EA685AF0E8AC9DF7A2AB6B168ED8254E9F9B72D278EE667
CC5C8F49FDDBF76B827E9EF70E7A19182A23A057CFF50C8EE1044E60F2844FCD3420C9AC00485CD8F141511846200600AA2952D37A03C6520C0A60B8
09454D100622B78F21F4D75AA6413D74844AE4BF08117F74FFBC01687EF6EBADC04345FCEB7475590D4CE07529E5306222A8C65343272E5A3AF0A0B4
82E3B40629AE2F55648416D5C1C025F749C7EFF90F7D1E3A7F65EEADDEE1791EE637325D9CC61DE93C7D6BF7EDC3AD3B3B0ED04A5238822ABF043EB5
EE68EE1FD2C3E65C40473BF7FBEB37A93F7E94CFDC14E814251F32AF451F3B3454BAE2502E90F3E3D2C4AAF397E5498169DFA66F31D48F34B37AD0F8
C3FCF2EA909FD34136ED42A886658519458851AD51629022FE49A601150EFD841D3F1485A80058AB8001B91EA5189AE385026F2C8C08046191EB9AC1
A45C5C431777DA0989A77F2244FCBE0BB095FF9426E14D65FAF54BE5FC8F945F8BAA60DC945E8EA06A23CD34DE183675CEF8E36F4F9838A45E49567F
D1C685661622DACD361EBBAFB9FAA51EF79B973BB55BF6E10E3D9E8DB2EE34DBD9EDB4E3A8991DBCEC7BB7C4BEE170B2F0803AC5B94DE72313DA8E7C
EE6C63BF79FF0FDFF24ED5DD0FF8C5FD70ED0F7DE275E5130E2B576FC9D117AE59125B13185893B2ED73C2ACDDA6924166663FE66D59592C3F73EC15
A00E3F9C0A7338CE303869546A2902AE2BCFAE06499EE294D0EDC3072DF915BC32C0B1DA26846D3D1C20001887D52082A90E6652451F105D2DF3CCF6
2E2E36FF8810F14701F07DF66FD89678B4FA631ACCD1B8BEB2A50542AAA5048696D7715CE1A6C12B476DDBB6BC8231B4E83049AE3A37EFE417BC668C
75978DCB7A5E89DEB7F0F9DD4E7B0BC6DA0C1C3BD073E40CA79F960C5AE2DCBDE3384D642C4753F5BB6BC39DCD86EEB1B5B9BDCACA6627B8FDB0FA5E
9474F6CCA92BD4577D12B49A1D3F400756D46AD2166F4B293A1D542CBBFA226ECA6530BD5F3BCB2DE94F2FD4169CBA2657A7055E68E6181667280250
3481305A95AFC7508AE2237E81FFC2F8B2CA374A0AA6219912E5EF86A51996AA6BC4600845E42F52B186A8164405D37BDC53C4EA1F1122FEACFF199A
2129DD9BD1EB6294729314A1291CD49A40D864822859469692234F8EEAEFB2E150FFAF94B61C4073322A55F1ABB3E98C3E365D57787B7D4AD833E9C5
22874F77AD6DBA8FF569BFDDAB7DD0C2311E1E1633D0BC520EA3EB967DD9DDA68D4F5FF3C9171CCDBD1AC32625BF09948778F4EBF32E657B810A0DED
117F7F7193F659E0C1F482EDBB0A148F8E3D1AF908BCE1DADE7ABBECFE9E8ABAE0501953FDEB4529C7014600561B942AAD26F36376636BC58F30F58B
2478EB4556C43623080C9B008C250856C809E2D9B5048EA27579E998F27915526C242A43FCF3C4D93F2244FC89FF6CEBF44F421B3967E88AD361F366
7E25080C3269310286088A6DAECFD172E75DDCEC375EF208C754650892F736038EDD5DC57DEAE13C6371EFF571E73ACF7E3C7259C694F6BD472EF49C
1D60332EC0D3C3C7663545D0344ACBF7E59C34B31FD179E80E872E96C79B079E6FD8F3F593B7FBECF557FDCFA8D539030EBD5CADC20F06BF4C2ADCB1
B7B4E1CBDC3D6B6F29365AB5F1B8597D695B5DC6C1F738517C278DE23084208D4D0A5829D3D53F89D07E9F552418009264705965B99406118A2258FE
459AA05986459A718204CA8B3490EA63115AAF26B4F55F03AA44FA8B10F117FDDFCA7F4CFF66FFF4E15BAECF729B1181602864D2603847F084315595
E822265958BDB9E61A8B480AF440D6A322ECE37E29F7A1B7ABBFFFE8DD27D6740E3962B7E64D37BBAE3DFBDA4EEB6F367C711757CB9FE02684C199EA
233527DBD8D8BB2DEAD3C1B247D1AA808AB5CF4B7D9DD76ED9B867F81A79E3C495F74655C14BD6943C2AD9B5225F9D382BE0F8C1BC896666831F034F
56E666EF4826C017A1551C4582286D68D018F569292653BD9EA2858185C2911F4610684B4D158680A430D39426280C155AFC102DCE52404EA1023794
94A26515389E50BFF407A378FC2F42C45FF47FEBF85F92D2365C1D3921BEE2FE83DBEFF20198C0888A388C132AE62ACB68E52157F7825FFDB2B55F12
11282FAA127E774EC1BCEED77BDCCEAD533FAD1C756C519F5F6E7939FB4F705FB0D8C92B2CB893ED7AF9C50A8E641477C08796CE669D87B4EDDA616F
42F7B4D7E7D2E6DA4F7AE8B7257F77C8F35501F7A65DAC1C3DB9F179CE9DB919F51FE62E3E742B78489B2E33F39A9FFF921BBEABD0F4FE5A3CC63138
CDB106A316ADBF73AB1AE35D3E4E5242D15F2BFF495D1D84A1424D90703EC05018AF3B10D840D21C949A8FE1A0C28823695A3CFB6392E71371F78708
117FCB7F219B86E91B32EE9F0D3BFE0D315514E5681152F9EEB18AE3608E832A693A3D60E8168FD9A9C9EFB27586CC8832EAF5492D73D1B98BFFFD61
3F7E9C3931686CBBE337BB6E3EDEB3F3962EEDD6D6FA9AAD378564F0F2A1E67CDDAD764E9E636D1DCD9DDE8DDA1C7D286B92FD80832B27BFFD70F39D
F7F623C32F7EEC3551FD3129D43FB1FAA1DFA2EC935307B619F2B25213B235F3E9B2FA8A73F754B430C98BE31099AABAE8C98D5A0641480CC685C9BF
FC2563388C424661701127EC3163599C601810501B508EC3B2D38C0C5AA7C5F0921C4479CFF461460D2DD25F84883F87FF34C3F001358113804CADAA
FABC796840920686944D30C3D56456D21CA6A74D522529DB61EF3AF4745EF69D7A489713514DDFB90035FA39747EF4B4CFD398E9D3577B589DDBEA73
6979DFC37EB68EF7A7B7DF497C2AA630B6FAAAE1AE95BDDF0A5797F69BE346BF7A1733D5D973F1B905F75A4E1E3DB8ECC6945B8F46F62B2D0C7F3E33
56BB7BE486DC82096E66D38AE0FA1DBB6B6E9FAC8D3E98893104C5702CAE5237E58645568042EB0F29E4FE29DEFF6B0BCAF43822D4FDB2BF6D316548
1A83642A80A65944528833A45C8D93BA441DFA29139EBF1F140FFF4488F89BF49F904FA3081C94E465165417854D1EBCEC42064DC120EF78CB2A31AE
2809238BCB50E8544F33CF84A677F7A50D86EA7B8DF4A1602465A0C5E42ABFDE31CFBC26AD68E376CCA6DFF561C3D7DA99FB8E32FB59FF309F86E982
BDB27D1DACBDFB39B4B18E9C37E54CC64107D751FB66054B52161F39B1DBED74C1F02E77B2CEBE9FFA097BD07F4163D6386F8B313550F6CC49F137F6
D65FDD5F4523284C738CB6F443664E5ABE81661094A021234EC204CD2A13BF962930102669AA75838950068403B2161D89B30C54DB4C72944C465268
5A0D5414A47ED6E92D25D25F84883FF39F6678E78FE9A545D9B1C13F4E597AFED0D38BEBCEDC9AB44BCE710CC1714499122E8DAD419A0AF4C823F74E
E39ED485BE2DFF4615DCC7B53F3EA4637CDC6FCCEFB0A0E044DF9085ED26FA996FDC39687D4F5BB7619681C8835C9AA46B4E1A7EB174B168EBD27645
CEC03319973DBA0D9BBF66C137E5BA5F24471C5697CC73DE9479FFC1C12FEAD7A3D6E467CFE96C31AE08891B38BDE0CEF9D28D7B642466C04906AFC9
2CCB298DC9A538C14A61044AA01849138D294D90BC520F63C21C004688FD79FE637A855E5001B042CF5FBD464DB354890485CF3E07268F9289EE5F84
883FF09F69F5FEC6E6F488C31BA604CC18EA39B2F7609F813D068EDC7B7EB4CF033E8AE6288E956999BA771FF4DFAA114368EFF6FDAE66BDAC48CEE7
E29E71090312D8376EDD86BB58EFAF98E67BD0ABE3F1AE7607E64FEEE360DDD36C9BFC76190C5215071AF7B7B1B3B331B38DDE32F1CD97005BDF5D13
7ADD00EE2DAE0CF71C597CD27A46F4EB390B138D91F3A7D7C7F5EB6131B55E73B8EBF2FCC85DEFD69DD1203A034A2150ADA452234F7B56420BF5BF34
66827123411292579514A0D0E1388591ADFBCB5A9318D226844229066C31301CDD2C6539BA380F4433F7684BBA9D15B37F2244FCD1FD0B71339173F7
C79E637F8DCAAB6D682CAACA4B7D7864C57CDF5E3DA77ADAF6BBDECC7B519AABAF23A1A8FBA91FF241D34D075BE72D61D5D09B06F6FE57FAD0A46A7C
AFC5B0D95D9D9FDE361BE2E7B87DAFF3FCB9DB07DA3A75B5DCDB729917E7A4EA45DD4D3B3B47A7368B631DCF27AE70F69ABBBFD78C9AECC1AF4B8776
BAF7CED6E77DDCB2E98FE51FF7CDCE6C19EF66B35EAF0C68FFB3F2CEFC5B4F2270753D4C1114F0394E5D9D97596DE0AF84BF56790B029130A8AD4E6A
C0781920CC01127A00BEF39FACCFD292044E9AE400CBE265353447A564E0B0F1542232B7678DE8FE4588F893FFE768C595D193F7BD540895B12489A0
088EE108D85897911411E4636DBF3E9B8F008CD25200CD7C1B7FFA1B623C606EDE7FF90193310C04CE24C9672F57560C711C31D9657ADC90B6333CEC
8EF90F1EF0CBE0F64E366D76366DAFC101521DAFBF66666BDBCEF38DC78017F7C6BAFF306390E72D95FFB49A79363BDE78F7BA91306A5266C4EB90F1
6FE0859D9C7653B2D1D6D71AEF6ECB6CAE522864BC8C971B529F667E4B2D6B16AE993700D25C158541F559C54A1D4B1204CFFED6ADA53CFD113D4A82
D55A0A810993DC4471446139C9E2B1CF218A7C731D7FD97E3B2E1EFE8B10F11BFDB9DFD4BFF2F2A4C98F6A6055E987D71FCA34C6C6DA2A9942CB73A6
F56D50D21C6BABD52D2C855565E8B0DA9A98902C936CB2F99EE8F1BBF0D27754E13559F9EC6D8A2FFDEC66F8F5B97DDF75CA14BB8015A3FB8D5D6CD3
D1BE4D90E9601D0E1399070D61E6CEB6ED2E6F6D7B217E74C7597B86740D287FE6F36A9FF5ECEB53C7BE7D32774E6ACAD688B1F3C1EB3636EB49D570
974FAA8397D544694AB38962F12634E763554E5229C1B5AEEDD09555990863557272A98221098AA178C52FB87EDE04A02D75261C6D5DFF831B11FE25
49196FDB0AC34C1C5B7C5E878C76A910E92F42C477FAFF36FF9B439F8D1D195A59FFE6E0FC192B8E9E7D9AAD84614D736159B5ACC960A050A1F6277E
6087CEA11089E6DD6BC29AE531A74BB0AF7D3F308F3A9FC41263A9B447B2C449876A823BD9F83A0C4C9A6F31D9D776E1740F97DDB6ED1DDB046903CB
10008F5FA20D6D63DD66CA1BB321E7AEB83AAE5B60DBE574EEA0D3E1F6B3DF0C1C59F3C567517ED2A2E4695D73923A9B4F271A7DBBE7956E7D819A24
59048D30445269F9E54F99E94ABA95FE6C69B214229A7312245206156AFF8519E0B4D0BB401210A8369224410A0541104250547D156F1CE01A2D6FC4
9E48B887ED776122FF4588F81DFF19167D38EE6E56F8CE9FD61CBB1F9F5D2B297817ADC3600444A52D2D4D693972CE84229C3E6466CF412F41F4E5BA
12A01A4EB9508E6FFF19818FF439A19494698AD36B13D61E58EEE236A47FA7431FDCC7AFEAF1E3BAE16E3BA79A75B4BC09DD91A2105170017B696E67
F7D0BBEBAAC7F32CBB8FED6F3F31EB9739C91386A6EFF5789FD57F42E5A711AF82DD3ED40D31EB27CDE932B1BE78432CA62F68C0C96A23F3F151DDDB
C4CC325010FE1C07C49790C6A49CD22603FF94F7F4F4F7457FC2CC7F12372138C93B7F610C20865324A994102C43C2828C913673E0E0AE8562F24F84
887F47FF0CC382C9130EE5EFFAE94A6A6979E4DDA0903B776E3C7F995C2597F2848339A2F4B38C865A343881962C755C2DC5221765D69721D13B4A24
B383F2B2774CBE056ACA53B393922F5DBFDDC7DCCEDBDDF5E336F399533D564DB2F10E71B6B00C6D79A8824C44C10DFC81595BBFC3167D7E4A767259
D0D1BAEFCD234332D7797FDC6B73BD76D8A0F477D3C3F7D946C8A79AF5288EEDBC092F3A584A54A58384294D09868796BD8AAD902074EBF95EE30729
472445D60913BC29A10008E779CE9042DD328DC12881F2FC47099214BEA794052049E3CD1A5E384030C71DB738888AD93F1122FE6D0018C6103171F8
E5986C835119757CEFE55B4F9E8447C77FC82C312200400AA532340AC1066DBD568763EF7AF8A722A967139F97112FB71765F8FBCED97624F832467F
4996461CB8AFBEEDE2D46F68C7F9498306058E0898E2657BDAAF8D7D787304C2FBFFAF5BF49BCC6C7EEC3E766CCC06F395EB6D3BFEF2757C486ADFEB
916E7B943F79C5254C0A3BEF1A5EE6D7B64BE23BDF7B5CE905195E91ADC7A409CD60D2CDC42F79528350CF4F33BCECD77386C26694C1585C2A415896
C285BA1FC13408CF791BF05B4F20C24B00A05A4FA2F28C2C1DC72206864BE9D8B75CFCCC4588F8BDFC6F39BFFC443300ABBF6E5F70ACA419D0012D0D
7905959284F8F41629EF377987C9E2040EC8642D0A004B9B3DEB70DE83F8BC135554E821ECA98DADA7FFEDA38F18D3EDD0EADB13EEEB8F7474EDE561
7DFE51EFD5E37C97F4B01A3CA49DD5CB8A4746CC44C4EF0496B5B5F3F5F4BA1C61E97FD2C97AEEC179678B7C777DE9B656B9AF5354DEACA8A5E6A945
7DDBD8853E9F99C2D5DD43E88C060C2B8DABC7A30E45BCC9C140002729CAA0A852A37449BA8C136AFC9AB315BCEA67708C128A83B59595798DC2FA0F
0415B697A3184E00121D8E4ACB257A0469D0109CDCC7E62826BA7F1122FEC7FD336CF38A3D8D58CB99C3CB662E5A71F0F9A963173F7CAE901BEB13EF
5E3B1D72FD531EC3B12845A210066BE5A0D268ACA98E7E74646B417188847C1201AEB4EEEA3820E5E42534E9587CCEF2E5BAF83EDE5E3E9D7A87F4F7
F61BB072B49B9D8B835B9422438F02C4D7135470FB6E637A2FFF34DCE6E8AA366EEBBD2CC22EF74C58D3FFED2EE7238DAB4E9FF1F99AE1DEC66A65E8
9A3C4EF55E4F16EB0965E857A4F6E5B9D8A4328C17FA0465284B6DD043E9973E9A68DEFBB32480B6B62CE0046F0B58FDCB6F001F0810142934059108
8651B84EA9AF97E9618C065BD43CED67DBF95588D1BF0811BF93FFFA9F8E62E0ED519DC69C3DF9432FF7D143270E9ABCFCDAE1D9CB9311DEE7CB8BDE
2656F3CE14C4600806E57203ACD7E038ACB9BEA0ACF94AB1F1E137E90F36EE9D17E6063F315C0F4CBC3921306B583B37DF2E4EF307775CE7D377A147
07D73656B150AD10FFC71E2077B5EFE230F4F38F96E38E395A0E9FDA6EE4DB9E81BB5CEFFF6CB7BA64F6B0C74B2445DDDA596C7ABC2F8F53449498AA
59267BD35B43E9D93BA905712A122519A221AF5E4F20252FD20996226961E2274BD3FC83C085833F4909254C02A0283EEAA7300227601C81E0C6F852
18011100E485CC119B1E6188C87F1122FE9DFD23EE2F418A17771CBC68C5F8E113AEBCFD96575A52569AF37CEFAC79EB6E16502CC3E89F3673248A43
A0DE640464D520A483101C8EDD942ABD5E093F2DAF1AEFE6E3FEAA29AC30E6F4B7ACD33FC4DE74F61EECDD7DD28FDE2BFCBAAD9AE1D8CDD6F945EE47
03029149C7991DF6BD5C82839CBA1F5B6B6D3F6794CDD943FDA3A62C7B3070F29BFD5EA117363FEB6B69B54176FF3D057F4D0341168E3B7AA6F1E1F3
6FCAF7712A21C18F5794C114AE686A36B666FA844E25EAFBA13F6F0C280234D1C2CE6F611920CDD0B81001E0180CA82432A91183619AE2121C5D0E36
D3E2E19F0811FFCEFD958D8F8DECEF32FEC0F465178BA5909061635896460CD2FC9029CEEB926896555E49222100408C069D515D5600A0B0CE8090B1
0BC31B4F27182312F326BA0FF40F7F19ABBDB5B5E9C3D0454D471CBA0DE9D8E9E7A15E3FF5F6FCA57F3B1BBBD7AFCFC3A0898ADF4F2CB0F55A7FAAB3
CB94AB9D1D2776315BB9BDC39DED5E91B33C6E06BAED6F59D0C3DBD636B832E409822786030497B1212EFEC8F92F86DCE4582941A0345A21C129B434
5BCE87F504C9EB7F9EED3CD999EF353F2A1067309C24849D20AD1301190C24486D7E69BDDCA0C5202D88734D7DEDE79420042D4EFE1321E2B7C25F06
3E3DEB925BAFE0D34B6E4871E1A5D607EF235986429AEEFFE0BEFE2BC635EC3D55D65452AB01153283364D86C37A0C4299D4F98F328FE401CFDFBF5D
ECEB3C63EB5575E58E7BD9C17E1FD54B3C27FB7499B7C475FD1C8B5F77B7B376CBCA0A4378EF9C7C820C70E872767A679F33339C87CC35F738DE6B76
4CCFBDC993F7C42D9CDD10EF33C4C139A836F8054697BCD471BA175B2232AFA4A8A213234BD57A9901C264D51083157FD5E28450D4C347F93CC7056B
25D4FE912D1ADEEBE3048E0BEE9FA4698EC6E5524093FBB55C0163186632A28C29C071603802631A444C008A10C17E77FF25FD8F4FEDF3F4C5926482
133A80D8DFE50559066F09DFD279F9178E6ED403F2D2B75135D2E2FA9AD85A14D01B5092CD5875F3FD9AB748D499B35B1DBA078E0932A4AC7D1F3A75
945CE2D7B787E3B080BE4336B51F71A04D07C7C4C29B6A0CA0638EE0231DBA0CF5769E76CEC9797A67B71F6D7CD246FB270CEFFFE1678FC7756387F6
B7BE95B12A01458A5EA8B8D2AB2F2EC52621E89B83E92AA5123221240151080A54EA28042170088130B275CB2FC52B0086416404CBAB021C865A4B7F
F8BBC09A1A8C8ADAF20623062108D06220A09FACFB5D961B3163AE51F4FF22447CE73F7E69E2B39E779EEEA96104C2B37F16081C257B38C579979E21
3104AFBD112657B6E8722FD7F0E10089B25CE9F9E843F353D0CCCDFB560CBAB969FC0BD5C3AB21D7FB5C664E39F97A0DF39A1CB0C4B7C79876761D93
73AF2328C42484B0E35DBABAF4F0BDE8EFD8BBBB7BC0D23E9FC21CAF868D3CF562D212E5DBCEBD3B5F2E9A1B831BB3C2CBB9F2E0C447BBEB39EAE48E
6AB514117AFD85399F0802A1388AE1A0B6590DF14A9F21109C97FF188081B0300390C0603E0420698AE5A0461DC1EA6B95301F3920A85E8A811B2D47
5D924226383942ECFE1121E2B7CA3FF9F44BABA77C3DAA155C3FFB6703D06A21E0825DDE3D2378D50D62BAB75FABC01663C6894C9A408D24C7A5BD2D
5C37BB928E1977A8EF98E8A9A3BEC93F6F3C7CDCAB2ED2AAA77B7FF37EC13E5E43BB9ADB5BC5C4FCA2379AD8F003D8044BF7EE36F34FD9790DB49836
CAEE53B4F7E6A01E8F23FC5C9FC987BB74DE9F39EB2B2E4D8F9352C9BF4AAE0D2EE7D4CB179518955A9224281A252845A55187C308084BAB95304AA2
18851A7196A26199112529A235F6472112E3C37F8DC2C471BAA23A8040508C008A94E841BB259F9B61902E3F9CC3FE210120DA0211FF58FED374C2F6
8A61C7232A5A9DFFF7D8FF4F6018B4FA4D2FDB79253806A0A8FC63AA1C80D2F67F45394C0DB01517CBAA7FDC0A33571606F5BDB4D17B73A5FCEEF0CF
EB26170CB1EA3CBA8FCFAE3EEDE7F8997772CB4DB886182036E6AA7E8253776F8F93E31C87F5EBFF83E58AAFDD077E1BE5973AADF36653B0FDB01FEF
2FFD424AB3721AC8A29F33A3DCEF734D13A7B6004A82403400CE6B80C6B7D528004388C9A435E0C2C46F9244018C66688620855980AD73C03161FB37
8543FC6DD4C5370A4D003889E717628FBD56E4E841989005BE80D83F6EFE120D80887F28FF190ADAF736BAD7EE38FA5FEEFF2F2100C7D1ACBE2A71AC
99DD6605011971B8FE4D058C159D3DDDCC21551C57755C993DEA24D93475A76F9F4BC31CEE352BFC9795041C7AE2DC75DC28FB9EE3DB8D9DDAD6DAA5
FCD9550204A98242D4DFCECB7AFA19DBEE5DBB4EB3989432C5EBD81EB7A3BB3DFB37BC77EC3B6DEEE678A62AA93C5F567322E251AFA35CD1F0F54A48
A5C161831CC049ACE171260C1B94CD3A3E04A0C856F6F3D4A61996BFF6D6C91F42FD3F414010C50815CB60FADD628644099C268A93A06F43977ED298
4C74DDBE3365CC1FE21C56E4BF887F2EFFEB76D60475B98033BF1F02F8973C01A16AF836C1A9837B0C1F03408429331E265BEEEE49E24A537903104B
BFE97E83D86733CEE3D05ED731C9A667DDDED6EC7EB3D1BB9BABB9ED6C8F8029ED1D3CEB6E3E65088C7B1F82CC73ECE674789B6D4FDF4183263CDCE5
B539A4D3AC5BDDFBBDC8E8D7C7B7CFD902DC18D9545CD3B47FDDB1E32734EFFCEF519CD6C090B03093902A0FCB2360634154354450BC8B1736FE180D
144E09318C70E24F92142B0C024784C30196AE8F7C50440ADD422C525180144D5EFAAADE64809ACFC53C6FF9E33DB2623250C43F35FFC790AF1E1A97
1D6A62FF9EFFAD1181702AA06B927EEED5778CF91198F7E33891FB0EA1E098A571FA75651C97F399DCD7F1639ADDE0119D8E4E330F92D4AD9D122FFB
FC63773B0B0F8B7E83BDC7B53377AB88AEA4018C8B580E4CB3B673BFDABBA3D7108F6E71BFB8ECDAD3C7F3E8CCB6DBF37D5C7A2D7E54ACD6BF6AAA79
AF7A39A0CFAB822BDB17E57398CA8421422D2F0B243EA9864865FCF31298129A7A79EF8F53800A010C3CC319A1EA4F20BB90FBC7308266A0CF4F33AA
215658F1ABAE92922513A6BCACD0B440B54195CFAFC1CC1FEE93150580887FAA0560F090CBDF46A550F46FB1FF9F7200EC6FC381384853DDB2DF73DF
FC0E63CB491C44A9DA7B6A9A4B5F1A79D6B79A636F25B78CEBFA71A75D80C74FC7EC0764557FFC75D306E0A85337D7A176B67EDEBE16D69E854A8285
49AEF00C3AC7C161EC515BEF419D7BEE3FE331FBC8740FFFF0D18BA227D977FA45A37B6D48CF327CCA8C0DF03EA0D8E8BC50CEA1C595460EE37881AF
79FAAE916CCC789F21A385C91EC24388F201443806145AFF7945C05004CED0044ED11410FFB85AAA37F16100D65C07B34563863E9034AB35C5CB930A
C6E4FDA6757E0B6E44F72FE21F1C0140C7CFBD0D96537F70FFFFF68CFF33149C3135AA0ADD3A87AD30733A8E9008C294DEA968E25237A61C59F81451
8743F99D277D72731B66BD6553EFE0F8579F6FDAEED34F72EDE5DAA9EF3A67AB8EED5D53544682C3B9B26BE86C4B8B05FE663DC7B4EF76BB97FF831E
C37AAFF51B55B4DDCC71B39E7D9651F7DCF8F1BD21C87A7BED8F562BD45C4D4E356630503C4BEB43936924E5518A94E2DD3EEFFA319CE009CF7F1142
A51FDEDAEB2B34FEA0B450155C55D45859A36E51E0348BC9D514973FACFB95427953E3E7A1478DE35736212CCBFCCDDD8A10F1CF8B0074BF9C5B9DF0
673EFC81FF8C5009C491448BE266A761D13F75301B5A41F1B17CC3CB1B1F88FC9BE9191B0EC38A38E86BAF3DBFB8CC761D9B7D7CEBF51D8F1316747C
983D78887BA73EAB07D9B8B7EF9A6E401986E43EADA3579B5B6DF2B6F01967EBBF68E4B959A3968D3E303245EA6439A1854BBF227FA62F7DD69CDE6F
A164056F65D0A22280E3D4C2D87159741347C53C6BE43D3B2684FE94C0FDD6F93E42AF1FF93D1940F0FA9F60688EF8FCAA09C69A6B5530C5F241011F
9F0CF13C5BDDD2927CDF67BFF4A7806B72A275399830ED9813F92FE29FEDFF754777CF6A6119EEEF0C80E027694A2805E6F9823132C931B785C51BED
6D9C1E510C823564EEBC4D361C2805E66EC79531C0933EA7C6771C3E33EDE68ED4356B2FE6FA7A35ADB3EA6ED5D67BA5A58D45F7EC320307115CDC06
7C415BFB052E76DE1EFD170E0E1AD9E7C89095BD369A0E5B385CE18CC7133EC8CAEE2934FEF30B560C8B01EBA2640C47E8858DA32569464E1B116360
18A1DE5FD8F3298CFA13DC3FDE4A7E5E001000C8AB029ABF8FAADB6F9514ACA96E42098C12CE013E79F7BE55559310BC6BF15AF4DE812DD7750CDDBA
E880E128F47B2F80680044FC63F97F6FD72EE3FFC67F4ED807F23D05C0E12C959EF962CDD68AEDFD3A9B6FD493D51570FEB98B0AEDB5E2448F034876
9271DB84F563572E5F386972D2BB0D5B225FBB1FBB65EBEE606F37DCC5C17AB02CA98A6309AEFC11B1B8BD8B974DB7899DA61FD83477C2C14133567B
65353B38595DD02445A4D448AF14824BC7A4CD9A23C1BEDC933138235418C11F62692CF7713AC16242A31F436018AF0364CA56F74F7D77FD24A133F1
FFB078C5D3C7192804C3CD4D08452038CDD1275C863D28BE7F6DFEB6373F299F47868DFD0809CBC1F8FB624D8A0A1542B162FE4FC43F150C274F3E71
0E64FF06BFBD8161A9EFF3015984A90DFD5C7EFA42CC4FA7E65BCEAC07F3D3F09690634DC8E994DB1D4F6259EF0AE6F7BE746B607B77D7E965472FEC
7CBDCB337B645B17E7B69D7B59B41B60BC94C1113857F2009BDCC6DDC2A2C728CF99C193065E9FD97F4BBB65F091F696630B4A2FBCAFCD0A2E322E1E
FD79EAEC16D3E7BBB518A230311CA7F89085E91E3FAA662902344024CD4B7ADE00A0CD8AD6F8BF95FD3C48142718561A1E5F08832004697508416028
C9554CB19AF4F2DEB26D5B0F4A6FDCCF3ED034C23F072169FE9E68BC2125AB0AC669D1FD8BF807C7FFDAEC0D6FC9BF70FFDF9C2078AAD09C9020A029
203AF07EFD950BCF77551D1893C72155854645C4FE74E9E5E8F58393C198A4C859775ED8B61BE0D165C7DD13B777662D7971CECCB34BC7FE7DDBDB8E
C3CE7CE3488AFB78949868DEC9DCCAA7F7FCA32B07AF5EE9366FA9C7F97807C7DEAFC8CC078DAA83EFF4ABFCB2570EAFC6226E29D0DA7A154D335549
4AAEE0C66703C5A2B4510BB7F6FB9334295CCF6FE13FFE3D0EA04920FB53356532E8747A18E65FE5C304FC6627BB1FEE5E1EBB7CCF15A8708DE244CD
69AFED552843F1DA5F1A1EF85981529418FF8BF8E7D29F618DEF97E7907FE3FCFF450A02A1484A58A6CDB024A97F7E6A43725A5AF0C9F8656B1A392C
5F52A78BDF73F37D68D4E8A195E89B1B1F2FBF9965D9D5CBDAF7CD84B3BF063E5E9038A6BD73FB2973DB990DD69F4C61418A7B761C1BD0D6CAD6C669
D8D1037DD6FE64B328C076DA3C779BDE29A69890E2865F9F19E68C2F0D724A5586BF56430A0964A2A1F878068F7E500AB24680445B27FAF2713F86E3
28224CFA1178CF339D7F4EE23A4D6572036C52CA654690E6B53F46902DF3CC07ACBBB06AD4FA8817906C634BC269C3349FD03A98A158A2E14A7034FE
2F5B27FE2188F88786FF2CF131B085F94BF8FF5D1BB40600881EE01D2627A4CC485A1515B831A7AEFEF8C20383B63553D29AFC5C63CC82BD9793AFF5
5EDE64FCF868C1A657CB7AF6F3E8FB22F0C8E79ED736878798BB5BF82EB4341FA2BF5FC19958EEC3237C60073B5BC7EEEB370D99B279E884CD3D7CF6
BB3876BC0AC8AEC5E3E74F356C18F5EDB8DD99E2A76F9444BD1AC290E65749903AE4BA9C35E801922284DDC4246ED4C03461501871026BCD0050388C
51242A8519D2D4AC52692196A1501CC5F1773D1CFCB7CC9BB4EA9D2643A5DC99A3DDA46C9EB33D4B86D0389A76F25431C389634044FCB3D53F2F0032
6EB4D07FE6FFEF1362AC49A921685E02500C8E33CD85C18B9328E59639A7FC7EAA44AB5BCA6334EF178C39557EC66E85461639DDFD65FE64DB5E2E9B
A256BE75DC7C6BE6A7516D9DCC7D2DDAF5D6BCABE470987B750FEADDC6C1C272E15E57C7053E6397F5B5BB1FEEDA7EB15C1A72D77065A7FAD71E1F0F
DB6D97E5E46286C6668506573FFF4C5604470080D180920C4EE342C48FE88CA4700E886108D46A007012C548B805E1ED93496984699AC40465A039D4
B3DF8FFE0357C7C8D988F0A693A96C5032B37FC9CBF2661C90C41C0F6FFE5EFC23FE1188F8470B00AEE16D03C372DCFF5E0FC3404D4A90114EDCF8D0
9B50166D5C25611B17ADDA35F48054534BE4C5196226773EA53CE9BAAD5E76D377CAB7F05EFDDC96685605CD1BFD7E71D8AF1DBAD8DB77ECD0579716
4E700417F50A1BD2CECEBAEBC931363DC7F7DB30CE611B34D7C23D1579BF4719B742F3C57DDFB961E7940A1D421AAA6108917DFA0AA65FCCA1181D40
B6E62085EA7E82C084563FFE6A300CE595BF300788A14950CE1B0204C36882C2789D42E377E78E5A386AD4867884831E87B6DC4AE79EDEE6D2A6DECE
A9526A1383EEA7A9C4A85F84E8FFF92F7D54C1DF57FFFCFB5D9426CF8892B8306083D6483F2D5DACE30AA7EFDDE47B166C6EC4922365717DACCF194E
985DD455AF36DFD830DDA6AB4FFEB6C1314E1BAECE7A38A0BD6B1B7BF31ED2DC53088B73A17790FE66B66D972CB6F41EE030EEE7CE56D912E7768168
CCCF6505E31BD47D165C1B76B249FF45A954579660B83C21D61875B18683742427F4F7F1B10A46B41EF7B7A6FD84E55E282CACF823412D6AAC2BA823
29DE422004CB6178E9D1E90B268ED99F03D240D18B32F4693CF7EE57025ABA3929A54EFE72C9D3B27A710088081182FF2792BF91FFAEFDFF5B03C0B2
4DB10A4470BB2489D5C53E9C721061E3169C5A343C5C5EA407A2AE655E74B0386838E4FE188CECE57E6EA7ADA777DC5BE788E5DDDFCDDDB2AEB78B95
93792F99349E43092EE22536B2AD658FFD768E63FBF79A36C3C775F7548B7EAAF2DD09C0CCB7E4869EC7E68CBA595922431034093056A54A0A5F9E6D
A2299867B4D09DCCD0BC7727D0D6633FAA95FF288CE2088A83E57992D2B402B075FB274D19B48DD2D413A37AAEC9C10843455C96818929E0A2179BB8
D0814F328B6A1E071CAA5750A2F71721A275EF6741A491F913FFFF1417F3D42BFC62E2836A1A8771EDE7F7A1532FD1D4CBC5214B66C5A67CD2809F2E
246FB036FF153CB5304676CCA6D7F6F9DD9C36347A8DFF342C68778FBDABCCADECDA7BD457C77038C646BDC0C7B56933757EDB2EA35C470C9ED7B1FF
4057F778C3DA77C0DAB9F833DB45F3463FD17F09D737A3A955CAB8DAC2FC13E7AA29A0454F30144E731CACC75A8FFB5061FAA7A001300403618C40A1
DA9ADCCC5A8CA6081C6308439904A492D7FFF409251A0B320BA41414FF8D2B5AD8CCD50CF9F17D6AE5B99F7E78676445F72F4284A0FF194E11AEFC73
FDDF9F13632C8BE516002481E37CCC5DF4E2D6D74BF114FC66DDB239F7F39EE6A89A23A3BF8E71B4D8652C0FBA913FCB6ADA9E65FE43A517CCAEEC1C
103C2DE0C88076F6D6BD548977399662CE1D847DDAB6EDD2D96EF8E099FDA71E725C31C72D8CBEB40D7CDEA94ED579F0CED529B0E4753359A66AD045
156B22021F481124A998A049866271658D9AA64808C0BFEFF615967FA0308461FAE2A6C6BA3A2DC2081BBF2842A791D6BD3970B10425E5AF1264721D
08964A68D3BE720E5FDBF5746CC2E55BAF7E9589DE5F84887FE5FCF0944284F96B01D01F820486D1C724E8115E6D83B83AEA641A7BAA1C6B4CD83B62
6E69D92B1DD470EDF3C7B18EEDD703F5730EBE1BE23E64C72BAF13AABE5EB7BAACDF67BB21D8C2C1B6477DF43B1683E92BFBE07E6D2C3B580D1DDF7D
B3C34ABF6E010E07B9C45965351E37D1796EAB16255339770D588A4E9BF8BCD674F397B74D24F0FE8389C2519AC1F2939B319C54A754F00E1F13CEFD
29A1E61FC7AB5E252B2104E17851801B104829C554F7A61C8318F2CBE93C0CAE85D0BC6F00753781E39E398FBD1B76F28DF1CC278CFB1B132742C43F
92FF0CDB90A9FEBB02A0DFA90461FB674D5249190CA31088943DBAC9E50613BAB8B45FFA6FAB7D71A3162A3DFDF5412F57CBE9920F4322C2BA76ED9D
32CFAE3CC2EACCF6804B037ADD0EB070F0ADCDAC6251980A7F804D6DE3E0E03C678CDFD0E1E79C867A4E80F5E33F209B03C023566B568713F9D78BA8
F4322AE94B9DE4EEAFC99806C8FA0CD13845D2604E8A92E03044ABC684199F18C590C2962FA32AF356224C12422D2083C100D028D3943C0DFAC472A6
D0EB3A422743E16FD91093F28CE4E43E5EFB8FAEFB809DDDDA4C8BEE5F8488FF99EE012556607F5302F89DF7ECF77731B42E412D89AB016103A07D7D
50C2FDFA90A80A2D38D0E9B4F4F6F526202E38F1667727EFFBE02E9FE41F86779C79B8CD0A7865FFD05EBB5699AFD96F6ED5BBA2B89EA230F2C5756C
761B5B33CF290ED32C670775EC645BC82DDB8B3DE815F6C67ED1D6B590243897884AA792A2D4A58F2F95193045D1272D45210461282E02380CD5EB09
96E24310FE159A64295A9B9F965BCD9B074C280CC20CB8A6EC4B4EE1FB6772962D7A998E938012967D2982C9DC4B7A96596BBB6ACFCCF7C493DE778D
CCEFCA9B4488F847F35F38576FC8D461F45F3A005B7B64FFFD8428BF96F7F25162BD516E908406B7D4EC5020A9378BF6CE4A90469F2901E3CE159EEB
DDFD1151E97BEDE5F8F99DFCDCAD926B875D9FBEF694CFA0591DED865436287012256EEFC0D6B671B3E83DD77BB2EFFAEEBE3607B813E3C0D435B3D7
779B707C52923428992A8E450A3E9A223726D743F9D509F5261A410812D6AAF8D05FD9A22384C41F81E22446511CD1F8AD5223B4FF5014C9BB7E0802
CBB232EBEB7572B8BCA4066169A50E694CAB820869B8946543AC671E98F595FEDC6759F2F7F59F22FD4588FAFF3BCD91FC1294FC4B11004343F2EF0B
C104154053FA872B74CAF4B7D9559296CC6DB3AA5EEDD6C3B1B1F5571F3701B72E01C097C0ABC7EDDD23994B8352FD47CEB21FD87E9C71DF9C2DFDAF
6CEAD4C5D96244457E0E891A88973F033FB6B1693FB4CFD04D03C775B2E94744B856A9279E7FE8EA7C66D831E3D6082AF93599725B533AF778554CFA
C574258DD3004EC0084530C6BC4A8CE078B293180E232443D647253733ADD5FF0401AAF52A449D145E09334DB9DF8A35FCBDF106C1D8D8A040307978
1345DCB69F7C694E3A1DBBEED685029C11DDBF08D1F7FFBBE68F319437FF0DFF5953B516E0FE551B40E365FB9E2A1444F69BE44F45B766CD936CBF62
D48725190A92CBB4D7AFC2C8955EC13FBB744F31CCFB78CF7BF944CF2EEDA3BE4D3EED7BF0AEB76B175B7F794A048183447A28BCD6DCD1767EDF55BC
59E8697E51D5E90A1772F6CA8EDEEB564ECD391A8AAADE9B8ABE1A25CBD7E54A1F3E8E4738CEC0E0BC5BAF2570494C23C95124DD5AE183E168E5B3E7
99268225508CE20D84418F996079713349190A24FC7FE363051CD1C955380A7E0D2E202B5677FA2162651AFB6175F6FE1BCD0C237EFA2244FC4BFF0B
43C00115CA307FF1FF8CA609FDBE138C7F50A43662666E4E24D51875E3D295E37E874B161FD565DFAD6851E7D6D4CF7B8A2A677BEC5EE430427A7E52
E34AFB2D63DDDB4FA89FB27469DF17CBACBADA8C51BEB8804230111BA29FDCC1CA69E2DCDDD61E03CCFB67F51C4FBEDE54B5D46DDCAA7E099702E1A2
5B458AD8D2C81B47AE25467F6DA0399065D1A67AA4510EE527A95AD9CFD228C9C70370FE938C1690C35BD94FE08D4D1862D411344DE194302F0827F8
874C86413A5D7D788A296DA0CD924F5B33B8E83DB25F5CBF1A19B1F45784883F2A00E2AF09409657D7356AE25F368226F096A005D5D7DF51AA840B1B
2FFEBA253769D2AEAC4FFB5291E2909AB8818958E9B0FEA736BBADAB9E7736A9C790C0EE3E7D0A4ECF3CE079E06D1757E7C98667277018203E1F578C
B5EDECD963FE0C8B81DD9D5E5FB6FB82AEBBB7DFBFDFB249FB234EC82A3E24EAC086C0310F6BB21469524EE8D2478B73618294859790FCEFE7E92F9C
05E075125D5915CA20144DD22C4F7AB252461AA406824618A14E9022190E8315268C524AF4AA16B8666D9F2D21C1995CDE05A4A9F34A152136FC8A10
F17BFFFFBF81A470A512F9ADF757D8B623597508B992416275218B0F1D3854F871C6C1924BC79FD445EFA93BEF554045BB8F3EE66FF738D8FEEA890E
CB168C72DD9B3864AF7FD7984D160EE355B7B7184180787D080C68D7D5BCEF1A5BE7EEED5600FD164BD75F7AD6C97BEA0F83C2364AD497EBB2B5CCE3
A9CBF2754FBE29488AE068A0A05241E3AAE87482C1318A86308A618C922F395A238AA20C87F3E1004B53608A0C073438861124CA1B0D8AC4B0E63A1D
0736B760D282D2C2AC8BAFBF1D2FE2AA1F19D8539E71382B0EFC1721E28F39C0EFF1FE9FC27F5E15D08CB6B609167A7F85415B3821D9155F7793C008
FDC39517C24F65676F3A54F165C661FDE733CA5FA6C8C9E3CECBCFF75F15E3D2F5F62CF745B3070464AEDAB0DBF650B88DBDBF2E2CD008A244F22D64
BE95A7E3C4FD361D9D2DB31A3CDFE48C493EEBE1EBBF72CBC338F0792A2ED3DE58BA250B7C7AAA8CE269CC4A8BEB598CA87BFCAC05404996C0601A54
989A725A005258FA2B10D998A5A309431382C2088112C25180518671A05681527A15006ABF34D67C8A5564473532D52F9BD9921E07F522F34588F8CF
14004353342BCBAC95D1C2BA0D9CC431B27E4FF9EB7B1C6CC4DF6FBB9D70343A6DF9DDAAF4BD5780ABD78CAB7693DA5976FB5739DE98D161FE2337C7
4193BB9D0D1DBEDFB3FB250FB321D27B5BF5262311775A33D5C2DE7AC3B4B61DDB4E32AD1A56B2E6E4F9DED396F6F13B1A0F86BE036B6B32662E0B6A
96BD2BA74913839715C11043A6DD29D491268426218AC1AB0A54901CA7119446495E1EA040A5AE7504200409E90518A3C062090268188A2314100364
15A43FA829C9B8558B7E7DA5C614937AE68AB9BF3F4580224401F0B7DCFFDE71CF5094E46B8D41E8B6C379FF8A510DE78A0FE7B12008A787C47E0EFA
56317B4EA5E6514A53B4D2B0399C2BF01EF164CC0F215DBB6CBBD0CD65DEE665890B37F939FA7A7518AE095BA1C35122FA24BEDCDAC9695B6F0BB7B6
2F9AEC2F452F89DBEEB4EAD4E8A985E0DD1306FCC2C35B2BCECB0D397A0E8768B2315BCF6048FAD3321CC584DD5E1487A567411485003A1287080651
359830A111986E3D05C470C26482752ACC08E01CA1A94639B6FEF9C3C47A6D5A62239A9F0013C0FCAE8721F1F31621E23F3400BCC866080289492F07
188C1056EEE008937E27656B3E63D213050732F3CE2A9F9B05C891C4B0EC34B87A7E0E17EEB43CC86ECB0F76DE49819683EE1FBBBE2870999DB57DBB
E1CA47731A30107FB2195B6EEEE4D0C7CADEACAF7CB5E78DBD1B82D64ED83C63F05353DCFC2AF4EBFAD06181F5860A9906D6324C733EC27BF4C84772
828470128318B6FE6305CA71A4C968C4548091D037A8605090FCC2EE5F02265108804CBC48D15394B6B4B881E34AAF7E9534653F4B024D79610D06E9
41DB89A5A2C71321E23FB400AC90F76331C8F42A2E1B15DAEE040140C0F732CBF6C7513C012BCE7E49391333BBED320D7A67464E0913B940CD5EEE72
7BE6A095833DD6DCB275DCF8F6F8A2AD579CED9D3A8C96E67E362120991488EDB676B135B7E9D4EE5AA665CFC709D3064E0D3AE2FB137075432A98B5
386AC392C6E63DB51CCCE1C646951E87CB6E461A6910A1698C17220DF1156876940446681832C206C06442494298E52B9C01085B415018076090C44C
95E57282D37E8C2825B336DF6D3202F1BBCE7F498E1FE31182899FB50811FF19FB7F030501B5D16F62480C23081CC7113A2348D3B22A8C864136FF7E
4BE983E7DE661BF5EAB99B8AD3D8A32B71FAD7491B6D0704CCF05C3EDDA9CBA5EBCBBCE367B67730EB535B9481A030917D92D8D7DEDAD6CEC96C64F6
689BF5095B47F41F73A25B40DEAD65356FEE5C3D367A7425B5F716AD6EA1EB4B9A41231415F24928F6A5500C26E9A67484CA7D9469204D104E995A00
1CA5519C113614D11C2B2401300434400486A86B25589D5A9595CED22F173D37A0B277A1B9B2C802C3AC619562C42B42C47F630318868115F5E59284
06B23502C0114A1DF4806EDEFA8AC551A6F11B26975CED68B516289DF0B5B491FE613741AC19D2C365EC8C459E8123472D3E15E27636C6CAD2D2B7F1
F96110C788F78BE1CD1D1CED1CED9D1E1F721E5E75DB727CC0D815FDE3BE0DCD95EECA0C5BB9B488F8BCDAC4A184AE0A30D2CAEB379A611204311AA1
4950AF361A720A4C2C63D4C00C6DACD7F2D743D22CC309EBC9F807515CA001791B63A8AF34E95BE4809AE3B01D3F1511205C1E598B29B3F114CF1398
C87F1122FE23F633FFE23F89A9B2B33EA8BFD6F0E29FF7FF10043E5A92CD682E4A380CA012EEA1B2C66BEE16BBD0B3831A3EA0E55D8E906196EEEE5E
DDAFCE99D62560DA9C302FD7D49FDB58F4948485200840444D05D6B7B3B7B7703CBE7F855350D9609B61D3BABABC968E7A645C5D583B687C3D91DDEB
06AAA3912C0989D79EBCA2214079490D8C90B4B1A8C154915F499274430140C3240161C23C42E1F89FA5855BD03CBA516C322088A2AE510DC90C3444
709295C7410C47735F975616C76832BBF52B65C4B13F2244FC7769408AC2A479F73EE34A0942E0288EC044F9AF6B9A58FD8B5A96220D4FEE28A4CAC7
132DF7D7CDBCD5F8817BDC71D72A6BF741FE5E8B3EEFE9B1F8D8B8AB0BAC7F8AE860D1AFEA4B1ECE8BF8DC1F9140332707D7ADE7C786F78BF0B71931
DACF732B736E04F2F81EB2B05B245A366AA996D0A28D35009BB3E1921A4055E91213411A73A39B80C21282A25BD27314388D0963FE847DDF34CB31BC
4FA7D1E217912A04359A545225A155E244C647A8E0701CCB19C88AA82AB9C924691C63770911DDBF08117F4BFEBF540230FFCA00D23485D57F3C9E45
33F26602877193B12671F92680938634635A027CF44D9B555C14D8E3FE2B9F92B4C7F0F67EE717DA771EDFCFF958D28AF177C62E0874B3DED4A58D57
5E42058AA264F224FD1BCF8E2EDB6EB804BF19FFB3ADFBC8616E3E8D4D3DEE7FF053DCB50F218AFAAD56D5A74AAB33003862F5D97A0C97471793941E
AD4B2C36D4EA0806C8F85028A3511823089A20710424855E1E56129F5968C2084067029510A431E186F47309F20A0D8798405EB92830FD5DD9B5B6AB
648C38F74B8488FF320E100C00587C6F7DA08193C5961A4D7A637141E4C4FD2457124DE2264C57A4D3A5471BDFCE8C39B650FEFC71E581B2DAC15D7C
FD47B83E897299F1F3D8637E969EBDDB0CAE28CE33113899EAAFCEF075ED3CA6F7809C93BD7B0C1E3EC0ABCF0B66DD34D39A7895CF5E4C3672274864
5542752D68F8EA08486ED4654B301444498D9C30C20C2C7995DC4C612846A2C226201A870092B75AFAE4D765384E9B0C461886B4321D4661D1117A4C
4AB094DE58942143F18A99E1C536BD3E60CCF75041840811FF21FF85C57F1445357DCAB93F3D9993DECE6C021B2BA20BC2FA9D6098841CC264C0240D
04F8E5361CB3EED5D6A7DA33CFE32E2B220738F61DEE302A76A1FB0EBF452BFB5A0D6DE3535A7CB59904C8FCA92DB7DB5958DB596DA89CBB29C065D0
300B5F28D5276BEF6166FE5C5CEFBF0735C52535C737A14F8E1463BA266D761D6F7A48A226B719658CE551311594B00180A2509C24189A14F678D3AA
CC842A3D69329A6010D1B6341A7134F26C3969AACA33A98C20DC40205472F7CBA6416DCFA828868044012042C47F950210040003E76670393F5E8320
795A438322FCF59743DDEF30707213AE32D417F1F176DC5963D4ACA7C72FE53F3A7B4D9EDEAD6B47DF414E332FCCDFEBBB74D31AB7AE6DA76924A71B
498C02B66BAE767076EA382A3AA8F7E34136E3C7996F978EBB9B3BB9E9FAA416E9F4C5009EFC0595D7A2F9218D5C7D1D2A53B128892171AF7434A678
F5A21AA1189412CEFA8541BF2445B314A547509A62F98080C40940ABD00354E5A563F990213FC740EA3958CDB05C6EDF47C4A1F693AB388A311AC4CE
1F1122FE6B01405355B14D9C26F043542690FB24233B2CF8E6DE01EF29F5D7A6662DDE52A4A09894FDD2AB33EFCD3E5CF72EE889725507EFDE1DFBCF
DAE3BA7B46CFB5D3BCDAB51951FF364885210C7AA2F65C3B5B1BB31FB2A6F59DE9D0DBD77C7CDDAC917583E35F387DAC183859AF8D48801B0A08D585
6CD850A91796F4C08CF4D83D03465745669978FF4DD2069822700A474892C610A55485D33CAF598E6149048421124A89066806526848B08106A42692
D3FEF49EFE64362C19275943A99216E92F42C47F2B00285A57F2262AB385A84D7C5A18FBE0D1C71DEB2F1F99938AB5945756E398A2564F732F0FC803
FBEF99F2441DB93837D7C36762D73E5B13C7F96EED7B728AA5535BCFECB0B32082B2E47DF92B7B2767EFAB41DD873A771CD3D729A9D6E9E6AAA0F79D
B6DE193CBF06781905E9CFA5E75D8AA19AD35B4CD542B55EC9E92700AD2E2DD0522CC15B225A6F124A7D2A2B719204753A9462499A126A00688A24F0
BACFAFAB288E808406210504B54838AEF95A1E13E9149081D0145E1CAEE4C4FE1F1122FE6FFA737FE802A448637DD28B9FCB098448B8939072F3E286
19934FFDEA178BC9AAE2AA49426B6861B90757243FFA4EE815AB5E3049B1D7BAEF201FCF33512316F759B1BE8D439BCEB961EB1A4994862FE63DB3B2
335F9FECDB6E9C43675FB34060C1989D53137BCEAFDE32A9A9F96438D272ECF6F5753986E62FB9AA321387692BEF44D1C66A21FEA72806C784459F04
81C30D5A1C46704A98FE2134020B6D8A3409D4BC7AAF660982A410C040E2350D2623F62DA28EB9ECB8518292385D7F251A60C4C1FF2244FC6FFCE77E
3FFEFF7FA683319411809018F4C95BC20819250F227393CFCE1F3C3578D1C878102C7EAFA428562AE188D0C8B4317D478C4C8EF6B8DD3CDC75F228DB
8197364D08F0BFD2D1A983AF34E394942158E47EE57D67579BE33B6D5D6678F7F5EA591BE7B8F2F8C3F1433E85CC4A95DE7BAF84438FC46C49C50D0F
9F22849EE33E3DFC22813431696A86F7FD2C51A3E4C98E800604E5C37F1C67384AE03DC1C70514CD542555C874546B550083A1185C518B70C8E38F28
7AC0E5701D421294FEC96EB9502A2CF25F8488FF4B00FCEE7B5E762B6520A1CC3CB1D468D2E891CAB484E28A88A56386AF1A39225DA9D14628798950
55C9E93E973DECECE73A56B2C3B532CCAAAB8F67D78187BB8C733DE96DD1CEABAC285C4663ACFE52E2654BC70E83C7798F761BDCB9EDA9BA3EDDAF05
AF7278F3657890E4DA4BADE9EAAEE70B3EA0BAA2CB596823C5151FBD5DC7963DCF40842E1FC4A4ACD1D218A66FA831C1380281BC221072810489E114
ABCABC97662285EDA01425D4290335251883BFFCFCFFD97BEBE828B2ADFF7B3D6108110816DC0930B83383BB3B833B0CEEEEEE101C82BBBB24248128
51E2EED69D7697F2734A7F5561EE1D79EEE87DFE78DFB5EA3343AFEEEA267456ADEF9673F6D91BE2BB9A1C50483B8674DC989318C77EB56AF29D9691
F97D03F0AF17FF9A10E4C84000AB8BBF5C46396D161B56529CA34FBC3EABFFB7ADD794EA48C317AB28AB425CC8B95D7A65C0A82E8B63276E36CFF2E8
F3DD988107A7F66F36BD91BB679B8CAC4B451CC391016951ADBCAB7FB7AE4B874E83EB74CCDD5179AD5FB391D3CE6F3EF1A5244BED546DBFB7FF1641
5039710461E56C679F2BED69CF733931C74714B9261402803BCA4A352845E3760C4AEB01A2E459D45C5E1EFCA19C61210160459B22BB4DE7A068A77F
046BDDD6CEB704079062CA96AF8B1103155ED6BF8CCC1FEBFFE7573F5D8019C9388A3A739F1829336236596D8AEC98C4E05B531BD5F223D5265A8FDA
05DE4A3B9EDD43FD478F6CB03366D2ABC8CE751B0F1BD97F5FA386931ABB7FD334E9F3CC2F3CC59371D9519DDCDD572FAE5CA973FD4AA7C33CBEBFE3
3DB3EB22DFC96F4A03A260E99633DB77E0245EFA51A5D1B29A1B77F3607CA45A10E37C454239C2428E2250CC6C163F41918024798E1643001C294CFC
189E85302C498AB13FA068A0523A5168D284C6D38A75DF9F2F23200E58C3D667274C82AC7F1999BFACFF9F5604C43FE68F4602A5CC119F6C9829578F
9516A3F9EF771D79BBDBABC5275AE72485D868C19CCDA02F1FA39FC7746F7C2E6FDCD2738DEAFE78BECFFCD63E6BFB56ABF16D6EDAA244F187592EBE
8FEF53A7767FEFCAB55AD4E99CF443B5275BBEBFD033EEDE94147B78191DB0E7C60EAD18A83F2B06AAB2E2FBCFA02A395D4A441C093936565AFE9346
0801404A8D3E799A17437D8EB4A9F36222B3CD242B9D4E2621A468862833D1784A4CB9960A9D3EEE4E3990360C0D5B5E3CBAC4B23CFFDB5F504646E6
3FEAFF97FDC1E9D46480DB1CB62707F458614809A12936D9620E4D7C73F3DB161982B3882F5B532E249630C6234FC195D62D5AA9C29B8D5AD7B0DDD1
A63F8CAE34674CA51A4D53D2F696B040B05E7B18D1D5B5469D2AB55BB9B92C78EF3A3CB0DEDE5E511F7B66C75F4EA61F2D39B1284760FD9F2A8B1501
978EC5A0788E060A0257945A4A8B7E1E1236BCA2050929DA01D102D014475B13E3CB35C506821193069210AF410E009B8D3607C73AA0F6F080692FB5
3447D19C61DED684A12AA9FE5FD6BF8CCC1FEBFFD757BE6601D64F5A803A70FDAD03796841889EB22B743812B07AFBF3F9030319D2C43D7CCDC3122D
2C9876DABAA04BBBEDF4B9BA0BBB79B6AE3366ADE7F2E595BDDBA617BC809CE8B2BF7C2A18EE5AA386D7EABE55BDFC565459B17A58D725E8F82BF62B
B125D14B362F8D6684A83D04967F68E89A38BCE8A351FCF74D51D908CD9180857A1364A55A6428695F6A4CC8147DCCB449B900459062F04F49B30178
924069DDE7020014F3BE5DF6D9CC42C071AAE5C314AB2EC09F279CC8375A46E6AFE9BF2205E0B89C601BB0D81DBA536710BC30A084C0509D059ACF0C
3D76A6FB46236B531ECB8416B3828CE87321B26BEBC6AFC9B9F567B7AD3DB0D66C0F9FADAE9E2DE2538FDB05966714A5451D5D3DFF67E425AFAA758E
366BD4BFE9A8A69153D6A067CE6291230E4DBECFB3F1276C69C561BD07F86349AB82384148B9A361041AD20CB4E6944BBDFE2AD40F685A9F5EA2C854
419624288220300A5038AAB36338C17089B7CA292A6694CF86704CFC9B1C9E3BE1DB94A02E893FBB7FD900C8C8FC35FD7F9D96C10B64649C8340CC84
EEE61383A5F8558205B7396D107E5A9915D276AA8167BFDC4B2801E94954F818BFB35DDB778D291930695187E92D8636AAB7C1C3AB5D72D136A5F8B3
2C2F4372BA7B7836BA35CAB3CEC07D6EDF0FAD3BEEC08A1E29898F8A3FADDE74F0344EDB37A5B07999CB7A3FC3A296DF72F27CD43387F8D718D1FDDB
8B0B6D0C0FA55A248622F56961A1D9464C9AFEC5009CA270B238B644596A4568C176FF7C54E2FDCD6DBE3B98EF60289AB2DC9FD8E02DBCB0C220350A
E20559FE3232BFA3FEFFA88E7FED01BE8C321298DD51743E0A75EAA26FC7D9811D01026AE563BB2DC2184752EECB24343804C6AEF15B32C6BB7FDEEB
A633BAD7FFB68BB7CBFC16959B84E56CCC1138CE72FD6A42A74A2E932F55F56A31E1BB862306B92DDBDFAAE0CE3C4BF1949557A6E643C3D60778917E
7FBFB38AE88DCF71C65A944A0B52C75F4091889DE4C5F49F61297D89C694166F86D2AE3E4543086802251C99491642FCE674D0A6BB96A0911D876F4E
300240012C7AD5885667194DAFDD38CF09B2FB9791F963FDFFCE9A20CF1BAF3E28B3632496FF22D04E1279412FD538E3C479C1C626CDD9868332758A
9FDEF13891CE39B3F0C48C26333407862DE8EC33AABFEBFC31557C22B5773402CB70294FF2867DE379BEBF57F5A5535BF6EEB4A8FDF4968FC098A319
4B677DDE19CA90B9172CF40DBF4917ACC1FBE248BB432B3A7F96010E9C1323008E956AFD294D5E5492CE6926780EC7514AEAFD2DC603A2FF2721CDB3
54E2AD3BE5E4A3D183F725DB598AA428F59DEE2BFB4FC190D94DFC59EEEBEF26EB5F46E677A3FFDF53BF88F6F8A55CAD997026AE581F82114462AA83
A5810D3762A871DD5100328C9FEFA2C68BE96CF18E596B07D63A9C3BB4E7EC462D96D46CDECDAB5950E903A3000013713A7F8C4BDB83EE359A2E6BB5
DCADEE8AEEAEEB984DEB8A67F40E3D7290244C2F0AC18739ABB71715DCCDA2D0220D2AB09C0D373BC50C409A40267A7FB32633BD946638A7D1A6D5A3
503A092C8D29231D0401EC1645591E10F2E7B79B14A0A5C5D080B6C52DDDFA62438B44F642DDDA3192FB97F52F23F38706E0D7B1C0AFBB8269EFDEFF
52A0362893D7353FA873E0A8B928C346A08C335E856D394D39531D112178D699143A7BFAC2D96D3B863FE93B7568A7990D2BD7A8D93E2C716586987F
83603FF58F2EAD3A7B359FDAAFF7AA9A8B47BAF571DE6B95FEB9FBADBCC971241D13407C19B37B778EF66929A5579A698A351639048491069032A2F3
271D8A92728211286B7E7EAE42F4F82C2DBE61D56A510C371ACA2D982028D6B4E87DBE58EA544E6A43C74F0A672ED67BC6C74E18DF2B45F8B7FEE503
003232BFB3F8C7FF621DF0B773415993FFDBF7E9692996F2D3ED675909339AFF2885362B61698083DC7117A8FD951F62C9B4B571E4EE766BC7579F1D
38BAC6A0AABD3B7CE3ED56EB6ECABA1881A668FF5549335DAAD41A527B5AABC9F5C6AD76FD569B59FB4C79EFCDC573B793F0F97EC7BBB107771E47B3
B2119BC241929AD27CA728768613A4627F1ACB8E4D2B77080CA6D03A184C3A0344731C4BE9F59AFC7C85CA890B82F648DBA69B633108086B5ED0DCA6
134BE84B0DAFB11973FC8B7ECC975A85F3BFB671323232BF5EFDFB752F50E167AF59D106C092FC2EF4657C8E45FDB2F726330B28B59AC04C0457164B
9167DFE01F5F96BD4B036F674616CEDA78A8AFCFDCA9ED86351935BE4AFD9A0DAF959ECBE738274C3CAF5FE652D9676CFBB5EDFA78CCF0F10CC187B6
31DD1D92E23F219DFBB42F2A68FC92971B94A6383BAA3552C6B4CF769C26154E690990651865D4BB5805C431850AE15941EA02C4480F020723EF665A
088EF8BCAD5BF3E9F7F5286E2E0BBB36EADB4EA72CE4B18EF759C5F81BF8ADB5C8BF7B7FCBFA9791F9833CFF2717293DD0C069361B0C08414BB5B3A2
089D79216F82D375DAF2F83923A329012B510BBCD944270690E49587D9E73293AF17A1BEB38BAEB6D930B351DB893DFA8DF05CDEB94A9D5AA7D3D7C7
0B80824977D4935D3C1AD55931BE7AEB0E2D5DE6E15BDD4F44F7BA9A38E00EF376434970CFEDEFBEBB05C28BE9A24CD41093AE662C797A0B8193503A
82F43136C7A9B75BD2B52C2B90B4687CA098FF730C2C7F7EB998E1D898D90DDA2EBD9D63722A430FCC6FD7B2DBAE1860D931309034CFD886A44F4B85
3FB97FF900B08CCC1FA8FFAB0960098EE1495B927FD8E7F78F6FBD78F0EE754444669905810E43416E89D184DA8AF6749A97CDA345C598405020F81A
8E5FBE1112E8789F48969F08538C1FF87C458B1E3397AC693CFF4C839A4D5FAB2F24F314C526EE0B1A5AB96E1D9F231D9BF69B54C7F3759CE7E08C39
1395E7A62B139787976CD8183C7396252F9B284EC1AD89A584353D3C1BA3753A9627B5F1D179A4BEA0CC21BE1244A7CF1328C3500C45957D09C84078
BE68599D3A539F95690C9F2F4EE858BDD1989B25943EFFE0AA1CC2B07A8D437DE3282AB7FE9591F963FD7F2D909506FFF2546E4C627E527CD0CB70B5
CEE2206D664D71CC9BC347B7DD0A2E353B340E0CA084F96D8F36E7392EFA3522A008F9F61EAA3C78F96681E17669F99E9D8EB0DABEFECD3DBACD5C36
B3CEDAC935DB8553F105101268F6A18F8B5DAA56DBF5B076AD3EFD2A4F8AE950FBC2874E875F760B281F7ACDB274FECB4DCB0ACC6FD49A4CA7985950
A637F1E510352A308E53C746DB1173C41794415816029EB659099AC11DB971D12A8C87C9ABEB78F459E51F1F757941AF2EB55BCD7DA5A72DB7A3135E
981DB693D34DA6ACA91FF15F6439F2ED9691F95DFDB3344538EDFA8FAF338A4C08701A8B15A5AA321DC5329C4D9B7A73E2A83DC905451A31F746EC19
7BDA6F5631D1CFB53C4E80DC9482987B371FDAA35E3BDE74BD004E74DFD9C2A365DBE5275A57FFAE49D3A7C89B5492C469CDEDF253AE55DB064DA8EB
D3A35DD7F77ED54686AD9895B0702D88DDACFBB8F4F08C898A8C8822524991095A9AF00FA320AD37908CDD90A1B0C0C48F69284DE300C5390151A382
BD2C2D215A81B1E88D5EAD3A4D5FB876D6D876DD865FFA7860D56B0452A95B1E63A8D966FE3429137E79B45A49CBF75846E68F56FE2A3AFD4803FF1C
164B765141960AC770A7C3AA2FCBF9F8F0C4B522CA69450840A9AF76ED762636285647D3A44DFD78FC9860B6E46C340D0880BEC8B4FA1DCB2E5F790B
3D58FD26B5B86EFB26F5ABF6DF3AA666AB3AAEC71433C301666735A7E20F54AE3262672DEF662DAAAD0D6FDEEDC6B9EAFBC2EB5E03DB77E6EEEDBFFE
C6DE737B0D5085D1FEBE0E42F5C84C235A2744D499AA7232E343140A58CC495AF400279D76D49C969CA64661FEA97EDE9D366DFAA175BD96534FA53A
AD6FA26C0005211B2238BD5A0D7306DCA6FD33769F31C85E5F46E6CFB27FA9C32E4BD88CFA84DB11B91693CEA028571A300677962B0CA8BE540B343A
8CCCD9D47CDC39DFB52FD53420CB23378EDAA7371D3A640608F5F1051DE7578C04AE4B47C6D58EB20F6BB7AA7F9D1E7B9677F9BE81C739E3A8208063
D0F1ECCB31F7DA3D1AD5ACD7D4A7D3A32D0DA685FD3838FFF2A042EB80E1218B87BD0D9EBE338FD1E2E0DDE9740C16A403E83091B6D498422AFBC07D
31EF27210E2914252D2C921C93926DC268EBFEBA955B8D69E955B5C9F8ABB9A2198BF22B16B30C7BC0C90268CB76D0B99BF7958545158FDB6693F52F
23F3E7FE9F632882D09AD4855684A101B4EA7273C25E446B309AA66C0E45B9D1A634E064F8B401BE7B470EBAAF8240137E75E09C7CF2DCCC22E0A41E
BD21536332A9B75B919C7A5DADC55D862DEF506BC6C1B6AD9A7FB3073FF0098861B9E1D0CB73DFD4F011DD7FED1AB35FB56CB96265832745EDFCC80B
CD4EECEDF9F8C4E8CF28AF4660F06D1D9D6AD59A0847B1B15495A970846C7D0B044E60710AD8702705B1ACC87803C6112F3A55FAE69B6A6D7AAEFFA0
63C4DF037FF4C006089525D6CF48A6971182E6C233F47398F564B39DA4BCE92723F367FA17FF834E822669A9AFF6D76C80C4F2762E5DB7F37E12104D
0065D21B11839327C31298C2959D8647DA59C5DBABD36666D30F272B21343F3483C04D0AEAF83DEA8ADB1AE67127DF85B53A1D1BD9D0C7750F0888A7
442352B429F27535EF26F5EAB56AE8E37BB5C5A8D503FB16AD1DA92FFDAE67F28223FA4DB182607901A26F036D929A4639DA5C66D718F3F35FBC2B10
F5CB98D2D2ED144E93943236D141F2F6D77D2B55F2E8B6315C8D54889B2B7B98C8403B026C5F0ACAB31C02AF08CC63882CE79DDADDC339B9EFBF8CCC
9FC7FF3C238DD86225F17315F1004753C5417EF38372324A080A02806218E1006CA186236FB46CBAA504965D3E3169CA17FAD16207CD94BDB66B4EEE
B76B3786E033DC1F83719D37D6F49C3FAF5A836F56E1471F020223D21786BC71FBA66643CFBAADFB9E1FD362C7B07AAF82BC5E69C7F73CFFC3DAA8CD
898200CF46E5DFB7E0A1F99C51CF9567A929FDB5659BA300CDE065B129597A848620EF996F32C63B6E74A9EC36EEAE429A05CEB190078E821284A709
1B466506A715439E33C4A433CEE872D3F78D8F16CA933F6564FE54FF5F4380AFCAFF79498036163EDE934E1A4A7255C56692020C4EB244441943C74E
AC36F099ADF0C69A693FA6905776D8193AF1B1C638F70E15BBBCB0A465C7F24FF5160CF01CBCC2BBBEEB06B0C68F26495A77A620B46A95EAF5BD5A0D
9D7FB76B9F039D66297BEED0EFAABFE8CAE4CFC7B6B302D817CDDDD3E1B7C220EE649C596A603ED46E5A224296A7854697B08003A024C8EF550264DE
76756BB4275F543E4DA30EA5C66E533B4517CF120ED492975962113F69CF2A27E8AC54E655F53D9735F2DC1F19993FD33FF78B1A807F1900C912B0B8
BDF0CC2D3BC04C7A83DA68401D182994679214633FD1A2C9CEE4A4FDAB7B0F4A04970FA22CFEFCAC3D64450E797799EAB9E71A7C7ADD455ED5E634AC
F6CD52E8FB90A648D272FBCBEB4AAE1E95EB4F6AB9E048B3C97D7A87AC1862FDD4A0D7D25E8F3E74B6086063ACE01742BE7F66B19670CACC721AD9DD
EF7C3E991F1B9AA760693BC9ABDE3ECFD692054947FA2C09D6083C4B420CB35BD466A662B01F8DABB47A755A090E0112FD2E544726BFC7CA5B9E3C3F
D7CECBFA9791F973FFFF8BD33EBFBC0809EDFB6D0104C00994B23A490A650494A611864F9BF46DDFBB617E235B4D2AA26F5D4459DBB9B7D8B5BD6A64
F70E62639D8894E6E3FA788EEF50CD75111B1C4B9390A062BE7CF4A8E2E53A715EED05837BCE6A72E051BBC2C27E3D0F4CEEFF78F36B8158FB42087D
CE85F9A32852565A88AA73AE5D480285CF038B2C2C429B6C58E6FB7C844958E233F3A1519AF8455024E540204D7362CEC2614E8DCE80AAB4760AA514
4F6EE5EA59ED7B3398BBC4D2663F252FFFC9C8FC15FDFF27F5F362984DE863F7CE4D00B81EA72B1A691062348F038A276FB4FB76D6ADF53E55E6EAE1
E370022A3629AD8B4FE0A679778CBDFA380FD51DEED9BE8B5BE5A1D8CD77348E033434A7A87EA54A43B678F90C6F33A1DEA82BAD9E80A555A79CE8DF
F7C73736E4CC4521ED1E9714ED04EAE43C4B79D2AAF1218C3378478015620EA3BAB43C390FE7B5DB1A0E7DE060C41C04B31B31D48CD290012C4523A6
7233CA201A2D2000991F546283ACF6851ADE1B8A3EAAE3CFC8EE5F46E69F1A006906284D5AF24E77DA8A50E672B5458AB53916D2ACD4902B655CE596
A3368FAE363C07A43F4A25828F191F760B61E2862745D4DEE71CD0B251A3E69EAE03347B9EF094A8DA7B57729ABBB43CD5C5B5CFAC417DDA9D9EB90E
DBDBACEB92353D57F9EF0A315C34E63F8838F41CA7F5F12566D3EE61BDAE22A5FB8EA75216525DA02A2C45118E0A1E37275874F60018F52446703446
71406A0FCA102425DA2582209C6A5556B189030C882886B923724D7DC717CAE9BF8CCC3F3300DC4FFA87A4B5F84EF30199C058189DA3C1041E5200B0
2C81D3AAD99EEE63239755EEFEC41C77E423797E77E9D26E39CCFE31F89E9A1FC33C7D3A79B856EE507EE82AA4290C9E5B9C58C7A577CF2A75BB4DF1
69B875C130E461AD4E2BA7351FA67C3C08CFD8A52E5DB7E785914E09C92473660F5D7001F5DFF404B56B8BCBB28A4A28CCFCE1D6A54BA19CC051A2B8
4BCD2CCDF20C4D8856085A2D6242C20B2C4502634878561694DA8545A600E3D2087243F3D326B9FA4746E6F7D5FFA701805417484354F96E48935B14
5417E59844B509945DD41DB001EBC15AAE3D0BA7B87A6F45D237675BF75ECBEC3FCD66EC7BD1D8A163FEE2EAF37B57AEDC157C7ECBF024C5BC3D56D2
DEB57A8D665D364E6B32F3C2B098B45E8D471E19D7FA153AE10373492184EFB12379299F7568EC9451F76E2485DD89A6E87CD2969D4F22E5D7268C7A
AAE5A453C862D0012D409AF6C18AB68063CDC96AC0092C238D0730451539ED1067D9CF9100DB72177C6CB6209894F52F23F3E7EE5FF83DF54BCD3F18
9A268D093BBAACD448B57602430B8028CEB6E10E0477DEAF5369446ECB4A358FE10527151993FD333B6E66DEB6CEBA5D7972925BAFE51E2ECDED5F1E
D162B400821E5ADBFF4F4D9F2A53CE36E87DBCDF697C6CE3B1A367B5D8ECDCB19EF0BD2D142D4A2BB997516EA2BF4C5DF9C5F77472729153A564CBE3
B58E42DFA15DD766331C4D320227BA7C146528C84096C51CB4980B307C4534C0B1D6C4124265231818F01C8727B6E2E8B899EBB3E4CD7F19993F3700
FF2113107E3E1750316DDB5916B070DC8B883056EACCCBF2A8325F6BC1EC18F9BE49A5CD3B5CABD63843A53D2142D6E9A33B3E61E62F32F6AFF278B5
7B072F976696A20F1C60692E6E4B4CF31A4D1A0D3CD4A7CE8A95BD534FD41BBC6D52AFBE799F47A89316DAD9FB6F9187B990A263E6DCF872FA4CA1C6
A124504B41B6C679A7EF806B69368622395A0CF41996C170064AC33D7802214928B52662188EA248759689B042CA1EE34F81BB13F4F0E9AA83FBCCB2
F79791F947FAE7859FDC3F2B85D7627C4D590BAE0E3F7E765DA4180100C0100AB5A934172561EE70F7416D1B55F2BC42870773FED768BFF6F1CA6FEF
5F7269F7CECBC5C3C52B2BE524C903924F5D18D1BCAAB7E7C645557ACCA9B931AA619DA9EBFB7BFBC68E8B508D8BA6FD7795DD8D1570E6D3EC07991B
2E59152976C6199F9946A42EE9735EC3429210137F5AFA2A1C697402C89090A6718C10CD8AE8F9391235DA74161CD849D2905C48E0EF5766C0821931
DD832959FF32327F4DFFBFED0456D1FA8FFFEAFE2BCE07415CF76C66A6FDD653B5188863A8CD82E83F879543DABEB1AAD794A695AADD662E6570CF4A
D935DF155F1F9BE053E5C10F2E755CAB15661C839C1830987CF533DD6B2E78EDD368EDECAE9FD77A379FFC43CD3EB1AB2E90870E50D6750FA23E0B90
7D3EFA55E2EA07C66225A0CB9F64E29A475336A6712CB0E11529BFA874318F102D0E255A002825151CE459C6662E2CD0012785E23A6B8156C7EA0EBC
5652C69D317B07AA1959FF32327F51FFBF8AFCA5B5418E974E034973B7A40880E700A67F39265B30049549DD413406A7D3947A239380C4618F9E2B3B
57AA95A31CA2D7CF51A203A7EAA65F9CED32FA8DAB9B8B4B46C1619C6438DEB25DBFDBADDAA57D9E9D5736DDFCB27E93317DEBFABC59FEA32DF847BD
7DEF9AA7DB8D82DD6FD8DD84258F1C256A1ACFB8F3C55E7064EE3D9A870EC26A07B4E8FFC5EF201D4F6218528CFF490A7050B43DCAC2D4BC6204A729
476A5C91A698A28D2B3711D0762E1F6B7DCA21A7FF32327FCFFF0B3FBF123D2EC47104C144A189F2139DAF43F77EDC6341D097A88B1D845E65B2E079
B7EE99017273FCD0EBC3AB76D51D1B5938683811DDDAF7D9E243353A678FAE54C5233B762DC9004600470A56BB7BB66F5C77CAFC366757D61CB8B445
D5AB65C313B09537F1DBC33F3C9CF63EF0E97CDF4F1B6E9B924832EF5D447CDED319CB52681EA57087159716FCA4D349B4188AD0D2F29F180100CC0E
B0E23CBD9DA2380C29897E98506AC368FBAE59A5D0762A85D9DF3A9994DDBF8CCCDF8CFF7F7A2AFEA16D65894F1FBF8F098DCF34DA099C0214665487
8F5D5D2010F60C05E9D4184BF4A4F5DC61258E154D1E7B7F6ED3A9597BE6CE1DB637E2FA80F39B4FFBB81C3EEFE25235377E8943EAD76FDB18BDCEA3
8677ED0E931AF6DFD5BAD1EC411E23EC7B766B7DD7C4F88D78FE68C249DFA0973B1E6C385868C08B1FC59445EFBB75624B3A4B21A2ED00144D025A3A
F1235A0069EA1FA428311D70E69B20A4A0407310457431891AA02169F3C3931A41B5EB03ABAFB1D04CCBBB7F32327F5DFF3FD5037C9D9747E6C67DBE
BF69EFB9C75667D1952DBBEE076614AAB53A8D2676C7B0F73C9F1183130EA4BCA00C279EEC539220B9F990B3DB078C0B0D9838F3C18F975FCD5BE43F
BFCA8CE4262E3532ACD7AC246459E24CC903EF9A5DA7F41AD160F69646BD96B5AC111E3EB02CFBD2ED17D36E9B8E5E319A152B4E3CDA9302D988DBC5
4C648B15EA522B0D085AF4F82C041810A54E5564FD508C002010431185866200276D0550168BD5261A062755FE2C9213B4E783586EBF572490C37F19
99BFAEFFAF2FA58A599ED1C6BF89D7E1A405B19E5F71CF4083E20CFF1BFBCEF805251528B28E7D7FCA5E1A9C67D1A14E4D419E1EBC3EA826C9B55E13
0EBE593A2031B5E7DE88015703C64EBBE9D530B1858BCB33F3817212179DF4EEC8579E6E7D0E34A9D36B458FB6B33BB94F2B6F39D37C2477CFF43BDA
8DBE762B38B1FAD1EA5008830309FE95F72C9374C6002738E0346215FE1EB7E3D06EA0384EDAF003146A32409A22694A1A10A08734413B282AFBCE67
8CCB7E50C433FA5A832DAC5CFB2B23F3B7F42F65FD8CA328FAF51B054EE1184AA20C0C9E3E74EA8EEB3641D006ED9E3FF1D09D84E2754D0E1384C290
5D4E3AADE6F852F0E94C0E9BD7E3BBC1EFD0A5C38A2F350E3C5EFDE29D6687A7565A55DFC5E541E9774990627864EEDBE83A9E4DA6B56CBE7843C3C1
0BEAD57D3FD5EBD963FFA4793BD57E93B34887EFBCD4B37B9D647224235C74FF11673054943D60EC250508434396166D086DB1325C457B229A90FA81
E0046429BCEC630103214193B6D48840DA9915A112007FAEC93B9A95A37F1999BFA37F29FE2772CFEC7AA32500455AAD7A1B82510288F55D3B79471A
1404B6E4F2FC76438F9DE9D62C8423F4A6A472548521054574C28D38FE49AF9E83F2B0614B4D2B3AA68C6FFC6286E781E60DABB9B89C30CC888624C9
525B82A23DDDAAD569D37DCE809ABD3BD759BFD37D67DC8C82C5C3921E0C4836665DFEE15DE896771F12A27966BDFB01315DA068CC46989D363BC1D0
62D4CF32E2CF20298EA9B00034840C101F28A7DDAED1301C6019738E31A7A858AB33A91915CCA83ED4CACAE1BF8CCC5F91FF2F06E471AC3928D84A01
5D5244767A7EB69AA4504CAAFA65CC45F91A5AFCA0FDFEDC76BDFB366D16C8381DC658B3210F305A235F7437143D33C967922AF0DB272543D6C5F79B
1FD675FAEEF19E2EAECFC8900292040CF7AEA0A0B147D51A0B7E68D4AC73FD469D6F366F697CF431B8D3E1C2194F49F3D4C61F14CB8A0DFB5358FD70
CF97020F1986C01C26AD8DFCBAFF283E883F03E05F6B7D21C5028A61ED369270885E5E320699715AF1ED220D46684A7505FDEBDEA4E5CE7F32327F22
FFDFACFF733CA30B3A5B4A94BF9E3BF25289DEA83321184E203846D19C402A0C084708024C585BAF76B5DA2F69BBE1CB17A230D4C06366A1F4F0072A
659EE7A284614D73223B9D78DDC4F782E7BA63DE2E2E0FB1774538CE01263437B34165D7E1D73A55EE35AAAAD784A1F5CEDEBB468E19913B679E957C
DD7447D2B637A9C732E1F5C65E71A2111258A7D1515268C1C88A33BE0C4B8B41800321A513401C235E1233FE9C408553CCFF6916678C2F1E3A781E4D
D643951A31B05175476B584E6EFD2123F3E70100C7FD34F593E3183427E049016E2EFAF8BE98229D841371E024899388D3E12458B63CABDCA0D7899F
0DEAE7E65AE31185E69F4AA7B22E46F1342114EE4E14727B7A4DDE5F736AF9B9DEC15B7ABCEDDB6E5CA52A2ECB927D4B1892A1F9F0644547D71A9766
B936FEA1AF57FD89B586249E2F4DEFF432AC531495D2D5D7B4E0327DFE6DC68CCAF592A5AF459B2CDA82621321EA5C9AF62B86032C8B96230C0B1808
1086447168513A190C253992808C26C12E70E60415B4E94894A0A63678CDB2F2B86F1999BF9200B0DCD7FE7F2C50264499218B63EA6C04B79A8C2683
D549A0A20120080766B7EA481E55EB53CB4841C0CE35AFE4F68072BE3EA2623407AFE30210628F5BB8006FAFFD8BDD8F1956F5BBDA69EB4EEF15752B
B98C4DBD9CC3028213EE5F527472F1EAE75E79F2D02AED07D6ECB2769DCA32665341FF3584F2FB85D8FE1B611B9C65ADBE19962FA9DF61B7134A0D2E
EDF653140415F13F4B39204DD1D06ED498AD340500279D09642005BE86FAAA08BD505E4E1256EE7CD5796AF192AC7F19993F89FFBFE6FC9C54EB4B99
223F62A2FA71D6141C6F50AA747AAD0525298C24288A42319433A5A11C4BE316254A0BBC6A63CDDA8F08DBF3C77ACE7AF5AC526085D01F53C9BDAEDD
CE0F6C971A5D6FFAFC196F074C18EAFACD5A2A255D3A4224BCBC64FEDED5D5E57F6AEFF3A9BE7D60ADED6BB75037FAE55EEF906F1D31CA7A6FE2BBD5
CF0A26B43B67E52159528009F602B3689248C42ED5FED10CCB887100943A7E510422DA015ECC000000346BD50201179313569341312A0541A35C4ABD
0EF104278B5F46E62F787F4EB0EB0B4D5A4D4E71A916F250E0085D5236816024A000991E8ED15FD7D1191CE5B04F598C182E60FA528774D0BE43BD2B
181A19A522C19BC3E1BCC06CEF556C1EEBD27F4FDD9191A35A9E6ABDD1D7C3C7D5652E880E81144F71AFFD94DD5D3C2B55FE6E44B54D332BCFD8D136
2DA0E5939CE6972C93469A6EAF4998141DD8DC3B8607C05E58CA6065E97AC0B12CEE404499D3A4C54943425AEF6770044AAB03623E4022044D3A8B95
1401048670941668C215046983F60195F7EA19D9FDCBC8FCA5E89F2F2CB125BEFA186BC348CC9CF6E67996522B46D75459D29B7B07C60CBF726EEFEE
8B1FF3A5CE7F3843463C314953C26D459FB3AD64F4A486C7AC58417AAC9D49597D86E2AC63E6E3D9CDDDD72E719B76B2CFF9E95D9F0FAF57CB7D69F1
EE9702C3D17C581039D1C5FD9B7A535BCD59E0DE60F58471694BC722EBC697CF1FAE0A5CFDE56DF4B55AB5CFA20C6353112CF9291A13188640310039
D1FB631A2B29E601A2CFA770AAA2DB8F98B250386048AB8DA1189A2ACC32283EC694439224D9B5557F0873C8C9BF8CCC5FD33F832695316216608FBD
BA73E5D41D5F2CC9E1D9E98FB76FD8B5EBE8C1D367FC8E6F9ABDFCE09991E32E18059A80392FE2509667C9D407E10AABFACC685F2766884FD2724533
F611746E8FE7E075A5DA2BBEAB7E6659D78B0DA6DF69EBE13638F9D8130152ACF0E10231DDC5C5AD69EB558B5CBE5DDFC4EDF0C64651AF7C924E76D3
BC9D9EBE726D6C8BE60F05DAAACFC505FB87389C17206947A53E7F2C0D20A0299212430EC810040BE98AD3C8142D26260C60595E9D9463B7E94D34CD
3808FE4695811B421959FE32327FC500083CC3C3B2E23C9D9D446C05053AAAE4717054D88BFDFE5A3B4262244E389D667D4ED2D9C18D3B1EB34A29B8
AA30CF29BADC822865619A36E1B8BF012285D17A2163D83900DFF455511B5D066CED3FF56CE7BD73EBBE9F5AC5A3AB2A3F9CE120CB7FDC48AE71A95C
BDC9CCA335BCD6AFF7F0491833D03E664D7AEFD7D9935E7FF068DFFEBB0C0E50504DF1B6D4149C878C34E343B434CCD7CC1F88F13F03A4BEC3B4740C
8865198A0500972A034BA34B305A6B6559C4890A9FBC5ACEDD6D97837F1999BFBA00C0F180C60C851A8DCDE2288CCE4058D44C203457717A4E3A734F
90769DDA96B6DFA7EF591B0D29682E2FB202AAC808F24EE598DE0739583227594DC70F7B4011CBA613E6DE95E72F6DB36749ABF3DD678E74736B9C5B
70DF2A0804F7659E76A98B9B4B9B3B3E1EBDA7B4EB7EEE5CA3AB5BBB27CF5F669EBFF761BFAE5ECD537806EA8CE58CA5A054EA2FC471B893A62C260A
32005292FB87D24A200D49693F90A371311780806189C8101D6BD59A501C313AB82F8D9B6EEC57CCCBAB7F32327F4DFE52AA2C0DFF63308B3AA7AC8C
E02A46668AD7BE6EA1F10C43028A41006BFF72B877FF7B005A8C9429AF8444B580FE74CA8C44C5993850145B0262163EC6945D0F81902AAD36CE99F3
78E29A63ED6B7B556F9582C4DA780E70D6EB9A0395DC5DD6CF761938AE7ABD9B098D7B97F5DC953135E6D4D22F133B0E6AFA160F8DD45A605E4A8103
B298E8F3318A048612ABA8760021C0494A74FC522950452F000EB749877F314D8A7F024562460BCE60263B93D0BEFEEB757BA1ECFD6564FEBAFEBF9E
F90524A4786927F0EB3050F1A9F8081C14A2B7198C06A703702581DBAB4D50D094CE4C66E603C244D26F9FE2485480059ABE7CCEC33236DDC1137C9E
500B5C3A4C6BB3FF46FBCD5B6AD5ACE6198A3E2A1744CF6D7F6038EFE232E2F437D5A70CAD54EFDD96BAF716764D9EB3D43271DDF30D2D5C6FC3E031
B75065444CAE95B1A30058490243AC6692E72180A4D94A0292FA5A07C854E40005A538C16061D7E30D8035A92D24868BF2CFEA5975DB8D06F972E5AF
8CCCDFD3BF1405D03F9DFDE77F7E87C31C06BD4E959118F225B1DC6CFEF07E55CBC67E0874984D399818A573F64779B42524480DCCC9EF1444FEB600
22A063529467D3E1DDBA1DEF34326059756FEF80C8315F44FD33CA59396FAB343ED7F19B81237BD4EC3CBBC6EC9C61178B4605166E9A31647AF5B394
3D3C992ABDF34A4DDB44995B111C6065854ED11AB10C4B594B8C9096D6FFA53EA41C0D580E4F2A147304557439C3118632036070878D4E1E5A7DCCFD
06BE505EFC9791F9FBFAE7BEF6FDFBF76920D18F12CAB25C8DC6EE4454A9CF1EDC4D573C5835A86E8DCE8114652BCBC159E0E4B2EF2278FE85756976
7D605039FE79692E756B91719C4BC3815EA3E757DD1ADAD2A3DAC3F77D42049AA01543A343DC6B3572A931B953CDEFA7D5E8726DFAA192A907B1853B
97D7ABB61368AF9450397E1F8D662729E0561B4533A5511AAE62002104989D62282065FB0CCB41234973E5114A80D8344E0E87D66C15244884A2DE75
A936B964E1288B9CFCCBC8FC7DFD57689FFF45CF2C517E84D62CFA79691238EB507CFEF821ECCDBBE9EE55AACC3640ACD44A73344D064780F2CB5326
5D5139A33EAAADE77659F1CB49EFAA371A31A8FBB976DF25AFF4A8F7CC7A3F5B8014ED98F1F98BA74BA54A8D17F6AE3E7376B3AB276A45F8B7CAC819
E5DBC7638693DA76C49173390AA18C902B8F2BA7795C19AF9222124E3400A0620B8061A4AD7F8AB0E9000355C53846234E92219CA54692C508485D6E
586B83E644CF0859FD3232FF40FEFFEAFDFB6F0B50D1FBBFA2DF8654C02BF5DE8F7A7FE86AE89A86555C5B47D0A4C6C2F2569BEA651E1D77BA77AF33
69CAB037411F4EC5E3416BECCB5C7AACA83E7B75D52D97ABBA5F643EE70834A0F0852FE26B5776AD516F589BCE03AAFF185563505EB783D143F72D76
FBDE4E5F58A04FBB144D19ED045712984160547E928997FA8E333C27ADFFD1152BFF246450A70DB08CCE4291368E6111142208074892516FAEDEEE26
79B6E10D79E4978CCC5F3500BFC900BE5EFBD93C540C00F9E97F866179A4F8DDE219AF92FD9635AF7E9084669C660C4C79603E91B677CC8AC7991F42
53434E3C753AD7FB957E5B6DDDCCE6E77B379C58A7C60BB03F5810030578F071713757EFB69D7BB59ADDBA55E8EA2A375ECD4BECD7EF569D7A397CDC
841CF3AD3894B0D1BC33A99016549905A4F42F8A7107B448A7FE389AA549820090E420C32AE20043106A2D20395AFC04E460647FEFC559D4A5FAEB9C
B2FA6564FE6100F02BFD7FDD17E47E1A02545173C3E258E6D6B147CD86338D3D97E09CD526E026461DAD474B370CF9E488BA1159BCBFF37D3467D0FB
D0EA2D1F0D19B9BC65C3AA6E3B6C2B1F8A2E1A803317953D5CEAB9F7ADD576A8FBE207AEDE4747DCF9B1E9DDEFBC228488C59FCB0F86E00A2D64F1B4
12812A4951330260291C70585E9E54F1074969BF9F91E67ED108929C0E71DAF42409301CC93202A3DCD871C64794BE536F48D1AF97FE64632023F3F7
F4FFDB77BFC6003F5900310D573E1F373FC17C75708DC5B8807102C970A69C523C7FE92C8533FCD091E0EEF55F13CFFAC45CF798B963D094EE7D9B56
1957E61BCC5200C037778CA3BFA951D7A7F5D46E2D4E77729D7879E6C8363B67D47D2FC4740A826F6EA34A0301B0BC621A4D4F344B96C75A8E319832
1FE3391E22668774DA0FD22C4F1300C700C91A131C3CA0A12038EF0E9F1882D0F0532B9F409A1764FDCBC8FC430BF0BF7383DF04001032CEA82DC35E
E289BBEB2CC004518B380DED7985CE94850F6CEFDF1EDBF6BADDD052626B8FA45555871F58D3BA6523D749B97B9ED200E220EC56F90417AF866EDD27
579D72B16A93890F36D619BFC4F38C90D1F299ADD42FB3280307BABC6CA08F8A475901D08EEC32278558699EA76982400125060190E3189AA1399A20
498CE2188AE79DAF0F9CBAA9E34926B9ADE76187AC781999FF13FDFFC2FDFF5BFF9209009028BA3DE92041DF6C30DB2CB030259BA2F48A1234F286F2
DD53CBD667811DF652C641A39247341F7F6B4BA3FAEE1B1DB7C3680A10E0CB2DD5FA2A8D9AD59B39B2F9C5634DC73E39FEDDC0D92D7E40F1B13BD99C
7731C6501BA98CD36BA2020B30D1FB130E4589112721C3332C415052E50F0B014132620400A58EA01401288E1362571D2962048A6383DB7AFF50F69B
96BF721D808CCC1FEB5FF8FD00E067FDFFCBFB4B6DB824295A92962D2C145ED41BA21460C1DDD71A0C4D29711EBF6E3D1F983BE1EA85862FE8A03ABE
016DEACFBEDBC4A3F20FCA371110C7081079A8C4B74A0DF7E11BBD3A5DEFED3DE6FDC0BA93EA7B25D2478695A4057ECE7A5B4425A5AB547ECF343410
1813826A6C34833150EAFA4733521910CD40C44931503A11483B1D38049CF9C1BE205C20204F1CAEEF31389EE485DFFC02F20D9691F92B3640F8E5EC
BF5F860692FAA5BA3B4692BFF4C8D10C83965E18F54608F5E992C582D2DBF7A31CC644FDCBEDCEEC43B657BD2E6F6A170DD774CB5BDBD07BEBC41A9E
A3D20F5EC5294082849365573DBC1BAC1ED168855FD30177F6741BDDA5C67EACACF34BCCAF202D3FD990998815AFB96A127D39AD564128FA7708696B
BE8162A439A415B34829A242FF0C2DA61D242413AF85A8053107E1E367D5AAD7FF35C2FEAFF845BEB932327FA8FD9F97C9FE5300C0FDB4FB2FE6DE40
D49E947D4B338139D2143CF79A90D7A745364D18FC8FBDB1A4FBE76DBB4186FAC54D681B34670C925EEBD0A37A1EFDB6B855EA9E74F4BC0D1018C83B
F4E5869B5BEB71D51A9D9F517BE1AB1EADC7D7E8AA708C1E90FB212E220E4D8F2A86F973F69A29074FE764D980E8F4091C572796435EAA3D90BE8654
058413507C49D911D21AFBFC731147E1E2174DEE51BB6EAF976AF257D54BB2FB9791F9FBAB01FF5BFF8CB4130F49968050EAC52966E0A214A12D6DE5
41129DD2CB0428CB33DFC7FAD08CEC833AE7892D8F7B6D4C1EF90CCC6EFBBA63B5B6BEF52AF5D50446911005D071332AA871AD164DEB4C7BDF7BD8C3
AD5DFA35F60EA4225BFBE5DC364732F961E54CE1C82388C34163A9310E41F4EE3633AA88D648B907277529E279E9E02FCE421A00D44180C298528AB4
5B05419531BD79E74E570C18CD0BB2FE6564FE2BFDFFBAFE47DAF973E8CD45B149372F0686A67CC928361218C5B280C4AD49ABE6E898D5BE80D2173D
59F2C6B43BD9188D44F6DC7FB3CDABCB03D4019E8737B8D59A55D7A555DE8B188A121D37792B38C2BB4A6D8FC64F96D4DA70A355875E2E93ACCA216D
AF3E4FFEA4D0F9ABF8CC3EBB095C0BD1B87827CB03BC30C78A29D552B6215900E93802CB1024EA74D8157A078D38ED24C08D24E70CC98AEFDDF87B5F
0D8E33FFCB7AC9775446E66FEAFFDF5EB4A235085B9A1453941EF7C62FE8D5E9ADCBB76DDFBEED7E82C14C41C2695385FD30349B092886A436EFF9E2
889B4332D07CF3AE66F7F7CCFC34FFB079F2F0C0E635BA7AB9D48B0B8D27485CCCE65F86C737AB56BFEAAAD0166D8F2E6CD3A7658BB4C48EB57F088B
CCFA881764F1F95D7EA4EC0E4A975D82F304B0E466EA088AAC18FDC30AECD741053C4EEA3372544AD4E6B43810D114986C65C179CA49EDFA9ED3D218
CE72BF98F9F19F0A1A646464FE7035E09706406A0E68FAEC7729D26A83787C5866E1DBB5EF703CFBFA820DF7D2F51862D2A43DE9D32A8CCE2EA5F3C3
8890FDA5F3BFCD76AA0B7B7D777FCA91475DAE6675B9B4BE8AB7878BCBC3C8835A0AA509EAF6838CF61E55DB048FF79AB4A861F7362E47B13995065D
791EF042834512DA81031DD040941668442F4F1B928BEC344E93342D75FB80342FCD20232C3866292DB05294416FA370C26A4548657144FFA6FD2E95
DA681285BF98F9F3D36906F996CAC8FC7DFDFFF48AD67F8A2EA2A0536F30D90DEAC2D2C8ABA7BE001644AE9BBC2CB8D4A0CA0B39D7B9D90B3C59A9F3
FD00F3037266F40DA06C1FDA0EBD772C606B5FFF7B8B123A54F5ACE4F2286D5B99E8C949FAF96B7CBEBBFBDCBBB5DA2C99DE6C50FD1E466DBBFA47C3
EF3D4D14D232F8A50D9338279E9089A334204B32CB500A42484031D360088A462CB4A04DC8D76973B4560A4311820098D34E500597B6B46FDE795B29
20315581542BC4FFCAFBCBFA9791F9DB0940857244ED5009D7020A735526BB0D315A6D4EA3D9AAFB74DE37861204FDB5394BDE67667FF0DBD3A2D967
504816CDF483B91F146B3B8551C4F3E647AECDBED7A55BEACA936B5C2A5771D993B72C9B4229C084FA73BEAE6EDF75EB3A7D60D71E3E95EF63EBDCE6
9DB874EF1D99BDC2F2A8CE7DC14CA644DB0185177F0E56511025A5F69E36A743A7C62946999A1F179618ADD5A200602405711C2381F1D4772D5A36E8
B32EBB28BF2C2C2F47CDFC5AFE72F82F23F38F0C8094FB0375ECB5A0903BE7AE3C8E4ECDCF2CCCD7DB11330688ECE7BEE162340E5F2D98FDE8E1F9A3
3B3A4FD20333CC9973CA167CDB7A70400CC7DE5915377FFBFEE6ABE2563CF4F2ACE2B2B768BB81B243928DDCA3B95AD5CBB3C9B28543FAFEE0D3D794
57A7A9FF86E057A9E4FC63E14D8FF1025716E100C07C7FBFBF8EC2481242862E0A549BB5168E16CAA2E2D2B32C98F8B520942C004E41C795DEDE2D5A
8C5A7BFD71CAA3A72109EA223BF7EF238CB2EF9791F92F02009A442CF90A0BE230EB2C9A7265515161CCC74823813B01CBDB5F1F4C103F8A5E1D3D6C
EBA223871A7C570CCD94F9FC41E5FADBC48121391C77E4D695EE37C7D57ABC7157171757977979ABB3708CC6B9F091611BDDAA551B3AA76DE72983AA
5F47E7BA2CDEE6FBFAB4DA6F54EC0F4729C6A18F2CA548C7D5DD49284660D2595F5AF5289791FA8C69CA8C7623410915ED40689CC4C5379F7F5FB559
CB7E2B9E1DD932ED467089CDAE2FA7D85FEE5D08FF8A62646464FE9EFE31274E712CA07956905AEE124E93C5A8FC70E17E9ECEE9A00541B569BF46FC
B06A55BB0EA38EAEADFC6D0669A7C827BED17362A87953CA84445FEBFE05DB7DC69EBEB5D8A59ACB74FB952812A3312E636D925F55AF8E1BFA34AE3D
A06D5FED67AFE6A767251F8F302E581FF48E8042564C818D75BC7C594E92769C02040D9CB9252C87B3F6D0370A4280024703C830340410D86E8CACD5
B4F7E8E9EBE64C1872F479668A9A24118CFD57E72F5E5EFE9791F947FAFFBA6F8EDA1D34CB408E97F6FFA5921F12C31C658FFD8A2D04856282401DDB
96E54004F65D77AFA6C717556A9647E908FCC4BDC8F9396563A72AC05C7FF3960FABAA4E7F7ECDA5B2CB04F4F1311C83186B7C5772DDDD7DE139EF06
ED7CDC57EB3ABB771DE97BF23271B55778840911B2DE3A13545464A093343970CC89E0E217003CCDD3256F434B180A13ED10454AC37E593275737B8F
3A3E3E430E4E9F7CF9695041A15E8193401A03CCFFDA8EC907806464FE410420D00426C9A9A2FF47850D601992B014873E5EBBAB945016C5E563B845
9F6BA605E3768FBA2796541AAA739AA17A67BEDFE2E2E8D1930A3EEC0537963FED30E4FC59B76F5CE601FF751401C5303EAA20A2BEFBD855CD5A4E19
D8EADD6A8FE63D66DF3FA1C8EB7F2D7D6BBE5072A958970BA223EC186A412D8A3C0B654E2D230446FF39361F88484D8868692E88FDC1C8DA95ABD5E9
D869E2D6C13F663396DCFC421BC502A92CF1171D8CE4054019997FAAFF8AA2DB8AF91FD2F3AFD300C4C09B5006AFEA3E304740B2F3721D1053592992
E742DA7AEDDC54631945E6610FD62A0FAF7A1F3B748D71637C8AF7CE99B527FCE855C96528FE683D49103CC0AEC6A436766FD4AA7E9F49DE3BDF7A36
F05E3C636D91B9F7D882F1974970F133F9C1A68FB742BD11B367A7EA1C6861324A1BC2E34B71121204C5D102236623F4191FB76FBC6AF6B999F9FAD5
B96D6AD694F1A50465002DC2FE6BE6EFEFF533919191F96B3900F7EF42A08A11E13FA5D40C56F472FE981BA8C0288A0AAC24C3B238C59AE6569DB0B4
E339A75565F57B187FFAE0B3A0098F1EED35AF5EE9EBD67565EB6A953A289F5F8590E459EC6C7C49EB9A351AD61F38B36FD49E3A0D3ADFFCE1997575
CDD7B7F761C29747BAFC7CA82229DCE03015665B48026558DDABCF6514603002D080A024352B06BA79B4FC7EE79A5B8981FB7D3F431645B3D21C94A4
7E9A61BFF62F177EBFA1918C8CCC9FAAFFE732E09F9EFFD415587ACA206541E3BB5C47196BBE5E4D503441A21471D1BDED24CFD9F63C7BD489B7B7E2
B65FDEB1287FE153EDBC956E3D9F7676F1708FB87F07C7044063BB03F39ABA7974AADDB6C3C2376DDBD75B346D1E72C36D5BE8E2523A7FE7F50FCF2C
7A1D8B14EBB5391A23915B8AE26977BE00C863040508022700070DFB6BB8379EF1FCF1D507410FDFA5B182C0F15444A29512DD3F849201E07E217D8E
932D808CCC3F4C027E6D10BE3E7CB50068CE9A3ABD0308DCAAB44B6702098404C12D6B776A768E85992FB2D6CD7A15F9627244D034D5B29E2DBA7D98
57C9DB3B2A3204123CE4A16F90A6AD7B8DCEAD7A35DC7BB66DCFB9615D430A1BF50ABD120D54B3F7449F48A0080AAA324B551A4A1F996245020EE7B1
344E9110400C27592E6F78F5AABD7E4CCC7EF74E8BD30C2FD019D14FDF3C7C59E00040EA0F0EA5CD01AEC206482B170C8E4359FF3232FF6D2E20F0BF
681320AD07D0DAB33D1B2CFF8C3B2C26889949CA8AC2E2BECD7B3488B619CF5CCE9D32C2F7DDE47EDA757B4E761DE7BEF94CA56F6A45465C4568120A
E0F8336B17D7AAB527F5AEBF61F8F0168FCFAC35F46E747D4F30A15FB92531E4A995447855546221461A42621D004F4C2280D529C6184E0483E694DB
4B862C7D0C784D3625AD4A32EA271F92C2D74EEEB72E51EF042445510010388223A89324C568807294292D1C27D701C9C8FC77D9C06F0203497CD6F0
DB5BBB1F7710084139AD00C370685EDA7AE85E54C89D1285FA8E7E163BE452C8EC7313B655EB77D2CBAD5946CAF41C5E9AE87DE529BBC4B56AF5D1CD
076CEDB87AC88E56C9D7BCB6AE780BF0874B1235C12601636D315F340EA6F4BE180A38EC46ABCD4A0006D1EA9D8E4F53BD5BBC011C87A695426919E0
F186436F03774DEB3776C1F94C1D02484001A7312F232537FEC3B3F84C4DB9D299916D3231B2FE6564FECB4CE03767EAA5AD005615FCE5DDB81546E0
44C8522D69B3620CBEA57D9B444138E0132D041CD6BFEFFF72C4C21D5D06D71854CFC5233C7FE427818150B8B183DAF58D875B5DF7FE3FD49BBABCC1
CEA45A13E6FD8840E554BF825B91108796A83415CDE63DCAC429CC5064366118CD501A05CA07B46DB9F0A903B716289D8260B93DB763EF855B46B5F6
FA76F2D1D7896A8BC96EC3AD96A2E8174177CE1DBF10A4B0610561E181B965089097006464FEEF9602FE6D01E8FC87FE45FBD7A4002761282391623D
C1DA56371B661552EAB48A151E6D2CBAFDFA448F8F539774F76854B3FA5DC3C1286902B010BE5EB3A49257951A35066F1D7768FA60E39426CB476743
F4CA81CCC48F669CD647A603C8A506E819C29CA8411084C170AC4CC170A1734EA5E2B8333F4B0759E3E161134EFA879D98FAFD94150FC2D3CA70CCAE
33162584DE3F783295CC79FCA114A7F2AE1CF78F2DD49443C8C9F2FFFF58FC28F3FF67F9FFB417C0733C1D3FE70111B0F00D61372914A8352E9BE38D
3F7A6E13E0348F7E25D8C60349C7125BAEBFD17753C36AD52AED72BECCE78140F3E94755872B57F6AA5179EC814947DB1C0CAC35A8C70B96BABD203C
CA2F0550F9AFF22D02171F8773ACE345980DC750BB03AAB3502174671824C9B2120D46D20F462C8F319A5E6F5B72EA696CA9BA1C75E07879E49B3D13
7E5CF0A48CE7701C02E38B854B365FBDE39F63B642563E052C5B0099FFB3FBF78BED35C47F741C889B760545D1723DB4BFCE1004ED7CCF2BC267EF3A
6B68E3A6FD1DCF5DE99BBE60DBB86F6ABB1D30AE0993CAF7395BA0E34B6317AFAA5DCE6D6B3C7151E2A0F683B6E2307EF6FD8CDBFE4632EA781A0D91
E06841E0CA1F8421C009501B81175BC9A05B668E35984C561D9E306FF41B5DCEED0307DF2728D4E5E55A0B4DA0259127E72E3AB967D58E0BB76F6F3A
7BEDF6EDBB1FF332E36EEF3D955248FEAB24404656BFCC3FBB7FBFEE08F86F18E7A5A12954F4E010D26E31A949C5D96C41403679BF137EF46899088F
B7E9543F69FA9C2B5D8FD676ABB25C3DE096C0009AD33C4A4F6CEC52D9A5EF9DF99E5D375EADDBA14322973C6EEF97E3A74B9D6F9687108CE6C81B41
2082AFC7409A702216C214E4AFD726A182D38990A4517B67F40FFEA147FD6EBC4E5129B44A9316455065CA937D5D9B8E9FDE7B9FF2C3C0630185A56A
2B414AB38AD5EFAE072A20C7FD5CC820DF4C1999FFCE86FF5C1CC4D2AAFD43A3A9B74B4D144E9BAD8E9483A201E01F8D8B55B6761BAF09F259DF75FC
A5E5AFBADE9A52D963B263C1358165206FBE1C1BDEB2B247BDF5C71AB79831AE75F701AB6C96254363020E2520790B1E52A064F355C0879D7C66E638
CC0E30D4E2BFE1BC8A112ABA8E3B1F6E5CFD7D8F5DFB1F27263E0ACD55166B4B4AF29429E197B7F66DDE6E7D64999351CC3E916763A4590115C5CB3C
A38DCAB1FF1400F0827C165846E69F89FFB75301B98ABEA0249134A2B30ADEBC49039C268A554907930546C8D89DFCAEA6DB21D38B7913BCE76E0C99
3472AF5B95AEC8CB6831FF07B46AD985DCEF5C5CBEF79F5E69F0A186BD67D5BA87ECEAF93C7CC973A37AF3391428D69D65D44FEFE50281612065D338
CCC5395A16D29C40EA42029EADEE307CF58B98ACE737CF8505BE7CFCE4E285FD8767F41F7BFE5E1E09AC98FA74FF1F429DD249C5AFC3023896E51DA5
06E4572581722E2023F3DF9B840A130039DDBD46EB28C2DFC288B2C24A9882ED1FC52C5FB33F64AEFB540BFE60F4C0A99DAF5EAC39B7B65B07F3D3AB
B8343D083F7A5239A8924BCFC30D5DD74FF3B93C61A825B7DF9EB2C3174D78D469139EBBF384DA189A490914CF19CB70458903486A26052E77F3009F
CD6107373CBB7963F782612BD7AF1D3477F2D4D12B568C6EE673152302B6DE52454F183DED9E5EAA31F8199655E82CECAFD42FEB5F46E6BF5BC5F99A
4EB390C588C861F56FB2065490BAF32AA2858239174994CFF24BEC522D02C9EC5A6746F39E1F0737E8E1D1B0F0D608A5409100DEDF18D5CDA592BB8F
7BCF999546EDF67EE558F843E695550564FEF15C67E0E11396D4648B80333C551CADCA2F71D20C016841D0FDD8A2D5F05EADEBD56ED1B3699D9A3DFB
4C99B8E3F19734BC30E3DEB88117CCFA191D9F39EF0F98347F573E64D85FE99F73185128EB5F46E6FFDC1E50984010AAEDDDDB7ED047575C00773E0B
91934F196CDCBB80970DA7DB1D8BAB2D19576BFBBB0E9D6A7927C59FD2733447D2F9F70B66B97856AADE7475CFCA4B068FD37D1EF4E1FFB1F7164055
ADFBE3FE7C7777D125482A82A058A08885DDDDDD858ADD89DD850A8A8158D8DDA888884A2821DDB061F7DEABD7DAEB077ACEBDE7C6B933FFF9DFABC7
3DEF331CE738CE88B336CF27DEF5BE9FF7EDA424B46A7A42CDFADE974BDFBDD793A416A979F1ACAE3007A17114A5A8ECC8F09E63BB877BF1D9F2611B
D71FBE959E555B479A492A77FEEC03B73390B2E9432BB1AD1E7D575D2E80C9DFF5FF7EF6872231831105FE0300FFD50AA031FF235A18AF57BE1FE2E8
F3ECCC2D9AC64D74DD8574F3F3B1A7F23F68CE3CEDCE8E86AE0A9A2D70B73ED95B2413263C3A036308096155A7EB6EC919C250BF8EF67EDD39516583
E6644F390EA9174CA93BA88854BDBFA126344ADA94F2BAC854DA78CF68430A7FE8613F776EB058E8342AB6567BEB4CB101C1D08F97D65D7AB6E7609D
A6F85DF1D9AD6AD374C781B169F52845FCD1FEC6CB0B081427FE417FE03F00F0FFDBFF86DF117AAD1ED56A1E766E3A41772099A6499C56DEADA11ECF
7C78E1E2A3777BF92EA9C6410E9B67598D1F67ED28B9767E42AE99C0512C6779C50309436167131A30B36BABA22F1D2FED5D5A03EF084BCF1F3D3223
EF4E058241FAFA0FB92AA8BAE13BD066326BA94D4090ABD86DD4856C04CA3DB0E76595AAFACE9A55D15733EB68F8EB9997A67747955F7ABA8D3F5D8E
9138F98F4B7DBFDD5EF8FBF96590FD0180FF4E2C30D3B8460DD7D5413ACDA9D62E77E8C7A5B4D908D1C5A74BE8C4D96967162C7CDB8FD1D57493DFF9
A0AFDF78BE9C7FEEC5F857146142B0EC858547190C2E43364D34DE75A37ECED0B3533F9309AD134B97ED48CDB99189C050557E5E6E1154031B319ABA
3DC6C6CED3D77DE8B12204878BAFDE2B825054797FF54513095186ACB8DBC5C61B7B8A52FDFDA7C6D69B49F37F04BC010400FE7BC580A9A8A85AFFB1
D8889445F8F6D7D34A84265098CEDC59409CDB91B4CA6FEF0377CE3AE340FBC3CBDBEDB4910A6394C7DE123802A3E52755710C9194D96B6CD3F19DF2
2B3A6FDD789DCCED3CBBE6E2B26475CA572DA2BA7839A7DCA82460B31633CD91DBB68B3A39F79C91A64D7577EF96231054A732EA6B75707959FE99DB
156469F415FDE3E623A3B79512FF517D50F60300FFD5FC6FAA5496692B0B6A94BA94D95EDB68B44647A246239D71A80E8D9977796EEBAFA785F23BAF
859DCEF75EEC2FE0EF7813F9420F633AB4EC58618A278BE5B5CE2A3CBC4DDDD6E6D356AB0C13DAA6DC9C75ABE85E9A01D6EE5C916528A933A86A7566
7C89A2FBDA9B95B939661A56265ECC304026A51E82D1DCD86B5F6BD3CFA7D1EA77E752D06BED8F3F9E924E51FF5E7B33B8050000F85F2C0698690AD5
94175566957C9A1DF49886B4104413109D73B4DEB87A4DCE8CA8AA29A2A08CE9D2839B3B788A45EB52A3928D086442AB56BD340E60303A4EE46D1E1F
1EE33B74EA1BFA8CC386CA4D5BF3AB3E6A11C3F18D55F575552852AFA5A84B61C363F21E7FA822D448D1D37775305A6F30186A72D32E5FA8408A6FBF
D614DC48519327BADF2D1F108BFFB1C707D91E00F89FFBFF2D0618F233DF3CC9CF9F105E4ED7D7A325040D916FD655D78F3B7575F0EB0F41C2C9EF5C
FBC7394B6D78BD321263F230A3CEA49DBD26DB95C1928BDAEE09EBDDB97FAF59F8D7C04E49C7163DCDBB5F89C097379568322A10754E0D6EDCD37CCE
975A5505F6F57371BD0943518D01D29466DC7B95A434A95E263C5315A6946A4C07867DAAEC330FFEFBD44FF0D1FCE70F0D00F82FFE4051FA8CA40FFB
EF7F1DBA00467484BE404FEBF147F38B3F4D4B7C7AF9C92E67CFFBBB5BC68E56B8F2BB575D084FC11BA7059DDBF4C19621647357EE918E5934BEF535
6AB6E2C0FB3909F0DBE730FA607D5D4E6A35ACADAA423F8C738D2C27609AAA4FFD8C9204A533EAAB723E3E48CA522BE1AC73B7D530F6F9D9BBD2B3F3
73B5C3C6D4829F6C00E0674091A60F59D973373E1E788636D6D1C49D64B30EBEB120FBD6E6F8D485D766584DFED87F54B4B55CD0F9CB16AF77A801D2
231F6FE947B3C46CBF5361E2C5A79C83EAAE2A06BE58BCA2E2D5B66AF445746DEECD0A589357A079DD4FDA2F8B345623AA323D45523AD8549B9DF1F6
6BB9B6DA5075EB42918EA84A7B9583A7CE2D20A6CFAA00591F00F83915004519B3CBB342275C5A9E4AC308AD5D184F42F0FD3539A7834F1FDEF82224
F8DEA96EE75A08C501997766E5A37A4887E49CD5AF66C9387DB70A9A1E9A275BF0A5438B9B27277F210FDE24DE1E85BF5C2C819465059AFC71BE51D5
3404D3A416D3923865D4161716285506B5BE222DFE8D0155D5969541E67791491F27CFCFA780FF00C0CFE822BF6FB1CB2A7CD5BCCBE9E832448FD05F
7B9D4120E3A9C87B8387574FD8BD3170C18BE1E33B7279DC84AA9DD9A841A7473E2DCC5DC49272023A70BBEC6DE67E7E9E7DD495F057F4DD1DCA9438
4DD19972A8E265A1EEEDC8FE8FCD667D0D8C18F43466AC576AB4B51A8DC9A0FA9C95F7FCA509ADAF5569D0E2F5338BDF7AF5FD8C93E03301007EDA12
80D9949477D7B7D3D6021C854D64F6C82D6AB86675C2BD09C9E767DD895CF0227A706B8E98139D1BF61C35198C70E9F20F47393291B587EFBCB58A91
77DC46648C998F668D4B454EE42157BE42B9FBDE21CF9AFBDFD1C12A8D812669CAF86CC7C51A08D5E96A3EBF4CF95850ADA92EAED1E8314DDCD868C3
97B6BD9EC2BF677FB0BB1700F8D1FA37560074F9D994FB11211B08D8001970E5AE6365D0E999B911E34AD7DCBFBBEEF0FD9E83B81CE6BA07766F10A3
D600E9629F9D6573851CD9C42143EDF64FF23EB7CEE961DDC8ABF8B143CA9B998682152795D17E33F3B0BAC2379F31584F6B4FADBA588899D406FDE7
D8EB5F6B35792ABD06C5F5997137D4F4ABF6E1F7EBCCBF5DF90B4EF703003F25005064DAA18BF0D35E0994C958AF47910FC93595471FBF9CF871C784
9723FB3D1878B83793B138A5477EE33501262C39ED8D4C20E5755AD2A2C3D8BB4DFADEE8B05EB36C2D96393EA33A05AE9CB3A56C6EBF7B3082A8523F
EA1ACC26BF3EAB6AF826B0DEA85556E8E1029DBEECF5B3531BF7BFD2D2E816FF898915C4DF9C07C77B00809FE03F4511EA8F7363E98C512F08A4F2AB
0642EBDFE9935696AF98F9B2FDC2E1CDA2E7F7DBC4618CAD38AD830D88D180DC3F9EE22410381C6CE5EAB9639BF5F2413EC52F6794C173F6A85E42DA
151B92676EA94611B8E8C147BD19A3CCE5B50DDF435FADD5D76A74A69A6A7551CA8B84ABA7AE559BCD49FDC336BFAC40CCBF0FF603F91F00F829F9DF
4CA0FA87FDEE9B73D71C2B557D4CAAA887B25F56EF898E773DB076F896162B1EFA048918A34AA26B20036442C8A76B2BBA0B84F3F6C97CBA9CF2F73D
D9AAE7B9195FA1DD034ACB54F0BAB177A625C0A8A1F6F18342DC509AA9D4AB491A2D7A9465C48C4695A6B2E0F6D1EBD5261427B1D79383D73DC9A835
E17FD71EF80F00FC28E97FFB9F6F6F00491285AA76B67B43E813538C70FABD94F7FA7DC9A973BFAE1A1ADB29BE87C3FD112C0923206DE41DC4A8878C
F88B515FA6F284830315ED06ADB29D39ACC3DCD66B35D743EE56BCA58EB84F9AF1043718DF5E495557D4251F7DAED698606DF2B5421C3318F575CABA
EAE23A935185D4DEEC2B0FBA9256596F44083398EE0700FCF000F0C7FA9FA4700CCE9D187C0F43610831E597277ECE8B573F8C7F7BE641BB29A75DCF
5DE03BF0FCCBA6474126036A244AE26B37F185523B8F3E73BA77DE133265DC94AF79BDB618F2F5A7FD43BBDC23F1FA5B57EBF4E58519D71F98883AC4
0C17A968BAE1EF8534750DD100C7B0F294E5F692810F4A543082FF76E29F06FE03003FB4E6A77FBB07A8C17E0C235093EA59B766F751A311818B2B74
E74D0FCE21B7D6EDCB39EA7EB165D35703595CC7F7DB067CD56A1113A14F2B3BC8178BEC83BB0DF65DD1B547BF9034D5A21E39A5B52F43825ADDA3D4
9AA45B46E47505FAE9419DAE12F9566250A8CEA881211821E157F36784765D792CB14883138D50BF9FFAA7400000007E84FDF4DF47697DCBFD248EE1
98495F96D8A97906A247D0CA543CE7957E5F0E34B2D931E3E46D5B85EBAED8CA5C9E3E1FF204D241305D3068DB436BB1CC2DA89BD7D035CE837AC7D4
15F6D85FADFD30755ED87A33AC2DBFF7B53CAB0EF91A9B02571948A2BE56A5D71A0C64C3F7A3903BC1CE8A807E07B394FA86CC4F36F41DD4DFBCA740
050000FCA0ECFF7DB82E896A8ABFD61BB4F53AB5C1A82DFD70A669BB5A44A3475FA49A2FBF2E3A5C74DF3DE04DCAB0534EF20DAD78DCBD799DA21023
04996B07CF4AB1E1291CAD0704B41FD5C9B97FAE32A27D5AC5B3F5C3A745E1EA72FCCE33DDE76C4C959AA6A956C2259F4AD14A0D46D3646EC282017D
C6F8F59C73FABDCA84E244C377FF7DD65FE31624F36FA3BFC0492000E007F84FE214AA2CBF1BFFF0D285BBF71EDD7D915C549EF562BD6C9409D2E286
F33894807EDD9715D374B472D5E6811C3727A962B77ED156D48442247664477A53BE936D93D96DBB460F0C8C55DF755959F971F190E56352293D9AF6
554FE42BD1FA6A54A5C2D479F5BA3A98A6CBB687DBF2050A8F889389992A1346907F9BF34B7DCBFCE4B7EB0840FE07007E84FF0D991FD2C2E9F792AF
94561ECE6CFC6DE593CBE7EFDD8F1D2F9803C34AB2F4345E136F2ED959BACFE9E68DCEEB7C396E72E63A64C73AD808A3B8F9FA862437BEADC07E5AD3
E0EB2DBA966476E894FA70C6AAC9361B7173417A366D7A9CA6AA4690CA32534DB506D2D2E4DB49CE22895DC7902E6151455AA8B1EF27A93F647F8AFC
EE3F680000801FE33FAA3715A42725D7D6D7649E8B4FCBAD3610B83E61FD92CD5BFA4A2351D444A75EACB8718BFE1457BE6C5ED1C2F5912E3E8ED29D
CAF05504642270FA58E7775D85B6B6EDC6B71CB9B6ED9494A5CEBB8B972F38DE72929E7EFBA80EC12EC568513D5EFCDC8894D71B30D3EB0142A1B57B
877ED33A8E7EA4692CFC89BFCFF837FF76139999066F0000801FE43F45E97575A939BAF26A4DDEE777174F1E59B7EE59B1515F70FDE0EC19031CD752
2842579E2C8ADA47DF8E4F1F9D98D0FE641FAE84BDE65DD3ED3884E2281DE39F3A8523759DDAD66B4B2BA7E507BD863C3D36FAA27B977A3A35B6C88C
9FBD8F965792B9670BF4E52AB82CA62B8F1B38F1C4D2409FF1D7CA600AC7FF50FA7FCFFEE6866EC0FC0FF712838F0800F85FFA6F36E9730A2A6ADE27
E7171616BD7D753DF642A5B1BC104691BC33C37AF8C6D0144D572427F7BD0C5DFFB8675A6678E4161B5BF6D8D7DD1F50261CC7E92743528E09044DE6
79773FE9D3F4E4AE805D55A3974CF4F84243A7AA69EAD97DCAF41932DCCBC2AA0CAA2BADD8EC36F1E987FBF8F6BF58A623499C20FFC1FFC61040E848
8AFEA3FFE0130200FE970180A260B8AC3EBBA4A2B416AAD755E7BEB8B462D903537E72BA164553C7B5F64DA149843E73E95EDF2FB57B4E0E4EDCDAE4
50273E7770FEA64B0614C370A2706BD269364BD1DE76CA78D7D6BBC6063FBFDA7289551C81ECBE6EA65F5D27ABDFAB54579E68F5FAF4D11226D77DCB
2817E7F6ABF3F4248637B4FE0DD5FFDF56FE1B1A7F9AC2F28A4992FEF32D80201C0000FFED0200D114D4E85098C430ACE13F63E9FBF39396A66B5EDE
BA5388954DF7092DA761AAF474E5BA08EDC798152B933ACE18C0134CAE5B36B0B22183633879ED51BC90A75074D8E410362B3C6C4DCD9811EDE61AF0
0713D574C1351DF5B98C2A48A9516B0F3567BB773D71A05913B7818F6A49026D5CF7FBAD00687CF5DF1004CCB8322B2DA79A00F91F00F881FEE348F9
7B036EC61BB23181630D5FC692676B43E25124E7692264886FD5BB92D6D06FB69C5CF45E1FB7AFCBC583AECD39ECCEC5B35B9510184AA0E6B815D779
2C3E7BF03CFE94ADB63EF1F7BCC35B17A05FC3EE53D93145A62F5FD479AF7394C5B3ACBACE4FCA3DDDB2D9D4E81CA421C810048521288A35EEFAFB26
3F89ABEF9C3A9D56558B82F53F00E0C7F94F63EA3C2DD5B814D708413424755495BED57FBE89202A3FE6D3E531E7CA61B33EF1C181314FEE27CE1AF1
A6BD4C21F2FF7C60BC92C4100CA64F4D4F73E073BC56DAD81C1DE9D8A760768F9EA711D3822550E1AED7E6DA0C13595C81D64E0FDD9D515DBEA449F0
9E37DAC65B3C29D2A82E2BD3EAEA351AA34EAFAE2FFAF22CE1F0FE3BD975451AD2FC87FC0F3E1E00E07F9CFF3195DAFCAD0BFF3D02A038549B15D72A
3C97A2910A6D314D146555110545379A065F3F75CC31610BD356E2F02A66621105E3B0898C9E5C15C8664E5CC1E8B1C3D679D7D7209FEEE5F0A9C145
057313341F9255BAF21A4A3967EE876AFDFB91AEFDF67E5253188121CABCE7CF6E5E49387EE0D28D7B37EF5FBF75F5D4CA4D09EF92CB546ACCFC07FF
41BF0F00FC8FFDA7281CFF7DFDFDB70880636875F6ED119EC7091AA7CBAE9174D9D39AEA472FFBBA0D3E921232F2BDBB48DEBEE676442E65C4201375
6B91BA1F8B35AEB562D350D7F6CF767BFA1D307CEA7E5317B1BAAEE4690556A931574C58AF842ACE766D3DE64CBA51AFAB2E7D1EB762F3A6CD9BD606
869FC8282ACE283012A43EBFB4A4A844877D7BFBF7B70000FC0700FEC7E9FFDBA61BEA0F1180683C0284287393E6590F2EA709BA3A29852A49D43F79
F366ACEBB4D45DBC3B2B99D2A0F731239231238498C887B32BC733D9D6C276535BCA076774B3EB5359307055CDCE295F0B8F275557A8E9B4905508FC
6166AF810BB67FAA5165BED83EA0E3867B3A8AA60B1E413405A37A93FE737CEC95A7A99975DFEEFCFE83FF407F00E07FEDFFDFF84308C0714C5DF876
A593EF8B861C4CC6DE228ACA4B2EE45D9ADE6F595AEBB17B987CC1B4F38B522193093391492BEB17B2B85C8E7C60E780F3C79C3BBC341E1C9C913DEC
B2F6D16DA552479E73598EE9F606872E1D392931FFD5F56D23FDE67D221BBE2F9A7B9BA0214AFFB1207D5578D8CCC59B4FBE2A453092FAC3F81FE03F
00F0A3FD27A9C61701188A1AABDF2EF3B03BDA100090C3A5F81B43F281356BE3266D3FDC71264FCC9DF575C1052DA23719A9A73333E6B38442EBD0B1
E2C1CFC3EC0EE8AE753E9D3B6949E5FDA39F2AB5D416E936BA627CF301F356F49F75F9C8E499ED3BDF33A118A17BF4B6A624BF2A2B23ADB4F8FCC9DB
0FAFEDD9B060FB855C23499ACD60FE1F00F013FCFF4309D01801105493B1A9A9EDC28600607C8B5616175D8EED323175F2A2CDEB9BF078EBE1E3670D
A81145C81BE33EEE96B079A16B86B7D91EE936B9A2A0CBB09C47231FD55D7C57AD23E25C63E98F9D9CC6AF9C3363F89471E3F7A79DCD420C1A2CED44
0941D33A38BFA80ED7149A28BC41FBAAC4657B1F56A1DFCE0100FF01801F1300E87F1300BEBD07C461589317E5697580C6CC3545F0EB2F0F2ACF368D
791638E3940D83D5AD20FB62266C40613C65C9CD682993EB37D479C85EDF80E7FAF94EC7DFF759901D7FB3B8067D39F325FDA4B9F3DC8903E6CEEC3E
363042A52931696A8DE887529A3651344590BA8A9AE2BCFAB27AB589208D77238FBED35060F30F00F0F322C0B7296004DEB813C0A4CADBD7DC3A8B46
A90A38E77D7A56FDEE6E297BFD8FB5624B3CD2F3A7DF44211C42C98FF1973DB85C81C867D766AF0DD0159B8185E7DBDDC988CD5199DEECA8C68F38BB
ED5ED9A1E3B6B12DC277169765556B0D0693112511BCA1C2C071BCC604E715951755BC7DFC564F90AA2BCBE3B211E03F00F0135701C806315104C551
1C297FB3AB6540154DD557A3590F8E15960677CD0DEDD58FC114DD79EE7F0B33A024463DDD96E0CCE47318214783BBA5558548A353FC967C5CF1A816
7E7344852D97B8DD8DF16F357356CB88D70459A5C64C4A585BA783608C201054AF46755FEA719D0E351A0B32DF6B49D2F034FEB91A4CFF00007E520C
68D41F235054A733A1A8BEB6F8EE6EE7803C1A5521C6F7873ED42E536CBDDC699A9CA7482C9EFA9C3241046E9A1591DC962DE3CB96ECF2DD02474BDA
DDEDD3AAE0E0F672E8F31E95693EAFF9B307411E5D77AD7B4DD13A8D01D3AA357A65A50AC7500846709346ABCDFA5AA96EDCF34F43B59A84AD65245E
FC464951C07F00E027F8FF5BFA3769D5865A0DA2CD2DCE3D31D6D9FD0BADAB8454D756A77F0E9724AFF7F563B1A2B14BCF1AEA7808AD9CB0AB68204B
C2B6393AC8715F713BEBA895E29833130AA1E2EDD9F062B6677A7473D77671A5088D6A7504642AABAE2C2FAE524346042149334E406849E2DDA79B36
C7E618699A822E744FC0F0D70F35DF2B00F0A100003F360234EE0426090243F51A557D8DA1B4447563C6C360AFAFE63AA3AE62FFFEDA0F3E53733B7B
580936654CBBAC456918D71CB8A05CC2967282CEB56A7AE8A65DAF2FDDE6164E4B8234DB9E2027D8B27B175D5DFADCA3081CD26046FDD7AC32A8A2B8
06D3180912AF7F7A2A3272F4981E1DB6A42E91F3ACFD4ED435FC1B2AD6AF317D5DF21A2781FF00C00FF7FF9BFE8DAFFF28B2A13B2F2ACC7A939BB9F2
C60341602DAAACD0A5AE3C509D283B7D5E622F9A95D4F58CBA9E4421E4CC900FFB5922D690D936FDD70C71DAB2DDF3DAF23D1AE3F1CDC6543B6654A2
8DD5D822BA14820D159545AF3F547C7A99A956D6D7289F6D1B1EE62E164A2462B14C303676F3E05632B1687976C33FE2F6552A7D4F0569A680FF00C0
CFF0FFFB390033891BCB9FBF4ACB8DEEB7B90D7F0609A975054F37AE554FE9983150288A40F6DED2A83112C65FB749B826143227F4F6DDBFA0D9A097
FDC6DD9892093F1E59A619CEEC77C9D17A458DB9B6D090979E539C5F565BF9EA9D12526BF4C557B7848676EA327AF9C64307634F9E3D73E3FA95836B
829D85B3CA683A2BCD9C1CA302F91F00F809FE7F2B00C8DFB70114DCBAFC32236AC07C07C163B246FDE5EDD376279214636E3B72C7C3CF2EA4EB28C8
84E90FC4DCE333186D9B856D99D86CCE8166BB86EF35A5F7BBA0BACAF539EEEB734A47657CAAF9985A5E4F28613DD2D05B6871D399899326741F30E7
7589568FC038643418D57A95E6EBB5A95E8AC3105D5F607EFA1426C0F91F00E087F7FFDF8E037D3F14642651B8FE76F4A93BEB164F1105998CDA9ACB
F7227B7D582C5D36963786CE38FAA58EC671187BFFE88D2B8321F259BBB953CF5B4B061EEDFDA16A758C26716AD3C8C0EEC968EDCBE7A5699935B4C9
A034D2660445A9DC110AABD6CBEE6414A3DFD6FDE9DF761C5018A67A1AE1D9E212856BB02BD9A899021F0900F0A3FDFF3E8BFB5B002049D8A0CD3D75
E0E0C06923B933514D5DCEB6DBF3228EDA483B885A191E4EFC6CC48C188A245EFCE0CF60305CB6CC944C7BDA64E9806D95FB9FEA576F9FD6DF655009
5C723436FFF92B549B5FA0A46918C368BAAEBF6D8FB3E9658D3347A8BFDF3E44346E39C420F5BDEE4DDBBCA7B1F757EA290AA47F00E0C7FAFF4F3B81
71CCA8AFBB7121BA67CF404EB4A91EBD7CEE46E4E1CE4207A9A7FAC59C1C94D263387AEA906A4483FFDD7674775D33DC6FF984D7091FE0E86E47DAF3
FA5720F96B7717A5BED7D4BD78A6A71108A6E9DC883603B79E293052288E350EFEA2BF971CDFE78012046C7AB3B37DD3444A77AB0807FE03003FC9FF
EFFD3F4E90A84EAFBAF6606F4853858F46A9FE147166E7FC482BB195F0F99B55292584D16842AE45D64C6630781BA7F186CC751B1D767DEF0BF375D9
C836EC619570DD912D059FAF1B6BEFBCAEC3092D45D7CEB371EC75F84B891123BE5DFEF1BDD230FF7D020149982A1FCE71FC447EB8AF062B8000C00F
0D00F43F6C036E0C00188A69AAD3AE3E79B0D846FC1E3354EC8EB9BFEF4C7FB98DE04EF2DCD43A0C8521AC6475C11A0683BB3ED475C5B49E5DD79FD9
4E2983DD8398239468EDE5D37995EFAB4BDE67E9C81A1D4D9CF07368BBFC5D55BD06FF56F093D43F9C3BF8BDE93055AC6EAFC3124A80FE00C0CFAA00
BE2DFFE33889E39032FFEE9AB3B53BB9CE2AA872F3D89BB736CC9109387B4E873DD019F50DFEEB573F39CC60B29BB4081B35A243A7C7D332E908DB50
71BB2A547BF9469E2623ADEC6D2E5A576DA0D38758B79872B3C460D419B0C6E1FF8D937FCDE67F3D816CA6749B4661C95938F01F00F8B1FEFFED3860
83FF8856AD34E8EB2BBF141626446CBCE5C3598DC0CF266DBD96B8CB45C4197F7841190A1B091CC6136E3EE171383C9F811B86071D8B5E43A7B9B4ED
EDFF92546FB96C28533EFFF4A542A7AB85A90B5E8E3D0F7C541B4D260823BFDFFDF13DE79BFF79E581A4A1C8ADEA17C07F00E007DB6FFEC33E604C57
97AF34284B4AD29253A3A7AF582F64DD406BD6F59ABDE941470EDBE7D9670C468C2806E14F4EA73BB1B87CC1883D2DC29F4F7A470E7768D764B3B97E
5EAC46A9BAF0B8A64053A7C60DABA54E0B9ED71A1104314126B8F1FCDF6FFEFF6B0DD01013B4D3A3EE22E0052000F07322C1B7B78018AC2C2CD71A51
BD36E34ED4DA53F3B8437458D6E8019D6EAE66B39C73F21E55A2108C63F8DDA8EAC11C99B0F981493E51F1338854A78EED43CAD08818285B7BED4A9D
B2B24E83E586F3DDB7E5284D8DC3BF481C42D47A8A24CCFF46FDDF8FFE57EEBEA306F91F00F849FE377E913866A8AFCE2E3352545DC2D6E907D75DAC
C9D43E19D67EE8353E53967EB2C365ADC18898E03363ABD672A5FC882BCD836FCC3B44CCB70E73E86288BDA37F69BC74A0B622AFD688DE70E3BB6C4C
2AAF31429849955D545B5D598F3476FFBF2DFFFF9B0A804C8E2BA3C0064000E0E70480461ADFC9E3A83A2FA500351B52A2229615C38FF6E6C5F6F3DC
EECEE07E383B274BAB860C263C6976D579A1B560CF5AE74589BD5375017E6D87BF4C398FDDD4DEDE515E925365440E71F9419B131EBCBDFCF2F5C7FB
374E67971BD48D4700C9DFD7FDFF35FD9B69DDD90F04D8010C00FC04FD1BD701BE89499A898686BDFCE6E14CCA78E3E8F07918BA3B227DBD475B7706
E3FA9D1D26588B680D947AFFC7EB52317F7A5F87ED1B07D4DC766C3FE03CB9BEE263E9F58DD5596FCB11740387DF6F55E49055B1A76F9E583974695C
59A5DE506682F0C60EE0CF5A008ABA170781020000F809D9BF3102507F7B0D00194B12365E818BCE5E0DBF4DABC66C7B31DCCA8EC5389B34F25CA116
C161BC60F7D74F8E12B18B6FC8E1DE11DAE9F2AEF9CA55A9E8DDD8AD559F1F9722D042B655F776DEBD623EA96075F297BBF1F74A209DB11C31E158E3
C5437FD2FF5374D17103F01F00F8E1B99F36FFBD30FF360E08D696BF1A7F087B1137B0691E7DBBD5D9C3362C0E6369DEF849D7949009253F8D7E9D6D
C717CB037AAC71D858D65C3EA8F2C93EC3B9B537739F3C2C818D53397C85EBC06B15A891A89A33F5D0B93B292A1845D4E52A0C6ABCEFEB4F020065AE
3D09760002003FCBFF3FEC03C251656EC2E80B68659CF56813B460FCA5415C39639AEE7042663D6EC429DDBA3BF5411285D7C425739D6FDF97F84F5E
FEA97AD3F2D759312FAB20FD4C362F64591684EA203869FA9C5309CF2A8C0446E068459912C3FE34FF37A0BF5F07D23F00F083DDA7FFB51527295365
EACEC08718BC5270967ADC26E2208FC5E852B6E6AA51879088997E19078FE2099DD66CF7ED92DF4BD029E67CF5BEC119B9910F9550D570AED7510386
298DC8D75903965D4A293551C4F75B86B5B52846FEA9FD6612CA3280134000C0CFCCFEBFEFC9C575EF1E2DEA1247A243830DC65DFDAECF6230BAAAB7
C4D4AB4D9499A44BE3CD514299F5AE2DF6FB3E08EC36956BF7347D56B3F4BACE9813CA1AF08140CB4C68C1EA96C3132F3EA8C171A271E36FC32FB8C9
8052FFB4EAFFC7FA1F2BD481F21F00F891F67FFFFAC746BCD17FD4947FEFEA68C72D18D43B9E78DB614C962DA363DDCE85372B309442E9D2459A18BE
5430BA9B67EA745E935750A27C7BC9E2DD6A38A3B7CB8E52E24B3EFA35D2CD67F996E389C5260C27BEEDFD6D0C010DF5FF9FE77F8A504160081800F0
136A807F598C238CBAF7AF627B48C6C165CB6ACB76B638D09AD113BBBD605F098D93045DBFB8E23A5F22F27098FBD24E3151FFCE33BC60CE9C72DD87
0B719FE0F2B4DAF28D4D9CC7AE18DCE1E07BEDDFEC27BF4D19A4FE506FFC73DC211B8700834F0300F889F5FF6FA7F349BDAEECE5F3554D851389E71B
55773BF8BB32DD4CE766C617D14443170F2D4F7A6FC3B711C96337F2DD5E66B41DF56142D74FCA9C8B05484D6176C59B492D7A8C3F366BC89ED4BA6F
978BFF3E5EF09F17FFE87FEC02280A4C0004007E4E1BF0CF37822246ACE6F3C5A1F6D2A3E8C95B99D303796CF683C4B5351A826EC8E8E6E83BF5217C
5BEB90333E820D65FEC165175B24EADEC7E7AB4B3FE6953C181C3A73EDE6478FD3F275D4DF077D98A9BFBDF9FBADEDF8430B026E000400FE1A1DC0B7
317C6612336A4CF5EFF707899DF355136EAE1FEACA6146646FFB580EE3661433C71FD28DE0C89BEF99C30E2B5DC83E991FBCDF98119BA5CB2D28551F
0B0A1CBBE95A8E4A03617F3FE0FF47FDFF2900D0BFEF3F3283F40F00FCB418F0BBFE587DB5526DC071535DBDFAF9B2D6B6CBF5FB46456F68C7E72CA4
EEC42711B8B9E12B69AF761E47E1BFD5DBFA6A96F5B067A36657171F7D6AAC2AD56A23ADBD7BAD3F5FA54711E2CFB6FA01C90180BF9AFDBF9D022629
43E9D7FCF48CC2128D4A5DFE316EBCB4795A41A765479D789CBEC4ED0987EB289440C8ECED05CB1912B9ABA46FF54AE1A1FB5E972AD71FAF492F3569
8689FD6647BFFE505D0F61C49FBDEC074F1B00F86BF6000DD53AACD6205A654D41FED7DA8AF40BB746385D81E2A69CEF60256847A54CBB5C4F63188A
D79FFC72842B134B6491A6915D8E0E3970ADCFB6D2D232D230694669D9B1F3E945CA7A94F8775BFD6890FD0180BF7007F06DADAEF1A43EA42EC82BC9
399592D26A84AAA4EFF6A95CB14DC69DAEB1B5460437184D979EC43159F602EEC203C1B7F70D4E9ABBE0CBE71C8AD877D7AC5FB9F4799D0146499CA4
FEA4F807110000F8ABAE007C9B03DAB867975015A725C71DAEBFED7FCBB475F9212BB17DDA9BE16775061837EA91970F1FD8F07CDD9B6EF0F5BD75E8
E08CC8ECDC6C8CACF852BF376C59961656C18499FCD7FCFFFD9C310000F80BFBFF1B286A7C167769C251F8E2BABC73417B87099DCB0CE75F54E9611D
A4353E3A931922B16932F54CCB1133272D9CFBF2D36D1342272F6CE1B3FA130C9960EADBB09F7FD30180A70C00FC95FDFFED4AC0862F0455C6464576
CC43BEDC3B37BE4BB05496001FDBF319223590567FFF5469578690BF70D1940B637CC61DCE29AE80E9BCB6361D4E949AB0C67DFEFFF8C20F74FF00C0
2F5400340600D2A44372B7AF54AC45900FAB8FDB4AECF993EB66CE7F98AB55AB55BA0BABD3873324B2AD7D075CEADB3B6A6EDAEB3ABA3CC069D0E53C
4D69FACB8AC6DBC47FDBEEF787897F3428000080BF7E0160FEEEBFD1A435663C1BEAF31137EC7A3BC8CA833F075D74F6E6CE2F5A75BDEEFEF6FC535C
41BBF8D10B66F759B9E35C6E320D77726AB127E7E1A523899930F52743FE80FD00C05F5EFFEF0180D0AB2BD59F93A25D3A1792A78F6E165BF327502B
8EBF3A52AC376A4CC5F19FDE5933BA0F76DBD579EAD019598F4DF464ABA67D9F3C3E9CF04947D27F32E50354FF00C0AF11001AFCC74C551F5F257D59
EF3529B3E84C46376B7E5F62CD8077F926834165343C4B2A6AC37215FBEF19B2714274B1DABC51D2DC37EEE2FABB6A8AFED331FFC07F00E09708008D
FE9398292BF156EAD3B52D8667BEBE3DDE9E2F2E3AED71A3D4A4D66B8DBA93ABD2BAB2146D7AF59D356976720519C3B7B59AB37DC72DB5F93F4CF802
010000F8350A80C67D00040E7D4C8C7EFC616EF3E9CA076E0E0EBC57452147B20C5ABDD1A43B392F6F3EC7AE47F751A7462F294393A50AE1D0E93352
B57F6E3F701F00F8C596007043C1C96D8989611EA7AAFB8AEC983BF503A6EE2D32E8B45A63D6C9828D1C856C68FC6EBF889ADA162C66F0DABE293AEA
CF077C83B77F00C0AF540134D6000456717BCED6A5611DDEBDF69232A39059B3FB3F850C7A0822DFE7C408AC452BD346FAC6EB67B2595D0E04C519FE
43ED0F1E2D00F0D7F7FF9F22004960F8E3E98BA7361B9815CB638ED487F7E9721BD2E951147FF4E6BE48219AB0DCCEF3FE1B05CB3F2B748D9AFC0F8D
BFF9FBDF0E0000FEE221E07B14F8ED3C2081531746EE59E0DAF3AA2FB38B266CFCE96CC868C271CDA9E45C179995ABA3CDB4D21D2CE18BA8FE25FFA1
EA07E20300BF4E00F87B27409128513C7BC7B32E92DECD989EAAD16B2AAB4C880ED7A72C7D58D74A646DEB19B8F37377F6D0BB7E6924F9AF39FF7BE5
0FEC07007ECDB500AAF120CF911137B7DB356FC774D11E38587BE723AC8575C5BB13A0DE429B0EBDC73D7DED6A3FABF536186CF503002C2D0090A4BE
F8FDB21EE1DBE7B671678A8A6E2EBD3AF1821E2508FAD646ED18AEAC53F0DC94A936AD84C37329122CF8010096E6BFB97CFF8CB4E2C5613B0739B2F9
99E9BDC7CC498370334E67AED446F264AEED0FC4B87774697EDB04EC07002CAF0120B19229F7E9AAF04EBD6C588C2BC5DD66C5A9499319A6CA17971F
E488B93E1BFC1D43140B8A70A03F0060790180220D7727ED30C4D8F6F01032F7AA4645DFCB2530023543DBCA9E4BA422CFB6CE437C7DAF1B29A03F00
6089058099486A36B666A86353097352FDC893BBF7D71910038E47A51738B104768A807041F74CF48F1BFFC08303002CC57F8AC48CE75CF6BCF2F4B4
61F6558D3B7C6177B28A82496AEF6B5D7B96C4B6437813ABF5B508F9F7719FE0B9010096E3BF9920959B6CCE0E155AB1BC3E8F5A90B5ED782D558B61
2B634DFD994241EFAE02E9192DFCB7719FE0A9010016E4FFB71EE0FD58FFB136B642FB8C633DF6CC3C01D1146D9E17699CC2B612D8598B9B24C33046
7E3FF60F1E1A006051FE3714009426A669671FAE987DE753B7E54B4F666BCAB35F75691BE1C3E4B2586CBEF8587E61458509F80F005862FEA71A0A80
31EE2D8462F6C1BBEE7D46847BB7F370E031190C26472061B1B95CC71E5D661D4ED3355600E0A1010096E53F499AB1AFCB5D5BBBF0B92D0205029E48
E0E818183676DEB623172EDFBB71754B0F57573BB7265D7794800A0000B038FF299234969F0B737096888422A1FBD82D271FBDCCC8ABAC505556E664
3CBEBA3ED4AE456B7FA9D469E46304F80F005858FF4FA1B82E75630B89882BB6EF3DE8ACA636E5CEBE4DCB97CF5FB470F1D050673B2E8F2B113BFB7A
89E501C7F52000000096E43F455058D5A949B64251B3CE51D1B74E9C3EB1B88D338FCD60305A0C6BC3603245B66D5A074A83979C9A6F2F71DF6E00FE
030096E33F49E0C567874803A62DDC7579C7925E1D3C3D447C2E57ECD929F2F8DBBA943D5B97EE3D969A7270FADAC51DAE9D68297439838100000058
8AFF246E78304CC87358B372547B3791586025B0750D1F7BE465592D82217A04C3104CA97AB7EFC68A43616D236E848942B380FF008085F84F11AA98
405BA9C8AB3D4F2092C8AD25F65D76A4D4C354E325E03A43BDA15EA357E71912A7DF09DD3437D0754C9297FC08020200006021FEAB8F788B9D155C0E
976F6565EF1AD87BF0FD8FE5B97AB4D6A4316AD57A440F61A6E2BA074B1EF59A14DCBC99E8FA41C1382598F705005886FFF0691F8154CC667383DA07
B8CC1EB2283BFEC1962D8FDE43A8568318355A1452D716193E444C3BE11B2AE07A32E7E6DBB5C802FE030096E1FFABB6423987C3EF3E7252339F90C3
83BB5CCAFBB8747BD1D33AB5529B968EE47E885A3F79FCE54DAE0376B97909794DD8B36A9B3B2401FF01004BD09FD6CD118A39F20953C7B971B89D7C
5A2F6A1AF0EC7960CFC391E7266D4E5AB8F6C9DCE95D5686BBF7EF28F2EAC667F3B862E6B0EB52DBA7A0FF07002CC2FF2FAD2482CE17D6D9B2397245
4717E9DA410EC11DECDCDC02229CAD67CC6CE12277F7EC1CE8EFEF2EB272E4B1584C8EFFE10122EF74E03F006011E5FF79857070F546365B2AF51209
5AB71F32C6D3CADAD5DAB5D9D6F14E6EB60C064BCEE3F2EDE5CEF64E0EAD460D6AD9EFEE276F7E70A9990237FD0100BFBEFFF04CAEDBFB87728EA327
57247791F9DBB5F0EFD173F9B4B95DBB794B996CA690239548C5D67CB950C491F658B4ED71ED418974AE09DCF50B00FCFA50E6BACEACF08270B64C28
B5B563D9F10286AFDB7AEFF2DDA8B1AE1C069FD10893C5E30BF97CBE84CB6432591DC636E5B67ADD38080C040000E057F79FAE6CC7ED9735B56D4747
1B070F45D765098FDFC69F8CE8622BE3721ADC1737EDD87DCCF895918B7B3415F078EC8600C0E6F19DCFA114180306005882FF359D996E9F320F6D5E
3963DCC1B30FEF5F593936C0512C1609F84E2D26ADBB9D53A7C320334ED45EEBE7E922E272981C81DB318804FE030096D0FED3868512FECEF2391DF7
BCB8FFFCF4A22057B15C2497D8489B8E8D8EBE7E29B3E66572E6937B979FDC38D8B5895D937EC39BF25CF76AA8EFB7FE82060000F8C5FD37E3E7BDF8
CED7922FC7CEEBDD5ECC6189EC9D4412A1901FBC72848F1D27304C2469225588ECC5223E8FC3191EE1C91F534080DC0F00584A0150BEC051ECB2F664
0887CD6709E4EE0A119B231271C5566CC6F7E53F2E4F2810F185421197EBE5C16BF30005B53F006031FE9B7322DDA4E28EF206E979761EAE5C269BCD
E171D91CC1FF31FEEFFF784C26972B120B1AFF90C59670A56B6ABFCD00070F0E00B08C0040551C1DE06CA5E0F39A0C1B68C561B3D86C1697CB92C924
C2FFB319C2E7B4F3E34BF87C3197CB53B8F03D6FA2E0CD1F0060494B0094F6CD926E2E5CEF9847BD396C269BC16430591C9E5DBBF64CDFB56CFEE675
B681CE02514300F05070FC930990FC0100CB0A0058F1A3E122EB9933645C3EBBC17F1693C917387B7710D80C750EDCB0B765DF008584C7759AEACB0D
FE4851C07F00C0A2FCA728E24D1BBE5824E47138CCC6453F06436C2FEBDCC466E6F20D334F0CEAEA28B16131AD5774E204A79B29F0C400004BF2BF21
0090E4451F85B8657B5E4304E033F80C9ED03B60C1984E875F5E9BF9625D3F2F3B47B6883367183F04F80F00585C016026694394423A288ACF1133B8
2C069BC30EEC1DBC63F2B8BBCF973E8AEC126EE5CC65335B84735B7F00FEFFE53F4EF00800FF9F7E5EBEDF00F2AEB34DE8091781882BB3E533ACA4AD
A785C78FF6B99279F2CB96BE4BFC7B04F199CE7DAD9ABF06CB7F0080C5F94F53946EA975EBFD1E2C86B78744CEB06739AE5E766CC2D0472F16BE5F17
B0AB9DDB4C6766EBFDEED63781FF0080C5AD0034688DEFF368733C94C50EB4617058EE42D7F133D71D3A9E9279286347CF87A344DB06300357D95925
90C07F00C0F2FCA7A8BB2DC5D343999CA60C868421933BF68C1ABB66DBEEECE7FAA34171F3ADCECC628A43A4C21D38F01F00B0B400409B293ABDA7B4
5D6B96B0B58CEFC61459B98E8A1F10F5F969E13BE37ED9DCA32D0F0C625AAD08152C8240030000585E00A0A89259A2EEA3988C117D78028684E9D67B
C3903DFA54551D7AD26AD675C7911358BCDDD384C3EAC1FA3200607101A0C17FFD5EAB292739CCF61E0C06C396E3357CDB82CD6BE7AB49222D283ACB
2F3CAA096FC23069F77290FF01000BAC0048F49A63978B4D99BE0246E3EC3FD982E373E276AE51E1A6FA415B9FD9791EF062F59A20F2FA08FC07002C
30005054667BD71B2399E3BA33B86C6B5B87C5FB873DD45560284A1CDB722BB4D581704EB7ED76CECF81FF0080453600957DF8113D19332FDB33B8CD
DB064E9CD076CFC7673004E14F565E5EEB33621CDBE3A09FE43C38FF0300582094D9384BD4AD13D33FCA9D21EB1D2CE9B9A0D3B6A8494A12C5AB125F
5FEFD9B91BCF2D7682300A03FE030096B80080454A164532990A015FDCA995A4E5FEF6132FEDA8C58C84FED6E3134B47B594F0D62DE24D3482060000
B040FF2962837C57B484CB12C92443FAF90F3A3776D8A5D7086220B59316C41FDE3853C19F36881F56035E00020096A7BF9944B649671C91B11CDD78
CC8065411317AF99BB784F09A6C59011BD6EAF9BB551CCEDD699E7F401F80F0058A2FFF855C7F6317E3C0F2B05CF237AF8C8454B16C66D3F8E621879
78F4E9AB0B263A08BAAF56C8EE00FF01008BF4FF9163D3331D193EBE7CB164799FC01D4306BF7A37578FC254CDAB8427A33AFBB39DA3DDB9FBC10400
00C012FD275E7879250C64758B90382A025B07C7CE9BF40A3A968B6104657CF27CC188413CFB43C379AB51B0000800589EFF04F5A185C38E8EAC8E09
018E7251BBB6CB764D3A539D918FA144DD8163491B3B0FE5B257EFE60D5103FF01000BCCFF545E90F3BA361CAF2D1E0A7B69488B7977173FFE7AA50A
87F157DD8F6427F65A1D289AB8DE36B412F80F005862FD5F18203DB094EF30D24AEC26B1F31B7DBDC3A5BAB3B91884959C7C9DF32074D12C4EFB418E
9E1934580000002CCE7F8AAA0A911D7C1AE0B6A19DDCDDAEE5CA050FFCB7225FAB500C4355CF6E646D5A1225EA1113ACB80A5E00000096E87F4D2FC9
DCF35ED647C7286C140353B6ED0D58AA5269600C83E1CBE3122F874C75F47CB05A340F070D00006079FE9B3573E53D226C84633A59B93884BE5D347D
EE88C34AC4D490FF919C599BD33AF4ED24ECB0C1A6AF0A14000080C5F96FA6F42B651D2E771078394A9DC45607FBB4DF3770665583FF386AD0DEDE7D
CDD9AAB74CB0B8B57B0AC8FF008005FA8F6CB5ED913554E46465656FE5B4627CAFD3476F1B1018C771ACE6C5B9C41E6E8B3D841716CB4F9060011000
B0BCFE1FDD651FF4A433CF3EC4566E2599B3DC6FE7D742150AA1388695DCDC747C997476286FF351CEA02A30030000B0B8FC6FA68EB9F85C0FE6B96D
6F29B5758ADADE7EF7D96989188CE03852F92EEEEACD567316F2C6BC92FA3E22C1F302002C2EFF53A79BBA5F0FE7392EF495486CB60C763FB9AAF5E8
020441F0FACF856F4EC68574DB2152C434B3DA0F83150000C0E2FA7FF3759F80678B6C6C1CA5B636F6637BFAC75D89DA70054675185C565777E7E1C2
2E077D380B3B59CF01330000000BF4FFA1AFF3D9B5B60A91A88942DA24D86AE29D63E56945B001278C857722A2D6FBCD0E662F9A2EF2FD00F23F0060
71FED31981B6E77658DB4B252DECBD7A8D085C13D93E31310B21481369C8DCB16471C8A0C1BCB1F112FB1B04F01F00B0B8FCFF255016BB992DF2B177
94FB45ADE87F746FC4CCBD5F211D45D214AD8B5DD1C16FB44072A68560250CDE00020016E77F5A3B9B7BFB399E915D1C049E93C7840C7C517AE88496
841A2A03B43E71FDA6A96D87BAB336EE12F5AF020D00006061FE93C49396D2F83D2CCFAD1DA47C5BEBE0F9419B3F5D3F5D4AA26633F662F59CCD8B67
FAB66AC56AB95628BF0BFC07002CCEFFA7ADC4E70E701C265AB3EC9C795D9E8EF68938FD454B52344968EA1053EC8CAE21DD1481FB85E2E510F01F00
B030FFD1C72D44FBB7B26C9776E429AC59219FA739F6DB5555852234499B69F8D4F0D9A3DBB8BB89F7BA724340030000585CFE7F1D24BA7255205DB5
C981672FF48A090A38B6EDC91733DEF0C718F5BEBBFFD225033DFC84F3970BED1E01FF01000BF31F7FD35A7CF9B68065D75621F3B69686F6687B6BE1
A6249C6A48FF18519F7C60FE9851131C98AE6F87F2D68121000080A5E5FFF43061E43EB180C7125937B777B05F1C79F9DCC12C65E3CE0012430C0FC7
878FDC1922106FDD20EB0B1A0000C0B2020089A676E48FDBCA168A18323E5FE0CF9BB07C76FC99CFF14A1CA261B4F8FCC28D739A847565B2BD83240E
57CD2000000096E43F853E09166EBB2090D85AD92A6C04C3DB8C9EB0FDC9C9FC63D9244163E6CAC311BBE77889FB4A79ED3634B199A501010000B0A4
004022575A088FEE65F35B852BDC06F9776B392060C3A5815FB75F37904602CD5BD339B8EB38F7306796F0C866B1D325FCDB9161F0DC00008BD09FA2
900BDECD3FCDE35BB935E72982BD15ACBE43D7C4073FDD795E0F9BF486CAAB3B3775EB662B134838EEC7FB49FB7C2280FE0080A5F84F53247CCAAD6D
6A5F8E935C602FE23B89845D5A8E4A6AB9E6F213C26436990A3E7FBEE0CEE708054E325EFB0562EB5D4A8C348383C00080A5F80F1D77B78FEAC575B5
11B5B196F7B2E6F40F1E5D3E7AD2BB8FA6AAEA5757660DEAE8CDE230F92C9E506035C35D10F8062540010000584CFFAFDBDD843B6F20476E2FF65528
0635977609F35CD0CA26C02BC8AB858DB3FBA0FE4BF79E8959E9ED1CE02A096F23E44FAE204910000000CB80A234DB1DF9872771AC9CC48E2E726B67
894CCAE50A6DBB6F48B8F8F05A728EE173BE892091A8619B577469E12812D8ADAA060500006031FE572D97F1A3BAB32432A15022E50B7856DEE3E7DF
2DA32812D6BDBA7B767F0FFFE15B2226B56AE2D4C446CC170B05363BCA400100005846FD4F5125E3C4CD639A73E572A995C44A1A1499574618A0D277
47E7B77265F1195C069BC96509C436E206F779129E80EFB4E43306020000F0EBDB4F9B49EA636BDEC4FD629E80C7910BFDD6977CBE3CB1671B1F8580
C1600819AC865FD94C3E5F2652489DED9AD84A596C2E87DFFF1A98050C005882FF04795D26ECEDCE60CB037C8787B5EED6C75EC4E3F3247C0987DDE0
BE4DC7EEA103068EE8191EDCCEBBA9C43770D2B2984D9DB85CEF3308080000C02F1F00CCA47E015BC0627AAF3F7D716A589386FE9F6F2555F0F91C3E
DFD3A3FBB92CAD418D200404D555665C9ADEC62E78CCF6530F56CA05AD9381FF00C02FAF3F45BD6FC66B3E386CFDC1E19E5C3E4F26B1B2920BA53E83
46AD3D9699918B62288A93A819A71A477F6A2FAD9FE4CAE3F89D3CEB22998D80A70700FCDAF637F88F5DB2E3755CD7CF99C56868EDB97CD9A4E9F6BD
575DABFAA4AF426A6A10C800C3288A6911C2A832163F4FD818C455F06CAF27C85DD3C1264000E057F79F34C5D8F3046C2683CB603085A15D3D1E7E9C
F73AFF46E5E3E2AA9A3D29080421188943307A6EF0E5EAA5B3B6F4B1155BD9F0825E75959D05FE0300BFB2FD8D2FFFF09C45121E97C962B2641D8207
3FBAED169D32E6C1A7D8D343663E7C3C351D43112833FB750EF5DE57F1E06ED70EC33BC9246C2BA9685B37D961E03F00F08B57FF94FE6C3B318FCDB2
F1EDD4B2BB9FF7B84793D66CB01D92903E41286B36665D56E1C7A32B42FB759DB5B3ABB5D3E84E2DDB0685F8B872449C81B71CE471C07F00E057B6BF
417FEDAD3031D731B4FD801EAD1A1A00E1A87787264FF4761B9D3C842F66BA05877AFBB0055227B14C2AB6B7B2B3963849DD7C9BB9DB76FD3082EF92
4C83DB8000805F3AFB9B6EF5B3725C72ECF0B4566C1E83CD100C5AD93768C8502FBFA9AD843672279E1D4F28B5B2B1174AE5F60E36622B895CEC2C92
F488B934842F9DA806F91F00F895FD27A9EC09F676DB5366B77014B1586C2183C995597B3A8478DA797B894542B9402E92B1851C8995B344E8E4EBE8
E46AC7E30B38A28DC13CB6D77DF0FE1F00F885FDA7280239EB61BF326DAC98C3B6F7E430592C2147266A69EFE6E4D6C4D5BE89C0DAAF55F340470F17
27B99D8DA38B8F95C8CACD5BCE97B11DED85B2481D48FF00C0AFEC3F49150D12F8BF5ACEE74A581E03D94C26A3698F7181AE62AED051E4D4C2DFDEAF
EF84F61E2EADBDED650ACFA612195FC057F879F3E43257B1685C1648FF00C02F1D00082AC557D1658B48E8D09423080DB255F82D9B35DE4BC8908865
3207AF165E016EDE123ECFD67A40A44B58905CC653F0C4625BA9402C759B9E4901FF01805FDB7FF32B6759AF053C99B5AD98DDE55264CB66A1B6229E
80C1E54985325B6B07378542E12EB3566C7AE710E62770E30A04626178A840E4BDA98004D53F00F06B070092FAE026EB70BBA5D8CA91E3B17EB9B790
C3E6F0D9120E9B2B93B93837F1740B6CE5D054D6DC6DFD5B8751C36D03B902B1CC2EA6ABB0C9DA7C30021800F8D5FDA728ED6891ED8D970B064D0A1E
D64FC2E5F3C46C0E93CDE4B2398A16F2A06E928E7D7CE7F5DED06FFA1B87E14B65CDD9222E67E2368170623E498109C000C0AFDE0090E6EB4EA2DE49
AF369CDBEC29904AA41C0183C71132053CA14D4FAFE1EBA4E3674CAF7FAB3C3633BDED8C130EAD5802BEFFD526FCF6C9600030006001FE53947E954C
E2D7AE59F4B486DA5EC0E7705DDA34E74BAC9932F120BFA0A34E53D6EDA475E8BB88D75D469FB5F5E0B0D9510BD9D647218204F53F00600101C05CBD
CC576E25E81C22918AE54E8EAE9E43DBB26C5C9C1DF85E7EDEBB07AD3C740DAF81F50792268EBCE6EC6DC36C765AC61D5D8903FD01004BF0BF210028
CFF7F5154BEDAD85EE53970FF097BA06C898E24009C7D1573179F9CEAD377095CA78F8C1ECFE773DDCC5CC8163588E179106FF29E03F006011250096
79625C90ADBDC465CAD11142BE3CBC25D7A98D83B0F50097B16B634E5EC35475C4BD4FFB267C191B6CCD0A73E074C921488A02050000602111802634
6FB7F6B1958ADBCD1AE82669E1E72FB5B30F683FC2BDEFA863270EE228865FD832A1D9B476C1768CA61CC9CA1AE49BFFE0C1010096120168326B510B
99C837FE948FA8450F8EC469F538BF96FD7AEC8FDDF0E8E9CD0DFECC468482E6CD39A23D2A8C02FE03009614011A84561E0952883ABE19C97269C5E4
C84705F23C423ACE88E8D6D6DB95C366B21BFC67335D9DB9DC43EA86FC0FCA7F00C0C22280EE689044BA6603BB65284F24140BB932858753331B315F
22B4627E87CD93487B3D4370E03F0060615054E5666771AB5B83DB78B35802A688C7654AF98E6E23373D4E7E7D7049971E1DBBB796B2383C6EEB8708
01EA7F00C0D22A002A6398B5F5DD438E021E47CA6431ADDA44ECBAF94CD7181A30B8A0F0594AC6B401439B4904C16F48E03F006069FEE3CAADAE922E
415C8E80C562DB0DB9F9B5AA0A45CDB419C6549AF28AE5C35EAEEB362562B283609212E80F00585A034062973AC82502A15466D56E5EC425E3DBBE0B
5570E305A166C4585ABC69F09399721BEBF94BADDD9F375E06044200006049FE13C69BC1D662A158E6DB6A7D5C3BD783FD9887F0EC7B108D55A4DDDE
BA6478AF1901525B81E7216FABD338580004002CCC7FAA748397D43B32D0A54393B0C10E9C959B67DC88EBEB1879644CA087B390C7E6B339425B6BF9
386FF1CEC6CBBF410000002CA9FF27EF84C91C1EDF756F1228B3F1EBEB171CE42065B158421187CDE430598D7301B85CB6EF2C5BA72B38F01F00B030
FFA12DEE9215173D9C3C026D1C46CEB46232980C165B289488245C3E9BCBE4B0394C8F49D183E5FD4ACCC07F00C0B2F43717F59777BCEA2C75B5137A
44F4153399AC86BCCF62B2B82CD6B7FD7F8D845E9AA8F03C8E36EA0FFC07002CC87FFAADBFE2FC11791399BCF3A45E7C168BC3E636D8CF17C9247C0E
83C96A8C064CA12DD77B77B5192CFF01009696FF5FB94A9F6F144A42D6AF77E60B042C4E83F26C06D36A60579B0E1DED1BBAFFC60DC03C878DA5DFE7
7F80080000580E94B9B09BD584A44DC73F2D9149ADC4220E9FCB134AF87EDDE74CEB70685388AF4CC615B244A231E50809EC07002CAE0030AD90CA06
6C8E682B122B147C3E8FC5B31288DD970D58FB34B6FE50B33E9E6CA9D841C85E0143A0FC07002CCF7F2A67B8B3442696C8E4326B8948C091BA487D42
070F7E5B5E89664DDE112295CB47B560AF2430603F0060810180CC5DE3E5686D2BB177118B255C01D7C65AD075F890B2E472F4C588F3DD9A0A9873C7
3086E9C0E63F00C0322B00CDE5F12E36D24923C512894824F00D710999B3CD90AD465F0E7BB1B993337FD91E516845E3E67F00006061FA37DE074016
CFB195AE9927548844225E8B69AE81115145195A24B9D7FD235DFB324744CB5DD3CD2448FF0080450600927AEE6FB57A85D8DACA512AF2EF2DED1B1D
5F9F0FA1EA356FE35A0D97F63EE7E8F68A02FE030096A87FE37D00F95DE5F3B64AC41ED39CC40E5D5C3A4F8CAF2BD0A3A6E8478F7AB6B677D8E56E15
8781E91F00806506008AD64F552CBAE8240FDAEE6EEDDABFD7ECF9CFF51002A1EFDF15ADE9D64E113D46B60606FA030016EA3F655C221CFEA683C863
918B20A0F5D4E8983A8301D1A185CFDE2E1FD594B370BFA05B8D1914000080853600449C6D8FE27552F769CD9A745284CD38A7CACA86215C5BF16AF1
EA4E36835EB8FB7D060B000080C5FAFFD4D3EE41A25C1012E4E1DBCCBB434CE69E97901E31C0599163078A6457FB48CF1004D8FF0B0058A6FF54C978
C5849532896B73AF517BFA855CCFB9988F9A1094806F4D9CE5218EDC2E9B6B02F7FF020096EABF6A81A2C774BE40E06C37FD7447DEECEA37454618D3
4254F6A0500F9EE7590FCF423309FC07002CD17F3309EF731C7D442C94398803C2149C3E59179E9BB488093668D6B50B10891356D9DC2581FF008085
FA6F3C62DF65079F2DE4296C9BD8B3A4C7AEE7A05A0C3753F48BC13DED854B97F323217001080060A1F53F72C1ADD73147BE542CB373F3933BEEBD53
63AC866A35791FB6F83BD88AACAD247E574C600F30006091F99FC22ED8845D6FCBE44BD922EB96EEBE3D3B4F0CF3F671EBD5A99DA3AD5C2614D90858
4E87D5A00000002CB2FEA7AEDAF5C91ACB9448857682FFC7DE7780495565EB7E5FC7AAEA486A32488E620401515410110C18471DB398C63C3AE278CD
3961C6881113064472CEB1494DE71C2B579D3A71A793EAED755AEF9B70C3CCBC77EF75FAEEFF1B19A0ABCF395DD4BFD6BFE2CEF70F281E323E7FE0AD
9F0653D115EF7DF3EE99FD6EFFF6993EFD5F9145055040A00BF2DFB2BF2C3A7BE78C4C7F4E66202BA7EF4DB77CB563C56E0BBE1E2F23FA6D7903DEC2
4F644D2C13FC1710E87AF4772DF3836E7D2705B2323333B30A2E796EED67E5090C2DFF76DA5AFED32BC7752FC87D615FDF925582FF02025D30FCB7C8
53015F766666FE8CF9237A3EF9CEC3AF1EA8C7BA6901FF3FFADD050585DD32EE6A1FE3FBC1155B400404BA1AFF1DC72C9D96E7EBD9ADF0EADDEF77CF
BCF0AD85EF4650F0B699B36EB8F0A21963064C185CE2CBFE745F51F7B522FF2720D0F5F8EF369C9B37FAD2B101DFE4BB4ECBECFECC772F5DBFE4B35B
72737D3E5FC017F0FBF20B7C99174ECB1EBB2F2DFCBF804017E3BFE3586F741FBF657E765656665E717EE1695386E7F7F5FB738B0BBB15F6EC5ED2A3
474EC1B0BC8CAC9C6B5A44FC2F20D095E8EFEDFF7293D7F75CB4B624272B37CB0727FE401630BB7BDF63FA1D3B7870517E6E6E46467E20D337F03DC5
11FC1710E862F2DF76CB268EDEFF608EAF28AF5B3738F53B2BB74FFFF33E5CF2C58A4DABD7BC7DEDF87E9919DD2EED1338F1736C09FE0B087431F96F
A70F9E3078C3E19BE75FD63B37273333C797993971EC83C96D07D269C7B1E3FB5F3E332FE3E1BBB28B5FD76CD801206C80804057F2FFE9DA192593BF
5C3875605E61811F4EFCCC3FADE88AE09795694B47C874175E797EC684F772B3EF8CBAE20C5001812EC67FC7315EEE15E8DEDD9F9DEFCF19D82B3B2F
CB7FC1B857CC523D6D534A08BDBEE0BCFE990FCDCC3CF6082C011435400181AEC47F1E01B45CD7AB87BFB82850D877CEE97D87E4F59E7FEA97F6219C
764C9331FAC4C0B367644E7E2BBBDBA730022CE82F20D0850C80E33A96B5ED9C8181B3E6E6F8875E7CCA90FE7DC69C3FE455BC3DC5E37FCBA2E4C705
B75D5E98B9685AD6A5295718000181AE26005C475B7C42DE6BCFF9F2C69D5FD26FC490F3AE1FFDA4592E79C78332BC75E9E38F4ECFFCCDEB59C30F3A
82FF02025DCE00586CF3A9812F9EF517954C281C3A75D8DC0706DC4A0EB4A75DDB35295D3E6FDEA22B7C051F0FF43F436D91001010E87A19807D330A
D6BC1028EAD7A7FF98F1731EFF68F6EF686513E7BF6312145C74F99657C6E4FEE17EDFB45AB1044C40A02BF1DFB300B6BDE394A2158FE68D38A37FB7
E1A3CEBDF7B39BE644424793B66531A21363D7D197A6E78C593DB0F8114554000404BA9A04709CF5137A6F7CA068F8F1DD4B464C9B74F5AAE7A736D3
9DF5AE693166A8B466E3AA6B87147CFB50E0B83D4CF4000B087429050069BE2543866C3A335010C8EFD5F39CC9D35E78A4FB12F5B35D36A58CEA3A3B
F0FD96B943732E39D8BBE879D9141D4002025D4C0138AF0F3E737DFF407E5E41DF7117DD7EDDBDAF0FBF2B71A0D9668C606418C9BAD7EE1E337054F3
0D79D3CB9923F82F20D0B55200E6837DAFDED4C71FE89E5778E2B1B7DC74EDA2B1135A28362921D8C01A335F5D7CD144FF9B6BF27BBD63882DE00202
5D8CFFE4CEBE57BD5E50D07D48F7DE63675E31E39ED2CBFBECC37109634C31239ABD6CC1C9837B0EAB995B7876BD25F27F02025D8BFFCAC57957BFE4
F31575EB575238E5EE47EE587275CE62BAEA3D0D234C88A9EACBCF3CA56844D6BB6B02BD3EA422032820D085E8CFF9DF7676E09D6FBB07F27A1E3BA2
EFB9B7BDFBE4EABBB25EA0BB7762C2E94FA869ECBFFC946967141DDB7E41E1054DA2075040A0EBF09FFB7F6BEBE8FC156FE4F8B27B9C3DB1DFA9C7CE
BAE8FD6B32EF32544639FD11C62CB169C72923270CCD7E6B4341DF0FA8E8011010E83AE2DF75D937C7F4FBEED2AC82DC1EA7F61D3EEDCC1927FFE1B6
CCEB5112518279048089FE55EAFDE1BFB9B8CF98D60B8B66378A31400181AEE2FD01F8A59E05AF8DCA0A6407C6F4CA1BFFE2ABB7BD744DE6C41623A9
318C1046067EF043F577673F38D9FFDCFB85BDDEC7B6E0BF804097E17F9AFCA178C4A28159F959FDAE382E30F4F65BAE7DF62A7FC926BBE6892D04F8
AFE327A78783D75C3975D890874EEB35B3C112530002025DC7FF2BD7FBE77E5D929593D5FB8EA1C5A3070F1F33FF9240E65396B1E8CE103110D2D5A3
B77C4B178C38E594FCEB370EEDFD0A11638002025D86FEE9D485FDDE5BD93FC7977DDC6B7307CC3DF1C617B7DD929779A365B1AF2BA8810DAC246A17
4B3F0C997EE9B009C17B7B4CAF147B000404BA0AFD6DBB69CAC0D7DEF16705B246BF7755DE71177DBCF9B1737B64CCC188EED8CB030083A07DAD9FEC
899D3CFEE263025FEC2FE9FFACE188084040A06BF0DFB14A8F2DB8EE86DCBC3EF93DE6CE1E73CE3DBF9B3865D30B59E3E5B4D3518F11428435B55556
5A8F1D77F9A41ED3571CDF6BCA3E4B080001812EC27FFAEDC06ED7CFC8CEF1E5160D9D5C32E8E253467F4E1ECECE39FFF3D5316468181B54D9F65D79
BC74E8CCD3FA051E79A767EFDFB4DAA20B5040A06BE87FFA56AF499FCDCD0DF42BC81B33A9FFC853A73E1379CA97EBCFCECCBE49473AC68890F56B83
1175EEDC0B87759B165B50D8FB0D711698804097E0BF6D698F145DB17C4276C9F88BAEEA5FDC6FE231A3F61CF667FBFC3E5F76CEDB1402005AF56AC2
DED9F8E1A4E90386FBEFD5E6E4CDAEB1C559C002025D80FF8E15B9B9E4E50D13F2BA17DC5C33C33FFA94E1D7C8CF65F97CBE9C407ED6A866860936FF
E56ACBF9A23AFEDB8B078F2F19F8F1F6BEBDBF3705FF0504BA82FF778E9ED9FD9D45DD020559C38ECCCF1C735CE0196D3AE77F4E8E2F3F27E3739310
6C6DF8346EFEF86D68DDBDC7E61D139814BBCAFF2A151580FFE07D156F81C03F0FFF378C3B66E9DD85C53DF3030BCECF9A36A5D78FDAB1D9FEDCDC1C
7F7E20E3419331C69DBD69C71796AB579FD463408FDBA2D3F25EA2620F8880B0B75D81FFF49D92A15F4DF317F8BA751B79D6C0F1C37BAD6D19989D9B
93939DEBF3659C27136A51C76E54DCCF36B49D5F32302F6FC33E5FE043380A40BC81E21329F04FCE7FC7D51FEB7EFAC6C9DD8E09F41F316AD4AC297D
FAEDAD1B9093CB2D406EB62FE3CC984D6DDBAA5A166C3958F5DA23134AFC13A397F9FA6CB36CF1191710F8A7E7BF6D375DEB9FB17B76EF41834E187F
DB09C79D50E85FD6323433372B37273B2B27737AD8C6C476353DFDF5BBF4D667E765665D9D1A903BB15E9C042220D015FCBFBD6D5AD1432B87F6EA39
78DE9C1B6F9C75FAC8DC27D155997981404E6E8E3FFB7AE49AB6092FDEB248DEFDFBF939D99FA8BD033725C4049080C03FBF01706C77D9D082EF9EF6
F7EE5738E1B84125BD069E9C7D5C746D8EAFD09F9DCBF1B86B93348BA2C4DEAA755FA8F74FF567AD59E72B7E55133B000404BA02FFC983F927575D99
37684861CF71B79C37F18C6306E67F6A5C985314F0F97D39BD36A7D3B64B43D1234B0D75C94F7B4F2F1CB27F9ABFDF5A66C1E620F1060A08FC53F3DF
B6B57BBA3DD634BC7BBF6EDD8B473F71CF79AFCF28C83E29B4293FD79F9D939D3D29EE50CBE379A2DCFC6866DDC999F3C3FDFC638F9816F87F610004
04FEB9FDBF13BDB5DB79E593F287F4EA3174E8B0EEC5132F3D7640E653C6ED99B9BEDCCC9CF71D06A33E9CE86A83BDF4D8E6DBB377287DFD178799ED
8A0C8080C03F3FFFD1333D021BBFEF51D4B35BC909934E9A3CFEECA953FCC5EB9A4ECAF0E7668E4D38BAE9BAD436235A4CFD646EE5C563B5B2DCDC05
BA25F82F20D0150C805B39D77F5DEAB6DC417DBB9DFFC28DCFBE7FC929971C9F757C55CD89D9BECCCB988D4DDBB5DD50B943DB77ADFA72E053EC565F
DFAFA8ED0AFE0B08740501607F57E2FFFAD098A292E2536F18F89BEFA70E197FC5A08C2B8DF579995927851DC3B299C9369531B66B77F0DACCA76B8B
B34E2CB705FF0504BA840170DDD8AD81B17BFED87BC49831373E7CE38D3D8E19F79BD93D7D8BF42B32B3FA3581FF67A6B257A1C6AA1632BB57FDD2BC
82CB52AE23F82F20D04514C08AE1FE4B1E28EC3EB2E794C58FFFF6A459D37A5F3A2EE3B8F8431999DD1B2D8332CCEA562469F8A36663DA55CE97FEA2
07B123F82F20D0552C40F28182924F2ECF1B3DBCEF9839F7CF9B3867E0DCF386176F5892913132CAF94F10AD58584A9A576365EA3AE7F240BF1F7E71
FF82FF0202FFFCFCB7778EF73F5C7752E098E387CDF8E3C5C553478F7CEB8AEC7F912666DECA18C686A187D6EFC5D5AB71684A5B4777FF9823AE9D16
FC1710E82A1180F162F781FBD60CEBF5DB472EBFE3F7C7CD9C37FAB25B8FBF012FC83CA32D25534CE248D7706305AD98EF2CCF0F5C1D779C9F0F0E14
6F9E80C03F3BFD5DC7A99B9E3BA3FA7AFFA56F9C72EC45134E9EDEE79847CE18D7F27246F1BEA3F54C56E3C850F18106F6C393E94703FE05BA23F4BF
80405732005F0F2BBA6C41A0F7BD4FCD39ED8C134FB978D23D574C2C5B9EDDBBA67C0FD3316A5C5F866B25F6D6CAF4FCFCA20FA9DBB9FE5FD05F40A0
4BF0DF4D5C1FE8B9E2A6DC8BF6FE307DD4C4E36E1B38E0FA714B4A737D072ABFC6BA217FB1687F7B858E7F5F4126E4F55A65FE7CFC8FF0FF02025DC2
0438E9D5A3F26FDFD0BDF0C95593CF7960C2D409E7AC3DFBC1A6BE392BDA171B443DB04D93CBAA59F3ADB4B65FC1F1D5FF7AFC9FE0BF8040D7900089
DBFD795FDE9039E4D613279C75D2718FBEB0F4A44BA3A3335F525FAD21B1F7BF8FB5EEAE23DFDCE3ACC9F3CF8C78EE5F705F40A00B5980CD6333CFDC
3832FBFA3BFA9EF6E08C8BC7CF3EFD62F5DCCC07C8275B70F2C0E1CAAAEAA4F3C207E94DF9FEDB35D711EC1710E85A0220757766D6DA6525F33E3EFF
9CDF9D34FEBCCBAE9A16999F3907AF793F26539338B52DE61FD6A55F0EF85FC162F45F40A0AB5980F48E7159A7359F9E75D9ECC209575E7BF3D9A71F
DBF844E6D850F9A2B64D3B656C0715E95F82CDC3FDC3B69982FC02025D4E01288F65653D7A5FCEF8392517CC1F39F6FA8B066F7B3DD3BFBEE3B19AD5
77AC5FD71862FB7E6FAECBF34F6F11757F01812E28008E9E94DB73D51F4E78ECDCB1736EBDFDDE6B46AFFB293367B1F162B976A86E67524AEF7D35BD
BF7BF14DB270FF02025D5000D07702BEA5F1E9E75E72FAB3AF5E3069CAB00F0EFA331FC68B561CDAD41A95C3FA471FA53717173F4E05FF0504BAA200
689D9E3F337EDB801BAFBEF884C231B346BE941C977127FDE1E5776F5F7EF4C743ED2FEE4BDF5FD47FA938F85740A06B0A80970BF2377E9B3FE98229
679C33E58A6EE7A9F3322E24EBBE8A57A6E48678FD6D87F4A145A34A85F71710E892FC77CB4F2BB8533E77C87DF75CFDE0E433FB1C9FB835732ADAFE
83194FC57645630F3737F72A981311FC1710E892F477F4470BFBD67F967FDDAD43A74E3F697AC9CADF670C95AADE8F1EDC890EDBD16768B86F60BE2E
C27F0181AE69029CED538A5E4DCD9DFCC425537E3BEB9CBC050B33FAB5447E32C2CD6E73F8AD07D2A1DEF90B4D417F0181AE2900ECC4FD3D87372C2E
7EF0CD8BA6CD9A9075E5964CFFF6E4E2D6F66D89F5B517DFEBEECCF37F690BFE0B08FC7A39FCFFC47F67C594E2EB37F71DF59B89C74E9BDCEFF29641
D9CBF147B5C16FF77DA55FFF82F3AD7FE06EC17F01815F29FB7FA6B1FB0FD2DF75ECE63BFA16AEBB3FF7F40B470F3DABFFB195B3325E662FBCDDDE7A
E0A3A367BD653FEF1B5D239AFF04047EA5ECB70D0D3BEE3FB69DDFFB169B7D373EEFDDE098B3EF98316DC6B8E32AEFC9BCCF7EE1A1ED91FD1BAB4E5F
61FD2E30279416F57F01815FA1F24FDBC155F7CD9DF7D2D60E85DACEDF6D01BC6FB0ECA66B7A5EAEDEE33F73D2C891837AAC793E73B6F1F2A237162F
5CB662F85A32C97F6D42A4FF05047E8D90979DD7BFB0A0202F7FF8B9CFAFAA962CF71F30018E4DBF1CDFFDD323FD0775CF2E1AD2E3A38FB327A285CF
BFF6E2439B569D582EF7F63F6A08FD2F20F0EBD3FE66D52DFDF30323BF297B7B72EF80CF37E89CB7CB74FBE784C0DF21016CB77E4EEEA8C882A251DD
F30767DEF95DEEE0D62F1EF9F8FB7B3F7D67546B852FF70DE208FE0B08FCDAE8AF7F3921503CF086FD96696B1D47B73C3DA76761F1F96FEE6FD3AC5F
92827F1BFF1D477FBF4FDE1775C75E3C77F0E8BEF7EFCEEB5DB5EAB5AAD2C57BBEBE8A96FA4B7E7245F82F20F06BA37FE4E142DFB0372B2D8B2948E7
C1BF436BD73D7E6E49AFFED31ED8D09CC29E1170FF16FABB163B3831779A74E509E78E1C35F686B26E8175EBEFDBF8E68A8E172F3357FA86EE17FC17
10F875B1DFB52BCECBF19D516EA3CA0F3FFAFEF0D1968650876C39B6193DF2D64525DDFB4FBD6A7159148201F73FAE0E7AEEDFB6A32F0EEAB1EFBB92
B32F3B67C259D527FBBFA87EF8C08AE515BFFF83FD5ECED92191FE1310F8F5501F386B2C1B1718F32CB66A1FBBF3BDBAB84DE2F571CDC4F14D1B246E
039AD73E726E718E7FC8DC3736B7CAC4FED904FC5B0BFC3A2F66DBDABE99793744468DBDFACC5E451BE7E73ED071DFD2F7DFFF66C87D6C5ED66F44FA
5F40E057447FEEB0E30BBB15DE506BA3E56FAF8C51C732438D86432C75DBC24547E438311DC78AEC78EAAA51858161173DFCF1DEDAA8CE2CDBF9F908
8F7F83FF96157965C0E08E4707CF3BE998C22F96F9EE7417DD76ED97FBA77D6A9DEABB572CFF1110F835D1DF6CBEB3A0FB93C4A93ED468394E2A1AAA
AE942C8726A514337559C60A6336070B6DFCE385430A737BF43DE5B6AF7636C414C3FCCB2E21B7D39C5896B1E3B4C2D70FF61A366254C1CD65B9B3DD
1DF36F5C5676C106765ADEB358A4FF05047E2DFC776C5A73ADBFF80D2BF2632567A69E6A0DD5B6ABD291D23DF5AA454CCBA628A1C4536DB565219BBF
3E5AF6E39B37CFE8559C37E69A9797967620F3E736A15F8EF3743CFE3333F8789F5312F37A1D3B71E889958327B11D93CE78FC8B5955F6397D3FB344
FA4F40E0D7427F07975F959BF7A1AD6ED942EC48D9FE2449234EDF437BEA11734C0747DA651AAF2A6F0EB725A209DDE6E4756DD4BAF9FEA97D7A978C
9D7AEB736B6AA2980703AE772D07927FB66D99A6BCF2B8E2753F168D185692BF785C6EF996FE43FEE5AD714722C306FF64DBE27D1710F895F01FEFBE
2CB7E02B1323498B351DA825E9B41D8FC531C696A5698D7B361C9630334DD3B2180EC7DBB61F66D4B680E2B8FAC8E60FAE993670DCA4532F7FE4ADD5
87DB9306A6D4B440FD9B8CA895F34B6E51670E3DFE84BE0B6FCB2E6D3C7EF2D20F8F6D3E5274FC5E57F05F40E057C27FB6EBECC2FCA5968989964A74
2EE66029A9391E558CA469A6E29AE3986A2A9408B5B72768DA6CA8B118570249CE74133C3D6A6F6D58FDD065E387F71B3E61EEB58F7DB47E4765ABA4
626428ED8BFB0FAAFAA0F8E433F26E7C23E3BBDABE2396DC34A066B56F72B523F4BF80C0AF83FF4EFD397981572CC6629595AAE53AB68B6BA2E934B2
A98E6D2F51E762433292523CA132663AAE6D29A188221BD14844D319B53A5FA4B41FDAB8E8F7D3464C185332EAD4EB9FF8EAA71DFB2AD6CDEAF6ECD1
E3279D77CADD9F673CD13AB4E8E4419353EBF2663609FE0B08FC4ADC7FE24EDFE0555A6D45EDE1A86A1A468AA22DDF6969EC98B24A4D274D49BCB131
6510D3326DC6A53DD7F716C64EDAB610D1316308CB06368D4E45EF9260F5C6B79FBCEAEC09A3474D9E7AEEEC09C5132BCFCB9B30EDF2DF669C1F9F36
78DEB8F3E46FFD97C704FF05047E1DF4472F160C396C773CF34DD0608621C56354EB90D24CC7EDF59AED6A49554F84646212629BCC22520A9996CB10
B39969D9A6C9183614AC63641AF5F5C8E27F677198B4FDF0EAAFDF5FF8E4B9DD7AED58547CDC80A1A7649CAB5E30FACBFB6E22DFFA9FB4A06FC071C5
01A00202FFC3F4A71FF4E8BDDBC6C146440DA9A5234991E1A41D29AC87E31A43A9F69095B61DDB54B091408C245A13CC429029D099CBFF4FA96B4426
535A12561AEDD91ED6514A8A2A8A2EEBBA69416B50C308FFBD5527F79C38EDF48C1191DF64FCF18FBFD55ECABEA1DDB0BD1AA1D73E2426810504FE87
E89F264B7B172C72CC688228F19A268520DB62E93449A13477F17A6375CA4D3B946144CDD0AEA44511B11DCBE4DC2654D15443553AA2989A4C65AECD
B8D9304C5555D48EFAB6789231469999BE2730E6C83303CE9D3F3BA377F37319773F710BBE2B376FF8258F2C5AB1BF392A23F64BFB9050020202FFFD
E2FFBDDEDDDE75F4A62347E3B5952D8452CA853EC106FFA28D901A55789CCF198F09A2049B8C72EE13FE3FCB765832164C10AEFF19330DE63243B36D
9B6A09991139812D1ADFB92EEA18E9DD45817BB68EE973C6E91919DB37F43EE7AC9B8D3BB2A75D7BEAF17DFB0E9872CEFCC73EF9F1489498F6CFDB40
5DA1060404FEFBBCBFFE66D1805D0EABDED6A148118428F7D8BA8A49DA4ED394C1A378D3D2348570F29B06FF12E35F278E65700DE0DA8E8965CD3429
D554B525E53AD44C63478FAB48D10D4AB4B69ABA08B2A8659FEB3BB9F4F2E2B3E764647C5277CDE97D1FD14FCB7A9ADFA466E5BDE7CDBAECDC9240AF
E1573CFCE9CAC3610D33CBFDBF864928020181FF62EF6FBC5194FFB6EB26AB52A6DC94C4F1584DF9DA5D32A7218E6E3E881DC6A9AE1988C42549A5CC
E0ECC70C1BC9A46A422580EA294CA8C5E42496256618AA65224DD6898E6583C592B04810232DFD8D3FEFEDA77326CDCEC8B865FF6DC3331E4D0DC8B8
91C990FFB7D369A7EACB0FFE70EB95578FE83FEACC8BEE78E88BFDB51D3C26B0FE7539B810040202FF45EC77D13BC58159496DE9ABADE4704542C191
EA75ABB706A9A925ABEB9ACC741AFA7BA0D79732C53019225C01D846439DC1BDBE6D7195AF714D605A1459B6AE688C302F7AE0164195B92F772056A0
D8D24FCEBD73FD88299704322EDB77DBC88C0F927D32EF6332366C42186C01711DDB343BF62E7AE9C17BE68EEEDEB7F7A853AFFAC3A295075B3C3BE0
FCA921480B532020F0FF87FBAE6BCBAFE5E58FDB6B7D7DC3FAD6CD2D88D1DA3A15A17844D1AACB648F69B66592A67D1DA6632AAAA6A524C2637D2316
3790C9794F9594A45816B308620C21D361049BC4C498E9714967F03B4275ACD92FE78E5839B3E4C2FC8CC93FCECBC9F8AEAD28E38F9CFF0E6D69414A
558B6428545615EEF26D87D62CFFF1C327E68C1FD4A7B8679FE3CFBA63E1F26D871AE31A35FF343D20148180C0FF23FD1DC731838F07BABFD4A4870F
B7E0F6286BAC6D288FE058690D52DA5B64FE1A2B851C1BC5AAF7ED6FD3313375A46A0634FB5886A498B66DDA2655310F032CAE0AB04A2C0BFA832877
FA90096448D119D329570D7665FF82B7EEEC7766BF8C93BFBE252F63536D41C60B4CA10E8E254DA931829824C75BDA5372AC26C5BCB12047AF3CBC6B
DBB28F1FB8E8840923C79F78DAF90FBCFEDDA6326E0788B770E0CF260D85251010F8FB9DBF63EBB5B7F5E8BEC2B519D559654BD38EF535F58ADCBE69
2F52EA2433EDA4D9FEAD3164E81694FB08D7F1C4E2BF508411340171D69B5CE9CB3237009C948CC7FC143113C201A80F30CBC4A64528E58201317456
EEEC270BC68ECEE8B9F0A68C8C2D470219AF580A9709B6CBC309DB3478DC802D47DAF3EDE1A44EA25B57D531CB1B22B458B47AFD474F3EF6D4EB0F4C
1F356EC429675DF7FBD757EC3958DD168BCBAA4E6CA1060404FE31EF6F2A3B2E0FE4BD9DB6F62C5EBDED68F9E155A5512A19EADEC6446D9B6A22EEDB
A53872699A599411DBB208E778B82EAE4B86EB300B4A7EA68BA56844336DC7B4305401B9EB675C02608BE904EBDCF59BAE6D19D861E985BE619F0C1B
3536BBC79B776564EE3D18C87CD7D1B08131D578D081314AD5B4CBF1702CAC1A9836ADDE9B326968DFB6C31D110A3B876053885953BAFB8BD77E7FE7
75B72EB8636CB7DEE3A65E75C78B2BCA82FC1E9D190261060404FE76F7EF58A12F4EF1F9466CFAE99ACB3FDC7368DBD2F57B0D4BB358E5E1D672953B
7423DA8E087F295514A2712F4E4DFE1F53C21A85CC1EE7A481B81890650D238A482AA672B9CF5D3D7F398F054C2392C20AE32103FFB3C1A543BA63A8
EFF1CB8B0765FA5EBE352363FF1E5FD61247B7DCB463634333341ADFBE7CD7FA5D614EE0585983EEB86670FD0F2BCB12366ADAB8624D43C2B2ED703B
363123306D686C7FFDE5AFBFFFFA8579C37A741B37FD8E0FD71EE950BD16A2B4300302027F9301602DCFF5CB1D78D79A6F16AD88A1D0674F7E1B524C
89AADFAE41866132B93594C22EA189705979824BF964591357F3A6056B3D2C53E57E1B116A5123A9D90EA3292D1995991E8C12263777688C9A488652
3E33F5B6568D39A6CDDC87737FB7B0C7D0C2A249C764E63496E6F936BAC4E171BEEB985C3390644CB568B23ED9BE7FD35189EAB1D21D0D1075047F5C
F8C6CE90A51FFDEEF34D6D929234529ACAAD8B050CE7614BAC72C3D2858F5E77C1F4D9B3EF78ECE36D0D0908087EF909C53FB280C0BF477FBBFC963E
F913F65A5CC96BB52B576DAB6B470DCB0FD0039F4A3643E17DAB2BA96BA95547DBE4A8CA098E6A0E8621E0B7A19EE75844910C6CC2E8BFA2316C5B90
F8430617068462C5A05825847B6A8B309318CC362D8B6AF68EDCC16F8D29ECD6A3B8382BFBF0169F6FAB834CA8E6F12740919A7A9D584A5B546DD8DB
66A9C18A860E4D952AD77DF0E3618358725BEDC677D7E8B6A174B407A3CD95DBB7B512832A9AC2A8033B866D664477BEBFE0C9C72E9C30E3BA375636
2474FA274543010181BF8AFECBCFCB2FB8B4C34E1EFA7AF19B9B1BF444DBD7CFFCFE9BA4A3B41C5DB7251E6A8D51433AB0B7D9B0D2D472284651C5C4
C964C2C0DE02509BE9A006B81687B880FB63CB74BCB660EECB4D8BC7F38C1B099800E416C3B6F8154C6C1A33B31FB8BEA024EFE49BF233577F9C9B7F
D025A0FF2D8C5B0E1C09C651E2D0FA307F342645CB632994900E6EF8E4FB6A8648D9EA65E55C64B8A65EBFABB2B95D575ACAC396491DACE18EA71FD8
A1297210DB3044E4B8A87CC5C2E7EEBBE637F73FFEFE4F554105798D0BE97FFC147301812E1AFC1F99D1EDA215A8E69507172D6B8A3225545BB17165
29624664F5F20D35B685DB5BCA0E1C504C578FCBD8B06C8B611DEBA91436B996E7EE1CB2019407F71473FFCCFF6032061E9F51E655081830DF64FCDB
88055100E800EBA5CC531FCDF7659CB66450E69AF77C8515AEC12C148D6A48D5A989D5444D8DE190704B5948A57AAAA943756DA2A0F6FAB2566473DD
A1C9A18AEA2892A49AA38C6B7B07AB8AC9CABED9A9986A4D6D9B81A98E49E779422C99A85BF6F0FCDFDE7AF3538B4B233C0CF9797D89F8871710F08E
F0649B8F2F382F54FBF23D4B8E6AB2DC50555ADFAA735EABF18E668D4BEAF0FEC6245611F7B14628C65C98FCE5313AC4FEB0D5D7856A3FE6F487BA00
B5B8D4B7A00180310C1D41460A3195C2EA5FA81844935024046B40EDBA41FD66E7E516F47D6F6EEE81AF03C774A431D705887281C0AF47B8FE37314A
346B96ADB496272DC7763A5A92C872E0B6D84C84757E4764A889A6A0994E9B75871348A336D4062C966A969AEAB67E7D089ECA42C8FB2151B875EF92
5B674D3EFFB74F2ED952934490BA48FF7276A1B00502FF8BBDBFF2EA808251EF2EFF707BCC4A361FDDBEA52A68315BD7514C42614EC6E6DA18A73CB6
F4B686A44A2D0B276AEA35E8E535A1591F086D70674BA11F80C692BAC6230010FF9CFC889946543122067F95E3724D9090F9DF33CE529B18785EC6C8
E33332331F9D99B3ED7DFFC8848B618DB80B6A823143E3EEDCE42ADF3612D50D35F17852AE2B57D2AE69C52321299C9499091D06DCB8C0DE91545D45
1CDB06435C8368B08284A5E46059DCE68FA034B433C7FC7907B14DDBEBF7AF7EF4EC09E7DFF4C6F7A5AD32B2FFB46B48580181FF85EC771B6FEED9E7
BEDA84CC2C2951B57B571B71EDB485DA9A1BDA53F170A8A641E1AF030D1F3A7A349C3451D3FEEA568930840CA235B4733F6AAA7189FB7DC6BDADBAA3
1EF1601F2680D5A80A09419350C3DB0EC449C9250501CB60541E5599C1DECF1C74FFE88CC09D67F877BEE73F564E631358CA8C70477363BDC15FCDC3
75BDB921A25848694FE214B26DD788486A0ABEE49A52631433EC591F6E5F1C93DF2C2525F88372E902DEDD4D5B9CF918727F4C6D690D218B759E4B62
EBAD956B96DC3B69F894336E7BB73AA6116F2B89E81C12F8DF487FBC7A52AF2B4AB98AA74ABC72FF7E05712F6CA70EEEDE5953839085234DD13A2EBB
39B974D5A429091149E66CB118D2A94DABDB794C6F5103638475DD604C6136C59CF408ABCDCD1A5802480CC0C640FE421B528798312D9CC254454DA3
B36EBE3337E3FC737CAB5EF78D4BB988D92429311E5B70EDE185170E536496F67602B940676A1AD8F6B684716BD256A7328BF057F3E8DFC52ABFAEA9
C513D4F4A4870ABD0696E9520A1685A1F8A145D5602228370B3FAF19940F6CF9E2C133E75D75FB1F5EFAAAB44DF3CE2B490B2320F0BF88FC6937BCA0
EFA01751DA4A2E7B7B5D4D9D9C4CD447D4C4BE9DF58646A5788AEA9C86ED552DDCB94783091C6E07850DEE9539B0E7CF73E990E137894591C130E1EC
05E1CF7466330312813F0386856DC7F2B27FDC2373234174FA70CE8459BDB3AF7C29FFA73702132497582CD626F12B72BEF39798E51147A7FC311917
02DCA4707DEFF05F39B801B0789800F9067E77D83BE0A2946661483BD891766E1D4C95D826E68A1F45B08DE1741247D6F8777313D1166A4F459312C6
9D471CC51B0F7FF3C4A5B34EBBE2EE3736D727B0D9D93BF8275345C214087455E5CF9DFF69FEBEDF3866EAC07BEF1EA849625CBBF950BC627B83633A
24DE2C4935E592657B2CB0A989F588645A709887473B9BF0EFC726E331BA63411D80621DE90CF2009649B8EEE6D29FF13F024B594CE2BF83C01C0683
18A1086B2ADD5C14E83F326BCAD2821F5EF48D4D3890FB031B0119441DD97A85EE1AC8E552DF51C33143D138DF6D02F546168F13B810741CF35B703A
1BC8F6AA8E5CC9B7B5B93087C4AD0CDC901EDD519F423A412684FFCCB183872A1A2B4A0F56CBF16444259D8B492D22576F7AEBFE99C74DBCF2DE459F
6F6F88C38A63F7170B29E48040D7633F783FDAF4705EEEE02DAE5DFDDDAEB8241DA93AB8F968329D26D84D63DC16344C2D1A569162528DC0483EA79B
05DC672D077808C0237A1639D84628D10DCC60AC87E89870156E31F0C314C67D6CCFDF834280823D74E6F32B60AE17B8318056FFA959A3CEEC71E2D6
D1AB5EF58F4B3804727AB6914AA654232E3B58495BAA4E34439364953A88406F50AA9A5F586E932DDBD31990E1636AB455E1DA010603B87680B21FE4
112CD84C6C595A284EF8E3687A22A9F0405F6ED7A84BA80B9D8646548B9757C70D994808E202F3D0C71F3DFFE8AC7E3DC64CBDF39D5D8D12B26C57E4
0504BAA6F377CCE4A6D9B9791754A4CD8A759575D51BD6D525EB9B719A33D54AE8B1ED3B13C8D02931E21A8ED734CA30C567737A9B983225C2A09EC7
48FBA12659235A5242DC0C208475836284A17E87B84386E91FF367967ABD40DE6961941B0BAE1510C12A7B256BDC694585DF9DFDD9A3B943C3368254
81A54BB02B48E75749AA504ED48924DB501530E178B1D60AA833DA2047BC8B11E8267229A6845045E74FE81517DD441BD7090A141BA0EB08EE4C62A1
7629D9BCBDCC76A98D89E34D3C58C6C6AF0EC542D5DB76D51B185756236EE5A46DEFBCF2D02D270D1E38E982873F3F12E95C38F22BB201AEFBE74649
D82581BFD7F37BA77BD6BD32327BC84ADB6EDCBA615BC3AED2AA30B538035DC34E2B074B771F514CA2265405F1D89EE85A2C1E89713F195589C90C2E
0428E4E3A00F8062876B70AEAC2DB00F945283783B0179A46EE9D8663C0EF056FF9B5EDDDFB4BC3E207E0D9D31640687149614069EFEE3D3B7E70E69
B71086AF806C20CC49DB69937845440B1B94C7F58E8D6442340A69078FD4F18620A2049A096D24A99CE12D4D2AF41E404622D28C180927A1E6C08509
042B8EC5144357EAC2902D94E2C44689BA186B5C5E9AD4B470434BCA54DBF66CDA51153389F7EEA4363D75C74D774EEBD973F2BC677F38D092429DFB
46FE8C7AFF3336C1FD0B03208489C0DFC7FDCE51FF2DE7E56677DF97B6DB763798D0A09748C912C149CCBFC6142EF86DA462DD1BE13728A8690DCB29
594D2084249DB397334A3140C73BB0CFD3B2A0BF077A7E18C69A464CC8C63384B93DC00C7C3FD7FD60046C88DA4139A888C70AD89A93991BC83DE9BD
5B66660F6A36550C3184062AC0F69A73A8D7416052C30001616109430B21A3281EE62221D9162526BF0861B088985F1C230AF38008A28DCE04A3E925
09741D36937063655BDEDE62278D527262C7B26D41AA27359318B072C036B9226A6AD264852B944873C293FEB1152FBDF4D84D03F27B9F76F1B3CBF6
B62434485D3ABFE07F283670FF0AE2832DF077F19F7BFFA5C764E78CDCC9153D2CDAE614D0F657E9965A59265197418DCCE1325D0B49E08E29520D42
6D07234655A223C3CBE43944C68EE7976DEFDC5FFE17FC5A3CB627A914853C01F403727F4FBC7900C6F539011DEE8047A6DEA40023CE07D981427FFF
0FEF9F9A734CB3A93162A88440BEC0FB507B26034C8C091D859D0A837ABBC5244925089AFDB867475C49783DC88CDF00C54D9A22702619839B50C810
5A4A44830A82770631940BD20E22C07A074214A2C005B8A8402E4E8565CDD62395DF3CFEF477F5C14428AA7B2BD1B63E74EBED174EEDD5BDFFD049F3
9F59712498D2A10A62C12C83FD6766401800815F7BC6BF735116FDBC67EEA017E3368F6EA983F96719B53720C75114C2596E30952AAA8C50A83C6470
6FCA38E3C1F763CE15E62DFA027E71A9EEB816F00B3A766D8711EE5539F1B9C5C05CB973DE73FABBDE10100803E81100BA41FE1F132FED4EB1D598EF
2BF4173CF6ECE4AC01CDCCE0860525603A807F533ACD74382B841148F151A8369A100E3084F8351D8822C030D85082847962C84D2A24568D4C1E2240
B2C2131F5003A03899A49CA9B662B8908B70ADD0D1185403BCAD445C20E824112CDFB6A5B625816D6C854AF787A8253794EEDF5BD9AA11D25C6580B3
A747577DB9E8F2DEFE826EE3CEFDFDA215071A3B925CC350A8173AEE7FB30D100640E01FFBD8A47FA6FFC745BE498D1045235B8E37D7D54B51CE6AAA
13CE23837F9829966409716F8CA17B8F00578DB844BC4A9E8A404D7BD975E8C507EA43470E671A33A0F716767E4082C0205852819F30F70F7D7F4C0B
C609970F20E80901923362CFF51715163EF9F505394338FF39B59B2A61AAD856AADB341DBB50E16326B1217CE00607231369DE3411786EE02ECC0AC0
7378B3858C52B98D72B3E0C0B2315000902404098094B8846953AD6ED969DBC54D871A1B24CDAB464AB10EFE83CAAA242763ADADEDB1FA1D2B5B2955
F564428620C4AE79FDFE6FA2163362092F2FA097AF7EE2A653C64C9A75FCE8B1675C7CF7C23507EB83705E91B783F0AFC8F9DF6700C4875BE06FFCD4
C05CEC3BDD7247B75B4A5DB54A9496CA86481861D3764DE8B4E5C1B42A31DBCB9C4182DDD4C206D4D6BCCA3D0CF4A906B08BF35C479000E497E31CA7
30EFE7D8B83302A05C00505D276A920B65AD23C6FF12427B1C8E3296D2B900300C15A9C4629AF39EBFA820FF8E8FE6660EA8A32A77D4C124A87CDB94
623A140B8986BDE081C52A0CDB41BA8534388BC86218D28490DDB7B14A3A7B90A03000194C4288377EC490A6C32A1113CE278847537A2A4E5C6E39B8
8821A98E445B3CAA2114EFA86CA8AC6E65906E28FF7EDB91ADCBBED8511B6AA8AB6A4D99B87DC3D20D0DE148F581D59FBFF4F2371DECE79E00261BF2
77B7DCF4DC2B0B268E183C6AEC5997DEFFDCC73FEE3EDA149274D0038EFBA796E0BF263F2F1480C03FCA7F577DBB77D18456DABA6D6DA5AE29B20ADD
32A8A33D54AF40550FB5EF6981C0D6EBE4B36054DF9BDEB5BCA59F9657C8F39AEF28C636785E4E52434EE95C5F13CFFB13352A43932074DE9AD07187
0D0ABE1AC12F8E4D3BD381DEB52D9AAE1F5C58D0EFF14FAFCA1958CBBC2122DB8B273887B83DB22C35963214299542A918B3EC60CA816C3EA4167844
C04D09641EB86AB160DFE0CF59027E0158090C3D8986AAC3CB602E89DB29062709B95C44E890E1E0960BD76F3F90623689EE5BBBAE416228DA1E52F8
53269A6A3CFB17DDF7F187EBEB1A6A766ED8B8BF26148E87154D295DB55732A0D1E067CE196DFB36BC71FB1993C69C306ADCE05143C64FBCECA177BE
DC70A4A12DA160EA0D15A4FFA48FD0FDAB6CCC3FFAEFF86F5401C4C75BE06FE03FD7EA8937FBE5CD8A5B91F61886195D6AF2B819D76DAC6CA8C5407A
A626A0850766626C2FBA87223B28FF44040A7DDCFB8205D035422112270441668F3B7E98A221089164D9BA6698FA61901EB32CCE6318C2017D80217F
07AD838ED712ECED06D72FF005FADCF4F1F5597DCB2964F9A17D1F6A09E0CF61CB8086B454288141AB13DC9EA23A857502D06660D2A0E2CD145818FA
0C19C4FACCFB2E954B05E6C0DE01FE8C5C0618200C6C66630D5A03288F5A2C1EB8730B106D550C146B899A04314BEB684552B8550A35347063468C70
5D6D4449B5B485E2314DD25372A82D1C6FA8DCB7774F433BA40E610ADAF1BAA87884938AD75596EE5FF5DAAC817D07F7E93760F898D3CFBEFACE573F
5B5FDA94D4A12BD9F9B36584FF97AFEE3FF6AFF86FB87FC17F81BF89FF6E6A619FA2D3D5B44D0C83E838DAD2AC712E6991B86913F0F610CB9A04027C
2948C0BF7A5B3C609FEFBE32CF71EB9CA6489775D3329081B097CBE7AA414FCA0A82B09EC7CA2AB0DB2BBBBB6002A0BC08A301DC5A6015C32A006802
F086F5CDF46B8182408FD9A7E70E2807AEC2AE608379860323A8217A15431E1030DBEBF62586B74B1C7A8821C3E8251F79F411AA0E429A10F600725A
421A825155220EFF2E8B290AC429FC6E7A44E1DC47C826E136EC8281B22C2595E26280DBA1969A9489A2B28C1CD7921A1A1BDA6412696A8E11CE7466
6BD1E6C6233BB6566B8E860C1E9434541CEDE80855D7B77644BDC606D842EA4D0FA9A170F9B64FDFB8FAA492FE4525BEBC82FE63CF9877FFBB3FECAC
6849AA705492E94D2FFCEB69055EE6C471FF2E387F022FF3DA3911210C80C07FCE7FC7951EE95638319C8EAC7AF6E9A73FDB74687F7390CB73154EEF
D131ACF2B33A9BF699A9EFDC0B83B4DE2E2F935B0A1CD7407973CF4E8C4454A1B0C99F47D2E047C19133434F48D8826AA2EDADFD418AC2B581B7FB1F
D80ACB8029450675BC5C9DE965014D94DE5EE80FE4E6F5C8E97784C203308A34EC95D618465E8FFF2F0B43E026A04BA82749740A9D45D8EB058263C7
20ABC86F160F590E56355932A816D33D99C1609388091C373131B594A2502D14853D44186906583D982FD4DA9B13321C700C1B88C3754164630E9876
500C25AA418ED34A296A32A52762B2A12A9226B557363444DAEB57FF5483B15A7F74E3D79F359A3F27016D1291A391C3DFDD79E189A37A16F7EED1B3
DF989917DCF9F09B5FAEDA5656D71C540C6F776267D5F4975682CE3AC29F90FB2F386F7B5390BF7CD11B5BF02A2F14CA1CA62D2480C07FCA7FEEFDEF
CDCB3B2712FDE2B94FB7D61EAE6809ABA6D4D2AEC7E5A86C986EDAE1A4320C4D27CCC00682BDFD9C1FC9F65434E62DEDA69D9DFC26D53418F335A0E3
DEF416FB78493BDA590B841E7F9802D2342C2574AF680733FA14F43DF0D1F6FA006D534D1A9C94DA99FE425F6E9F1CEEFF0D025184092D45B6E7DABC
44A3373FD4B93810560D432D819B1E6C9B0A82DF7B670D589DB53C96DA7F185B48679AC2CD0C0315C33CB3E1D52BA004C864A862C076001E38C04F61
24EBAA1B822DD1A4C1B876B05134126A695560AF11FF8F2B15FE4446ED9A43C8E31EB3DD64EDBE83AD06E81AFE23310587DA2A2AA334A52BF58B9FF8
2AAAE9B1FA8D2BCB19E87B6E713CC572E48E332F3A7BE639C7F87DB9FE82BCEE23C69C3061F6154F7FB272DDE6AD87AA9A5BDADA62899424C3AC64E7
3092671B3B950D9856385E9975663C30D6351D696A4A9692B144321E6A6FAADAB763E5F7DF2C5D5F1325220810F8CFE97F4F5EE17D2CF9CEF791B465
1AD196849CACDE7508CA72AA46DBDB30F084EA84FB2706A7F8602EA571A223164A72074EBCF91CB373460F227CE2F920F6737FAF572D605E9AD0263C
A4F0A6825504DCE3C23F24431200599DC701599D8DBD606F9CE77D79B9532FCDF5FC3FEC0E01CD01DE9AB3C184441EB8756AFFE2EDE0FB94880ACB06
91B7610C4C83D3B97184C70B5E3F00B2BC33846C68414618860F084CF2D820D511756DD39B2FA6282973BB846525A9688C39448EE96AE5C6D208F2B2
858E4D52094557CDCE6922D3D282AD296CA83A174536A638DAA62554B57A6B1B3C17E435356E327452BE6DFFF66D4D4D47D66F5CF2E1D31F6FDFFAD4
5DAFEC0BC1F3EF9A79FCE5CF3C3F332F901BC8CDC90DF87AF7ECDF6F50FFD10307F4193775C6B4D9BFB9F98EC75F5FFCC1072FDD77D7B32F3DFDE083
0B1E79E4E9E716BDF3F282EB2EBCEDF1E71F7BF8FE3F2E78F491E79F78F1B97FB9F5EE6BE7CEBEF0DC0BE69C3D75DA8C49830695F4EAD5A7FFD4ABDE
DC8F452BB0C07F28FE9DE8BD79810B4DB2B3CA928287372E5EB2ABB56CFBA1B6A0A147EA628A7278791D0AD735C60C48A943268E28D8F0E25BEECFA1
9AEEE5033D25EDED0000CFEB8037F7326F9E6F87BF714DE8F985E51D9603793BF0DB44F2D67E421700FF2B9C0025E138FCA14CA7EE989C6E8F2E281A
54C6BC13446D6FBB80D5F92B4412707EA8E95915E834849B60AE5520AD0FEE153485B72F108A9360381C7824EFF1BC40038E20D5F8B3C81D70902085
EFE3EF03D565FE04861C8FC70CDB9B0534CD4875A386756F23186C1854DB9BAB823CFAB79957864CEAC1A628A73185C30D6D5DEB58BB3945B0DE5893
E23190CE9C587B074AC4621D3195392C71F8C8A1FA03870E54EC79FFC54DCD91B2CD8B9FBCF9856D5535EB1F9D7FD115F33FF9F0F927DF7B7FF9570B
AF3CE7B419271E377CE4C809E34E9D79D939137B762B19347478DF013D86CFBA64DE25679D71D26933A7CC7DF7B397A68F3971F2A973CFBFFD81F3C7
CE3873D4C90B9E9871F6D5575E7BE935375DDEBFFB696F7FFCE3AEBA75F3F2FA0EFDDD7E4B5800817FDFF9D3BA79FEDCBCDBF7FEF8D9B29F3E58F7D3
AA8A9481DAAA6A8F54C60E7CBE2E496C5B8B050F6CDD561F0EB746C232C688F210DC4104BAF638E34C2F63E67854043F4CA0838678753C48B19B04BA
7A9897B947187670C0F42FB4E880D8678E91225E34CF158411F5AA009CDE5CCBE32B330AEFB8BDB8FF612E104CAFD3C0C42A0810062B85140A3381DC
96D89071B7BC7E626FC108A306BF7F627F1052151AB7051881ADB1BD2481D719E4FD8729A7238F8EA387A320A14149535589D7554A26520C0DC3F31B
6AB83D26499AE9806D30A91C0A2674A4C0AE51AFAF580FA52CB067AE0BC94B461515320AAA6AA40CB93E644AAB3EDB56515AD7D4D1166C4FC65047C5
E15D4B97EF6E28AF5CF5E58E8EC32BDF7960DEC5973CF2DECAF75E59F0D6B2A06671FB60DA60F9B4A8371B85714A5521CFA91DDC1BD10849B51F6E34
5D6E9D8C84CE62516E9C82755145C70A0FCD229A1155D32E867398B891537F5C9EA4941C7DF7C2D13DFB1CD3FF98073B440C20F0EF787F2BB56A5466
FEF11F2FB9FFE62F379556D57CB4626BE5AE5DAB77D71E3ADA5E7E2006277C85DA54DD225493134A381593B867E49F51C8DFC1A08F165531D215C4E5
3FD56D10DDAD090C53B9544924156A999ACADD2D810A82376D03D9692F6DE778A70153450697EDD5FD616720080587C1FCDEA75999136FE8D5F708F7
E19013E3BE98E8A0C14D5D8D440D08DC2D0466C1311D6AC0354192308BAA29429295ADD04604DE1C210B0A6D602B30CC05C07DA9974C90994B35EC42
4CC195389342E0B9415CF047C07A34140C0793DCB5BB3CDC37E34D15A59BF6C5187521E166199128B291624159C12B35280A310D64F0EB1004FDCA8A
A4C64349955969AD85B194B46FC38EF204570852C3B6CDAB0F1CFAE8EEDFBD72A0B675D7BE0D3F2C6B9558C5774B6A12C18AA06E461A771EDE75B4B2
AAB45E22325189A44663472B3BC29A45B92CB21C9D3FBBC1ED9145A0350B4E5B85592C0A7397A6CEEF6E90484B82C08E95DDD7E5E7FA02FEFEDD7BF5
EC795993300002FF46D9CF75F4230F15175DF6CDDA0F1E5DB266C5D36F3EF1DB0B2F9D71CDF50B1EDDA63A0C49A950CB919DCB36D4F148584F45BD28
1E247C5A6DD3616B266760B0B452215AAC3D018D0088711E85158A60ACD6E8086970FA17D1B1C599C83DA89BB6191CE405B4F7F259D830090FC411A6
3039E34D043362B95E4ECDACEF9BD16382AF70B7055B429D34EC168773C4617310E1A284F2DFCB3A33E5244E43B2FEE70C23BFA6860D4610316CAC9B
183A83A17CE17AF305D00DC4C38C2002C5A1EB5CB0BB267F1902FF8FB0977637616F901C6CECD03A5703459AE21AF7BD554984F94530D678906F2B49
83071534ACD85E4BB1855BEB5329CDDB664A1045B29E680CC3B145694B6E8E3635D5B5D6761087C5767CF9D2F3AF7DF1E9AADDA5A5559BD67EF5E9EB
EF7D505AB1EBDD575F79734B484AEAEAD175DFACDDB675D5A186C3DF7F5F6B18D1B283BBF7EDDEBEA131D8D89842DC8826A5287FA7A1C30A87E20C9E
174458C791762EC87435966AE9E868DCB227A930EBC01D3D7272FC05B3166DBFF9A6DB471DF306129D40027FCD7FB3EDF3537DBE11E7CFEA5F386ED1
47BFFDEDB5E7FFEEB9851FECAE69DDF8F5EBDF7DBF756FE5DA771797B5A634142CAF48E98996A09E922CFE49B291978E265889C7BC413EA8A8A72224
D6C148637BDBFE46CDF24EF6E61A8142DF1D5218BF19D4E15DD0CEB0850B4684A0AF0F32738C7476157B17F56A5CD02C7C4166414156C1568B70A301
1A9F799EDBF436804276CFA25C7FD886C60500B63AFD1B6739C40104EECB98A4C37311AB73C2977AA7823A268E4725E844E65684FFD935517345B3EA
4500B0AD9024A3314543BA46BCAA1B6BDB5B1D5164583B02C79BE1785DA30C618EE5506AE398AC35D719A06808344327AA9BE22D074326D2F5688210
4D6FDCFCE5F79BAA8EB6852C9A3ABC6BE3A18A504293D558FB9A153B4BD71FA88A6B8EE5A286B2AA04E20A235E51537FE8A7553529EA69248A553921
AB48AFAFE6EEDE6088C5933236926073B88ED249DA52936A420D870F94AA16B191A1A8B178044B9A93B61615E56417FBC7CE5ABC35FEE6DB0F0F18F6
A325FCBFC09F72DFEB509356CF2ECAEE36646449FFD94F3EFFE17BDF34AFBDE4868DBB1FBE71D1AAA5DF7CB16E636D847B1CCDB26958AEDE5121251A
362EDFBCE3306C02704DCBA08C07061242DC758327C234D9D8910A5217D76CDF5E9FF06A025E3C0E3377CC82113BD8ACE3BAB6D74608DBFEA04B1F82
78BBB34B17BAF3BD3D3D303CC474FA7866C09F99B7D1E6729AC1843DF30C0597BDD0DD0BF93EFEAD707617D723AE99E652A4A52981A0AB00D205B686
70588607807C21570F869134BCC141CF8EC0294426415071645CB0634475031A1152A138F252FB2E5822224B12540DC09269864134A2CBA0210C9998
5EC6920BA09437F1E8B8FC354D3BAB12B1FA24A41392BAA932A33D9C30A5849B46154B3F5BC3E326E3E88F2BDE5FBE6B7D698C7A5190AE4889147F18
83C45346A22D1C6BAC5761EB604CE1C63699C49663385A5D08F1C0C268A9684A3A69383F8DDB066ABB5CE9941D0AA928D41CC2A6A1EB0995C00987FC
67B51A6E0BE404863E7AFF1F0237BC3667F690C0F84F0D417F81BFC8FB39ACFAC121B9831FDCD0D051B1BF49AD69AAAA8954ACAD74935F7C549BA8AA
8D9B8958E99EE55F7DF4D3BE98626A71440D2B5EDBA66B1677E5F069E7C4E14A93CA299D7A7DFB3C4260B68A194CF0DB94A95827525B9C1B87A46299
9A9769E4944983908F460CE69D06026B023A3303909F839E414810003718A23FF8FC79B9053B6C04F93A93331232F528A943C0ADC6B8FB85D181CEA5
DDED8BFE78CE09FD0B8B2EDCA252DBEBF92598B30283ED8153C6F99F74D2B045F66A8210C374363079F57F882A4C1E8120292253A452B73343C16583
DA1102F643A723E6B7A28E77F488ED6D0EF65612C1BE63A82C98562AD61A8D24115437A8ACE8C9AABA889648B234571472C7FEEFB637739352F7EEDD
7F7CEC831FEA7930E51898D0444794DA289A8AC5A358D56187497B87CEDF84485D433C1A4B12FE3AFE7C8A0A1B15981E0D22FEAE721BA28175E30229
4D43BACB958796C258A5A908FF41EB9B6D236DFDB15F41FF6E535F786CFC5DCFAE3C36B7EF84179A441BA0C09FB1DF5BF0BF614677DF0D21CEABC8BE
1DFB2BA26AA82ECC03584D852D9869E3F0E15D879A5A9A237515472352A8A335C5D9DF148E6388D0136D1D30930F653E0BD66B00525102077F993C2E
4F212AB5D51E3DD2DE72A432992AFFFCABD2B0A113EFBE6933AA302D1C43588763BAB8ABEF6C06A49D9E99C11A1EC8A831AAD1A37D3303BEBCCD36B2
C05E792DBD9CFF6ABC25C5F9AA621ED8337EC1D0B6379F9A736CEF92C23CDF318F55689CAA04C918C99CA08E098BC7B8C497641563D70C218804F8AD
BCF624A80FC22C017408A85CF223A92505BDC7367F771817FCFC213455D53442E04E2681DD8496EE250499B7BF901B4048736264ABE1CAFA362DC5B5
874A8CE0C19A9AD6E6A4C4B50EB158F39A9F4A532C9D669BEEBD69F111351CD2B9B28920100E92412395ED1427BD8624157A25D5BA8A30C4378EC562
AD2DEDDC1A04DB6197B203C302D0EA60A5F6B7C1144372F51EFE93A5D3C45B6310677A75A544B7DF73D84C93FB0AF24EF8DD35F75FBFE0B5E79E1E95
3BEC8E4A3B2DF82FF017EC7713EF0D2EE8BD8439EDCB5E7FFDF3AA140F71E3299BF24F3B9CAB91D6AA6A5A132A33A15FC7C452B08D986D3BDB98858D
58AB946A694DC622D4EBF081295FF860B27082FB56AF6D1D1D952D23160E7548AAA168AA166FDEDD14DCF2E1375FEE6ED6E4780B8205E0063288C155
2C44D510C97BBB3CBC2E0197FF01E99CAF589E9599EF2FD8E4F097A5613800D28D291ED24B4909BE0E1B00EB17CF9F72C6C9030A7B772BEC36E599C6
74DA6B4A60C45B40667A2EDA464752DCD4B8302CC0D585ED78F5446F94C00146F3678F1F3A1AD1B8FAC7DEC6106CC49B2218B886131231B1DE5656AF
432FA11A47103778BA1F536CEA3044A448869E082755981E72600A22D9D1A884B9DB2FDBB8A2E248D9B6A34D26ADD857B968D6DCF71B303282AD314D
56343DD152CE5FA4B7B44310641BB244BD0E050BA918CE2DB0343D129110A52D6BEABC1392B8E6E0D6C6A04647559CB058CB91D507F95B084B46B4FA
16359A8C361C0CC9AD75349DB8A1A078F29BAB1F98747EF392910579E7ECD0D3AEA0BFC05F24FEEC869B0B03671CB2D8C1777F4A72EF9B48D6563707
5562A64D9490597CD75EAE05EC542219AAA9A86DAC95128629B5CA8E46D5C6BA785232CD68426BAA513BD7E578FB3B11F7AAAE038BF0ACCEC96166E8
1A41482394A5E470E5EEBD6BB7AF7FFBBD36CB8C56B7C3A62C4F4BF327819E01BB33F2775C0AC2DA3452090D29E4F1CCFCFCC215B60A8BC16CC6ED05
D563B071D3EBFCA1D2CAFBE6CC3DF3D207DEBCA66F5E9F539F5A15735DAE8F2D642185EB042FE90743C1B651A13BD017CC992D53882C601180371E00
413D1C0D628425C20C0AB57C8769446B6FD2B9D1634ACC801E3EAAC615CA658DAD4884ABF6643C898957A94444C5D2BE23BAA5AB5E159244EA625AAA
21AC69513DB1E593BD9587AB0D9D36556EDA9968FBE8F3169B26139246BCAD62AD6551C5949A7874646B3C704A4423C4F1BA9B2047A234B5A414E29D
5FC27F1CC2B4A40CDB0CF93D601B83EB24EA1A3AB8EBF7CE2B32AD508D8AEBAA8352874C2C275D714661719F7997CFBDECB3E7661407863D5F6DA7C5
14B0C05FD17FD729B9BE5B5566C422887BDF604B4DB8A53161729D9AD61A0E34CA2A8CBC602E9551BC392813C38BF02D8A4C58EA6B91640A59385979
A88D130CC3342074DA5B5ED18F79836B36E1DED1E054C57066A78E894E1CFE2994CAB73748A9839F7C7A486ED8579E32F92798695E5A90C2B95C8E93
260992FEB9D28FCC7545F9DDBAADB4248CB89F36104CEDEA0667ACA587ABBEB9EAF89173EE7EF5893B7F337EE895AF1F4CB8DE6642AFC32752AAF3BB
7B5DC19D3B81BDE502D003A462AF8FD05B19EE6519610300D412BCE926EC8D29526F6B118FF98992823401F7F11CBA09E71ABADC8668DCED6B5833B8
3A6206B783411EB163FE73EBE1DA9A069D1A512D7A34948A1BC8B514E446B7BEB664A3E3053E72127B1D489438A602458364876223D56491A4E570C1
C0F90E55516E11DAE23C70804DAA305D60A17890AB039B4721C43BF12CAD36C80EE1AC269299764D8324EA1B548DEB125D77E233F20A4A2E7BE1E29B
6FBC714076F11D9B55D1F927F0D7993FF587C1FEBC97098D26DBE596DA0D9BAA63DE217D5EE7AD86B9AF4C9B54430C3E7B06F4CAC0511E491E43831B
741C18F033890433720843A40CF3C1DC40C493B0EFDBCBE4C3393B06E2EC8771400A9577CE2B4EF0B45C97446AB2B9357260D5C1FAF56B8F1EEA4822
853B76C8A7C33980404718232040B7D6F1FEA2EEAB1C83A55D02BB8188CE898C925B9F9D39626449E0A46BEF3CFFB86117BEF8DA7EE8D465B07180E2
98045DB9308D0802804B66C79B91C10624FEF8FFA9DE29A45EA701EDCC3B7855074AB16EA8BA64A4407F385E0461A430D091133B09534F48E53F8105
31028F209295470F35B5AB500430648D470A48ADDE7948A134A946434D4DCCC54AF98E75ADA103EB2BA20E4E1D68D211EC4F548C48B3CA43132D69D9
8873DBE0D206B215AE4B08864D24500D85F2A6EDF5249A5ED6C3D415FEAEDA899A3AC58237C7B6B8F5301C143AB8A7D93453E1444D1014528ABFCC7E
38AFE79087BE7EE2952523F37ACDFD2062A745DF8FC05FF29FB53CD23D30780DD9F1E467CD65ABDE5CB43EEE4DED90844C592CA9F0C8D64D9BD291D2
7D7B0E762479A44D750C33B250BB73BC333BBC0439F7EE5E46CA9B9F830440AC4E818173FEC986AE5FAC636FACCFD464449261AEA45D38199CF3C785
D43977EF94B4D5D4562B542ADD76445563714E512F3901E375E08B75855CE72B2CD99086AE41D82A8238831B97DE7C5AF7ECCC9CDE3D7D03C64E9FFD
E477112E835D383A04431E9129AAD93910E078568907F85CB6D836ACFEF27600A1CEBDE3A0FE75CDAB3842FE0F2A7F9AC6550A0F3ABC915E6FB0D8AB
459A54896B146206D3D3E3D0916B3B0CD6948714D830424C2519920955093653085B6A38C4686CC7B6D2986C25DA12B0AF486B4C407F23A172EB8615
41D876423BA7FD650D8E4AE0EC66B025995B4FA539CC78ECE2B54FC233425C012B8FF87B4664E3E7BDC216C210871886491D2D1152114534B8BE8EFF
8CEF16771B7CDCE4A2A2CB6716F47BA28609F60BFCB5F8471B6714E44FA9B2623F6D0A6135DC215307132559BD736F843B292B559F4A33CB65329282
0AA7109193DE420D4E0D25AA03AF404053D7E641B502C780C2E979F029D631E402ECB4830D9369B1B0667ADB796133A71C859382C14D23C6A5854D89
C160ABA8EBC0062023D4A404372FDFDE126E8B9BDE33526C50AC26C967FEBC92356906D546EED4F5AAD72F1BECCBF2E5E4179DF1FCA20F7757219486
997A03567F5A8918F3EAF016F4087B0DBFF0A4D08B0B4BC77800617BA34186AC106F75008575E4A06D88826132498F2BD02E80158CC20A1C5A00890D
1E39107E092FFF004D448EA6F2B8C165900EB0201ED1702A1589B4851A938A9EE491414A96B927AE5FDEC8CD92A31969D7D1A17BD1B50C851A105051
0421917752B9CB623A17140E7F33F85350D7746247CA136053E0A908F44C38D06B010B973AB782C04A53D3D1FF0F7BEF1D675759AE0DFF7EE93D2120
A0205590224D410405111044512922A2A208221EBB880541A94A1535A0F41208BD26101202E9BD673299DE6776DFAB3F7DB5FDDDD7B3E3797DCF7BFE
39DF9F87D9933AB367EFB5D6ACBB5FF77585C59AB7DBB003A3E80DFA5EB8F76D3AAD173F3075BFAF5CFED5EFFCE293930F7BAA3ECAFF33FAF86FCCDF
B963EF997BFEDCC950ED1A51335C3841A175C9DBEBFB835863F0E7E8067A6376B09E99CC1F72AD88875F096B030EE6FEBCD5CFA80E560116E573DC92
60F306B83E46ECC6F23F65C57510F637594001DB0186468428017010F6FE06EF0DEADE844CD82B0F0E57075E79E0DDDE8EFE00381BC1C3BA1C3C72CA
BE6FD9E38EEE3CF38263A78F193BF393A75FF8EDBFBD39249B311CDE2887B498084AAE2503C8297092C5347779C14E8A611FB0C2A9E50C9031652796
63430721F8876251EB1C408930BCABE0638DA852F11D8157CA54DFB0306033C27012B0A4065DAF0874444930145ACE81981B197AAA5E1E742A3BB754
A91E4A6268A0875983631A98D8E57DF26754B9803500DB14017849E9BA3472AA40A05B8CE624788CC2C08BE254A8265F211D588DAA95CCD293D07981
1DAD52441730611E5542014F544F4FD420DF515DFA8F85E474D7EE37698FCF5CFDC8A6F99F9A7AE95A3D6AFDA38FFFD7FC1B23DF993EE38C6D8D4C85
B1608E47152F946E6CB60E1E9EF64E8A955872374D7E6D150958115030CA76C3E879E5812411955A60EF75F4EC03AE3CAAF663406D53F4021886E34A
C2E2E9061794DA22790F7A2B14DE93BC010C1FD6F4291B477EABA97C16C6AA80B27A6DA4500B0225C9A89969DC3079F60BD1F2C76EFBDDE7268C1933
66EC94537FF1EADB3D2E65ECE580011108178516BF0A0B81B63861B2FD862D7BC9ECE01D528B29B267402F1F49A931B2A0B446A2E8A7E8C938250274
EEAEDBD40D471B24859260A247DA2B5C091E29CB35924AE9916710149F335DA853ADE106B0555EF1C03760EA45053C63B1635D31269F245A3B18280E
B01CC4182809195DC288823D90CC944CB00A457F96249C418C280B076B5CE754F73B25723F5227515F8F1B7846548A55F28861ADD85DAA38826A19BF
528F79E7DADE8EA15E89044876AE6A19289BDEA3C74F3DF3EC69E7FEE4A00FFDA5061FD8F86FD945471FEF5FF36FE4DDE74F9D7945DCC8BC540CB50F
57A5D3B27C24CFA8FA875896703B960DE7A9A4B08A2A1FB372ECDE21965A1D0DCCE9C8F4A384EA5F8F2234BAEB14E1756474282CF52E08A71AF01221
16F0535628717816AABA1526FE5007A63CD60DA8029622C4263160BD466203082D787B985912D66ACC73AB15FDEADE7BFCFE577BC2F6C74C38E4C8B3
EED8D851F3B456A1A874EFAA473106DB8945FD8B50E358E9138895E465C0ECD1803E2858389A8A405495A09396C149C1E280A385D3C86C11BE9B5F00
5CA6A02DC35FD830D0AABA6BE7305314D4EB355F581D415C0FC119F33D1515BB9D24D5CDE6422258777FC1F38CA24BAD46CAC676F2F2249498473A3D
5549EE53D8E1A3210FEAA55611D55399C975582A45197A8F69B133442F146A87F4692AB23CAA4042DF8FC8F8752D49CA23543CF15259EBCC195EB6EC
C5373B82ED6BAA525F3E76FC51AF5FF7C10F1E78D13A39CAFE3BFAF86F1ED18BC7CF9C763DDD153D0B8ABD0BD60CF5148677B5700A5226625EBD18A8
18C119907EC6A198178AE6766D4405B0AD4A79A4E9B3940C534C141E40328D44791ED983CF348475C2D050799FE84AC1E27BA97E1520DC15115800EC
F62D90879A72024115848C6DD77B777D9080252F73CBE5E1B9179C75E905E79E78E219A7CE9EFD91E913268E193BFEB3CF0FF616159E01555F2A4202
518DC27A40419B92101372AC108691CD276464A2D0408524032B4930D0E535F53F954955EC973CB08720A3C7FE9C9DBBDB4D04F0949117040A31897C
6507114A064EB1520D43EC3132CC11401F860F9D02E52C86BB39B622310A290FB58E8429928F8445063946644938C937614D907B22A222DD1DD459CC
B52D01C095E8D5229998EA90175B92A324668E02CF8102D361185B7F05D74957DA715A573B617B0F1D44CBC6A172E6AD7FE9F1B75E1E347DEB86A274
DEB8B193BEBDFD89C30FFE4B35DBCD2B3E5A008C3EFE2DF8A71D974C9DBCCFEF5583D7078AA2BF9715BBC18F49A1D11F6EDBD5562A380A88196C9EA1
176D40F847561A38F51A8BED667E1C46B017DCFDD2466D03B43FB668C241971C44243CDF2D2B70EDFA8A53869023EE4A2AF639268071F22F319C0610
0288A4D0F3CDA29D15C1EAD59ED6B615F77CF733A79E76E2ECB1F4983861CAE46953A74C9D3CF3F49F3DF6D6A086D85702268E34B6DE02C09D405A5A
51480061A127F0281B6898C86511472E2F19743F825A8569E013A01A424F574DD902F0836A0080AC1F40F6414EC3C4562B00BF29FA0741A499DD30D2
78939862B72E0D702087C1806889C8D029888B3B0606DA7B860B22CD31C48BAA354E6E95AE80204F65945307DD185D09CC2FEB4DE6E4049A2309D53D
CC8D13BFAF0A6664509CE1EDE88AD1BB57477A7AC9EE335D6EEF18A152ACD6D3B1AD2772C2724F6568C51619F0CAB040BFB4BC384CF3EEBDC78DFBD8
272F3AF5F027F9A8D98F3EFEEF9E9F7D844F1F39F198AB5F7DF4ADA16A07EFDBB0AB566494DC0E8C70E1B5EE1C0EEA8C0C9D43D39737A178E85A67B1
A5ECA0DC5F5B750F8DDE37D6F7906F8B00581CE4E07142412B56A0F14D39A3B4552BC934A35A7B37EF6DCC8BBBFAFACA1D6B1EFCE32F6FBDF1BA1B6F
BEE6D2CF9EF5F9AF5CFCBD4B2EBBE8B28BCE39F898630EDFEF83B3F7983171ECD871E3E8F724FA07FD39F384533EF5BBC7DE8DA8A000FACFDED8800B
41339492880C5BAD262337A434E736FFA754DA32E2426D90B211972B160618479249C99807946B23270055309D80E333AA41306D6701D97F0DB05FB0
8769CB26AE6BC36E18EC660520C788734A8DEF48A58372A0E1E628D5E16EBDD65FAA6ABF8A3A5EC4C04734F5859561BE0F653355ABF92CB29B0EF6DD
244F764F4E5361627A8E913E68C0D19D8C780A89444E8E2BD2BE0B1850D8BEB5AB2842D7359EC45C8275166A4EE80EB4BAA1E7E58D9DCF0EE8063F6E
FCB8A9FF71ECF4A39FD5A359FFE8E3BFDA3E95F8DBAE9A7EF4E36E2CB7FB2898BBBA59C3E84ADBF2965A1261D394E252CF90C7AAED0E859F88974706
CBBEDD704B3181A24AD96ECA35E976ADEE1FE5CA32CE77535437D2DD2C9F1A2C9FB14DE8C98AF2B4D8B672E982E76FBCECE4C367EF39FD43D3A64F9C
366DBF030F39F4D39F3AE3AC4B2EBFE9C539474E9A3071C2F8C913C6C0F0C78F9B3875C69469B0FD09877EE2BC2B1F5B5E2A76E946CC870B2DDB97AD
EF74BDB47956500CC69C2C06933FB6FD552005C668A85684A12A42C12560D3275309E7C24844FAACC95368578E81588A5C8EFEA096F590C1E745816D
D763161FF8940B80D9505BBD9E145D416C28E9A6484154F104AF166B9A393DBD7D91A2744602970C4103B8CA1CB34FF40FA8DEC76C2F16D8E2017600
1D49681623A9E083454F94FA5146C948584662E607E81A88E161A6B3866D5A4238998A96AE91B4C1C25851CDA31335B0644BCDD56E5F92AC7EB49F7C
E037264E1AFF81E3A79FF6663C6AFEA38FFFCBF8D1DFE26D771D30ED9B0513B73FD709884D6F09B3F881F64208D02E743A42E9169C24359E0F1D6E23
65401981D5C2C57FA1E9637363CC9F617621126B6D752A506B5A941DE6806401B68357EA78F7D11F5CF2AD2F1EB2CF070EF9F041A75EFEDD1F5EF98B
1B1E7A6EE1EB0B5F5BB573E5B2A13A594671C34F8FDD77EFBDF69EB5C7CC1953A64E9D3A6DC6D42913268C1D37EDD84BFEF66E7F0D2BBF8D465CF3B6
DCF3605775606BFBF6CDABD72F59B66A5B2183B60FF91F582236F6642099E494065059DD30401035F00C41211A84A294975044A6EA5DA529A50909E6
194ADA1D05AABC35B8B40C248E2C7B98D53A01E6D1585570A9FC286B58F4939D092639361D201FE8F9BDDDAE92C00BD32534425B8480F6EA3EC78A21
E5282945F038A6CC022556AA2170284C06963272004041155A2B5C958B31D4156B9440A1B692A5B2ABE5487B19B518A702077D864485413548D5485F
6DDB3F6E5AE89AF6C7DE728DEA7A6A6BBAF159A71135E64D9E34E103274EFCC452359AFB8F3EFEDDFE2D696DCFFD47CEF8D09DD238EFFDFD9945CB97
0E752EE94958B5A308140BE6EF16084BD531659B64CA755727826C416806523E10FBC5F606A7F21F392B20002E103468C06756E056A30780A89F3574
AD5278EBA613F71D3776CCB1BFBCEBCE791B06DA37B914C52408EBE324DB7D6CD57F7EFDD069E3264ED963D69E33664E9B3669F298B16350F41FF4D5
3B9654634BAA4747441EC54BA2813E37A2F418DBC2D5812DEFBE1B35129FFC8C4E592CFC4D7337A3BF40A932195BFE5FCEDFFE05B81278B260D28262
BE8C65E42188E2FFE40DBC6A9662D6B97BFF985C820521C74DD651CFA3BC3C8C73643EDA70DFD7E0142043C59030314DDE31CB5E826228AA5A0DE31C
D950EC0C8A38C6C282254B369295077CD4533A6239500FE4E19A5C28E46F3DAA5220669CD487BAAAA59A4FC70EF0222453A4E021B807305561BC3AE0
276CDB0E299DB5F73DD23F329C34C8E77D7CCA940F9F3BF3934BF5A8F98F3EFE3DFAA37535FCE239B3277D6951AC0B6BD77971656088973A6B51207C
E4C294DD2761A878CAAB2163896A56DA520ADF2EF048340429B0470A3CBF512D680ED2235799545618506C74AF4BA7BF8AEE409C955FBAE6EC138FDE
EBA8F3AFBF69CEFCDA6E63C75EAB163CE2A1E42AB442B8A56F53DA3F61FC94C993274C183F7EDC9831633F78F9CFAEBE63C17B6D8C7270A1C0E00174
0C384390A13BA1C37423DD7D6ECEFC5BFFF2D39F7CE78A2FFFE087DFFCD887BFBEA573E9130F3CF4C6EBEBDA8A55A712F2A0D8B3FA95B92FFCE9573F
FBDBAB6BD75460B72816B889043698CB8EB1364B564FA7275450E6185DEE8EFBB1E20CF33C0007A1660266C0C0AE3827F5AE72B1BDAD8E4D61617B90
E87E6AF20EDAF2196278E886D016A0AC00509FD8DAAF951A37A8FE0D6701A67FB147E9536A588AA5672C2228E131BADC19169262EE6B915B05D51C2F
1BF50F453C554E2192E805D2CF4C96BB231EBCF78FEB1EEFA2BC85C7A671CB94297B7D7CF251AB44366AFEA38F7F37FF3C29BCF6C5BDC6CFFE7364BC
B2EB535D5B4929B2D8DB44FA35261251DAD4A52863C63A3B598396960F472B2E0D583992B8D8E2DB9898690B45030B674A9933DDD591B0A8760A6C6C
DB23FD6442EE43274D9B75DC39B7BC352C6D9D4C773645386575818542461E9BE1B04176AC76CEBDE91B5F39E9A31F3AFA9CAF9EF2F1CFFFFCAED565
C0FE81C3173E95F2AE23C16E1B16FA197DA7EF7319D9E01E05ED8BEEF9F5673EF2A96F7CE592EFDC7CDF532F2DDFB473E3DAC50B5F5FB57ACBDAD50B
5E7DE289BB6EFDF59F6FFBC111B3A6CD400BF1B3CFBFB803CCE35A734135825440E171DBCBB08304459E29D7CC8A0F623000EC3DF44CB4E52C4F63CF
CF53BA2E7992A7CDA60765F9CA661AD1409425997103135421466E3720A50BAFA17CB7EA59C673CB6CB45BAD109E2063A10F0AC29869992A2C1B59AC
1F53117D9788420AF4E4191AAABCB3B54BD075CA62D5B1B9961AB3EBA96D9468388538AD0FBDFDDAF2A16D6FFCF1F1775CAD7A1E7F3369BC336DDA8C
D9934F5AA94617FD471FFF77F1CF565EB0CF9459E70F1A53F30AD5CEE1F64A90B0B066925C46C39D5DD57ACFEA5D3EAA78D1D351102272A8E0942145
490FBA3EC62AFBDAA5793054716E7413EA924622007CDEA8FA5048B77DECEF1830C9AE7327CEBCE6B55E87DE3ED68C5E80694E09823256B1434674D3
CBD81B86AE467322A94397EE722E11D7D3D08BFC88F190BE4F79C335AA0172DED7BDAE405EC8AF7814887BE65C75C585971EFDC10F7FFDB667D655CB
6180E53CC515D6839596AC73C93B6B5EDFF8DEFCB73B5A9E386DC6C419FB7CE6F2B9157A6110F8512E2180474E35F0779479C79A53A19EA8881B136A
B26F55F314EC38CEC96007DBB0F343D590965958B43A7A4DF5418B74962AE14EB52304BA4829CB32A29B9583C02E410A7E131086512897D02041A6A0
804FA2AC20F2C81544A9F4E9F7304BEDA481B22DAA6062DEB3A1BD14F22CD6232D3D9ECB2BC3237EEFD61D9ED1A5A58FBF19E469FB338B077ADE5EF8
C6B6ED2BD6B6AEA9C6F59ECE050FECD4F119D3674E9E75C23A2C688DDAFFE8A3695BB0FE64E4C683677EF8470FFDFE86DFDCF7C4BA671E797C55FBC8
9A671F6FADFBFD85A12ACBD0AF7270E3963774F50ED5846411A873D6ADA8BB75F0F65A8C0C2ADC9CEE5B1630B2366D11757487FB8EC6E68F64315472
2C77CD5563A77CDFA6FC718E3C9A59892E40D82D2B4830588C105A85C300D1C5263BEE56F06D658944BE6104BA785222362375A104427B117AFD1C3B
3D85B9D7FDF677BF3D61AF0FFEE8AD0DEBB7AEDDB9B3676747FB605FFB3B6F3CF5E2F34FFCF49CC30EF9F445577DEFCCD38E3B66AF591FBBEEB1956D
2CD9ADAE0DDEB0819D9B5EFDC35BCD6D5FAE98D2A25295AA6717CF7415B861558D521C1E8C281D2E620C98D8EF562C8F05A7B2C4F54CDC543D8A4D6D
84FB4033A596BB14915EDB991F654CD890CCD34C69CA34424BD29F68DF379631ACC947DA10A16B642203CAA9D035D5BAEEE28BDCF1C1A710BA650021
9D65EF2D5ED352F41331B065737B9CB3CD2F2E6E5BB6A98FAA81916EF25AA6F3F5ED2C12322B7D68C6EC69A76E001473B4F93FFAD81DFDC98282F967
4D9975C5AEBCF6D2838FAE6FD9B1626B77B954DEB460692DD361A2290AA771B5ECD22DC8DC10B97545245263E1DFF53DAC98C656783B499ADC5B391A
7C60EA276B227B4C101F2DAA8E5E804764BB66C9A7A71EF06C151B3EE42142E9A954C4BC1A014440B72BAB3A82ECCE24BE6329C4E91D51DB83021476
44F581905ECF50534F104E09BBC4D0FD56C0E46995D85E9EDCF8C8EDF7CC79F0AF73EEB8EFEE7F3EF8EC7D7FF8CE170EDB77C64107ED3171E67E679D
7EC55FEFBDE9BEDF5D7FD77A965951A0C419DAB2E295C76EBBF6EBE79FBCFF3E471F755340EF12521920A972CE92C8AD788A4768E8E362C0F438FA0E
E871D21BE6748C79D6C8D2C88492EB9A8BB51C464910667A39C03F59DCB4E8D8CA9A623E48190ECEC90E516C4B0076AF6ABE1538440F02934115D50D
E8C9534B242874B4AE1B3D02909F49B7568BC1C350DB3214D7029DF07AC9AD832D64F3C221E99660E019AA2EE7E5277B298D094BF91313A74E3D617D
6247B1A3B7FEE8037642F782DEFAA399930F789CAA49306C1B24EC5AD70406D391C11E8EAAF60E17FD6A472736F91B2AE31E45ADDD027AC8A8AD2227
7A5439F832333BF84F8D40598C3B1E6875DCC0F4C1A4CEC91D24234F5D72CC59BF7CD98DB195A2047D4819726DC0FA6F40F267F5409B9081BC11734D
954313138885BDD0F143D745D601852DABEB95404D9822B6F53FF01AD97F9E22D078226C7FEFD107FF74D763AF5C79C081471C76DC01871F72EC495F
B8E2E165DD6FDEF5CB6B7F7CDD8FAF3AEBD48F1CBCFF7E33A74C3EE4D00B1FEC28074230B4D56D5AA3C81750F52E18C4C3B2464A671CD6980E1D897C
1D4AE00966A5CD16668E2B02D680C102D405C0229436B0EB689A1AA4312547B85606EE036AE591ED9498D46A8E487BF689ED1C2067C8A443F93F644F
63681B073BEBC01DDA1D8580D3A1356255ACA48D2C5351EBCE91C40F9CA1EE92CC0237E53D9D3BDB7BB9E97BF8F17AC2DE7B68D8A4A74C9874E0A2AC
C92830EA00461F36FCA723F71E347DF6151B309CCF531DD507AA2385C1AA032C0C8B4CA05336D4DA4DE56AD85FA21C3D6086526DCBC2AB418643E5BE
65B807881F3DAA88372C7A266EB6FE2C85AE8CD1E48684073D45D9D81EABD697AE3AFE4B3F99DF679D90325192675403A760EFB61B441AB68C995996
885003F29E672015CAC9CAC10A889C037E068A7BC01451DD2F258CCA20125B4DCEC4EADDC33F602901CAE45AF6B57777B7F7B6EEDCDA3B50AA28567E
F337D7DC74EBCD9F9832798F19FB1EF2E993CE3AF7C8BF82C34C7A61A53652DAF2CE6B2BFA2C75215E807C6392F030B55DFD4490DD4AA808EC5E816E
A481D7BD63CBA0E23CA15484E23A649031FDA3B36AA42143A3135BBAD226F8CAEA88847E18451AE78BC2DFD85F8008591E61ACF3AB9248E8122B067C
85AC95B4166892C45497E8240F5A7A05647B63E3B675A8A0AF6DA47544E466646D7BCB7B8B362D5DBBB4ABED8DF991D9F8F583E7C4D9F313264C9D27
779BFFA8FD8F3E60FE7CD169D3665EB69AAAEDDEBE4AA97357C171DD81CEB61259BA007A2541B8A3281E61D406744A0C7649D0F18696ED46EF78D1B7
14D7584FA5381855A16F6DAB5BD8B12683660A503F68EF66869587DD04945FF8DCB247FE7ECDF9BFB97FF530054FBAF3E99654D8065616B38F8DDFD4
5A373665B0CD1ED7AB220C2C8E5F226B6EE26C8445CD60BE00C20E8BFDB77198C7D6FE13ABD445F6E20BEE3308FF827B083B7F496EEC7A00B2EDF296
753BB76FBDE68413AE3D6FEAE9FF7874EE3FEEBCEE57577DEBB3279E73DD91E3A69F7DD33F6E7F76D9AAC57D71CAE3B43E8063C4241E2F9D5BCC73D2
B7FCC9B9BF3EEBA3074C9DB6C7599B72930AB0A4C5C2F28190C72293338E80066906F30768210AA141062074C401FEA7184FC985A22BC6280F8AFB3A
01EF8761839408F042156956A0D84E471B75B78ED4A9A8CADC4D3BA206FD8406FA850C9278A8DF44D52CC9659D7983599215B6BFB078F988608B2E39
EBC9385D366DD2E41BC22CCD47ED7FF4F12FF32FFF7EFF99A7BC4EC1996D7D72656DA475109D7B6938875607FAF2991E68AFC712D1157CD9D0AE9794
DA67997603A3B9E4BDCB87A4D5C5ADB47708BAB343DF20FF45B54C99B7EF753DD3D220FF8154994B19D5458C10A8ECC64B9AB3ED2FFDE40B9FBEFA11
4C16D1E3D64C30CBA62522D5EC7891555020CC35F00574C8C8A0C34037B5BEC946784992A7012981F15962377600C36F24A530A7B00FA6610D45F128
F499A0145A08C0FF186A11836520C15918726301039D7EB5F0FCDDBFBFE6CAEB7EFE9B2BBF71CED93F7FF0F75FFBF084699FBCF49B5FFCF99C97DFEC
B79045EB5C12EC05C6CEE0BAB71FFCDB6D77FDE267177F64CA94C993F738E6CC9FDFFA484B42EE01230CF26491DDF943B6CE5566D94FC003068F552B
36F70801A862823C5019A353881D79AE16496997B29D425C698006A5AA0CD53D3F4AB5D0A27D5D994A923A1343C3598325FEA63699F9B16AE9CE8342
A311165AFDD2B05B5DB5E0A54DBD7E980EB52CDFE225FA95E91327FFB0FAAFE83F6AFFA3D64FC5F6F64B667EE0865A9EAA5094CB209EC1361F4556AF
7B98228A9F5026ADAAB528023E057CDAD86867B2A138582F3472624A1262E0E7A50E8A0509ED598A751A98006154A5C844FFFC6223C3542C4B984070
CEC9F61436FAC926AC2ACFE0C3177FF86357AF12766F20B6383974F3B0012F409F87282A7D669A3D7ACAB84128801E01F806040333382690A0F204E0
DDC92D1716F20700ED9087C4DCF10CC8439944DFD2C296319FCBC8B28C12125C5E1917FFBA2E3802A3DCA0FF823163264FBD60C99696359B57B56EEF
DCB1ABAD7B47F7DA3977FCF64BA79CF885638FF9E0DE93A74ED963EA8C0913F69DFDB14BBF7BFB8650610C60D71A808F22BFA4E18CC0C503C5031DF2
04FB9168F583F8C0128FA1374AAECDAD25F475100E6BDBD880DF63B6F988F5DF6ACDEB199132B5A2A6516F35CF65A12F848F1679DABA9E374C57BDD2
CD3C378F6BCB17BFB76A656876AEF3F19E516947056C21EB664FDDF7A7158BC11E35FFD187ADBACBFF3870D2E16F420B83B27606BA7A505A085EF79C
003B29C2726F640D2ADBADEE2E19752A29B58D45A0ED0E1A55DC0C12DCA85D298105EA8C72F5B4A1FA0B94A586B23EA29A4932D0681A5FC7525E6679
7E357EE9CC7A80FC4F1F98BCEF126869591703324D63A1B24DCE6DBB27870D606CE20AABA4D34CAD0191859C06368C28A9C0283D2A8456012783A7A0
1720EF618457A384864BCA2D225F72CF57B137C2419D2B18747BA2C0650A190A0689E475C8070A4E6EA1E72FB73D75CF591FDDEFC04F7FED8EBB7E7E
F6E99FF8DC7147CE3EEC8089E3C78D1D37F9C4D32FFFC19FFFF1D2E3974D1E37EEC827DED809807346257FCC02E0F481EB836FA48CC438E4BA28E949
B41F020A81393FD939B62680E305D699CE9353EA443981C4396538FD3C5799AC014F9C889ACB22CA55408D062403796356E749007E62EDD4581E758E
C4DC2D55F2BC77555FCA4A89DEDA9BA4E5751D7DCB36B48BB867585D3C75C615A5A6076D4E7D470DE07D6EFE6AD9A7274CBEA484F59E2D4B8683B54F
BCD61F05BB56EFA8049E1B82C137935E71B84821CCF6F7D0ED37517DD801E65D59AA5C54B81AC14B55A3B85A46AD9A29A4E831991AE7219301E87E12
F4A753CBAA69C1AE1919878C28EA72ABCBD568044BCF9B39F16B05406B33196A4B9E8B59D96EBEF0B8B95B60778D0D4796427E42039F430E258040A6
9D182881BE03E5175993B70F4D84300C354B641450062D06FA022139953A519C2A47625D21F4B9609A47919DB66596271091974A1FE3473863C58BDD
9BB76DEE1D79F7B75F3CF4E0038E3FF92FCF5E3167C5BA15DB2075A036DC71DA47661DFBEB22CEC2A88829CAF8434A5A18A40014E756115455141535
BB171F25FCA8B474E13868AB706CEC044009CF431F05C76E624E2E8EFC80A72CE9108FC1846254A0E39C32A2340947B0440CEED24C53C51F92BF0EEA
9954BDDB6B59EEC685F5DDDA7FE9DA67DABB36154570EFBDEAE549B33EB36B34EA8F3E761B7F232BDD3C7DFCE1F7660DD1DDD5D6D3C3586FCB802B6A
8B1F59C7C12863F3833881C6459A48A86EC368C9527C854116457EDFDEE5680AC4A9096416D68CD64144CF45C4B2CC3789AA4532768A6183EE79CBEF
49C1DDBE3DA5D7E88F5365EFEF78EC579F983179EF6B4A76552F5595D00AFDC8A63A97B1B5BD96B1D5F5C6729C6C2E1B634680467888B53AB0F3DA79
81E5DBD5018440A9EE4FA417441A74FF0C70A1A10A791D382C3A113B6F4F293F078C48DA86BEC00A2DC6F9D01B8E43487BB3907C98D27A378D4871A4
A04DF3F829C32F2FB8FA13077EF2478B5A7803DC624924A220948C8580E9A1AF27B1A96F5F33F54B20F087740F4578941FB1ED9DDA5E8B153EB01B00
54010108053502E515073CA912BF465E829C58027D84FE751BCAD036F202B73842E503EB786F690FF9896285E5496128ADF67A232C2F6ED855DDD666
D28D7FDB90E8C11DB2FCBBDB03FFB0C91F7E6BD4FC471FFFA7EDFFF909D37E399888A59B1998E7634AA0B38861031E76E755EA9CC120A80A601C1968
D55136FFC7E6891FC89C577822D0E24B2D060F0A40144999E599696461BD6DCB823757CEBFE5CE7B1FFFD3EF6FFBEBED77FDE9EFF7DE7FEFDFFE3AE7
AF8F3DF7FACB2FBDB068A0BCE1F9BF3CF2E48F0EDF67DA9E475FF3D476305E90A5E4A9B0F13F81C4AF65D7B30C1B88A2100CB2DDF6262728F67931B1
04FF5D3D02261F3C5E2999AC0002D7FA02AC1202686B737AE570A015307A6B6EE9D8780F2A1DF4E0E226454F02B690143D79654338F480B96D1C92FB
C33CC1587E0F32EF87F69D3AFE9025FF0766D0B06C4340F170A14319AA48A9540EDBCD623ABA5031CC04CD6E55B1DC2E4442A5070D027B66C00C6BFB
1A8CA90C7E80DCA50EB8C70321A0ECC94A2375EEEBB86F559153E5A4BA9F7EE8F565D5DABA55CB4AA9DAF8EAC0E0525FF961B8799750517F5BEBBA5A
12AD5C16745FF7824A2F9F34EB376AF4C61F7D34F1BEF1CE9FCF98FCF185715A5CBADC43FD4F5970260D0F43ECF353A92CDCAE96960ADDB34148D106
F577C0B5006F05F0FDB5BA349A03BD0E9E2ADBDF6E325D25417DA86FC5DCDF7DFBCC333FF5A5CB6FBDFDD6BBE6CD7BEBF5775E9FFFFCBC7FFEFDEE3B
FE78CBCDB7FEF69A8B2EBAECC203F6FBF8C7A68D193376F699DFBEE3B10D653A2C49A9032A7AE08363B4BC252CD3DA3F1A0CE14815ED32BB128F6380
41632C086EEC84977C6CCCA3B790DAA9218E48C37D71CD360590C6231FA25C66D0958B75C8A1D58DF96093841C22A0B2D9DD64453FB158E4D884056C
2271810FA8034B85C520720CE44F38D3F3668E1933E91B7F7BFC85A7DE5C3DD2BD76654767CF8EEECDCFFFE3D13B2FFDF245C79F74D285DFB8E2472F
AE6BC9AD8E52D25CE5A7F3A4BA26CFA2FEF6B5AB176FDABCFA9DE5EB563EFEA79BFFFED8E24D7DFDC57244497D8E8945026430D8D34CA2DDC8AB504D
82DC807E723229AFDC21C157AADADEEAA063EE786D839F88960D41F7EB4EBE79DEF68A9BC672E5E2C15D85B8E58DB5B27AFFF234BD7EFCC413FB47EF
FDD14773DD87CFFDD8F47D6E2DA651D9D769AEF950DDE4C6B8D52A0B2433CCAD498B0A4EB927EDBEAAB2B3366BEDC6F7B0094BF149201747CC6EF2FF
C7FEA6779FBEEBBB179CF9B9E33EF1A9ABEE78E5EDED03D2AF7B7693BF89AD07CE164098C18DEFDDB9EFF80F7CE1374F3EF1E4E6E6A26E1C28C143AF
E0A5361CA2FF9E06C358BA4F6D44A58C03DD72C07F20C4252527030F84CDC48130F2EBC8A565B3D988797FD2CC13A084DB05F12D20EF62CBBC835C9E
127ECCDED03E04638F65F1B31538BD5348A1D5B29227BA2A417640B1DF5639E006131877C01D712ED5AEBF7DF5B089E3274C9A3465FA51271FFFC103
F63A64FFB34FDF7BCF634ED867BF4FFDE88F3F3F72E6F137BDBB66CD70B1BBA7BBD0B16BC592F9F3E73D3D7FCEE91F3DE2D4130F993E6BCAACE9FBEF
B7DFDE471C73FC41871D7EFCD1A79C79E9A5DFFBC92D7F9F73DFB32B3B8638D85390E684E814A8FA8047A503A523318F95DFD72A72CA93E80BBE8A22
59EA6379103123B6B8C9AA7B96F4B338ABBDB9B22C073A56BCDAABE3056BD2FCD17113A73D998DDEFBA30F3BFA71AEDB7BFA9777A0590D3C8F196A1F
D1591E76F5F9911D5B316D016AD86CC50C0DF367A8E3688E619FB7B1025E4C28DA20A68241DBF8239BFF72E5799FFCD0FE879E73F57D73172E69E935
B0411EB92208C0B8D364EB4E937060D99337FDE4AB9F3AE12B37BFB539327669C02801BFC2C1C45918168955B7066386EE17967113F1DE1B51D825C2
E755D1C3F63D026324419107DC617504075BEC0649516A893B8DED215AD53E0DAD0D9B3A40963006BD8700035842793848C9AD1E09AA7E54FA76DA6E
1B0B46D6A80CC72EB2C01C0F40E638934A2692612C197190180FB56FDBBE6ECBEA0D9B772D9D3FEFAD45BB06BBDA0AC5BEAE216EF4D3677CE93BA71D
B4DFC1C7CC9C36EBF093CEF9E89E7B1D75E0211FFDDCE7CEFBFC17CFB9F807F7BFB078C1AAD6AEF6E14269D809B9E7BBCE50DFE6556B9EB9E1C67BFF
74DD7D9B82527F2110E816A070A9F407465639269D4C8474F5F38CBC9EF659D9E526A1DCCD98EE2D3C59FE5CB7ACA4D9F6275EAFC79B17BCB78A25EE
DFDE4892DE3D278C3BAB3A7AF38F5A3FE6BEF1D6F3A74CFF25E6CA6407C3A5AEFEB28E33D5B6C1A35A5859582F424B6E4A9D556D77D60CA5A21A0CBE
6841557781C0074666C57B7B57FCE96B47EDB7CF07F73EE2D3DFFCE7866D033C4DEDFBC0AFC063D8DDF82CAB6D5AF9CF5B2FFFD609871E70DC67AF7B
78719B009918DDC98285DC600F8E42B2899495D5CD9AA8DDAC015C2F9AFD58C90FF1E614C3C950022E1D0FFAA139BC909408955EDDA04611680AA07E
B0E5BF45DE228EEB7FCD03F4EEB122DC5B6D10D11DDA5A069976A16649BE2C022149A33A9EC7E067C007A42DA21FE5063618A57515428461A4D3FCDF
AF6E66398F70E6A1CFF8D0B6B79E7AE2F51756BFF6F0BC552DC3DBD76CAF3B41B87B9F79772966E722CDCD1FCB0D4EA99613686F64C4D3DC75228C36
B011159BBA5B8D42E606F5CA50C0A42FE12359DF5024459C448E11EDDB7B4C65F1BBCC54E27CCDB58B18DF30776D18C75BBE3F47E8F4F68953F77A37
1DBDFD47ED3F6F6443F7EE33E9C3736596BAEFBEBA69F982B55D9CA9D00CBDDEDFC8C3B5ED52B158964BB0C39A8B5B12E1D02EF21A6C9F4B6DE7FDA9
ACBBCEA6C7EEBAE8C43DA78C9F38E9F02B6FFFDB732FDC704F31AC9251A451A5C7EF6F5DD1DB39D0DEBEB977CDE33FFCF18367ECFB910F7EEED667DE
DC3E1C3541361885513A412FC95528B23C125E64ECC00FEFF66F538A0CB878C6390FB4A4A23B0CC946C3C192E4E4C0504D039A9FC4361F812027B68F
4C6C6180AC6EA9FC4D224D13149818CDB1970C114F93D67685311736D7009F482881B1C3566E2CC0FC2F3051008247D15B496155B9B00B0419417A41
1E711161E110B9080EA1116B5F801E49688B5AD2E43CFE751E999D4F02EB831F00B6FF313E889D42108AC011D87E022282AE7024A17362CBA534470D
028702FDB3447B3AE5837D75E1FB942E70274AD8F61D82B22651DEDCCFB7AC71587DD10A19778964F18FDF34993BB7857E66EFDC329F7299F0E00913
FE2846EFFE51EBA7CA7FFE519367FD7C384FFDA8B86ECBB0176AAE8ABD25C342130FBFF5541B7664CA1BDA2538F5E25C029E6E7BE9E890239B6E727C
1AA767A8FBC96F9C3275F21EC77FE1DA07B705220A2B1B1E5DB47D98C9A4F7966F7EF9079F3CE1E8B38E3CE5A8233FF5DBAF9F7AEAF7DF796FD99A75
7D460C8CB841B118663170BC966D4358DB5464600173CBFDEDC3E0D04DFA5ABA786DA4B7E25B0BA24417DBC4619D93BDC54C44F43FC62826A354C881
B041C415A8668C55BF4B4CDD8100776240509E26B2A26C59409540441535DA00062823CBDC496510F915380A7C16BB821AC05B2CE225D61B781E823D
04FF5C9935A25D0EFAFC12FC0300FAA5D8FD873C295D40ED27BABD828D5E308AC6D0FAA2C441870C5C25A1A2774E185720168046601E7A3A534E89CA
0FCC3AEC0EB0C2AE14F211EC2D68ED16AD542178545370FE53DE9F63329A29D7644EEB081FEA8EF8C840149602D6B2A15BD7B657E2F685DBE91B1EFA
E9B67864C516898CE4EE09E3CF288CDEFEA3E69F27033F993AE9E4A564528C525F8A956486EED01043261C44D5DE028F7C5E1F19897418CBB0E252CD
59E4885C899DB6DB8888DE99425750D40A1B576DAD2B1BA79BD0B988A25E525FB9F0D9C7BE77E1B9D7DE78E7F577BD5D767CD73E85F2F11D37FFF9ED
AE75F73CB85D56372C5DBA767D5F5FC9A777E1010B7AD6B495DEB8E2DB3FBB6DE1DCA71F7DE59D371E9D7BFD05675FF8C5AFFDF8C9252F3EFEE4CB0B
9E79E5A5B7D6202B212B05322F5411A832E145B2E63EAC55F852268B9BDC39D8033454A4639F4E16BDD8628AEB21250054FB0BACDAA18E27834370B6
35028C0F007C9D81B5B359290081A0ECE6B2DDDDA7B7F0B6FA195920E54102D982A04A065B017CC821D7E38BB8EADB02C0AB579CBA685634E4B3D037
8D229830900739F44F21739C35F5D22CD2013916552C1A43C634B3B065E9176AD83402F6210AB5554D53A9AED515288CA3E1A111D7A77384D4C0A6E5
DDCC19E8722C476092BE73DD524F77F50375B8F3E713C7EDB736191DFD8FCEFDF8B2D3C7EFFB67DE101B9F5EEC0FEC6C1BD8B5614D7FC04553F58E22
A70CA922976E144177423B6E02185C13FA02EB07A907EEFD66E71CA105F96C03038034C3FE2D907376CAD588D56E1E9B346E604B47535C658EEB873A
EADC55966EEB8615CF3DB670497B5078E1AF735FEF8DB63CFCF8968EF9CB7AB7AD7F6AEEDF1F5AF2CAF5373CB7E8F5379F9EF7CCFA8E57BFFF9993CE
38F9EC2F5DF98D9B1E5AB4A55C62DC8494410BF08C34C08941D581A90FB5ECEC0AAD82367A7EE8F867CDA6BF3096AF0FF25AC10E4A7784B60A6554C5
90DF412F10454ED6AC3D2C8749860C22B15E0F0C66281CC0E54D21D8F204031565B781551830088AA7B94521E44023D099DAC9A26661B532341250F2
A41465F0B2DC57A1BF5DD72216C08488742307A89959BE74E8F8A59669106C40102405F1B0BD7EA98800A472E904204A94CA08E04B725B3C54292B6E
DD51671DDB7645ACBC7ABD4A9DCD2F2CE6C1C6450379DA316FFECAC5B75F3C7BC2F8498F8C927D8F56FE71FBAF674EBEA8BFD108D6BEB87C7060C596
A22C6C5857D114FB4AEDD5042BA62C7485CB292BCF72CC9F3182B72A3749EC5723EED6183A00C6AD1ACBBB01643D0561BC3AA228C4012DB6869E83E8
8BDD406CF018C6F172828216987B80F747856BC9B6751CB5AFEFEAF72849F63C3F086C5F1105417530F83F07EF6FDFF6EE9285EF3C75CF2FAE7F74C3
C33FF8FEEF6EBDFDB737FCE60F37DE7DF7A38FCE7DF886F3CF38E78BC7ED3F75EA1EC7DE9E524210739837F4742D5A091CDC76B3304E585583D1C026
2F2CC66680B06A9E48E2B5889BC3062811E2C956DDD37A81269D898E312C4C2C17685386D02961530A58854866CD9C1DCF437056560104BC5E7825A5
C2FE2EAA0F12890B0BC60FA503D72E30B14A3FF0445A526EC54260A9BC6E073808D0A8A9085BD0D2F1B127CCB06D85D7874E30793FCA30EAA91AEEDA
E5542B0519BB8B5FDDE407C5F5AFBC17145AFBF36CF8A907DE6C7DFBF5F9D7CF9834ED3A77D4FCDFF7B97FF8E49153F6FA8725D881EE24F6EF6DA73D
12FD6F3DBCA44E0145864202941FE9063431B0B583BD5FC4C59805924C54583C4EB5CD45DD8DEF3794ADFA7D11FAD856010403339B4D6BC3A10CA864
4A453B9606F16F085C190EDADA0694B48CE5BC40B7CF6EC025941F4BB0D5A1AE4E73AB97D15CE547CC6D928DF81E2FAD5FB062C5F377DCF1A7479E5A
B468CEE3AB5E3C6FE258A8FE8EDBE79C6B5E43073E91569B4B185D8EA0C4A3E326F3089881803A3656EF4FD8CA41C7CDFE1F4CDDEE1A995A31F20B58
63B63A6636F7B180044C41C8FAB1BD9458BE30108D81A3832AA76A604905B270488373C420A2A3A880F780A882B6688384BE283C552D732C360596C7
C0F0A06E29FFC96918CF811BE53E720AE425C00AEB045B90B69D09FED04A7B55699168B7187AAE9691EF6821791AAE7E75E1B2A25BEBA6FCA6D43DD0
68941EB8FB6D95795E9ECF993AE59ACA28D9F7FB3EF71FBC6AD68C2F6F242375BC7AC6358BA35A5B6B81EECC7AA55809F31C0BC08335AE180B06FB1D
E65B8E1911597D1BAB646DD56F2D2F9789507B6756EB8FFEE10D84746BC75E248A0165B1DC152AD458B28F6A82850A7B7B10AB4CD0FFA28F2AB3DBB9
B980E63DF86FC93860234CDAAE029882A0C745D5448E24A4E12660FCAD843C0C251AF12ACD1A4DAF90A925179EF9E5C3A68E1937ED908F5FFEB74ECC
0310CBA91C305E2934E45060C7144DF3466A357AB19167B1FF5AD885610DD764C1064DAE62940B11AB0D5A086014420910404018B2EFA74D6C4162B1
C996B19F321DC97924501B206B091C5B78002601F74ADED3A43A84A3C05E03BD6126225B5134130DBBAB9025200A49BC9E3260C57052F0A46018B173
06E810C149495BB5F01DAD4C914F319C43BB1C94C335A72AB3C2D3EF7A2CA9B22475DF58596DE41B6E79C149D3979F75F49553265ED29D8E86FFF779
EE2F169C38E988177816F30D4F6ED8B1B595027DCBB2F6500302033C2AC5C348457EC80345C96A606CAD8BDE1414E914B2CE880147A32C3B9DC507A4
E4324A94B182F98EEC42D403B0D482F78371305599D41D02176F8AE57F4CF2C1744DE6EE512D9A0234C49A936F65197B6CED0D440196720DA6793115
C4C016F6D4C1E5530B05672222AF02AEAD4CF56C7DF99E8B8E9E3A7ECCB8A9275CF6D8EA62D00421DAEF27B3E54EDD6E2E5960BF8A658245C3C43216
2A2BEF83CA9B6C0CE1DF320BA1C7CE43D435F82E9C3534CD237258F40FC6B2CC0D6D819FA07717071C5BFA5009A0ECC698CC6A00A3002037672B0CF2
9D76973A6B04EDE484E266FD915A8DD4B42900A46C7260ACC78863DFB5134EEC306036902541DDB227FA0E7DAABA238AE35ABFEF0F54E35CEBE1CDED
3A15453A2A19B46FF45323FC211E97237AEFF6C7D6A579BCF955F28543F7BD19379E9D3CE5CCF664D4FCDFDFC13FEDF8E5B4891F6B49932009B6B5F8
81E38B62170539BB3C1F6761E855877C54DED2AF3B11B4AE54D36F6449B9A7226D0320088C0A7CC9382C163684EE388F02091C200F1905712CB864D8
1FC0D68CADA72992593C2F4A681EA0BA0090C642DBC1C1155BE10F8BE6C3323FCA0404418F4353845B869D046C020DA88A62E15F496770E3C287AF3F
EFF8BDF79832EBC0F3BE7BDB9C3737D6393A8D99CDCE73F268920377948076531936E29B20701D4F2A2AA9C1BC438706EF1263DFD6F6D3E8F854009E
5FA7A4ECCE31E5358CA194D141DD15D810245F12AA0CFB45A01E41A51ED121BA11223B38D2D0BFC7CE019C00D60DB984EAB0B4CD81118D3203EF026F
9BA5CDBE2218406CF281C11EE8C9720AEA54ECAB2687322504259D35786167379DF5827B2A89D9F46E953949E4948ADBD6178DEE787BE76034BC6165
1F0EAF3C5CD85508ABFD4BDEEDCEF2B8B743646CF58D8B939C1D3FF9E045A35CFFEFF399BF7CE3631367FDB88E5D3A29296FD48DD8EB1DD4E863D1FF
A2C07387EB0E546B293E052E7661EB4EB1ABA0B4EFF352D71093E8A5C1DE292B0F387A7C914F55A93019F0E9E0DC91A1409DEA977C8B62D1ACBD95C7
293796F3366BA0A31E5358D460D568EA5FA16B0EDB337E8D43F7C728A75AAC968549C27AD145E61D59064DF3AFF348824AFBDAB9D79D76E084F1938F
FCDC65D73DF87A4B29DDBDBB8354BA798B83430C7D474CD420AF9D2B5F82A7808E55A631C837445F55353B6918F2230483644F925983AB1C14A65649
97738FB1B038548AA0140ABC12D516E487ACFA1E39088F0C950C0FE41E961CC1F7EC3E22529D4642275503B2483224E9CD6681D510B408E6A6FDC7DC
77A8741032742CB5B785FF35818B89453071FAD1EC2C0F57350F066B94B8B4F72569D0B5A1B54F9037D2220AA9A0EAD8554C1311C7BBB6943CC9FA46
8ADA88B5F337C9A1A5D7DFB58BCAACEF4ED9F34F6CD4FADFD7D61F8FFC6EFAD4D31665E9C647B700011B7011FA32B268B9242353A71C9EA5DAAF813D
035BC036BF67BD034A3B25DF17E41E9A65691A92CB100140F9CC379906B806F07AB4B812B08150ED2E1865DB94A657478250E3953119A49B973B1150
6E4A013A8FDE980694167D76FB1518240FEA60E1308A01ED1762611FEB078D862EAE78F195676EB9E4B443674E1C3366CCC419976C65BBCF2F61C56A
C5EB68D9D2DB5B18A45CC4950918BE04C55F23299E229846519C9B4646460C518238ECAE42B630F005D93974C133CB872DD07497CCEA18809A370A18
63328C640A86128BF0B1590C06F59002AB820B985EA30916A22B12526A04841FB0844A18818501213015C4970D1741D8C4165BEB4F8C5B638C2A9590
7B81A22305CE97B2112F02CE5237950C93521FD4907015C8434426E5FDED23D8BC4A59292477E1F7B5597901BE7503834C0279D3DEB6CDF3B67B834F
3FFCE048A6E3EF4C9E7141A1316AFFEF63F3CF8285674FDCF7F7E554BE79D58FD6268D70C38BEB39592E575148A61F0E7607C0F5C7BCD03958F32C3A
0DD4FD6874DB499C4A2594A533ECE314E8D63391B6152A3C8BB24336FA00076E9AA7E0AEA51A14B53FB9823A455BD81D18BE4C5865E8BDB3802C0A7E
030CFD40BF801318BBF931B485333B5B270B2C0D23536094AEB83B9EFDC5C5C74E9F3469C2B8B163279DF3F7BF9E3A69EC84D9A75DF4C39F5D79FEE5
575F7BCE511F3EE51B3F3EEF33275F78F10FBFFB1F5FFFD6B7DFE97AF886773924B67224F5749EBC465E4562829FA197A9A03460A4703C74EE415D92
3580B49564E9E01A4DC01E0EB71452D260897A77A3F5733AC104CC26503440673EB54B44CD6E086672B1E75A67120B9F23B107CD2FBE03FABD8024E9
D8B2A200868478EF06765711F37D054650831CC82B520D80F101F4036DE98453B1E8A02472238791C7A68BC587EA3AD543BDD89DA08BDEDD2B1B89D3
313CB071E5A61A1DF78ECD83A9B76667E9EB93F7F8F8D651E37F1F9B7F663AAF9F3DE5BC355956E8DA31E2EBF24877CBEAC1ACC18A4E50ECDEDE1344
0124A5C2C81374EF8317D3B2ED61D3D5B2EB1B4DA134B58853AA538705A5CF02BD6CF4F912CAB2A160933461F139B8F84D02216ADCC314C65543AAAC
D65741518DBA9ECAF2FEDE90928308E83B95A6101841B41528C41303CE3EF23B52244BD79A40BD7CEE79E79C78D81E1FF8C897BFF48993BF74CDBDF7
3FF75E51EBD639BFBAEAF3A79E7CE6B72EFFFC99C77DE6077FB8FDE1071E7BECA55757ACDEB2FC8DFBEFBD7961D796656FDFFE9D8BAFB975E15AAFB9
8994FE27D35D4A75BFE5F10006D8508111DBC95F6E337EE66A06D6A01C89911660FBB0369BFDFBDE8C4DCC41340AC3679E541C33108B06A2EA228414
89A888A692AFC500409BCCB20C19AB1BBA7BBF01EED16E26609B2806CBA762A909597DEB486AC8ED2946298327A1F797A4C035A35E704A82AE11AE54
79D84B536F757BCD8868C7AEC82D3592FE655B77B4BFF8DE4095C7FCA5B7DD42D4D29B6F9C316D9FD746753EDFBFE69F65E1C2CF4E9C7133CB750DF1
D588C100F57416D7AB2A837A0526E02E65F44139B2F36FB283BA87186815ED13ABDA4177BCA8961D0AA1E0AB93B0D98C1C45A643863D16CBB69DCB80
02202514742B23BC6387D0EAF228A953D0F2A5E87D69054A3BCE40BE81BD3E74E3A86E3082553C4BD49BEE66BC58BDD144FC1F279F7EC1B77EF5CFD6
5A66A210C81BA36510683BAB6BEA6D68663BFAD84AA4842419F61ABAEBE547AFFDFC9E93C78D1BBBEF67FFE3B6B7576D59BB7AE19B6FAF59BFB5ADE4
932D93216A90ED6B9166542368A76D48509240CEA9E2A8A0CC286B68A4436E06B7077EF2D4BE4FCE763CF5C0DB3B16CE6B6F58F611AB35A0C34A1915
8F0E4B1EA607D810CE28AB1135A876C0CA130D7C242F4619C40FEC7C31B31CE7B660B06443F865E05092C8A9B811EF0FB05B00A673B5632B156101D8
42EDE2B184FE518AE64BE839F4BD232B7B53A69D55CB6A7E9A38CB170FE549A9A3EBDE15A274FB6D4192B46E495A3F3179FA9D66D4FCDFBFD13F19BA
6DFF0947BC4D26D6B7D98F9D5AAF8358A6829DEDF5CC56ED9407A3CB4EBFC918A2A6F074AD06D94B8B9941988AB398057EDDA9A265AE7D04A6849E5B
F7A896472C42AD40CFD9B899D27EBFE2F852F776591D6FF22E60B9A277149CAC845C4B8A991E56DD78A49AD64036CDB9E70AE597B9ED2390D9F9648E
5EA1165832DF38FD175F6DACB98C582842E88EA2CF97DBFD60CA26988CB86022EBB9F2B4CF1D396BECD831D30E38E9D3679EF6850B2FFCE4491F3BF9
D4F3BF71E38DE7EF3B7DEAF43D8EBB6011D5F59C83BA23C5FA832F2A5D557A052EB5EC2B1A2096218C291948F8EDF640CFAE95CFFDF61BA79E76CCBE
077CF996EB2F99C332909103B160793E7914D339E4C80998EF4A78502A8A8035A073518229BFAFD7035D2AD000E0F8454D225994682732A8A3506B09
9E262209FB5A6B505EC07695E7CBE2AA0D9EDD4EA2942453CEA04F47C4D17FA98F9402BADEEDEBC961EBDADA4D26158DEE67DA589C75BED5BE625924
E6DE524AA3477EB6FAD95913C77F371A65F97DFF4EFDCCC6AF4D9F72C120DD7DEDEBFB3A5AD77478B1E85BB97163CFA0C64C5ED729350F8CE5968D75
48B7B52F2C511D6ADA26D62735159EC85AD8D4A7A3F81C78711C327A6E48B71F303B95329AFC14F913B2760BFE156D23B8FF21E4A5E9E60B05460054
59080E0511851920E043BE97E479E7808922ECC9A61C7823DDD7AB5226A89260557A1E587429194F2C9B4096519A402F22144FAC0E48069F80A51C3B
19246BCDB3B077EDCA7796AF6C691B1C2816EB4EE47A3516374CE1A14F1DF491433FF9BD97DFEDB12C3E706BA9525CD936416E32640492D21246594D
D35ED2B86FFECDD75DF0D5737FF3E3F3BE7EEB9CB7D7B40F4BE41B26562AAC977D17C4A7B1F44AD28A6C2ACF17169B0B9706EDB154074C26BCDCDD07
14B305FA59588FB773D8B214D4FB913BD9DE096565CACE4139BC0E5A0B26107C78106902F48262A7D2D7D7CF535D2D9B74F8E5275B751E57D7B4EB4C
C5234B7AD035D9745F2BA5196DAB2B7E5DC543F33ABD47BF7A6D5FFDA089E34F2D8E46FFF7AFF9AB05274D9E7E83CC63E1862A914E4839F9AE073741
BED2542819AEEEEC6675A7B8A9D557EE88C728F23A029D2C6D7BF8DA96A986E2B91B35D9BFD0908A0602B24B30D4C596135CFBE500C96A88D218543E
9CEB009CD7C2E58C4A004D660BC48D6461C2C944A8FED650E0D0A25A3349BEB10D5BAE90C6907D9ED15151E6A8C60585757A9EC2BAAFCA2CAC26678E
A078AB98174940E752ABB89B50A49661085671E0F8ED9830B6E2DA947108C534EFDDF2D60F8F9831FBBC47E7B7547231B4A37FB067EDEAF54311A6F5
769F9FB2A4262B39588629DB4858D7BBB75D73F629FBEF35FBDB739E78B9C0E3FC3FAF6873F59F333F0A85C742A624E39EB18D05BC06655098126019
220B7AFA7D2962D429098054C25399E4963EC024DC0D00AC56CDCD6AE3FAF6DA62CF9FD231F21CDC05765AD0754BE9A7572C0C0382AD1DCFA8D697DB
A984E8595ECA559A0E3DF9AE89AB49F9AE558DB8BC751DAB6D73E291FB9EEBDCF8E8629DBD3661CC47BB1AA3B8DFF76DF2CF1FDC67E211F341499380
CD4B6054C57B7B701FD7DF5B3642E935E5AE42D677F54451FB164AFA4349A13A3199718A0E94FF9266148FED5D6A6CB99A272EB3DAF19955DC05E81D
3B7E76E986CC16C8760E75CF62212C73065980DD3CDD80BAE156A42769C6EA94F2338CF6620D2D10C8618A6A64F56EC92A4D16E760C2114600EC9B27
0185DE46658053EDC0B1CF1F5529EA5A4D706013955F0B588229436A5707842FC812055711B28725BFBEEC63C77D709F8F7EE9DC2F9DF9C5D33E7DDC
19279C7CD0878E9EDFBB724BB1B3B7B3C2AA232CDF8D2FF08637BEF5878B8F3C70F6DED3F7F9E8352F758C480CDF7423B1DBB809251C60EE8D424C15
C0FFC771053330004074105CE81E1D088BC0F699539D91C49E570971B1C883295750F2D3E4FBC9A40429782623C0FDD33812B6F7411E3522BFA01D19
1623408AE8DACA0AA3BC82B20045C790D3D1C09FB87D0CCB53432BB7415C356D6F494DB1D0552D747A49D7ED2F63EC92A6C543261DB67534F97FFF36
FEEB374F9D7255B191471D2DF57234B461D54838E0739762A9669EA77220CE621ED86532E55785521E7038D896096B819D6AD91D786CFFAAEA10B70A
1CE8CF37A0E5D18831DCB26A3D06F81F00637D67B08E9EBA91AC503728864365B9C01BB2B461C7DA05FF9C73CBF5DFFCE259175EFE8739ABFB4245B7
343004969407287BABFACD406C4546A2A192C184F218C0060A3E4CC7A0E3D6CA64AC6E8D927C1A8BE03AA4E77349713A37A91D3E2AE86B18F2169106
EA3F74572E5EBBF1AD57E6BE3AEFC507FF72ED6F7E76CFAF2EF8CE6F6FFEE3B7BE79DE0F6EBCE29BF35EFCDD9F1F7BF2E95BBE75F807C64F183FEBC4
5FCEDD39C8811FA2B30ABCFF725DC1D42195DBBA600BF4114346CE2F4F28C25B2A13ADFD01A758B32A40CD0D27EE614959077615905ED0160794EE68
4A82122009D1808D2CC2DA8A8340C408D71BE948827DEBCA763FA1A25FEB5CD3C5F672FA26608F73BAFA7D9B38E62A41A958AF74067E77678F341BEE
DD48AF94D613E7E409B3178F9AFFFBD6FCD39E2BA74EBD8BE268B8F089453DBD431B56ACDCD551084CCFE6CEEE3E910093D350AADCE769612536B1C9
AB38A8B56C819C5B80AE6D59C759ECEF5C3362613E98806799DDBDE1819D10A06368C06501C14B0FD33553DB584828BE25C38396C6ABBEFC813F9F73
C4AC430E3EE680233E74F0070F3AFA5B7FDB3CACACD80E1D2785778B884F9B0A438960F4CE8015B092872E3B0BB10A9B6598A5030D887DB88C72142A
DB95C354B362065A472476C780BC8A8C59CF0017E5B25056DBC3CE3228F1DE8D0F8C0DD7DEAE1DD55AA1B3BFB350ECEF59FBB59953A74E9C387E9F43
8FDE73CAE443BEF7C71B6EF9E3DDB7FEF3FE1F5D73EF7D3F7860F3CE52D51DEA5EB76EC9DBAF2D5CB2F4E5571FB8E58A8B8EDFFF0337FA91CFF5EEDA
2065833D4B9F9BB76ED80B90E664069DFF062A1B0ADC74F03EB902CFABF506541BA45958A76B5E0B61D57475F36239430FD51207DAE562EB79C15F4C
DE37A86273001801D3300DED76B4B22C028A30ABB6D679171B5AB2D45DBBA196C41BD7F83ADE3877981C744FAB09CF9C34FDF651C68FF7ABF5E7C9E6
0BA74FF92BD2737F90C96A1FE5E2458AC9517DC7F2DE50DA4046C6247C46C6EF0FBAA5029506943B7B3C667186EE1FF6DDC1A08D343466B500B60E20
AF0E2987A54C14A18F4558A51390C24E557F2FC72C81227461C950C255DBFCBFACCBD58639571C3C63C6D19FBAF8EA7BDE68ED1CEAED6E1F2CDB7D78
09722DF8934C8185C8E262A46956FAE1CE629C5B8BD71C0297563B134C7C09B34FB16C3C5416C02508EC0E90790BCA4FC228B193461307352907DB03
C555E04AE9FB6110851158FCB0DDAC988A9B0CE3A6B4F681DF5E74D8CCFD8F3CF9AC3F2CEEAB54D6BF72FF83F7DF7FD7DD7FBCFDBA5BBF3D6DD2ACBD
A77EF8437B1F7FFC69279E74EC672FB8F48C933E7FF6270E9C3AFD8053AF78B313CBF8860583C5CA8ABB7E7BD1273EF6E9AFDFF372EB8047D549035A
C5491E032D9C001F617D9DD8D15EAB568536BE0E5D4DE5821FC914CC24C950880500B457B81F63A91A0069EC41A626B4522C49522F4A25851F7BED5B
6B28B4CAA56478D1C6BE75F570F5D262B18DD29E6D8F77C4D5D71644B15E79EB9C5A7CF1846957FAA3E1FF7D6AFE8D68F1F1B38E7A43C77E570D9365
904F51A0F6A2DE8D2B4B59C3386D5B36B4452095019B757198D7EA66377307D25E89DD19CB2F03FEBBC4D8C19FDD644F7323B4E5B9503C4D3930324A
51C0D5523865AF1601784BF9A9371C392B37B47B4BBE78F03E879D35677DB9CCA9B4C877D3CF836D4F3249656F146AB007C596AFDB4ECF31D54BD106
B0F3090CF7E08A80D0CFA4C1982EE82F82BB0FCF845A0E0261B69B3DD784D2B2E7324CD6B5915C30898D3E154641E493652A1EA74D445FA2C3C2A33F
3DF7F829E3C6EC75EA4FEF5AB2AB1E91ADA8E62A506E59878DD7BA69F38AF34FB9ECEA1B6FFCCBBD0F3CBFAD6DF9ADDFFCDEE91F9A75F0FD2DED5B77
BEF7ECFDF7DCF2BB2BBFF3D9333EBBEFE489E3277EEEDD32B84CCD7F89B919AA1106CA01E94864371586655EB82E8E7D030B10441595CABA4CB96B90
0B19BB26487F805D54E9540DAC1AD03AAA55CB9DC3982D46EF3ED9EAB46C2B57A21C8BD32371DAB2F2CD96EAA679EB387BEE9A4BDED4D94F26CD38B7
6FD4FCDFAFD1BFF4D70366FCBA68D8EA1B1F2824D5F56BFC58389C955BDBBB5D8A89B2B879FDAE6E1F9A91A6C96D97DAA57E6862A560A0F5438B5C69
82FAE96EE47643CFE2F35244377809EE41AB475A8A0A66D88063EA65CE5475286214C7EA4D55CF25C77CE8B093BF7CCF2B5B5AEB1E59860149000835
0C39000D2D5BD08036281D0802AE7D66F35CA41C39FD6967EA064130B65BB5F47901931E294244CBB20F213F8E59776977C71F52A0B1A9750C0BE019
2241053A34823070F482508621D70140F2196BF9F1C9FBCF9E30FB2397FEC79F9FEE0A2C87591E732DB90B7401D97124A9B440C1D033088921F251FA
9683A64D983465E6ECE9538F3964C6D871E3A6CE38E4D8B3BFF4AD4B7FF8CBF32EBBF9A6A30EBFFED69BFE78C71D0FFCF5D975CB9F7EEEA93977DCFC
8B9BFEFA56ABD04224762E48AE4940E80F4C1F94E9905F480D5D64039F8A654AE9282BE801F232D3642EF2EB8AAA7DF27BFD654AB822721164FE7479
3BDF1B48E0E454CE2919DA32C4E7DFBBBDB0FD9DED45BDE3DAB31EAD66F93313A67D7CD3FF829DBF5107F6FFCBFCBBAFDE63E69C34F18AEFAEA82BA7
B375D0F3B58EDD2147820BD3F5982DE093547A837D1EEEC3A695C500F3480AF6DC92E0475EA06CFDCA1426FB362251B12EC3244F44AC23AC0A28E36F
585AC6EE804CADBEDED030104514BF05E5B45963DBEF2FF9C241B326ED79E0695FFCD5139D09E5C3CAEEBDD6CA9E36E0AF30F50414DE8E63A87288E2
381CAAB8023D45AB86D7E4184A54D9B5DC9D50E7B105726AB14B29043CE2B0B5D722EDC01540A7104A54325A834FD8E73EA8872C213F7C8BC99DB6F9
4F5C7DC1E1FB8D1BF7E1BBDB4A0EC2BC91B194CC0A8FC6646D08CF681E6A83FDDDD85E27E8F0A97F9E7DEEB9979C79DECD77FFC7D5D7DEFFF81B6F2C
5DB5ABCF7783AAC92AB53C9EFFD0F3736EF9F14FBFFFB3EF5D77EB05DF7FF34B93A65FF0F39FADE8AF65A8E29DFEAD3DED0BE78F049167F5CE311FC8
AC2A397AA726B700E8D89E5BCE79855236068FA18B8E6532076B88B15C0669E49B3C55350D4CB3E9EB5789697FF2B5EAE625955414D264A4E7BDFB57
D0EB6D9E3DF1F8F7FE3714FFA3F6FF3FBE6074976C397D8FE97766B96879A7ECB0DE353B93FEF515BAD3AB9ED596D4CC1874F013E177B7753366F1BB
5A01A5037E7F4062B800FB5CE884DAA51C3352D8F8D786EEBD4437DC256FB7A439F0A9B109390FDB57B5FB51E42AAAF853A4E70EA3D702CA87829D48
1A54A40F6E7AF9FE5F7DF9B8FDF73AFEC73B3343157C26DC626988494EC9BF2E70BAC17984D69E87149D813707FAE0293A659673C8E852DDB61F9003
58A6FD7FF16700650FC610A063055218294332E930240B52BE00EB96D18E5F1ADEF5CE93AF3E3DF7CE3F5FF5F1BD8E3AFB86DB5E5EB4A54439855101
AB78A14FC75B1B0EA94417E43EA4EF42D3CC1284A4D01B4495414ED24A7B80201C5989D52DB2968BAB0E767EBB1A8DF78B96FEF10FF3BEB3C7CC4F9E
7CD41957FCF64FD7FEE8FCD30E9E3D6ED2E4B193CEF8E657BE72DD2D5F3BEB07F7B5B4BDB10DAE95B2A10894E531A23C048BD238603E950841BDA653
515128050006D699D242871EA363A8F6D5634A1558576790A9ADAF6E5132C833B6B93ABC60193973EE65BBF69F70C092381FB5FFF7A1FD535C79EF88
E933EF6E642B160C0FF5AED8D1E6F8860721EF5BB851C7D148FBE6AE8A0AA2DA48CBA63ECE138A2226726A5437538A60600149A69465D3B603FB5853
ED0A7C3A85426F57854A02E6CB4ECF72631ACB716159F064CDA2D4A1D201F90B1BE339543738A001396C58ADFAC894F107AD48C9B278A52C58C9AF7B
3EBA79A9D28E4775B90E7C65F943C173955A7401D51AB16E928FEEE6E6DBBD1D1FC7B9D51F82FA98AD5110C8A5DD4D00934864F9F4050FDA5EB9FBAA
AF9D74D407F6DA7BDAB869D38E3CFCF35FFBFD7B35DDBC5010D8965EE009A7B32FE2F59A4A32031664439F328657BD4491A7C49C01FA02D6B819B697
C348862CE40C4804C9622A0DD400E7BEE040030BAB1E561CE9DBFCCE3B0BE7DDF1876F7FEBF4933F79E6C9679EFEE50BCEFDD60F2F3BE98893CFDA73
C2F8B1071F3176FAB58F3DB2B212FA8E156385C6896558D3419654DC847B8500EC0A7466ED1B5D9DCAB05C2E453CD30EFDE446449AE8819E1A181DEB
063B89C13B4F3DB6ECF9ED51525D3E92F61C3AFE438BFE77F07D8DDAFFFF74E89FF305074EDFEB854663FB83BD718D559C9A25EBE8E95DF6529B3BB8
63D5D69D5595387D4E54E8ADA03D458925F76D3B9AEE6F118711853C8BA0A3CA335171B31BA7A099AD8C57A788CA018D6DB6B65243352B790259E81A
52023BF219B6033DA1ED2C9F53892A850F8E2C7AED9E5FEE3371EA79AD09E3E81780DC8E279117D96DD88051A9AFEB0EA47728D980F26E66A3BB116E
60B0F986E21EC54BC8AC0BB0C30974C998B06C1AE80368CA531C41BE884A029E280085577C7AFA3E7B1F7AE8295FBC6DDE9F6F5CBC6D70A869FA74C8
5C727A93D04A0368C9ECB61E7D501D60ACE08F090BAE6032D7326D682B0A001150CEBCC8ABFB21248015F7992915A9FECE558F27E8731193E02C8A93
7FFB8924C605AB11B83D801CF6ABF57ADFCE2D8B16BEF1C79FDDF58B634F38F7CB97FE68EE8B2FBEDDD25F7640909242FB23160CBD12FA13AC617E4B
85CEC9A9F9D85830CE40299499F287DB87A878CAB593255CC8750B77182A7632F9EE9DDD4974CAC49973B2FF1DA6336AFFFFB3C23FCFA327F799326B
5EA3BCE8ADD624F183CDBB862B7D83BD6587AC83FB91D5BDF3870A16F106D14D4AC3AD022E3AEE221680C1D273A29ACC2804A20985CFFABE414C8559
0B9BF01A33B06DF1EB0FDFF7938B4F3EF58C733E73E2478F3EF9B433BFF9839BEE7D61719B47F6CE2021AA29990E42AF2E296EAB6597ED316EFC69F3
646EF764A1E00B825FF41D94B16ADF94D947DCAA7030F41AD1CDC3D715D7D6BEA1B341B62E1C24C0A965DAC67A416A29FBED261D2CD477A5A5EBB5C3
3443B65819AA783CB0D265266E16C47633875E0FA47D16DDA7ADA2486A784C39755F7F848A9C7C5AC8A2186A7F5E586C2BD1B3B02D20A15648975152
798359A22277094F84DC80636150A07900ECA2015160C336F629B3B72E45585DAFD4B204E3279599E296857FFFF3B7AFFAC977FE3FF6DE3BCCAEB25C
1FBE2E5226934927F42A554247511445388A8A5D8E62C58A881E1504BBA8800D4414051144E9254020A18490DE13489D4C26D3EBEE65F5F596D5F7F7
DCCF0AE7FBFD79FE4E662E4599CCECBDF6CA7ADEA7DDE5D3D77EF59A6F3D23F33F8B538027D21C6D01B9B094E2DE83147318FB554C3902A50C0B4082
DDCB86FBFEB3B056EC0DD26077D11DFAED8F0692E4E6A9B3BF220E90C89988FFFFDB4DCAF2AF34A9DC366FFAE18F98777DE1D94577FE6BE7E287960E
ACBCF7B9DD7B3D7AFA128D674E05D6AE379A10FDCB1A3DE3C3DB46C2849FD7480FEF3080697562DF2E8D5135A06D3702472EA3BA80329B62170C9BFE
7FD7E3D77DF5DD474F6F9FDCDE36E784F77CF4EA6F5C7FD32D3FFFF10DD77FFD73EFBFF0E4932E59198439678F12590025DFF2C6DF7E7866DBBC4FDD
53C586812214D0DE48407F17EABFF8F7087E1954508772AC0676010C81A1E6C7D865961DF31D2AD7812C46926C8E557C992B86C2E8839D80D1252B90
E7D8C44BD139835388BB8708B580542A4C7BC7005D48353C8E23E09F8D7A40BD37F6F3A187EA5F39011F15F40B8E60247FB3A9AA0DDF73A196023D30
F819382EB8053262D314B60CF0FD4078BE436D824E52ED47196B06D0BD833CB8AA876C2ECE6C67F4294AB1483A6F30510FD57B963CF5D0F3D538D751
888A0D0602D1DD57F00B0A42BA16F891A4EED6115610881A1515668D57160FA6236B2B41A3DAD8B6BE585C7BFF5F0B49BA6B5EDBF9C503257026E2FF
FF9EF9A12433F899F6E9473C3874DB8DABA3BED1758F3FD9333CB6EA75979EF2A17202C59D984AEE66D1E68D1FC5B91F69FAB5C0AD392A71F6AD1E92
50AC0D540409FEA8DE35E8490394D65604CD4AD1F3C2ADDFF9FC173FF08ED96D6D275F78F6A5D7DFBF6AC5B23DEC0D889E1CA79037BEF5E9CFBDEDD7
F7DF7AF7FDCF2CEE162D555EB5EC899F9D3F7DF2DCB7FF7623151772BF33109A87FDE89F7C474E71AC90270367D8817706D0AF09F48130ECC30110F8
768021820FBE803356A1A0C1083DA1841BD6C61D36E981887E1E614005B047592030A3C8310A5AC7ADBE4EF433C97EAE3E45ADE1B24338F7187186FC
CB2E00111F1BD84FB4B4C0F85FC2D024F05D201DC2D0F7403B34CB0A5347DE3DC0C514AE42984DEA920360630B5BCB38A3579170F00B6D6C5423B61C
C49015CA22C027D2BB33E001AA85ACEB1D3B7D0EFB1278E658D18BE1A0E0790C920CF66D68601A39F4EAC6061D1CDDDB752B88C5D0F66183FA82ADCF
FDFC9FCD34DD7AEAB4E3B61E309C9F89F8FF3FDCA2FDB93F4BFD4517B6771CFF7069F36B8343DDCECEA7F7967AEFF9ED7A636CA0E23A4E71D8ADEEEB
1E1F2D072AF0A4EF1996A05A38960D51EBAF86C2A5FA1ACEB830B5A118913A680CD63DCFF150EFEBFAE04B3FFFE45B8F9C316DF2D48E051FBAF9A13D
E3452B417052A2961E255757A38B48E837979F79C373BFB9E5B6DFDFF0890F7EE1F24B4E3EB17D46FBBBFFBDAFE8D546F78CF7F5F78FF7EF7CF0E7BF
78F0AF3FFCF24F6F7BE0D737FCE481C5EBBAA92D014538849036221DFABA98FBD1C160363C6CC3517E53FCFA05537B1016E10304DEDE9E395AA8024D
0BCA60CE51D414571493D2C2FC00F019EAF443D783FF46C6BBBC301F5962DE0E7EA3A6B65C700B81A88B58CF2BE48A2002C70980214C09F01FDF12BE
CCE7FC14B7D68019B88666E03F12B51476A5AF81F57ED6A238B71AA873D230B74C08FD8A0D2D4245FD0CA8D321C6FE0ADE6A3AE0FA05EF00A0B3D910
1A9020ED0A6909AA6AA45937B8491AD935ACE2C8DCBB7871277566DE9E3D5ABA5DBB5F5E55A7C3A17BCBD67D6B8675F2C4DCC9872E4927E2FF20CBFD
1964AF377D6656C7517FDCFDD887BEF5EB47BFF0F56BAE5F3A78E7BB3EF76253579C66B56095AB6341A358A122DBB6ABF650D79856A26ED646861D21
E34C31EF47F9553F4B45BD09354B28E8D2B32DCD6D4FDD74FDA74F9F7DE88C696DC7BCFBD635450BCA8054F82A2514FD145AF7246931432EB12A5B5F
1963CFBACAF573E7CC3874E6EC691D534FF9E2B51F7BEFB9EFBDF4B357DD70C70F2F7EFBDC234EFBEF4F7DF2FD577EFAFA1BBFFBFEB71E7FDC65DB29
80D20CC87E6610C33A4CE55D74845D9C82C7068E80C82A35C1F3CF1974144D76C392E572944AAA86C5284F0AC04AA0DAC62DD990E880371F7A0E58F0
068C5F00CB3634EB801A0579651147024ABEB93BA88472511CED670C50370E622FB6121CC6A1707AC778DCC0BF1AB10F18EF2251361587EB2684C871
12A21D89192799B0B22FD71FCAD781ED61CE8A23054820843CE4FF2360194208A15A0D9566F0138B59994945F6DE9D1549958EB36D73254A03BF6BDD
5687927FEFCA7D6E34BCFDB52D4374B9DD2FEF49B24C16C297E7B6CF7DF80012FC9988FFFF5BFAA766BB78CFF13366CC38ECBFAE5A30EBEA871EFFFD
83DF7CC65C72C5D5EBFCCEADEB362F7D68612DCDCCCA1B5DABFAFA37AFE93155A56EFB62604F73B45F02930713BE487991ACF9F42F835D9ED240A6D7
363FF9D0F7DF7FCC9CB9679F7A4CC7EC33AF7B70CFB0AB233800FAF40C23BCB8E685976ECCD8DEC8A3275AD577AFB8F5CAB7CC9D3BAFBD63CEA5B7FE
F69A2F7CF08BDFFFC73DAF2DBEFEB6FEFAEE378606AA4D68104BCADA5E6FC16B56299BE95ACFA04A3DA754EEEDEA81894D3E3D878C06453FF5DA5206
A6E3783E8C4AF046C2547150DCD5ED61584F850B15C702CE3F8D1A048E23DB829E8939E4323B887A778873004A13436BC0B3D155034208D63E1D0F90
FF0B0C1D879E6B07B9D379C49E2714CA745B580E2D46146BD3E3ED03B39DD9FD34411D31528CB5025419DE7C3CBDE4D961AE6FC427108FFEB8C68873
9445E0413A54B1AE27CC97711834478C4A0DA8A738CC82C8AD5B747C894209A44C6FEB1B514B788EC60224B0766CED312A3DAFAD1AA38BEC79FA8912
9D536651D7CE9CDA71DF81A4F73511FFFF875B84DC3FF2C8451D6DEDB38F6B9FFABE1F7FEBB6EF5C73FBED975E77EF0F3FF2EBFFDCF0D6534FFEC6DA
671EEFDED457AE74F70D96B62CDB323A5EA978633539D6D7305CCA25D4EF9786FB0B02DADE019E792A98E91FAFDF72E1FCB623E7CE9975F8A1F38FBE
E8DB2F8DD3E32EC60673AC0E3ADE50C2DE42944CF0E25B2D7BCD6D3FF89F5F7EF78B972DE8987AC894E9C75DFA9DDFADEA13E8B7D17BAB48F4EF1341
2E02D082E69FF3D2E7DE76C10D5B4DBBBCEB8DA5ABB78F6CF9C3C7DE73EE39279E79DAF95F5DFBE24D37DC725F3FFD14C58AD2C2F7A9CBB05063EFDD
49C784CC8419C67EA01D47BAC272552BCE5A96C4ACD0171C61E0154BED39B91F6914982E9524662DCD63374736A7663581F137E606102A305040409D
047A4451E41703D4FF29FA01BFE6E037199A4435BDAA99403E00A19F303AA84CC78DC6B01FA70D3AF690994D3976001BD2607F3D91C01B30AF3AB4E2
956308CF0130A0A5D6AE1B2A9FD5C8A3D0AC0CF65B2C158A032DAE8C24B007813741A96054C6EC24B18A237692969FF9E37A3B1E7CE2EEB54AFD77FB
8C1FBB13417370C53F3D65C6924BDADBE69DF1C32DBB9E7A71B03934F4C6DEA717FDEDF7CF6F7AE9F9A7FF70DD4D372EAFAF7EF199DBD677EF1BAFFA
61A5D8B9E3E5357B87779585169E875ADB29D6AB76DD021E5F78490A7D4D6FD7ED1F98316952FBAC638E3962FA61EFBD739D49F5782C3D591A118168
BA9110512004BBD45A56C3AA544B0B2F993763CE8C591D734FFCCC4D77DCF5F4DA1E0F5797682112086E9A2C0BEC61420E054D30DBAF9D39F7A8991D
A79CF7BEF32E3AEBD2779E3D6FE68C23E71D7EF48293CEF9EF5BDF77C62F9F5C6DE61EA15451A0CF904AD922ACBE61B6CC31DEE347B164A93217B81D
2A16D88F08FF4D191C1005E8D5C151D049E8035CE7D46150CC2E5D10E34C9A23418EFF8DF2323EC71887C2303146F460718C7D1DC8C50D2F1F2D4097
379082CA09964967C59E10CC2638A060A909EC14DB81307D2AC9340402B4A772EF2FBC182A095C27F50EAC51007742BA1ADFB781F7891B8D28A1EF15
4AAE0D4846044575151BC5A01584BEA5D268EFAA42998E2BAA86E83ED6D6DDF9684F24BB6EBD75A397DC39ADEDABD69B24AB89AF8324FCB3B8FF86F9
EDB37ED669F0C05C7B55DB4FDCCC70A25698014FA2E4D6A7477BCAF18B3FFCF7B2871E7CA5D83087BAF7560C5F794D67B4B34A19DC55610A2D3A30E4
B02537167DE2C829933BE6CD9A3B6BF6AC73BEB5B44C6FA4FC66A98941BC56EE58D9F1946A1655CCAB3EDDB979D3F0E88BB7FCFAF6071F7EEEA5CE1A
1E7F4AF3592B082028125046A52E2354ACDE01977B200BA924187965D19247EF7DF2A9579F7AFAFE3BBEF1C9B3CF3AFDA28F7CE33BBF7F60CD7F2EE8
F84277B1A77B4FF7DE81FE62D53285658D97E11C14C66F3EDE18A0DB0A6B40153AF55AE09429FFB35251AAC76DEA1B4A56C8EEE42903F991482376D4
4DD83F1B737EF80E30F007138384B97FD17ED1831C81C0B542C6BD3E370FB546A005E002CC3D88D8D32308F7F3A7E0CF61C81C8F90E45FF489CB2624
BC81AA4878ECE840D80B7C65C07DC09674AC26631AA18C061DD54AAF43BD825FB1138614532D8113C02ADAA06A65A928EE2A4732F525B53FD1E8A6A7
56F444BADE3BDEA42CB07EFEB4771626B2FFC136FD534B2F9DDE7ED8120A34395A1EED5EF2AF7B57AC7E63CFDE61BBDCAC36C7DDFEB1F2784FB96757
2B79F42FAF1BDB57ACABF95EA3DEA43894CD71591F37231B052EB665142C3E7D7BFDAD97CC99366FDEFC5973A64E9EFE856787A8170716208E1B6544
AF0A02B3D170DD62D7184FD02996E4D8B0B73F2E811D548E2F022F08816995053B320A2A56C0E8492D284274AEF403A55D9E8061470F559CFAF0AE9E
316152483E75EA191FFEFE4FAEFFD1F537FDCFF7AEF9E697AFFBE1F59F39EB6DB7EE6858100D4082EBAFE00460D4AFA731D49301380B92D2ABA4CB29
1A49A4461BE022833D14F95462B3E60EAC3272947E10600319E67CE7BC38072D11FBF79819B9AC46021794BCD7E76FD52BD02C43E4A27EA033C56EFA
9AEE199C0B51CA43B1807901CC4B407C07185C22B2F1197DD873FAC2F1D095C0890DAE1FD22805F9EBD30B6141E0C5524696057A90ABDD52933A0DED
969B746A5135242BFD324DED5DEB7AC370E8950736D59268ECB19535EA30D2E2826947EE98C8FE07DBF4AF78F3A1ED53CE5A19159FFEF5FD0BB73CF9
8FC7D6F795B76F1BAFFB66DF86E7966E18192C56ED309349CB1AD73645A46D773DF9D2405FFFDE9221043C7E13A0D354A89DC218766BB1FCE9CCB6E9
D3DBA74E997CD8D93F5E2792D41366BD6462B787F613DA1B8192CAF73DACD563A87066215C6C580133846997C0943E516ECD5316BDA80F19BD10BD76
1C4A7865629E0ECCBE40312158D137CAE1797014F7C3C640B9325E1919EADCF3FAD6D52B972FFEE9F4430EF9FCBA1D4363F53AFA8E50ACDD5E59B7E1
C107EF7962AD13469ECE6F460425126C0EB1B144DC2B8FBDC9581F488B80F93C583566A1029E47B1179FCE77FF8C26F61C8C1759D528C83700310F00
140F11724EC49B46C5DC692037E325A195986B1C321A08D400AE209288D51151C787803DE39008641EEFD86BF2AB607ACAD548132C4798308A374946
C1F6D74C1CA2DAA8D51A5A97F6558078F2373FFEE260BDEBC5F58D24D875DBF7FF5D082227752F9D36EB3F132EDF0759F827AFBFAB7DE6E71EAC1A8F
FFF0AEA5E3A6D73FE05292AB0877346C658E97482789DDD8345AAD913D41D65F4B4CCB2FECE9ADEEDBBA6C9F936A9906A60E8D7D4D7A648346D15030
EB283EF2B58B16BCFB0B7F5D3650A7D7B23C1F343E707852BF06461C14B091F2B0438FA0059482F482675A284D69D887BD8E4A336A30948FBD1D3DC9
02767AF9968D9EF9A613E148A02A5BD209C0DA7FA10062182305A5B9CA4F5BEC044EB15859F2FDABCE9B7AEE57BEF4CD6FDDF4CBBFFCFDF185CFF778
1B3E74D8FC49871C72C2E77637CD12E246B116A1F4B049CFA0231225816F53BD5134D9AB47952DF6D484EA0F84F4A82BA25CCCC4223EFD38DD2B961A
06AB40C74A84F9BC3E853801453A5486A039C0C61D8048711590D7FAA12350002830A8130CFC8157CE7B0654185C7344BEE0797FEEFCC1E740C84803
940DBC54B4AAF47EAE4CB413283F49FC302CBE3C98C652FBD6C0AA95437EA57B284EACC4DFBDA55CAAEFADBA495C7FE6DE578A3AF6BCC4FE50DB94AF
7813E17F707D05CF9D38F388FB29D56EDE67043AAEAE7BA137F087EB5093A5FA98E2677CD4776D606465A1926414AFB64315A9E7D4DD8A01DDFEA47F
939FD84B361AA88803F850C25A23A1461643FA40B9B0AAC0742D83E24E6081AFC7A240ACBA45E56C4D60489E818C02DF0ACFA22323F4EA50C48541A0
A027D9D160D3877EC354F96F277B7A8078433A9550040FFC06142F02D60687EB98A09AC0538AFE14305CA8F9EEB8EBA77FFAF35F7FF3CBEFDFF0835F
FDEAEA0F9C3F63CAD4C33FF4C5FB07C66CD7D54984AE97123BCC7538022832039EC329B7EE81C54039198A8454A6500DEE29DE9AA082C71C40877003
4B791E40819861A3182401332160DE19A28C4797C0BE7FF827C53D2A1C561D01521F28DE88B9044AE76862F6F9C44CC0B063E61869AA6BAA0E638172
3E23F27F6CD7F35968047D053AB0A849A9F9D04B1779A530504AD1ADA8F29E11AF95468651E82D36478A0DD1D70CE8CA2BBD056A477A2AA63F72C9E4
F6F3FB27C2FFE0FA2AFC70DECC639F074F5EFA3D3D2BB674AD7A62E3D0F060736CE76E2BAC16EC91CE3DC581AAE42C16A5DE9895BA48B9BB37575DE9
54BB7A6B9A4ADEC81FDC35D094AEF6ED9182A0E717C83BC86452CDAC94F2A81B8DF2C518CB56A528693598B0E8DBFD1187B76D90F284C59DB0E87CD1
8D1168DB20DEB447B1AC5C1F9917A70657C8C9683185225E22A859C0004D37041D1BF483123D3DB52086E3366DEA8D8348191EE537CA9088B7407BAE
E9EC3867CA94495326BD7FF1134FAE5DFAEA0B4BB76D5F78EFEFFF76CF8D1FBCFAA6C7D63E74FF43F7AFA8A6183E004F8C0945C0E29F78633AC8429B
F577B1360923219A550FF9590128042620255F6BAC6641B44399B6666B72150AE509A8012BDFD5B90D7A94E39201FD63C74FD8070442C134941D0DC3
3CB5534984AE08778A3EA0C7138B28572BE34588554D5878891710DC5F44B0138C23CBAE06DAA8C08F994EA56633686571A17FCCAED46D61C89E9D16
BD65E7D3BD54CE742EDC668BE8B629879CF17A3C111107D5E0EFA5B33B3A4EA4F02F2C5DF2EF3B7EBF729F9F26B5E1821F84EE9E1E9506AC021C3941
68F57736E84194459545AE11D89D0556C690A2A981B28B6D2BA234AFBCCA40DDC32E3A0D90BD53EA6DBDA6A1B9B4A6B097F424668CA1032E8D9E5EC7
827F4096468C9B8F4C870A85487B7430B826BA7CB00A43300A42C13E5EF44E801B608F8E9E9B9235054BE42BA73664D34FBA4D9F11359473D58809A7
208A03BBE40217E783412CB40F6A5092ECF9DB2FAFF9D695577FF59AABBF7AE5672FBFE2831F39FBCC53DF77EDF72E9B35E998CBE7BEF733677CBE8F
A9C8FFFBC5AF19E792A629B0B950D7841B8F14F59248E2B0D063E4683C4CF2BD5A6D98ED07EAD5400D7795203BEE33D2010B8C30DAAF4444C50A447B
D80121F74750A8ECE10ED4A829E620C6AC54C211CFF8E17C9B98A389228610078C5C8EEC064E3716F565FDD3C0AD578B0D13AC9F4435B6ADE991B2B1
67DDBA22C406D2CCA9E2A2CDD79FDCA153FF95A747E837972D9872F89A89D1DFC114FEE9F88D73DAE75CBEBDD5B2FE71C7F67EC3548618D8B9749D1F
AAF12DC3F0F9AD99DB5E5E57A5C7BB7B754FC9EE1D1201D5F2664F13BDA962240BBCBD43C6990A5DEF2918B6362C8B67E05858538CBB951A44428043
E5A71D132F54ED3E1E6E25724D5096C58975536694F290F0D9CB8A4E8F8815FA35C5AC842207BD53C348780DC6627D3C3173EAD6584F3F9513884919
877C068475414709307B9532D516C0C16357AEE107247508BC73AC85E9B9C2F1A922F0DCC0A743CADDFEE4DADE6F9E71CD0BFD752CEA62AB38D8FDEA
BD4F57E12890D3FBD98230D5400F535EA51A01B869BA8CB151601FC2DC0A9C4A1FC15B40A34137A0E901ADA75DE13BB6E33127195748D16E8BC0055B
304A0BDD71CEEAE13857A5C1BD0DCC3630E4E43D022A008DCE06011E878EA38114507056007A32D0E5AE92440122586239CEBC51DB1B29786EADD6B0
7A162FEE155ECFCE9E9A434757D20A1A55C884D43A6B3A4B6B2FAEA41334FDDB9449339F4F268AFF8324F1E3A98D367E747AFBE98F21C99BE3D5D076
C787B7EFD9B06EE37818F63FB67478F7A0F42CBB6BF1A2E1C82B97EA61A5D85F483249CFA461630F17B1961E1E5EECBBB09CB2074B70DBF31C1587CD
6AA35E03B025C4CEDDEDD9ED3140960A03A8E9C3730393B9FD9B6E2AA41593EEB99786022E2A59B87345989DA3D24507CD678556A80260298448A00A
9C4A05FC19E54E4C1054A29DA805006294D7D0DA66EF404FE28861360F13668492BEAF7CE14B1FC8E38C27859142F3DF58BD6AF9D2DFDFF0F9AF5FF3
87CF1C77CC490BAEBABBF4FA9F06DF2448BE891DE0355918D1EFF1F759785FD9742406FB1768313509F4D2C214ACFDA3A534862BEC4F8043CDF26390
9B63E9B1FEB0E766B902515EDDEB84D93DFFCB12E05A3FC4BE1EF3C380BAAC26331A584B2CAF0D5C87D540E18E98B422A762C651C94CE4785FD1A879
546BD8A660ABA456281B2A4E65A56C7BF4968575FDB1F6B2856D6DD37FA42636FF0749F8031FEA3C78F2BCF62BEA2D7FC04806F76E58BE675FAF9330
2B471B063D31A6BF7F4028AC427743F805CAB3A61742BF23F2A9E5764194A72A1FBA76619E91850C2C4F35EB94B0ABB558508909043CFDB4685A0182
18A47AC94E7E4E906FBA10DDCA192FD540A74755AC2B4DCD48178808C9D069FAB9A04D00AD2CCCDBF24E1C55724A2744246CF4CD90BFE5E343D0E103
781ED70770C1B1EA1EBA61386925110405ED808934CA971E657F0F7BC840B24BB90A7CE97AD8FDD98DF172D7960D9D8BCE9BDCF1F1FBFEF9F7BB6EBE
EFB1871EF9CFA30F3DF7CCE2279F78E8E1879EFAE59F7E75E763E38C014EC1D3C5323F34CA02B2E7D0FAA51220633C8F7604CC8480CDA5FF45A70E49
1E67C48CF07180DF57C032F101C29B8438FF196C498523715D504F82DA01A4CAE9A32817F3572D42DF854521A8D82ACCCF081C90411C28C774A80BF1
E166826395FE8EEA453B054CB1D2BDFCE5A1241A5AB6D9A53B69F6EF36ACBA95EE9C3579EA97AD89D9DFC113FEF1F8CD87CF9C76B5D4CF7FED9B8F36
E2B854AE48E6C426144A4C26495A709AA27E1178BD9087EC89638B1842F4B6638820B44B252AB259FE1BBC57FCA2A7B43B340A410D7EAF84C53790B5
C0718122382571E8F122A473A86C14DA43050BFA5E116ADDD81DB6F2293B102E5E60972950928061B33C038BD85D8C5E4006189B45D263703C7662E8
C7230FBE4429A3F4231616374418EFEFDBB1C7B706420A498935A3E742579F8A67AAA5552E631E3A36F079419AB4C0504C876EBA60FA21FBBF261D32
69D2310BDEF38E77BDFBF20FBFF75D1F386AF625BB18C683FA852D7742CF4645FEBF5574126BED41BB4C51B101E61E75ECB23A6067B81BA84E7896A7
F379406E77CCE53C2088511A0925E8B88D03CEF01A6C099EA2E0648C636CFE44CD876619AF1434FF2AF614AD963B64272C719646A890E8ACAB0DEC1E
1569A0B2FACA67D7F607D1BE85AB5C557EF9E145955096BCC43CB3ADFDD2D244F63F58FA7E7AFE767E75FE8CC3EE0EADA717BD31EAB4A81C1C293261
8E4A6E51F18CC0AB56C23474E9D9A6BC0377DB4073CF49A57FE852B0BA9ADA5077CCC232CEE2BD38AA6D495D2E350B8E4096866F6C9657F9E0E28347
1BC17227CAD21016746C5AE152DFAE585D4B61784EE584E5DB1EB218C42CA48A42DF76B40E80FA9760E2019F07DD1FC5EFC1A61D1A49340C303547FF
C01C7DE078318F439C698CEDB30033C630F22B0152B2D01AB53F6565B01063541168A613EA569C04DAE4811FC2C93BD595EE4D9B366C5EB37AD9ABAF
AC5AD33D6256CCB1BE7F9F71E459DF5EB82FC7F240B624718C62B5EE3BCA575EE033E6DECF180798A37CA859F795D442C5B26183DE9F02BF08223F75
24EC8D12E63D14437614D7460213069611A61F6D7881742CC566C3091F1A743B25F47FB0D3C019064B136C0D93963BAA23B89C69F669749535D453A8
6BCF8E3255EE190C5AFEC69776B8FAB5BBFFB3CB481B6B75945C33B9FDA29E09DCDF4112FF546ABAAF9E376BFA47BB52D9D3A46C951437ADDD3446D9
766CF386353BCB1A79BC5A74029165F1D0832FABB485A6C0C1122A36BBFA1D4F43361FD9564B186B24E57190CAC096A17C0ADD6E17A89EA0A15A60AB
81BD16B03736E4F35B098A82CC7BA3016EBB5F3579AF15C4AAEE60FE45BFEB1B7698C01D3CA5CA5C47523A6336F7BF9895870316255D0A0DEAE881D0
0959CE170A18E8A6F19D2886AA56086E30F83CAC4C1E64494C9D09F3F41C09263FB8B6802A05B2095C32EFEBA99B810A010392E01546DD395D33DFB2
3437FBC2B99465379D78F4DC9973BEB166F913F7DD7FCF7FEEB9E7E7DFFDCA972F7FF799177CFCA61FDFF2BD6F7FEE1BD7DC7ED72D37FEFAD63FDFF0
A39B7FFEC52BAFFEFAB5773EBB714FD585DA5F1473F905DDA156683A7447A992A1C22681E8173EBF666F0FA60D53D5444DBE6609CFC8A62228A859CC
FE8D7835C06823EA8A802252F90C00C76B2BD5350B9C20EA37B248D7CB9E6E8CD72365D6A9CF18E9EA97A9DAB772242A2F7962AF4AC25D7FE84BB335
53DB4F583981FB3B48BEE839AEDD7EFCF4B977C6ADB47B87ADAB43CB576CAD37C3A8BCFED5B59D3B46CCC6E0AE1D7DE389BBA34F67B1B56508167251
5CDE81C15928FB0629F664284486ED9EAB528A88428F42056A15EA0650B29A1A538A9BD8667DFA40320B4DB3064E28134CDA92B8B9AE98F1A8C1A6DA
D8AAB881590B02E16BB00840B00BA09C373E1E4B388341F142DA028A7CC98891B369C45815753DB0B371DEFB1A302043285998F5091706238023C3DF
2F706CA0125C7A1D4B30FAD01C7220A0876620E25101D43EA90040EF1EB72C3AE51A36FDB1CF2E1E38A3287FFB10418C5A8B2E7DDF0D5FB960C1FB2F
FFD8873FFFD5DFFDE3EE9B7E70F33FFFF5C4BF1F5FF2DABA9D9BD62E59BE79B0581CB6ACF1D55B77BFFA97DFDD7DE76F3E3E75D2F4132EFCC5089C79
A928822F012C907D41FD071D477EAE8D48A7A35BA91926438423A8F551AB45573DE40089486759166ADE0952D5C01F54D7EA66DD87463A5DA302C737
ACF7370CC34A94A0AA4D955C29AACD4AC3605DC138A90FED1B185351A5C78AE3DA9E5A229D45776D89B3E2696D33EE0DD289F03F48D27FDCF7D54367
9CB0266B3597FCF6C117573DBFBC2BC85A6EA55CAB05F0E6ADAFBCFFA94E274A1A5D4EC680FA44F6F4F8F1F016805C14D375E9E155367BFDE900DB26
C7A6BE33086CA3D2A0AC29297F637A47F57F905ADDC32E450EC29009AD5E8097640C6C12B5122076F1944B47790298D9FFBF08A59FDBD913F3460BD3
3D6151F10CA650BE6D84E220F68F212FBD796A0E013210EBA039CE323E6012A0AB801F96ED614AC71C19CAC25463DB451F33429EC8E5583BB7EAA263
863B796614B5360CC52A7F746269C7F44306E8C8900A82D8CE82DA58B96C19FF7BC5AC6416E688FF7CFE974B4F50ADA37A9F79DBA443264DBF19530E
6CE8B34CD1AB04D0480B1B6F1443989FCA00082BBA89662DD436D4445890946E93EA77811966A5807CC5CF9821380DB9ECB2C0C58067C36E7D6CD14B
A336D50DF4A163BFD85582033AE32FE000E88C8F96E9E7E99EC40E85BB947B1F5C4F2FFDC723264FBDDE9FC8FE074BF7AFD75E34BBFDEDDD69B0F799
3563AA52B6A00DED56AA35E056851728C716F4F0155677E33966A30C4C9BA4410F8E0AB151A7AE33F1CA4E04A16D74AE823213C07741C4356E656F4D
8924037B2D6BC95103F0F45004416350C0AE2E6DD18F53D4C03C1815B807556C25B0854B64B373C5E3FFFEDB03CFF76A8AB2A85C66ED4D5F4BC0ECA9
468E3806F08202643A00E439F619FF1663D2CD541A047AEC98124B82109B7BEA0720CD9925186DDA2635DD786770633173C358DD83C95FC47C7CDBE4
A622B7358C18111CE114C21B7BC38AFE2DCDA230CBF90AF4B9AB0ED040B26EE7BB84FAD67D52619E11825434FCE20D274C699B7BE6F79EDC2B85EB4A
0F3A8454DD04ACD041B16BF954BE5B3EAA0CF64BA77B2DEA386CA0A5E207FB970191887962007D330F1A21B866B6FD61BE55430AAF49B78A0A29BC33
F8946A6CC7A0AFBD40570CBC639478A3FDA04E3B1529E05322776DDD3C1C07E9EF27B74DFDBC3931FB3B58C2DFFBD731B33BAE2CC67EFF88C862BB19
0973DB4B4BBB5CA8F886FB91A640BD37EB10DC51968B051735A95A094FCA50FAD44B53164420511EA53F0F53F4DD09A519897A592ACF128E2734359A
2A659C8DEB337B6EDFBA7A1CB4322FA40AD8A393A5250C3BF225D5BF680FAA8B963FF7A3CBCE98D33E6DDAB4F6F68E77BD8667564B6001F99D206BCF
1C796C1A42660DB30E80028A9DF7E48CD3A38EC0F3D13DC7125344F60103C008F2C2141A81AB230F753D736AC09CF358AD4F3B2ACA5F002B05760FC8
11F809DD33E9E6733C7A5B00205BF821BA263EFD64DCDC38DEA223C21C7022462AF63EBD153C812C885AC16BDF3F7ECEF937FDF989E565C6EF80F1F4
FF045ACEF98F43CF876AB2820408B5003E1629B83EFA63A124C39C32201398EFCC5A228C64E65921CB05C7F5AD46C4CC49F0A8E820B3EB96E5556A61
A282B83912E234F487876DF418D2872B69AC5575B71B876EB2AEBDA3ED8AEA44F81F24E1DF6AFE70EEAC39BF1570E349FD466F7FD1E97A654D774F3D
490A1B37F69458EF5E8B8631540D5A81F20B0335C477A429F67D4F947B474C1186D41798751D513DECF55730B8D749586BC618A861362F84254D2F54
32A5EE3553DE401DEDB36379AA6EA5DA90316BDE43EF338C58C3B6BEFCDF3FBAE898B30F9F3EE3F0F947CC9D397DD6895F7AB81EC3394C8295CF7D3E
F7C40923DC43780647AC892965C38484660E7C8FF7C365183F907B62E7DBC2440BC8730150A751FD4B1C73B0E89485512895E674BABC0FC09915E4D4
59A45EBA56614197C7298F9A014372F6EF1714B37C1865000B3219F0725EE6A03CAA08B2EE8F76CC9C7BF22D5D3DBDBBF66C7F7DD9CB8B96BCB87CE9
C65796EEEB5EBE78DBF01B6606DDAF5068708C054CD253EDD57B0CE801059AA10098EC2B161CC5D1A1B8250843DF67E230CB0E3002CB19B0007F40B5
8F69871ADCD86BCB084AE1B1BF6994CE8EB077E53E3B898C5DEB7ABD206DFEE96BEB5AFEAE9D35150F9C34B5FDE20997EF8325FAB3D10FCE997DFEDA
28B476AD5EB46AEDB86C25C5BDC33AF30D69EF7A6575EF7061B83AD6F7C6D66D6323C32355D4F00190F6D8E1839AA7CD8241C70375FF56CF48426580
F22AB6D63C1F930163E878251708BFE1792173F133E5D47CED3BAEAD84ACD40DDF8F858EB2A895779CF1E8E3DFBBEC8293DEF1D94BCF3CE3A813E6CE
6A9B7DC287BEFDD85EFC09A2CBAB7B9AB5EE73DD7BE08E85CFFA586808522A9C15B27DACD4FEFC8FA939D2364A6BC1F16E357CEA090495C423324BF6
97CC302962B62098735C60B36F38F7E00CE661DFF010D83E609785EB05C27372DA7DCADD7DE80AAACA255849114BF0C6A87FB4EBCB8051BC61988EFF
FE7DE79E72F29187CF3FEEB8B71C7BDCBCA3671D71ECD90B4E9C39FFB04367749CF181B72F69B17B07F882395C20009FD9944940150F9F348972F995
33E61043AE5852BB809968AD813E8211C9EC39ACA1BB843D0A3C4C20258EAD67626EDDD6DD43ED5B5ADDDC13ABFA8ABBEE58EB25E1E66F5EFAA05FF8
F57FFFA42FD9B5A0ADEDF2A16CA2F93F38C23F5E76CEF4D95F35E29147FEF4C8C62DE36E96E95249A7C2B219B392267EA36E378DBA8BAD5EB5426522
F5FF1E822E111E3B685105EC24D2B4B5AB2159217C66CB533FE917C647B62C7BF0CE1FDDF4A33FDEF1D3EBBEFBDD1F7FEBD6FB5EADB4FC2049A4DBF4
0C613B964BB5412C285683FC8A8CD57FFAF8A9477EF0277FBFE5971F3D6FCE21538E3CF623B7ACED017488A22EA612435B655B6ACCB69823E0377D3A
125C6ADDB5E7B126878C5B290885947DC3FD4C1DF6BFC1B84C385C240B3F60D460A60A51C6A91EBB442E0028D326A011329E90F17771AEC8072621D4
7A21DF8FF960007422AB68C26B438AAC05BE3E2316801580421890C591F67DF077B1E7A0F7A443A23232B267FBAAD75E7CF5C5654B56AE58B77DF77D
171F36ADE3CC8FFDECD5BED1B4953647DE58F9C47FFEFAF747973EF5E00A8CE9B0B810782D89757E6C0E5A285752B07AA0AF1E2B2F88A8630A4C27C4
FA30624631039A63209AE818A6CB105C39C8D0D9BDA6C7C111AC0D379603DBB70CCB346DFCE7778BCDB4FF4F7F7CA9102F3FA6ADEDFD2313E17F7044
7FCBBA6DCEB4F97F0D92EAC247F6C2DFC2311B5AD1B3E58B5C2B3B0D85C4863B27C137064C6CC183A415B7E22C741C78D868AAB963ED5BAC5901015A
CFEB79EEDACF7EF3A6AB2E3AF5B4238F3EFCD48BBF78FDA72FFBCA75BFBAE9D69FFCECB357FE75D91F7E789F03855C0DBA0A656D15E4FE3DF57DAFDE
77D582D9B34FBAFA8FCBFE7AF1D469874C39FCA6C53BFB244F0C90EDA87485D33D683D140CD844649900FD1FD951519E55885E9E2F80F216F916F03F
21C3DE986118F174EF7F75B63396184A11C394F76152C2193F4CB462F70CD6EA0B72EDAE0432BD51947BF9309A70FF945DD0ADA24353604F02C09156
2EECC3520600820701E55ECB0E791B8F712478167CB2B25F88ED3D70FCCC53E67DE5C1DB36EB7078E9DD3FFED8FB17BCA5A37DEAD4A9ED53A674FC07
56C5F48A3C6891CC8DD6AED21E407C29EBA0C6B97E80666D50201FD92B88733D40CF9A8F2DC96848855301B7A7158656DDF03CCB426193EC7CBC2B4E
F5C8884F2FB8E788F65917764F2CFE0E8EAFB8FB93ED1DA7AC6B6591A7243D2995113B85BB9C3152AAECEDB3A98655B5FEA26337ACC248A1DC5F7128
372258C3944755E148255063639EB4BCFCB4080A2BFE73CB6FBEFEB1A3A64C9A74C8FCF3BE71E35DFF7866E12B8349B67A1DC539667A9EEAF9E43B3E
74F92B03BDDBF776F6ECEDEE2E54449844852DCBFFF2CD0B4E9839637AC7311FFBC6F51FBFECAC0F7EE7E1875EDD1AF0F44E6117E005ECE0C9D27798
8C3132394B3586EFC8DF123C79CED7D00E623740CF62221CB84374B1D2CF4174099B63A01BC840D8675FE250ABFA9089E13EB70B6128D9928B397661
2E2518E588227A052141EF63FD6C8A2C0D4C1D8C7119571B7A4D3FF7DC665B103A34FC5CDD90DA11909840E767C5518A6661D73B7FF1B3DF1E3369D2
A419A7CCEF38FDB28F9E7DC6395FF8C1AF7E76F7BF5FFCFBBBDA264F9E7ECCA2161C42E91357BB1A5249DB03BF023B007ACBE58F6E1A42A3B19FF7AB
A9154A7C27D26960FA713E0BCC12416D8A43EF27EB4547A03980AA83E728136AA3005B24EE2BAF5289E174B9AD50C4AD1BDAA7BF65553C91FD0F82A9
7F2B73EE397EF2B48F3BAD60BC5E6D0CEFDDB9AF4489A3B46BD3A0E9066EB152EA5EF5EAEB2527A0B8716AF542B16487B229A11B1D05CD511333FBBA
457FD4F09356545FB762CD839F3B6D323DCC87CC3CFCED57DDFAC0867D3AA1B0886C7A64FBCCCCD156D3D53A0BEEFDC9BA277EFAEDCFBFEFD273CF39
F198A3CFBEE8CAEFDE7EED02447EFBCCD9C79CF4DEDBDE58BFF1F58171057C4D002E0EA57A47B81A7A1F82753DA5328ACD20140A60E25021D462AFCA
0E1B31CEA68C71FA6C8FCB0EC079DAD6AE8329370B6CD051270276CC4D01EB8B03C7F2ED8603B63ED7082957F5CC604A2367BC015941AE240072A43F
63C41485B7869A40C27205798D92849E1D413D393F33B4B49A36EFF7A1DEA9598744A12F0154D2DFBB716870E7DE17BEF9F1AB3EF7F7BFFEED5FAF2C
DDB47677BFABE880D9FCDBB74C9EFFB57F6E198BC191001AD81D71E82E384D27806E0870BCADDE3FFEE4C3DFD999EC670486828E17281247A16D242D
FE1E5DA04DB1EEBA812EACDF650491E71A810480A895897D9D43745459DB1F7C833E9CB9C349A0CFF2C0ECE9473E134C44C7811FFD5916EEF9DAD4A9
47FD530EDF74D927178D6E59BCA2AFD2D9B961F3EAD77B9CACA5D8A5229261D4AC0D8F1A209A5146F5966EF0A534CC20B460E99360FA4671D2DCB5F4
3BE71F3AB763DA94C9B3DF75D52D0F6CEF2E82864605B3EF0585023DB2D566AC4DD3B328F741832FB1EA85CEAECD1B5E7CFC470BA64C3AE49043DA0E
3BEEA8D3BEBD70F9AE3E2CC8D94D33DEAF69C35F683FD0316896020C15ACAE1A651E834B175E19CE8811307D10633A7A8B28F64C2757C00E59510B10
39F6FD544C59182D2BB089E8D2946128D5ACFB14174A828BA025C685C8E1CC3E0A3D03381C3613C77C20CCF97DD8F4879ECB4E2052A1D24F72EE5292
86591CF83ED002FBB5FAD019018FC863516086C0924E75BDE2807D0C9E4292BE9972533F297C6AC6B4237FE7B20743C6E82C1D0796928EF06CF43E61
9AECFF6BBCA7E3CCE7B06B08624925BE0A6484D6C8A5432C4B2C2A7E9467BA3E153EBE3BDC55A5AAA9D95474764187C9DBD1597323A37FCBA621BAEE
FE65955866C9D2CBA7B71DBA289C088F037EE64F3573E5DF274D9BFEE5BD51DF93CFEE280651B367E7E297570D942C2FACEDDDB6B7CF05C8274B7D7B
70CFE66DA336C7B91E5C354E2DA4E731932F4B5A9465ADED2BFEF5C563664C6B6F9F71ECC77EB3BA979293AD0405A8A7FC9A83ED78296925D5B2B06D
19C1DA4341BEB6057FE0AD4FDC76D3A78E9834E3D4AB3E72D8A1C7CFBF748BA839E5B1C125BFF9E75EA36FF3FAB1C2EEE5AFBDF2D2C245CF2C7961D9
D2175FEC1AAD97A800108295B64376E38EA16AE342E5C3553CD1A3075FE3D5BD3027D1B20E81634151877BF624B14D74119EAB91F725FD1CAA0AE94B
0F1460D6D462FF9D98E1B7B1F69D807D7D5266DF826B07CB63D40DD49748DE0F26B5F128CDD58A80CFA73F04C11F0B409CA148FFAA68278937ECE78A
C0401E7115028F31E07D7D2505A3FA7520022F1AFB5447FBBB9F5CBC78D99637469DD8B66C6A371AA3D4C1AB861B46C2A9ABE6504FCF8ED757DEF7DE
B927FC063A6A6C118C755FAC835C51D8AF0D1A742E50F407D658DD17D02AC726360C1A8DA6AA77BD3E9E665163B421E8839636BCD2A0B6A1F0B169ED
B38E78229ED8FC1DE8B99F9E6767CD95F3A69FFA6C1C375D1044B2A8B27BF7CE711137BBC76B1B1E5F31606A287638769DE2C4879236A591B181A2A9
2862222DF63B6E3737FEF0EC29871C32694ADBACF3AF5DB42B0049C6B7463ABB1C4AF481A02A1DD81C5E54D3EF600F862D029D2AA5D54F7CE7BDB326
4FEA38E2138F6CEA9563DF3EEEF093DF79C199A7BEFBB29FBEF0E9B6F613CE3B72FAD4134F3FAA6DC6DC938E3BEDA4B79E7CEA82F77FE25DE7BEE3D4
055FD82BE900005B47410B008541687900B15195ECFA6C771BA43ECB0D0957313638A89B016F2CF3D3807E48F1388FFD3B5965109C19461589465551
1F8FF2C32983698B3D47C2EA2214EF4C310C234F411321D2948D317F07959E9D3E723520CAF774F601E9982BF6E1B742634404F56E974A013A0A98A3
085D73C803C3AD836E96E4DF54EC81985597FCE6D2938E3C72CEB1734E38F743175EF09E1FFCE2C68FFED7777FF683EBAEFAC16DD75E70FA290BCE7C
DBF9175F7CF9155FFAC9D261ECF87111A9AA580C0050F44FDFAF565CF0FDE05F50EEA7134F0AC3F762FA97D79EEE147290BA7DED366D95AAF286654B
76C9C8967D17B41D7EE8BCC727383F077CFC53EEEDFFD969B3A77F6C2C8B363C5F0FEDB135CFBD5642315DDEFCF8CA8A6F18690BFA90B6836E5F42CD
0BAA72F6404549D307755EA32C6DBC7CC30573A7CD9D35A3EDC80BBFFD2C944162AAA6B51FB8BD3D7EA8A5564D195245DDA2DE14BB2738E0A8D1BEDD
6BEFF8D607CE3DF4C819275EFCA385FB0A2EE029DADAB17CDD86C79EDF582C87D9F83D9F3EE5B40BAFBEE9C9675F7E6EE996AEAEC1D1F2D8E848DF9F
8E9A34E990E99FE9A358413192600306150110E75DEA08FC5094C65D6A82811172D34CC79E2D94E32AEA9A35CF00580D0CA706A47728190A9E1826B9
0A79201183A1D4692BDFF807224C8149A6E48BEA3E8E038AEAA051F1916D294A623A21A1DE21D17FA73849A81CA06B825629548D227E95D06F34A0DE
4D9F916EA20E626D85A1098E2E2601F03D08D92B0C9F48DA9A8E471D42175D3407FBDED8F4DC1DD7FEEA5757BFF5F0B9279D71D1F9179CF4D6E3E7CE
3AFDBFBEF4C56FDEB96DB856B565CA3A89CC2516CA6D80C1CC1C405401D86A68E606A12588DDF12618814ECFFA12FEA2A89569D4CDC81C58F9F2A61A
3B1E26CF1E75DCEC335F9A18FD1D04ADBF78F53DB3671DFA0BA7D57CEEDE753BB7FCFB6F0F753658AD823A7BAA05C0928F83A6E94817DB39CA4F8E0F
78BD1F0B61BB98BE47DEEEC7AE3963467BC7DC995326CD3FEBF723FC4C513E4D90F094E7404D2770F635A0A311E6DDAA347B97FFEEAA734E3F7CEED4
99175CF1A95B9E1B30317C539EEF42FFF64D924F4ABD45E67A223FA912D6C9DAFAF7C7FEF5CE29ED57DDBB709CE7F1544328BD5F969635AF2998B116
0B242BE16BAED4A1F71303D82F282EA8EDB06A36DC87917901ECCBED7960C4ABEAD5BEBE7D659C09602260308AB8614D4D447EC6D19FD08F518083B6
043E9076EABE62C60D63EE216104C17D08F4B83A02AA172B3AC6DFA9FDB043B007582CD8AD05C02E823B1445B98111CF28A417676CD91D7B0A2709D7
1C692BDEFBD0037B1AD2716AF5EEED7B2A6009602E90B622DEF3413641940CD088B87FE1290578CB0AD821D6098AA551AC00B351B7E9F55A51E658C5
5A5D04C6D6F53D54421945CF8B779F347BC6859D137BFF8360EEEFDE7DEC9C8E0B9E4E526F707BD5DFB3619B1166B15DAA0E5BD0A94B9AC395C2D0C6
471F1EA342DFE9DB396699E37D25C84950454FA700D87199FECEB4B669EDD3DBA61EFBA5077AC61D3BCE6931D0A3149472A1681B044E4F1D3DAD36CB
1B1EFAC935575C74C4544CFA8EBAE88AFB4BC2E7ED35D50710DDA0EE57040A68FC3468C52978A7381A34D5E5BEA2F0BDEE9049F475F81DAB5E5AFCD8
D3F7FCEE47BFBFF577D77CFC8A3F0CACFEDB731E74482256B9C9F1FE54AC03E5CABC004CF8E958829AA7E735EA5EC044C3963D6AD995A18DAF3EFECF
DB7EF2B52FBCE7EC63179C74E261E77CFF75414755823BC06881002B3ADE10A4BCC5578E073E0E5ED71DB5ADFE211FB21D5E9331C64914E4BE3E4992
9F3C11AB9145B9DE0710468E09AA80300C3764D7B0BC5FA07FD48D28920D33E0DE814E9BC88853380C8A4002D698E4E75C0617C560BFA24F8C4A82E2
3C79530E5C1A0ABA8B4002BFE9100005056D545D8C6747461B0DFCD5EC581F6151E91855DF31DC60EF7D5BA1AD5C1FF564BCEFE8E96DEF1E9ED8FB1F
F0B9BF9595BE337FD6FC6F8DA751608D9A95C18A57B34BC3AB1F5AD85DB2EAC65857490491EF54C7476DC756E33BF63475204C5746AA2E232D3D8FF5
285A9D5F3A7FC129675EF4B3413CF5F4ACEAAC36E2857EB3E91B85BAF6257D47C1998AB25F5879F4A2F933E7CC3CFCD8F33E70DDAD2F74356B20AE53
6C5358E94201DB70E8E0843C54876657904A6A46A88E0FE11B427F3CF8AF3BBEBEE098B9B3674E6D9F3277FA611D87CF3FF5C463E69CF3B5B34EF8B1
AD6DB3522DD62BE37BF60CF66CDDBCBD7FC468365D8AE31032E580ED534A46278C5E415677EFD8BAE85BEFFDC085A79D31B7636AC7F4A38E5E70F1A7
3FF9D90F9FD63669D2F4174298FEB5124750059004BE2303CFA603814A21561687DD06768CB13D24215412BA7E648F29AE1594CF8E1DE97E7FCFC8F3
586904461F4C1C8CA525316D14ECF00DD122A8F840C97FCCA23B65186C04848AC3EE0325927E57392823E88278C31F28E5BAAEF0613E10D0DF029D97
00441A054386CC4D0AB9F2616B20742254EBD8E5D1821F45FECE65E30A0EAB23FFE8CFE296DC55F4051D614630B44950CFD1374467C48A396DED9F1D
9F88FE033FFEA36D97B5CDFAE4FAA49528A57A7A0CC793BD3BD72CBC7363401DEDC647EF7A55B6B0DB02EBAC6E0BCBC47C8982D18B9AAB7B95DFAC79
0238782AE985DBF401E777A8EC45411ABA1465541E6BC7853F06FB6E21338749165406F676EF191AAE5A8CB29541B318A42A895AA93FD690545230B0
2DAE94A902F0EC100EBF92AECE329A3EC47DA113141823FB766C5FBF71C3F62D6F6CD8DD57E81B2BD47F3773DE59EFBBF86DE79FFED693169C7ACE71
C7BFFD6DC79C72E6E5DFFCCD8D5FFE65B585D3247F9E75D7F6A1DE6D9B5E5DF4C01D579F3B7B5647C7B47927BDF3D2AFFFE00FFF59B76D78C4485A2B
3E7BC1D193274DBEB03382924F2B9310218850EA27E61E2A3094E986ECE56BD48016A4B385A53919E313BA8AD13EDA80DA1677F80C0DB02DEE2E02CB
63D91E9C1A09F72519BB1E2491AC56352A038A53C95A5DA1E332D65858981ED08FDBDD820E0EFA0C80F9D339202D578358456FD668426E097B09BFE1
48D888B926B51D9156C690C15E610A2547035B88281CDA45F75AFBC1685798B6EACFEE49A99818ADD15BB97BBBBB7A46E950BAB16DEAB4ABEA13B5FF
011FFD49E1CFC74E3BF43E7AC007770EECE8ABD3C3299C86E7F9C38D286936EBE57ADC8AB02CF385675B1AB10E8CAED1530C4481AAC6A6A5BD902BCF
FDAF19C430D773FC1634C3901D519DFACC4687A236253011ECC793C470DEF285E5699ECF51EA83F47EA43D57BA00DDA7B611B35E263DB4BE6C510C86
58CD01ABE3F942E7F87CF6D06617F056CB7DFCABEF39EDF8234F5A70F107BFFE3FB73FB272D78E8DA3157BF42FBF5EDEB5FF498EC5D6A52FFCFCD2F3
CE3B7DFE0CB410538F38F3CA6FDFF5EC608D9DF402C07493F17BCF9EDE7ECCA99FB8D766C802EFD64360FC294BC762D44F63843B068E49693897EF63
4640CABA205A791AE63D81EB43EA1FB41F5EFCA15C807C9FE3E7D801F44748E6215C8A9937A040730CB128300D3087C35A19EE0840252A9BCAF9C42F
A1E0C73820900C1FA01395255723051DA358C95C132C0439C0337DBADF5EB96788677E01A47F58F30C8EA9365551B186B45064AFE98500E3D8F26218
340676952C25A2A1B74F6E6FBBBA3111FE077CF83B8BDED63E75CE63ADD6F83D373D323C541FEB5EB7EAC9E5BDD4D7D636ED705CB7CE8EBA69229AB6
4F75B9EF28B4A0DA2F0CDAB0AC86F45C14D8354A2B9ABA74CAF7F4E04B97225930AE1EFBB2981D313CD3054E87B9B600BB61ECD6B727C8775D0E2CAE
31444F81E381224F185A4375A6CCE6689956A438C4131F937185E980802EBFEB8371E7C20D884A5C4C0B8555AAD6D0A7A304914B87E4CA8F9CF8E597
9E7972E19A276EBBEEBA0F9E7FFC0973A6B71F3AF7B893DEF3C99FFC7DF9D661836F45421DB490AE0A8C7F7EF0F863DF7EED9F378F61171827AEC900
FF3CC383542C20C3CDD99E3EA0F2D21CF7C7F20051129810310A731DB390C785E8F6B1F5C7E990609A907300D015D05159AD85FBF944FC59A1F61BC3
C09BD54A700A72BD1447AE4299103056A2396EF1FC31CE7B7AE9096E1C125DAEE1AF2461A6130E98546BDFF1988D081852E2D7216F4A57A1E8BBF477
54291915A7E98354D818DEB7735DD1854971108D9E3575EA8CEB8D89B5FF011EFEADB8FBAA691DD3173C2BBA7EFB891FF4B8B231BAF4EEBB57771606
06FA4CED562435F12DDB71A264F4D575855A9165B00555EB40C5A6F85FCA28548C57BB9B5AD0339B50568EA15A47557324598F5A8143831340D1F951
AE02ABEB824F4BCF68985061AB9CB1A2040A7F74F98862549CA07E56D03BA9BE5D2E55F909DE8DF269D0B0D34CC56ADC8878B78F1AC1AD56D9055804
AC641BC0639B6A7C407381088854ACF794BD3F1F31F5904993DBA675CC98D13E63C679EF7FCB291FBB65F9FAEE319F1932A8A6B94E41E1ADAD87DF3A
E5FC5FAD1E4A99E784D29AFD7A6C3F14D47968EA652C2F8D318B44A59F4BF267ECC42D4361AA48F416A34485D484B3D116FB10C44CCD41EC637981DA
9FFE4547C290D42BEC19E6DD0A70CAAC53CA86E5385B428E7B2A0EE01540B73384F711948300EBC9677A80336BDE74305731F24B36CF18E3480A6094
F84CC2011AC6CDBD3E9568D53EC10E8D4168D7751C8C770F554C6C575AEED8701035C6554A27864A1B974E997AD89D5E6B62F67760477FCBBDFF9469
1DE73F3DA2BB1EFCD56F1FED5EF1C8137B065D68D918554349A3AE9B1BF7980E1C3D87FA9A6EB3E1348C5281422E7E53612BE00E177A7794ED513A83
F4023421AB6823B36965DBAC9143CF76628C52587A86CF3A1AECD4190973A4E0292CEB2467C1377DAC58C33BC83535B9AF4E28D9B319A0A0DFA53204
8877D928DAF4EC47857E2CCDB21846DB143B8A0A5BAD6B22433C5128743D7BFBED3FBEF9A77F7F6DF9C2D5EBABC5DD6566FE871026B0078B411A64F4
76AA32EE8A57DF33E9AD2F18D8B6435ACB03BD8EC2CE97D4512B97AA6A0A1AA80123C7E6DA008E47199D593D680E28B8A3DCA9CF69EE671920B2632E
1FE807234FA6F00C67EC60506BD04DF054EE50862242F82AC17093C1466C15C48A85BEC3606744B642051043F324A46FE076803A00E7518809C0322D
4CA46D2AF6228DF6BB026172B36F990D09F682690AABDCDD6B38741048EC3FE84F5B69792C600F72BA485F6795F74D9D7AE4E3D144F41FE04BBF74E4
5B6DD3DA2EEF0FCD9E71DFD7FD7B7776F6289508BB2961C0551EAE46BA6FCB485D78DC99DA9E6B54863AF78DBA0A7237DAF3A94E07DC35665AAC1610
D0838A260666D8DB315F3682C736983B3C2E47F30E689B08E8DB3145B765539EF69D862313AE8B13E903B51BD82ADAAF6F1304B697E74CD6F48C2122
18344BE86D23B6B4D72A6E0E336626C94455F268CDF7D3D4F553CAC4C0F152219CFDBFE75E964ADFF554E007A23E5CF40007127EEFB32B8DDE774F9A
747D90BCA9D519462C78908550C8A7FF074D312ACA25187CCC1D8ACD6DC3100E4C342E1B2504AE32F71482EF26607888717A2F19668952A68F2E012C
6448A62816EF119021A28B16868B5A89EE299F82395101325EAEC58043AA2B14FB92E5AA1F401F074E038E288CF505F220E565BFB4839435C3127AAD
285700731B906B8C42EA154C63A8BF00E9107BDC8EC2BA99A6D5061D4DD250981304E9EEB3A64C3DF6C509D4CF019EFC33F7D5F3A74D7BF75352D6FB
4A63EE68B3E1058932B6EF31B1EF5655C7A5DC17632657DA3106C5ED7C8287A95B6D5F3DF4AA65886106A131546E4A4F5BBBC6284A3415C949E21975
C11CF438F7E1C6CE2A4FE3EC4013F1F43E64BDAE7A839E709C07A828D254D8BE025357D7A1CA11C6D5317A668D1A587880F5A5E88E35C59004AC97AD
C258DD133306EEB263C59A5EF8FFAD8CB76CC0D4228F4BA87742D5837A5DEA1812D19440FCE110A133C417C22995ACF1AF4E9A34EDC3FF7AF1951716
3FFBD4C38F3CFAECE20D6FECEBAB56AC660E4682CF30DE5141E90B4E62DA010B0A5D8E8C8512901506AF38E1ED1BEF307010B86E945AFD544785BC98
43CB0E59A018E023CD7626100A735CD45498CF513D0FB6011F1E5CE283CECC3B03698B3CA3C78C41A2D280A72991865D78CA68E3BC7B80AF789AD80D
C90ACAF8BEEF721D45FFDEA253334B32E9372CD735036A74E85ED53D85C1619CBD7274DB94E35F9908FF03FC2BECBCFAC86973EE0C226FD1D3AF15AA
A5C1F1CEC1A65B19AED183AE8DDDFDCCC10FA8BB16815BB2026152BEA70E93DD32DC712B8290948B9C52D8BAAFEC347D6A7E11853A67C13876D90872
976AC6A28680CAC68C8D83AA1E76F83C92C25A9ABDAE630C9E426D3663E06354901BDA44CE909D33EE63D41309635B8388A208365F980244CAA30A3D
66F71C54B0097BE1B2FF650A5D5CFA09092720499D061034516C0BCCCC6305DA9170A8A2A64B860717BDA650435F9E3269D2F4234E7BCB8253CE38ED
8CE38F3FFCA879471C75CE7917BFE38A0F5DF34A73C0CD890B12FA7D50F488991C986B0750E0057C7D98412429B527B51D4E4077294E7DDC545FA658
0BB0EA4FD2F20D60942C27E08E865543030E4E3415BE4FE91DACE050B868AE5484C5239A09E563A7C86355F63CE6BC0F5274DDE5B9BFC42E81E782A9
6727319DC33E8C14E91E47565D267407E81C8432A1944044090F062C495AEF1C6DC6B128FA2A5B3E67DAB4052BA289D1DF819DFC8DBB4F98D671E2E6
34EAFCE32F368FEE7E66930A8AA667554454DCF8C68A1DE314705ECD52DA182C47F45CA794B2FAD7F6DAC8D42A663D9D2C4DBD8697EB5DC1799B7DE9
816E89E08A973A36DC2664EE544B95BBF674CA23710A5BDF74432ACBA35C3103193CCC3D6D2924B316BA508482CE197DF9ACD1290A8A3DDBD2FB35EE
73263E4F07825CB7176B8224E49860B7AD0CD5301F32B9EA275D5A501FB79D71931DB429BB26E1D8B0D0F0D161DCAC908E17C8D5375E72F105EF3AF7
C2B75D7CE17BDEF98E73DE7FE57F9F336FCE9CE91DD3665CF9DCED9BA0E743E709B534BE4293EDC26638959640194FE74AC0C26201DC0AA2A8BACA80
9AAF629D4F1C145EE87304D2674285805E9F8E078C3AB85F8776392B0AB2743A9438222AE5A9FE696284CAA8C1302FF4D398EB7FD632CF8F0021993C
2855EE954C7F45F562960A3FF2C1220417D8F5F8DC8A980CA932E65C43FB877A8D3D0BFBE9EA57DDB335C8068E9BD17ED1F66422FC0FE8A55FB4E3EB
733BE67FBD331E5EF7D4AE82123D4BB73B14FAD6B68DF5F2AA97FBCC2C69696FFB068342C49654D3FB669284BA5617CAAC98214CA87214AB62344F04
BC0E55F5DA814B352B60311F3708CD26C6D7D09FD5B92D37D8A85ABB54DD2B3CD5C89B8AFDE8501120B572DD1929E8607950AA44099226892EEEA226
40D6A86385D268AC60C51B53D94DA9D30BD34C3839E697BD35E884F02CD4C7212B6C42082BCA708CC4D0C465CD5E3F34DE58D75F2E9B401986702CA2
CE868A72870E327378D796C5CB97BEB4A273D7605F7568E935EF38FDF8B7DDF88FD776D42C9FAF33CC158028DC95E501DB2B6A8614DA16DADC5385CC
A167D63DCAD34AC0395C400F48518D0DB5B2086C026C1B33BE0EAAF135BCCB38FD07DC2D05F8685816C0382068DA747CD81EAB85F0B671FF50856E4F
0A3D45BEDB8CA874058E19E894F23C26969DFDF005C377526D561A957A9DFEB61CAA7E3C298C42D46A19A33DD46D8C0E47FD4F8D25F1C6DB1FAF24AD
E492F669970D4FD4FE0776F89B0BCF9B3BEB7D2FB6FC653BA9CCF6CC6AFFF2F5CBB6369BFD7B8B9AED7D5BCADFF56AA70BBC7D0CA778A76F4467AC5B
07861BFD53E756D4B9332F3DB6667F0D2AF9BCA9A3B6B45452889237FD6619250BBE2CFD12E8B870F160A00FA30AE1C62D78B4051A3125A890FAE516
76581CBDA1F21DDBA1AB10485609FB0107AEC0F913713C6BFAB6F26376D960F64C04BF2C1F6714277E94CA296C7DE9DDA07D1F620A49B164722F4D6F
2AEBD822009B4011EBE1288BC3FDE21BFBAE3C71EE617F6B8CEB376F1E34079D9EB50FFFEB2F5BD969948A691C7D41E82601151D4AE03D00F3F784EB
3B2EDC7C73A3DF48373DC7D8BDB7EA177ABBAA568859026FFA4DACFA8DFE4AEEEAC9DB03E8106659287475C0A7035386AC2E4A9FD8B303861D418F34
96616E36C8187FAC4F98DDC0DD12D547F5374AAD4CD41DA029ECD1A166DD77446D68FBEE5D55A51A4614F56E5DBAB53334D63F5A4A5E5E1995FE76F3
1E8C387E3573F6797D13C9FF400EFE2CEDBDE1B85987DFAA21D79726CDA14D3B3B4B8E5F176EB1AEB1891650EEEEEFAA878C19D5F59ACAA2C6A01B83
BFA392C08B52BF5830DD28436E51F9AE4E576A013D6814D31A23BA6A05ECB7847B02049850C2854C965D353DB7D444AE0F13E5478CA14D23A76E79D5
2A50AF8107D7600D49FA2C41AA83C3EEF0A840474D9990393CF84FC8AB006CFAA09DC7423BF9B01D33413A29A8C365BC9D66E39B8419BB5A6A5F459E
6215BC0CB38498CF17A86905F0BB62E73EA57CDFF06C2912AA425ABF68A7AF93BFFED3BF3C7ADFFDB7DFF5A7DBAFBDE67BD77FF753671D31E7E4CB56
E29051D4CFEB58419217883FEC0A7037285E311FA08F53DCB2F1B525F7FDF9A7D75CF1AED3DF7ADAC933DB661E3167CEEC53CE78EF6F1C5EA046BE49
1D95F6AA6600A8BE5596F4D9802EC862B3261AD540B1AF3AB8CD18A72ADE29D2516115FC801143001DD107C32421A71508C02BE82F4140B0A02240BA
761D0DF6727DEF86675E5B5F49EC4D9BDDC147D6C6EEAACECE51116DFC9F7FAEBDF9C7153AD9DC1B3B669EB26922FB1FC8D19FE925971C3AF313AFB7
5AA9B96DCDF6754FAE18805104ACBBCDB192577F7D7341AB7245C1D19B6A7669356AAE4D25AE19C60972B36385A1701C67B4819A34C8B3366FB051F1
A3900E0086A1342525FBD607AAD4EDC54ED987D7DCB82D1A5DBB7B7D8E6B4D498C65F3A8C60F1B232CD8ADC0CB815198301A22CCFD356494A5FBB97B
BC4D8838EB21E67DAAC8032A8F598B97DD7E8048C0F40F1A004CB205EC0D7530F5092634B815A369B87FD9CF8B0B198E4FFD864571EF3A56D3F5249C
34E993B49A2B17FEEC9357BEF3A2B32F3CEBA233DE79C125EFBBE9C67987B49F7EDB8A7D45C5C863D87BA0A781A0273379DFBCD5F6C0B287FFFAE72F
9D3DE3988EB91D7367CE3BFC8413DF76C9472EFFE0672F3B72DE61B3DBDA4E2F66113B0D82B597B3F3E81A1C131AE1EC230C4451A0B8C600532964B6
128E0C8DD1E5D05E99E6FA62747F1C211DC56E7E69AC5C9FF109AC76ACC2160A26CC56B569D2CF99DA0F92C117966C1BA68E6EFB009526A3BFB9FEC7
7F59B7A964CBC2FDA7B4CF3AE28509BEEF011DFFA59B8F9E73FA03416BE0E1ADF575CFAEEF6F72BFCCC8D3E2F652D018698469656F230B62BB305872
BD7AA0A5081CE0FE78D4ECB908BA24764D204C2960614BC953FE28819AACA0EC272CEADB13BFE283F1E3C98609CF2A6CF79187633DD2E5482D8DBA47
89BB6AB26E1685A26B52A2132C8EA3A1876576F6C948FB8219765C26C07A1B1D330A746808505C3B2ACA35B7E9BF0ABC81003D35E36760C71781BE1F
854ED3E08DB88EFC9AE458E300CB59B1D84DE49EC3C20D946F49CF879B28A00CFA4D51BD16D5199A7A6F88FAAEBF607AC795BFBAE75F4FBFBCA9B3E2
E004834147AAF3D00F8DD2E6E7FF7AF3372E3FFDA819871E75C65BFFEB337F7C7CDDDA5D3B06861B904FF346565D3B6F5ADB8C53BFF706CB7E51B2CF
D0A020BE196D8C357F9CDB1824A86E80D35591A8D5E8E6B0DB27C8968A2729BCEFA4CF14D425A399782CC0BA1FAC519844D294B047565B77D9EEE8FA
2D75163FA0BF533A8BD3AC36B0A149D73B78D5D7B616C7A93768FCF3E4C9549C2C9C58FC1DC0D1DF922FBF63DEEC6FF4D1C3BF7161919E1EE074C2D8
1EDCB4BDB3B7A7D2543E7C6BA9C975ADC240D7C0886D0AA99DBD4358B855860575EF709F77FCC8633B8920560C75615E4F08609CB77103D59E210E88
34F2353BF2B25426B5B33249819B8DD1695364DB562822637DB1159802F85BF6F89580F360BE08041135EBDAF0C2D0875C30951E4C810B38E29384BB
7468F3A22A4858574B4B8FBA0CEDFBE0EC87DCFC337A3E627E0DAF061B0599F0F5409A3FE1E545A461004A3F8CB2995581E8877D8FDD00E9D85250D1
445F2334DB6CC9A8F6F25D375EF19E33CF7BEB99A79F7BEEE51FBEEADA5FFDE0861744EAF5EED8BEFEB60F9E76F29187CD7FCB82CBBFF6B757BAC72D
13BBCD283F1AE2BEFBBF71EED16D1D277CE0A697467197416ECE583790BD49053C3E79D31F67115802747DCAF7201618BAA6CAF1FE701FE7F97E1CE5
2AFF511A36E38C37812AC07C952EBF5E07930228237A1B7FE1E221CB7601BA14D5D54FACB77162796B168D20FCF77EE2E7CD30F22AF4B7F7C88CE973
8EFA573011FE076CF8B7E27D5F993FFDBC455448429F5ED46B226C510C3BE57D6B57AF1D284969577C50D062A7ABABE65894880CDB97D557365914C5
8D311FE09238904D13BAB2BE29C1F0C131E0189663DA3E256BD1DF07D3AC10AD3AC5A4422685859E0D992CB75117A976E00586099C82086D480F68DD
0D750AAD5D387BE7627349C688790C17E2C86C04D0F795BE03EA6180B8C8556CB0FD439DCF9E3CA89323A6D9531C0B500458F952A910F8DCC0F4819F
F720FB89988FB274BF0960E4DBEC39C2FC21B00BA9B01786691831D44CB1ADF03D6AF46D8F8A6CC7756D9F378F81E79BF5A13D6B562E7EFA0F97CF9C
7EF817BEFBE96367CDE8686B9B76D6BD6B773555EE181C401901F3BFAE17EEB8FA9C638F7ADB1537FCF9D5ED1695155473D0A920B1A6A05380DED795
9CB781286235C638E3922680A18ACA590111FC556418EC3730646B25F45D259C07F4E1AAFB0453949B6B477136C289204E8D9776A5B1396C86A12B2B
AFAFE955895BDC5528748D51AD521C7AEC29977E7FD970D834BD5FCC9AF5072FCB26A67F076AF6970BCF9C3AF35BE560E4F917560CEEDEBDA1B31999
567160EB08EAF224F2ABA6F00339DE5D0D6A4E8AF934AA7BCAC68EC5590C1BEB8822361554F56AEA3923651AF09FF24AE386742DB0F831BA62EFDD04
685D0ABFE68800362F5F284B99B642AA2398BE135038A7510BAC98200E78BB4D7159DFBCF8918796D4403D61507CACB12E48422D3C19AADC41CF1A0B
78EF4DE58829C0F24B5C9747826845A8C14F28D199BC1B87D84894F09CCCAD293A90A4DB7459D01FFB088414BC6E5D2ABDED8215458225F7AA740C78
A2DE641B4FFABF1233411184B62D7C9BBA144FAA3737808139F2FCCFAF38EBC8A3DAA6B5CD3EE25DD7DD70C909F3DF7EDDFFFCE00F7FFAF763CB966E
DEB075A0BCEDB5479FFAC927DF72F2824FFFFA951EAC1779678F1C0E7B601C052D4D35BD52D245C992C474262B8C457213C21848C77C8E42079FAA55
82C030F47E801F54BD922C6AE6CE64B11C015F5989C2486E7296C52DAD1B068EA266B9D0B03C80AAFCDAFFC7DE7B78D95996EBC36B259064924921F4
261D29225504C1020A1610ECFC043BC8010F56548EA082221640407AEF445A2040208414D2EB24934CA6973DB37B7BFB53DFBABFFB7A76CEF1FB17C2
4ACE3A2A64CA9E3DEFDDAFB26EED904E5A116F890A457BD7EAC6B651A5981B7C6BFAF7ACDDD5FF031BFD49EE7773671EF1BC0A2BDB37152AD67843D0
53532DDA6CC269B574E88C2C5D53A561DD1B1A6754AC6BF9BAE35B810C8C96A5D9692360B4E3E8D887333694EA74AED70373C6AA7B82B5992800A140
60072EBF1438DAAAE24C95B55B60F06A025F892AA33A1AB0C817716BE7889D3AC52DCB163C77D50953F798BCC7B4678C99462B4995EDF0C0AC01719A
403B4F99A806FF6148D9A5B1308A37916D1BB5EB48D6DDF6D29049B331343A97695B850BEC621978E58A6B90F7C6BFCB68732921A4F68A5628A1B115
F3AA2F0397BE1EA371452648815C096EFB7E40E51FF70DFA242FDFBFE0C62BCE3E61D694497BEC39B5F3F46FFEF3D5BE714625B53E3AB0E289DF5EF3
C54B4E3F60CFBDF73EEFB7E71C7EF83E471CFDA59F3DBBD1A52F473D8FD9721A3B60891A8F669C3B585CDABD0D416D12A38C279ACCA097B11D94F815
0025452922E24DEA4F7C90828DEB182EB294448C7AB002C3824620261CDF6024920847566812419D906628283706DE441D2E2D493CBEC5C674D0F37C
91BA2E95F7C5DD538FDC96A5BBCBFF0734FC83F91F9B33EBCB83099039E817292679EF760AFD96BF667D215F1ACAE59AD2704993CC55D5FEFE72DEF5
94B9BF1B65B9081739B3AF820786363BF6CCB8CE01B51E341AA0DA14C7395C303D160AD7A1A202F66A487306E580B03AD4546973382F7CCF092C40E8
225C9C23BB7BC5E3377DEF3B979C75E09C29F0FD9834F592054E8AC5041CB65D19095307391381C34238FD32CC0D40C1A546D18766056810E09A1749
1C1DDA1EDF8632039671DA660464001BD19C1B00F3DB66E11AF5EFC8207AA51F80FE12F834E12BA8E642A33482F626361352528D06A3D86B8C6D1FDB
7ACF55679E72D6BE9326D3FF1D79FAD77EFCC0B2A651E4034FA8FD7647BE9A58FBF686CDC5C767CF3EE0A0E38E3AEB0BBFBBF39EBFDFF7FC804B937C
D9D186A81853F79E04B00602B387BA28B7EAC818AC82D8750CB1B86D3FEA053E340653C082A041881F50196D4553F733344206F68C1E4D2143191551
CD71024D0D50D2BC034238BC5C64467D59B75CFA6D84838B87C0647C6D65A0EE9BD5F180DE4DF9FB40FEA16769F0CAFDE7ECF7170AE16A2F45A01350
291CE8EA698699F2B7BDBA62DC4377AFBCB1C13A0B6D01C7AC24D3427A137668807BC0CD18CB6BA860A033C5F23A4BCC160B577E05A49FF246C778AC
5B71986953C015ECFD5479FD0E9FFA856A238D32AFA1C0C9A31A28797170B474E785A79F7EF2172FBFE5A1472EDF7FEA9E93F6FEFC1597DFBC165921
D4B8390AD35368E393235CABC624CD1D9C85F0ACC3D1A2EDB067D8F161BBA2873E1009860D6088800005187A2D7CBCF54E6C3CCDD3DE78C3D5C6D70F
86D8D2F6031648F08AB10C748580436F443F457B73176CD9D6BD66FE33B73FFEE7CF1EFAA18F1E3EAD83D254C7D40F1D73D1EF86DCF6874093131E02
DAD7DCE447D309297768FBFA0D43DD255B078D91914DCFDC74FD17CF3DE1C3E7CF67CACC563291510BAF49684A6BDC533E856D4069515BD0096DF4DB
116564A729715D898CE7A8712E022FDAB1DBA10F55925445A968824A041513608EB5D5574AC32CE9598F6E00D30325B1C085E3B7844243FB58911B69
C5F1E08DBF1BD34FCFDAF3A2C67FCE97BBFF7CB0167FEF9C3E77EEA94B5BADA0E7A1DB1E5CBF6D62A0FFEDFB5FB76335D2FFFECAED0247787A407861
435715CC79E6F9D435329DF081610A36E1BA8D468D4A9271B4A17AA5CCE49C8227686A126A2435FBA0D524F091D386634B9DB41372B7BD3B88FFF368
F9FDF31EFEC5D59F3BFFA8834EBCE167BFBE7F432962DBFF7CDAF4E9077FEDB66D46AA063EB5E01C58135ED8D8310806505B6513333D0E62C0F38175
DF8A61048C5B9EF1D80A694C97622706D0D47C9C1BCC6601CB740062E29DCE78CA7730511874A3D47EA549433DA51A2897503E60924911436FBBBCF8
5FCFBD78CF155F39FFCB975C7EC5155FBEE0905933674CD9738F59875F78E503F3C7AAD2BCAE160D17F453975892C9953BE0C14D51CA64C0ECA08DDB
8382702BADF6DCF9617438D38E3BED277D2C08B0D1FFDFF704F59B01302451BF231A3C6059DA2829642DE173EC2E283FE0E61FE32000EDBF8AC41923
963C025458F6F5A229901EC70D837E97EB9BF40E2D7BBCDF5C0592D81BF1608A067C54AACDAF23C1774F42F6CC834DBDA273EA8737A5BB7BFF0F62F8
83ECB3CF9CFD7F5168B59CB55B0A4E9526D9895CDFD66A12E269333D5F14F2A0E1D0CC1EB250B37CA1520FCCEEDAAB50A566F59A638F97750AE61836
782396061D95CA2FE8F9B1E9BD2964A4C3A8CC497BDB40BD501CB3B953F73DE1FB42B7A7FF4054763C7DD30FAF39F5D839279FF2FDEB9F5DB0669DE9
9B071FFBD281475D75F79A121A674A22D2CC05C26190190D69DAD55CEA50B5B92EF51D4DFA9EB8E001D22B0A2AF43045839EC08535DE084175836E2E
EAAFD1E749DA19C010900DD4C1A0E7E1D6833505A67FE5D78406F79EBE07753E3AD3CAF506D6BEF1EA03BFFEE2674F3EF7B3BFB8FDE5B716FFFBE15F
7EFDC45934ED1FFC896BEE5831646408A9A2FA6DE3B3987A06858385F1F833F208ACE9D1BB2959A07C29D2EAF33FF9E8819D9D734FFBE99B79B08275
003955D67CF76FCFBEF4D293CF3EFEF46D7FBAE627CF082ADEB1E9F533830A4ACD0990F20253B84F726DAE16B85152F66501AEA2A2D9C48786D53E27
51B1181AF2B51234AB803994E45EADC6F85C9A1AEA1B4B5022A1AF9234DE1C69B5A82368A55A185EA1F43F3175DFC7A2DDCBBF0F64F8CB159FEB9CF9
9515F4BFCBEF0F1A451EAB58A588111C6E312956ED14419C7B56A04564D78B2E8BEC0A831E65DBB116EEB00937C21AC09B063E1C27B1BF86E50F3D53
C60D03DBBE46CD6BFB03093B700DFCCE90E0656DFBEB8FFDE0ACE33E7FE1D72EFCE90F7E7ECFEA5583456DF610D1D0F29BBFB8F77E5F98576BF7A3F4
874A267D5F1EB85AC3873B35751DB0F8D0A9B93CA835757B778EFD609ACA888F79E00B264604A3AD328CF55802DA7CAC3DD1662A2699F1E18C84329E
7D062FA043A352229942EE32E284A6128F6C5F7CEFFFFCF31BFB1D7ACA67AEFFEB6BEFAD5DF8D6CB0F5C7BDE0153264D9D7DC299573FB0D5321FC503
2ACB9A7EE85835EA22A4F7391481D4262126B2E25247647CBE6972612A6CCDDBBB73CE41FB1C78C22FE7BDF5F2FC679EFCFB6DFFFDDBEF5F78FC8907
CFD873DAB4A953A64FDBE7B4EFDD76F33C2E23E7ED0DACEDA5A830C2D00B6434138035686EA2060969044362FAA634FDE840D0CFA3B04C8884760B00
55D3E4069EB6D6368FA5C153C954E4DD84B5191B8B9FA851839126ADC6B8C48E2388DE98D2798DB53BFA3F90F15FFC9FFDA7CCFA131E29BDEE856D91
CB9AC5BEB111BB09B49B71D588BC428DDADF088C7ADBB11B162CA91223DC4DB33A8D0340B8070D8C0409AB3B8C53487941BD0E3C205A4894542AF22A
40AD8B22E37B6D748597BDFCF643F7FCE68F7F3A7FEEDEFB9CF2F1AB9FD93164336A386962D7B6608AF7FCE51BC7CC39E05B7FDB1A0272446593E25E
0AE6F8326642B6C04282DB9E847F8E0AA5E3EAB6A96564DCB81443714F12C652EA5E044582A12350EFDBACD12CAF8CD996E000D27283B405312E6000
18D36443A1A4344E8740CA41D793AA3975C53DF7DFFCCB6B7FFFF89D772DFADBE577BEF8C40BF7DC7CC5A70F9D3965FACC0F9F73DD131B7A8A0DDCE0
7928031DD22C1DBA358AAEFA8E11DBE24CD84D6176A5344888AA9F68509328D942CD5BA4A33F3A7AAFD91DD3A74F9D3C79D21E9367EF3175CAC1871D
76E0C99FBCE627E77DECE617163DFF9CE90A42B774CF41477CEBA70FDF7FE78383D0158833F01239A51069F25798421FD1A08015C38A80D57D8C4C52
61354B6DBF6145731F970C466F0AFD56E89539752F70AD6A3D94755B0A6B5D29B32A9418791DFBCAD80EC5697B1ED3BD1BF7FB415CFBEB459FDE63DA
A98BD2565477AB132A4E9CB1F5055734FDA61FA5B13F3432B27D43DF902B43DF2D37CB8DBA8378A728894CC8295F506053ED17D213989C991B266146
3D40E0E0449648EA4779E0E4B7D624E46A70B5D7CDEA60FFAAF71EBAEAC0C993A6D0B37EFA4FBE77FBA6FC4E286D1233EE483F143A587DC9F449530E
FADCBA0C71CE5BCA1234262414AD558F8534DD7BC500CAF9484A63C32E2EE64652107803016E9F0CE2D4F36884068EDF07F02FE28E4075F730CDFB14
DD11487421AFD741896196305CA5380C1C6E0E6846F912F7F7E1256FBEB770FE8377DDF7CF8B0E3BFBEF3D7D1B57BE32EFA6EF7EE5BC93CE38FE9C0B
7EFCAB87DEEC2F1B486ECC2BDCA1DC023B31DBA3E23ED645DF494377CCB78DD89FF13280A52E22D12093A571E288B2707CEDE279AF2CF8F7FD0F3FF1
E67B1B17BDBE68A8303CEC0288DB4B6F8D1E185F73EFEB2FFFFCD47D27E3063279EA2167BFC1B5E12D1BB55565D8CDE69811FA812E57815AF655ACCB
ABB6EB2C73AB003970C3010A8DCF5718A65012804E4124CA63751E447E90EAE11D1E350C2C7255EAAEDA1AB59A569846E96FA6CEBC35DC3DFC7FF0A2
3F2DFD63CE94436FAAC4D11BD73D58CD983B5EEEEFDF5E8378B6630B95E45F5C51729A0A603DA73C3ED15F68B29A2F4C6F4E352EF24A458A774EC534
044BC7E876A7A968E6CCE3862B3AD59718518B857B905FF2C4EDB75F79EE3187CD98BEE75E7B1DF791737EFCEF37162CC81B607A06C5DB2CE59CE16C
1FE59EF9F89449277D7F45BEDF8B9808A8624E04AD0CAD81C798A8533689AC212FCDA07EC3DDF1511B2773A31C4EC370A58101001788DE1EB311CF32
E0F8E905F910DF6E0BF3C8F69C10295CD2329868542CCC050618677CF7E8EF125EDF32EFD7179F70C09CFDE6ECDDD9D1317BD67EFB7DECAB5FF8C819
879DF0F1EB1F5CB9BD9AF7DBFE6034E7C0E7B7B0C5839789C3639EF7681C6FD27F423FD3DCE23348F6468030078ABB0135250D1FFA2691024EFA3F72
9AE94EDEB4F107D651CBA8FE8CF63CF7897D3BA64F9D3179F25157DCF0C4F2B28DF5A4198A58E2D727EA554F29342EC6AF5020CFC53687EB6140BD4C
B50E6DC46D3BE09104911F7A3DAA19606F427FA093AC5356C1FB5F18910D27A40CAFFCB79E6B86F5F1A899B6EE9D32E5B4DCEEF0FFC0457FE6AEB8A0
63DFEFAFD7BAFEC8ED8BDD341D7EE6B989320B3CA739B46283630064E6B7DEDCDE9D6F68A05B84E6AC38345E1D9B60F484496A07EA56ADE4D17C0982
2035BE0653E72ADC9C92F602993BBC965BF2F06FBE75C2DC3DE80F95AF1987FFBFBF2DDC5CB5BCC860DBE1E4C529EAE1A1819577E00EFF71FF43CEFF
CB7C9B6A637D827A7C9E64128E20F49706CDAB23C693765025A6F26161666CECA1F945F5DF4CF471D6B2CA183F80216C7FB051C0302C401A110C1A18
7702DC3F29171A9D3CE3850BEA8CB96586CBAE38704647C7B4A97B4C9EBCE7F4A34FFDFCA53FFBC323CFAF59D9B3B9ECB6B57361151A1A2263469FC9
EA4C409D346E3B6C28233A0ABA9E99A30CA058635E82472AB3E0E0D5520DCFA87BE22A0FC5412E94A677045D7D3BCF62EBD8028DC1DDB6A277786CD1
E57FDF91779D66DFDAA5EFAFE95BBFE1ED1F7CEE8B47EF37EBC8E3FE46F31563B16CFB7A85123B50739FA4942CA12C120E8C47101AA21FD31D1EDBD4
6DA7D80B8421B41513DE354CBF8BD2B8926E9CD6D7AD741A13326D48D58A5B4F4D9FB1F7A3F1EE88F9A045BFDAF6ABA3A67F74090D9085EE7E0AD5C6
E6052B46B4482AAF7657D6BED5C32201BC07D4B6ACFE2D230DAA0938BD69AF325665F9A126A8AD4A88D258D982822EFA4F10F052D4B15809D7A95507
D63CF3CFDBBFFFF1E3F7868BE7A4BD3EF9F5AB7EF4BB9FFDF2D96D79A3DC035D2B2DBCC0F3C70673DD83E652D66A89E55F3D74CEF1DF799FA64EEE54
AA01D301A4FB43CA24A9D124375E38947BC067011B46190F3BBFE185345E97AAF0FF04EA0524380AE0C8B8E28650D74558E3DFEA36ACDF443FDCEF0C
D31F1C80A82D9DAF77AA6782BAA77A1EB8EC9C733E76C649A77CF39FAF0E5971F41FCC444C530F8D0A92229D2AB431DF100904B918A34C482F4F4498
C98DCB094D2D0A0E5EED570B1415D8FF499C85D22DD670880035173A8450E012B8C7E3D62295A64462AE9D29BD30E35CA037DEF085534F3AE698830F
3DE4B84F7CE1C25B8E9FBA27E5D403FEB2300F053373AB8961D793C4CAA80925403625C62AC0AAD70A1EB48682159B6D0B5C00D96CD402AE5C567DEB
2DCAA8D6B015E758A6D7AEED1692326FD08AAA69DF01D33BBE58DB5DFE3F507FA8E2D5EF3F6BDFCE2FF547F1C84B8B73B9BE25FF9ED75DA022C5476F
7FACAC2598E2541A1C2A21661F1ED56B9410FA17AC69F8DC105C28D812A8FB854C475C4AF87668DFAA34EC891DDB5FFBD3B5979D73C69187754EDD63
CF0EAA9DD33E71C5777EB1BA86EB60060D7A4D5332CDC982177240E6F14AA3365A01B330ED7EE0D333A67DFEE1220B1B76ADE64937F068BC000E8662
0B7D7A8A8D1F9876D433F0188E1895020DEFC27729D0A25AB785BD3E3E92860F387219EFEA16C567DC16FCC01FE37C41354F59029241800A9BBBB90E
A98106430E83B492710B2AC1ADD09301ABE5CB91E9158C2E9E210A281F867CF03582EC10969F66970EBC00F6224C61E7686217EE07CD208BDA5CC3D0
1FAF19C403EE92341505C2F81745D247DB601409E8A7642A62D8B34239CD6078909BBC30907E14B2AD0B1F7BEEA5951B7B0B5B1EBEA563EA949933F7
FAE1BA2D9B56AE5DB1E2F17F7EEFB3271C77E9B7AF5D5A6CF39E53B81CC3519D7E5BF9BA5B75613A5A6FE0942F837C6E70648212D8E8E62ECAEF135E
855E7D2D4EF3C3D4D1F4DCB619F75BD1BAAF63F6C10B770BFE7DC0AA7FB8F15B7BCD9CFD8391963FB8FCFDDE950F3E39AF2B4703AF0EC6962FF5B147
A29A99886D8BC72786B7D7652885C39D66A3EBADB5657BA24A3196C09982B28029AA3ADC7EE3B7BE71F1A59F3EE9A0A38E3F68DFC30E9C397DDF838E
3EF48823CEB8F42FCFBDF4E8CB1BA9504BCA2D2C624E4388C0633E17CC0F93DAB04ED3D00DB4A374166FFDD379B33B3A2F79516709032E88C1B7CB30
6280E3FBFF25AF6427331FFB7ED9B0A4749A712BA219A2EE45926BAF81D3588B5297C75AF83C9AB27128A3699A89D8D80218AB3DBF061600BDAC0087
CB90B2D8F804C48A3854B93009709964A8F8D032830589442862A66985C53CF4C01C61B4F59A3C5199E6296ABE113A05622F31AA9E515BB713B464EA
3234B4BBED5C601410B13005F701FF4FDFAFDC4C8C9E37FA02F876036FC06CC76BB73B948A54AD4AFF1D601BD2DE0ED04B7172CFFDE8B2ABCE38F6B8
CF7DE50757FEFDAE3FFDECB73FFFF577BE786EC7941917AD1E653B6D0D50FF6315734BA6F47BD5A95DCEDA6796AAD5A4DF61DE95BE480A9BBA8B41D6
2CD9BEDBDBDBAFD9A6DE301DBEF5A646EBECCE19DFD9CDFBF9609DFC5BF6834774CCFCF4B3549DFC3E2FF4EB4DDCADB457DF3150D39886257CE703B7
961F1CA67696A6FE304A291C8123B7776C2D53EF5FAF701EB51572E8512C3F7BF31F7E76F5D7BFFAB96FFEE29A5F3DFAFC1DF72CEF1F1F182F38DA9C
FAD244878C9A5B9525A201CB5FE9D56A06DE4FEDAE5FCA971C16965FB878EE9E538EF8F9D2183CA0B6E60748B7B5C2B21BAE7DFCE99B1EDAB072D5A2
456FBCFEF8630F3D73FFED7FB8E3BE575671E3B5030030347FE20CA6B770012E0CDB51665A74D9EED7130C24A8BC8D3A435C1B456E6C0BCC946FBE17
10395A376C807F389391413187D47243355C2805F6AD51193748DE56C890218122D63A09FAEA9014F5B4C91BBA6DD71753F6301822AC4469BE885458
838E377507D2648890537A01C23889AC80D28FAB52E30496E8443A944B6810A90F578CE629EEFA1019D0984A6820603EA3D90367049315DDBCE78135
84EF8B192969FEF7B4497B9D79EC85D7DDF3FCD22DB99A60A303D42D0988B7B57823E621780BA91E9D80DBE0CA3B9B71EC36BB561574943A92BB9CF7
9663971A9BC1FFBE694D6B68CEDCFDDFDD5DFE3F48B5BF956CFE52E7F4A3EEF5D216EB7DF2D9F7B72FDB966B8E8F742D5D531D68707AFAC686B76FEF
1A2E957AC6C6DDACE56DDB30D694CD72C14347C995B07DC527064A126AF908D2581B3B2AB0ED032AE3D1FFCD1849AA94AF6C4F292EA18D1106550F9A
3C4C7A752CC6A3CCDD51B10AF9AABBEDCF1F9E3C79EAD17F2B9BF4C4C697BDFFDE3B4F3D74D30F3F7EE4F4B933264DA2BF9BD6D931A773EEDCC38F38
F19CB33A60CC3BF9840123D34DCDBF8E27F2898CBD920F0080B44B5C87A51C4B58088DDF80B2142F56206B5B778CEE9ED41C9780B61F6F6254430DFE
3F322A23E03F19025068DCF614BCC6A9F86B8EC925C9EC227D51CC476DDC3CC682D0901C741C488893292344A295570C5340FEDB8E1B5A783D762200
9754619B83047B5E54F27C573335FE82917131804601D4C538BD96B02D739E41BECC2809F92E7D79164AE3AB44B9A905315020A32823280EC72078A1
8EFFF1BF5F5C72E74FBF76D1A7CEFBC4673E79D1A5271DF7A59F6F6621B417730E357F86D8D0BDC1A304585DB8D8E241AEBBB7E4258CC9A45849336B
225F97B1FDFE387DC73B66749CBFBBFC7F90A23FF51E3BAC63DFAFAC0EC3ADAFDEF5C83DAFBDBD78CDD6256F3DF5CCE3CB26844C1C59DDB8A3D9F01D
A1DDA21DA4B23C98AB0616ABE6AADC0E74220CD14731A3726D743622A39C457901DBFFC8943F3073A16F4DADB555E9CF730CC79E1DC5F5ED65CE40D0
038A00DBC55ABF4DF380F38FB993F6D86BCE5E973EF5E0DDBFFF9F6BAFF9DCB19D7BCD9C3375AFCE13CFFFFA753F38B6F398AF5EF9C003CFBEBE7163
4F5F756278CD77F79C3C69EA597F5FDEE0BE0B909D943A6F4592B936050490BA781D8571A0645B995F9C1054114B01A6708A796D0FD3EBE19E0028D8
F426463B3B35AE97A629088D63A79107C0AE1CF22311F7813A8451772A6C786FC4B2ADC8679C44F0D9C0E219CB91A87DCF840E8167F0C4865E4CFF9C
305C1F64A5AC8D8B2F04BC1B0DD388A4D4F198DB639A8ABC87B4445903DC01EC24F01E87D4BEB4423816141C7369118C799CDE736C3760D4415D4F08
1826FD88B87566347E80F713CBA0D0BBF8857B2E38F2CC276DA884A6EE388C8FA00EE00CFB692646E9D5383D1BC60B79F3AAE0E24663DF808AA38145
0D6496F366CD7D2BD91DFE1F94F06F25F6AAEBF69F71FACADC9D37DE79FFF21EA3E8AB45B134EAD92A0A94B6ED9A4FF36DA62BC3754635ADBE2567D4
B410C2C3AF2D1E6B288F1ED92635CA0CE69252C1F33A86E9B7E640A2FB5219A75A28EB50DDA787152D37ACC37C7A802DD77120E1058430906882DA4C
26E45F664C9B3299AAFC1E93A7EC3175FAC99FFEF217BEFB9B5797766F737514FCECCCDF2D652A81115ED5195C76FBC5C7CEDEF78BB76E041A11AE5D
E0F5614A8F5B18BE0356AED49B150AD12C4B54A42618744A12BBBB00565F9448A5ADD19232730BB83FC68DAF8D1D3202793A34CA98F417C1B88890A4
34EC75126530753BBD430C5F36644D1E997BA2712683AEAF8ADB5AC1519B374CF308E625EC168CD569DBD04F36FCFFCD9A104188DBDEA511350EB866
249286746A26442458D4F645377E8036AC4A90AB286DA95079BE749A81A269C0EC3721EA0FD7D4B69B524C630F8CCF6D2993A81DB995A146140A8896
597E2B5589AE345540DD495C188BB244974A928D5B317E2701E362FD7246FF63733551416BF5CC69DFF27797FF0FCAE09F398BAFDC7FAFCE137AF4CA
BB5E9F08325E5CF7C693AFCE5BF6F696C1FAC0D6A1A0F4FE92FE6242A36EF9E5BFFEBB0A462C15718A524CCA14BEB5BC570FFC66D7D2B7B73BD4F40B
28DD606B45B36EE8D783083D004D91B258019F3489DB3EBE46D3C65CE15B89F6951ADD81D336667783D6A5FFE12C7BE4864B2EBCF6377FBFE3D1875E
59569740A2B71007816F95E929D620C934963E74EFF99D1D1DFB1EF2CD859BB6AE5EB3A66B5B6EA2D4BB63E5EA794F6DAFF72DF9EB1B4A353DCA0BC6
104342A8CB5A940723498F7615A18B274369DCEE8CF5251842AD847BD00985DA1622356E7B886167A1E1208AA24AED3418373BB90189210E623BC0EB
CC8CF906656CE435CD97A63C6097EB8DA66A411D10A53FCD6FACBA8132EE1DCA3026DAAE0346A3D3F896525D97A2CD60345D00E55BA87DB7B79CF8BA
C2C730220CBF028A04348AD04FC22C902D332C11A3C86D005118B6298C8974CB6B4723A3C79809009981F6E5014E1458AC94D7F6E3FC17A74D4DBF12
B03386B67BF8CCE2369A3AB696E2383F805F5D1A9C36FDC055BBCBFF07A5F8876BBE337756E78C0F3DF3EEEDF76C2FD7973F79F3EDF73C7FF763FF5E
F8DE86C5EFBDFA6C9F1A1DAAC8C67835D6F5BE1D551ADE4518D860C5E3B9A7C1364D339564B16B7BDAF329F8B90D2A7BAA2C87EBA0E1877649A58A87
EE5879A7A18F9E180CA0E6192B7A2671DF8A4520C7377BA1AAF6D725A24180558F5D00F3B431BF4435CFC25612BA1CFF9EA25631CE686AA0D9BB55FA
FBE7AE79E0DC1953274D9DD2D971C851679F78DC91C71E7AC0AC833E7EF6B1C77F63030DC0D47530A3BB6DF43F43AE53EAD4856733CFAD4C8C355904
7770A301821D1C35F21232FBA6F67B018EFFBEC4D38E164619A50F7A7989C2EACCD86C213E0D6ED0909A3194A7E834503A134A6ACCB8938FBCD2ED72
E3731223B612560B02E5B03469FBEF68A3E8ED358C4F6F6CBE141452553B1B9AFC92F87E3B51989E83DE1457EB580C3453E39F44938CC036D2AE41A6
25644656C03136CAF824493D42296763FB402D83EFA277823649B564F69ED470F14A7F93FE4512D5ED3802E6DAB77AFA20A3547D6D887A321AFD4A5D
6E9C046EFC93E97BDDC87697FF0F46F467F5BB0E9D3D6BAF23BE7BFDD5975E79EDCD7FB9EFD7776E7682D80BBC9A0AC7B76CC9297870A48D2D831C22
783E28AF4963B8EC3562A7D8DB6F51B7CA559A81BD9F248147653A08B469636DF8FA50850A6C054B4AA3336D84F3A3DC000473A9F76DC13704950EF0
17059C1B36636DD74A0073215AEF0BCE2CE62382B1774F522F2763CC0FB1A2F80733BDD5D259D633FFD13FFCCFAFFFF9C492A56FDF7AC58F2E3AE0B0
B3CE39FEA4ABB7EFFC3971C937CE98B9318B0283FA5B70F5695A6101B4BF6263AB61AC40107B388586C687D42F953CC179719081AF988028CC0457DA
AC058CDB1E6BD3F5E39D4EE4A68627700240A64953EA3002737367767BC636D53D3153469270239102ECB1D11C0EAD224775375F087B7F2D2084A88D
6F411ABB7E649088ED3D61123601AC6470508FA4A7802980569F410F4375318390B90EB053A584AC93B8018524EA0D78A9EE53DA8FB50E2A231EBDD9
7CD086A02137FE2BFE906E81E458777C87A61DED6DDA2483C1DEA682ED6AB87A07BB75C6ACD377EC0EFF0FC6CD4FACB868AF19B3CF7A6053F1BDF77B
2BAB578D7B54BD33468F282F6CDBB1F4F5554EC41A35B75129B934C772BF414588CBA0D0B77DA058B54A16364B709044B72C20A2A9985F6D52EC8852
35841A06050AA7648279380CBA1797E3F6790C43A934669D685055C470D6834A45689EDF54C3ED2B61810C24CFD7985BF3A9E12D5165922D7F4CB722
CCD41C24BF50C49A9A716DCC7CE08A832F9A36AFFCF0478FD86BDF137FF6B7879E78E5A9C7EFBEFF6FBFFDDE65975EFE952F9F76F299DFBAFC82AFAE
4F50C985C1D767101D44489BFE9962CE1B2B81393FDE5F63562D083533B6E050DD36970D4C31AA8DE0A7C120D16D45F0C49C0C00F8C33E4F30CCE1D8
FFA36DC77F6BB3F04BA0CD11CB14EF5A08C8916EFB92E2FB81A2A7B01B88939D0E5D9166062D88AF052D0F695C3F8113306B00C7F8126A899E03321E
46943C358EA0065B6C8CCE81954006D13BB7927196548ABAA500960CFABA034AD5F6DA570B91F114F74267A06B82C92C29556C81AD9FA3793D101BE6
57E81BF53776FCE196AA3A6FC6DC7FF0DDF1FF41F8130E5E7FF8CC39273E3822FBD76F1B1F2C6CE8A152AD5D56B3464773760DD45CC5C75E7CF4DDDE
124FC2DCDB4BCB9C3ACB40C7BC52A6C9B1857A029080C55B99126D4F0C85CD51A4BDAA0F173E006E115A095839AA3962476D352AE3C9C32C38F2B535
3AE0F7A9BD66C375200884EB7B16C606EDDA5C32A2246CF798032FAF98621ED6B682D59A92B920B8A1FBD6BE6BD98ED7AC573D271469E23736BCF4DC
DF7EFD83EF7EF7CF7FBEF4E0FD3BF739E49413CF3EEE988F5E78E619277FF5B39F5C04CD1B1A0D8045504960E5CB2C82994675D570287AB6B84AD7B6
F4347D9A6E7428D3F6E5329385BED54BDE5C38625CC3B0DF4FDB6E5A713BFE1373D3A39F469ABDBD59F387695B793C5446BBD77C2C9A204CFA6D8F03
740E46CBDF5C12DA8701446A080722693E16DE66A9912341D94F0C8B11B23DC0259A4D6562B409B0763078447A3F21DA6BA19948DB56C6CA408E0DF2
424BE01BC166F04B224D9557AF5B9272060D0A826D5B3A94462DFDFE12370BB615783D561339BB39EA44B1787B7DF7F91755D2E20933CEAFEC0EFF0F
42F40FDDF8E1CED91DC7BC9D59B9EED1DE91C5CF3D3EFFD5052F3C71CFFC85AF3CB56644F04899E2645782B0D668B2D2BA01A168E40E35B07169E82B
97B900D4C6821B0629D6DBA9744A8E6ABB7B6896AB9AE9524466CBDDB6E18C95F4A1B78D195F43884F9BA3378430946B7317F76FA8D9C15817F702B7
E0726CA04381D55CC22B25C7169247A256ABD79BA3398FB24DECB89C398C352A81654915C0C4EA3F7FECDE75AF2C5C5F2C36A869915C521BBD792028
700D69A2362A41D7EABE049729F75E037E430993584706A2ED6CAA2A5BDEBDF9BAABBF7ACAE1333B264FBA42C008BCBDE7FBCFFFB7DDCC12A0FAA8FC
9A2E06781A9A9EB80F4C24205451A20A7EAA117D504331FAFC49521F6ED07BCAA8D308B01EC5850ED025980BB44DC8A527E2D4C899A4C67A0CC6A5D4
3F407BD810058D91790C3834F68E21BD3C7C8679BF817D823D81ED52EEA0A91EFBCA2433A40D657696A2112035C0FF90D73CED592C91ADD197F22D5D
1E110D27AEF79525E761E4AFDBA257FCA137E45FEF98B320DDADF9B9CB9FFCE3DC0D474C9D39F5433F7BE9C69B6EBEE98A9FFEE886DFDCF9C7BB9E99
FFCEC68DE3BD3B7AFA36F62F5FF0E2C2954D852581F495A7A962C4F434056D5D3D19B722D7771A35C61968C1692BC942B734D1BDA160D682A0AFC492
531D973EA8F489510485AF0E920164EB927639041387C28459B2ADB103860C05A6B0B19886A3B7E40EE436D3A066B73DEF62A1554A7DB7B6037BC7C6
929FA5DAF7E155E179C219CA43B6B61A44A9711CA7D99F2580FD69ED0986F960E28D658B5F7DE5E9DBFE79DFBC17DF5BB4687B736CB89961E9E7092B
601EB4BCCC3D0D2A3CD4C76C7EF177FFEFCCA30F3864EEBE9D1D202A4EFFF8FCF695027E634AEDACFF1113B1710F4C136D229BDE2DDC07001A2C57DB
7A21601BCAE18AB9EEE1A880291D2D44D7A202751ECC88FB1AC35E2F10060CBCF3B210036504065262F4C9801F4ECD222508E2A82D629A185D738005
13596FC074153F7BC443DF8296A216AE80F67F736C475D63039830D77833275699614E80ED52E40DD19B98B25A49593A93AE1DB51CA153A7A78B1A89
F1D70BE1F80A270C6F9FD179F5EEEE7F970FFFB4F9DCE9533AA7CE3D775BD27CECF90DAB1E7AE6A5856F0CE79AE8C51D6BC42E145D36DED7B56953B9
B77BEB784939658B35049421F24B37D6AB4DEC7F4378715B133541B55EAADAC88E62636CA03851E598914325D0C327F4443348FC29F04F71FFD7AABD
E10F2458391CDCE0809A6F0CEEC80D6DDFEF102CBC56067B5DD7F78136E6D447380D7C4006376B733760258B6A13E07540F8B1C02D5482DA40954285
737A75D80CB6327ADA331E0ADFB62C4B947F72C1E78E9C3E6B6AC75EFB1E7CC421C77FF8ACB3BE7CCB5766771EF85F5DF55C53865CD672BDF644B16F
C79AF7E72F786BF1FC979EBCF5F31FDAFBFCAB3FF3A1FDE77476ECB3CF319FBE7D1BE8CEDCF08A13989846C672835AA1D8EC00D16B8BD0EC321273DB
8BA34680C490181DE1D400798C8B68A8773A9262EF68148753E37FA21C3B88C2662334446118F8B611C2692B04FD97995B610299FE5A6127162831CA
64DAB8ADC486CB6CE042DC1681833E46689C55EA3BB6545C30A668A2F01C2082A27CCDE01624505CAAAF97C63DEDD76AAD281B2A365A9CC9561C0E6F
2C467AD5E33535FFAA1E21270ED9EBC4A16CB7E4EFAE1CFC14106AD30FE774CEBDF8C67757AD59FEF2634B4656BEB3C38E26762C5D3DDA3FFFEE5FFD
FE4F4F3CB5DAB6863CB8F28E6D5BBA6C6DEFBBAF77D7F2C591922D9CE2486EDBF6ED13155E2BF10CBB3B5CA763D6FDCEFB394E1D356C7169FAD41E8F
CD552C046B8E7A7E0FE2974A526048147E532603E48958D72BE61828E236284056E9EBB462D8F9A51239C22C0AC3C820725135E969E6200C35276A1C
C18F73630A086CD3722C2874EBA8854B7B2BA44EBE15F99449A4AC979814FEFC3F5CF3839FDCF6879B1F7BE595C5EF6DD8B8A57BC98F6777764E9B36
E3D07D8F39E792F3BE78D6C74EBBE0F4530E3B64BF8366CF9DBECFEC03679D7ED1678F3DFEB8BD0F38FDE2BFBFDC3D58F33534CB04350B01C07A66E6
0F0D73096C222007D1A953B94EDB4EC206D14B89226D19B4436A04C8A1436EFC07DB1903A861B32980005ACC8300FB518F21920D2808F390515BA941
05C5DC47D00B8471B980462B32431A5571D88F99059F31093552EA306347B7059360AF54A44E29A3772F6D9F40D296E783191C49D713A1EC7DBB116A
1F12EA6274C55AAE5B92BA306778827E37134B2AF2A14BDE555EF8F28C7D1EDBEDF6B94B37FEF4C815FF75DCECE9E73D459D617E751FABF1346A6543
3DFD5B9FBBFB91E79E5B34EFD18796AF19181FEED9F0DAFB795BB7425DD8DE35E8C4D22AD4CB4D47D2A3A4EC6A4D70687F88280D771A67C0EB074773
2365A5AC52939E728A5E2123AA41013C68EB450F6B262A6B5906203B5A6C57499B01DB02A70DB4CF91BBAA44514F432A187EC67D9B4BE0D7A0638589
99718F49ED95190D0E34B5A25708B190A78940D49B528BCCD875B76822F09B2A55636E8BD28E2F2CD793C65F0F69A88DCDC95AF9EBE66D9BFFAD534F
F9FCA73EFBA96F5DF6E31FFEECFA5BFFF48F5FDEF66EFFF8E6CB67CEDE7BCACCFD0E3BF0F44FFC7403247D0CCA2650C20F38CB8D09C0F9A8A23B5DFD
AE54660BC7FA2B58CCD1FC8EB5BB4116E1CAE8C9CCB0073253FB1328A2E39F6229203696B4D103C6A90342068A661BD31A50F0F300E2603033A00E68
BC9AB5CC0511EA1FF4E5B5658560FC50FE83A53AA71EADBD89845FB9B109A62FA000C854BAB0AE8C3718848B4098F56BD8CA8AB95696646234CF1475
58DDDD61F1BDA188B7AA6F2FD9D16885B668398B1EEBF25ACD5583D1B65FFC7915FDB2FDB33A3ED3DC1DFEBB70EDA7E877169E3F7BFAACDF046A7CDE
82F7AA425A767FAD51E3342BCBC06E04A90882D0E715FA8B756FAC6FB098279E63D1782F12A52CA7E65A2E366EC0FE29810D5468B6D678B6215211C6
D4B4D383C61A35886F2AA84FEDAC61BA51E43178ED5423798351190A4A0C073BF0EF42C6686E67F4D80A07572F688C4BE5E5739689D7F6653E0D0513
013327359A7E438EBC810A08381BA489609325E3A805154F41EDACA417E6B9145122048C95DBCCF3ED80D56C9F9E79C954A25DECFF849F6A55A7C424
8CCF8FF2CBEFFFF9E2FDA6CD38EA133FFEF3E22DC2F4BBD4B128FAD6825746A9EDE01CF73BFA3EA957836C09BD06A144C3C50D200AAC7C811B5C0F16
FD22508971D6CB0C8D3FC3A510D3B87002CC3B6DC45FC4192E8CB502750ED81DC02E25A6CE42856DFC1F1AF604EC5EF086436372AA207004576441BF
192D0C2822445791219962C780864B8010D558B24AD06BF01B1248C2FA08471BD2DC42E5BFC5BAB7DA501C8CB9575BDFA5531D0DBCDB6889C40D5AAD
A537BD9B24E38F2ED11377BCC554538A6B67ECB328DD1DFFBB70F5CFFCEE5F1C31BDE3D8F9A17EF5D2AB9F786CF5D217EE78F8A135CF2C62F4F734AC
F26A79C8C6A0DF603C0C2D4F377DDF63A935E268C7AF8F0F57989046AD121B6A11023E634466E9B9D29ACA152A36B7CD7A9F35AB3E9E661867F830AF
03855D68D0DE552D6781A9EA0EBA40FA81F356DF3EA042EEC9B6F62575BA19C5366F54724D1E3030E64D0290AE0F0F30EA3312A702A13F694045F422
EC008E83D2E56986095C3B36F67F40FA52DD0EE1331CB10A55EE463D50CEA8AF41E8870E50CA3CAE80F8C95CC7F391B1FA5EF8F6B153264F9AF6E93B
170CDA3BDFB994E242447ECE5352DA96A9F7901FA248C05802CB33E000E80381DAA19184F1F6A22E34E87B49C940A5460600D8DE0C373AC580DFD321
2512FA4FC53C4F817E3BE4D3DB222405796294CB8C1E3AFE43D18F0299C2481BA7AE04B74FFAEC464099C387BF373E481B750F4AA9319445A0C74C2D
1284CAF31696068107ACA6EEE9C5F2B5DE2FE8D547956D9ED9A670D72934B234486A25CFA3E12D6291FBF6A2B0C596AF81FA471CF637C2A7A6777C47
EC0EFF5D37FC5B72E45FA776CE3CE64F853858FED7773CEAB31B8DB2F41A033C1959B8BA54EB1A66627478B850A7426D3168F832451D762B5CF36AAF
EB8EF7F694A078818EDFAFA29357A8976DE12CD3AB1BD78C24E8E9C742396CAB54A5D269FA80F651F153210D02C09E06062317DACC91D0AA0D632F5F
14924215BB2A6B22E7D1E01E19346044FF26A6C620458F41ED6E42B32EBE89570B8D7027E65FEAE487B663272F7C296B65A3FB6936F0A658655EB3A4
62A559DD8D14EC747D07C733C45D942AECEA80E54BDD268D1583377E6AC6E449938FBEFA0D9AB6330392822E0E8843F4330492BB826A2EFA72CA6C94
078D5B9E361167E21FCE82261F0150280D9A893263D4F236348CFB708AA183A67F01851F13DFF840C1394336A22E49A9A4C51875083B7B033873603B
82636596D06B5406379048BB62DBCAB7E83B6B169ABD9FF935D04FAEA2C001E01719875EA3C12BAAD060048C67D028407FD6B8840A08654658FF39BD
433EB7A238089A9205CDC8CDD1DFE6368E0576B194A57547C763EB74FDC419276CDF5DFE77D9F06FA5A5073F336BC68C1F94D3982FB963ADF62A3B7A
2AD2D3D6E008935DFFB8A3DB19C98988F9E54D2FADB7E9810C9C7290C6013D3F29CFDB14894123D0421996EBE89B3D4C2BE643B54FD51C3CAA914383
2EFC25E3F12D40FF84CA0CB278A6239ABE151CF7826653C40A0F1EB5014C2AE3CA09BF60B4C489DBE3C4B1CB5CCB66CCE712AA7D40B64B66BBB131ED
3542B6D89EC12A030E7540C0A1594E44802E5C7AA1DF95CB528543B811F8CD52270FF74EC54C5208B19D54DA7CB60A3C412DB5F242A888A6754B45FF
9E3169D2A4CE597F74CD5D8F79D64EB25C66943F308A731980D718F2800B370EBC9032832837E845F0B180062C6114BBA8D10984D50C8D11327DDB4C
4F08A06E6848D25E001291794530EA89FD7CA0904FCC56CEBCEC866324BA128821D24B55E62ED816438FCC9784EB89640C8D1A8D614141985D23229F
DAA9260BBDA66B1B4090E8CF33A10DEFC9589AE26D0924FD4A5881B7B2288B030B3A6D43CFAFD369A1CF2DD4D2BADD3F1E57166D4D5BDE922E9AD006
982C0DC7F1F0939578F5ACA9BF0F779FFE77D5E86FA965DF9839A3639FDF4462DBAD3FFBFEA3EFE4D6FCE2821F3CFEC2BDAB6C1A5E71ADA74253DD34
1AA632E8EF1EF50220E643556EC0841A029EF440352A7E9B131BB19A43D58B07921EC8665F039527E8B7A17AB713D1877E14CEC0069F97C4CCC6C63C
E20DD8D7D32C4E1D67BD405F220AEA10E18AC21645477EE1386701087960D74681CF8D33A7AE9564CBA0EDB02B68D755F4FD583798A207942D058E32
B4425F665846E2DE263C15BB132EB80968167C6E98B9D2079E403A0DD6AEF0ACADF62BB9D4BF9E3C190EBDA77DE5928BBFFAA58B3E7EC647CEF9C24F
6EB9FF865F5CFD8F7BEE7DE8DFFFF8D73B5B8B34524B288208838C808F1F7C4AE9E5CA5C153314E6762D3D3B04BE10AD12BA8CA49D437031F0BB4670
33817C08FDC3B609961F10F1CEDB20D847D46A54192592DA98CFA4E1F600C6EBD4DAE3467B1A81079179973DC7E2ACD2CFCC5A107E24942E19A81550
0CC30960FD98D9759ACE8C3268E834607144739E32BA85FEE235BE4EDDA503F45507D6171BDAED5DBBC98FB6BDDB14F670BF522363A3559EF168E8C9
C1987F79FAF143BBC37F576DFDD389BF1ED2F9899BFFB523693E7DD7036F2EC90F0E56C7F39EAC8E099C93A18F11BA75AFDC74AC3ACD9ECC19C2D018
070305F0DAD089A7B231365A0FBC7C1EA256E8343577C0D0F31DB8C545BCE29A663C35207540D3B5BF65243670008A48E0D3CCCC90647E9D2A3BB76C
1D046A629585E9D6F88BB768BCE6428399A6431665AA0E6A1EA33617FEB386814319C32BF9588FA976A419F67DEC7831566501D3EDE403ACADC12125
40EAA3EF4D42B7EAA21F575E00255DC5235352332A88145BF4BDB950DEE6E76FBBF2B24F1E457972CEEC638E3EF1D0C33FFFBB5F5CF7B56F5C7EC3EF
3EB3EF7ED3A6EFB9C7E1F78DF1CA6895791AEC9DC4F89B43082090D8DF87569AAADAE010BD8BF0240B425E0F347441E0BA69B0B9B12A796982E1C3B8
0D3006AA83B12745BCA65039531E124630E180B768A0C1A13B3A1628BCDADABAE53BB04854666717C27994D33C220D2D484A002341E235900170011C
F82319D542C33D6A2E59CF76F6506033F8EF3E374193420D566ACD9523AE27729B5694D9D0CA829BB3860AFEB6F536640F46D6BF301E27FFD331F58F
BBFDFE76D5E8F75E3B7BF6CC6B9D2C8986976D61C279EFFDA2E09AB903454B4EAC7D7DC1FBF9BE9A44A8467EC48427A457AC0E0F34225EE2AD36D785
A677CFF27CEE4E8C53FC29A3741752B4C2C3A7CD58C1B41F1B5F0FC32945F16B8C6246C6949C61F24555A6AA530E29DA7CE54B2FF0BC826BA4BC77BE
5443BD05763696D001E0AE21F4783E748275A308F88B74C01A04893657574A7BB64A9300FB76CA2942F988EBA83958309771736E4F0DA3C64C0F9438
606F0356239204BD0122C64E1D959CDA156E260A51EEDBD13D3CD2683673B97A2BD18E1FE8916F7FE8F803A6CD3EEF37F3BA6AC5A213B0A8EDCE91E6
AB345100114839AE951BA1169D592EF48CB4E651EC167CB88DB76206688461039871253610DDF65503FC3C74FE98422058886480B545062C6464DA28
66557DDDE2ABAE3DB173FA41B7D7B5E7E01D0E7CE009E22C82E682B1F880E29F360AE351CD327C24E00F940C8CF64F1AEB809B7F83D79424C1C62575
EA391A4308EBE19124B51D5D28C451713C7423F1E682C161BF158A81FE897C2DCE16764E3B7E6477F8EF92AD7FC6565CB6DFF459D7C77163F333AF37
6A85F51B87BC281E7E73FEA6EEC25829B7B577B46BD3EAADCBFBE2AC11D0D057E819AED765541F2AE3961F5142A066516B73F58B80D735338001C1C6
A6E4E388A561A8297CCCBE0696622875146E0AE05CA804F8159A39D1BB332FC13ACE5479835F4DBD4ADFB6D52BE6CF7BF08107E63DF7FA7B9BFB4B36
8A33D0FADCA370ADCB8C8A9E5BAE424DCF70ED11CF95754521FD5C41432A93810DA3541C502F9125CEF69C51CD6E8B6D19CC5D5BD907F41CA554226B
8DE6F858F71095ED169549BF3EC271DF0F1C069E7EF27F0F7A4685D471BB1F3B77FFC3CF38FE2B4B407FD0C6B1AB5ECCBFF7FA5D8F0E544A260FB4C7
7387B580D009B9EF498E436928B084A074E50AA41C835F42E556F50A33D7011005A81BD1467F4B276C6CA48A7EAC95B45D3B69A4978A92B11441CF93
971E3EFBF033AEBD735991DB1EA3F451AEC6B19104C146867A292E68AE829F29BE8DC7906D0C7BC82DD60DB20A2EC349A89B96C1155177569669A4D4
F63CBDF2FAD638B3687212A9AC0E94ED245A74F3D68052E39AED636932984B0BC74C9DF2C0EEE5DF2E19FE51EF4F8F9A31FDD30B5838F0C6832FAE28
4DD8E52257B98DCBD70FD1931028DB0D639B857C60CBE6AA6D4DD8819B2F09E11A214C9D89500E7755B441B0266D0C1E0E7629AA97D2463C2635E856
AAB5CAAB15EAA1F193C1104A8F2F1362FB064B7B95AACFFBB7FBB859479101D378F5D1FE558BFFFDD86DBFBCEC8A2F5F78E6D91FBBF81BDFFFE6B7BF
F1BD4B2E3AF9F8634F3DFDF33F7F5BE122295C918809DDD2819689E9DF51DD4C7F2F247211088371E27BC01A422C37CA4C358F4D601825BE76EF8B3C
65A61693972A432A810621F524F406B99E57617125C7250FA869C778001A7E0A341F68C66F7C69EECC99D3A69E72E3B32FBEFED4D38FBFBB7CCBDF3F
7FF1259FBFEAAA4746BD8105AFED18DCF6FCDDB73CFEF4DA35BC1521BF7006DA4290220F52B41BDD9F35BD66C24F30B94B01512233DEE35DE409F0BD
29A505EAF84D27419166D6258A4BECFF9C2D8F5CFB992F5CFDD43BDBEBC835F4F73E9C8AF186A494FAF085B54F2D49208C5FB768FB15C5EDDF97AC94
18B888C62F55273CC79059A00C828B40CA1BD8C8E402602E125EF3CB851AE5A9FAD33DF4BB779EDCACD26CE95B91FCDAB4A95FF4774BFEEE8AC5BFF2
D7A3274F3BF9DE20DBFEF66B1BB60510A4A88FD49455AAC7345A8B286A5486AA4EAD4E05A2E9BA52850CFDBE62D2B651B43046360A4D2C0934DCE435
552409999B989E73BB54B304547B135CB114959F4AD3DCAB35B5B80C4EB64214C77C698F5561E5A7E0AC9DF9B5B155775DF3C9934FF8F8C7CE3BEFCB
5FFDAF1B9EDC323E34D497F79415349867EDD8B061FE199D33E6BCD032483595D24C9241363735EA7AD43A23C3607F66FCFD90699839E3A10E2A1E1B
2F13A8F963E16F146EDB9DB021D369C1390CF120D89119252DF87CD3FCEF4E8C15380407183C0AB3D00A70F74F52EC0595DBFDDE93375E70C49419FB
9DF9F91F7FE71B471FB1F7E1471F7CD899A71FBFD7ACA953A6CC983E73CE9927EC357D9FAFDCB7165711C159C07CB017632F8C5A0688D02A36330360
80D397F14B6456801D09A623CA0871E64E687CBF88720094CB8163D614FBCA5D71FFF53FF9E5B34BAD36119946777AB902DBDAD448ADFBA0516342C3
CF8CEB0C432C1B28A0A4AF288A837E86B92B0D335C071D178A015BC70C2B29F67C97D2717390A7481AD5F284ED0061517CB7877EE11B7EBB8672FA23
CF4AB9F543D38EEFD95DFD77C1ADBF5E78F61E532F7CC3CEC4AA7727EA8D916ACFCA7F2F98808933F7546D70EBF0D0A02503AE84135459128E8CD1A4
9EE1894CD4E880ABC686947918E959A3169CDB0D9F6600CCF8A8EE70C7E35620617E633599E501E906CC4A9A14C6E929F72B452FB0B8F1F26BF7C7DD
2F3D76C38F3F73EA71677DEAF26BEE5BD43D6803BF4B0D3476E9D49D8BD0713D1A9C575C79E8DC034F7A95BE0F37932DB686F4A3488065A81B8639B7
0F51DC18AA3454AD38335DB0B2FC5058DAF43C46CC8BAB34F4823665A62D9C6710316D43BF4C4B65C08B3E8F321DB92E204DD881308EEBBB0F2E6008
910E901E8D7259A3A7BBAFE969597FE4826953274DD973FAD4BD0EBBE0278FDCF2BD4B3E7FDA79271F73E40577DEBF4630611C088C2E8FE09EE099F6
A87D09A927916EB3E230E12915D43D69573CA3CF6BFB3438C5991C1A106D953F7400F442B02FB46BB5A77F78F9ADF39D3600D970FEDB69C074593813
9A118C5EA64ADA1784D0B89619D5B0140B0DCF0AB05BCDD2725F79ACAC23380E88050F36D00FC9FE65E32A6B35DEA852AB25AD66D9E70D9B12FDF093
9B43E53F7DD1F234C9DDF27412649BE7ECFBE46EE0FF2E18FDEBAF9A35F58017A923CEBC05DBF88E971E7FEEA9271FED4BD2DCB2377B1C9B296BE4DD
6543F5C1255D96EF36AA23F6F8A81FCB91F56354DD53C0679B9B4B716288EDA971BF2C8DB9E0EC070E3DE4F0A102E54DA88835EA5EAA0583AB070A54
565D5B8F95A8571DE1231D649E5DDE72CFB78FDC7FEEE95FFFCB3F9F5E30013F6FC47BA8FD10D14509887A62AB6C29ED2FF966C794E947DDDA43FD44
447D74AAE8EF021A88A3E6B8A3831A3782778E4D5FD52E53ADA5986856C28085CACA595015C040928267CFDC20923E43A652810C65AD2970BDD3E646
19D1CC6088B32EB52CE0CB1BFC2F25ABD02ED2E06ED8384A3311695FFAAEEDD6B9346EE234DB33D977FBD557FD6BFEFBCB7BCA9497A2E76F7CE0B52D
5BF2796DEA70C80A1B96BDB96C5BB51109C878A59402A18F44B599AA79E8875A589ECBA02C6633C164310F04318CC9525087041C1001AD6EA9DC4397
5EF6B7BFBC380C0F31B806008C6132612B650CF38E94DAA87B8136E542582C31CAA3F4C3520A4A2107D8C0EE159A1F91CCE7C60A0ACE5FC2E97A61DC
581D762F1A8004D0CA15600EAFEF626136FAFB95515C7866BDA70A2BEF5D1F47EFDFFC4EEAA5C927F7FC646177F8EF6A9D7FB8ED577BEFB1E7875EAC
2E7BF3D117E7DDF1CF179E7F6943AE1ED93DAF3DFFF0C26E2F061A37F22AC313EBE7BDB37DC57B796553016FD4AA856DFD16103D2828C21654A3B82B
8D046E2BE5D0042F977C1E26CA9D00DD17B5261498AEEDD10A4B8DDE3C7D1C7027C00AAAFCCA154FDC78C6E1FB1D7AE099BF5EB8A1CAD2B63A6E205C
1EF82C11DAB3DD30F0808BA1D01AB9EBABFBCFB9E0A679B9FFBD5D800868D082AD54DB9E0A2CC95C99253EA7C0EB1A518A81636BE9B21D69D7018E21
A56C87B563ACA984E372090EADD753D2BC612434CD5C4CB555C2E42FCE62D7E55226700A82072F7D046F5A68975B82D25BC038BA700E8A1F250FC1E8
8F4B23827119318ABD328EADBEF7AEFBC6377F71FD2D7FF9D7DDF7FDE3B66F1E356546E7C1879DF489AF7DEFCF772C156D7B20334B1897328A6F2619
7506AEEF709FF1A0218D4F420B0AC534ED233368D5D8B8F4E61F5DF6DFCF1584A4118BBE3B26A7A4E9E3AA2F80D403A41A72A806F407552053FF99AB
12D984E9890EEDD17240F9119B45B42219D6262AD285CDF93A3E39742776306A0B86ECC1A015A7C3AF15323D70DF5D5B824AD74094EB5B35C6D2C2C3
770D622EF841C75E0F87BB436AD7AAFD71FEDE933AA69F707D7FFFDD373CBA3697AB8F56EABEAF27BADE5BB699AA35F76502CDCDDA78D30D28127B37
6FCD87118B1CE6B7A036279A833949B5C4AED77DAFDEBB230FEB78FAC28107051F3049636595470A456843027D1605DCF6A8E0969AB1892B205EF4D0
EB577E789F534E3CFF9B772E18B1A081A302CD45D3D5560ED21E037D8CD583A059B7DCF1020F1C357471C7DC7D2F7B7F6C7C626BD7D2A50BDF7C75D1
7B1BB6966A121B3D6C1EE17215945DC807A6191B1189E0946528B2FD40F8DCF520AC23A9C1AD94AA9E4F1FD3C6E3A818727938B6A9B65C7F6A1044B8
C927DA76428A3D440410FB21A710A1FA9984F561E81A39258F33EA51201016FB1526A92B1154BD7DCE5D7C539AE0A9836E4D2CFFD58587ED3F73F6AC
C30F3BBC738F1973A6EF3969F29E7B4ED967FF732FBFE6C777BCF6CECB0B5E7DFB99792F2DDB3836DA3F1CB66DC92D6A24287DA2B337BFB194BA7F9F
F95286CEA6E71EBCF9E1252569DC37F2BD5EA41399B4E2CCA55F42E407E00B025EADC0A836FA5EE6E64AF37DB5E071CE146E34D2B765E205D47550BE
607009C8F06E045DFDC068C4D5D204A7B417ACEC6161925913DD23CC59FC4E31899C1C35269B776C09C3E213AF17E3F1D0BFB463FA678ABBCBFFAE75
F1AFBFF9F5BD3B8F7BD0831264DBF5318D6D51EA5FB6B64433A2DD181CCC95DC60C7AA575ECFB9AE767C158E766DCBF70E0C43B19A9E759797FAD7F5
79788CFCBAEDB8A57C9151D564B64A8D16A62199258215281A295FF82AC13840B5CD76DA7E5F72E5DF2E3E76D62167DDF0EFFE5280EA91D81EB329FC
79E2A9B8B8A94921502D39CD669800565F2B55190FBB3F3977BF39871C76D431A79EF4B18F9C72FCA91F39F9E88F1C77E6473FF50FDD8A0C2236C20E
40E1AE96826323281D19F12F8A59A140816F851E8FB1BB349C1CCA53BE59481A4B9FF61DC0B8F3500E30DEE43185072CC45B092F05943BCC6DD32F16
99DFA88B54C6D273A8F3A7E087EA7F7DA31B516B4EC9D1A3B0A7E8E5D8E4D14FAEE91D89AC426FFF48A958EAE9EA9AF7DDEFFEF9BE87DF5A3E32917F
E367A71CB5FFBE9DB3669FF1EDDFDFF2FCE2B5F77EE5DCAB7EF0DDAF7DFFEF5D9E14DA2E6CDEB466D3D67A93D7C024048E429556DE7DDFBCEDA5306D
855140834050080068C2AF907A2B287BA11730E08108CC006E976DE88D85CAA16E059C0CB89648D01AE8F78E2E2F8C6A5E2A8109A4498D267F95B1C1
AD0E849AF5509FA2E1C0DB36E227232F6DA5372708051BB39DE101A197BE1915B7A9F896A9B30E9E17EF8EAA5DA9F5B7DFF97F074D9DFBF34696CA4A
C92EE707C74647366DDCBCA3AFD80C7994EB2D366B858D5DA5FCF6D56FAECB35464B6ED396D6C840BD383C6805D4ECD35C49615E1D68A880A64C140A
6A2B2DAA564E435215071C3702029F6A16650FAD9937DAC7A1C2112B03AAED7BF8D2233AF63DEDD7CF8C3A34B74A2E78A5BA7E615E79C6BC965141B3
0219281AA83DA174C81D873B3677999C58F6C27DF7DFFBF80BEFAED8BCA96747A1501B9CE8BB62FA8C197FA0A22E65DB790F4FBDE1E5D098A1031608
8707814F5F3BA0F2DF0216CFF8E61A5F8B4C0CFAB1711E828019787242B5AD7E7062875A2ED4EFB03013137022930239004BC334D6F03344C742FF43
60961055CA3C8233CEB9DF70F0D1F04341EC51B649A1BD01E1AEB6A307300499D1CAF1EBE5A1752B5E59B265B477D59FBEFEA94F1D63BC0A279FFBDD
9F5EF7AB4B4E3BE4C80F1D7DE8C7CE3DEFE39F38EDD4CF5DF7BB7B9F79F68F7F7DBD4ADF30F441DFA77EC4D8A5606469512B83306EB581175832EA90
313FF0A11912F14ACD9394CE426A8E225968D24BA797146561581DAF860AC6C661B562C5324DBDA5AB253054435B9A611672EE3593C0DB3E1C2501BD
69C56D8534FF66CD5FB3D5F246EAD2FB48C7CCCB76FBFDED2A951F53BA5C73F5E1D3F7BB72539656DF7E69E5BAC1BEC1EE75EB4727C68ACD12B5AC4E
637480A75257CB651A3F1DC99C529522587174BDA1905EC875BDBF8C039261B14BA7490F0FF6CAD2684D1BBDDA9DEC34E641A49A8B30696EAB534B4A
ED6963C563579DBACF47BE73C73B459A3D23CF1FAFF811D7B6EF04F4F11C8E5275007110C7349B67C0148986E08E0B7D3ED462031D4CFEF7472ADF7E
E4F48E8E5F32254015FC5FB191C8AAE6EAF57AB99E2B957CC75C092839C48A3ADBA89518DC7CAC450C77331E2BA9253712DD946AEC00F8421D28E6E3
3E4E211254EAE03553E3CC542CA5591264D8ABB599F60691481F6F58F9A1D1D3E2D2CAD77DC7566DF19394225384001F32DFE50EBD294AA10E5376A2
D9A8FD63246AA4AF77E983FFF8D7BBAFFDFE2B9FF9F437BE75D979177CF6F2CBBFFD9DEB7E7ACBAD0F2C58F0C03FEEFAF36993264D9AFCCB42C3A697
2284A7BC86DCE9410609215C57D0F3501796B4F5809304DF46E29E8FBD615325C6B00418CAFA7000AF013C0E8D5CEFA8CEA26468481727226890E7DF
D918431AF1ED97CBB8EB8830766B1B9B5C57EC5AB5B77F602C5103EBABE57E916E5C1F259B3B671EBD6A37F07F17097FFA45A5F95B4FE8DCE7170349
A673CF3CDD173816147642FA0B259A4EA9C4B56583C8E68A28B20DC1940A9E640158BFCC484FDAD5B1A56B2CE0F5F138B146AED146CFE0AE8F9BB5F4
6CB786153E1B1BC3794D2AA6C0218A44F713DFF8D0CCC34FFAF16A8123BC602E6F7225537AC4C0945132A088F0BA87A973701B546963C1271A2DFAFC
BE82B2BB7B2BDC97018C2CFD801285275BAD78D5B7F79DBAE7D42FBFA6A8D8D353DF1C59B9F091BF5F77C5673E76F471C79FFBB14F9C71E6474E3DE1
F4D33E73F16F9F5AE7285F8688DB38CE00F2D1D56288B295505AF12536E7297D27A1E12312CB3AB5B80A6282F4A3B97E94C850B936F4B681B9C31D5D
9BDB1ABD52C3EFD73C80733075DF5AC1AE102905861B3BC57063AB2E2226B550825EBAAF0DEF5FE0948EBD5EABA5A84E376BB040423C8792FA18E5F9
AE1452FC9F9C5E9ABF64C61E93261DB54A7ACC8B654F35F26C8613092EAA2DA8894AFC7A0CF639315EE3D5219A214CCF11D64B4AD63C0D1D751E8632
C0D50F87CF840DBFDF53011FA3EF99BEC0963414E456764DE059197BE9B526B548F941DFDAB472554E95DFDFD0BB6D6C20EFD59AAAB96AE5986A6DFF
D13A96FC786AE76FF8EEF8DF55A23F2ADE7BCEACE9DFDC1EA7915BAFD445E0E4470787C69B63EBD70C8E3703EA9275647BA9561B5EAF51A8C07D278B
45BD319073B84BA5034F9CD36056958570BB047E5CFB680B8D4B8E7197A1672F148DBA8371522429B605519AB0550FFFF0F469530FF9C9F23A03BDA4
EF5D7AFE7C0654BB67E8A7DA77A2763FC168C62C9659140AD7EB99374A5136BCC9938DF16280F8748541D0E9A425DFBF6476C7DC832E9C1FC68D9E35
FF7EE48E5BCF3F72BFA30F3DF5ECCBFFEBB7B73CF7DEEA25CB5F7FF3D5571F7AECCEBF5C736EE73E4F273B4309F2E4D20928A5A1FE6778C9029735AF
02CFDEB6026FA6B1F1C72090182E40621CB3DCC0778DB827140695E1111A9E43046DB3282A0F83EB231B8DB60D28DA1428A152306394D09027A27C8A
EF050E8E52604361131B52C5F69A0CB40A2003E967347D85995100F6C18B1EBAFB9A13274F3AF89BBF5DD1648CFEB6392AA0DF8D6F62948C13CA7F0C
351DB7FC94661CACF2709C35BA1FAAE6A671D309F1CF7042F1393405A9D1F77D8B26B044D9E5AE891880C170E3B3457A9F92E1556F0DE8342ACD5F30
C2173D5F724465D50ACB6BC4EEC0981F672BEE1F4D5AFD3F7E5DB636ED3DEDEC1DBBFDFE7685F08766D6A6DF7EACB3E38817C2548F3FF27AD36D3ACD
CA48259FCBBBB5F56F6F2C5946B0C76B7ABED4B977724150AF7A116B58AE571D1C6EDA55C7868817846370D7D75E8F63983F896C04F6780D5272123E
7A1A187BB7A9C336A8941E8E60C75D5F9B39698F93AF79793BFD13639247E5E52E65128F0284377496A830712B00A6612EA69ADA36E0F0B356712B94
7FF23DEF57E16A15F8BE5001247613B996A27FDFC3BFFDCB7FFCF3A75F3CEBE473CFFEC297AFBFFE96BBDF5EB169A38BBADACA4C11448C6559EEBB1F
FDFDB2D7DF59BE7ECDBA3AF5DA5CD2544C3F2B4F8D3297427F1CEBB24D3301E4B1320A3B9ACD0D00DF089147462DC3B583D23883680F7602C80EC687
D7686C52A9562C30A285703A07038F5B01583D957105C43E754B418DD11BC41C6A119A23B538A1FE42C4694B8DB15437271A5C5382B01B9C536A62B0
0A66528840D180C0D6FFEB9CD9D3E61CFB5FABA12DA405BD398E8C9541231969511588566CF60D6DBB8FE1312434A0069407D692A0C60C3B03A82B8A
42CE87465B14E307874A59C4F30E9881DA2EAB91D7C0539E58B6A6E1B7226D2D784B66DBDF2CC591B360A5F2769447BA8A2ACBFA6E5E1CB6C66F58DD
CAEA6774CCBC5FEC8EFF5D62F28F477F757447C7F4CF0DA4CAA9BC74D3DB3D6B16ADD956F5213B0D659A4896D7AF9C90AC7F1414DDC4B3BCFACA97DE
1BEF5FDAE727DA711D279858B6C1D5B6E5D1E48A8D195FB88C85B8513B7556DBD24D6549C2798F060A0ED11BC500BFA39663C5F5671EB8E7D157BEBA
CD4F707AA7DE1AB497580570A1A60F532DE360237D049A21E4A16849A85E00FC17820DCF87166CDCF472176C6FE0EC672DFACB3973A7CEECECECD8EF
E84F5E72DD2D4FBEBD63AC59F25C3BE7F8F9DEB282AA10CDEDBE966C74D1538F5E7DD211A77DF8E80F9D76FEA7F7E9FC741E429A46111B973D6AEBC7
9BA889340CC4D41CAB4A456A67B446C9125C1FFAA87CE5FFB8725AF986B1D8660D45C6D1482B8A56C809622030CCC7C818713256B7649CA598F0CD5D
81DE19C01159BE2A4256C7F4928AA247E3BA43D55AF881513F1FDE08FD30B8885260271E5361ABF287738ED8FBE04B6F7B6F4C2AEE716C3A3C9A83B0
68841309661AC503BB6748B4FDC1522ADBE32E4D1AC65F4828689F6A63CA02D783A8D69F17615B938D452D9A87C2D0AED127F2B0399053C6CFC8D96E
51445353531F0CFF3FF6DE33CCD2B2CA1ABE3E9A4ADD4D374D92218A690414C3E8CC600045511C15510CF33A23A8A333E80C73A988A3A8E088206004
03064492C406BA699A6EBAE99CAB43E574AAEAD4C9E1C9E14E4F3CDF5EF76967DEEBBDBE303F0D75D056BAAA4E78EAB9F75E7BEFB5D7CAC5A89D86B2
B6AD2A4A9637DBA2A267E3BFADCF3AC1F70F753AF1157DC75D515974FCF86338FD69F5FE4BFA8F3BF33FD7083937317CA454698E6CDD38385E9C6C94
4CB8DAC492FB931B7698CEE87C020908B076A46332ACD804127EB17130F3C2BE6ACB31E1E69508915A3396571B1F9DA9B658E4945D09A27DAC02CB93
748752AE83D9C6E8CD6F3A7ED9CABF7FA209360C53AE1FC454C6C38D2ED62A60A9EA645006C7B610DC77A1409E61AB8D2A60885866F4BEA4ED523DDD
9A9EA77C4CE836716F7BCD92638E39B6A7FFB88BBEB471D0F02316B3A0DE68B5174CA31D70374C647DC11786A8FEF45F7FF4F0F57F7FF5A7EF5A7B70
DF4C59753A37BEFAF2E156797CDFDA67370DCF784E0CCBAB56CD63026D082A7878CCA1DFDF6EA5FEF639B88C4A363283CE9E76009329B799F6E842FA
24B4EF399456EB0D20762D72A4B4C2280500A7E603021058D03AA0A0F6208CE9F1A14D050184CF712423EDFE417125A0A286D041E0C2B693BEEE9BA8
737890CCBCB9EFB8935EF1DBBCBBFCAC95BA74FC02D380690D40DD0B4D3DA1B7FB12DD9515D02FF2A24E8EF51F8234BC45E502F49305BC84A08E4075
1A331D890B5E998B530969E56D9EE61AA82A240F08019933D0F79E6674595ED83EB37722729A9EBBE5CBFF3694749C3B9EEF745A572EED3F6563BCB8
F8F7877FFAB3E6431F59B57CE915B39DD479E8274F4C8F37AAD5BAE3515E37CA4393A31316B4AFD0C3477F4DD8335B364F0978C76B6640CC58B5C453
9EC4C15CCB8908957A9E82D61E9D5A7401320CCDC09C8F08E56A6D1FA8ED6489F3FCB7DEBAA2FF9CF7FC6C2882AA8596F8EF603B3FD2776A063D00CC
A05329E00986DA42E55AFE2BA203493919BB31E89599A6B6D5C47B997D6EDBAFDF72CC92638E3DE58DEFBA792B0F45A930519D3566B73C352104EFE4
984A441485C2241532D8B76E98A2C87FB39EE6FFF537BFBEF5B65FFDF007D7FFCB45AF79F5FBDF7ED9D7F6CF79DA19F77FCB6110D38BE2D99D1618BF
38F97A5B10C69A702E11C951657FB0E898AF3DC853BD2D8387360BA5F0D9B2C1D3276C9DC3538B3EA96163F697682DCE4465C2F722ADF4AB5D0BB4AB
172A27A925BD516B883876551EB1F4BE25BD2B4F38ED8AEF7CF7EEBBD78C5BC38707E7E6A7F68E14DA4DC7A4481BA0F9D0C195C95341354D68D36F11
866A84BCB84CFD820DB1D0448CB5D254B8220938A180587B0C0635D89CA5D6EC98854D81C41C994243342C6C2881739467D3E37E164F6DB2E9773778
382CBEB06BEB6343F5FBAFFA492D977BBEBB5AA8EAA7FB572CFF8C9F678B9B7F7FE06DBF2C5CF7FEFEE5677DF801CA38CDBDFBA627A647468BE5B97A
2300AD3C4ADADB7EB5663AC09A1C6C6BDA55A33C3A320D2A3A1DF19892881156A65A8C418D17AD6ED99ADEB5BF08ADA828C2B327294500E1B67DE437
018B5B3AB7EDA73FFD8A9E65E77FEDB9AA4347CA65AC55769DD9960CAC062C7FC154251041B95A24A25CA5BCAB188FA0DE918681432F31E751E6A7F8
801C5C990D08F547795CFAF1EB962E3D861EBD1FDDDD1A9FDC3F7A68786862A6CDB06DA40F7814B9568389D8374DAA026486B8834D564ABF946C8FDC
B52DA8592CC149B38BFBB73EF89BEFBCE9CC377C66CB9E2D5BC646E707B7EF981C1A3C60D0E9CF45C2589A404D10AB0550CFD3CD792023A5BD3C630A
0D7401FC50767D77BB6A1DB0D985D51F361CE9BB15708DEE1E26AE0774046D8D48838554571070FCC594813BF40A54D80B46553A41FFE644D6494D89
BDC6EA57DF76C1192F5ADABB64C931C7BDFCD233CF38FD6D6F3D67D949279C73E6B9AF7BE73BDFF6F6EBDCEE7420A7604485886B86F45B83D20FD40C
43052343F004425FF25639A0F785BD6BADCFC08C80305CF1D7AB5B1499B2847E9A429533B87AC3414961A9A3869E1AB6843D544AB3647C734D547DE3
C801B3B971779C361F5E5BA0C22FFDF18A81F387D3C5F4FF07DEF6CBD9E88D272D7DD9F70A94866A63AB9FA807F387462A8E5F29ECDAF2FCC80241F1
3C69154A8D46DB03DE566E3BA0928E070EE1E850E430D0EE48D793CAE3CA6E871161656E56CA216AC8D8AB364242AA68F5311908EC9D51A6F7777FF9
FCFE1597FD705F09679283A98AF698729BB5D94A33802F45AA37F4716F72CC1AB1F72AD28871CF73238B2B67BE199A6107848250393655C462F47B97
9D7FFC09AB962F59D277E117B6CE19F387B735C0FF87B526BADE74A33687C60AC5169DD47AB9BCBF9DC20B08A8068E415927691585C6427A222E656D
D723CFDD794ECF9225AF39F745FD67AC3865F9F2E3CE3969D5B20B76622A40E05CC65D513EC892A1D7C6E97F5C97514E877927D7021CA1E87A976939
842E7398EBFCAA3DF9E860D7061B7A204FDFC3BB9E3F78B790144581CDB18817430011B6DDA6A1F53964DA1114843A8CE550EB4FA43BBEF389877E78
D743DBD76C7A62DDEAA7EFFCE6359FFFE475D77FEB969BFE6A49CFCB374ECE9996EF83B3A8D90069D73B28EFEA13233C4BAD671CA56816409395DEA7
59A7371C810951D86BE9E224093D9689D8A835E14542186CE8E909FAE8469B8BC2D808C56885C9693CD7927261CF782A0A14323EB0F4F85F6BBEE7E2
31FB03C6FEE9E42DAF5BDE7FF6DA3431773DB6BA50991F1F695B5365BA4F29DDEE7DE4D9AAE355668C4832CB6DB4A6C6C61DBAB9256563DFF68557AA
70AF50EF7462BD1D1E05C35BE784F005564A61669DCAD64C857127A0F295523FF3E71608255098B8EF0357DE351265041D825827F6D0EBCAD4EA03C8
9A7508E56ACDAB9870294FB5DF0F38A8CAB27CDBF642186A4454796A7DB14885B547AE3EA77FD9A927AFEC5DF5EA7FBA6BB0E00647BBCE49EAB53DDF
94F41474AB33ADFE09141D8721BD634CD4A00044CF32F1D43C9ED371AB33DB1EFBF6A7DFF092A5033D3DFD279DF5CA335FFBAACB3FFAB5EBBFFBCBD5
D79DD6DFD7B7F4071D1C41F8752639567D62EECB300C184313A249F10FB5B8D61786050734B5E25C2BF8E989A098B13075A368164454D808D3D5EEDA
11A436195CFF22CC14352888EC8AD0FA4871E4D0658AE90F89020AF66614211C23EDEAF6E911601A47BFAF52726D221A45C6356FBBF4EDA7AC38E9DC
57BFFF639FFBA7AFFCEE87BF2D4C3CB1656162CF965D5B77EF397C70BE31DBD8BEBE16685E33BDD3C477AA15C09678A69CEA18E636584E9F169CC966
8D473E34068194B2C6E06020F3B47170A255D93A1F67CA0953A6DCD9416E4F1C91CA5A578CB25BFB7B3E5ED5C77F3100FCE13E9CE73E74FCAACBBE5F
B0CCE6E4130F3CDE02472771F73DBDBB156B5749C5CCB036576CB92C4898539B3E54F1E87E81D28C104C7ABEA8EE9FF6021539CD2617DEF881BA08A1
9295A0125610FB413D0B361C0E4052DA5F0CF4CC0DBA00F4F7040FE00DECD4AD068CAAD1DCA3C313382E65663D7BC2BD89BDBA58FF0312AA14BEAF38
FE89B5C3765E79E1CEDB3F7BC9C9CB579CB86CE9AA8BAF7F6E8EBB61BDE9A7DABB9E3256D9F38F0C251DEEEBD93D94FBB56E3E2CEE52E62639A7E396
B73FF7E2975FF69EF75FF2EE37BCEAC415C72C3966C9CA4B3F7EDDBFDFFEECE82805BF3A53F442C6EF2E1838F55377ADAB4116144BF4DAB35C19B612
C22DCD623C91522A859C698ECDFF08D38058EF104541ADE4279A5B1033EDBF9B31BF51F2351B4F6B90A077E034422004ADF58BF1641C1852BB005281
AEA913840824B6A762BD14288596F050011CD023AA8098AF04180A4A708A2ED2DB37EC1B130FFFF84B9FFAE0556FFFAB979D7FE56B97BFECE52BFEE2
4DE79D7AF685AF7BF599AF3CFFBAEF7FE4F39FFBE867AEBAE5E7CF1FB0A0BCC6FCBA839E4BE251DDD7499B0BB6DB8D6569EC992C8F22CC087C3AFDD5
9926EFF8F6EC542DB0933075AA0EC419659BB3D13535160D7D7D7394EC3961E905871725BFFFC0D37FF56B7F75FC1B27A1EF12FA54DBCF0DEE2D55E9
ACB7CB33D576E00AE8768067E2B89EE00E83EABB3F335F834D9FE2A1C543DEB44DAB5217A27A7077A16A417D06527552D4660D16042D5B48684E50FE
0B5B2EF4B7DD2E110DEEB6118E781252FE097C2A2B4A0D097E4A1210C2A73CA998E0BA918E6C86846DF30C7CFB445034D15FA30F10EEFFE67BCF39BE
A7E798DE9EDE638EBBE0862D619651F561B1D1CD4E967A53E39596E0B0B2CB336DAB23346B87556B6008C4F0DD80BC20018BF98F9F7DE629FDFDFD2F
FBC48DDFBCFD969FDFBB6B5C2650D1C29E3C8A70E36797AE7CD1453F1CD2BB91B203D2AEE563F6E04D3430E333EA5C0B883B5CA5AE85E31FCB6255CB
71C48953272090471E4048ACD5C7C183124690C3E54B690FE48820158B33AD38865183826611CA1E744C215C9448EDE6E304A8B7F28E02AB2F0DDB61
97B34BD10CD22158D6826702612D65D7552ABB94A6CC6935CC606EFD93AB9F7A7AFFE6D5EBC747773EB1EE99D1F9C347660F3EFEBD1FDFF0A3D53591
B0D939476568DBEA65E04ED4B223ADDC4ED8A05909524ED12EF23DE9362B6E8CBE2B43BDC43A59E1A11D3E3BB0F6C82183951E3914996B3EF66D5F85
17F52CFF195F3CFD7FD8C73F79FC152F3EF96D3FFBC5C307E6B0259678B5C3A35BB68CB4A92A66A19CDEF4F4B04D45374B44F550214865C219F658AB
7E583EB07BA45CAE4FCF8E951CCA618CEA77D3F794A342DE2A16BD4855A6EBA1E286516D2A6E9B8D8079F3A5E1D1AAE4B0DD92BE0A18A7544709DC6D
41668E8E4418604A4EA882E96979E4DB6DA999C3DA0F201683453AA71918E8F0EC24007FE8FEDBDE79EAC040EFF2DE638E5972CA25FFB111DC943872
7DC862561C9D83BBDC1CECB86336869DD7301441BB45F7B18A7DC31366A097E03BAA5898397C78B8AEEDB635571E184772574AD57EF8B2D356BDEDFB
B3749A6217F27E74F23AB1C5B2B06230165A2D9F4BA8F8A6A93F11E4B968EAACD93C6476697E7177C4103702F081A4A414CBBAD2A1DAB507F61F2D5F
4B0E694603648705617DCD1CD4D3124262CD52A47D3E11182822866815526D13360CA5DDFF12CFC126448375529160C0A7C027468D93622B02D22212
5CC0ACAB2398EA392A94CEF069F5ECB5B57B4B0B2B59843E22BD7F94283DEBE848D9DCBD5F0195D84557B08450496CF9896C4E1FF11B22F30F6F7EC1
CC8287EF7B786316EDDCD3DA75E37BAFAF0A75535FFFA5B5C523F6871E00A66E7AEB4B4F5F75C2B215C7BDE3FAC747A9F8031BC4B3C346158D2063FF
F34F6CD8DB0E1C6E8C8FD65CCF96F0DEC636BD5F2A508D2EE8C64F8CD980F02974EA83A03E5F6C8F3EFFEC42E83A3E1CBACC4AB9ED34460E1DA95161
91BB85BD4DAF5DF702A18216654C28DCE6847E25F66251FF2BBD5393221D822EA349745DE44F479AF02DF342AE53A635FDF017DE7ADAB2E503037DBD
3DC7F6BDEC033FDDE1C3DFD76C37E00EA02D2BA08B9D2111C78A2AD81453733A3A942721F62FED402A67A1297CCB0955C2C0A3EB2A8B045E182A3BE4
CC86D52FE19C83DF7CC7ABAFFAED21AD0C1609D37739588BD8DB53CAF138F79B0B2E175DBDFDF62C8F326F1E995B3A0EFC08A0AA7D748927E96412EE
3D097A8794BEE104864A4429BBE6E0D4A12F87F5DC34277890275A35113514A1FFC0D38AA438FB2A93A1ED76BF3F094219B8719C2BD345015F34E3B4
933066D36FD2A3D70AA14D8267899860A6C343693916B413E0FA1B8934320C8EB8CB252F4E3BD0FF8A34F397B0850AB13448F1D4B6EB731EC50A36BE
65CC4A7868FBE1FEC767C2E77FB3EED90DEB5BE51D4FD5D28EFDF44FD73F329F4EFC6CE3FDABEF7D62D6E78FF7F59CBC79B1F5FF877EFE3BB9B5F777
DFBDF6E213972FE9EBE9EB7DD5F319181E04BE6FBA3F03FB2613336314F43DA6B85F2B57CD804E044BF4666A06CF6EF0F1A69F9A025D8F19AD568562
44125161AD8DE12077E331BA39E30CDE14D868352BCD6623A89B6D53322F169579B3EBE90B6A1ACE44AC75A8D082EAFEAD2E7BBBFFD5FAE1B6A362E3
A677BEF5DC538F59B2A4A76FC5292FFBBB9BEFB9E7B15183CEB86597AA746ABB5D74ADAE03ECAAF1B30C530A51DD7F47CF0C6EBF0A745CC51CCAD14C
ABF00BADCDE10484BC293F2371DB070B2AB879D5D2B3BFD93C7A2367514067A773745720D57D3EFC5C0077DF44DB20C5997009414B4AE0A1004B17E7
0FAE9C22ED805F1F69C7802CE984471A5DD74D3ADC4CEAFE0865EB1C9EC729F618BB3E24D84CF420280C8F62B01F11C532F8FAD1558D0B4704C54E5C
2970FD2342550EC535F082033FC4F7285F7632ADF61199872C6EB86E80BF3F4A4BC844DBD592DF4065710A6B048AE361DB92CAAC495419ED62B90E31
D4C49DDF3B361F85EECCEE7DB6181AE2936B36D73D6FF2F0FEB2D7892BEBD6EC7258E6CC8D1DDAEAF8B2D1983AABE7B8DBE462E3EF0F3F0020EDF1A9
87BEF19E57AEEC5DD2FFAACF7DF2EA6F1C8EF27C742EA7BB09FC71CACA2EE67154B15BAE33B8762872AA2C4C74698C5C913627EB5423FBE3A31537E3
38229D5C09AFD2A29C1819769AB14C4F8CD0D2EAA8A6C7AB45A339D3B04224A39A0F7DEB44F3D615F603F5922C183A899EAAC32C38D2FD3A8D0112B8
7957BF70F947DEF1D6B77DF473B7FC64C7540DE63E5114389ED570E68B14676C1310413B6368274DBDCCC770CA250875510A8B7111616199C7B153D7
2E3F949913AA6E2C93A00358C14950F222164C0E0EFED39263961EBBF4AC4B3FF9EDBBF71C9A9833027BF70B13DBB7CCECDF3D599D9CDA3B3469D9A0
F36B830C686953BE4764CC82B6541ED53258E901979E20879EEA75A7F104F43BC6581B3A23141E0C9B8A1A4C08CB0591042DBFA3193B99C20620A8D2
842DB44C175D931CDB7C2A26E810678908AA265A0C84112880F298119868371CF89E49D7429483D70A867C304F72165CD3E6BA034A57B8350F2781D1
05380685547171746452B72A6231B9A3E079A6A07BC37B615B1BE6885165EC70DDEAA4A1513D3C346254C6E68BB391CAA2D69E117AB3E9D63507A5CA
BDD961CF9D98F1C2E942FC604FEFBB5A9DC5F4FF473101D43B35F6D80BB7BEE594638EED5BDA3F70EEAF1CFA0234B93934B859281786C7DA82F2B257
7B61DB6CB9C5EA539375976E2FC94DC34F42DF4F82593F87C10634F1933CABEE1975A5DF2AD74D97437BCB3C304E3762307EA864D64B6DCFB45AE004
D2BDE83391FA15073437A86EA7BABC4D94F602A683051C2191C8B04D0F2E0C30791433CB8A8F9EA5D8B73CE971375098A72B41483968CED7781CD2FB
83F6059DB0908A11CAE95C6126AF1DB59580F280EAFA64275A2790224688091CC19B2809E60D088298E53D1B1FFEF94DB77CE51B5FB8E98B5FBDE197
F7FCCBA597BCF9CA4B5EF5576F3EEDE4D3CE38FEC4E5A79DFAB257BEF57D37DD3F383151D36AA621F69EA11DEE50D4746BBE82F25F8495461C862830
5D01E4009F81943952AFD932A91588D3D825341298F0214BE91D064D57A691DE1F06232FCF94F6438CB330945A51411204C0B27094619D073B81549F
09ED5EA8B45049AAB9BF8101BB6411840CCF9CE9B504D998F2097854B9F04C420AA68BF8C44B479A69CA0B55FC12E80E683CBD1B225F53870BB38EEC
60F850B5DCD41ADB39372B93DC2EBE30582224547AF2293BCEE491B55591559E2C453BEF7517CEEB7BC9B38BADFF3F9A2220D7520F517DEB572E59B1
746059EFC0495F7A9E4AC090D3F90C093C8A8543FBD7FE6CC3A415A4717BB6997582C9A1A9B1F905CB6E2DB41B66A5EE84C5A2CD589A852DC798DB35
5AAF9408543B46696AC1F58A07F68F8E6EDF3D6B316BB6D4A4B88199BFF42A2E0CC0BC46D5F74C4AC974871272D5B298DADA3AC54947FC813697D286
3C20D9713A4F146762CC1FE1A34D9996A0440C5F4A8A19A108409B978DB28F3E010CAE0127A6E6B03F40953B178974F1F30D5F69C671A429EE1427A0
8D8F044B2021A79A5E29194416033F2ED10E5FD0E08852BEFACB77EE2CACBE7BC3E35FFDCC9B2F78F5E597BFF884979C7FC9C57FBDF2F4734FBE7033
948D34EEA0534B85451804350CE7126CCDC4A979E4975FBEE2A2F35F79DE270F61BE47D583EDE01D767784414D883B39666D9891A8C6E0361FFD7F19
31863FE94C8F4C77F9F7211844E01748D88E01DBA471E8C2A51C44642D4CA4CD41D020809AB9537582487A3693DAC500E75F280531638A8E2C602A53
2EA3022EE2162314D36EA6F01E89D4C29AC134975179DFD62AC5749567CCAB3B06AF148D4C75A2C2333B07C3B8E36EBC7B0B4BE3855FDF5748B3A16B
1F0C567F613478536FCF55C162F6FF6342015AA13B8F9D83BFF9C0AB56F61F3BB0F4D287589E4B21F8EC58AB16E65934B37FFDA136976912798A3BA1
6CD5A70B4D2952CF65E0B9273653EDA9C90937750F6C99B6712F0AE55861E21507870B6DD729CCB4A5EC24988EC350DB7AFE4090C83CE1AED0AE572A
ACD4B03D4FF567C3811E7FD380D125580651E04A6CA342485F69242FF0DF38092C82DDD010E254A760A1185E4219FAEA5489FB28BCB563969CAC81D1
47E98FAB84B7E15B34E6A6BA328844221C2F4AF87FF910664175CAE4AE810D03E953D6F4B50BB69E22680503742270C502B3A9C223236D43F2C6C607
7EF886975EF5C0E30FEF192E57CC858546D26D5974BB99B1D8F29F9FBAEA92BF5C3A30B07CE5ABDE75C78CC62EE845EAEEBF7EF722049D10B126CDBB
C8C6AE85D0FA54F305A57DCA93BA91E7ADA119139A7B40F67AC7005C0BE8FE6924A1421DAEE06DD0FD2C94CDB5F2107D0278FB2698AFA67CC1D53127
42079360173AFE2C009C209492C936D351C29A2F9BA9A270506321580C5950DA3BED47DCB3E006AC0EEDA88432CFAD271F2C6589F3F8B79FB39364FD
7B6F3C74D3D58564436FDF19C38BC5FF1F510F4083803CD34D3A61ECFBD9075ED4BB62F9D937B729DBC8D6C6D537DF69E83991B46CD312094B02B752
E1AC61517E1669EAAAFA64C99541B33A65E96D93989BF333E3F5BA634907D3713AF04DC7B4A9D6A7EAD898B715A7845F812127665040B5A8F07D47C2
F12A26340F5E2F6750B0085CACBDB84C6ACC0EFE1B7C46091450010F678F8C37D1BECBB4AB17768622ED92C5AD60DF10D47CA19783B3286540053815
D3DC41DE6F4BCC1DEC90CEB85376E3D40E91BA41154E42933265FBC8E4741B15504CE93684043E3C8C99CDDD90600B26ED70FCA3C259C2F933A5CBD6
2A9746D7DDFAD58FFCE3073F76F179AFBEF63BFFF9C0FD77DEF5D51BBFF2EF1FBDE2A357FC455F4FDFC08B5FF9AEEB9F390C73C25C190D99D3F38630
FBE0CD30EECAEDA4512EE2DFEFCB6435138027725A52F3951216669D6DFF32DE55EF540A362CD0FC8630B71F690962B72C75DA7779AADD8F32A80D40
AC2BD6CE9E6823123699DBE7C5E02572745794EE568896AD18E62DB9318FEE6DA426B7DB143D6412D6FC2C15B3F3D2DBF6D81E022CFE811915676A6A
5AFB22A99D8FCD7792C377FC940A027EF7853FE533B7CC25D1BBFAFA7FBC48FBFDA3EC03E47A0090D9EBFEEEACE5C71DF7A22F8DC1F551D67E7CD37E
09E55CC9E61F5CE7A25AB6DAA59989A2030B304A26E5F56BA76776166228CDC441A966D9D5C3F7FDE4AE9D053A8D9287CC177ED8A89A1CAC92037BEB
04BDE10298459126BC283AA2DAFBBB5525A88F1A556A3639F25B63470D82419421E169C729CD2106A01788CD36FA29CA4E0405E04FAFAB6C8C1BB034
6FAA4619552FCCB4A556F4D3B2245A3D3BEAC2E034A2CC1F6B3F9F54F12CD17A3D89E6B80AC5DD2A7D31991BD61C3ECCC3115F986D5238F2100C9236
B4BD59C0E1B643E7AD7B1513B33C3ABC6FF5FB5E75C53B2FBAE283EF3CFFECE38FEFEF1F78C5A59FFBCED35BC6E799AE23B2B8C36A2D9C762CE24223
34D2DBC11043A02B98454D074F05B0C1BA9B058C4E7FC68F4CE6040AB00F4D153F5D58BF8DD582242E4F845DD581D46E6AB1E2A850D79209887DDA65
3C39CA3AC2FC43854D1B972065263A2488C05057923AA265AD350BF87DC8F96D25AA0D64E2CEF9792CE27653D427CA2A6772E8B7834162EC9BA37291
7E21A2495F9E5EB3D58B44F3F6BF7B2ACE46F646D9637DFD1FB017D3FF1F5B0F20EF3A697585C0726BD70FAF3CE9B8FE151F78065DEDA43669BB36A1
EE5096F78CCED68B3CA7FBBF52985B181A1EA9D3F1150E9BDDB8B7EE34F76CD93DB46358736DC5E8DE832D6BF260450AF8D2B42BE3737612F92643B7
8AEEEC5481B69E6B701B4BCA51496BCE979463A9E245154D2080224EAD25851F689E1E5A811C7A1AE0A960A72D221C81B64182D5FB721D133E68056A
69712D3A9628DB4544511884E94922B4F9503F7703808AD4EFCF86EE3F6253075EBAD0C5D3F203737B243A6A60DBE10D868EE59854441094187FBC21
E908313A24022205B9763BCCB4D0597CF88EFB7EF0EDAF5F7976FFCA97BCFF9377EF680954170475A89A019F39816EBF8844B9C2858099A765B63180
CB2229B3A43CD5263CC33C6C0428293CF07AA1BC9132141E988F761153B1A23F41B25001D941D1A508581452C448A5A5BA6E21D8BBD2FE6BC8F794EC
5D37F04035C250D233438251E014A4A90A03E514EB3C6A14B52D70D870E17E108763631867E4CA42EC561E6F1E18B2B85FAB609A93EACF94B64B4E14
5AFB7F7CEB5496EEBEB994375ED473CAAEC5E3FFC70900FEEB5F08C437375F73526FEFC06B7FDE464088EBA1F9C2ADBBA9BE669387EE7E769C616BB4
5C3EB069F5E100A4B49813362E3FF7C0FA5AEA5081D03408CB5B75CFDCFBBB239E68D7E66B32F142871B2E657734A0E2D0B5F4FE2C6654DA7C37D528
5571E978862F304AD72E7550C6D00E4150B202F60700E8C410C803EB2839EA034E38192B877448F344EBE3496458CD9047AFD05AE05D4D9EAEE375DA
F5F5D65D77ED778BD32BE964231BA2470EA5F0548BE621CFA529B37DFC756C0C5715D0B7AC96A80EA05044F14746BFD71BCEF8D8F6C79EBCE7C637AE
EC5DB2E4D881F33FB37A5203832C95F01F84CA16211F3AAE7441A50A4B453AAEF0FCAD549A947F33AF4ED554E28B0E1C76441408BDF5AF381DFF00EB
0339FAFEE85706A64EECF0188E355F318D9A35CE1D85730E0E3014CBE9F384A057683622AE4412B68C50681D0504C4906A215F538E658B82C7E0FDFB
807870E558238822CF6BB3A01652644BE4F86CA6189A020E6CCF04BD87D09AD833D830EBE387F6D6C6265FD8B46642A4F5876E2C64D995FD277E452C
9EFF3FC6F3FFBFFDD6343DD4DEFA95972EEDEB3FF973EB19546A9B877E7778FAE975D3812A3C77E7966D07F616CA4110F9AE23CBEB1F1E2E515AF2CD
20298F4D39D6F043BFDCE961CF258F0A8546BB55DCB6B9A6285944961380ABA24215B4EB734D2A2533552FD478AC5BFC31322B1D2FCB37C7A6458293
A67334A8F55A73338769A56609609CAF197BB172A06C89C55D3059234414A18D7065A2193A740E22DD54D7A140EB904647DB7D52D38FF0CAE0D708DD
664704407B4C29C7F74A369719FA8ABCD6D6CBB9ACC931E1234C114B4E950DFC6FF20EE795B10DBFB9E1967F79D719CB96F42E39F698932FBA77FFC8
7470342A50C4E391B29B542C09DF2240E5BA7EDCC1AE3FF69B59A56263273F528DA6F6F2CA5145C1D15486D674D11C7E70538BDE2FD090A6343A0BCD
569B6204A040DA75EBD4AC09FA749C9E0EDDCF5406B870214891B03D820B0106065023C405A5C815D8BEF2AD0072E5D6F622C504DF870AA01421AF16
663C26FCA983F5109F9C8FED9E525E9AD88663268C79BE91C8CAE3F73CB46FC15B7FFBB3F307363EBC7DD0F593997537FF6A2ACD1EEE5F71D1DC62F1
FF273113C87336F2D34B4F5E32D0F78ACFEFE118009BC3BFB875F508EBA4CDA1E70FEC9E356C42C08117B3F2F0CC13BF79A15469998DC94AB5305E5F
187EE4C97DF589C167360C1951C6333EB5ADAE7CAE9571905D19BC374D03D89C558B1E256CCAEE16782EDA6327119549A18F21ACB9E1C59345BCDD2C
CC53569C1B71080C63753987097D12D95EAAE7179936EF8C39B87174AC44B7BF7D547C23030D58E844881303D65C06550DCD0A04E14875AB700A0B5E
3BC0D43098AFF91E50083AEA19EA03ADF319019068AF3D9CB7D1E7F73EF6C54FBDFBC2979C76FC09CB96F72C39E698634E7BF535B77DFB9CF7FF6ECA
3E7A25A16014F8461BCA06201D06CC774D06998376A1E814C767DA3ED410250BFE9B32C723C3228C648D8E97F7DF75F7A863292C0CA781972593CD14
4EBE0855689668A9F1A40B7DA0DF2BA5660963962AE12B121C2E0642B4E75D5C883CD2A2C44C04DCA3DA45E9E5E3C2027A05D078C48BF8AE1552FD91
448147BF2ABFDA9C1B6A6117EAC8B32516F8A33BA70941D807F6FA69C69FFBC5AEE2AE75134DBAA07B1FBBFDD11D74358DB3969DFD48BA787AFE18CF
FBFF09DAB436505078E493E7F70FF4ADBCF8E179BA7342AAC877EEFBF99D5B6C828FC1DCBA5FFCF6EE9DF3AD5A94C6D59D9BF76DDABF7DD46373535B
468E38CD62B3353D366AF0148CFB549941085BD014EA5F219760A0638135E73E08AC0404125E34A36E9F0A1416D707B735A03021224400FC70E1084B
84D9A2F21F2338DDD7C237471874EBCD393A13618569BC003FA1B05886418156D582DB67084C11132CC6BAA366D0282DD31B7575B530878F98D53402
C5B4710FB036CCB48FB29199CB7052B1ADAB2493100F906EED9E0F5EFAD6775EF2914FDFF0AD5BBF71EB1D4FBE30592490FFC059BD275DF0F12F3FF8
C29E11AA95A95481BF96529940274F508DC29B758FB5A767DB86EFC6E8C4C39B238BCDF2D0F65F7FE7EA7FB8F2EACBDF73CD10DC3B6452DED48ADD48
403AC4ADD09B01778087582A02328AB4C8B0DEF0D1A2A5A847846EF669CE5496F8F05CB74DCD1D42D0C346015498145D8634EE245CE4A9F6359089AA
CCF86040818AE442C531F54C8AA4E0464CDF3718A58DCD2F0C5B4150DA7690E751E99907370D3D7400DCAAF233F71F2848618BEC96DE659F30F36CF1
38FDE974069539F4B3AB5F7AC2C009E77EE9DEAAAF18739D6D0F3D53730D5B78D5836343E5FD4F7FE7D67BCB028AF15E7161FCC14727EBB37BA727DA
B136FD8CA39601576AAD959534276B50C194EE96DD5A1B40787EE02A4D7FCF9240A512444074AE312248B8E7841A8AC71AD9432D37EFFA5D26AD8A19
44DD8E3DD27312FB14472021C0AD0806F79A31183374EB6242DF3815510229223AFFDA162BEDFAF9E9BD386450D40A680E0ACEB1FD07F57CE87B6BBF
6FD41C51C216E61C4CFCE8AB84E443ECDD43225B47CAA83B46CD30B59090D2997FF2FE1F7DEB6357BEE982BFBDF3E1BB7EBAF640952A00CD01047F01
3336E607BCEB1808BC638C3E72F7377EF6BBCFBEFEF56FFCEC97BF77EBBDCF3DF8F92F5F73E383549F30ADE0996554D3A07399627F4861EA9730C2FA
F47912BB49F9D977E034A85585BA8E03D80AEE6AA874B7AA0063227A49A67548BB6A43E895789E9626CBB5E6985D419B237117E6F58A21FDEA0B758B
6741A5B46DAF54F5435354FF4B6F623088E3CAD31BE70BC30D8A35F673DB16DC9CCFD16F67EAC465AFD8BBC8FCFB53E30710E0BBFB8A937B96F5BCEC
E6991A1A4AA908CA8FFCC7DDF74E41832F15E6C48E9167D6AD2524E89B53C34DE541FD8679127A72D0D386922C651E6BB6519BA953592B84B37A17C6
78A9DF7262EEF911465861ACFD7508C65A0653AE2B19E37EA884E798011C82BDAAADF753E1689159F375A6D535B42F9FF2E0C4AD4709DAB63B95CA82
5F01217B74EB8101B0D2133812391F143A9C09FC07FD778DF123DD1FC0A9C12A0234BFB1E887C68326CA8142CF187A04296FCCD8C2670C1E3ED0164C
C0D9492558493C0C0216C2F2BB8B81F3FAC803DFF9DC67FFE386F7BCE72BDFFFC60DBF9967FF7559A1444AEFBBBEEDAE6B3F7FE5652F3F63C5396FFF
D1F6D50FEF382274CF209A6391F7E4AF872886265AEF2307B7103B46D860CA29AC89D9DD134D6CFF0706BDA79A05219154EF496939114FB730D15D15
47271C999C9EB05CA67C2F2C2C68C6B50C837645BB7DA08ED21C6C7A78B5051FFEA672BE592BB459262BFB264B9E8CBDF18544B98A1D5C5F4C44E1E9
AD4E96843EE7DBFFFD57228F4BB73E1E760AE7F40EDC2E17CFFF9F2043483536FCEBDF9ED0DBF797E7FEE32F9B98DBC7DED0BE27361CA1C42A70AFC7
957DF7FCF2E6BBEF3B3C39E1F39A105E1CCDECD9B340E7B81E2469EECECF342AA58207E91F1531AE9538D0EA0B34732DF283B6A01B549A8D9A1F29C7
1171E0B3002EE39C19CD56E0CB203466EC0E3609BB8AF2FA3C6879AC98E2061C7633ED24E4E1D0D3D9B28AF3D874A3323DD68E77910A02D789BB8A3C
60CBA4747CE1C89144DC17A00D6AC18DA4BB6D8440843D22A00F99A6CAD5FA1E10C774EC5ACB1642728EED1E9952E647A751730FE09CCDA0160E5A7E
7674FF35676124B6DEF2AFFFFCDEB79CF7960F7FE3478F3EB5F6C9350F3FFAB3DBBF71E3173E7BF1EB5FFF9A0BFEEEF3F73E3D38D212FF059B25256A
25FD483CF383EDC558B306C2B66FBACD00C95A7FFC4874A43B33CBF4CAA45E0ED6DB535402840E4743405B2BD2C7C530B55D867A67E2D7C32C8FDDBA
CB5A5C3B1A5B65DB0C7427C56FE23360B790B52A6E0486A1BD65CC6750289ADE5A4D65C6E78FD43D142AEEE1D93475F70FA94C0E6E66D6969B7EEAA5
F1B3573D11E5ADF3FB96BEAFB6D8FCFBD38C0139AB6DBBE52F7B96F40DBCECEFAEFED4534CEFA7361CBBF2F8C6C9A79FAB11F60CC68E6C9F39F4EC90
0D527D9C30A3D808D5F496613310F642D9E399EE57B3A0BB2103E44BD937F0E32C4A734EA9B95DB5ABC5E9B61386A16586224CE988853616712242D4
2973189ADE59E82A2867E9C4CC39F43C1A0C52427A40E759545E7034C7EC0656002951535DAA9BFC9C5E0B5338DD2A87BE00D448F47A31D65F45C42D
B7FBB60823F0508BF5D25FBB7EAA327B810A132645E8FA2D23D0D4026BD6661EBD3CBD07E46BBA1C4CC684FCDD16240630FB00050107265181CB99E3
B49A836B9F7CE4EBEF78E9E9279EBC62F9B1C7E8C7EBB6B842571960E7B0343ACA23F2211D94A57CFCDA0B2EF9C8AD777CEB6D679FFCB2179DB3EAF4
D3DFF39EF77DF9BEF543664761D29974BDD451C320733342EDD26FF3AECA470CAD15D00D22B7CD40994AB441AAEFC77AF699C656C9F0A82889F23C58
6863B7104A639E1164D8C42A6E9A4B04D50FF6C1ED7692C9FAD3DBE605A12ACFDE3444D768F0F956163C7ED3E6DAE167664532FFC82DFB08047DBE6F
E0B4FD8BC5FF9F6C7F90B27CF1EEF79F7FCA406F5F4FFFCB3EBD46661926D1F6ECD3D77E7DC4B11A5E9484A25CF6A467CD97164CA1D5AC8285E9F9E9
724D8FAED2D8B1ABF335CF82D68F76EAA52422A03811E7316F4C8EB5F4660F2568C3924EDDC78610A77308AC000B8EA883A6BFDD927AD386703C23DC
2F83F2B40506017278AA38565928A264E882418554A9D9612C30836108484C557D4C353EDE4F840E58A4333F7D2F1B1DF27106A03B8217955D5B9E90
3E06A33341B1219401C7C98E83667BB61C4A0F9EDF1012140E84CF74D6D59A5C32D1BA1B426F2A30CE84E3B861204367F277EF3D6D697F6FCF31036F
BCF9F1A76FB9F8F4733FF6E51BEFBCF9B6EB6FFFF66FF67285031A345D5FF96EB535BCE781FBBEFF37FDCB7AFB7A7B8FED1D78E55FFED54BCF3BEDEC
55A79EB5A267E56BD775B0C583D2258B43FF68151432D07D2822C4A03B49A5671671B74892D86CD4FA8A69DA5D6B128D6682614CC0C69FB0A0D1064D
13AE5795B3E4C8A63655516965CFAE393BB59DD1353BE929C5811F3DBE795C66FEC66F3F343AB16FC74479E38293355F78784276583671FCCAE5DF4B
16D3FF9FF04490FE64CD9DB75D76F6C0406F7FCFB9DF1C89CBB340A798E2CD1E2ECE557D0BAC73B3B46BFDBA49C7879064D6C995EB9B41AC194679C6
B9EF378BFB8BC840614007C61B9BF5C17A4F6CD743F9AFF755E8EE9C1D378474B4621F1CEA092AC49EE09DD86BD51838F00924F103B4F4BC8617504C
C16660ECBBD0E1449F1BBBC3A94EF5513BC4185153FE745AA7DB1CE29D9986FAE80082158C6C2DC0B8A7EF05B310AC7C680B7B364504680F32F4CD21
8C3BFECC50CDF6230EE920C20D3003C3382497C576E4431D19BB85042FA007EE070ADB7A3EF7DB7BEEFDEC5B4F3D7EE9AA57BCEB43D77EF3910994FF
8DD1E15DABD71E78FEB19F3DBD61D3EF0E58C6DE7B3EFB8E575DF8EE1BBFF5D50FBCE355E71E7FE2B2A5FDFDFDA75E78F5676EFEC5F67DC3BBC7F60C
ED9AF5EE3E7179DFABDB394BBC799780BA53B64124207C24300101472268B254B804B6A0A904BE6357565C7B81620E80EF8B8C16E732B5E76C51DACB
503C65BCE1D0AF84D0CFE4AF9F17194BF8AED5C36627982BB9153B4A7964AC7FE44831E1F3EBD6EE5C981EACB303234E279D7864D4A6AA224DDE3130
7089BB78FCFFF45140648FFCE4A32F5E756C5FDFC005179CF6E15FB6B33C611EA50E47FAA1274527932633B93DB17ECB024CE5649485F5B1998AED06
1E7A7311DD607BE77C11DA300592C591A2D36C868966DED3010CDBF3064F941BA421C8F549180AA7D260E86CD1CF4BAFE523E5F2C075DB26438E85B3
10E1EB36DDFF49F5908F5DFE38D37E57F151F70DCCF7BA937E708034FB8FDE97D00370BD02008B024C02A12F18517120EC724DEB8F27CA29D5AB50EA
C8B43346E886B585008C40CEB8A88CF859EC37EC1C43B2A8BA6B21E65150AB430711CF0E5521F40B9FBCF56BB7FDFB87CFEC5B7AF1D757AFD93FAB57
73B1E7E4B24857F39ABC947AFF74D185270EF42D5DBAA2B7EF1858981D7FD65FBFF1E2377FECE7636DE5B0C037FD86E5D187A87DACF784F3BFDC4EC3
8C152936514D42293B55A1ED449A01819A1DE47CA9659430ED8F34BB214ED1AA803A311DE14A1CD97E42459B5103B550A085C9F71E66E04E0523EB86
283E24FECE8DCD0E812543D2CF04CA39B2BBCAD2D02CCE4EAAC0192A95C6DC4E1EDEFBE3D9345A28CAFCE1BE8173C6178FFF9F050AC8636BE4F16B2F
3C8EEA80A5037D67FEAF3515F052A97A0CD3CAC88EC3D37ECB756C1158AD62B1120A48F9C57E6568AABE30D694945B29053B8546DDA79A5866DA44C3
9882204D1431BA3999F09CC242918A6118D8C785CD3356AB16A4A88D0343A6AC69716CB728BF7268C44689A0B5831241E9384AF8A4916A9D2FE6855A
FD074371406534C364E02BA552BDF1AF120A569A397C941DABBBEB9194A83F54E8378C2E455028DF73859F49C797917D70CB90056B2E3C59C8CDCDCF
984A34AA5E18261EE34E5B482F704A15D37362AA9D3BAA65DBDB1FFEE2852B962E59B264F99B7FB8CBD5D3CB089A2A8CC2E54C99098EC97BAA0891C8
6F9D79C64B2F7CFB072FBBE835AF7FCB55D7FE64ED7899BE0F1ABE32B09A95FAE41CD48CA7BEF792BEF73E57ED74C2A6254061429D12854AD8163A8F
5CC73EECFC687101894D81083D496C5BA661C8A39839F585430BE055D0C509A94CA3CBE410640AA8FEA1CBE0B542819961EC95661901A9A84D4FCCD2
A4B5E7C074CCC062622D677AD71064DDF2C1DBEF9C13C98E6F16F2C6A9BD4B7FD559ECFDFF393CD0D9CE335979E1979FBCF0CC9503CBFA969D7CF937
D78EC79D0E9A72CDDF7D75A754E5390F3DF0B85D9FAB5349ACB5FD734A96755B6AC63F2A75CA2D814DE71E8DA9BCA3E89B146545ECDA09C1BAFCFC24
6AD5459705483F2FB24E0476BC6C3743E5B6EBAE50B271A488BB9A10702EC07D4757C25311A77B3DF07D6D6609FD7CC622C1A85CAF99883390144169
4071A36265DD093F5834582B84507808C840EF0B2A7E2C9159E8352AA1A4274063015BB4CAA95B8117A699046F2682DA7F98A554A104B0EDC9A09D69
6FBAE6DD577FEC952BE8EC2FE979DB8F471D4C2A8036A480B201218B6A5B2F32C201904A28E9954A0DEC3670DBF58F1A77E2812D1E9F5BA62958B4EB
4B27AC7CFB0375701D7920ECD87323BFD66A43F4974BB88A48D836E4E88622F92316506DC3A49001BD4C63D28B15211B576F394A4215116625543571
AA0B9C9AA4EAAC5AF4B5CE50505B7062AA015AA30615114938BE77968A27C78BF3F684D39CB382B0525A38B8668867EDB5DFDE17C51FEC5D7AF522FA
FF3342018801096B8D3CF1B5CBCF5AD9D337B07CC545370F79FA8B560B7B2BBEE3829092B1D2DEC32D09A84BE7C6B5594847D30D389D67D549AAC3B3
2C167AB32DCA990CCCA20165EB3435662A7E65F3160B28DCB3E7E6161CB73154D1A3C88007B34766ADF29C1B51B8A90E2E60FEE7794C41C196321FF3
1C3DE5771A3507B2F9119751C81023A8486E3AE6CC7C10692D0F6C0624B68B2DC0EE22302668015331733139836A08F764247C6B66DE0F93140E9D1C
1AC8996ACC5056E4000F9D4E1CB716389D5A0E5E0FCFB270626CFB7D777CEB6F562C3D6EC5F2E5A75CFCD51F3CB4D5567A7B51697931CC3FD0ACC046
6366D729CCC4A951F544486F08DD47A1DDC4A81267A54A97102D3C93C9FA2F3EB4E2C4773C5249B9AD98110129C172D4AC5618FC5561B21412A8E9E4
CE8C9D8A900E76462027687950344597B35646D513A08B910661E03B2248E4CCF602BD979CB5AC8825F6C65D2EB414DA9305D3483C8253EB87D24E1A
0CADBE7B9B195B16CF13733E28EF9FE1EEBEBB9E9A68071D3EF1AB4766D2ECCEDEDEF36717B3FF9FD344F0E820BE13DBA38F7EFE5D2F5E35B06CE9F1
E77D71BDD49DF5C049C2844D0C05D0AC754D163ADC4307995223DDBB8E11807B2BE9EEB4B1834679BFB95093CD96A35AC5195379E5D6FCE44C61F7DA
03C26F94CA63F34DCFB02AD59900E219D373B2AB7EAD557EB1F91AA3A3AD2C1F8CBD34761BACBBCF0B9D11A07A699A61106A4DAC0838DB0FBB4A5829
AA0302218EA13497184724E3D530C3EC003BB6AE6585C2AEB69CD67095CEAD02974886E5824BC93748F4CC214F45B31A2E4CD994109504FFA6FAC055
E79DFFF2B7BCE9E4BEFEDE65A75FF2B1277CED0F1C4819FA7015A1CBA32C7A8B09B8068458DA3310F68D9BE5D05DA8CF9543C658A0B02110E741736A
8EE24B794F030B810F5CD47FE667768499D56E5B8EDF16CCCAE8CC3B58FA099D28A082825350750CD9DEB2D12584A2FD8E9CBA53F608F3837CDD9122
07CD27E18A0277A0529149595CBDD58873D5DCBE6E4E26CE63EB591EA7D5B58F4EDB49C01C736A78B8C154F3C8F6E105A8082771BAEF31BB315AE551
E5A96D2D95C787BF7BFB06BA0E13ABFA4E787451F0FBCF28FFE7DDBE95568BA73BDA1CDEF09DF79EB3AAA77FE0ED3FD83BE482EBCB5564EEA9435F06
63B6380B9974EB0EE53688F4A4A159A9D6CB0B0C8E35791E18A57D4FEDABF82C859756E497AB5E4A15B68A5B1B9E3A6C195008869FA0084465789A47
DACB82106BD3913EF8EE0A8E1DB3CF0ED3BFB0D86B38188CD19DDED5FB4EA4E3C0134069C552C8E0455117FE4BAE75B0842DBA7FA72B80003C7CED33
202DA7EDB0A9CDAD280E155C4A1437AA4EBB64453965D644C4964D98A53C61BB818D75A0A8B4EF8EFFB8EEE295CB579D7EC6DF5E78EACA377EF2F143
ED36C1065F523D235CF4E961F141559083155C6C26D39BB7608C0E453FAAEE196CD4296EF9B53A81854CBBF9E5BC6629E3B1F7F52D397F14F18F8B50
98518EA5A77D86DE5D10705D05928F2125C6B14C45950668D1503FA09AC39DAF5359447580560CC30AB26071D6C14B192D3AFD3234A6671C1E4F3C09
3B978535CF15B46F19F78D48E58A357944BF6ACF1059FBB1078B539341144D3EB0DD8FB2F017FF7CF76C18CA674F5EDAF371B678FCFFDC6601795739
A41BF863737AEDBF5D747C4FCFB265AFBCEAA92A48ACA9F015DD8C953D8E6699A924001F3F21BC49A5BFC73CFC7F9574F86C39CB9A7B4658D2981CB7
08E16B3E8C33BE6B786E7ABCAEB28897470E4F07316F1EDA336C441DCCE7C0AFF30E8D57E6F4D0109BF0E1C894C724CBE2F6E894AF99814AD3FC757F
4FFB08811E0C1EBDB6E5D1271E7AFA94B523688AE25912DDA4A0E301A0AC22E589C89E174AAFCC8B84CD1D9CF622DE893A5EA3EA0489516B51C52D30
6193A5E76FFBEEA7CF1B5876DC09279F72E249CB57BEE6B67D1257270CF4C6BD60AEE9E9D7A0D8926A3D90A8EB3D9ED9300CCBF5DA0E3486A0E501A1
00741FB0054CA7993ECA63C72E593270DAB5BFD9716476AE36D73043EFC86CC77DF7F507760F4F948ACD6AABD53C305C234425C02DC6F644E6B7BD29
96A321C85AB386A41ADED53BCF0C6BCF507AB2F73640FBA577A2E8ECE62C60AD300503C8CF53D52E1C19AF517DC0E3F6EC8C3137BC69FDEA2959DF31
68B76622696E7C6A42C69971CFBDD52C9F9D7FB0AF6FC5EB671797FEFFDC2A80DF4780DFF7033AB9A86CFCFEE5AFE8EDEF5FBEF2A3F7D575A73DA6AC
B16DBA5C0BF55C6D76D617CEC4DE31CBC31987AE56E8478DC3AD8EDE4569CD17862A4CE8623571EA33C55A9253F519DAE35B371DACDA8174C2A89353
BE2B3EB3C38DEC06E3A11D5141A1DA630DBD0D042F9B54D92653F2A8F716FE0EFF0F2EA2313303788A6BCD20C86325BACB972874E321B54FF83F41EF
2E3455A602FA114CCB19532C9041E073C51D58637662DE82827994664AB743F9E4CF3FFBEE375CFEB54B57FCC5E9AB4E38F1157FF34F776D2D53C013
895B6BF99464CDA1098E7D273A9E9E5D9A67501DF1EA3CA32B4031B1C4E8C4866E2C1AE5B64F95BFEBE8D925F67C5CDB12A8EC95F1838B4F3EF6D89E
DEFE8133CE7AF505EF7AE765EFFECBF33FF5F6254B7A7BFEE2F4D7BDFE439FFAD2A72E7ED5F9D797095B119E12B114A1AA5660A49662121B321E8793
254C41159889B1CCA3CC5A3D9661DB013123CB79687888C77168A50422FCA95DCF0FB7E1283A39198AD89B797ED798F0DC5678D08CDD0D8F6E23A0E1
6DFC6D21ED640BC5E4C3C70E9CBA6351F2EFCF350CFCD7FFEBF604BDE2C3FF7AEEAAFE81A567FFCBA6A9A6C0343A2C6F5B3B152661AB54AE79AED532
6CBFD16ECB249CDA332BF1F56AB94DD57998249DACA3D2DC3EB07DB445F766A6BDB3234770E955EA9CC31D3852A912D6CC42108E1DE149EAFB4E2CDD
E14DD0F8F50D53A529B6F6B0AF9328DF07B3B7EB08443548CA27472DA87C482D04A2E60E33906EB9D6D9157444308DCF3316A5868D1161D695197486
2CFA81C01578E584995387A7246167D880D3C33CFCDB1B2E3B6DF93F6F7AFE17EF3D9ED2FF059FFAF5E034DAA0894824AB369C44C8B65DA86B2733C5
2597CCB5A157DC30A0578A5981D09D46C2E766035A4AD2C22E0FBD56B5251863A04B4878178CFDE44BD75DFD8FEF7EDB79AF3871D9B255E79FBEE2EC
979FF5E67FBAF6ABD7DD7CFB93A3B3379E74F24BEF8BBACA685004D2C48262238FAAE5109ECC49FDB095A76EFBD91B6FDBEF5308537973068A8309DC
9A33C59BED300726822C28E6A24A33865579C7A8EA748F76DCD83174A0A9D491A73658BCF9E40DFFEB7EB7930DDDB52DB55FD473DC9D8BC4BFC5C751
1C90C76CE1B9EF7EECEC15BDCB57BDE6EFD79434BBCFA3346A6D7AA645B8DE0FD3841B164BF9CC9A474705167F5B56BDD6F4FC04478AE02F9BD9B167
72D68D33513D3265C47EBB6A3B741E63E7E00B9B8F8C5B56DD89E14CEC83DFEB19747209D2CA40592D07ACBA546FF646516058D8F7D33AD920F3A7CA
C7A49DF0BF125E28F8C47E0A1DB16FB8B1F61F80F761482578A556ABDB29B6D933ECFDA6C6DA56040721FA6A248C4A7DE410C1663DB1F447BEF7C157
AEEA3DB6F784175F7ED9B9C72C5972DE5B6F988E32EDCFAB12C18B831332C2D001EBCB10DEC19E12817FAE8DCDB133045D12B885B9DD5E265622208E
422F15A791E1A38311090A3CBE6F4549779D28E1EE9167B715D4E4C1567BA14211229461F8C487CF5879C2E9EFFBD4E5FFF0F16B3FF9C56BBEFAEDAD
6383231D554D3B3EE711638E94BE1DF3B0D198F8C18B5FF4BE07FD44299EE754A055A190DA1C2F383956A9CA07065B9A544CEF8DE04E7B70EF8482BC
317DA6C6BAA71E1C0EB270C7EE1237EEFBCA35BF98CCF3D2B7BEB15D45EFEDE9FBB770F1F82F3EFEBB3748E09715B77DFF53670C2C5D7EFAA59F7FE8
9085B140223DCAE491DA7DD36F8CAE6217D5DF66C1B082246B4FCE0ACAFDC2982F45BABBEF2FB8C6D4EE83874AD5DA5461BA69B9CDE6C4FE51DBB3E6
1B1677B1B9933BC5DDFB3D4AEB1E7677E187A17504C1E203025030D5C0569EE761A93FC30E20BA0210C809E17CC75C4DF6A57A033581E0C26DD67D25
3DAA93054508CB69799CB0435856D0DDD30AC15230ED9592399B7E79FD27DF70DA71279E7AFA29272D3F1E33FEB3AFB96FAA42E5894F2FC3922CF7DA
263C0DC0C2072338EEFA07A4601B200868A5416D4C0C2DF12C4B8F7A9D411905C9173EDC703CC7BA7E0C974F57C92C8296A082022ACC40A0EEE5F89E
17B09B4F3CEE82379DFBD7577DE883EFFDC87BDEF3DAF3CEFECBD34F39F313B7395D4D02FAF83C94016697F40E6EEAED5D76479CF08A95E6C19E1D41
968989A9081A03993CBCE9B083A2A75A0B282036C7E7094F550101E23DDF59375ACCF2F2035B841A5FB3A9196559F8DC779E5759F289BEBE772D7AFD
2E3EFE8FDE00FD99F2E29A6BCEEF5FF27FF5AEBAE4AA9FEF965A4A2B8F33EFA9673CAA3AE7A6BC986E3E0BDA5F0163228D647D76B250A450E105BE61
3466860F4079D0A93111CE8D0F0E0E6943FB58187EA36A296F76FFE4A1612FE1C25360A84916A3DF4D776F00AB10CD73472E8B245C40D052F0AD44CF
FA65CC6D41C03FE4426FC410D2E701F37CA6A28C33400B4E6532B767676A5C0A3A3A52F98D86C702C2FC146C529EE6F7AFE8A35A7CE5AA138F5BDABB
EA452F7EF9077F3C628056CFE867315E8776314F44D58DB09D07971E6D3812EB717E24A1A100760E1D4DD0F78E2A9A4BBD17190ADFD71E685AAD980B
302659A8B58A09B02FCC387EDB6A33680C849C6A09E1F9A15001D5280121F50ED5174E73E18E171DD3D7D7B7E2FCBFFAE79B9F7EA15887E60A678957
DB72CF75D77CF4CD2F7FFBA6D4E36DC712AD7D15BCBAE6257658DBB23C9904EED8862D8D246EEC1F875F60F5B15F1C8AB3C3F7BC20284AD576ED606C
7A6B5D44F681F59B9EF3D224BBAEA7FF75871693FFE2E3FF613E400FD9187EF4537F7B7CCFB17D2BFFE1E6DB7E7EC0D41458AACF03736AFB81DDD373
A5E9A12D230B66DD2C56F76E189549C72B1C189F29971D9E28F0019D20515EB3DD706CEC06A5A26DFAB6397EA010DA632F94924C26514EF9358BFC56
CD744BC559132BBE14587219C23997107520BAA4FFDCA9AAEE46402A2C989D07F4AF6160406A1B8B089176C703C7C774D2AE133154C444E04B7F7AEF
4861868040CE0B1E2CBA8DCDF75EF7F1775FF9C12B3F7CFDDD87A70B25BD44C462BF45E79AD703C2F8301AC556334A6CF897E34CEB5D3CCADABE561B
876B01876A31257A6B966965C218D608810F2493F965382B439B108DC83809CA2C9CDA31AFE2D0B7D2D0EBE4F4B940739481D016C19A4325A16A280E
FEF4CEDB3FFD6F979F31D0DFD37FC289A7BDEBA35FBFEAEFDEF30F97BCF335AFBFE0F2279C591B33968E6C69F933A99DBE3BDEE6DDAE627118B6DBA5
5692B71FFFD1263B61B5E12345C20662A2DC893339B4CB8FCCB5DFDB1DB477DEF1602DA0D221BBA177E0CC471715FF161FFF6FB500FD1999871EFEE2
DB4E5ED9D77BECB2D3DFFAB1DBB71950B122742B9ACD05BF3EF3CCEF76D96169E3F63DC3061D3B6B7E74AE688BD80E2CB339B7E5D923A17F68DB649B
8E276B1CD9B9FAFE2354B0B60C46777CE82C54A6DA32C7962FB7DBCC0E199D66BFD1A2741786213834F0D7365DA9B93A69D81210BE030F08233B1904
3C6C950A0B13339E9E08428D17767B76D9C432905E0D60501A47DACEA12B4207A6E165142CF0E14408EF0CC8FBA181A804E5D31C9B85B266C3BD8033
FA170D3860BE15DAEA2800201C52B18566E683E727B1389C2A4BD0A7C0B442D1F1F55CACEE32C23D82D15B45884893B05A309DB6012AB3346DA78CC5
FF50B8AD660026BF90F04B8656B8C75C079E25F40EC391AD0F3CF89F17FDCD159FBEF10BDFBA73EDE6F139C723B41588A6A2DFCDF4AEB0AB184438A3
33F993C3719E0681D576B027C8B63D3A11B1F2DEC122DCCA72291342250737D5556BDD4F361BCDFBBEB62732174A89F7C19EBE9537B0C51B7DF1F1FF
D111D4FF1B1476DEF5E9D79E75E2AA15C71FF737975DFE897B4C68D2E9AF61D6EE9BF30DB3343459681838E97EA532B3EB81C7272002E8D5AB4DAAE5
65D870786C97A62B33934DE471CA77F56AC53703472652C42A439A4C126955C23CB64C46505E6FFDD159005B20850F76DC550244CAA572177A7881E3
55EA2283E40FCCC258E88AC0935893073740626CA8D7E5A15088853905AF8F24A757833862877006961E62E597CA3CCB7308FE675A4103233CAD298C
798356DDD323497A87419B20470E8D7FA59D4D74BF3E3E5AFF27B046D61C5DAAFD43A952CD6E3066EB961F4A68A2E23D311139611AF9327229EEF108
6CA280400F048D157D2CAA02D4EF7F03704381736AA2DDD66468D9EDEAAEDD4DCE860709D8D4DB759FCAAAE621973041242A0D093B9746BD19A55179
B820E93289E98DF7EDF0736FF76EBAF823CF8DFA0BBFFDE70D713864C483AFED59BAE2BD8DC5BB7CF1F1FFD30DD03820E39503BFBDF1C3179E7EE2B2
653D7DE7BDEFF24F5CFBE8509369C68AE358ADC9C2BE3D9BF68F1B96DF2ECF1D1E99D857F0E9DE265C4FA75688238FFEF2A09F261DDE1E7EE100C7F3
B1405174E1CC94B93282C85731E664698A951DDB85A69FD6FACF310C030F398D20DCE5699E0DF080DE23D63262946433C1A8CE9661C0E8BCC1F19222
0733308F87A03E188219A4F49C30E098B1E3948759045BE07078DCB16B8ECDDCC0F76D066C90E651573D10D335787C28280FA8587449C83104379344
A0490091CE1CDA83984B300E365F945022D79EC67A6A48DF1FD69C1811C36301F323A75012AAEE2ACFF6BC80D74686EB36F345C04205F53FE55A2C90
6927931E360014F6939C463B0C85A4B8EA06AD203506F7CC35BD9089C6F0CE83FB4CC26106260B1D3E3E6272651FDE5FE7843704260F311B5FFBF496
2777787CEF8641DB94FBB7D50E1F58B7A932322E58ED8ADE65CBCFD9BBD8FA5F7CFCCFA602DA6634ACEF7DECFA0FBFEE2F7A7AFBFAFB962D3DE9A557
DCB1714F154C3DCE8D7AD51A7EE44B575FFDAD3D4DF8CC74724627B3B87BFD73E30BC3E33BF7ECAC53E990A571B359374DCF345B296E5B3A9AE69C61
D989EE9F6B860BC875BAB5AE37E2F5F6610733814470A99D30E91FCA913A53D3A9846A50BB1D255EC3E05986934B68810288D7ACB46D199854E01322
501ECC3F4D4B124E3600D8A9840E3BA9355D6E1850FDF26C9B85CC74CC506B6F278281930B614E7DBEA1450437C218AA3C084B5A6A84F0449A746D4A
09C20B0EB652D7AA105B3AB599293F818A699641B55708CAFAB13BBAA7E23AE87B0654211853B396D12CD393A3C5E0D8BEA2AA810A95DA8293691492
D447F6CF42083DCE19BE45D5270D282AA542869E21BCA0F0CC2ED88E74C285719F7BAA3C6A41CD2C495DA7D9E64ECD866FE3C18DEBEB5E58DFF09317
1E585FB102B72C9B6DF10F032B4FF879B478FC171FFFE35AA0CB138C596DF0B79F7DE549CB067A962E3B6E6069EF8957FC667315DC55C8E5BB8D9985
9A5924703F553159A558D8BF7964C17353E9D4C2C6FC1CA300604F8CD5449A969E9BD242FE9CFBC20E3C907B9193A5A45CCDBBBA7A10F1CFB55950E4
FB4A85814CB4EF702E0C0B8DB52810BE50B15B6BB128F6DAAE5E218C82968B6C1D346A66CB682D9403886607C2E3018C45E8E984CC73C17D934E1B6A
7D900B21B2414FAC9C9A51F7E14BC05CB06FB984081FA6FFACCB3DD2E2FC5022C3C49D1049D4D50C4A633819766B8004365B141EAAF3D56AD7BB37C1
17A9CA575AB738E00CD4018A2B6EC0286084816F574DFAF8D2975C3858534EA0259828FABEB1DB575B42F851A65C1E4351C5854F732E390BED488AEA
ECE123364CCAFCC9712B92BE0D91050ABE7E79C60EB0B09053D859D83962C6CEE0DDBF7D74EF01CF0FD950D937B33D2F597EFC35F6E2D6DFE2E37F3E
15CC8F22011C94C6D0A3FFF9F76F39FB9C93FB7AFAFBFB8F3DE515575CB7CD135A0333751A85E93D87F66CF9E94F1F7E6E7351D08F24540BD021F1E7
4AB530D6818470FA4C8109DB26DC4F593BEE7AF804AEE31A742A54AA08B647C2F75CCF0EB01863FB7427A3052F3D37947EDBCFB17F400823094DBAA7
334DCE81D370A6AAA3051F223D9A8E43593FD3467CB1A07CCAF11DFAF35082073F1872E140F8DA472CD16A3E8EE7639240A1C1280E977D0CFF38ECFC
3C2CF8A6A01943DF33D6FE4611CC4652E9B16EFD6F38DA5913BC016C2EA3048820E549CF11B56679AEBB05395C85226EB7B8542115115CF1E2B61221
A14EEE1E1C1CAD96666B7A564A71240B0ECFA17CA26F9998CD3A84873C9F82AC364282136B75DA11A954AE5B2D956D7F74C35C1886787D77A1A93FED
C2CE4376D06A3B4938FFD8EF9E9FABC5627EE7CCD3FBE95ADD73C2D2E5574E2F1EFFC5C7FF781ED87519EA328472BDDDA68CD11DF77DF18D2F5A726C
6F6FCF40EFD96F3FEF4D9FFEF653FB0C4CC82565C8A1F583ADA9A611B44A6EA25CC9E8000473D37B86CA22834E207D1BA71C97A468EB73AA9A59DB34
7DDE76A9F28D3C834E9EF23C4531C175CB94204364E3B0D59C5F30A1FE2B1242F2C8C3CD19C3ACC00D177D3EA803853E644851F66BBD4C1ED2398B09
CF4B16343D46B88212269CFA0859F05061EF9F90BFD49B8D5A0894B72C1E734ACF75CFF76D4E67FCF01363919CDB5A802C21F889890A42095E00BD0A
F610BC2052604735EC5423100231702E47098EB92058C0516D5AC2B21B7B8A4AB88C3567A12B40DF4AE57FD8C44A7227364B0E73860E1A89C06A739A
2C0CF348868123FDA1ED210508B57D17CBA0E0AE122EDDA6571EAE3359A97A49A55A712AB7DDDD62296F8E1F79FE609845D6AEE707D7FD7CA35D1F31
02F7E0534F559B52F20D4FD40F5E574DE3CFF6F6F6BFF6C8D116EEE263F1F13F01009DDFFB0C1F4502F893B0BA31FAF82F6FF8C8EB5EBC62D9B13D4B
7A7A06969FF9A1D55E8ED113843D8CC6D0B607D634B323BF7553AD0414D7C6B64D408A57E1C0814D90A1A3DFAEDA81513631C7E7B68BD91A6A79DFF0
A4D3A8B56B96C8945F6D796DB07A20D24BA83C8434519233B35D6E10A848252C3115BC77B54E07403AA48AD0BFA35221F167E74D83C93CAA87477B1A
AAABA54DE0DDA398A3A93FF01E4E422F909A5D4FC50A547F54FBE9E758CCCA75E07D680DE6B1E344D2B342050543CD4F92911E50E45A781805837637
46FF008188C769E0A7DD76664CB5FBFFDDDE7BC04956956DE2BF9D8ED5DDD39398214A16541050111041D105410151C0F84912454140415194110109
921964089318661826E7E9E999E99C73AAAEAAAE9CE3CDE9DC736EDE73EEA0EEEFF70FBBDFEEA722DE47E986AEAA5B55DDF5C6F3BECF433A092A7019
8D75A0A8648A91EC37EB36511625A34248D2A1037B26A12400C48B5CA48CC37E65CBE33316E1E9C29751902AF1A4C3A9564466AE67FF98C60D964C9D
EF58B76AF7E696CEF0D0B677A779904F61F7AB96DA8ECC16A08846577543634DBBA1DC8073B6B3DBA0B7F4EFE13FE500DE5F18B2FF722E407A022413
3071421E1B6C79E9F64F9FB97449436363F379B77FF5BA3B5649AE8215999435F5B6160D01599089FC95A6C8E4C8DB3D3F2702761A9BCAE7CA0008C0
B05432D79E62651EAA0A53C9958B946C01EC4C70ACCE6625928D2BB8920792040B09E25170B5AC2BE4C4AE14C9E705320B20030870C54006762C2209
C4CA843518E94A9E82242F37952C3021CD928A81ECD293E304E2257497718BB0029B44CD48231BFE64FE0897F596CA21778A409645E46A76F202845C
093F9D829317A2867654CF4473330E96C3A50B24B38AE4C1E514994F726B0C77AE51E125229DAA234575C5C2016025BE9C13A0AE39A4D948EA7F271B
E20D990D17897E8242B3B84E9139A998A490AA18522242D40E4831A422A64C695C6E6220C9F36401A16FF74C4E2FEF59373693E32C23DB3992D795B1
D600473162E0D5E7FCBAB972B96A6EACAE6D58B216A74CDEDAAF87FFEBD3C1F7FB82845548E3934756DF72F999272D5CD8BCC037FFEBCB973FFCF4D8
D109130000536663715613F00334435614FC59B79898645BB01CCCD148614A2C39ADD72C0D186A3697CCE400B474A2998973025A5258D5266CE206E0
751C16C54C4A255CC0065034CBB041304ECCCDA06281A4462606DC1E3CBE5A717832AA1CA509C781DE22B7315945A1684E266205599C4EEBA66D6107
81BD012B1B4769CC4DB269006449217D3AC3E518279AC42AF140C49F212232209193419DC809BB52A6E4BCCF321D952A08BA715424853C50A1C9AA22
224BC1F8CDA832768106E143C1050F990F809222413E38CB9BAE4212E12FC06F994A8AAACB2826804A9AD725ECD344DD9D0750B9EE15EF150C4B8CE5
75BD50644849427EFBA4D3A097E62AB6C695CBD889E28C6462CD7ECED0C32BF7E75BC74BBB7EF956D1927FF56DDE30AFAFAB3BEE69E6A8C891E7013C
FC571C0F5A6E5F9E34090C3131BAFDCD7BAF3DE3D83A5FADAF6ED9B9B72EDF78C06FE04C1FAABCCA4885EC5C9E74DC25B2435FD9D3A558A6854CBAC0
0B495C1213A56BCBD1190929841947E2A1E1D8BA980C944442BE09132996C79641CEDB20A10891B9685296294898AE49333E3B331042D84F1854816C
276BCCDC5830C201A4AA9C4CF687B123B00D4478B6150898489226645B847358C4A158CC9454C2C1035D8D02B2436BBADB3C64AA87D095918D44C39D
2724CD491CAA8932A12A5302B90DBB0A4D25FC80926C918121851C101A4490837094C9D8DDA9A6CA0B0A74341B673D8A4A367B0DB747A964C6F286EE
FE2E05A091AE9C454E50394D2F8E264A154E504D0008078BA6302A6068840B1C7AB03D56E688D401D9CCC095181D95A430329502AE954C9137D4CEB7
C7784BE97EB373E8405FD17F604CD727AEF926675A7B71F67F6F917025D8DEE2BF87FFC20EE1FBE703AEAA2E37B9EB4F3FFEC2A74E5C34BFB1B1D1B7
E8A7C3D1B4883FAA3866E677EE6735D3AE740E95354B881C1E660D852693331A8EBF2EE9A70DC2649A181A5CA542058762647C07E7050A728C4C3050
210BC247353B4531974C5640AE6CBA0AA03691F02D09B66DD83015CAD338803332D9D683BACC55CADC51A2601CC75DDD4C0BAA1C20BD34A29E47B8B5
549C4063AB96B0590A1CA9FF09EF296904E846BE4F24E29C2EFBBF49B67EC94FDDC9201CF46577740969CAD8044ED40139B937219919206D7A327B80
CB0E55B10C528C480E210C562A127984EE727A41DA1F4A97643A33D8359CC8A9A429A7D29A43E61D94784454089BA95E2810212484C88493E5EAA80A
6501620FA595630C1944D407BA290559A220D96C994586D87EB0626973BB3A621A050C1003A0F2CCA54FABB6DDD358D77047CAE54BF0ECDFC37F516B
F07FAE088E320A593A147293EFBD70C7674E5EDC30DFD7D474FCE757CF12067F1BD18A8AE3FFDC708A960C331F4C4B82A498AE8687C62B203822EA94
68C9D09D7E83B9991251C6B164C5B20883104D0EEBB9D181588566B294E8EA71482225BBDC2116B02CC8230B22A84815AEC4E3A8A9E908875828103A
605C819B50540945B82BBC49B66A912010BD630D8A1AF96F0597F6467A5CD11456207783845E440ED04474C02415074EF56189C1E1DC6D126AE4B001
1169540472ACCB2C4AD800540C85D420D8BB10693140A6F51143ABD8EB30B9A929C124A789D88E55808A7101F2B2CACC054B90AC1005BADEEB552C68
5B06D97A02F88A622AC9139E0396032A2093CC0E5F019648885262DBDB8A8665173B6624FC1E2A15BBE20F538641F50C6B7A6AEBF6842DE23A838A8D
A7B26F3C3E8C7DCBE6FAFAFA6F844DD3B63DFBF7F0F72D0A70568A9D407A62D36F6EBFECF4DADA5A9FEFC2FB5AA22A9118C5C683EFC5C16CC920EA23
166072A5BC84A3BAA69703926350F94CA4773CA5408B706710AD3FA4B87CFB0E245BB74AB12028228F0C19C74F41E68850B6AD17C3AC2BAA47B27177
3D8EA25544C6720154C8793CF6040CABCA22C297933477C21FD7F742364B93A4DE282A381DB7482281B36845151981A886110A72DDC1C980264BAE66
1F9915803C9114C686EFEEF6EAA4DC20AA9E36AEDD11CFAA22472A0195CB97881E0719FC4196027928231B5B335018D1DD5042B1095AD7DCB500B2C1
487E6F82800A03C3219E78439647744E312C796EBCC89BA8321D2A3324E350B552B02F0CE5826D727DFD51A89BF2DCF631D386D104928BDDDBDA2B06
2895E8ECF49101A8F3D0665AD70559501AE670E2F098AF69D98D73A4F5E799BF87BF7F77D0FD170B4985D94D0F5C767C7555554DF3E7EEE9286A6452
08C7208749CE4D65A02CF202950F87FC45402CC8049224C8D4CCE8ACA8B0D8EC155BC964B8A24AAE25CF144D77F70649C552AA42944A893EA1C5E555
988DD09A219529573113C74B3E97657091C009003F834AA406D45251C6693BBE2DCF932D3EC217EAAA0721CD821962FF8446CF20F648EE450EF7717D
8E337C590448850A242D3D1117F022436E01AA8420A1244315D6C057D28F1E38980AAB1A4053F2892CA7698C4044CE20233110D71DB6AB62EAF60AF5
D08B2FC50DB261A0286462D17420CD11E515FC36916CE6BACAA85230B4C264248B9D1A4A1C196580028B89221D0BE62D87E5011D9D1008EBD8D8DA43
38390A8F4970684FCB00EB98B0582E2BC5E9383D20731D2FEF8AB0793D59142CEAEB750B96DD18B2DE9FE6F23EA21EFE615E80B00BB63C7ED971B555
D5F33FF6E56F6CE68863C0B933B7799B8804DA20E56D211B4C978188B360404EE74CA2D9A191FE16174AC89AA02B8C5AEE8C9A3A12CB94C844D3894A
25341214918C840C4F26732CDD110BBC63F1C932CFB392886D5D5578088B5379CD9DC601C0B53D532FE224D9C0D9B43B8E4872760B313AA1EBD5145E
E01557CDC324AC24BA4A97699E93805B1D909E3E4F55F85211E88A485443C8503099F8D77192AFDB24E3B7B0FD11EE32EC384CECCB6880339302995B
200B888E9A8E940A25C1328050A00975AACEE583199CA560BB972B38D721FB0D26F60EE12957D290931487E81D2336595415AA3217C6AEC7F59E0EA8
A864CF307F64A71F7BB3E0915CF8C0D6EEACE958C9CEEECCD4F647966FEB18CD0FB60E8AF4E49EF621C60A9F55BFF8C4EB429EE17BF82754036E6350
49B73E79D3B9CDBEEADAE37EBCD5E5B8C276A72A9A904A292A826268B03DC5D196EC0FE1C4D8D481ACC6E74A3A400A54D978363051000C83A03ABDBE
B3C24BD8D8A4F4645A63FDE3DDD34817753695541549D4ADCC91390025083441A48022AA95B932D9E2418053C9AE3D4EDD81EC2AF948C070A786355D
45AE7696A04281A664EC37A0CB4768E0DC815788D21F708FF9103409A12891DF5254C2F647BA7C84AC575344D18DFFA459A1C938AC8B245B20279648
E779EC2DB02701B8D851382A51A0140900E2016D2DBC6F3246F12A3955E000515FA0CA32E07454E0702583CB0752451992C4A84065189563E151C925
155188E8A0C0E49E5EBF1C1E9F9C139943AD318E76C0584BB432BEFD607B3A3796D3B84CCF1BB1D2864164074FA95D76CC3531C7CBFA3DFCD30E09C8
2A716E60D35D9F5B5C5D73F22FDEDB3727E00F3E4406B7E3B975E3D8B44149144A7465AA63B460593202B4BF7746B490827886850A39B9234AC385E1
A802F85C8693C892BB215005D25F9767DB0F8CCF8AB8AE90CA105B8DC624829158B924937E1C49D3A59959111944BF880C05622F400E0488C9113931
48744678DE308120929B7532B5431EE54E1299EE43884E10D901B6C92491AEE3E2820116D91D42B6E91EF85BA8C4DB960909B107247506693F12C123
53566C9545E4D49D4C1C40571D94EC5294C612F8B5F18A20CB8A842D5BDA7F200F11C357D27919E740CC74447694B9908AE40A008A121C97090BA9E3
40CD741C63762A3E3AAE9686B64FB38689E87C5997533D01D98C2524FCA6C771A5DF7EEF86A030392259D327FB969C70A5DF337F0FFFDC4CC0CD03B2
876F3973618DAFFECCCB1F7C8F233C1A7A3E5494502994514A529E85F448B0CC43728E60E1FA980ACCCC29E49CDB642305D310A0228B6278FBAA9DDB
46F326D926522860996A3E3F37332BB943430E8206E48B6549E0295CBB2315D702101449FB0FE00A4342809C2FE00A5F967490CF48EE762EF62DC451
104D118DDCA27069B23EE80EF7907300C2CDAF8A3C42D01526D6544DE609DF2019FE832EEB88A1E68AF8A18AA0C864059F6410BAA50313A70AD81F11
261357AD9C7004BA66887819190639E4ABF00C8F7D42B97592F43F8A89249D832635934CF2969E1C4DC96545A0E25466FB21B224AC0358111C27B57F
6726ABC8C5C98128D95DA0A66292323D5A70CA11CA1041E1D010E2D7DFB7AB0C0F4F21736469ED314B3FD7EF99BF877FBE0F205F5165E4E96F2C6DA8
F1559F71D7A0EA1E198A6A64F59A18A7A8060EE7B0208812AD2099479A56CEE4CB09569098E0481E005E10145A1471129E4F47F28950482ECE66D512
ABD9AE3A0F89CE964976EFB1EDE2606F2850562542E143B8B814504A96186CB6447BCBC0B987C4D22CD41599E58A9C480A7C32866B102F6090053FCB
9517C5B19C54040ACEF2654323ACE48459C424598289449657705D8033094D5714E21864527D1082102D5322838E38912067830E290E4C13200192F9
1B315FC02E0E18AE7F317579B073768EB70C5DD6D41C6F706A2651213E62A02387848A9A9B9A2AB3A2AA8906D43856B6EC645B98E42CC9FE90A21B40
9FE90EF1037BFA353E57206392B1AE0C525A37C5ACD068C134B5AB6A172DBB709FE699BF870F8E0F90C7D7DF7696AFAEBAF1B2B7C6257797402E548A
81A972BEC0D09221E686467A36B7C9A44D2065C391A1A9880CF96896A7B21CAFEB6484C548F6EDED989091CC50F1A2AAB98C7F384F5058A8B1120F81
A892D120322320015D21353BD2C5228B33691D91CD1CD2C3C7FFCD9733997070A60834A294E9AA97935E9F452672C9F99E3B6C88FD08BE0E21FED188
36183922D04867028A0220AC0096E6EE0F11250E455488AE1F2AE76632C050246CB1145428A54409523119245945797AF4F004E7E8B645320ECB311C
F5AD27A39685AB0A4D2D96E5B2A8B13A769376E6400856A092999E869A8CAB0FC8FA43B2E5089BDB4D076712639D050EFBB8CC4040E5075AE62AC538
B09CCAE0C66E5A4FBED983A8B6DDB461657E5ED7BCF89C77A067FE1E3E50ED00538A6CFAF9F98BAB6A1B3FFDAD570620911D13E3BD331526074C3247
687265462824A2B81826BA1682AA40918E948B228E9A8425C476757149C78D30FD113D72727E561A98480982C42B44B71B6253952880DD005038A152
2E940036D264B05411A1A2A8D80F48504CC4784525FC1E48034C91235D3EB269EB92FE9A8490486418B2768CF37E15B83B7E8A488E0E558E1C0010EE
013767209E03419CA203B69213B9EC68583200F60F9651C92B723A1CAA14D2431D930564C8A10019D47774901D4F92EFA66D4462A6856B15284B2546
E290910E6A384B604AF869D0EC663F47CE1FE4D0C481F53396DD79EF2B84F94B9B3C9441313A397438A89BF900A70E4E423371E0C19F0D98F0DD6FBD
33FBCC35AFF09675E4A43ADF82B357F2DECA9F870F9C0BB03436BCE98E637C35F58D17FDB69F706820196AC0A0FCE3291CE56D931EDAB53FE15AB629
3151FFC0E42C6B18D4784F90C7553C47F4360DA1AF25962DF250E1A70E0E94E9B9D14C29972FF3B2A0D2F9A03F931898524C243040D7805A2E924940
2AC80101E7FFC8DDFFB154D55DDD47640E586658512E6781E52EFB93793F0400974A0B084AAAE4467E0980128DEB702D171488423261FFC80A504626
221A45A228A583291AA83299FD27FA685016112A94155344B97054925891038E81F86CA6303E5E7137F19163C9792643950BFEB98C402323B7A98730
07DA3612F5E80B7DD867C042DFFAFD63C3FD62EEA56FBFE73889E9C4DC742133144DB61E19873013AACCF9D3C5ECE4EB0FFC62C8D156FEA81516D66D
532C6BADAFA6B1FE8C75823BF4EB7DEA3C7CF02A01C5BFF5A12FCCF7D5FA4EFCE1060ADB9CA0087D3B47A934AF221110D5EBC17777CE94A16188C9E1
91128ED989705A827A65A6ADE3C0E89CAC31B1CCF09E4351860D1D3A18C830865D9ECCE32A5C108AE9428513285A903851260BBAB641D676B037C0BE
80ECEF69441D88183F99DF21C502D9A851899287ABC56792A560A810311F7C5F44EA7B0D489424D192A96820196008F327B6FB4A5EC2D7273B3D9224
4A722E89D3034BD379A459F85593B24300544E45B6A2E4B0EFE02192658B4F555419FF02207EAD164ACD4C8453F9F19980BF2B046074CD5E01DF4278
076DB0A5D380A012EBEE985264FC220B41D671E47856A523A1D1B09C8B49155AA8C852986363D3B16059B60B7B77309AC4C89A817EE3F335349CB551
B14CDB9BFAF1F001F501B649B52D3FABA9AAB666E99D93806C099B8E25B19D3BD6F51770021C1BEC6EDFEB4F498A6168597FDF74C17A9F97383AD6B6
BE274DB4F4A8A444CEF32C5B473C98DB348D645E565D7E31C754B38914AF1196611CCC454DD3154A000A525D6E70936C0E111501B7A8476EC7801CDA
11EE7E728B56CE93597BC2D703148527DAC5B81C505C825FA4AB383BC07506829A48BE0803ED6989E3456CB61251020625915745407310B0C5329444
3A5FD15850E47062E0209DE4401A333B339C83F4E4785166444A11FBF6A655AE7DF59C832F4CAA20BBD4032D8D4C449980CF2A12AEF06DC2A2EE943A
0E4E0B7A1E2039EC6775A97F78BA25CADA9693ED1A902D364C51489BBBBAB6B1B9EEA37B54C3F4467E3D7CA07D000CEFF8DE92EAAAAAC60B1F792B4E
AA5B5D0A05FDF9B27F7787409A7DDC48DBCED1AC14DC7720671AB9F6B5EFB626148397592A3EBAE9EDD6E9503C46B322D43885D385409268F0164766
58A22568C8ACE0EEF12BAC0A159DACEFCB3A8EE482480EEEDC2E9FAB35AEBA722206591B80447E4C07504288A675C2E2819D0704B2448E160D882B01
5106644CD1225BC14496442393C7C1A7F60A2A84B2E67208F0AC2C2A39493505C9C4358D5198280024C25CA6826B110765454773F4DCE86C3AA502A6
20EA1A2B41A4A44A48537239CD222AE3081836CF193A3D7C682E2F66150825972508FFBAE283014A96704E216C599755526FBFD2E3CFE89A1ED8F44E
C5B187FB19C3846B4EA96958B4F0E3DB5DDD16CF0178F8007B00FCF934D9E197EFFCCCB1D5F5F5C77CF3084FEA6A0399D4B617876440A527FBC395C4
747B14476D1B32A37BDE59F7DAA1D92C835468E892A46B22930C1C3A3891A44CB25B4F48F9D470DB68916CFED88456174769757A4A459CAA41019108
2D1E5CD96BB8277C22A9071051E82486AC024576233EE1FF30492580A0AEE4922C50C8A632212E863C5B9680CB61E6720292CD3C2423831AE16D4D36
249E55111F9ECEAA501589DC0FC08581A9F64F13CA73FF6C84D169D5CC8F46F11BAFC405404A20681B94805F41A922AB0A0D75240AD0AE847B87B319
3F6382C0CE6D0331197B9B625E055C9ECE87C292869307CED443AFBC320BC4BE5E91AC240E3FF7873E51EFDB32A19A76E01B7535F58D0B3EB68D98BF
ABF0EC99BF870F7443D0467CBCFBA99B4E6AAC6F3CF6ABDF394034EE74D574E4FE179FDB39C4E348AB646662D1DC5C9C2544195A36DADBD9369E1320
21DE3521AAC4FDD3D3897036998A1409F32622F4630ECEC309AB2F42209635B08D4BD81790506F44B7872CD2E5372A59BA5CE418EC145477051090D4
DC140817211924200381E2CC9448223CA121D60C0BE23CC2B20D8B109510522F972F84F09BE1670240E425A8C8A5B2C48669EC4C80A92BBAAAE46754
C7163BBA0B92AA29F8F5C4634017918193003915CC23056972745B2FABB208AAF9544553554128436E2E66B9C2C48ECC95E33D49BDFCC4775F9CE86C
890476ED1CE36C637C6D2BB00CA6FDC8C4A1C3FB573CFE6EC22EB51F62707EF0C231750D0BEB6B2ED927987F59F8F71C80870FB0FDBF4F272417FA56
5E7F42555DDDE21B1E7D63388D6F31F8689A05D9E0743831B9BFFDD04024331B8D14CA8CC464A2D3C140361C0AE4744B5575395B4AA5E35478E2C878
5E75B936B19D8AD95819E2F82E4372BA86244E430C402AB65640067E2CD354795EE14B491A1C65F2C25FF94A81A145C0E7B29144B442B2039A824813
2805A4D2643787706C231EFB0C7EB63D41C881090F924EB4C1088DA046EA76DB023875D7728C01087538C8A98E612747B139CAC876E8F25C5844C831
34B19CC88668136719A3EB4788EA2992E8B4A463EF62E25F4ADE1FA8C8B66A5B309E9F693F48E9D38FEFA47519C9A9AE3D5D15293F25E1DF4F7AC3DB
DD832DFB7FBFFC002E3D38D21C085E5557DBBCF4B8C6AF0CA2BF26FF9EF97BF8E0770248AD0FF2A36FDE714643554DCDB2EB5A894ABD26A65213AFDF
72EB9FE3AAC194D3543E383278A86D881C185A502D95670F6F19644132CC93797D720A2F9524B2ADAB692A120A490A204566241EE078ABB2DCCCC6A4
3BE943147A54B60420D91E8464F34FCC2545C2132467421CB667933097A68B84EF4BA015047229514C5206E1F2866EBDA0EABA50061A524C134886AB
CA6B7002210A23E4FDC8E544534D721098894C142C07099A85B0CF419599D01C243CC9A1A9A4881CA48B4679AA3B44260F0BEDADDDA37499CFC9A218
DC77309599CD80349FA10D8715244349D1D8AB51AC625B8E4615B027E113FDED63845658DEB36B8E016406427CF9D8FAE6F94B8F5B764754FFABF17B
F6EFE15F240F20D38052A2F5E7172D6BAA5FBCECFADF8FB8BA7E303B35393BEAF72713A1391EC84C251EEAE81929A9EEA1DDF8918048A5533C8EA39A
00490B40424A3E275BEEB91789AEA40A000A92CA95E8589108875AAA20694A2129BBAA82A60E24B952ACA8A601CDA3D988E60A911910A81A2C95C908
30701F872C8D4E615762B82A619AA1AA0C951D4DE23BE31A401C296974AEC2CA14796E5B1254D5D0633B3BB3C98AA0296406D840EA6C5F5113718650
9A1CE780633A7A3E520CC7F39A2AA1D07B87C3A560E79A3D134C64CFD6014A2BCD6512314A30E43CE5A0F2D8748A9E0B5574EC67A4784274721BD78C
F82BF89DA5E72238B391012E14365FE6F33534CD3F66D92F0A9EF57BF8D7EC0590A85FEA7FF66B1F5FDCDC70C9EAD75E1DAC1086EDFCC8A1ED5D2189
965C1111D57FE0F555DB42F251656F9C315772C160A440CB728542D874794127EABB00168343694D973895086D92BD5F32EF671BAA691B64C39ECE08
48CC93406E13ABB72CA4338582842B061D20439171798FA33A39D3D7140D17EA96EC8FE357C3CE165D1D20688A954A45B7C9819E5D985220CBF29A09
744D8C6732EE10723125DB96AE3BBAC564437C291C560C5B1533851C0F1DDD7114AA9C6218830C29C9FE1059EC0D87690DC5DAB31237FD6E27A71A16
979C91743372643097A06819A1D86CB9C25ADAECBA96989CE7C1E01F9E9EC2EE4C443AFFD6E71BEB7D679C78CCE265BF669DF739196DEFF4DFC3BFA2
0F30F8C8BEDF7EE6385F6363E3A22F3D34A091895C44515CAAEDD0E1031D65DDD04132341CAB30855C5EB1B19986FD73655E55FC33B205353A4FF15C
2A95A1CA99F1E9E8DC5058350411EA4023F33980ECF093B301959A1ECE89B47F3C43D6755D861E81A7E6A6E2140090117880881A0061FAD0151C85FD
3968C81588AB87A911C9D004221A48E6745C2E5E492EF2263014839C4BC895B942A510A23980330C081D0DA6773FFE7C3BF64AD06273F1EC6C46B44D
5CADABC5581A3A50E68BC01041859D0B27C87A62AC67726CB018EC4D615F07C7061547EDEA2971A69A0AB2F181881098AD984267D8A9C4BB27E73ADA
4A08E06C86DA704D4343C3A27B9E3C69E9925F08EFAB3779E1DFC3BFAE0BB0616AF74FBF7C72735D7D6DEDC79F4BBDAF5C25A4DFB9EFA908D10E33CD
623CC0A406760CCE1588FA0FD465DD10B1819A0615CEE99A0474978F90A518555309AF1FA1F0AE70858A42CE1DA5705C22646086942D91C3406C29AA
886CA4CB84ECCF3000EB92FA12050F9C3770B39383451CC771AE01E3411CD44D50917186A091093BBDD03AC4C99AC5852613195986F827A6D07D2022
6A32C2458015EC9EEA1CE74C4BB792C323D19C2A54E2E9BC520C764E03DBB2D98C681A0AA530C9B08434A5942D14A822726C7C154588C6118AB64E23
07E587C612E9E9C4647F0899A5CE906633934305C9B13459935B7FF8F9A68686C6A52F6E3C73D1090FF35EC8F7F021E80812E56171EACDEBCE6C6CAC
ABABB9E257AB5F7EA617FFD0541504D4C448EF6C54203DB6DCE49E9E191CDD910254D95469B2118C10050C9BE8661A966113AD5EC0E4455D03B1C9A4
82B0CD4AE1E19EBC6193E93AC304244E5B50A652B429138D1DA8E95CA26F82312C4540EE8E5F6E244859B8CED744958A2A36BE6025CAE2E7D31481A5
863B2624D39061653A942B086986880E89913CD10FB690250C1C9E926C5C6D584274603C6D620F939D1CCBD08199BE9463B3A18E11DE514C5899E82E
6A9A34E0A7840C2E6BC87C51643CA3236DA02F8E9D416264BC84B8C8DC24AE06529B3A412A3D1B078E86FD1BECBA6E7EB5CF776CA3EF3BAF1CEF3BEF
25C6DBF7F7F021E908922F4A6CC33D179E50555B535F575BF7F915256CD148D2A54ABA5CE270B0772C838ACDB08C44367505A95864933338022A4ADE
ED1558A60C44206B6CA827C2E9201F2E533223E4B3695E545C220FEC2B00042031D495940A710127EC842BC8600BC1FE846E517322A1F2946501E936
36CF549A0C09114D503A877310293B393498CA137571A2F1079165C00A9F1C8A2A84694821730C74EB9E02AE12A005F854A484E33D2F6B5291829AE1
5091706C747292D6B1A752F25146D3F94313FA6C4096CAC8B6B5FDBB45D58EEF0B6147967C7B3765A07C415101DDFEE757231A28322CF9FDB03BFF78
E5D2EA9AE68B3FBEC077FC798B9A2EDC217929BF870F970FB0D442FFF24BE6FBEA1636D4559F7ACBBE0A99E0C7711D82F0748656208EC35C7466205D
01BAAAE854F848D70429B03541D63499AB64F2F4ECD050B2802449B2549635548E571553158A953C631CE5E6D12436071C91D2DC9A5E678B92ABB961
D34959263B4396A18BB39D319151DDC95F03157179A028D4DC4485689C90757E11C996E5CAF971B3E394A80043D58A13B15C6FC2D6185DD7F2E1BCA5
D9E4CEC53CA552F8DB50CB1819F205914436D49F438656F1278470923AB2AE8BCBB6BDDCC2D8FAF4D651C57662EF6D1B61293ACDA940F2FFF9F5395D
A4F04F1D69EC9E739B1A6B7CD527DDF9CD8FDDFABDF39AE65F7458F64A7E0F1FBA4AC0B6F5CAC0BABBCE6A6EA8AEAB6D38F5975D007FC60D512F05DA
1EB9AF17C982C464C766FA2238BAEB966330A3A1304EE4654114224909A96A7636A50332CB6FA8482A8A48D1784EC1593AC3F08AA1E428F3A89B712C
22D52374EE19C9AA38B3701449229BFEBA51C95510B7E58068B9D281A60E2A7922770E05E86A73139F41CE0B755DCEB0D07235100C2065C6771E4A03
9A6625CB66E6A68BD081B6A108E51443692A60832195A80C70C36BB7B50673BA30E12F8A164DF95B0FCD21B0F58519CB520A91A44E33C550D0E0B371
9697294902846110BF2676FB5D171EEF5BBCECD853BE79DF5357DDF8F8AB572F5AF0AD5EC5F2F67D3D7C085D0091C096223B1EBBEEBCBADADABABAB3
7FB88BF064D986917CB74BA6C602324288CD4C4DB674CD20878CC61526DBDAB2B65D0CCA96221A1A5D103509E0221D4DB49675810524F337717580F4
4AE7BE89AC3F2011178063BD3EF56ABFA44C0511D7792021A9B87E3F38353CC39A5A0512AE6F8D0D1C99A64B9A8D2CBE884C5D99EB194CB30A2F6A96
CD8DCD2618C98038A11054369BA669D314F3FE08A7A9259AB665C92CEE78694F8C464C747AA88706168C6DDF15CD9590E180AE3FBD1D35D974BE7060
47816411A487A0AB260469BE6459482613CD848984288C4B93EBBE775E635DDDC265A77EF4C1DDE557CE78BAFB6BCB9AE6DF15D44C6FE1CFC387D303
9001375B63836F7CE7F805F50DBEC6E3EFDD4111511B13495C32160B04A20581978B91BEC91491C775503148013D7A78BA282A9C5CC8B3B3938C6360
FB65096B904262A8459475145DC826B2B10C700C8A217CFF858A6EE882A864A723962DDBC1B77B38C20A4660C2543CD1D75EC8B5F9155BEA98C097A3
0281BCCCA58EEC0A6B7066828216B6599DC88E4A32206AC4722E451995146569A662A8E99EC37195F3BFF7CA4BFBA75139D3FFD83301A28A649476AF
ECC5E69DCCA89AAC12965157A884300C698A010DA2904A60996A6AF7A3577FB2B9A6AAAAE1B4FF78FBE12FBD08833F3BF7C0F8271BEA8F792C459402
3DFBF7F021EC021CD518249F6EADD2F3EA8DE7CF6FF635F94EF9FA2E892867A8402BB63CF45457061086006E7C36A9DABA0432D343C95C314F7803C9
A860A920B00807488BA882E59229114252D743B2248F9F4393122371D332CA0C820A690A10D66E1CF255422CEA2A02A733B12DDB0ABA65C081F50392
96ED1889A6123CD1015198743850A0E624C765103190215444031172420BE00770F991C9282D68127E27063B7E78FDAE6E56F6AF7E75E51ED171C4B6
CEE2F4247E815329CD31C83401BE22D40028958BF2512E7272EC91EE797BCD733FBB68E9BCAA79F3AAEA8FBD72FD8EC0C82D6F4CAEFDE81F0E6C39D3
D778DE7AC1FCDBC68FF791F1F0E1AB01DC2C001B30AC0437DCFBF1B39636D4D69FB57C7F896CF94194CBA8862A89644B9EEF1D1C3EB87BAA8CCB77CD
9273B3C96832592EF0E9420551AA4AA67D109D8B16336C21319B91B06913AB2C1985146B08D3EF2588069881CDBDABB7A36491B95F6CC8D8B6B8A1C9
785C345425178EB338B19744255F88B31C1D9B8D28F80A9A6A3AC8828960321F27824208AA2E8FB0A15B8E35D519512C47D1B842FF789C12C9D8FEE0
9A2D5943CBAC7B6A4D1CEAA613D83CE5108E7FC71A5DFDD0ED977DFEECF34F3EF9F8732EBEFECAAF5CF1E52F5C77E9A71634D5D5D7E0B85FD57CD2B7
6F593DD0F5D4FD2F5C71DBD0571AEFDD72467363D3858715B2F0E3D9BF870FB703C02035AECEF90F3C71E5471B6BEAAA975CB5812656A62323BC694D
4F50C177365916589AEDE6EC1ACF44875BC74265151705257A2E3415E58EE6E62637BEF3E0A1086799E2704F1A1956251DDA7E44C6E93FD7379A9504
8B37340701130EED9823DA40440AD81499A2A8E342025A64BBAF54981E1B8E09866A1D553A31E35B77BDB7ED300D2BE1B48C7899E635D7284DA20326
E49331FF70C9204AE25C81E191563AB46B076D3A3CCCF4B4CC40F7E4BE787B7DF53C82AA7947BFBBFF5EFB9153CE39FB98C5A7367CE6DEB1397164DD
DDB79EB07CDF7D8F7EFF82175E3CA6A979D9CDED84D6C87A7FDFDF337F0F1F6A17E07ECE7531D3FDF4351FA9AAADA939F1FE7D8C8E73742A91E742E1
F1B86EE3ECDD56403634363693D7C8E40E8EB8F9C4EC1CA5AB5C30A95BC8948B65A0E92A1D29CA38630813665D479779990884F6AFD937B72F488364
96D1B1FD6AA9DE2830F8049D2FA66906D8388E1B0E9308062B8C88140088E6972DCE4DA745C7E0720A4501DD5407F7F5CCF44DC600613340BA8A9393
74B2C44A2EFB1884144F2359A70A1CCEF7CB63ED39FF344F98C11CF8E7531A162C59F48DE5AB57BDF1D5B3171CD774FA19E79CBCECFCBB775365C3FF
9D5BEE5D1E28828EDB4E3EEF82754FFFE8BB37DDF7FA3B2F2C6A683EF6E72164FDCDFA3DFBF7F061EE031C9D6E27CC4150C81EB9FF5375F3E6D5569F
F5E3C369771D485526A632E14956D51599A392E16045B5C8A69F282707366D6C2911E64E1CFAA1BFA5636A385E82B663499934B05954C94AB661EA82
4A0FF68EF4B7C7C73B4274282D180A11063040FFC696F6D19182A1B3F879D8507FFB962191280400155012A07BB6EF6C0BDB8417483575A0A777ACEF
9A8902263536D633C309A97C2ED1DB5B3465C7E462955034A342211C0B40C3A45B07B245153A40C17EA8FFB2E6C6C575A7BE028988D8DC586BD7914C
451E5D396128C6E486A70E3FFF5AACDCFFEC17AEBD7EE7F2F5379CF6D6BDDF5AFF699FAFF9D21D45DD5DF6FFABF97B0EC0C3BF433A408432C5F096DF
5EBAE8BFCDAB6A38E92BF7F74A6E160ED552FFDB9B0E4FF338281A06CBE4FD2196E4EB9A14CAE5E954ACA89940505450628B1230A47C1E015BC9F5F5
657441735983740D596ECF1D998A56E2299A62E6BA262A34009A842D934F0DF70E248A8665C9850827CB92C2E7274BAA635AC51423035CF183C1BD31
558360FC60CFD440D6B4187FEF507780ACFE273BDAE62A9A307B68283ED2C36BF9BE29DDD659C192F04B7FF4B8854DC77DFB8FDDA2D01BDDF48711A3
902D81C2F4FA0E46D979C7552B873A7F369E79F2C60D2FAC59FEFC776FFEEECAA796341F575773F2935122206CDBDEC6AF877FBB7C807CB1D462DFD3
377C747EB5AFB1E9AC9FBFDC91203F64FBDB426C59A221AA8C4CCC66F290CCF59AA6AE73201B4F88A1FEBCCBF56DE8145776F3061370B266E8720602
2014C7F2FAFB1B471C239426C622E9304E0368A888B229E7C626A6A188FD8FCAF46D7BB77B22C70B8139C5311C2137D41F934B62BE980CB30642ECC1
6DD132A90EB4D4705F54C457D3C7F70D957443EF5D338DE88A04F3B31929C4C6C70C0B27227F5AB874E1975BB0DB181AA4BA1E5C9FDBFBF8687EF8E7
8FBD99645EBFE1893DB9153F8D0CDCFBDBF553775CF6E61DD77CFFBCDF7CA4C657537DCD11E56F81DF0BFF1EFE3D6108A9AEDF9DDF585FDFE85B78CA
2D2F078189344DAB5492E9E1AE8C65DA9A14E83E323C3C93CA67786401980D0E07429958AA2C96F97478B06FEFD635FDA58A3FC888453A1E3A72F8BD
CEC32944A72B7383DBDA1374AE48D4BE101084D9FEB6E0D0885F838E81A864FB68A838BABF2F9F2F272AA2511A6B497148532A7234C6EA08AAF1FD2D
D384398093F2111138D051D21DBD7EA459881A8E6920A151636DFD01452D0F500644D66B1F3DBDE901C931038F7DFBC8EC0F7E997EEE9ED1C9DBBFFD
5C2E1D5FF3E09FCB7FBC75F3E4D6EBEE9D7AEDAE87EFBE6DE9F75A2FBC7C7163DD89BF2FDB8E6D7BC1DFC3BF796390F880CA91DF5E78624363BDAFBA
FE135B64D5E5F7D123AB7A5D4100892FE503A181D14CA15886D0E4464603E1BE8110C031B9144B4CFA076289EEB7D60771A52F4B1C94CB83E1D8A6E7
5E3B3CD93D13677004D799D078A863FB4BEB5E3F324ED3B6CDCDF9A7E6828A6980ACAC9615599AE9DA1BAB58406615A44BC8C8AD5AD13F7AA80869AE
34369C17451B96A36CFF913ECDB0742BB6A3579A9C94D243733D255517B382A53AE88645272F78D67162CF7FF1E7470E7DFFBE97EEBA3354B9E12BFE
5460F9DDF7F63DFBDD475EEDFBFDCDFB9FB9E1B5B71FBCE2EAA75ABEDEF4F1450D5FEFD78FCE47FFAD2FE27D183CFCBB7603F01729F8DA2D27D55757
579DF5BBFBEEDECDE3D4D8E0039962302C1AB6619B66A52CC70E0EB3150E6A44F35367684E4CB2EEE3CD5447201A48974CD3550235355160241B67F5
442C08706387C7F2B3B1497F81E7352807BB371ECCDA96058BDDDD899C0AA30777864B96460DF40ECC2424D5E056DEB62183A0A68A5AA56F54B098B1
96B12C9D15B9CC1CDDBF72ED64363B27E707E84A09899D23D0D04CFB76DF69676E74FCBFF8D2BD7BCA2D3FFAE9CCF6EFAFC83EFC8DE1F82FAEB86360
64E5F79F5C75EFE30FEF1EBEE281C3D79CF4D86CF0525FB3EF846799BF1ABC67F71E3CBCAF2BB6F23F4EAB9937AFA6BAF6A4DF4D90D1393EBEEF10AB
10613E4DE04B89AE3D5BF68E520C428686B293C381BE113F655B50A533C1C0A1681900C3D0C831A0425C8AA95BC010A5E11D7B449C1A5896AD0C1F8E
986A341529588E95ED1D383274684AAA74EF9F54C462B492A6A7FAE62223C1C87EBF6D9B427234C64ECD0A7AA17B7F80A8F7C0E1965C7847AB64F907
CAD111263487A22FBDCDEB82666DA83F66C15AE7F58FDCD2C6F7DDFD89DB9EF8D965CFAFFFF125EBDEFBE93D2F978ABF3FEDB25F9EF3F8B695377FFA
8515975C75F1C5779F5D5F37FFE651CFE83D78F87F5602A61CD97AFFD90BABAAE7D5D57DF1A5B0817F02351EA4860F768539DD1443FD69C140840918
FF9F91C1ECA63D9D33C3A3E3F93825E310BE3F26425A2E098C9C9D1C28EAA5AEEE9D2D2DF13249094A51B93C1B93445EC74F05E8B1B182A66BACC871
62A17D5CB035CDB225046737EE2F910545980EF84BD954BC5C1E9B546DC772A4A1298548891667FBC3459A1746FC93DBB69775A81883C72C3AF11BB0
75F14FA40D0F9D7AEE8D3B5EB96F7FE1D987F6EEBDE9B7017EF3054D5F79E2E00B0F5E73DFED4F6CFCFC4F86AEF4F97CD51F7B53F0FEDA1E3CFCBFFA
00DB5253ED8F5CBDA4B6DA57BFF0F39B44B2DA67B0E55084AE8C4EB2BC62999618696BD93534C3D926ADB3819EE9A85466C9B42EE27AB6BC33323BDB
1D2E49F9B1FDFD19A14271803734894B6D7C72AF6411CE1DC7307427B9FF20671B8A415404A9D9318E9C17122D213DDD3E261A7C68D7D6115543AA32
DB335254716EA2F6EF8E93658258D7E1691A3FAED23B31D6457148908DD98FD434F87673E77D67EA37675C71D56B8F7CFF4B6FECBCF9E7D4F897CF5C
DEF2C259675EBFE5B6CB8FBF76DFE635FBAF3C71D59A139A7DC77ED76F7839BF070FFFDF598005F9A997BE797A5D5DFDFC93EE7A634226E7FF00C2EE
ED152044BB87C6C3A1543A33922CF0C876F7832C9D4DCE1EECA411A013A57C4EA25AB6262D43721C5D02169F6EDB7EC0EFCE163BBA0561858B460239
558764432F55100273BAE3EEE8D94828482234CCEC747F5EB3B06F6027662AB665337DCFADCDE1D7A0B6ED99086BC0CA0E6D6B89AAA552461210FFC5
79CD359F2BFEE6DAD8A397FDEEBBAB567EF1E6173AD6BD9CEEBCF2B2675F38F747EFF4FDEACBF36F7A7EFBF253CFBCFE98C69317F81A3EF926ED75FA
3D78F8FF4F0270A4D6D25B7E766A537D636DC3C756E64C320D048B8C927BF5B7EB8610AEE63521F0EEBEA1D6314E352B9160767A62DF7BAFB714802C
73A12C150C25C20959E5DA57ECCCE6270ECEF1F8B2E5C891AEAD1B5F7C61F9EAEE740A28728CE7260E1D0964290D188E998F459313BC68208DD50837
E778CBFEF675FB2615204FEF7A7E550BC2AF67E6D58DA30223F807DA075A52FE8943834056B5FBE7D5CEFBE4DE67CE58FF9B2BDFF9E97D8F9EF2DDE9
1D4F66CD81732E69B9E7B2DB4746EF3AE5AEFB573CF1A9EA05775E5C535755BDEC9E09E891FC79F0F0BF7201AE0A889AD876F7E71635D4D69DF0B51F
1F944CD39444B1C2A9E9ECD4689271E4C1DEE050BC3839391B89B1AC86A07F3C37F3EE7E964FCAC5E9FD3D3C98E8191DDC7578556BD1B2CDA9AD6F8F
F973ADEBDF1D1CA7752BBC6DF74C0146BB3A29902FE0D23FB3EEEDB1BE503C3DD8BEB375D54009195AB2BF6D722A4D0933EB576C1BB76C35D2FEF813
2D8A88CA7DBDC1B0911A6BDD3F24AAD0189A3F6FDEA563B3573C3DF0F5BEEDDF68FBD1E5D3A92F1EBFF9E0273EDBDE71ED1F523B3F7DF21DDF5EBAF4
C40527BEB4F394DABAE37F7280F3FAFD1E3CFC6FE500D8071076CD8E3FDD74414343ADEFF49F3DDB9BD04D0BC9C5E1279E19A5F3340D0D8BDDB7AE3F
A9C9B85827EC5AD2E88ED1C8C0CC5458524D5745C8A1E7BA6663E5F8543F4F460591A5E9B2910B8DB4B4974C5DD5B9709848FD99E14E9A3C69E5D0DA
3FBCBA3182EB7D8DB016E9D970A86D6B2FB2731B5EDF78A0A5A81BC98E7D51C9AAC48706FB2951134CE9A279F32E2F755F7647DB0F5E3D72DE5DFFB1
E48DC4D7CFF9D97517DC31F3D8392BE65EF9F4A5BB5F3FFEB23D0F357CE6BED3ABE77DE21DCEF682BF070FFF99420097E552EAC0D3379ED6D4E06B5C
7AE5F251609A5A312B97A3E6DC8BDB33FE4449871A12673A06A6DABBC24C29D6D5C3B133417F7A26162BC68A2254D598FFF000ADE753D3A168689C92
F8A9AE293144E48560A853B450AA278FFD8726E493BC2C1745D37198D5AF8541BEFFAD37236ABECC1A959171D93498BEB7D73CB77253A124470F1C29
42032AC8BC6BDEBC13A6273FF3C9CEBE1F3C72F515CF7CF59999DB16BC1CBCF217ADB77DEAC5773FFFF137B6DFBBF89BB1CEF39B9A6BE6357D7F54F5
2A7F0F1EFED39500610FA3C7DFF9F979C754FB6A1B3FFDE8DE1CB624315B4AACFBFD738725644031D8BB79EDFAB6C33D07DB0FCDE9AAE0DFB276E5CE
D5DB56767286A64BB4480753134756B514F25DBD3943C3B62E9706F7B407FB47FBA22AB7E7E15FBCF95647EBE1971F7FA76012B53F4B1DDED6199D79
E589BD45598E719A1801B665F2DB9F3FE40F993A4C95A2D371D18202B2D65555CD5B153EA379F98B5FF8F28AB59BEE5B9178FECC1F3C71F9A5B72C3C
F9DD572EBFEB99F517362EB9EEFA93EB9B6B6BBFB4B1E0F17B7AF0F07FE0028E6EC79B20B6FFA7E72EACAF69F01D7BC1BDDB68A20328A65921D0B171
6F4722C391443E1F4DD286924BC64B5945419C34F0EB3B97AF1B1F6C4B1AA97C9CC25EC488F51DF053D803405498E9658052E96D0DF56E3E28197A3E
2B5B7271D6DFDEB7378A730053C1AEA37B181AD34FFECE8F530EB1246916D429FFA1B18CA901808019585A53F525EACE85BF58F7D94F6E58F1ADAB97
7EFCFC8FFE7AF8B19F3C70CEA3FED653EEE8BAF198E64F6CF871E3C2FA4597BC9A40A667FD1E3CFC1F5702644B5ECD762EFFC219F3E7D7372FF8C4AD
6FB74A86A9CF6DDDD7361393CBD9C3BB06A622995D2FB605FB9319A01908E492FD9B37774F8CAD5EBB71D5ABEFED4C8EED0D94C233D3D1FEC9844A36
8591529EA5F882C118B24A7847506EE74B7FDED3910ACC4DC7896677646CA0A8A6F6EF3C32D9BD7BC580A4163233AD938111019848172453B8A4AA6A
F1D0C6336FFAF5CA5FBE3C70EE47EFFCFAFD575DFFC8B7CEBCE99E1FFE79F9B917FEFC8A450B7CE75C787CFDB2AFAC9E336D2FFA7BF0F07FE3038E6A
63A2F2F4A6DF5DF9B1939A1AE7379D75FD1BEDAD211D72D1049D7AE1D7AB77E76D334F2BF978C7DED5CFFCEAF9151D656C74B83C10E70666C25D53D1
884A47355D8B0F6FDDB03DA41885AD9BA35008A70DB3FFE515939CCE64C7A7081F187B78F96F1F79A39F922C31E18FE67535DFD7B5EA9D1ECA3F3619
A00B62A5BBFF08AF22EBC1AAEAAABB074EB874EB9E5F9D7EFD8DB71E39FCF8969BCE38EB84F9679C776CD5C28B1EDFFBE8B225CDCD0D4BAFDD90F4A8
BD3D78F82FF202643C880DEFFEE1274EACF7D5D5D72FF9C696514E264AE3BAA54B635D1305153B847447E7C1C9C8F45B9BB7F60618455771F62D2514
0B96FD81A129594B7676660DD0B52D0C2618C731E8F1C06C6A7C20A0EA2654270EF6A62A6AB16002FF70FB0C633299BED1AC81C470ECC8BE3DE1E1A0
3AF1FCCAB58C60EDAFAFAE3A75F7C71A7EF8F00DC79F77FBF9CFFEF184F3BEFBF5EFFCF8FC47B67EBDE98A87BE78C19FBFDFD8B870E1951B2286B7DA
EFC1C37F612D40BE212AB0FA27179EE8ABAFF52DB9E41B8FAEEA2E591664DE7EB26F7CDB8B07F250B54CCBD12677BFB9E18537FEBCEB483413DC7D30
CF75ADEC2A0C7687CA8A6698A59D075425653A4AE89DCD798B8D8FCDA5D9C891F57BDEDB19832CE06323533D239483CA2D2B57875959CE1DDC309A4E
AB7C269D6E6FCF6B8EFC899A9ABAFBEF6E3AF3DC7B1FFBC9AA3F3CBCE12B57B76DFE8FD79EB8BFFFB1A55FFED357EB7C4B1B6B175DFCECACEA59BF07
0F7F171F60B2B31B1FFCF49286DADAFABAC6B35F8D125920C7E4F7AF7864F3BEF58FBFDE352D123D102614E8DCB9AF7FFB906617A261DE3034A4CA12
9C19E60C813584D6EDBD958996090014B1B4F7F13D744691C26FBFB8BAAD4723CF30FED6CA60969A1C68DDB86E9A4C0C17A6428A6AE39BEEA8ADADFD
EE330D9F7D6AF9A6DF3FFEB53B82BFFFC6E8BECF9E78CFCD8FDDB2E49B7BBE54D3D458BBF8532F4D6A9EF57BF0F0F7F401F4C4DA1F7C647E7D93CFD7
78ED93A3139DD071406F4FEFC8D627FEFB256B441D6ABA501A9F2ED862F850574C1C3D7460702CCFC058D7302A1FDA958ACDA474CB48F228B4797749
084CAB1A87CA543491B334DDB0C26F3E13953421D5DFDD55884B6CF9D08BED9C42760F9DA7EAEB6A6F9BFDCE73899D775FB9E8B20B1EBC76D1FDBF3F
7DE9F9672FF8F49265979FB7A8A9E9D8ABD7278E127B7A7F2A0F1EFE9E3E0086B7FDE8D34B8F6968AC6B6A6EBCAC95C8081286FFC9B7DE080426F66F
79674FC834B4D2D4EE036BF795B35B56FCE6C927B70D954423D93FABEBB6298DCC44675A57DFFFECBE1D6BDBB7EE50D5AC6E9766828737EE1BDC36C9
D87A7CEDB6696E7C78BAED95158FAC891A86CA42674F637DDD1DE1DB1E4AAD3EADF9E21B1FB9FAA26B6EFFF529CD8B975D74CAC2F94BE6CFF71DF3D5
1579B2C4E06A0B781EC08387BFAF0FB0F8D175B77FB6B9A9C157EBBB784D4B9825FA8286520A478767625CA5D0BE61A4CCA254285E60F9B8FF406F8E
3BB07D0EC96CD6D684BC14DEB66AC01FE8082743FBE36C4836E2FB66C6F6F7F815E4502DCF6F1D53A8C8A64D9189EEACAC411541C18C1D5F5F7771E7
F77D1FBD60C1B2AFFEE0DAD3175E71D9F99F58D678DACD379F75C2E263169EFCD91792C8F1143D3C78F8C7F9004B4C6CFEE91527F91AEA1B1B8F3BE1
6B0FF74690A16B3AD4B9E896D75EDC93D091AE27B76CCEE9BA1658F5CADE49DEAC1CE9DAD5D92729C53C9DCDC82AD4E8A2888AA635F1D0130315C3B1
CDFCC84F7E7840836CCFA181C99222C1CC10A7C9B2265F515F77EA9F2FAB3FFD8EDB2EB8FCDA93EA8E7FE0B767FDF18AFAEF8EDDBF64D1A2E693AFDD
18323C425F0F1EFE811EE0E88C90941F78EDAEF38FF3CD6FF0D537D45DFCEB8886042852695AD215A9D0B175AAA25486361F9CCAE7D8E9B1BE48C50C
ACDA5FDEF3FDA7A6F4488E4E66CA34A36B70FF55774E8B6A21DDE5DFDCD99D93A27FBCE2CAAD9290E4C5C99B5F5555C5306EABAF5FF0C4F79A2F7AA1
F59E33CF5A7AC10F9FFFD5C24B7FEA6BF8DAD71636CF3FEDA6BD39FD7D5E6FEFCFE2C1C33F360BB0352E7CF0B7177F6461E3025F6D75D3AF4B477582
919EF9C3553FDE97139552FF3B0309838F6E5CD18A6C58CCF252CFAF7EBD25912EB53ED3C2C274FB2B2FDCFFD5BD8EE97F79CFAA2425075B36BEF6C8
F2770A8009275B6F7B28A30AC0FA5D6D5DC3FDF7345DF2F8BB972FB9E8EAD3FFFB4FAF6C5CF6C099D8DFF86A8EFDC9A1E25FE57CBCBF88070FFF580F
E076DC4C3933F4C20F3EB5607EBDAFF6E23FFCFC7B6F95F16DCCBEAD21B5B8E1579B815A7A736D4B3ACF53D3EFFCB18596A021F6BCB9252BB0A19577
3D7CD3D7BF772031333DF3EAAFBAE5A2C60F3DF46AB2229ACC5457281D3BB09F839260BC535B5777D54B4BAF79E5D79FBDE87B9F6DAE3DFDF4058D97
5E3EBFD1577BE60FB655FE46EAEFFD393C78F8C7E7004452CB71746E6EFD83E72EA9A9AB69A8F57DFC876FE5496EC0C43AF775EEFDD5A3DD0349CD76
E4DEE75E7EB6CB4209CADF95D215C8FFF99A3BFF349C746C7AF5AFEFDA094C7DF0F5B6237909A623EF3CBD3544E7525015387D724175CD79CF9DF9E5
5B6F3BA7F992AB9B4EFDE15B97D5362E6CACAD3AE5D6ED69FB7DF6226FD1DF83877F6216605A8E63A9A9C34F5C77C2B2E625CDF37D277F6BC5BBBD31
A8E962A075C6344DFAD06BED92ACBEB46176C37D6BC212974BD2D199A4026D08A9D20C2BB31B37CDB40DEA349D2DEED9365C8070FFE6B22A4A5AEABC
79359F5C73E327AF5BB0A0F6E49F9EFFD9D7B7FCA0B9AEAEAAF6FC7B8F942DCFEA3D78F800E400AE07203280263BBEF7D16F9DFD91658B9A1A6BE67F
E2E91EA01986214266F9E76E7A38942D0653D3032D339238F5E6EF36BED83B974DC9B1916D539C63F4FE6CE5AC218747770F85F2D0B0CCE93F8CE98A
A2759F36AFEAA47B3F537BF185758B3FFDC8034D4D672CA8ADAA3AE9DBAF0544FB2F84258EA7E6E3C1C33FD9031CF502380DB0617962D3835F3C7161
63437575F5F93F7837A05B8E321BE103D3EBD6C714C971D4F1150F0FD2B10D6FFD7EDD9A577FDF9173ACF496FEA42116F73DFBB31594693A865A1CCF
1A00696B7DF36A161E5357D3F899CF7D7CC117AE985F5F53575575C2CF3A44AFE3EFC1C307CBFCDD384C7C806D412EBCFB81CF9C5857575B53356FE1
033189C177D375AACC4F67A6F7BC974CFBCBE9194615E0DCF82C7E14DD3FAAA972EF8E3DB919D6B27461A0336B1AAA01FF5035AFAE667173FD273E56
73F637EB1B9AEA6A9B965DFAA761D9B65D196F2FEA7BF0F0417306EFBB03839BDA7CEF45CD35B5F53535275D70DE7FAC67B1D5A23247AFBD679D6E23
42FAAB53D1E148D250A6A765C799DB16453AF61D088D6DE8552C1D99891BAAAAE65555D7CD3FE53B1FA9B9E14BBEFAF9CD173DBE27E6EDF878F0F081
F702F81F438AEDFCD545CBEA6B6BEAEBEAEACEBD9B90F15B0850EC70E0E0FE646AA86D782C666AB1CE9CEE48F4C40C340D4560CBD934BE9B6AEA3B4E
ABAAAA6E38F3B48FD6CF3FBEB17169BDEFD47B76A5A17D54BFD4B37F0F1EFE05A0E6C7DEFDE5A58BEB1A1A1AEB4EBAEC9977C29A5E60935C5608EF5E
BE8E32D5F86C3A2F43DB8A56AC4C6BDC0419221BA61A1675AFAFBABA76C1E7EEBAF3DB27CEAFABAB39FE86B7038AF3D715240F1E3CFC0B3406C8379D
6ABBEBDCF98D3E5CBFCF5B706B4C43BA8174492B1572CC447B4683727436521474AE6B44944D4D530C5B3F78617D5DF5C28BE657D5FAE657D7367DE1
8909CE337E0F1EFE057D0051FBB3E88135377ED2E7ABABAF3AAFAB67652F511702A020E4594D1761EA9D4EA00122362221D5B2CDA13B9A1A1736D55F
F985065F7DEDA22FBDD0C3BF9FF67BF0E0E15F340F50A36B6E3EB66E5E5D535355D3050FEF05966D99D0D06411425AD264001004D032D9C08B4B6B6A
172C6EF61DDB54D774C6AD07DC5142CFFA3D78F897F601F88B12D875DFE58D3575BEEA9AAA731E181189DA80AE28A22A8B2C5010348D816F9FD6E4AB
AEAB6B5EBCA8B6E1AC1F6D98646DCFF83D78F8706401B6038B071FBF78597D6D5D435DE3590FAE79697309E9BA6198A66D19807EA1795E757543BD0F
97FD677D7FDF9CE4ADF678F0F021CB02CCD2D0A63BCF3FAEB6A6C6D7E0AB597AC37B23A3BD8199EE476FFCF205B5B5BE86BAAAFF5635FFB677133211
05250B869EFD7BF0F0A1F100A41FE8207A62F503977CA4D137AF6A5EF5B205BE450D7535B575F5B5D5B5354BAEFAF15BA38AE5B8A3846ECEE039000F
1E3E4C8500F101B64627FAF73F77C7851F6DC03EA0A6A9C9575BDF70DE8F5F1FA4C994CF5F383D3CDBF7E0E1C3D709B05CA561C7D185ECE4CED71EBD
F5D2538EBFECF6156D294056882CCBFE4BE4F7AA7F0F1E3E945EC0799FAFDFB6758DF2B71D4E32068EFBD6D1B4DFEBFB79F0F0EFE0068E5A3CFE9F43
B884FE62FB9EFD7BF0F0EF9207B84C82FFB3E97BE6EFC1C3BF950FF81B8F88B7E1EBC1C3BF9F13F002BF070FFFF6E6EF39000F1E3C78F0E0C183070F
1E3C78F0E0C183070F1E3C78F0E0C183070F1E3C78F0E0C183070F1E3C78F0E0C183070F1E3C78F0E0C183070F1E3C78F0E0C183070F1E3C78F0E0C1
83070F1E3C78F0E0C183070F1E3C78F0E0C183070F1E3C78F0E0C183070F1E3C78F0E0C183070F1E3C78F0E0C183070F1E3C78F0E0C183070F1E3C78
F0E0C183070F1E3C78F0E0C1C3BF2AFE0756E173B3
>
}
\edef\QuillImageObject{\the\pdflastobj}
\pdfrefobj\QuillImageObject
\newbox\QuillImageBox
\setbox\QuillImageBox=\hbox to 2.7in{\vrule width0pt height4.05in depth0pt\pdfliteral direct{q 194.4 0 0 291.6 0 0 cm /QuillEmbedded Do Q}\hfil}
\edef\MakeQuillImageForm{\noexpand\pdfxform resources{/XObject << /QuillEmbedded \QuillImageObject\space 0 R >>} \noexpand\QuillImageBox}\MakeQuillImageForm
\edef\QuillImageForm{\the\pdflastxform}
\newcommand{\EmbeddedQuillManuscript}{\leavevmode\pdfrefxform\QuillImageForm}

\makeindex[columns=2,title={Subject Index}]

\newcommand{\printsubjectindex}{\printindex}

\allowdisplaybreaks
\newtheoremstyle{bookplain}
  {8pt}{8pt}{\itshape}{}{
  \bfseries}{.}{0.55em}{}
\newtheoremstyle{bookdefinition}
  {8pt}{8pt}{\normalfont}{}{
  \bfseries}{.}{0.55em}{}
\newtheoremstyle{bookremark}
  {6pt}{6pt}{\normalfont}{}{
  \itshape}{.}{0.55em}{}

\theoremstyle{bookplain}
\newtheorem{theorem}{Theorem}[chapter]
\newtheorem{proposition}[theorem]{Proposition}
\newtheorem{lemma}[theorem]{Lemma}

\newtheorem{corollary}[theorem]{Corollary}
\newtheorem{criterion}[theorem]{Criterion}

\theoremstyle{bookdefinition}
\newtheorem{definition}[theorem]{Definition}

\newtheorem{problem}[theorem]{Problem}

\theoremstyle{bookremark}
\newtheorem{remark}[theorem]{Remark}

\newenvironment{chapterguide}[1][]
  {\par\addvspace{8pt}\noindent\textbf{Chapter guide.}\enspace\ignorespaces}
  {\par\addvspace{8pt}}
\newenvironment{keyidea}[1]
  {\par\addvspace{8pt}\noindent\textbf{#1.}\enspace\ignorespaces}
  {\par\addvspace{8pt}}
\newenvironment{warningbox}[1]
  {\par\addvspace{8pt}\noindent\textbf{#1.}\enspace\ignorespaces}
  {\par\addvspace{8pt}}

\newcommand{\Q}{\mathbb Q}
\newcommand{\catz}{\operatorname{cat}_{0}}
\newcommand{\Clz}{\operatorname{Cl}_{0}}
\newcommand{\sClz}{\operatorname{Cl}_{\mathrm{sph},0}}
\newcommand{\Hnil}{\operatorname{Hnil}}
\newcommand{\Hnilz}{\operatorname{Hnil}_{0}}
\newcommand{\sHnil}{\operatorname{Hnil}_{\mathrm{sph},0}}
\newcommand{\nilh}{\operatorname{nil}_{h}}
\newcommand{\hofib}{\operatorname{hofib}}
\newcommand{\hocof}{\operatorname{hocof}}
\DeclareMathOperator{\im}{im}

\DeclareMathOperator{\Harr}{Harr}
\DeclareMathOperator{\id}{id}

\newcommand{\Com}{\mathrm{Com}}
\newcommand{\Comash}{\mathrm{Com}^{\text{\normalfont\textexclamdown}}}
\newcommand{\bs}{\mathsf{s}_{\mathrm b}}
\newcommand{\ls}{\mathsf{s}_{\mathrm L}}
\newcommand{\qiso}{\xrightarrow{\simeq}}
\newcommand{\ontoqiso}{\xrightarrow{\;\simeq\;}\!\!\!\twoheadrightarrow}
\newcommand{\qis}{\xrightarrow{\simeq}}
\newcommand{\onto}{\twoheadrightarrow}

\newcommand{\indec}{Q}
\newcommand{\mA}{\mathfrak m_A}

\newcommand{\mC}{\mathfrak m_C}
\newcommand{\mM}{\mathfrak m_M}
\newcommand{\Disk}[1]{\mathbb D(#1)}
\newcommand{\AugCDGA}{\mathbf{CDGA}_{\mathrm{aug}}^{\mathrm{con}}}
\newcommand{\ComAlgNU}{\mathbf{ComAlg}_{\mathrm{nu}}^{>0}}
\newcommand{\dgMod}{\mathbf{dgMod}_{\Q}^{>0}}
\newcommand{\Indec}{\operatorname{Indec}}
\newcommand{\doi}[1]{\href{https://doi.org/#1}{\nolinkurl{doi:#1}}}

\title[Homology Nilpotency in Rational Homotopy Theory]
{Homology Nilpotency in Rational Homotopy Theory:\\
Cell Attachments, Retractive Towers, Divergence, and Stabilization}
\author{Paul-Eug\`ene Parent}
\address{Department of Mathematics and Statistics, University of Ottawa,
150 Louis-Pasteur Private, Ottawa, Ontario K1N 6N5, Canada}
\email{pparent@uottawa.ca}
\thanks{This research was supported by the Natural Sciences and Engineering
Research Council of Canada (NSERC), grant no.~149521.}
\subjclass[2020]{Primary 55P62, 55M30; Secondary 13D03, 18G55}
\keywords{homology nilpotency, homotopical nil-length, rational
Lusternik--Schnirelmann category, Sullivan model, Andr\'e--Quillen homology,
retractive tower}

\begin{document}

\frontmatter
\begin{titlepage}
\thispagestyle{empty}
\null
\vspace*{4.5pc}
\begin{center}
{\LARGE\bfseries
Homology Nilpotency in\\[0.35ex]
Rational Homotopy Theory:\\[0.65ex]
Cell Attachments, Retractive Towers,\\[0.35ex]
Divergence, and Stabilization\par}

\vspace{3pc}
{\Large Paul-Eug\`ene Parent\textsuperscript{*}\par}

\vfill
{\small
Department of Mathematics and Statistics, University of Ottawa,\\
150 Louis-Pasteur Private, Ottawa, Ontario K1N 6N5, Canada\\
\href{mailto:pparent@uottawa.ca}{pparent@uottawa.ca}\par}
\end{center}

\vfill
{\small
\noindent\textit{2020 Mathematics Subject Classification.}\enspace
Primary 55P62, 55M30; Secondary 13D03, 18G55.

\medskip
\noindent\textit{Key words and phrases.}\enspace
Homology nilpotency, homotopical nil-length, rational
Lusternik--Schnirelmann category, Sullivan model, Andr\'e--Quillen homology,
retractive tower.\par}

\vspace{1pc}
{\footnotesize
\noindent\textsuperscript{*}This research was supported by the Natural
Sciences and Engineering Research Council of Canada (NSERC), grant
no.~149521.\par}
\end{titlepage}

\thispagestyle{plain}
\null
\vspace*{3pc}
\begin{center}
{\Large\bfseries Abstract\par}
\end{center}
\vspace{1.5pc}
\noindent
Homology nilpotency of a simply connected minimal Sullivan algebra \(M\)
is the least integer \(n\) for which \((M^+)^{n+1}\) is contained in an
acyclic differential ideal.  Unlike homotopical nil-length, homology
nilpotency imposes a rigid quotient condition on the fixed minimal model.  We
develop lower- and upper-bound criteria, including methods for
rational cell attachments, and construct an explicit finite-type minimal
Sullivan algebra \(M\) satisfying
\(\nilh(M)=3<\Hnil(M)=4\).

\par\medskip\noindent
Complete strictly coordinated retractive towers identify the obstruction:
they characterize homotopical nil-length while exposing the additional
quotient defect measured by homology nilpotency.  We then show that this rigid
defect disappears after stabilization by
degreewise finite wedges of simply connected rational spheres:
\[
                  \sHnil(X)=\catz(X)
\]
for every simply connected rational space \(X\) of finite type.  For the
Sullivan realization \(|M|\) of this algebra, the stabilization is explicit:
a single rational sphere suffices, and
\[
       \Hnilz(|M|\vee S^{17}_{\Q})=\catz(|M|)=3
       <4=\Hnilz(|M|).
\]

\par\vfill
\noindent\begin{minipage}{\textwidth}
\phantomsection
\addcontentsline{toc}{chapter}
  {Declaration on the Use of Generative Artificial Intelligence}
\begin{center}
{\large\bfseries
Declaration on the Use of Generative Artificial Intelligence\par}
\end{center}
\vspace{0.6pc}
\noindent
During the preparation of this manuscript, the author used ChatGPT (OpenAI)
as an assistive tool for language editing, \LaTeX{} preparation, and the
verification and exploration of mathematical arguments. All mathematical
statements, proofs, references, and conclusions were independently reviewed
by the author, who assumes full responsibility for the content of the
manuscript.\par
\end{minipage}
\clearpage

\thispagestyle{plain}

\null
\vspace*{2.7pc}

\begin{center}
  \begin{minipage}{0.78\textwidth}
    {\Large\itshape À Andrée Martineau\ldots\par}
    \vspace{0.8pc}
    {\large\itshape\hfill\ldots\, mon étoile filante\par}
  \end{minipage}

  
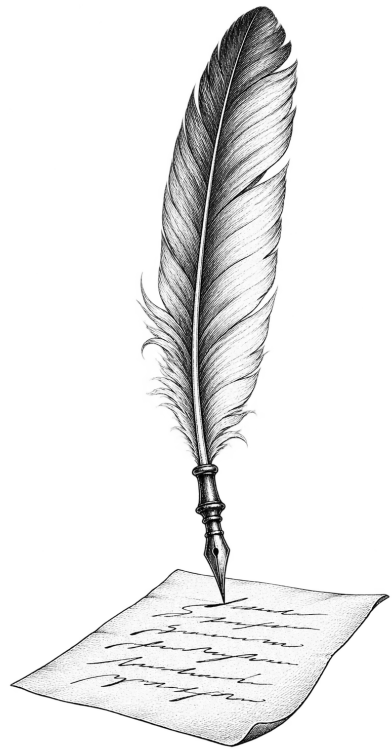
\vspace{2.8pc}

  \EmbeddedQuillManuscript

  \vspace{1.5pc}

  {\large\scshape\textls[95]{Citoyenne du monde}\par}
\end{center}

\vfill

\clearpage
\tableofcontents

\mainmatter

\chapter*{Introduction}
\markboth{INTRODUCTION}{INTRODUCTION}

Many useful invariants in topology arise from a tension between two ways of
simplifying an object.  One may replace a space or an algebra by any equivalent
model and ask how short that model can be.  Or one may keep a preferred model
and ask whether the same simplification can be carried out inside it by an
actual quotient.  The first viewpoint is homotopical; the second is rigid.
This memoir begins where their expected agreement fails, then asks what the
failure measures and how geometry can repair it.

The classical setting is Lusternik--Schnirelmann category, which measures how
economically a space can be covered by subsets contractible within it and is
closely related to construction-stage invariants such as cone length.  With
roots in critical-point theory, it now belongs to a
wider circle of ideas involving sectional category and cone decompositions
\cite{Ganea1960,CorneaLuptonOpreaTanre2003}.  Its rational theory is especially
suited to comparing flexible and rigid constructions.  Sullivan's
correspondence replaces a simply connected rational space by a minimal
commutative differential graded algebra, unique up to isomorphism
\cite{Sullivan1977,BousfieldGugenheim1976,FHT2001}.  Geometric construction
problems then become algebraic questions about differentials, ideals,
retractions, and multiplicative length.

Two such lengths drive this memoir.  The homotopical nil-length \(\nilh\) asks
for the shortest commutative differential graded algebra among all those
quasi-isomorphic to the minimal Sullivan model.  Homology nilpotency \(\Hnil\),
introduced in the ideal-theoretic study of rational sectional category
\cite{CarrasquelVera2015}, imposes a stricter demand: the short algebra must be
an actual quotient of the fixed minimal model by an acyclic differential
ideal.  Thus \(\nilh\) permits a flexible replacement, whereas \(\Hnil\)
requires a rigid quotient.  For a simply connected rational space \(X\), with
minimal model \(M_X\), the basic comparison is
\[
       \catz(X)\leq \Clz(X)=\nilh(M_X)\leq \Hnilz(X).
\]
The middle term is rational cone length
\cite{Cornea1994,CorneaNilLength1994}.  Every rigid quotient is a flexible
witness; the converse is a genuine strictification problem.

Minimal Sullivan models are free, acyclic extensions are abundant, and the two
lengths agree through the first two levels.  Yet multiplication, primitives,
and closure under ideals need not strictify simultaneously.  A boundary
harmless in a flexible model may, in the minimal model, force a primitive into
an acyclic ideal; multiplying that primitive can then create essential
cohomology.  Homology nilpotency must therefore be computed independently.

The first aim of the memoir is to make such calculations possible.  We develop
lower- and upper-bound methods for \(\Hnil\), including a quadratic-shadow
bound and criteria for rational cell attachments that replace an infinite
ideal search by filtered or finite-degree tests.  As a first test of this
calculus, Part~\ref{part:calculations} applies the essential-direction
obstruction and an independent quadratic argument to a cell attachment, where
the two bounds meet at the exact value four.

Retractive towers organize the comparison.  They repeatedly combine intrinsic
short quotients, relative Sullivan extensions, strict retractions, and
pushouts.  The existence of a complete retractive tower characterizes
\(\nilh(M)\).  In a strictly coordinated tower, the images in \(M\) form an
increasing union of differential ideals.  The cohomology of that union records
whether the flexible object obtained at the limit descends to an acyclic short
quotient of the original model.  Thus the tower formalism locates, rather than
merely names, the strictification defect.

That defect becomes genuine at augmentation length three, the first possible
level after low-length strictification.  We construct an explicit simply
connected minimal Sullivan algebra \(M\) of finite type for which
\[
                         \nilh(M)=3<\Hnil(M)=4.
\]
The strict inequality is forced by a boundary that has a unique
primitive, while ideal closure forces a nonzero cohomology class.  The same
mechanism occurs in every complete strictly coordinated level-three tower.
The construction therefore explains where strictification fails and why the
low-length argument cannot simply continue.

The negative result leads to a positive one.  Although the rigid defect can
survive in the fixed minimal model, it disappears after a controlled
stabilization.  If \(\sHnil(X)\) denotes the infimum of the homology-nilpotency
values obtained by wedging \(X\) with admissible degreewise finite wedges of
simply connected rational spheres, then
\[
                              \sHnil(X)=\catz(X).
\]
The stabilization uses a specifically chosen sphere wedge.  A fibre--cofibre
construction associated with an optimal Sullivan representative of rational
LS category---the LS root used later---produces such a wedge, whose homology
nilpotency reaches the lower bound
\cite{FelixThomas1988}.  The added spherical directions provide the room needed
to display homotopical simplicity as a quotient.  The proof combines this
geometry with homotopy transfer and Harrison--Andr\'e--Quillen theory
\cite{Harrison1962,Andre1974,Quillen1970AQ,GinzburgKapranov1994}, connecting
rational LS category, commutative differential graded algebra, and the broader
problem of rectifying weak data inside a preferred model.

For the Sullivan realization \(|M|\) of the algebra above, the repair is
considerably sharper.  Its three-short witness contains a hidden contractible
pair \(r_{16}\mapsto xt_{17}\).  Removing the primitive direction exposes the
degree-\(17\) class and produces a three-short quotient of the minimal model of
the wedge.  Consequently, a single sphere suffices:
\[
        \Hnilz(|M|\vee S^{17}_{\Q})=\catz(|M|)=3
        <4=\Hnilz(|M|).
\]
Thus the negative and positive halves of the memoir meet in the same local
mechanism: the missing disk direction responsible for the unstabilized defect
determines the explicit sphere that removes it.

The chapters follow one question: when can a homotopical simplification be
made visible as an actual quotient of the minimal model?  They move from
calculational tools to a universal obstruction and finally to its geometric
repair.

\section*{Acknowledgments}

The author thanks Daniel Tanr\'e for many stimulating and enriching discussions
over the years, both online and in Lille, France.  He also thanks Jonathan Scott
for introducing him to \(C_\infty\)-structures.

\part{Foundations and the Two Nilpotency Invariants}\label{part:foundations}

\chapter{Main Theorems and Proof Architecture}\label{chap:main-results}

\begin{chapterguide}[title={Chapter guide}]
This chapter is a map, not a second introduction.  It contrasts the flexible
and rigid invariants, states the two principal theorems, records the
dependence of their proofs, and isolates the exact limits of the conclusions.
\end{chapterguide}

\section{Flexible replacement and rigid quotient}

Let \(X\) be a simply connected rational space of finite type and let
\(M_X=(\Lambda V,d)\) be its minimal Sullivan model.  The definitions in
Chapter~\ref{stab:sec:alg-prelim} give
\begin{equation}\label{intro:eq:basic-chain}
 \catz(X)\leq \Clz(X)=\nilh(M_X)\leq\Hnilz(X).
\end{equation}
The middle equality is Cornea's identification of rational cone length with
homotopical nil-length \cite{Cornea1994,CorneaNilLength1994}.

\begin{center}
\renewcommand{\arraystretch}{1.22}
\begin{tabularx}{0.96\textwidth}{@{}>{\raggedright\arraybackslash\bfseries}p{0.20\textwidth}
  >{\raggedright\arraybackslash}X>{\raggedright\arraybackslash}X@{}}
\toprule
 & Flexible problem & Rigid problem \\
\midrule
Witness
 & Any \(n\)-short CDGA \(B\simeq M_X\)
 & A quotient \(M_X\twoheadrightarrow M_X/J\), with \(J\) acyclic and the
   quotient \(n\)-short \\
Invariant
 & \(\nilh(M_X)=\Clz(X)\)
 & \(\Hnilz(X)=\Hnil(M_X)\) \\
Permitted change
 & Replace the CDGA representative
 & Choose only an acyclic ideal in the fixed minimal model \\
\bottomrule
\end{tabularx}
\end{center}

The last inequality can be strict.  The memoir therefore treats \(\Hnil\) as
an invariant requiring its own calculation, not as another notation for
rational cone length.

\section{The two principal theorems}

The first result identifies the initial length at which flexible replacement
need not strictify to a quotient.

\begin{theorem}[Divergence at length three]
\label{div:thm:intro-main}
There exists a simply connected minimal Sullivan algebra \(M\) of finite type
such that
\[
                         \nilh(M)=3<\Hnil(M)=4.
\]
Consequently, its Sullivan realization \(|M|\) is a simply connected rational
space of finite rational type satisfying
\[
                         \Clz(|M|)=3<\Hnilz(|M|)=4.
\]
\end{theorem}

The separating model is built from a finite three-short witness.  Its minimal
model contains a boundary with a literally unique primitive; ideal closure of
that primitive forces an essential cohomology class.  Chapters
\ref{div:chap:witness}--\ref{chap:strict-separation} give the lower bound,
while Appendix~\ref{app:exact-Hnil} constructs the acyclic ideal proving the
matching upper bound.  The same dichotomy obstructs every complete strictly
coordinated level-three tower.

The second result shows that the rigid defect disappears after a controlled
spherical stabilization.

\begin{theorem}[Spherical stabilization]
\label{stab:thm:intro-main}
Let \(X\) be a simply connected rational space of finite type with
\(n=\catz(X)<\infty\).  If \(n=0\), set \(G_\sigma=\ast\) and
\(W_\sigma\simeq_\Q X\).  If \(n\geq1\), choose the section induced by a
strict retraction of an optimal normalized Sullivan \(n\)-LS root, let
\(G_\sigma\) be its homotopy fibre, and put
\[
 W_\sigma=\hocof(G_\sigma\longrightarrow X).
\]
Then
\[
 W_\sigma\simeq_\Q X\vee\Sigma G_\sigma,
 \qquad \Hnilz(W_\sigma)=n.
\]
Moreover, \(\Sigma G_\sigma\) is rationally a wedge of simply connected
spheres with finitely many summands in every degree.  Consequently,
\[
                         \sHnil(X)=\catz(X).
\]
\end{theorem}

For the separating space the stabilization is finite and explicit:
\[
       \Hnilz(|M|\vee S^{17}_{\Q})=\catz(|M|)=3
       <4=\Hnilz(|M|).
\]
One sphere is minimal as a number of summands; no claim is made that degree
\(17\) is the smallest possible dimension.

\section{Proof-dependence map}

The scientific dependence of the memoir is summarized below.  In particular,
the explicit separation does not use the Harrison--Andr\'e--Quillen chapter,
whereas the general stabilization theorem does.

\begin{center}
\small
\renewcommand{\arraystretch}{1.18}
\begin{tabularx}{\textwidth}{@{}>{\raggedright\arraybackslash\bfseries}p{0.18\textwidth}
  >{\raggedright\arraybackslash}p{0.29\textwidth}
  >{\raggedright\arraybackslash}X@{}}
\toprule
Block & Principal results & Role and dependence \\
\midrule
Foundations
 & Definitions and quotient criterion in Chapter~\ref{stab:sec:alg-prelim}
 & Fixes the two invariants and reduces \(\Hnil\) to a surjective
   quasi-isomorphism problem. \\
Calculation
 & Chapters~\ref{calc:chap:acyclic}--\ref{chap:why-close}
 & Gives the quadratic bound, essential directions, cell-attachment criteria,
   and strictification through length two. \\
Towers
 & Chapters~\ref{rt:chap:foundations}--\ref{chap:tate-selection}
 & Characterizes \(\nilh\), constructs the comparison defect, and solves the
   finite Tate-selection problem after a retraction is fixed. \\
Separation
 & Chapters~\ref{div:chap:witness}--\ref{chap:strict-separation} and
   Appendix~\ref{app:exact-Hnil}
 & Proves \(\nilh(M)=3<\Hnil(M)=4\); independent of
   Chapter~\ref{chap:aq-toolkit}. \\
Operadic bridge
 & Chapter~\ref{chap:aq-toolkit}
 & Types the augmented/nonunital passage, the Harrison bar--cobar resolution,
   derived indecomposables, and the mixed edge map. \\
Stabilization
 & Chapters~\ref{chap:construction}--\ref{chap:stabilization-proof}
 & Uses the operadic bridge for the general open-disk surjectivity theorem;
   the explicit one-sphere calculation is also verified directly. \\
\bottomrule
\end{tabularx}
\end{center}

\section{Exact scope}

Four limitations are part of the theorem statements.
\begin{enumerate}[label=\textup{(\arabic*)},leftmargin=2.5em]
\item The tower characterization is existential: one complete branch is
      enough, and failure of a selected branch is not an invariant lower
      bound.
\item Tate selection is exact only after the short witness, branch, and
      retraction have been fixed; an arbitrary acyclic quotient is not claimed
      to generate an asymptotically Tate tower.
\item The general stabilization uses a degreewise finite wedge of spheres,
      which need not have finitely many summands in total.
\item The priority statement concerns only the augmentation ideal of a fixed
      minimal Sullivan model and is deliberately qualified; see
      Chapter~\ref{chap:precedents}.
\end{enumerate}

All CDGAs are cohomologically graded over \(\Q\), their differentials have
degree \(+1\), and structural morphisms preserve augmentations.  The complete
algebraic and geometric conventions are in Chapters
\ref{stab:sec:alg-prelim} and \ref{stab:sec:geom-prelim}.

\chapter{Algebraic Preliminaries}\label{stab:sec:alg-prelim}
\begin{chapterguide}[title={Chapter guide}]
This is the memoir's principal repository for algebraic conventions.  It fixes
Sullivan models, finite-type hypotheses, augmentation length, short CDGAs,
the two nilpotency invariants, the quotient criterion, and the lifting tools
used later.  Cohomological degree and augmentation length are kept visibly
distinct throughout.
\end{chapterguide}
\index{commutative differential graded algebra (CDGA)}
\index{Sullivan algebra}
\index{relative Sullivan algebra}
\index{augmentation ideal}
\index{augmentation length}
\index{short CDGA}
\index{acyclic differential ideal}
\index{indecomposables}

\section{Connected CDGAs and Sullivan algebras}
All vector spaces and algebras are over \(\Q\).  Gradings are
cohomological and differentials have degree \(+1\).  A \emph{CDGA} is a
unital graded-commutative differential algebra.  It is \emph{connected}
if \(A^0=\Q\), and \emph{simply connected} if in addition \(A^1=0\).
All CDGAs are augmented; \(A^+=\ker(\varepsilon:A\to\Q)\).
\index{commutative differential graded algebra (CDGA)}
Throughout, an object is said to be \emph{of finite type} when the
relevant graded vector space is finite dimensional in every degree.  For a
rational space this refers to rational cohomology, and for a Sullivan algebra
to its space of generators; when a CDGA itself is said to be of finite type,
the convention applies to its underlying graded vector space.
The indecomposables are
\[
 Q(A)=A^+/(A^+)^2.
\]
\index{indecomposables}
The induced differential on \(Q(A)\) is denoted by \(Qd\).  The symbol
\(\Lambda V\) denotes the free graded-commutative algebra on the graded
vector space \(V\).  For a graded vector space \(W\) of finite type,
\[
 W^\sharp=\bigoplus_k\operatorname{Hom}_{\Q}(W^k,\Q)
\]
always denotes the degreewise graded dual; no completed or unrestricted dual
is used below.

\begin{definition}\label{alg:def:sullivan-algebra}
A Sullivan algebra is a CDGA \((\Lambda V,d)\) for which \(V\) has a
well-ordered homogeneous basis \(\{v_\alpha\}\) such that
\(dv_\alpha\in\Lambda V_{<\alpha}\).  It is minimal if
\(dV\subset\Lambda^{\geq2}V\).  A relative Sullivan algebra under a
CDGA \(B\) is an inclusion
\index{Sullivan algebra}
\[
 B\longrightarrow (B\otimes\Lambda Z,D)
\]
with the analogous well-ordering condition on \(Z\).
\end{definition}

We use the standard model structure on connected CDGAs
\cite[Chapter~4]{BousfieldGugenheim1976}: relative Sullivan extensions are
cofibrations, surjections are fibrations, and quasi-isomorphisms are weak
equivalences.  Only the following elementary lifting consequence is needed.

\begin{lemma}[Sullivan lifting]\label{stab:lem:sullivan-lifting}
Let \(j:A\to A\otimes\Lambda Z\) be a relative Sullivan extension and
let \(p:E\twoheadrightarrow B\) be a surjective quasi-isomorphism.  Every
commutative square
\[
\begin{tikzcd}[column sep=large,row sep=large]
A \arrow[r] \arrow[d,"j"'] & E \arrow[d,"p",two heads,"\sim"']\\
A\otimes\Lambda Z \arrow[r] & B
\end{tikzcd}
\]
admits a lift \(A\otimes\Lambda Z\to E\).  The lift can be chosen
inductively in the Sullivan order.  The same statement holds in the
underlying category of unital nonnegatively graded CDGAs, without
connectedness or augmentation hypotheses.
\end{lemma}

\begin{proof}
Suppose the lift has been defined on the generators preceding
\(z_\alpha\).  Choose any preimage in \(E\) of the prescribed image of
\(z_\alpha\).  Its differential differs from the already prescribed
image of \(Dz_\alpha\) by a cocycle in \(\ker p\).  Since \(p\) is a
surjective quasi-isomorphism, \(\ker p\) is acyclic; correcting the chosen
preimage by a primitive in \(\ker p\) gives the required value.  Limit
stages are unions.
\end{proof}

\begin{remark}[Derived inverses and Sullivan resolutions]
\label{stab:rem:derived-inverses}
No strict inverse is intended when a quasi-isomorphism is inverted below.
This use of the homotopy category amounts to making the usual semifree
Sullivan replacements explicit.  Indeed, let
\(f:N\to C\) be a quasi-isomorphism, let \(M\) be a Sullivan algebra, and
let \(\eta:M\to C\) be a morphism.  Factor \(f=p i\), where
\(i:N\to E\) is an acyclic relative Sullivan cofibration and
\(p:E\twoheadrightarrow C\) is a surjective quasi-isomorphism.  To spell
out the fibrancy argument, the augmentation
\(\varepsilon_N:N\twoheadrightarrow\Q\) is a fibration, and the lifting
property of the acyclic cofibration \(i\) gives a diagonal in
\[
\begin{tikzcd}[column sep=large,row sep=large]
N \arrow[r,"\id_N"] \arrow[d,"i"']
  & N \arrow[d,"\varepsilon_N",two heads]\\
E \arrow[r,"\varepsilon_E"'] \arrow[ur,dashed,"r"] & \Q .
\end{tikzcd}
\]
Thus \(r i=\id_N\).  Since \([i]\) is invertible,
\([r]=[i]^{-1}\) and \([i r]=[\id_E]\) in the homotopy category; this
is the only homotopy-category equality used here.  Since \(M\) is cofibrant,
Sullivan lifting gives \(\widetilde\eta:M\to E\) with
\(p\widetilde\eta=\eta\).  Thus
\[
 j=r\widetilde\eta:M\longrightarrow N,
 \qquad [fj]=[\eta],
\]
represents \([f]^{-1}[\eta]\).  Hence localization at
quasi-isomorphisms introduces no additional algebraic assumption; it merely
suppresses the choice of a Sullivan resolution and of its lifts.  It does
not, however, turn a homotopy equality into a strict one, nor does it
automatically preserve surjectivity, kernels, quotients, or strict
retractions.  Whenever such point-set properties are used below, they are
established separately by lifting or strictification.
\end{remark}

The following decomposition is the precise form of the usual
minimal--contractible splitting that will be used later.  It is included
to make the disk argument independent of any textbook formulation.

\begin{proposition}[Minimal--contractible decomposition]
\label{stab:prop:min-contractible}
Let \((\Lambda Z,D)\) be a simply connected Sullivan algebra of finite
type.  There is a triangular CDGA isomorphism
\[
 (\Lambda Z,D)\cong
 (\Lambda Z_{\min},D)\otimes
 \bigl(\Lambda(U\oplus\widehat U),D\bigr),
\]
where the first factor is minimal and
\[
 Du=\widehat u,\qquad D\widehat u=0.
\]
Here \(\widehat u\in\widehat U^{|u|+1}\) is the generator paired with
\(u\in U\); the hat denotes the partner generator, not an operation on
\(u\).
The minimal factor is unique up to CDGA isomorphism.  Moreover,
\[
 H(Q(\Lambda Z),QD)\cong Z_{\min}.
\]
If \(M\) is a minimal Sullivan subalgebra of \(\Lambda Z\) and the inclusion
is a quasi-isomorphism admitting a CDGA retraction, then the splitting may be
chosen under \(M\) with
minimal factor equal to \(M\).
\end{proposition}

\begin{proof}
Write \(D_1\) for the linear part of \(D\) on \(Z\).  Choose, degree by
degree and compatibly with the Sullivan order, a vector-space
decomposition
\[
 Z=H\oplus U\oplus\widehat U
\]
such that \(D_1|_H=0\), \(D_1:U\to\widehat U\) is an isomorphism, and
the classes of \(H\) form a basis of \(H(Z,D_1)\).  After rescaling the
bases, write \(D_1u=\widehat u\).

All changes below are made along the exhaustive filtration by initial
segments of the chosen Sullivan order.  A successor change uses only the
subalgebra generated at earlier stages, and at a limit ordinal the changed
algebra is the union of the preceding ones.  Thus every substitution is
triangular for one fixed exhaustive filtration and the resulting composite
is a well-defined CDGA automorphism.

Proceed in the Sullivan order.  Suppose the preceding generators have
already been changed so that the corresponding disk pairs satisfy
\(Du=\widehat u\) and \(D\widehat u=0\).  Replace the current
\(\widehat u\) by the exact element \(Du\).  Its linear term is the old
\(\widehat u\), so this is a triangular change of variables; the identity
\(D^2=0\) gives \(D\widehat u=0\).  For a generator in \(H\), subtract
successively terms involving the already split disks, using the
contraction \(u\mapsto0\), \(\widehat u\mapsto u\).  The resulting
differential is decomposable in the remaining \(H\)-generators.  This
produces the asserted tensor decomposition.  The last displayed
isomorphism follows because the indecomposable complex is the direct sum
of the zero complex on \(Z_{\min}\) and the disk complexes
\(U\to\widehat U\).

A quasi-isomorphism between minimal Sullivan algebras is an isomorphism:
induction on degree and the Sullivan order shows first that it is an
isomorphism on indecomposables and then on the free algebras.  Hence the
minimal factor is unique.  In the relative situation, first replace every
relative generator \(z\) by \(z-r(z)\), where \(r:\Lambda Z\to M\) is the
given retraction.  This is an invertible change under \(M\), and all new
relative generators belong to \(\ker r\).  Choose the splitting of the
linear relative complex inside \(\ker r\); the substitutions
\(\widehat u\mapsto Du\) and all later disk corrections then also lie in
\(\ker r\).  Thus every change fixes \(M\) and the resulting decomposition
is an isomorphism under \(M\), not merely an absolute isomorphism.  The
linear differential on every remaining relative indecomposable is zero;
together with the minimality of \(M\), this makes the remaining factor
absolutely minimal, not merely minimal relative to \(M\).  The inclusion of \(M\) into that remaining minimal factor is a
quasi-isomorphism: after tensoring with the contractible disk factor it is
the original quasi-isomorphism \(M\hookrightarrow\Lambda Z\).  It is therefore
a quasi-isomorphism between minimal Sullivan algebras and hence an
isomorphism.  Composing with this isomorphism identifies the remaining
minimal factor literally with \(M\), under \(M\).  This is the relative form
of the standard
Sullivan decomposition theorem; compare
\cite[Theorems~14.9 and~14.11]{FHT2001} and
\cite{Halperin1977,Sullivan1977}.
\end{proof}

We next record a normalized form of the surjective minimal-model lemma.  Its
proof is included because the normalization is used both below and in
Chapter~\ref{chap:why-close}.

\begin{lemma}[Surjective minimal model of a decomposable quotient]
\label{stab:lem:decomposable-quotient-model}
Let
\[
 p:(\Lambda U,d)\twoheadrightarrow(A,d)
\]
be a surjective morphism of simply connected augmented CDGAs.  Assume that
\((\Lambda U,d)\) is minimal and
\[
 \ker p\subseteq(\Lambda^+U)^2=\Lambda^{\geq2}U.
\]
Then \(p\) admits a factorization
\begin{equation}\label{stab:eq:decomposable-quotient-factorization}
 (\Lambda U,d)\lhook\joinrel\longrightarrow
 (\Lambda(U\oplus Z),D)
 \xrightarrow[\simeq]{\theta}\!\!\!\twoheadrightarrow(A,d)
\end{equation}
such that the first arrow is a relative minimal Sullivan extension, the
total Sullivan algebra is minimal, \(\theta|_{\Lambda U}=p\), the relative
generators may be chosen with \(\theta(Z)=0\), and \(\ker\theta\) is an
acyclic differential ideal.
\end{lemma}

\begin{proof}
Apply \cite[Proposition~14.3]{FHT2001} and then the standard relative
minimization to obtain a relative minimal Sullivan factorization of \(p\)
whose second map \(\theta\) is a surjective quasi-isomorphism.  Surjectivity
also follows directly from the usual inductive construction, since
\(\theta|_{\Lambda U}=p\).

Choose a homogeneous basis of \(Z\) compatible with the Sullivan ordering.
For each \(z\in Z\),
choose \(a_z\in\Lambda U\) with \(p(a_z)=\theta(z)\), and replace \(z\) by
\(z-a_z\).  This is a triangular change under \(\Lambda U\).  Moreover,
minimality of the base implies that \(d a_z\) has no linear part.  Induction
in the Sullivan order therefore gives
\[
                         \theta(Z)=0.
\]

Let \(D_1\) be the linear part of \(D\).  Relative minimality implies
\(D_1z\in U\).  From \(\theta(z)=0\) and \(\theta D=d\theta\), passage to
indecomposables gives
\[
                         Q(p)(D_1z)=0.
\]
The map \(Q(p):U\to Q(A)\) is injective.  Indeed, if \(p(u)\) is
decomposable, lift its factors through the surjection \(p\) to obtain
\(w\in(\Lambda^+U)^2\) with \(p(w)=p(u)\).  Then
\(u-w\in\ker p\subseteq(\Lambda^+U)^2\), so \(u=0\) in
\(Q(\Lambda U)=U\).  Thus \(D_1z=0\) for every \(z\).  Since the base is
minimal, the total algebra in
\eqref{stab:eq:decomposable-quotient-factorization} is minimal.

Finally, the short exact sequence
\[
 0\longrightarrow\ker\theta\longrightarrow
 \Lambda(U\oplus Z)\xrightarrow{\theta}A\longrightarrow0
\]
and the fact that \(\theta\) is a quasi-isomorphism imply
\(H^*(\ker\theta)=0\).
\end{proof}

\section{Augmentation length and short CDGAs}

If \(A\) is augmented, write
\[
 \mathfrak m_A=A^+=\ker(A\longrightarrow\Q)
\]
for its augmentation ideal.  The notation \(A^q\) always denotes
cohomological degree, whereas \(\mathfrak m_A^k\) denotes augmentation
length.  These are independent filtrations.
\index{augmentation ideal}
\index{augmentation length}

\begin{definition}\label{def:book-short}
An augmented CDGA \(A\) is \emph{\(n\)-short} if
\[
 \mathfrak m_A^{n+1}=0.
\]
Thus a one-short algebra has square-zero augmentation ideal, a two-short
algebra has cube-zero augmentation ideal, and a three-short algebra has
fourth power zero.
\index{short CDGA}
\end{definition}

For an augmented morphism \(f:A\to B\), we write
\[
 Q(f):Q(A)\longrightarrow Q(B)
\]
for the induced map on indecomposables.

\section{The two nilpotency invariants}

\begin{definition}
\label{div:def:nilh}
The \emph{homotopical nil-length} of an augmented CDGA \(A\) is
\index{homotopical nil-length}
\[
\nilh(A)=\inf\left\{n\geq0\ \middle|\
 \begin{array}{c}
 \text{\(A\) is connected to an \(n\)-short CDGA by a zigzag}\\
 \text{of augmented quasi-isomorphisms}
 \end{array}\right\}.
\]
\end{definition}

Let \(I\) be a differential ideal in a CDGA \(A\).  Carrasquel-Vera
defined its homology nilpotency by
\index{homology nilpotency}
\index{acyclic differential ideal}
\[
\Hnil(I;A)=\inf\left\{n\geq0\ \middle|\
 I^{n+1}\subseteq J
 \text{ for some acyclic differential ideal }J\lhd A\right\};
\]
see \cite[Definition~20]{CarrasquelVera2015}.

\begin{definition}
\label{div:def:Hnil}
If \(M=(\Lambda V,d)\) is a simply connected minimal Sullivan algebra,
we set
\[
\begin{aligned}
\Hnil(M)&=\Hnil(\mM;M)\\
&=\inf\left\{n\geq0\ \middle|\
 \begin{gathered}
 \mM^{n+1}\subseteq J\text{ for some}\\[-0.2ex]
 \text{acyclic differential ideal }J\lhd M
 \end{gathered}\right\}.
\end{aligned}
\]
\end{definition}

Throughout, the infimum of the empty set is understood to be
\(\infty\).  With these conventions, \(\Q\) and the one-point space
have length zero.  Conversely, if a simply connected minimal Sullivan
algebra \(M\) satisfies \(\nilh(M)=0\), then \(M\) is quasi-isomorphic
to the only zero-short CDGA, namely \(\Q\).  Uniqueness of minimal
Sullivan models gives \(M\cong\Q\), and hence \(\Hnil(M)=0\).

For a simply connected space \(X\) of finite rational type, let
\(M_X\) denote the minimal Sullivan model of its rationalization.  We
use the notation
\[
             \nilh(X)=\nilh(M_X),\qquad
             \Hnilz(X)=\Hnil(M_X).
\]
More generally, if a simply connected augmented CDGA \(A\) admits a
finite-type minimal Sullivan model \(M_A\), we write
\(\Hnilz(A)=\Hnil(M_A)\).
This is well defined because the minimal Sullivan model is unique up to
isomorphism.  The subscript in \(\Hnilz\) records that the invariant is
rational.  The rational cone length \(\Clz(X)\) is defined geometrically in
Chapter~\ref{stab:sec:geom-prelim},
Section~\ref{sec:geometric-invariants}.

\begin{remark}
The invariant \(\nilh\) is Cornea's nil-length, for which we use the
descriptive term \emph{homotopical nil-length}.  It is unrelated to the
Berstein--Ganea homotopy nilpotency of an \(H\)-space
\cite{BersteinGanea1961}.
\end{remark}

\begin{lemma}[Direct representative of a zigzag]
\label{div:lem:direct-zigzag}
Let \(M\) be a simply connected minimal Sullivan algebra and let \(B\)
be an augmented CDGA.  If \(M\) and \(B\) are connected by a zigzag of
augmented quasi-isomorphisms, then there is an augmented
quasi-isomorphism
\[
                              M\longrightarrow B.
\]
\end{lemma}

\begin{proof}
Start with \(\id_M:M\to M\) and move along the finite
zigzag.  If an arrow points away from the current algebra, compose with
it.  If a quasi-isomorphism \(A_{i+1}\to A_i\) points toward the current
algebra and a quasi-isomorphism \(M\to A_i\) has already been obtained,
the Sullivan lifting theorem \cite[Proposition~12.9]{FHT2001}, applied
to the quasi-isomorphism \(A_{i+1}\to A_i\), supplies an augmented
morphism \(M\to A_{i+1}\) whose composite with
\(A_{i+1}\to A_i\) is homotopic to the given map.  The new morphism is
a quasi-isomorphism by the two-out-of-three property.  Iteration gives
the required quasi-isomorphism \(M\to B\).  Since \(M^0=\Q\) and every
morphism under consideration preserves units and degrees, the resulting
morphism is automatically augmented.
\end{proof}

Thus, whenever the source is a minimal Sullivan algebra, we use the
direct-map formulation supplied by Lemma~\ref{div:lem:direct-zigzag}.

\begin{lemma}[Surjective-witness characterization]
\label{div:lem:surjective-characterization}
Let \(M\) be a simply connected minimal Sullivan algebra and let
\(n\geq0\).  The following statements are equivalent.
\begin{enumerate}[label=\textup{(\roman*)}]
\item \(\Hnil(M)\leq n\).
\item There exists a surjective augmented quasi-isomorphism
      \(M\onto C\) with \(C\) \(n\)-short.
\end{enumerate}
\end{lemma}

\begin{proof}
Assume first that \(\mM^{n+1}\subseteq J\), where \(J\lhd M\) is an
acyclic differential ideal.  The short exact sequence
\[
                         0\longrightarrow J\longrightarrow M
                         \longrightarrow M/J\longrightarrow0
\]
and the long exact cohomology sequence show that
\(M\onto M/J\) is a quasi-isomorphism.  Since the image of
\(\mM^{n+1}\) is zero, \(M/J\) is \(n\)-short.

Conversely, let \(p:M\onto C\) be a surjective quasi-isomorphism with
\(\mC^{n+1}=0\).  Its kernel \(J=\ker p\) is an acyclic differential
ideal by the long exact cohomology sequence.  Moreover,
\[
                         p(\mM^{n+1})=\mC^{n+1}=0,
\]
so \(\mM^{n+1}\subseteq J\).  Hence \(\Hnil(M)\leq n\).
\end{proof}

\begin{corollary}
\label{div:cor:basic-inequality}
For every simply connected minimal Sullivan algebra \(M\),
\[
                           \nilh(M)\leq\Hnil(M).
\]
\end{corollary}

\begin{proof}
If \(\Hnil(M)\leq n\), Lemma~\ref{div:lem:surjective-characterization}
provides an \(n\)-short CDGA quasi-isomorphic to \(M\), which is precisely
the condition \(\nilh(M)\leq n\).
\end{proof}

We use throughout the rational-model convention fixed immediately after
Definition~\ref{div:def:Hnil}; it is distinct from the
presentation-dependent ideal invariant \(\Hnil(I;A)\).

\section{Two elementary surjectivity tools}
\begin{lemma}[Surjectivity on indecomposables]\label{stab:lem:Q-surj}
Let \(f:A\to B\) be a morphism of simply connected CDGAs.  If
\(Qf:Q(A)\to Q(B)\) is surjective, then \(f\) is surjective.
\end{lemma}

\begin{proof}
Induct on cohomological degree.  Given \(b\in B^k\), use the
surjectivity of \(Qf\) to subtract an element of \(f(A^k)\) so that the
remainder lies in \((B^+)^2\).  Each factor in a decomposable element of
degree \(k\) has smaller positive degree, so the induction hypothesis
lifts the remainder.  In fact, connectedness and positive grading suffice
for this induction; simple connectivity is retained in the statement
because it is the standing hypothesis and makes all positive
indecomposables start in degree two.
\end{proof}

\begin{lemma}[Products of cocycles]\label{stab:lem:products-cocycles-prelim}
Let \(f:N\qiso A\) be a quasi-isomorphism from a minimal Sullivan
algebra and assume \(Q(d_A)=0\).  Let \(a\in A^+\).  If \(a\) is a
cocycle, or if
\[
 d_Aa=\sum_j x_jy_j
\]
with \(x_j,y_j\in A^+\) cocycles, then the class of \(a\) in \(Q(A)\)
belongs to \(\im Qf\).
\end{lemma}

\begin{proof}
If \(a\) is closed, lift its cohomology class to a cocycle in \(N\).  Its
image under \(f\) differs from \(a\) by a boundary in \(A\), hence by a
decomposable element because \(Q(d_A)=0\).

In the second case choose cocycles \(\widetilde x_j,\widetilde y_j\) in
\(N\) and positive elements \(\alpha_j,\beta_j\in A\) such that
\[
 f(\widetilde x_j)=x_j+d_A\alpha_j,
 \qquad
 f(\widetilde y_j)=y_j+d_A\beta_j.
\]
The required Koszul signs are encoded in the following decomposable element:
\[
 c=\sum_j\bigl((-1)^{|x_j|}x_j\beta_j
        +\alpha_jy_j+\alpha_jd_A\beta_j\bigr).
\]
Indeed, the cocycle identities for \(x_j\) and \(y_j\), together with
the Leibniz rule, give
\[
 d_Ac=\sum_j\bigl(x_jd_A\beta_j
        +(d_A\alpha_j)y_j+(d_A\alpha_j)(d_A\beta_j)\bigr).
\]
Consequently,
\(f(\sum_j\widetilde x_j\widetilde y_j)=d_A(a+c)\).
The product on the left is a cocycle whose cohomology class maps to zero;
it is therefore exact in \(N\), say \(du=\sum_j\widetilde
x_j\widetilde y_j\).  Choose a cocycle \(v\in N\) and an element \(w\in A\)
such that
\[
                 f(u)-(a+c)-f(v)=d_Aw.
\]
Since \(c\) is decomposable and \(Q(d_A)=0\), both \(c\) and \(d_Aw\)
vanish in \(Q(A)\).  Hence the class of \(u-v\) in \(Q(N)\) maps to the
class of \(a\).
\end{proof}

\section{Cone length and the category lower bound}

Cornea's nil-length theorem identifies the flexible algebraic invariant with
the rational cone length defined in Chapter~\ref{stab:sec:geom-prelim},
Section~\ref{sec:geometric-invariants}:
\begin{equation}\label{stab:eq:nilh-cone-length}
 \nilh(Y)=\Clz(Y)
\end{equation}
for simply connected rational spaces of finite type
\cite{CorneaNilLength1994,FHT2001}.

Let \(M=(\Lambda V,d)\) be minimal and put
\[
 A_m=M/(M^+)^{m+1},\qquad q_m:M\twoheadrightarrow A_m.
\]
The F\'elix--Halperin characterization states that \(\catz(|M|)\leq m\)
if and only if \(q_m\) admits a retraction in the homotopy category of
connected CDGAs, or equivalently, if a relative Sullivan factorization of
\(q_m\) admits a strict retraction
\cite[Theorem~4.7]{FelixHalperin1982}.

The first inequality in the next proposition is also the specialization
of Carrasquel-Vera's general sectional-category estimate to the
augmentation of a minimal Sullivan model
\cite[Proposition~21]{CarrasquelVera2015}.  The proof records the short-model
argument through the F\'elix--Halperin criterion because that form will be
used again.

\begin{proposition}\label{stab:prop:cat-lower}
For every simply connected rational space \(Y\) of finite type,
\[
 \catz(Y)\leq\Clz(Y)=\nilh(Y)\leq\Hnilz(Y).
\]
\end{proposition}

\begin{proof}
The equality is \eqref{stab:eq:nilh-cone-length}, and the final inequality is
Corollary~\ref{div:cor:basic-inequality}; only the first inequality remains to
be proved.  Suppose a minimal model
\(M\) of \(Y\) is quasi-isomorphic to an \(m\)-short CDGA \(B\).  After
replacing a zigzag by a morphism from the cofibrant algebra \(M\), choose
a quasi-isomorphism \(f:M\to B\).  Since \((B^+)^{m+1}=0\), the morphism
\(f\) factors through \(q_m\):
\[
 M\xrightarrow{q_m}A_m\xrightarrow{\bar f}B.
\]
In the homotopy category, \(f\) is invertible; hence
\(f^{-1}\bar f\) is a homotopy retraction of \(q_m\).  The
F\'elix--Halperin criterion gives \(\catz(Y)\leq m\).
\end{proof}

\section{A low-degree minimal-model lemma}
We use the following standard consequence of the inductive
construction of a Sullivan minimal model; see the construction in the
proof of \cite[Proposition~12.2]{FHT2001} and the relative existence
theorem \cite[Theorem~14.12]{FHT2001}.

\begin{lemma}
\label{div:lem:minimal-extension}
Let \(P=(\Lambda V,d)\) be a simply connected minimal Sullivan algebra
whose generators satisfy \(V^q=0\) for \(q>N\).  Let \(A\) be a
connected CDGA and let \(g:P\to A\) be a morphism.  Suppose that, for
the same integer \(N\geq2\),
\[
 H^q(g)\text{ is an isomorphism for }q\leq N,
 \qquad
 H^{N+1}(g)\text{ is injective}.
\]
Then \(g\) extends to a minimal Sullivan model
\[
             P\longrightarrow P\otimes\Lambda W
                 \xrightarrow{\simeq}A
\]
with \(W^q=0\) for \(q\leq N\).
\end{lemma}

\begin{proof}
We use the following induction invariant at cohomological degree \(q\):
before any generator is added to represent the cokernel in degree \(q+1\),
all classes in the kernel in degree \(q+1\) have already been killed by
generators of degree \(q\).  Thus a closed generator of degree \(q\) is added
to represent a class in the cokernel of \(H^q(g)\), whereas a generator of
degree \(q\) with nonzero differential is added first to kill a class in the
kernel of \(H^{q+1}(g)\).  The stated isomorphisms and injectivity
force both groups to vanish for \(q\leq N\), so no generator is added in
those degrees.  Since \(P\) itself has no generator above degree \(N\),
every element of \(P\) in a larger degree is decomposable.  The standard
minimal Sullivan construction can therefore continue above degree
\(N\), with decomposable differentials.  Adjoin the new generators in
nondecreasing degree, killing the kernel in degree \(q+1\) before
representing the cokernel in that degree.  When a generator of degree
\(q\) is adjoined, no previously adjoined generator has degree \(q+1\);
its differential therefore has no linear term in \(W\).  The resulting
minimal Sullivan extension is quasi-isomorphic to \(A\).
\end{proof}

\begin{keyidea}{Reusable lesson}
The quotient criterion converts homology nilpotency into a surjectivity
problem.  Every later argument asks the same question in a different form:
can a quasi-isomorphism from a minimal Sullivan algebra be made surjective
without increasing augmentation length?
\end{keyidea}

\chapter{Geometric Preliminaries}\label{stab:sec:geom-prelim}
\begin{chapterguide}[title={Chapter guide}]
This chapter translates the algebraic truncation of a minimal model into the
geometry used by the stabilization theorem.  It reviews normalized rational
LS category, \(n\)-LS applications, the chosen Sullivan root, and the
fibre--cofibre construction.  The final suspension is shown to be a wedge of
rational spheres; its degreewise finiteness is deferred to Part~\ref{part:stabilization}.
\end{chapterguide}
\index{Lusternik--Schnirelmann category}
\index{Ganea fibration}
\index{n-LS application@$n$-LS application}
\index{fibre--cofibre construction}
\index{rational wedge of spheres}

\section{Sullivan realization conventions}\label{stab:sec:sullivan-realization}
\index{Sullivan realization}

We use the standard contravariant polynomial-forms functor
\[
 A_{\mathrm{PL}}:\{\text{simplicial sets}\}^{\mathrm{{op}}}
       \longrightarrow \{\text{CDGAs over }\Q\}.
\]
For an augmented CDGA \(A\), its Sullivan realization is the simplicial
set
\[
 |A|_q=\operatorname{Hom}_{\mathrm{CDGA}}
          \bigl(A,A_{\mathrm{PL}}(\Delta^q)\bigr).
\]
Thus a CDGA morphism \(A\to B\) realizes contravariantly as a map

\[
                         |B|\longrightarrow |A|.
\]
There is a natural evaluation morphism
\(A\to A_{\mathrm{PL}}(|A|)\).  For a simply connected Sullivan algebra
of finite type, it is a quasi-isomorphism, and, for a simply connected
rational space \(X\) of finite rational type, the natural comparison
\(X\to|A_{\mathrm{PL}}(X)|\) is a rational equivalence
\cite[Chapter~17]{FHT2001}.  These conventions explain every occurrence
of

\[
                 |A|,\qquad A_{\mathrm{PL}}(X),
\]
and every reversal of arrows under realization below.

\section{Rational Ganea stages}\label{sec:geometric-invariants}

For a pointed connected CW complex \(X\), the normalized LS category
\(\operatorname{cat}(X)\) is the least \(n\) for which \(X\) is covered
by \(n+1\) open subsets whose inclusions in \(X\) are nullhomotopic.
Equivalently, the \(n\)-th Ganea fibration
\[
 F_n(X)\longrightarrow G_n(X)\xrightarrow{p_n}X
\]
admits a homotopy section \cite{Ganea1960}.  Its fibre is the iterated
join
\begin{equation}\label{stab:eq:ganea-fibre}
 F_n(X)\simeq (\Omega X)^{*(n+1)}
 \simeq \Sigma^n(\Omega X)^{\wedge(n+1)}.
\end{equation}
The rational category \(\catz(X)\) is the category of the rationalization
of \(X\).

The normalized \emph{cone length} \(\operatorname{Cl}(X)\) is the least
\(n\) for which there is a sequence
\[
 *=X_0\longrightarrow X_1\longrightarrow\cdots
   \longrightarrow X_n\simeq X
\]
such that, for every \(1\leq i\leq n\), some homotopy cofibration
\(L_i\to X_{i-1}\to X_i\) realizes the indicated step.  We write
\[
                    \Clz(X)=\operatorname{Cl}(X_{\Q})
\]
for rational cone length.  Thus a noncontractible sphere has normalized
category and cone length one.
\index{cone length}
\index{rational cone length}

For later reference, we also fix the normalized sectional-category
convention.  The invariant \(\operatorname{secat}(p)\) is the least \(n\)
for which the base of a fibration replacement of \(p\) is covered by
\(n+1\) open sets admitting local homotopy sections.  The normalized
topological complexity is the sectional category of the free-path
fibration, so \(\operatorname{TC}(S^{2k+1}_{\Q})=1\).  These notions go back
to Schwarz and Farber \cite{Schwarz1966,Farber2003}.
\index{sectional category}
\index{topological complexity}

An application \(e:E\to X\) is an \(n\)-LS application if there are maps in
both directions between \(e\) and the classical Ganea map
\(p_n:G_n(X)\to X\), commuting with the maps to \(X\) up to homotopy.  If
\(e\) is a fibration, it is an \(n\)-LS fibration.  This is the notion
introduced in \cite[Definition~2.1]{ScheererTanre1997}; in particular, an
\(n\)-LS fibration admits a section if and only if
\(\operatorname{cat}(X)\leq n\).  We retain Scheerer--Tanr\'e's terminology.

The algebraic construction needed below is summarized in the next
statement.  It combines the F\'elix--Halperin category criterion
\cite[Theorem~4.7]{FelixHalperin1982} with the Sullivan realization of
Scheerer--Tanr\'e \cite[Corollary~7.3]{ScheererTanre1997}.

\begin{theorem}[Normalized Sullivan \(n\)-LS root]\label{stab:thm:root}
Let \(M=(\Lambda V,d)\) be a simply connected minimal Sullivan algebra
of finite type, and let \(n\geq1\).  The projection
\(q_n:M\twoheadrightarrow A_n=M/(M^+)^{n+1}\) admits a factorization
\begin{equation}\label{stab:eq:root-factorization}
 M\xrightarrow{i_n}P_n=(M\otimes\Lambda Y,D)
 \xrightarrow[\simeq]{\kappa_n}\!\!\!\twoheadrightarrow A_n
\end{equation}
with the following properties:
\begin{enumerate}
\item \(i_n\) is a relative Sullivan extension and \(P_n\) is minimal
      as an absolute Sullivan algebra, that is,
      \(D(V\oplus Y)\subset\Lambda^{\geq2}(V\oplus Y)\);
\item \(\kappa_n|_M=q_n\), and the relative generators may be normalized
      so that \(\kappa_n(Y)=0\);
\item \(Y\) is of finite type;
\item if \(\catz(|M|)\leq n\), then \(i_n\) admits a strict CDGA
      retraction \(\rho_n:P_n\to M\).
\end{enumerate}
The realization \(R_n=|P_n|\to|M|\) is an \(n\)-LS application; after
fibration replacement it is an \(n\)-LS fibration.  Thus there are maps in
both directions over \(|M|\), up to homotopy, between this realization and
the classical rational Ganea map.  No identification of their total spaces
is asserted.
\end{theorem}

For reference, the source of each clause is listed separately in the following
dictionary; the proof below explains how the clauses are assembled on one
chosen factorization.
\begin{center}
\small
\begin{tabularx}{\textwidth}{@{}l>{\raggedright\arraybackslash}X@{}}
\toprule
Clause & Source\\
\midrule
Normalized relative factorization
  & Lemma~\ref{stab:lem:decomposable-quotient-model} and
    \cite[Proposition~12.2]{FHT2001}.\\
Strict retraction at category \(\leq n\)
  & \cite[Definition~4.6 and Theorem~4.7(iv)]{FelixHalperin1982}.\\
Realization as an \(n\)-LS application
  & \cite[Corollary~7.3]{ScheererTanre1997}.\\
Distinction from the classical Ganea total space
  & \cite[Proposition~2.7 and Theorem~3.2]{FelixHalperin1982}.\\
\bottomrule
\end{tabularx}
\end{center}

\begin{proof}
Because \(n\geq1\), the kernel \((M^+)^{n+1}\) is contained in
\((M^+)^2\).  Lemma~\ref{stab:lem:decomposable-quotient-model}, applied to
\(q_n\), therefore gives the factorization
\eqref{stab:eq:root-factorization} directly: \(P_n\) is minimal as an absolute
Sullivan algebra, \(\kappa_n\) is a surjective quasi-isomorphism, and
\(\kappa_n(Y)=0\).

The algebra \(A_n\) is simply connected and of finite type.  The
degree-by-degree minimal-model construction
\cite[Proposition~12.2]{FHT2001} therefore shows that \(Y\) is of finite type.
When \(\catz(|M|)\leq n\), the F\'elix--Halperin
criterion \cite[Definition~4.6 and Theorem~4.7(iv)]{FelixHalperin1982}
asserts that \(q_n\) makes \(M\) a retract of \(A_n\) in the homotopy
category.  Since \(\kappa_n:P_n\to A_n\) is a Sullivan representative
under \(M\), comparison in the undercategory gives a morphism
\(r:P_n\to M\) with \(ri_n\simeq\id_M\).  Choose a Sullivan homotopy
\[
 H:M\longrightarrow M\otimes\Lambda(t,dt),
 \qquad \operatorname{ev}_0H=ri_n,
 \qquad \operatorname{ev}_1H=\id_M.
\]
The evaluation
\(\operatorname{ev}_0:M\otimes\Lambda(t,dt)\twoheadrightarrow M\)
is a surjective quasi-isomorphism.  Apply
Lemma~\ref{stab:lem:sullivan-lifting} to the commutative square determined by
\(i_n\), \(H\), \(r\), and \(\operatorname{ev}_0\).  It produces
\(\widetilde H:P_n\to M\otimes\Lambda(t,dt)\) such that
\(\widetilde Hi_n=H\) and \(\operatorname{ev}_0\widetilde H=r\).  Set
\[
             \rho_n=\operatorname{ev}_1\widetilde H:P_n\longrightarrow M.
\]
Then \(\rho_ni_n=\operatorname{ev}_1H=\id_M\), which is the strict
retraction asserted in part~(4).  This is the homotopy-extension
realization, on the chosen Sullivan factorization, of the strictification
clause in \cite[Definition~4.6]{FelixHalperin1982}; it does not assert a
strict algebra section of the quotient \(q_n\).

Finally, \(|\kappa_n|:|A_n|\to |P_n|\) is a rational equivalence over
\(|M|\).  Scheerer--Tanr\'e prove that a fibration replacement of the
realization of \(q_n:M\to A_n\) is an \(n\)-LS fibration
\cite[Corollary~7.3]{ScheererTanre1997}.  Transporting the two comparison
maps across \(|\kappa_n|\) proves the last assertion.
\end{proof}

\begin{definition}[Sullivan LS root]\label{stab:def:LS-root}
A factorization satisfying Theorem~\ref{stab:thm:root} is called a
\emph{normalized Sullivan \(n\)-LS root} of \(M\).  Once such a
factorization is chosen, its retraction and realized section are part of
the chosen root data.  The root is called \emph{optimal} when

\[
                         n=\catz(|M|).
\]
The adjective ``chosen'' is essential: neither the relative generators nor
the retraction are asserted to be canonical.
\end{definition}

\begin{remark}\label{stab:rem:not-classical-ganea}
The last assertion does not identify \(R_n\) with the classical Ganea
space \(G_n(X)\).  F\'elix--Halperin
\cite[Proposition~2.7 and Theorem~3.2]{FelixHalperin1982} compare the
classical rational Ganea total space with the realization of the truncation
only after adjoining a wedge of rational spheres; that comparison need not
be an equivalence over \(X\).  Accordingly, \(F_\rho\) below denotes the
fibre of the chosen \(n\)-LS representative and is never identified with
the classical Ganea fibre \(F_n(X)\).
\end{remark}

Put
\begin{equation}\label{stab:eq:root-kernel}
 K_n=\ker\kappa_n.
\end{equation}
Since \(\kappa_n\) is a surjective quasi-isomorphism,
\begin{equation}\label{stab:eq:K-acyclic-short}
 H(K_n)=0,
 \qquad
 (P_n^+)^{n+1}\subset K_n.
\end{equation}
The second inclusion follows directly from the normalization: a monomial
of length at least \(n+1\) either contains a relative generator, killed by
\(\kappa_n\), or lies in \((M^+)^{n+1}\), killed by \(q_n\).

\section{Admissible spherical wedges}
\label{stab:sec:admissible-spherical-wedges}

\begin{definition}[Degreewise finite rational sphere wedges]
\label{stab:def:spherical-wedges}
Let \(\mathcal S_{\mathrm{lf}}\) denote the class of pointed rational spaces
rationally equivalent to a wedge
\[
             \bigvee_{\lambda\in\Lambda}S_{\Q}^{\,n_\lambda},
             \qquad n_\lambda\geq2,
\]
such that, for every \(m\geq2\), only finitely many indices \(\lambda\)
satisfy \(n_\lambda=m\).  The empty wedge, identified with the one-point
space, is included.  For a simply connected rational space \(X\) of finite
type, set
\[
 \sHnil(X)
   :=\inf_{S\in\mathcal S_{\mathrm{lf}}}\Hnilz(X\vee S),
 \qquad
 \sClz(X)
   :=\inf_{S\in\mathcal S_{\mathrm{lf}}}\Clz(X\vee S).
\]
\end{definition}
\index{Slf@$\mathcal S_{\mathrm{lf}}$}
\index{rational wedge of spheres!degreewise finite}
\index{spherical stabilization}

\section{The fibre of a section and its cofibre}

Assume now that \(n=\catz(X)\geq1\), and let
\(q:R_n\to X\) be a fibration replacement of the realization in
Theorem~\ref{stab:thm:root}.  The strict
retraction \(\rho_n\) realizes a section
\[
 \sigma:X\longrightarrow R_n,
 \qquad q\sigma\simeq\id_X.
\]
Write
\[
 F_\rho=\hofib(q),
 \qquad
 G_\sigma=\hofib(\sigma).
\]

\begin{lemma}\label{stab:lem:fibre-section-loop}
There is a natural rational equivalence
\[
 G_\sigma\simeq_{\Q}\Omega F_\rho.
\]
If \(j:G_\sigma\to X\) is the fibre inclusion, then \(j\) is
nullhomotopic.
\end{lemma}

\begin{proof}
Comparing the two homotopy pullback squares determined by
\(q\sigma\simeq\id_X\) gives \(G_\sigma\simeq\Omega F_\rho\).  By definition
\(\sigma j\simeq *\), and
therefore
\[
 j\simeq q\sigma j\simeq *.
\]
\end{proof}

Define the fibre--cofibre space
\begin{equation}\label{stab:eq:def-Wsigma}
 W_\sigma=\hocof(j:G_\sigma\to X).
\end{equation}
The preceding nullhomotopy gives
\begin{equation}\label{stab:eq:wedge-splitting}
 W_\sigma\simeq_{\Q}X\vee\Sigma G_\sigma.
\end{equation}
In the terminology of F\'elix--Thomas \cite{FelixThomas1988}, this is exactly
the fibre--cofibre construction associated with the map \(\sigma\): one first
takes its homotopy fibre \(G_\sigma\to X\), and then the homotopy cofibre of
that fibre inclusion.  We use this single fibre--cofibre step in Cornea's
cone-length setting \cite{Cornea1994}; we do not iterate it to define the
classical Ganea spaces.  Notice in particular that
\eqref{stab:eq:wedge-splitting} splits the cofibre, not the total space
\(R_n\).

\begin{lemma}[The spherical factor]\label{stab:lem:spherical-factor}
The space \(\Sigma G_\sigma\) is rationally a wedge of simply connected
spheres.
\end{lemma}

\begin{proof}
By Lemma~\ref{stab:lem:fibre-section-loop},
\(\Sigma G_\sigma\simeq_{\Q}\Sigma\Omega F_\rho\).  Henn's
loop--suspension theorem states that, for a connected almost rational space
\(Z\), \(\Sigma\Omega Z\) is a wedge of rational spheres of dimension at
least two, possibly together with ordinary circle summands
\cite[Lemma~2]{Henn1983}.  In Henn's terminology, an ``almost rational''
space is a connected space whose homotopy groups in degrees at least two are
rational vector spaces.  Both \(R_n\) and \(X\) are simply
connected rational spaces.  The homotopy long exact sequence of

\[
                         F_\rho\longrightarrow R_n\xrightarrow{q}X
\]

therefore shows that every higher homotopy group of \(F_\rho\) is a
rational vector space.  Moreover, the section makes
\(q_*:\pi_2(R_n)\to\pi_2(X)\) surjective; exactness gives
\(\pi_1(F_\rho)=0\), and the fibre is connected.  Thus \(F_\rho\) is
almost rational and in fact simply connected.  The circle summands in
Henn's conclusion are consequently absent, so \(\Sigma G_\sigma\) is a
wedge of simply connected rational spheres.
\end{proof}

The algebraic invariants and their geometric comparison are now fixed, and the chosen Sullivan LS root supplies the fibre--cofibre model used later.  Part~\ref{part:calculations} turns next to calculations inside the fixed minimal model.

\part{Calculating Homology Nilpotency}\label{part:calculations}

\chapter{Acyclic Truncations and Essential Directions}
\label{calc:chap:acyclic}

\begin{chapterguide}[title={Chapter guide}]
Homology nilpotency is computed inside a fixed minimal Sullivan model.
The problem is therefore not merely to find a short CDGA of the same
homotopy type, but to obtain it as an acyclic quotient of that model.  This
chapter develops the corresponding calculus.  Cohomological truncation gives
the first general acyclic ideals; an adapted decomposition identifies
directions forced into every such truncation; and a degree-by-degree
elimination argument gives a rigorous upper-bound criterion.
\end{chapterguide}

\index{acyclic differential ideal}
\index{homology nilpotency}
\index{augmentation length}
\index{essential direction}
\index{adapted decomposition}

Throughout this chapter,
\[
                 A=(\Lambda V,d)
\]
is a simply connected minimal Sullivan algebra.  Thus
\(V^0=V^1=0\), \(dV\subseteq\Lambda^{\geq2}V\), and
\(\mA=A^+\).  The power \(\mA^p\) means augmentation length at least
\(p\), whereas \(A^j\) means cohomological degree \(j\).  Minimality and
the derivation rule give
\[
                       d(\mA^p)\subseteq\mA^{p+1}
                       \subseteq\mA^p.
\]
Consequently every \(\mA^p\) is a differential ideal.

If \(A\) models an \(r\)-connected rational space, put
\(q=r+1\).  Then \(V^j=0\) for \(j<q\), and every word of length
\(p\) has degree at least \(pq\):
\begin{equation}
                       \mA^p\subseteq A^{\geq pq}.
\label{calc:eq:degree-length}
\end{equation}

With respect to the chosen generator space, write
\[
                         d=d_2+d_3+\cdots,
             \qquad d_i(V)\subseteq\Lambda^iV.
\]
The first component \(d_2\) is the quadratic differential.  Although this
decomposition depends on the Sullivan coordinates, the assertion that a
specified linear map \(d_2:V^j\to(\Lambda^2V)^{j+1}\) is injective is
well defined once those coordinates have been fixed.  All later uses of this
condition occur in one fixed minimal model.

\begin{keyidea}{How to use the bounds}
For an upper bound, construct an acyclic ideal containing the required
augmentation power, or invoke the quadratic-shadow theorem.  For a lower
bound, exhibit an essential direction or a primitive forced into every
candidate ideal.  Under a cell attachment, the category test and the
\(\Hnil\)-test are logically separate.  The hypotheses in the criteria below
are sufficient and are not claimed to be necessary.
\end{keyidea}

\section{The quadratic shadow}

\index{quadratic shadow}
\index{rational Toomer invariant}
\index{coformal CDGA}
The identity \(d^2=0\) implies \(d_2^2=0\).  We call
\[
                         A^{(2)}=(\Lambda V,d_2)
\]
the \emph{quadratic shadow} of \(A\); its realization is the coformal
shadow of the rational homotopy type.  For \(n\geq0\), let
\[
 q_n^{(2)}:A^{(2)}\longrightarrow
 A^{(2)}/\Lambda^{\geq n+1}V
\]
be the word-length projection, and put
\[
 e_0(A^{(2)})=\inf\{n\mid H(q_n^{(2)})\text{ is injective}\}.
\]
This is the normalized rational Toomer invariant of the quadratic shadow.

\begin{theorem}[Quadratic-shadow upper bound]
\label{calc:thm:quadratic-shadow}
Let \(A=(\Lambda V,d)\) be a simply connected minimal Sullivan algebra of
finite type.  Then
\begin{equation}
 \Hnil(A)\leq e_0(A^{(2)})
       \leq \catz\bigl(|A^{(2)}|\bigr).
\label{calc:eq:quadratic-shadow-bound}
\end{equation}
Consequently,
\[
 \catz(|A|)\leq\Hnil(A)\leq\catz\bigl(|A^{(2)}|\bigr).
\]
In particular, if
\(\catz(|A^{(2)}|)<n\), then \(\Hnil(A)<n\).
\end{theorem}

\begin{proof}
Fix \(n\) such that \(H(q_n^{(2)})\) is injective, and write
\(A_{[p]}=\Lambda^pV\) for the subspace of words of length \(p\).
Because \(d_2\) raises word length by exactly one, every
\(d_2\)-cycle \(z\in A_{[p]}\) with \(p\geq n+1\) represents a class
annihilated by \(H(q_n^{(2)})\).  It is therefore a \(d_2\)-boundary in
\(A^{(2)}\).  Taking the word-length-\((p-1)\) component of a primitive
shows that
\begin{equation}
 Z(A_{[p]},d_2)=d_2(A_{[p-1]})
                    \qquad(p\geq n+1).
\label{calc:eq:quadratic-exactness}
\end{equation}

Choose a complement
\[
 A_{[n]}=\ker(d_2|_{A_{[n]}})\oplus U.
\]
Equation~\eqref{calc:eq:quadratic-exactness} for \(p=n+1\) says that
\[
 d_2:U\xrightarrow{\cong}Z(A_{[n+1]},d_2).
\]
Set
\[
                         J_2=U\oplus
                         \bigoplus_{p\geq n+1}A_{[p]}.
\]
This is an ideal: multiplying \(U\subset A_{[n]}\) by a positive-length
element lands in \(\bigoplus_{p\geq n+1}A_{[p]}\).  It is
\(d_2\)-stable, and the displayed isomorphism together with
\eqref{calc:eq:quadratic-exactness} proves
\[
                              H(J_2,d_2)=0.
\]

The same graded ideal is stable under the full differential.  Indeed,
\((d-d_2)U\subseteq\Lambda^{\geq n+2}V\), while
\[
 d\bigl(\Lambda^{\geq n+1}V\bigr)
        \subseteq\Lambda^{\geq n+2}V.
\]
Filter \((J_2,d)\) by the decreasing filtration
\(F^pJ_2=J_2\cap\mathfrak m_A^p\).  We use the cohomological convention
\[
 E_0^{p,q}=F^pJ_2^{p+q}/F^{p+1}J_2^{p+q},
 \qquad d_r:E_r^{p,q}\longrightarrow E_r^{p+r,q-r+1}.
\]
In every total degree the filtration is finite because \(V^0=V^1=0\).
Since the full differential raises word length by at least one, the induced
\(d_0\) is zero, and the first nonzero associated-graded differential is the
quadratic part:
\[
 E_0=\operatorname{gr}_FJ_2,\qquad d_0=0,
 \qquad E_1\cong\operatorname{gr}_FJ_2,
 \qquad d_1=(d_2)_{\operatorname{gr}}.
\]
Consequently \(E_2\cong H(J_2,d_2)=0\), and the spectral sequence
converges to \(H(J_2,d)=0\).  Thus \(J_2\lhd A\) is
an acyclic differential ideal containing
\(\Lambda^{\geq n+1}V=\mA^{n+1}\), and hence \(\Hnil(A)\leq n\).
Taking the least such \(n\) proves the first inequality in
\eqref{calc:eq:quadratic-shadow-bound}.
If there is no such \(n\), then \(e_0(A^{(2)})=\infty\) and that
inequality is tautological.

If \(\catz(|A^{(2)}|)\leq n\), the F\'elix--Halperin criterion implies that
\(q_n^{(2)}\) admits a homotopy retraction in the CDGA homotopy category.
Its map in cohomology is therefore injective, which proves
\(e_0(A^{(2)})\leq\catz(|A^{(2)}|)\).  The remaining lower bound is
Proposition~\ref{stab:prop:cat-lower}.
\end{proof}

\begin{corollary}[Coformal case]
\label{calc:cor:coformal-equality}
If \(d=d_2\), then
\[
 \catz(|A|)=\Clz(|A|)=\nilh(A)=\Hnil(A).
\]
\end{corollary}

\begin{proof}
In this case \(A=A^{(2)}\).  Combine
Theorem~\ref{calc:thm:quadratic-shadow} with
Proposition~\ref{stab:prop:cat-lower}.
\end{proof}

\begin{remark}
The equality of the standard rational LS-category invariants in the
coformal case is classical; see Lemaire--Sigrist
\cite{LemaireSigrist1981}.  The perturbative content of
Theorem~\ref{calc:thm:quadratic-shadow} is that the same quadratic
upper bound controls \(\Hnil(\Lambda V,d)\) when arbitrary higher
terms \(d_3,d_4,\ldots\) are included.  Although \(d_2\) is written in
chosen Sullivan coordinates, the isomorphism type of \(A^{(2)}\) is
intrinsic: the linear part of an isomorphism between minimal Sullivan
algebras intertwines their quadratic differentials.
\end{remark}

\section{Cohomological truncation}

The following elementary construction is the starting point for all the
acyclic ideals used below.  It is the cochain-level truncation appearing in
standard treatments of Sullivan models; compare
\cite[pp.~146--147]{FHT2001}.

\begin{proposition}[Cohomological truncation]
\label{calc:prop:cohomological-truncation}
Assume that \(H^*(A)\) is finite dimensional in every degree and vanishes
above
\[
                    N=\max\{j\mid H^j(A)\neq0\}<\infty.
\]
There is an acyclic differential ideal \(T_N\lhd A\) such that
\[
                            T_N^j=A^j\qquad(j>N).
\]
If \(A\) models an \(r\)-connected space \(X\), then for every
\(m\geq0\),
\begin{equation}
        q(m+1)>N\quad\Longrightarrow\quad \Hnilz(X)\leq m,
        \qquad q=r+1.
\label{calc:eq:dimension-connectivity}
\end{equation}
In particular,
\[
                            \Hnilz(X)\leq
                            \left\lfloor\frac{N}{r+1}\right\rfloor.
\]
If moreover \(m=\catz(X)\) and \(N<q(m+1)\), then
\(\Hnilz(X)=m\).
\end{proposition}

\begin{proof}
If \(N=0\), take \(T_0=A^+\).  The reduced cohomology of \(A\) is zero,
so this ideal is acyclic, and all assertions follow.  Suppose \(N>0\).
Write \(Z^j(A)\) and \(B^j(A)\) for cycles and boundaries.  Choose
\[
 A^{N-1}=Z^{N-1}(A)\oplus E.
\]
Then \(d:E\to B^N(A)\) is an isomorphism.  Choose a subspace
\(H_N\subseteq Z^N(A)\) representing \(H^N(A)\), and then a complement
\(C\) to the cycles, so that
\[
                 A^N=H_N\oplus dE\oplus C,
                 \qquad Z^N(A)=H_N\oplus dE.
\]
Define a graded subspace \(T_N\subseteq A\) by
\[
 T_N^j=
 \begin{cases}
  0,       &j<N-1,\\
  E,       &j=N-1,\\
  dE\oplus C,&j=N,\\
  A^j,     &j>N.
 \end{cases}
\]
It is a subcomplex.  It is also an ideal: a product involving an element
of degree \(N-1\) or \(N\) and a positive-degree element has degree
strictly greater than \(N\), because \(A^1=0\); all other cases are
immediate.

The pair \(E\xrightarrow{d}dE\) contributes no cohomology.  The map
\(d:C\to A^{N+1}\) is injective.  For a cocycle of degree \(N+1\) written
as \(da\), subtracting the \(H_N\)-component of \(a\in A^N\), if necessary,
produces a primitive in \(T_N\).  In degrees greater than \(N+1\), every
primitive already lies in \(T_N\).  Since \(H^{>N}(A)=0\), every cocycle of
degree greater than \(N\) in \(T_N\) is therefore a boundary in \(T_N\).
Hence \(H^*(T_N)=0\).

If \(q(m+1)>N\), the degree estimate
\eqref{calc:eq:degree-length} gives
\[
                      \mA^{m+1}\subseteq A^{>N}\subseteq T_N.
\]
The defining acyclic-ideal criterion for \(\Hnil\) therefore proves
\eqref{calc:eq:dimension-connectivity}.  Taking
\(m=\lfloor N/q\rfloor\) gives the displayed universal bound.  Finally,
if \(m=\catz(X)\), the inequality
\(\catz(X)\leq\Hnilz(X)\) from
Proposition~\ref{stab:prop:cat-lower} supplies the reverse inequality.
\end{proof}

\begin{remark}
The strict inequality in \eqref{calc:eq:dimension-connectivity} is
essential to this construction.  At \(q(m+1)=N\), word length only gives
\(\mA^{m+1}\subseteq A^{\geq N}\), and the nonzero top cohomology
representatives in \(H_N\) are intentionally excluded from \(T_N\).
\end{remark}

\section{Spherical normalization and adapted decompositions}

Let
\[
                 h_A:H^+(A)\longrightarrow Q(A)=V
\]
be induced by projection modulo decomposables.  Its image is the rational
spherical homology of the minimal model.  If a class maps to
\(v\in V\), it has a representative \(v+\xi\), with
\(\xi\in\mA^2\).  Replacing \(v\) by \(v+\xi\) is a triangular Sullivan
change of variables and makes that representative a closed generator.
Proceeding in Sullivan order proves the following normalization.

\begin{lemma}[Spherical normalization]
\label{calc:lem:spherical-normalization}
After a triangular CDGA isomorphism, the generator space admits a degreewise
decomposition
\[
                         V^j=(\im h_A)^j\oplus U^j
\]
in which every chosen generator in \(\im h_A\) is closed.
\end{lemma}

\begin{proof}
Proceed first by cohomological degree and, within a fixed degree, by any
chosen order on a basis of \((\im h_A)^j\).  For a basis vector \(v\), choose
a cocycle \(v+\xi\) projecting to it.  Every generator occurring in the
decomposable term \(\xi\) has degree strictly smaller than \(|v|\): a factor
of degree at least \(|v|\), multiplied by another positive-degree factor,
would give total degree greater than \(|v|\).  Thus all factors of \(\xi\)
have already been treated.  The substitution \(v\mapsto v+\xi\) is
triangular for the degree-refined Sullivan filtration, is invertible, fixes
all earlier choices, and makes the new generator closed.  Repeating this
construction degree by degree gives the asserted normalization.
\end{proof}

We now isolate the exact structural hypothesis needed for the lower-bound
argument.  Let \(X\) be \(r\)-connected, of finite rational type, and
suppose that its rational cohomology is bounded.  Put
\begin{equation}
 \begin{aligned}
 q&=r+1, & n&=\catz(X),\\
 \tau_n&=q(n+1), & c_n&=\tau_n-1,\\
 \chi_n&=q(2n+1)-3.&&
 \end{aligned}
\label{calc:eq:three-thresholds}
\end{equation}
Assume
\begin{equation}
                     N_X:=\max\{j\mid H^j(X;\Q)\neq0\}<\tau_n.
\label{calc:eq:top-range}
\end{equation}

\begin{definition}[Adapted decomposition]
\label{calc:def:adapted-decomposition}
Fix \(M\geq c_n\).  An \emph{adapted decomposition through degree
\(M\)} of the minimal model \(A=(\Lambda V,d)\) is a decomposition
\[
                              V=W\oplus U\oplus Z
\]
with the following properties.
\begin{enumerate}[label=\textup{(A\arabic*)}]
\item \((\Lambda W,d)\) is a sub-CDGA and \(W\) contains
      \(\im h_A\).
\item \(V^j=W^j\) for \(j<c_n\), and
      \(V^j=W^j\oplus U^j\) for \(c_n\leq j\leq M\).
\item \(U^j=0\) outside \([c_n,M]\), \(Z^j=0\) for \(j\leq M\),
      and \(V^j=W^j\oplus Z^j\) for \(j>M\).
\item For \(c_n\leq j\leq M\),
\begin{equation}
 dU^j\subseteq \bigl((\Lambda W)^+\bigr)^{n+1}
       +\Lambda W\otimes\Lambda^+(U^{<j}).
\label{calc:eq:adapted-differential}
\end{equation}
\end{enumerate}
\end{definition}

The existence of an adapted decomposition is a hypothesis, not a formal
consequence of \(\catz(X)=n\).  It must be established for the family at
hand, or checked directly in a specified minimal model.

For \(c_n\leq j\leq M\), put
\begin{align}
 B_j&=\Lambda\bigl(W^{\leq j}\oplus U^{<j}\bigr),
 \label{calc:eq:Bj}\\
 K_j&=B_j\cap\left(\bigl((\Lambda W)^+\bigr)^{n+1}
       +\Lambda W\otimes\Lambda^+(U^{<j})\right).
 \label{calc:eq:Kj}
\end{align}
Then \(B_j\) is a sub-CDGA and \(K_j\lhd B_j\) is a differential
ideal.

\begin{lemma}[Obstruction-space isomorphism]
\label{calc:lem:obstruction-isomorphism}
For every \(c_n\leq j\leq M\), the linear map
\[
                 \partial_j:U^j\longrightarrow H^{j+1}(B_j),
                 \qquad u\longmapsto[du],
\]
is an isomorphism.
\end{lemma}

\begin{proof}
Suppose \(du=d\xi\) for some \(\xi\in B_j^j\).  Then
\(u-\xi\) is a cocycle whose image in the indecomposables has a nonzero
\(U^j\)-component.  This would give a spherical direction outside \(W\),
contrary to (A1).  Thus \(\partial_j\) is injective.

Conversely, let \(z\in B_j^{j+1}\) be a cocycle.  Since
\(j+1\geq c_n+1=\tau_n>N_X\), it is exact in \(A\): choose
\(y\in A^j\) with \(dy=z\).  There are no \(Z\)-generators through
degree \(M\).  Degree considerations therefore give a unique expression
\[
                           y=\xi+u,\qquad
                     \xi\in B_j^j,\quad u\in U^j.
\]
Hence \([z]=[du]\), proving surjectivity.
\end{proof}

\section{Essential directions}

\begin{definition}
\label{calc:def:essential-direction}
A homogeneous element \(u\in U^j\) is \emph{\(n\)-essential} if, for
every acyclic differential ideal \(I\lhd A\) satisfying
\(\mA^{n+1}\subseteq I\), there is an element \(\omega\in B_j^j\)
such that
\[
                               u+\omega\in I.
\]
No closedness condition is imposed on \(\omega\).
\end{definition}

\begin{theorem}[Essential-direction theorem]
\label{calc:thm:essential-directions}
Every element of the adapted block \(U\) is \(n\)-essential.
\end{theorem}

\begin{proof}
Fix an acyclic differential ideal \(I\lhd A\) containing
\(\mA^{n+1}\).  We construct, degree by degree, CDGA automorphisms
\[
                              \phi_j:B_j\longrightarrow B_j
\]
which fix \(W^{\leq j}\) and satisfy \(\phi_j(K_j)\subseteq I\).
Before the first degree in which \(U\) is nonzero, take the identity.

Assume \(\phi_j\) has been constructed and choose a basis
\(u_1,\ldots,u_s\) of \(U^j\).  By
\eqref{calc:eq:adapted-differential}, \(du_\ell\in K_j\).  Hence
\(\phi_j(du_\ell)\) is a cocycle in \(I\).  Acyclicity supplies
\(x_\ell\in I^j\) such that
\[
                            dx_\ell=\phi_j(du_\ell).
\]
Because there are no \(Z\)-generators in this range, write uniquely
\[
                    x_\ell=\sum_{i=1}^s a_{i\ell}u_i+\xi_\ell,
                    \qquad \xi_\ell\in B_j^j.
\]
After applying \(d\) and passing to \(H^{j+1}(B_j)\), one obtains
\[
                H(\phi_j)[du_\ell]
                   =\sum_{i=1}^s a_{i\ell}[du_i].
\]
Lemma~\ref{calc:lem:obstruction-isomorphism} and the invertibility of
\(H(\phi_j)\) imply that the matrix \((a_{i\ell})\) is invertible.

Extend \(\phi_j\) by sending \(u_\ell\) to \(x_\ell\), while fixing the new
\(W\)-generators.  The identity \(dx_\ell=\phi_j(du_\ell)\) makes this a CDGA
morphism, and the invertibility of the new linear block makes it a triangular
automorphism.  Every word of length at least \(n+1\) in \(W\) already lies in
\(I\), and every term involving an earlier \(U\)-generator now has a factor
in \(I\); therefore the extended map carries \(K_{j+1}\) into \(I\).

Finally, apply the inverse matrix \((a_{i\ell})^{-1}\) to the elements
\(x_\ell\in I\).  For every \(u_i\) this produces
\(u_i+\omega_i\in I\), with \(\omega_i\in B_j^j\).  Induction through
all degrees of \(U\) proves the theorem.
\end{proof}

\begin{remark}
When \(\dim U^j>1\), the invariant datum is the whole obstruction space
\(U^j\cong H^{j+1}(B_j)\), rather than any preferred coordinate basis.  The
correction \(\omega_i\) changes under a triangular change of Sullivan
generators, while essentiality of the direction modulo the earlier
algebra does not.
\end{remark}

\section{Eliminating a long-word ideal}

We next give the constructive half of the calculus.  Let
\(A=(\Lambda E,\delta)\) be simply connected and minimal, assume
\(E^{\leq r}=0\), and put \(q=r+1\).

\begin{lemma}[First-slice elimination]
\label{calc:lem:first-slice}
Let \(I\lhd A\) be a differential ideal such that \(I^j=0\) for
\(j<a\).  For every \(a\leq s\leq b\), choose cocycles
\(\psi_{s,i}\in I^s\) representing a basis of \(H^s(I)\).  Suppose
there are decomposable elements \(\eta_{s,i}\in A^{s-1}\) satisfying
\(\delta\eta_{s,i}=\psi_{s,i}\), and suppose
\[
                         (a-1)+q=b+1.
\]
If \(I\subseteq\mathfrak m_A^3\), then
\[
              J=I+\langle\eta_{s,i}\mid a\leq s\leq b\rangle
\]
is a differential ideal satisfying
\[
       H^j(J)=0\quad(j\leq b),
       \qquad J^j\subseteq\mathfrak m_A^3\quad(j\geq b).
\]
\end{lemma}

\begin{proof}
A nontrivial product of a new primitive with an element of
\(\mathfrak m_A\) has degree at least
\((a-1)+q=b+1\).  Hence, through degree \(b\), the quotient \(J/I\) is
spanned only by the residue classes of the \(\eta_{s,i}\).  These residue
classes are linearly independent: a relation modulo \(I\), after applying
\(\delta\), would give a relation among the basis classes
\([\psi_{s,i}]\) in the corresponding cohomology groups of \(I\).  In the
long exact sequence of \(0\to I\to J\to J/I\to0\), the connecting map sends
\([\eta_{s,i}]\) to \([\psi_{s,i}]\).  It is therefore injective with image
all of \(H^s(I)\), so the prescribed classes are killed and no lower class
is created.

Every new product containing one primitive and another positive element has
augmentation length at least three; a product of two new primitives has
length at least four.  All such products start above the range just examined.
Together with \(I\subseteq\mathfrak m_A^3\), this proves the asserted
high-degree length condition.
\end{proof}

\begin{lemma}[One-degree elimination]
\label{calc:lem:one-degree}
Suppose \(I\lhd A\) satisfies
\[
 H^j(I)=0\ (j\leq N-1),
 \qquad I^j\subseteq\mathfrak m_A^3\ (j\geq N-1).
\]
Let \(\theta_1,\ldots,\theta_t\in I^N\) be cocycles representing a
basis of \(H^N(I)\).  If there are decomposable
\(\sigma_i\in A^{N-1}\) with \(\delta\sigma_i=\theta_i\), then
\[
                         J=I+\langle\sigma_1,\ldots,\sigma_t\rangle
\]
satisfies
\[
 H^j(J)=0\ (j\leq N),
 \qquad J^j\subseteq\mathfrak m_A^3\ (j\geq N).
\]
\end{lemma}

\begin{proof}
Products \(\sigma_i\mathfrak m_A\) begin in degree at least
\(N-1+q>N\).  Up to degree \(N\), the quotient \(J/I\) is therefore the
span of the residue classes of the \(\sigma_i\) in degree \(N-1\).  A
linear relation among those residue classes would, after differentiation,
give a relation among the basis classes \([\theta_i]\in H^N(I)\); hence
they are independent.  The connecting map in the long exact sequence of
\(0\to I\to J\to J/I\to0\) sends \([\sigma_i]\) to \([\theta_i]\).
It follows that \(H^N(I)\) is killed, no class in degree \(N-1\) is created,
and lower cohomology is unchanged.  Products with one new primitive have
length at least three, products of two have length at least four, and all of
them occur above degree \(N\).  This proves the remaining length assertion.
\end{proof}

The decomposability assumption has a useful quadratic test.  If
\(\delta\sigma\in\mathfrak m_A^3\) and \(x\in E\) is the linear part
of \(\sigma\), then the length-two component of the equation is
\[
                              d_2x=0.
\]
Thus injectivity of \(d_2:E^{N-1}\to(\Lambda^2E)^N\) forces every
primitive of an ideal cycle in degree \(N\) to be decomposable.

\section{The top class and acyclic completion}

At the final degree, arbitrary complements are unsafe: the unique top
cohomology class must remain outside the ideal.

\begin{lemma}[Acyclic completion at the top]
\label{calc:lem:top-completion}
Assume
\[
                  H^{k+1}(A)=\Q[\Omega],
                  \qquad H^{>k+1}(A)=0.
\]
Let \(I\lhd A\) be a differential ideal such that
\(H^j(I)=0\) for \(j\leq k\), and assume that no cocycle of
\(I^{k+1}\) represents a nonzero multiple of \([\Omega]\).  Then
\(I\) is contained in an acyclic differential ideal \(J\lhd A\).
\end{lemma}

\begin{proof}
Write \(Z^j=Z^j(A)\).  Since \(H^k(I)=0\),
\[
                      I^k\cap Z^k=\delta(I^{k-1}).
\]
Choose a complement \(L\) to this cycle space in \(I^k\), and extend
it to a complement \(P\) of \(Z^k\) in \(A^k\):
\[
                         A^k=Z^k\oplus P,
                         \qquad L\subseteq P.
\]
Then \(\delta:P\to B^{k+1}(A)\) is an isomorphism.  The top-class
hypothesis says that the map \(H^{k+1}(I)\to H^{k+1}(A)\) is zero;
hence
\(I^{k+1}\cap Z^{k+1}\subseteq\delta P\).  Choose a complement to this
intersection inside \(I^{k+1}\), and enlarge it to a subspace \(C\) for
which
\[
                A^{k+1}=\Q\Omega\oplus\delta P\oplus C,
                \qquad I^{k+1}\subseteq\delta P\oplus C.
\]
Define
\[
 J^j=
 \begin{cases}
  I^j,                              &j<k,\\
  (I^k\cap Z^k)\oplus P,            &j=k,\\
  \delta P\oplus C,                 &j=k+1,\\
  A^j,                              &j>k+1.
 \end{cases}
\]
It contains \(I\).  Since \(A^1=0\), multiplication of either new
degree-\(k\) or degree-\(k+1\) summand by a positive element lands above
degree \(k+1\); hence \(J\) is an ideal.  It is a subcomplex by
construction.

Below degree \(k\), acyclicity is inherited from \(I\).  The cycles in
degree \(k\) are \(\delta(I^{k-1})\), and the cycles in degree
\(k+1\) are \(\delta P\), so both are boundaries in \(J\).  If a
degree-\(k+2\) cocycle is \(\delta a\) in \(A\), subtracting the
\(\Omega\)-component of \(a\in A^{k+1}\) puts the primitive in
\(J^{k+1}\).  In higher degrees every primitive already belongs to
\(J\).  Since \(H^{>k+1}(A)=0\), this proves \(H^*(J)=0\).
\end{proof}

\begin{criterion}[Acyclic-truncation criterion]
\label{calc:crit:acyclic-truncation}
Let \(A=(\Lambda E,\delta)\) be simply connected and minimal, with
\(E^{\leq r}=0\) and \(q=r+1\).  Assume
\[
                   H^{k+1}(A)=\Q[\Omega],
                   \qquad H^{>k+1}(A)=0.
\]
Fix \(m\geq2\), integers \(a\leq b\leq k\), and a differential ideal
\(I_0\) such that
\[
 \mathfrak m_A^{m+1}\subseteq I_0\subseteq\mathfrak m_A^3,
 \qquad I_0^j=0\ (j<a),
 \qquad(a-1)+q=b+1.
\]
Assume the following three conditions.
\begin{enumerate}[label=\textup{(C\arabic*)}]
\item In every degree \(a\leq s\leq b\), a basis of \(H^s(I_0)\)
      has representatives with decomposable primitives in \(A\).
\item After the first-slice enlargement and after every later
      enlargement, a basis of the remaining ideal cohomology in each
      degree \(N=b+1,\ldots,k\) has decomposable primitives in
      \(A^{N-1}\).
\item For the final ideal \(I_k\), the map
      \(H^{k+1}(I_k)\to H^{k+1}(A)=\Q[\Omega]\) is zero.
\end{enumerate}
Then there exists an acyclic differential ideal \(J\lhd A\) containing
\(\mathfrak m_A^{m+1}\).  Consequently \(\Hnil(A)\leq m\).
\end{criterion}

\begin{proof}
Apply Lemma~\ref{calc:lem:first-slice}, then
Lemma~\ref{calc:lem:one-degree} successively through degree \(k\).
Condition (C3) permits the final application of
Lemma~\ref{calc:lem:top-completion}.  Every enlargement contains the
initial long-word ideal, so the resulting acyclic ideal contains
\(\mathfrak m_A^{m+1}\).  Definition~\ref{div:def:Hnil} gives the
claimed bound.
\end{proof}

\begin{warningbox}{Logical scope of the criterion}
A successful elimination yields an upper bound for \(\Hnil\).  The occurrence
of a linear term in a single attempted primitive, or of the top class in a
single chosen ideal, invalidates only that branch.  A lower bound requires a
coordinate-independent obstruction, such as
Theorem~\ref{calc:thm:essential-directions} or the retractive-tower
obstruction developed later.
\end{warningbox}

\addtocontents{toc}{\protect\pagebreak}

\chapter{Cell-Attachment Calculations of \texorpdfstring{\(\Hnil\)}{Hnil}}
\label{calc:chap:attachments}

\begin{chapterguide}[title={Chapter guide}]
A single rational cell can remove one indecomposable direction while turning
its differential into a new cohomology class.  This makes cell attachment an
effective laboratory for homology nilpotency.  Essential directions give a
strict lower bound; quadratic injectivity supplies a one-step upper bound;
and a strengthened retraction controls the case in which \(\Hnil\) is
preserved.  The chapter closes with a projective-space attachment for which
the lower and upper mechanisms meet at an exact value; it serves as a
controlled benchmark for the later separating model.
\end{chapterguide}

\index{cell attachment}
\index{essential attachment}
\index{quadratic differential}
\index{filtered retraction}
\index{homology nilpotency}

\section{The minimal model of one cell attachment}

Let \(X\) be simply connected and of finite rational type, let
\(A_X=(\Lambda V,d)\) be its minimal Sullivan model, and let
\[
                              f:S^k\longrightarrow X
\]
be rationally nontrivial.  The rational homotopy class of \(f\) is a
linear functional on \(V^k\).  Choose \(v\in V^k\) and a complement
\(W_f^k\) so that
\[
                      f(v)=1,\qquad V^k=\Q v\oplus W_f^k,
                      \qquad W_f^k=\ker f.
\]
Put \(Y=X\cup_f e^{k+1}\).  The Sullivan cell-attachment construction
\cite[Section~13(d)]{FHT2001} gives a minimal model
\(A_Y=(\Lambda\widetilde V,\delta)\) which may be chosen so that
\begin{equation}
                 \widetilde V^{<k}=V^{<k},
                 \qquad \widetilde V^k=W_f^k.
\label{calc:eq:attachment-generators}
\end{equation}
Thus the direction \(v\) is removed, and the image of \(dv\) in the
attached model is a cocycle.

One way to see this construction before passing to a minimal model is to adjoin a
degree-\(k+1\) element \(e\) satisfying
\[
 e^2=e\Lambda^+V=0,\qquad d'e=0,\qquad d'v=dv+e,
\]
with \(d'=d\) on a complementary set of generators.  Cancelling the
linear pair \((v,e)\) produces
\eqref{calc:eq:attachment-generators}.  The minimal model is the one used
for all word-length arguments.

When
\[
                  N_X=\max\{j\mid H^j(X;\Q)\neq0\}<k,
\]
the cohomology sequence of the pair \((Y,X)\) gives
\begin{equation}
\begin{aligned}
H^j(Y;\Q)&=0 &&(N_X<j\leq k),\\
H^{k+1}(Y;\Q)&=\Q[\Omega],\qquad \Omega=dv,\\
H^{>k+1}(Y;\Q)&=0.
\end{aligned}
\label{calc:eq:attachment-cohomology}
\end{equation}
The nontriviality of \(f\) is used here: it excludes the null-attachment case
and identifies \([dv]\) with the cellular top class.

If \(X\) is \(r\)-connected and \(q=r+1\), then for every \(p\geq1\)
the common subalgebra of the two minimal models gives the identification
\begin{equation}
       (\mathfrak m_{A_Y}^{p})^j=(\mathfrak m_{A_X}^{p})^j
       \qquad\text{whenever }j<k+(p-1)q.
\label{calc:eq:word-stability}
\end{equation}
Indeed, a word of length at least \(p\) involving the removed generator
\(v\), or any generator first appearing in degree at least \(k\), has
degree at least \(k+(p-1)q\).

We use two standard category facts.  First, a one-cell attachment raises
normalized category and cone length by at most one, so
\(\catz(Y)\leq\catz(X)+1\).  Second, suppose in addition that \(X_\Q\)
admits a CW model of dimension at most \(k\).  In the range \(N_X<k\),
rational category then does not decrease, and consequently
\begin{equation}
                 \catz(X)\leq\catz(Y)\leq\catz(X)+1.
\label{calc:eq:one-cell-category}
\end{equation}
For completeness, the lower inequality follows from the rational skeletal
comparison theorem \cite[Theorem~1]{FelixHalperinThomas2002}, applied after
replacing \(X_\Q\) by the indicated \(k\)-dimensional CW model.  This model is
the \(k\)-skeleton of the corresponding attachment model for \(Y_\Q\).  Its
inclusion is a rational homotopy isomorphism below \(k\), and all possible
homology comparison groups above the skeleton vanish because
\(H^{>N_X}(X;\Q)=0\).  The upper inequality follows by adding the cone on the
attaching sphere to a categorical decomposition of \(X\).

We write \(\operatorname{Cat}_0\) for normalized strong LS category after
rationalization.  Cornea's theorem identifies it with rational cone length:
\begin{equation}
                         \operatorname{Cat}_0(Z)=\Clz(Z)
\label{calc:eq:strong-cone}
\end{equation}
for every simply connected rational space of finite type
\cite{Cornea1995}.  Thus later occurrences of
\(\operatorname{Cat}_0\) use a theorem, not an additional hypothesis on
the cell attachment.

\section{Essential attachments}

Retain the hypotheses and notation
\eqref{calc:eq:three-thresholds}--\eqref{calc:eq:top-range}, and assume
that \(A_X\) admits an adapted decomposition
\[
                              V=W\oplus U\oplus Z
\]
through degree \(M\).  Let
\[
                        \chi_n=q(2n+1)-3\leq k\leq M.
\]
Suppose that the functional determined by \(f\) is dual to a nonzero
direction \(v\in U^k\): after choosing a basis,
\begin{equation}
 U^k=\langle u_1=v,u_2,\ldots,u_s\rangle,\qquad
 f(u_1)=1,\quad f(u_i)=0\ (i>1),\quad f(W^k)=0.
\label{calc:eq:adapted-attaching-functional}
\end{equation}

The category part of the essential-attachment theorem depends on the
following rank argument.  It is stronger than the assertion that one named
generator is essential: it uses the entire obstruction space \(U^j\).

\begin{lemma}[Relative-rank obstruction]
\label{calc:lem:relative-rank}
Under the preceding hypotheses,
\[
                              \catz(Y)>n.
\]
\end{lemma}

\begin{proof}
Suppose, to the contrary, that \(\catz(Y)\leq n\).  The
F\'elix--Halperin criterion, applied to the word-length projection of
\(A_Y\), gives a relative Sullivan factorization
\[
 \begin{tikzcd}[row sep=large,column sep=large]
 P_Y \arrow[r,two heads,"\pi","\sim"']
     \arrow[d,shift left=0.8ex,"\rho"]
   & A_Y/\mathfrak m_{A_Y}^{n+1}\\
 A_Y \arrow[u,hook,shift left=0.8ex,"i"]
     \arrow[ur,"q_n"']
   & {}
 \end{tikzcd}
 \qquad \rho i=\id_{A_Y},
\]
where \(\pi i=q_n\); see
\cite[Theorem~4.7]{FelixHalperin1982}.  Put \(K=\ker\pi\).  Since
\(\pi\) is a surjective quasi-isomorphism, \(K\) is an acyclic
differential ideal.

For \(c_n\leq j\leq k\), we construct CDGA morphisms
\[
                              \lambda_j:B_j\longrightarrow P_Y
\]
such that
\begin{enumerate}[label=\textup{(R\arabic*)}]
\item \(\lambda_j(w)=i(w)\) for \(w\in W^{\leq j}\);
\item \(\lambda_j(U^{<j})\subseteq K\);
\item \(\psi_j:=\rho\lambda_j:B_j\to B_j\) is a CDGA automorphism.
\end{enumerate}
Since \(U^{<c_n}=0\), begin with \(\lambda_{c_n}=i\) on
\(B_{c_n}=\Lambda W^{\leq c_n}\).

Assume that \(\lambda_j\) has been constructed for \(j<k\), and choose
a basis \(u_{j1},\ldots,u_{jt}\) of \(U^j\).  The adapted differential
gives
\[
 du_{j\alpha}\in\bigl((\Lambda W)^+\bigr)^{n+1}
                 +\Lambda W\otimes\Lambda^+(U^{<j}).
\]
The first summand is killed by \(\pi i=q_n\); every term in the second
has a factor whose \(\lambda_j\)-image lies in the ideal \(K\).  Hence
\(\lambda_j(du_{j\alpha})\) is a cycle in \(K\).  Choose
\(t_\alpha\in K^j\) such that
\begin{equation}
                      Dt_\alpha=\lambda_j(du_{j\alpha}).
\label{calc:eq:rank-primitive}
\end{equation}
Extend \(\lambda_j\) by \(\lambda_{j+1}(u_{j\alpha})=t_\alpha\), and
send new \(W\)-generators by \(i\).  This is compatible with the
differential because \((\Lambda W,d)\) is a sub-CDGA.

Below degree \(k\), the models of \(X\) and \(Y\) agree.  Since there
are no \(Z\)-generators in this range, write
\[
                 \rho(t_\alpha)=\sum_{\beta=1}^t
                    b_{\beta\alpha}u_{j\beta}+\eta_\alpha,
                 \qquad \eta_\alpha\in B_j^j.
\]
Applying the differential and passing to \(H^{j+1}(B_j)\) gives
\[
              H(\psi_j)[du_{j\alpha}]
                 =\sum_{\beta=1}^t b_{\beta\alpha}[du_{j\beta}].
\]
Lemma~\ref{calc:lem:obstruction-isomorphism} and the inductive
invertibility of \(H(\psi_j)\) show that
\((b_{\beta\alpha})\) is invertible.  Thus \(\psi_{j+1}\) is a
triangular CDGA automorphism.  This completes the induction up to degree
\(k\).

At degree \(k\), the same argument shows that every
\(\lambda_k(du_i)\) is a cycle in \(K\).  Choose \(t_i\in K^k\) with
\(Dt_i=\lambda_k(du_i)\).  The attachment removes \(u_1=v\), whereas
\(u_2,\ldots,u_s\) remain.  Therefore
\[
                   \rho(t_i)=\eta_i+\sum_{h=2}^s c_{hi}u_h,
                   \qquad\eta_i\in B_k^k.
\]
After differentiation and passage to \(H^{k+1}(B_k)\),
\begin{equation}
                   H(\psi_k)[du_i]
                      =\sum_{h=2}^s c_{hi}[du_h],
                   \qquad 1\leq i\leq s.
\label{calc:eq:rank-contradiction}
\end{equation}
The left-hand side is a basis of \(H^{k+1}(B_k)\), by
Lemma~\ref{calc:lem:obstruction-isomorphism}; the right-hand side lies in
the span of only \(s-1\) elements.  This contradiction proves the lemma.
\end{proof}

\begin{theorem}[Essential attachment]
\label{calc:thm:essential-attachment}
Under the hypotheses above,
\[
 \Hnilz(Y)>n,
 \qquad
 \catz(Y)=\operatorname{Cat}_0(Y)=\Clz(Y)=n+1.
\]
\end{theorem}

\begin{proof}
The attachment model omits \(u_1\) in degree \(k\) and retains
\(u_2,\ldots,u_s\).  Lemma~\ref{calc:lem:obstruction-isomorphism}
therefore identifies the cellular top cohomology as
\begin{equation}
 H^{k+1}(Y;\Q)\cong
 \frac{H^{k+1}(B_k)}{\langle[du_2],\ldots,[du_s]\rangle}
 \cong\Q[du_1].
\label{calc:eq:essential-top-quotient}
\end{equation}

Suppose \(\Hnilz(Y)\leq n\), and choose an acyclic ideal
\(I\lhd A_Y\) containing \(\mathfrak m_{A_Y}^{n+1}\).  The induction
in the proof of Theorem~\ref{calc:thm:essential-directions}, performed
only below degree \(k\), gives an automorphism
\(\phi_k:B_k\to B_k\) with \(\phi_k(K_k)\subseteq I\).  Thus all
\(\phi_k(du_i)\) are cycles in \(I\).  Their classes form a basis of
\(H^{k+1}(B_k)\), so at least one has nonzero image in the
one-dimensional quotient \eqref{calc:eq:essential-top-quotient}.  The
acyclic ideal \(I\) would then contain a cycle representing a nonzero
class of \(A_Y\), which is impossible.  Hence \(\Hnilz(Y)>n\).

Lemma~\ref{calc:lem:relative-rank} and the one-cell upper bound give
\(\catz(Y)=n+1\).  Proposition~\ref{calc:prop:cohomological-truncation} gives
\(\Hnilz(X)=n\); hence
\[
       n=\catz(X)\leq\Clz(X)\leq\Hnilz(X)=n.
\]
Thus \(\Clz(X)=n\).  Cone length rises by at most one under a one-cell
attachment, and Cornea's equality of rational cone length and strong
category \cite{Cornea1995} now gives
\[
 n+1=\catz(Y)\leq\operatorname{Cat}_0(Y)=\Clz(Y)
       \leq\Clz(X)+1=n+1.
\]
\end{proof}

\section{A filtered source of decomposable primitives}

The upper bounds require a way to choose primitives without a linear part.
We record the exact filtered input.

Let \(A=(\Lambda V,d)\), let \(n=\catz(|A|)\), and put
\(A_n=A/\mathfrak m_A^{n+1}\).  The standard semifree
\index{semifree module}
\(A\)-module resolution of the word-length quotient has the form
\begin{equation}
 Q_n=\bigl(A\otimes(\Q\oplus M_n),D\bigr)
       \xrightarrow{\simeq}A_n.
\label{calc:eq:semifree-resolution}
\end{equation}
It carries a filtration in which a word of length \(\ell\) in \(A\) has
filtration degree \(\ell\), and every element of \(M_n\) has filtration
degree \(n\).
Its acyclic kernel contains \(\mathfrak m_A^{n+1}\), and the filtered
contraction has the property
\begin{equation}
 \left.
 \begin{array}{c}
  w\in Q_n^{>t,*},\quad Dw=0,\quad t\geq n
 \end{array}
 \right\}
 \Longrightarrow
 \left\{
 \begin{array}{c}
  w=Du\text{ for some }u\in Q_n^{\geq t,*}.
 \end{array}
 \right.
\label{calc:eq:filtered-contraction}
\end{equation}
These are the filtration statements in
\cite[Section~29(f), Lemmas~29.11--29.12]{FHT2001}.

Hess's theorem identifying module category with rational LS category
\cite[Theorem~1]{Hess1991} (see also \cite[Theorem~29.9]{FHT2001}) provides
an \(A\)-linear chain map
\[
                              \rho:Q_n\longrightarrow A,
                              \qquad \rho|_A=\id_A.
\]
The usual retraction is initially obtained up to homotopy.  Evaluating the
chain homotopy on the unit and using \(A^{-1}=0\) shows that \(\rho(1)=1\);
\(A\)-linearity then gives the displayed strict equality on \(A\).

\begin{lemma}[Filtered primitive lemma]
\label{calc:lem:filtered-primitives}
Every cocycle \(z\in\mathfrak m_A^{n+2}\) has a decomposable primitive
in \(A\).  More precisely, there is \(\eta\in\mathfrak m_A^2\) with
\(d\eta=z\).

If \(n\geq2\) and, in addition, the resolution and retraction can be chosen
so that
\begin{equation}
                              \rho(M_n)\subseteq\mathfrak m_A^2,
\label{calc:eq:R2}
\end{equation}
then every cocycle \(z\in\mathfrak m_A^{n+1}\) also has a primitive in
\(\mathfrak m_A^2\).
\end{lemma}

\begin{proof}
Let \(z\in\mathfrak m_A^{n+2}\) be a cocycle.  Viewed in the algebra summand
\(A\subset Q_n\), it has filtration degree greater than \(n+1\).  Apply
\eqref{calc:eq:filtered-contraction} with \(t=n+1\) and write the
resulting primitive as
\[
                   u=u_0+\sum_i a_i\otimes\mu_i,
                   \qquad Du=z.
\]
The filtration gives \(u_0\in\Lambda^{\geq n+1}V\) and
\(a_i\in\mathfrak m_A\).  The connected resolution has
\(|\mu_i|>0\), so \(\rho(\mu_i)\in\mathfrak m_A\).  Therefore
\[
             \eta=\rho(u)=u_0+\sum_i a_i\rho(\mu_i)
                    \in\mathfrak m_A^2,
             \qquad d\eta=z.
\]

For \(z\in\mathfrak m_A^{n+1}\), apply
\eqref{calc:eq:filtered-contraction} with \(t=n\).  A coefficient of an
\(M_n\)-term may now be scalar.  Condition~\eqref{calc:eq:R2} nevertheless
ensures that the image of such a term under \(\rho\) is decomposable; all
terms with positive-degree coefficients are handled as above.  Also
\(u_0\in\Lambda^{\geq n}V\subseteq\mathfrak
  m_A^2\) because \(n\geq2\).  Hence the second assertion
follows.
\end{proof}

Condition \eqref{calc:eq:R2} is the precise strengthened
filtered-retraction hypothesis needed for preservation.  We state it
directly as a property of the chosen resolution and retraction.

\section{The quadratic one-step bound}

We say that the quadratic shadow has no spherical class in degree \(j\)
when
\[
                      d_2:V^j\longrightarrow(\Lambda^2V)^{j+1}
\]
is injective.  This condition makes every primitive needed in the
one-degree elimination decomposable.

\begin{theorem}[Quadratic one-step bound]
\label{calc:thm:quadratic-one-step}
Let \(X\) be \(r\)-connected, of finite rational type, and with bounded
rational cohomology.  Put \(q=r+1\), and assume
\[
 n=\catz(X)\geq3,\qquad N_X<q(n+1).
\]
Let \(f:S^k\to X\) be rationally nontrivial, with
\[
                    k\geq q(2n+1)-3,
                    \qquad Y=X\cup_f e^{k+1}.
\]
If \(d_2\) is injective in every degree
\begin{equation}
                         q(n+3)-2\leq j\leq k,
\label{calc:eq:quadratic-window}
\end{equation}
then
\[
                              \Hnilz(Y)\leq n+1.
\]
\end{theorem}

\begin{proof}
Set
\[
        s_0=q(n+2),\qquad T=q(n+3)-2,
        \qquad \ell=T+1.
\]
Since \(q\geq2\) and \(n\geq3\), the lower bound on \(k\) implies
\(k\geq\ell\).  It also implies \(N_X<k\), so
\eqref{calc:eq:attachment-cohomology} applies.  In \(A_Y\), begin with
\[
                            I_0=\mathfrak m_{A_Y}^{n+2}.
\]
It vanishes below degree \(s_0\) and lies in
\(\mathfrak m_{A_Y}^3\).  By \eqref{calc:eq:word-stability}, \(I_0\) agrees
with \(\mathfrak m_{A_X}^{n+2}\) in all degrees needed to compute its
cohomology through degree \(T\).

For every \(s_0\leq s\leq T\), choose cocycles
\(\psi_{s,i}\) representing a basis of \(H^s(I_0)\).  By
Lemma~\ref{calc:lem:filtered-primitives}, they have decomposable primitives
\(\eta_{s,i}\in A_X^{s-1}\).  Since \(s-1<T<k\), these primitives lie
in the common part of the models of \(X\) and \(Y\).  The equality
\[
                         (s_0-1)+q=T+1
\]
allows Lemma~\ref{calc:lem:first-slice} to kill the entire first slice.

Continue one degree at a time from \(\ell\) through \(k\).  By
\eqref{calc:eq:attachment-cohomology}, every ideal cycle in these degrees
is a boundary in \(A_Y\), and such a cycle lies in
\(\mathfrak m_{A_Y}^3\).  If \(\sigma\in A_Y^{N-1}\) is a primitive and
\(x\) is its linear part, comparison of quadratic terms gives \(d_2x=0\).
The injectivity assumption in degree \(N-1\) forces \(x=0\), so
Lemma~\ref{calc:lem:one-degree} applies at every stage.

At the end, there is an ideal \(I\supseteq\mathfrak m_{A_Y}^{n+2}\)
which is acyclic through degree \(k\) and satisfies
\(I^{k+1}\subseteq\mathfrak m_{A_Y}^3\).  Suppose a cycle in
\(I^{k+1}\) represented a nonzero multiple of the cellular class
\([dv]\).  After rescaling, it could be written
\[
                              z=dv+\delta\alpha.
\]
Let \(\alpha_1\in\widetilde V^k=W_f^k\) be the linear part of
\(\alpha\).  Since \(z\) has no quadratic component,
\[
                              d_2(v+\alpha_1)=0.
\]
But \(v+\alpha_1\neq0\) in
\(V^k=\Q v\oplus W_f^k\), contradicting injectivity of \(d_2\) in
degree \(k\).  Lemma~\ref{calc:lem:top-completion} therefore completes
\(I\) to an acyclic ideal containing
\(\mathfrak m_{A_Y}^{n+2}\).  Hence \(\Hnilz(Y)\leq n+1\).
\end{proof}

\section{Preservation under a stronger retraction}

\begin{theorem}[Preservation theorem]
\label{calc:thm:preservation}
Assume all the hypotheses of
Theorem~\ref{calc:thm:quadratic-one-step}.  Suppose additionally that
\(X_\Q\) admits a CW model of dimension at most \(k\), that the filtered
resolution admits a retraction satisfying
\eqref{calc:eq:R2}, and that \(d_2\) is injective in every degree
\begin{equation}
                         q(n+2)-2\leq j\leq k.
\label{calc:eq:preservation-window}
\end{equation}
Then
\[
 \catz(Y)=\operatorname{Cat}_0(Y)=\Clz(Y)=\Hnilz(Y)
          =n=\catz(X).
\]
\end{theorem}

\begin{proof}
Repeat the preceding proof with
\[
 I_0'=\mathfrak m_{A_Y}^{n+1},\qquad
 s_0'=q(n+1),\qquad T'=q(n+2)-2,
 \qquad\ell'=T'+1.
\]
For a cocycle in \(\mathfrak m_{A_X}^{n+1}\), the second part of
Lemma~\ref{calc:lem:filtered-primitives} supplies a decomposable
primitive.  Since \((s_0'-1)+q=T'+1\), the first-slice elimination
applies.  The shifted injectivity interval
\eqref{calc:eq:preservation-window} makes every later primitive
decomposable and excludes the cellular top class exactly as in the proof
of Theorem~\ref{calc:thm:quadratic-one-step}.
Criterion~\ref{calc:crit:acyclic-truncation} gives an acyclic ideal containing
\(\mathfrak m_{A_Y}^{n+1}\), and hence
\[
                              \Hnilz(Y)\leq n.
\]
The category comparison \eqref{calc:eq:one-cell-category} and the general
inequality \(\catz(Y)\leq\Clz(Y)\leq\Hnilz(Y)\) now give
\[
 n=\catz(X)\leq\catz(Y)\leq\Clz(Y)
       \leq\Hnilz(Y)\leq n.
\]
Thus \(\catz(Y)=\Clz(Y)=\Hnilz(Y)=n\), and
\(\operatorname{Cat}_0(Y)=\Clz(Y)\) by \cite{Cornea1995}.
\end{proof}

\begin{remark}
Every degree in the intervals
\eqref{calc:eq:quadratic-window} and
\eqref{calc:eq:preservation-window} is used.  The degrees below \(k\)
force the inductively selected primitives to be decomposable; degree
\(k\) keeps the cellular top class out of the ideal.  These are sufficient
conditions, not asserted to be necessary.
\end{remark}

\section{An exact value: a projective-space attachment}

We now apply the lower- and upper-bound methods together.  Let
\[
                              X=\mathbb CP^3\vee\mathbb CP^2.
\]
Its rational cohomology algebra is
\begin{equation}
 H^*(X;\Q)=\Q[a,b]/(ab,a^4,b^3),
                    \qquad |a|=|b|=2.
\label{calc:eq:projective-cohomology}
\end{equation}
The minimal-model construction, carried out through degree eleven, gives
generators
\[
\begin{gathered}
 a_2,b_2,x_3,y_5,z_6,t_7,w_7,s_8,r_8,u_9,v_9,\\
 q_{0,10},q_{1,10},q_{2,10},
 \mu_{0,11},\mu_{1,11},\mu_{2,11},\mu_{3,11},\mu_{4,11},
\end{gathered}
\]
and differential
\begin{align*}
 da&=db=0,             &dx&=ab,                 &dy&=b^3,\\
 dz&=b^2x-ay,          &dt&=a^4,                &dw&=xy+bz,\\
 ds&=a^3x-bt,          &dr&=xz-aw,              &du&=xt+as,\\
 dv&=xw+br,            &dq_0&=yz-b^2w,          &dq_1&=xs-bu,\\
 dq_2&=xr-av,\\[0.2em]
 d\mu_0&=a^3z-yt-b^2s,
 &d\mu_1&=z^2-2b^2r+2aq_0,\\
 d\mu_2&=yw+bq_0,
 &d\mu_3&=xu+aq_1,
 &d\mu_4&=xv+bq_2.
\end{align*}
Here and below all degree subscripts are suppressed inside formulas; the
remaining subscripts distinguish generators rather than degrees.

These formulas satisfy \(d^2=0\).  For example,
\[
 d(a^3z-yt-b^2s)
 =a^3(b^2x-ay)-b^3t+a^4y-b^2(a^3x-bt)=0.
\]
The remaining identities follow directly from the derivation rule; for
instance
\[
\begin{split}
 d(xy+bz)&=aby-xb^3+b^3x-aby=0,\\
 d(z^2-2b^2r+2aq_0)
 &=2(b^2x-ay)z-2b^2(xz-aw)+2a(yz-b^2w)=0.
\end{split}
\]
\subsection{A finite-rank verification}

We record the exact linear algebra used in the minimal-model construction and
the spherical calculation.  Order the displayed
generators by
\begin{align*}
 a&<b<x<y<z<t<w<s<r<u<v,\\
 q_0&<q_1<q_2<\mu_0<\mu_1<\mu_2<\mu_3<\mu_4,
\end{align*}
with every generator on the first line preceding those on the second.  For a
set \(S\) of generators, let \(\mathcal M^N(S)\) denote the ordered monomial
basis of \((\Lambda S)^N\).  Exponents of odd generators lie in
\(\{0,1\}\), and exponent vectors are ordered lexicographically with respect
to the displayed generator order.  Thus the matrix of the differential is
completely specified, without any choice, by
\begin{equation}
 d(v_1^{e_1}\cdots v_r^{e_r})
 =\sum_{i=1}^r(-1)^{\sum_{h<i}e_h|v_h|}
 e_i\,v_1^{e_1}\cdots v_i^{e_i-1}(dv_i)
       v_{i+1}^{e_{i+1}}\cdots v_r^{e_r},
\label{calc:eq:monomial-differential}
\end{equation}
where \(e_i\in\{0,1\}\) whenever \(v_i\) is odd.  All the matrices below therefore
have entries in \(\{0,\pm1,\pm2\}\), and ordinary rational row reduction
gives the following rank verification.

Let \(C_{<j}\) be the algebra on the displayed generators of degree
strictly less than \(j\).  Put
\[
\begin{aligned}
c_j&=\dim C_{<j}^{j+1},\\
r_j^-&=\operatorname{rank}(d:C_{<j}^{j}\to C_{<j}^{j+1}),\\
r_j^+&=\operatorname{rank}(d:C_{<j}^{j+1}\to C_{<j}^{j+2}).
\end{aligned}
\]
If \(\phi_{<j}:C_{<j}\to H^*(X;\Q)\) fixes \(a,b\) and kills the other
generators, write \(\rho_j\) for the rank of
\(H^{j+1}(\phi_{<j})\) and \(\kappa_j\) for the dimension of its kernel.
Row reduction in the bases \(\mathcal M^N\) gives the rank data in
Table~\ref{calc:tab:minimal-prefix-ranks}.
\begin{table}[tbp]
\centering
\caption{Rank data for the inductive minimal-model construction through degree eleven}
\label{calc:tab:minimal-prefix-ranks}
\[
\begin{array}{c|rrrrrrrrr}
j       &3&4&5&6&7&8&9&10&11\\ \hline
c_j     &3&2&4&5&8&12&17&24&35\\
r_j^-   &0&0&2&0&4&2&8&7&14\\
r_j^+   &0&2&0&4&2&8&7&14&16\\
\dim H^{j+1}(C_{<j})
        &3&0&2&1&2&2&2&3&5\\
\rho_j &2&0&1&0&0&0&0&0&0\\
\kappa_j&1&0&1&1&2&2&2&3&5
\end{array}
\]
\end{table}
Row-reduced bases of these kernels are, in increasing \(j\),
\[
\begin{array}{c|l}
3 &[ab]\\
5 &[b^3]\\
6 &[b^2x-ay]\\
7 &[a^4],\ [xy+bz]\\
8 &[a^3x-bt],\ [xz-aw]\\
9 &[xt+as],\ [xw+br]\\
10&[yz-b^2w],\ [xs-bu],\ [xr-av]\\
11&[a^3z-yt-b^2s],\ [z^2-2b^2r+2aq_0]\\
  &[yw+bq_0],\ [xu+aq_1],\ [xv+bq_2]
\end{array}
\]
There is no kernel at \(j=4\).  The five positive-dimensional target
classes \(a,b,a^2,b^2,a^3\) are already represented, so there is no
cohomological cokernel requiring any further generator through degree
eleven.  The standard inductive minimal-model construction consequently
adjoins exactly the generators listed above: each degree-\(j\) generator
kills the corresponding displayed basis class in degree \(j+1\).  This
verifies that the generator list is complete through degree eleven.

The morphism from the displayed Sullivan algebra to the algebra in
\eqref{calc:eq:projective-cohomology}, fixing \(a,b\) and killing all other
generators, extends over the later generators to a surjective
quasi-isomorphism.  Hence \(X\) is formal.  Its rational cup-length is three,
so \(\catz(X)=3\).  The cohomology quotient is a three-short acyclic quotient
of the minimal model, which gives \(\Hnilz(X)\leq3\);
Proposition~\ref{stab:prop:cat-lower} gives equality.  Therefore
\begin{equation}
                \catz(X)=\Hnilz(X)=3,
                \qquad N_X=6,\qquad q=2,\qquad c_3=7,\qquad\chi_3=11.
\label{calc:eq:projective-invariants}
\end{equation}

Through degree eleven, define
\begin{align*}
 W^{\leq11}
 &=\langle a,b,x,y,z,w,r,v,q_0,q_2,
              \mu_1,\mu_2,\mu_4\rangle,\\
 U^{\leq11}
 &=\langle t,s,u,q_1,\mu_0,\mu_3\rangle.
\end{align*}
The displayed differentials show that \((\Lambda W,d)\) is a sub-CDGA.
They also give
\[
        dU^j\subseteq\bigl((\Lambda W)^+\bigr)^4
             +\Lambda W\otimes\Lambda^+(U^{<j})
             \qquad(7\leq j\leq11).
\]
It remains to verify the spherical condition in the same finite range.  The
active generator sets for the successive algebras
\(B_j\) are
\begin{align*}
S_7&=(a,b,x,y,z,w),&
S_8&=S_7\cup(t,r),\\
S_9&=S_8\cup(s,v),&
S_{10}&=S_9\cup(u,q_0,q_2),\\
S_{11}&=S_{10}\cup(q_1,\mu_1,\mu_2,\mu_4).
\end{align*}
Using the monomial order and differential matrix
\eqref{calc:eq:monomial-differential}, the corresponding rank verification
is recorded in Table~\ref{calc:tab:adapted-stage-ranks}.
\begin{table}[tbp]
\centering
\caption{Rank data for the successive adapted algebras \(B_j\)}
\label{calc:tab:adapted-stage-ranks}
\[
\begin{array}{c|rrrrr}
j&7&8&9&10&11\\ \hline
\dim B_j^{j+1}&8&12&17&24&35\\
\operatorname{rank}(d:B_j^j\to B_j^{j+1})&5&3&9&9&17\\
\operatorname{rank}(d:B_j^{j+1}\to B_j^{j+2})&2&8&7&14&16\\
\dim H^{j+1}(B_j)&1&1&1&1&2
\end{array}
\]
\end{table}
Row-reduced bases for the cohomology groups recorded in the final row are,
respectively,
\[
[a^4],\quad [a^3x-bt],\quad [xt+as],\quad [xs-bu],
\quad [d\mu_0],\ [d\mu_3].
\]
At each stage these are exactly the differentials of the displayed basis
of \(U^j\).  In particular, no nonzero linear combination of the
\(U^j\)-generators can be corrected by an element of \(B_j^j\) to form a
cocycle.  Thus all spherical directions through degree eleven lie in
\(W\).

It remains to extend the displayed finite decomposition to the whole minimal
model, because condition~\textup{(A1)} is global.  Apply
Lemma~\ref{calc:lem:spherical-normalization} inductively in every degree
strictly greater than eleven.  Add each higher spherical direction, represented
by a closed generator, to \(W\), and place every remaining higher generator
in \(Z\); set \(U^j=0\) for \(j>11\).  The newly added generators of \(W\)
are closed, so \((\Lambda W,d)\) remains a sub-CDGA and contains the entire
image of \(h_A\).  No condition on the differential of the new \(Z\)-block is
required above the cutoff.  Consequently conditions \textup{(A1)--(A4)} hold
globally for an adapted decomposition through \(M=11\).  The direction
\(\mu_0\in U^{11}\) is therefore three-essential.

Attach a cell along the rational homotopy class dual to \(\mu_0\):
\[
                    Y=X\cup_{\mu_0}e^{12}.
\]
Passing to the minimal model removes \(\mu_0\), while
\[
                         \Omega=d\mu_0=a^3z-yt-b^2s
\]
becomes the new degree-twelve class.  In particular,
\[
 H^*(Y;\Q)\cong
 \Q\oplus\langle[a],[b]\rangle
 \oplus\langle[a^2],[b^2]\rangle
 \oplus\langle[a^3]\rangle
 \oplus\langle[\Omega]\rangle.
\]
Theorem~\ref{calc:thm:essential-attachment} gives
\begin{equation}
       \Hnilz(Y)>3,
       \qquad \catz(Y)=\operatorname{Cat}_0(Y)=\Clz(Y)=4.
\label{calc:eq:projective-lower}
\end{equation}

It remains to prove the matching upper bound.  For \(n=3\), \(q=2\),
and \(k=11\), the quadratic window
\eqref{calc:eq:quadratic-window} consists exactly of degrees ten and
eleven.  In degree ten,
\[
 d_2q_0=yz,\qquad d_2q_1=xs-bu,\qquad d_2q_2=xr-av.
\]
These three vectors are independent because the displayed monomials are
distinct basis monomials of the free graded-commutative algebra.  In degree
eleven,
\[
\begin{array}{lll}
 d_2\mu_0=-yt,&
 d_2\mu_1=z^2+2aq_0,&
 d_2\mu_2=yw+bq_0,\\
 d_2\mu_3=xu+aq_1,&
 d_2\mu_4=xv+bq_2.&
\end{array}
\]
The five columns are again independent: each contains a
monomial---respectively, \(yt,z^2,yw,xu,xv\)---that occurs in no other
column.  Thus \(d_2\) is injective in both required degrees.

\begin{theorem}[A projective-space calculation]
\label{calc:thm:projective-exact}
For
\[
             Y=(\mathbb CP^3\vee\mathbb CP^2)
                        \cup_{\mu_0}e^{12},
\]
one has
\[
       \boxed{\ \catz(Y)=\operatorname{Cat}_0(Y)=\Clz(Y)
                    =\Hnilz(Y)=4\ }.
\]
\end{theorem}

\begin{proof}
Equation~\eqref{calc:eq:projective-lower} gives the category values and
the strict lower bound \(\Hnilz(Y)>3\).  The two exact quadratic-rank
computations above verify every hypothesis of
Theorem~\ref{calc:thm:quadratic-one-step}, which gives
\(\Hnilz(Y)\leq4\).
\end{proof}

\chapter{Low-Length Strictification and Retractive Closure}\label{chap:why-close}
\begin{chapterguide}[title={Chapter guide}]
This chapter isolates the positive strictification theorems.  In the simply
connected finite-type minimal setting, every flexible witness strictifies
through length two.  Rational category is recovered from either invariant
after taking retractive closure; Part~\ref{part:stabilization} will establish
the parallel spherical result.  These conclusions delimit the first
length at which a strict quotient defect can occur while keeping pointwise
and retractive statements logically separate.
\end{chapterguide}
\index{strictification!through length two}
\index{linear contractible part}
\index{retractive closure}
\index{spherical stabilization}

\section{Removing the linear contractible part}

The next lemma explains why one may study short witnesses whose
differential has no linear part.

\begin{lemma}[Stabilization of indecomposables]
\label{div:lem:indecomposable-stabilization}
Suppose that
\[
                    B_0\supseteq B_1\supseteq B_2\supseteq\cdots
\]
is a descending sequence of connected sub-CDGAs, obtained at each step
by splitting
\[
     \indec(B_i)=H_i\oplus P_i,
     \qquad dH_i=0,\qquad P_i\text{ contractible},
\]
and setting \(B_{i+1}=\Q\oplus q_i^{-1}(H_i)\), where
\(q_i:B_i^+\to\indec(B_i)\) is the quotient.  Assume that every
\(B_i^k\) is finite dimensional and put \(B_\infty=\bigcap_iB_i\).
If, for some \(i\) and \(q\),
\[
             B_j^k=B_\infty^k
       \quad\text{for all }j\geq i\text{ and }0\leq k\leq q+1,
\]
then the degree-\(q\) and degree-\((q+1)\) components of
\(\indec(B_j)\), together with the differential between them, agree
with those of \(\indec(B_\infty)\).  Consequently, a nonzero
differential on \(\indec(B_\infty)\) from degree \(q\) to degree
\(q+1\) would force
\[
             B_{j+1}^q\subsetneq B_j^q\quad\text{for every }j\geq i,
\]
in particular for \(j=i\).
\end{lemma}

\begin{proof}
For \(k\leq q+1\), every product in \((B_j^+)^2\cap B_j^k\) is a sum
of products of homogeneous elements whose degrees lie between \(1\)
and \(k-1\).  All these homogeneous pieces have stabilized.  Hence
both \(B_j^k\) and \((B_j^+)^2\cap B_j^k\) agree with their counterparts
in \(B_\infty\), and so do the quotients
\(\indec(B_j)^k\) for \(k=q,q+1\).  The differentials agree because the
underlying subcomplexes agree in these two degrees.

If the resulting differential from degree \(q\) is nonzero, then, for
every \(j\geq i\), the contractible summand \(P_j\) has a nonzero
component in degree \(q\).  Passing from \(B_j\) to
\(B_{j+1}=\Q\oplus q_j^{-1}(H_j)\) removes a representative of that
component.  Thus the displayed inclusion is strict, as claimed.
\end{proof}

\begin{lemma}[Linear reduction]
\label{div:lem:linear-reduction}
Let \(B\) be a connected augmented CDGA of finite type whose augmentation
ideal is nilpotent; thus \((B^+)^{N+1}=0\) for some \(N\geq 0\).  There
exists a sub-CDGA
\[
                         B_{\mathrm{red}}\hookrightarrow B
\]
that is quasi-isomorphic to \(B\), has augmentation length no greater than
that of \(B\), and satisfies
\[
                 d(B_{\mathrm{red}}^+)\subseteq
                    (B_{\mathrm{red}}^+)^2.
\]
\end{lemma}

\begin{proof}
Let \(q:B^+\to\indec(B)=B^+/(B^+)^2\) be the quotient map.  Over
\(\Q\), split the complex of indecomposables as
\[
                      \indec(B)=H\oplus P,
\]
where \(H\) has zero differential and \(P\) is contractible.  Define
\[
                 B_1=\Q\oplus q^{-1}(H).
\]
Because \(q^{-1}(H)\) contains \((B^+)^2\), the subspace \(B_1\) is
closed under multiplication.  It is closed under the differential
because \(H\) is a subcomplex.  Hence \(B_1\) is a sub-CDGA of \(B\).
Furthermore,
\[
                           B/B_1\cong P
\]
as complexes.  Since \(P\) is contractible, the inclusion
\(B_1\hookrightarrow B\) is a quasi-isomorphism.

The operation is repeated because an element that was decomposable in
\(B\) can become indecomposable after the contractible summand has been
removed.  This produces a descending sequence
\[
                         B\supseteq B_1\supseteq B_2\supseteq\cdots
\]
of quasi-isomorphic sub-CDGAs.  The finite-type hypothesis implies
that, for every \(q\), the descending sequence of vector spaces \(B_i^q\)
stabilizes.  Put \(B_{\mathrm{red}}=\bigcap_i B_i\).  To compute
\(H^q\), only degrees \(q-1,q,q+1\) are involved; choosing one index at
which all three have stabilized shows that
\(H^q(B_{\mathrm{red}})\to H^q(B)\) is an isomorphism.  Thus the
intersection is still quasi-isomorphic to \(B\).

Fix \(q\), and choose \(i\) so that
\(B_j^k=B_{\mathrm{red}}^k\) for every \(j\geq i\) and every
\(0\leq k\leq q+1\).  Lemma~\ref{div:lem:indecomposable-stabilization}
shows that a nonzero differential from degree \(q\) to degree \(q+1\) on
\(\indec(B_{\mathrm{red}})\) would force the stabilized degree-\(q\) piece
to decrease at the next reduction step.  Therefore the induced
differential on \(\indec(B_{\mathrm{red}})\) is zero, which is equivalent to
\[
 d(B_{\mathrm{red}}^+)\subseteq(B_{\mathrm{red}}^+)^2.
\]
Finally, passing to a subalgebra cannot increase the nilpotency length
of the augmentation ideal.
\end{proof}

\section{Why lengths one and two strictify}

We now make precise the low-length phenomenon discussed in the
introduction.

\begin{lemma}[The square-zero case]
\label{div:lem:square-zero}
Let \(M\) be a minimal Sullivan algebra, let \(B\) be connected, and let
\(f:M\qis B\) be an augmented quasi-isomorphism.  If
\[
                         (B^+)^2=0,\qquad d(B^+)=0,
\]
then \(f\) is surjective and \(\Hnil(M)\leq1\).
\end{lemma}

\begin{proof}
Every positive element \(b\in B^+\) is a cycle.  Since \(H(f)\) is
surjective, there is a cycle \(m\in M^+\) such that
\([f(m)]=[b]\).  The differential of \(B\) is zero, so
\(f(m)=b\).  Thus \(f\) is surjective.  The conclusion follows from
Lemma~\ref{div:lem:surjective-characterization}.
\end{proof}

For the cube-zero case, the entire obstruction is linear.

\begin{lemma}[The cube-zero decomposition]
\label{div:lem:cube-decomposition}
Let \(B\) be a connected augmented CDGA satisfying
\[
                       (B^+)^3=0,\qquad d(B^+)\subseteq(B^+)^2,
\]
and let \(f:M\qis B\) be a quasi-isomorphism from a minimal Sullivan
algebra.  Put
\[
                         \mathfrak m=B^+,\qquad W=\mathfrak m^2.
\]
Choose a graded complement \(V\) so that
\(\mathfrak m=V\oplus W\).  Then \(dW=0\), and the differential on
\(\mathfrak m\) is determined by a linear map
\[
                              \delta:V\longrightarrow W.
\]
Let \(A=\im f\), let \(A_W=A\cap W\), let \(U\subseteq V\) be the image
of \(A^+\) under \(\mathfrak m\to\mathfrak m/W=V\), and choose
\(V=U\oplus N\).  Then
\begin{align}
\delta|_N&\text{ is injective},                                      \label{div:eq:delta-inj}\\
\delta N\cap\delta U&=0,                                             \label{div:eq:delta-disjoint}\\
W&=A_W+\delta N.                                                     \label{div:eq:W-sum}
\end{align}
\end{lemma}

\begin{proof}
Since \(d(\mathfrak m)\subseteq\mathfrak m^2\) and
\(\mathfrak m^3=0\), the Leibniz rule gives
\[
                         dW=d(\mathfrak m^2)=0.
\]
Thus the differential is indeed described by a linear map
\(\delta:V\to W\).  If \(u\in U\), choose a lift
\(a=u+w\in A^+\), with \(w\in W\).  Then
\[
                         da=\delta u\in A\cap W=A_W.
\]
Hence \(\delta U\subseteq A_W\).

Suppose \(0\neq n\in N\) and \(\delta n=0\).  Then \(n\) represents a
cohomology class of \(B\).  Because \(H(f)\) is surjective, so is
\(H(A)\to H(B)\); hence there exist a cycle \(a\in A\) and an element
\(v+w\in V\oplus W\) such that
\[
                           n-a=d(v+w)=\delta v.
\]
The right-hand side belongs to \(W\).  Projecting to \(V\) shows that
\(n\) lies in \(U\), contrary to \(V=U\oplus N\).  This proves
\eqref{div:eq:delta-inj}.

If \(\delta n=\delta u\) for \(n\in N\) and \(u\in U\), then
\(n-u\) is a cycle.  The same argument, after projection to \(V/U\),
forces \(n=0\).  Therefore \eqref{div:eq:delta-disjoint} holds.

Finally, let \(w\in W\).  It is a cycle.  Choose a cycle \(a\in A\)
representing its class in \(H(B)\).  Then
\[
                         w-a=d(v)=\delta v
\]
for some \(v=u+n\in U\oplus N\).  The \(V\)-component of the left-hand
side is the negative of the \(V\)-component of \(a\), whereas the
right-hand side lies in \(W\).  Hence the \(V\)-component of \(a\) is
zero, so \(a\in A_W\).  Since \(\delta u\in A_W\), we conclude that
\(w\in A_W+\delta N\), proving \eqref{div:eq:W-sum}.
\end{proof}

\begin{proposition}[Cube-zero strictification]
\label{div:prop:cube-strictification}
Under the hypotheses of Lemma~\ref{div:lem:cube-decomposition}, there is a
quotient \(C\) of \(A=\im f\) such that
\[
                         M\onto C
\]
is a quasi-isomorphism and \((C^+)^3=0\).  In particular,
\(\Hnil(M)\leq2\).
\end{proposition}

\begin{proof}
Set
\[
                            D=A_W\cap\delta N.
\]
Because \(D\subseteq W=\mathfrak m^2\) and \(\mathfrak m^3=0\), one has
\[
                             B^+D=0,\qquad dD=0.
\]
Thus \(D\) is a differential ideal of \(A\).

Choose a graded linear section
\(\sigma:U\to A^+\) of the projection \(A^+\to U\).  Then
\[
                         A^+=\sigma(U)\oplus A_W,
        \qquad d\sigma(u)=\delta u.
\]
Moreover, \eqref{div:eq:delta-disjoint} implies
\[
                         \delta U\cap D=0.
\]
Consequently, a class represented by \(\sigma(u)+w+D\) is a cycle in
\(A/D\) precisely when \(u\in\ker(\delta|_U)\), and its \(W\)-component
is defined modulo \(\delta U+D\).  Thus
\[
 H^+(A/D)\cong
 \ker(\delta|_U)\oplus\frac{A_W}{\delta U+D}.
\]

On the other hand, the decomposition in
Lemma~\ref{div:lem:cube-decomposition} gives
\[
 H^+(B)\cong
 \ker(\delta|_U)\oplus\frac{W}{\delta U+\delta N}.
\]
Indeed, \eqref{div:eq:delta-inj} and \eqref{div:eq:delta-disjoint} show that
the \(N\)-component contributes no cycles and cannot cancel a
\(\delta U\)-boundary.

Consider the linear map
\[
 \theta:\frac{A_W}{\delta U+D}
      \longrightarrow\frac{W}{\delta U+\delta N}
\]
induced by the inclusion \(A_W\subseteq W\).  It is surjective by
\eqref{div:eq:W-sum}.  Its kernel is
\[
\frac{A_W\cap(\delta U+\delta N)}{\delta U+D}
=\frac{\delta U+(A_W\cap\delta N)}{\delta U+D}=0.
\]
Hence \(\theta\) is an isomorphism.

For completeness, the resulting isomorphism on cohomology has an
intrinsic description.  If \(a+D\) is a cycle of \(A/D\), then the
preceding argument shows that \(da=0\).  Define
\[
              \Theta:H(A/D)\longrightarrow H(B),
              \qquad \Theta([a+D])=[a].
\]
This is well defined: changing \(a\) by an element of \(D\) changes it
by a boundary in \(B\), because \(D\subseteq\delta N=dN\).  More
explicitly, if \(a'\) is a cohomologous representative in \(A/D\), then
\[
                         a'-a=dc+d_0
\]
for some \(c\in A\) and \(d_0\in D\); both summands on the right are
boundaries in \(B\).

Relative to the two displayed decompositions of cohomology, the map
\(\Theta\) is upper triangular.  Its diagonal maps are the identity on
\(\ker(\delta|_U)\) and the isomorphism \(\theta\); the possible
off-diagonal term records the \(W\)-component of the chosen section
\(\sigma\).  Therefore \(\Theta\) is an isomorphism.  If
\(\bar f:M\to A/D\) denotes the
quotient of \(f\), then
\[
                         \Theta\circ H(\bar f)=H(f).
\]
Since \(H(f)\) and \(\Theta\) are isomorphisms, \(\bar f\) is a
quasi-isomorphism.  It is surjective by the definition
of \(A=\im f\).  Taking \(C=A/D\) gives the result.  Since \(C\) is a
quotient of a subalgebra of \(B\), it is still cube-zero.
\end{proof}

We shall use the following foundational short-replacement theorem of Cornea.

\begin{theorem}[Cornea's finite-type short replacement]
\label{div:lem:cornea-replacement}
Let \(A\) be a simply connected augmented CDGA of finite type, let
\(C\) be an augmented CDGA, and let \(f:A\to C\) be an augmented
quasi-isomorphism.  If \(n>0\) and
\[
                              (C^+)^{n+1}=0,
\]
then there exist a simply connected augmented CDGA \(A'\) of finite
type and a factorization by augmented morphisms
\[
                    A\xrightarrow{\,g\,}A'
                      \xrightarrow{\,h\,}C,
                    \qquad h g=f,
\]
such that \(g\) is a quasi-isomorphism and
\((A'^+)^{n+1}=0\).
\end{theorem}

\begin{proof}
This is the augmented-CDGA specialization of
\cite[Proposition~3.6, p.~107]{Cornea1994}.  The quasi-isomorphism in the
conclusion is the first arrow \(g:A\to A'\); the second arrow records the
factorization, and the middle algebra remains \(n\)-short.  Since \(f=hg\),
the second arrow is also a quasi-isomorphism by two-out-of-three.  This is
the orientation used below.
\end{proof}

\begin{corollary}[Strictification through length two]
\label{div:cor:through-two}
Let \(M\) be a simply connected minimal Sullivan algebra of finite type.
Then
\[
\begin{aligned}
\nilh(M)\leq1&\Longrightarrow\Hnil(M)\leq1,\\
\nilh(M)\leq2&\Longrightarrow\Hnil(M)\leq2.
\end{aligned}
\]
\end{corollary}

\begin{proof}
Start with an augmented quasi-isomorphism \(M\to B\) to an \(n\)-short
witness, where \(n=1\) or \(2\); its existence follows from
Lemma~\ref{div:lem:direct-zigzag}.  Since \(M\) is simply connected and of
finite type, Lemma~\ref{div:lem:cornea-replacement} supplies a factorization
through a simply connected \(n\)-short CDGA \(B'\) of finite type, with
\(M\to B'\) a quasi-isomorphism.  Apply
Lemma~\ref{div:lem:linear-reduction} to \(B'\), and denote the result by
\(B_{\mathrm{red}}\).  Let
\[
             \mu:M_{\mathrm{red}}\xrightarrow{\simeq}
                    B_{\mathrm{red}}
\]
be a minimal Sullivan model.  Composing \(\mu\) with
\(B_{\mathrm{red}}\hookrightarrow B'\) gives a minimal Sullivan model
of \(B'\).  Thus \(M\to B'\) and
\(M_{\mathrm{red}}\to B'\) are two minimal models of the same CDGA.
Uniqueness of simply connected minimal Sullivan models
\cite[Proposition~12.10(ii)]{FHT2001} gives an isomorphism
\(\alpha:M\cong M_{\mathrm{red}}\).  The composite
\[
                 M\xrightarrow{\alpha}M_{\mathrm{red}}
                    \xrightarrow{\mu}B_{\mathrm{red}}
\]
is the required augmented quasi-isomorphism to the reduced witness.
Only the abstract isomorphism \(\alpha\) is used here; the argument does
not require the two minimal-model morphisms to \(B'\) to commute
strictly through \(\alpha\).

If the witness is one-short, its differential is zero, and
Lemma~\ref{div:lem:square-zero} applies.  If it is two-short,
Proposition~\ref{div:prop:cube-strictification} applies.
\end{proof}

\begin{remark}
\label{div:rem:first-threshold}
Corollary~\ref{div:cor:through-two} explains why length three is the first
possible strict separation.  The proof of
Proposition~\ref{div:prop:cube-strictification} uses the fact that the top
layer \(W=(B^+)^2\) is annihilated by \(B^+\).  For a three-short algebra,
\(B^+(B^+)^2=(B^+)^3\) need not vanish.  Consequently, adding a primitive
to an ideal can force new nontrivial products into that ideal.
\end{remark}

\section{Intermediate witnesses and retractive closure}

We first establish the intermediate-witness characterization needed for
retractive closure.  For a simply connected augmented CDGA \(A\),
\(\Hnilz(A)\) denotes the invariant of any minimal Sullivan model of \(A\),
as specified in Chapter~\ref{stab:sec:alg-prelim}.

\begin{proposition}[Rational category through an intermediate
\(\Hnilz\)-witness]\label{stab:prop:intermediate-Hnil-witness}
Let \(X\) be simply connected and of finite type, and let \(M_X\)
be its minimal Sullivan model.  For every \(n\geq0\), the following
conditions are equivalent:
\begin{enumerate}
\item \(\catz(X)\leq n\);
\item there are augmented CDGAs \(A\) and \(C\), with \(A\) simply
      connected, and morphisms
      \begin{equation}\label{stab:eq:intermediate-Hnil-witness}
       M_X\xrightarrow{f}A\xrightarrow{g}C
      \end{equation}
      such that \(g f\) is a quasi-isomorphism and
      \(\Hnilz(A)\leq n\).
\end{enumerate}
The CDGA \(A\) in the second condition may moreover be chosen of finite
type.
\end{proposition}

\begin{proof}
Assume first that the second condition holds.  Choose a minimal Sullivan
model
\[
                    \mu:N=(\Lambda W,d)\qiso A.
\]
Since \(M_X\) is cofibrant, the morphism
\([\mu]^{-1}[f]\) in the homotopy category of augmented CDGAs is
represented by a morphism
\(\widetilde f:M_X\to N\).  Thus
\[
                     [g\mu\widetilde f]=[g f].
\]
The inequality \(\Hnilz(A)\leq n\) gives an acyclic differential ideal
\(J\triangleleft N\) with
\((N^+)^{n+1}\subseteq J\).  Put \(B=N/J\) and let
\(q:N\twoheadrightarrow B\) be the quotient map.  Then \(q\) is a
quasi-isomorphism and \(B\) is \(n\)-short.

Let
\[
 q_n:M_X\twoheadrightarrow
 A_n=M_X/(M_X^+)^{n+1}
\]
be the word-length projection.  Since \(q\widetilde f\) vanishes on
\((M_X^+)^{n+1}\), there is a unique morphism
\(\bar f:A_n\to B\) satisfying
\[
                         q\widetilde f=\bar f q_n.
\]
In the homotopy category define
\begin{equation}\label{stab:eq:intermediate-retraction}
 [r]=[g f]^{-1}[g\mu][q]^{-1}[\bar f]:
 A_n\longrightarrow M_X.
\end{equation}
Then
\[
 [r][q_n]
 =[g f]^{-1}[g\mu][q]^{-1}[q][\widetilde f]
 =[g f]^{-1}[g f]
 =[\id_{M_X}].
\]
Thus \(q_n\) admits a homotopy retraction.  The
F\'elix--Halperin criterion
\cite[Theorem~4.7]{FelixHalperin1982} gives
\(\catz(X)\leq n\).

Conversely, suppose \(\catz(X)\leq n\).  If \(n=0\), then \(X\) is
rationally contractible; take \(A=C=\Q\), let \(f:M_X\to\Q\) be the
augmentation, and let \(g=\id_{\Q}\).

Assume \(n\geq1\).  Apply
Lemma~\ref{stab:lem:decomposable-quotient-model} to
\[
 q_n:M_X\twoheadrightarrow A_n.
\]
It gives a factorization
\[
 M_X\xrightarrow{i_n}P_n
 \xrightarrow[\simeq]{\kappa_n}\!\!\!\twoheadrightarrow A_n
\]
in which \(P_n\) is minimal,
\(\kappa_n\) kills every relative generator, and
\(K_n=\ker\kappa_n\) is acyclic.  The
F\'elix--Halperin criterion supplies a strict retraction
\(\rho_n:P_n\to M_X\) of \(i_n\).  Moreover,
\[
                         (P_n^+)^{n+1}\subseteq K_n.
\]
A monomial of length at least \(n+1\) lies in \(K_n\): if it contains a
relative generator, \(\kappa_n\) kills it by normalization; otherwise it
belongs to \((M_X^+)^{n+1}\) and \(q_n\) kills it.  Hence
\(\Hnilz(P_n)\leq n\).  Taking
\[
                     A=P_n,\qquad C=M_X,\qquad
                     f=i_n,\qquad g=\rho_n
\]
gives \(g f=\id_{M_X}\) and proves the second condition.  The
degreewise finiteness assertion follows from
\cite[Proposition~12.2]{FHT2001}, as in
Theorem~\ref{stab:thm:root}.
\end{proof}

\begin{theorem}[Parallel retractive-closure theorem]
\label{stab:thm:parallel-retractive-closure}
Let \(X\) be a simply connected rational space of finite type.  For every
\(n\geq0\), the following conditions are equivalent:
\begin{enumerate}
\item \(\catz(X)\leq n\);
\item there is a simply connected rational space \(Y\) of finite type such
      that \(X\) is a rational homotopy retract of \(Y\) and
      \(\Clz(Y)\leq n\);
\item there is a simply connected rational space \(Y\) of finite type such
      that \(X\) is a rational homotopy retract of \(Y\) and
      \(\Hnilz(Y)\leq n\).
\end{enumerate}
Consequently,
\begin{equation}\label{stab:eq:parallel-retractive-closures}
 \catz(X)=\inf_{\,X\text{ retract of }Y}\Clz(Y)
 =\inf_{\,X\text{ retract of }Y}\Hnilz(Y).
\end{equation}
Here ``retract'' means rational homotopy retract and \(Y\) ranges over
simply connected rational homotopy types of finite type.  If
\(\catz(X)<\infty\), both infima are attained.
\end{theorem}

\begin{proof}
The equivalence of (1) and (2) is the rational Ganea retract theorem of
Félix, Halperin, and Thomas \cite[Theorem~28.5(iii)]{FHT2001}.  If (3)
holds, monotonicity under homotopy retracts and
Proposition~\ref{stab:prop:cat-lower} give
\[
 \catz(X)\leq\catz(Y)\leq\Hnilz(Y)\leq n.
\]
Conversely, when (1) holds, the second half of the proof of
Proposition~\ref{stab:prop:intermediate-Hnil-witness} constructs a
finite-type Sullivan algebra \(P_n\) with \(\Hnilz(P_n)\leq n\) and a
strict retraction \(P_n\to M_X\).  Sullivan realization reverses the arrows,
so \(X\) is a rational homotopy retract of \(Y=|P_n|\).  This proves (3),
the displayed equality, and attainment when the common value is finite.
\end{proof}

\section{Why retractive closure does not imply pointwise equality}

The preceding equalities might suggest the pointwise statement
\[
 \Hnilz(X)=\Clz(X).
\]
The statement is false: the example constructed in
Chapters~\ref{div:sec:witness}--\ref{chap:strict-separation} satisfies
\(\Clz(X)=3<\Hnilz(X)\).

\begin{warningbox}{Retractions do not preserve acyclicity of images}
In the optimal-root construction, the retraction
\(\rho_n:P_n\to M_X\) sends the acyclic kernel \(K_n\) to an ideal whose
image contains \((M_X^+)^{n+1}\).  It does \emph{not} follow that
\(\rho_n(K_n)\) is acyclic.  An image of an acyclic complex under a chain
map may acquire homology.  This is the precise descent gap that the
counterexample realizes.
\end{warningbox}

As anticipated in Remark~\ref{div:rem:first-threshold},
Part~\ref{part:divergence} realizes this length-three ideal-closure defect in
an explicit \(56\)-dimensional CDGA.

The calculation tools show that the flexible and rigid invariants agree through length two and identify the first possible ideal-closure defect at length three.  Part~\ref{part:towers} now organizes that defect by compatible retractive towers.

\part{Retractive Towers and the Strictification Defect}\label{part:towers}

\chapter{Relative Sullivan Factorizations and Compatible Retractions}
\label{rt:chap:foundations}

\begin{chapterguide}[title={Chapter guide}]
Part~\ref{part:calculations} showed how individual acyclic-ideal constructions
can succeed or fail.  We now replace those isolated searches by a tower that
follows every stage of a flexible short witness.  This chapter supplies the
relative factorizations, strict retractions, compatible homotopies, and
pushouts needed to build that tower; the degree-one refinements are treated at
the end.
\end{chapterguide}

Throughout, we use the CDGA, augmentation, connectedness, and grading
conventions fixed in Chapter~\ref{stab:sec:alg-prelim}.

\phantomsection\label{rt:convention:quantifiers}
\begin{warningbox}{Convention on the tower quantifiers}
The characterization of \(\nilh(M)\leq n\) asks for the existence of one
complete coherent branch.  Such a branch is built from one chosen short
witness and is not canonical.  Its image ideals depend further on the chosen
strict retractions.  Consequently, failure of one branch or one retraction
does not imply failure for all branches; only an obstruction proved for every
complete branch yields an invariant lower bound.
\end{warningbox}
\index{quantifiers!retractive tower}

\section*{Tower stage notation and arrow dictionary}

Fix a minimal Sullivan algebra \(M\) and a truncation level \(n\).  Every
stage begins with a relative Sullivan model of the augmentation truncation,
equipped with a literal retraction:
\[
\begin{tikzcd}[column sep=large,row sep=small]
 M \arrow[r,hook,"i_k"] \arrow[dr,equal]
   & E_k \arrow[r,two heads,"\pi_k"] \arrow[d,"\rho_k"]
   & C_k=\operatorname{Aug}_n(A_k)\, .\\
 & M
\end{tikzcd}
\]
To pass from stage \(k\) to stage \(k+1\), we first factor \(\rho_k\), then
form a pushout, and finally factor the full map to the next truncation:
\[
\begin{tikzcd}[column sep=large,row sep=large]
 E_k \arrow[r,hook,"j_k"] \arrow[d,two heads,"\pi_k"']
   & E'_k \arrow[r,hook,"\lambda_k"]
             \arrow[d,"\beta_k"']
   & E_{k+1} \arrow[d,two heads,"\pi_{k+1}"]\\
 C_k \arrow[r,"\iota_k"']
   & A_{k+1} \arrow[r,two heads,"q_{k+1}"']
   & C_{k+1}.
\end{tikzcd}
\]
The strict identities
\[
                  r_kj_k=\rho_k,
        \qquad \rho_{k+1}\lambda_k=r_k
\]
are the coordination data.  Hooked arrows are relative Sullivan extensions.
The maps \(\pi_k\) and \(\pi_{k+1}\) are the chosen surjective
quasi-isomorphisms, while the other two-headed arrows are augmentation
quotients.  The vertical solid arrows \(\rho_k,r_k\) are strict retractions or
their chosen factorizations, and the remaining bottom horizontal arrows are
tower transitions.

\begin{center}
\small
\renewcommand{\arraystretch}{1.18}
\begin{tabularx}{\textwidth}{@{}>{$}l<{$}>{\raggedright\arraybackslash}X@{}}
\toprule
\text{symbol} & \text{meaning}\\
\midrule
C_k=\operatorname{Aug}_n(A_k)
  & The intrinsic level-\(n\) augmentation truncation.\\
M\xrightarrow{i_k}E_k\xrightarrow{\pi_k}C_k
  & A chosen relative Sullivan factorization, with
    \(\rho_ki_k=\id_M\).\\
E_k\xrightarrow{j_k}E'_k\xrightarrow{r_k}M
  & A relative factorization of the chosen retraction; its pushout along
    \(\pi_k\) produces \(A_{k+1}\).\\
E'_k\xrightarrow{\lambda_k}E_{k+1}\xrightarrow{\pi_{k+1}}C_{k+1}
  & Strict coordination of consecutive stages, with
    \(\rho_{k+1}\lambda_k=r_k\).\\
K_k,\ S_k,\ J_k
  & Respectively \(\ker\pi_k\),
    \(K_k\cap\ker\rho_k\), and \(\rho_k(K_k)\subseteq M\).\\
I_k^\alpha,\ \widehat I_k
  & The map-kernel \(\ker\alpha_k\) and truncation-kernel
    \(\ker(q_k\alpha_k)\), with
    \(I_k^\alpha\subseteq\widehat I_k\subseteq J_k\).\\
J_\infty,\ T_\infty
  & The union \(\bigcup_kJ_k\) and the comparison defect
    \(\ker(A_\infty\twoheadrightarrow M/J_\infty)\).\\
\bottomrule
\end{tabularx}
\end{center}

Chapter~\ref{rt:chap:towers} constructs the towers and identifies their
existential meaning; Chapter~\ref{def:chap:comparison} extracts the
image-ideal defect; and Chapter~\ref{chap:tate-selection} studies which fixed
retractions can remove that defect.

\section{Relative Sullivan extensions and factorizations}
\index{relative Sullivan algebra}

We use relative Sullivan extensions in the sense of
Definition~\ref{alg:def:sullivan-algebra}.  We write one as
\(A\hookrightarrow A\otimes\Lambda Z\), call the chosen well-order on \(Z\)
the Sullivan order, and call the adjunction of a
collection of generators whose differentials lie in the algebra already
constructed a \emph{Hirsch extension}.  Thus a relative Sullivan extension
is a possibly transfinite composite of Hirsch extensions.

The adjective \emph{relative} is always understood with respect to the
displayed base.  Thus a
homotopy relative to \(A\) is constant on the image of \(A\), and a retraction
of \(i\) is a CDGA morphism \(\rho:A\otimes\Lambda Z\to A\) satisfying
\(\rho i=\id_A\).  A linear splitting of cochain complexes is not a CDGA
retraction: once the values on generators have been chosen, multiplicativity
and all Koszul signs are forced.

\begin{proposition}[Relative Sullivan factorization]
\label{rt:prop:relative-factorization}
Let \(g:A\to C\) be a morphism of connected augmented CDGAs such that
\(H^1(g)\) is injective.  Then \(g\) admits a factorization
\begin{equation}
\label{rt:eq:relative-factorization}
 A\xrightarrow{i}A\otimes\Lambda Z
   \xrightarrow[\simeq]{\pi}\!\!\!\twoheadrightarrow C,
 \qquad \pi i=g,
\end{equation}
where \(i\) is a relative Sullivan extension and \(\pi\) is a surjective
quasi-isomorphism.
\end{proposition}

\begin{proof}
First make the map surjective on the underlying graded vector spaces.  Choose
a homogeneous vector-space complement to \(g(A)\) in \(C\), together with a
well-ordered homogeneous basis.  For every basis element \(c\) of positive
degree \(m\), adjoin a disk pair
\[
 |u_c|=m,\qquad |v_c|=m+1,
 \qquad du_c=v_c,\qquad dv_c=0,
\]
place \(v_c\) before \(u_c\) in the Sullivan order, and define
\(\pi(u_c)=c\), \(\pi(v_c)=dc\).  Tensoring with a finite collection of
such disks preserves cohomology, and filtered colimits of complexes over
\(\Q\) are exact; hence the full disk extension has the cohomology of \(A\).
Its map to \(C\) is surjective as a graded map.  We regard this as correction
round zero.

Starting from that surjection, construct an ordinal sequence of correction
rounds.  The round number, not cohomological degree, is the primary Sullivan
ordering parameter.  At a successor round, first choose a homogeneous basis
of the cokernel of the current cohomology map and adjoin closed generators
mapping to chosen cocycle representatives.  These generators may be
well-ordered arbitrarily because their differentials are zero.  After this
cokernel step, choose a homogeneous basis of the kernel of the new
cohomology map.  For a basis class represented by a cocycle \(x\) of degree
\(m\), choose \(b\in C^{m-1}\) with \(db=\pi(x)\), adjoin \(y_x\) of degree
\(m-1\), and set
\[
                         dy_x=x,
                    \qquad \pi(y_x)=b.
\]
The representative \(x\) lies in the algebra present before the kernel
killers of that round are adjoined.  Thus every element occurring in
\(dy_x\) precedes \(y_x\).  At a limit round take the union.  No ordering by
cohomological degree is imposed across different rounds, so the resulting
well-order satisfies the Sullivan condition even when a late correction has
smaller degree than a generator from an earlier round.

For completeness, we record the degree-one invariant.  After round zero the
map on \(H^1\) is the original injective map.  Closed degree-one generators
added in a cokernel step map to classes independent modulo the previous
image, so injectivity is preserved.  A degree-zero kernel killer would be
needed only for a nonzero kernel class in \(H^1\), and therefore never
occurs.  Degree-one generators used to kill basis classes in \(H^2\) are not
cycles.  More generally, if a degree-one cycle involved a linear combination
of such generators, differentiating would give a relation among the selected
basis classes in the kernel of \(H^2\); all coefficients would vanish.  Hence
\(H^1\)-injectivity is preserved at every round.

Every element of the union is a finite linear combination of finite words and
therefore occurs at some earlier round.  A missing target class is represented
at the next successor round, and a cocycle whose image is a boundary is killed
at the next kernel step.  The limiting map is consequently an isomorphism on
cohomology.  It remains surjective on underlying graded vector spaces because
round zero was surjective.  The transfinite composite is a relative Sullivan
extension and proves the result.
\end{proof}

Neither the relative generators nor their Sullivan order are canonical.  In
every application below the injectivity hypothesis in degree one follows
from the vanishing of the relevant first cohomology groups.

The following refinement is essential when the truncation level is one.

\begin{lemma}[Degree-one-free factorization]
\label{rt:lem:degree-one-free-factorization}
Let \(g:A\to C\) be a morphism of connected augmented CDGAs, with \(A\) a
Sullivan algebra.  Assume
\[
             A^1=C^1=0
       \quad\text{and}\quad
             H^2(g)\text{ is injective}.
\]
Then \(g\) has a factorization of the form
\[
 A\longrightarrow A\otimes\Lambda Z
   \xrightarrow[\simeq]{\pi}\!\!\!\twoheadrightarrow C
\]
with \(Z^1=0\).
\end{lemma}

\begin{proof}
Use the correction-round ordering of
Proposition~\ref{rt:prop:relative-factorization}.  We maintain the two
induction invariants that the current source has zero degree-one component and
that its map to \(C\) is injective on \(H^2\).  Round zero introduces no
degree-one disk generator: there is no missing element of \(C^1\), and
connectedness leaves no missing degree-zero element whose disk partner could
have degree one.

At a cokernel step there is no degree-one class to represent.  A closed
degree-two generator, when needed, maps to a class independent modulo the
previous image and therefore preserves injectivity on \(H^2\).  At a kernel
step a generator of degree one would be required only to kill a class in the
kernel of \(H^2\), which the induction invariant excludes.  Killers of
higher-degree classes have degree at least two.  A degree-two cycle involving
new degree-two killers would differentiate to a relation among the selected
kernel basis classes in degree three, so its new coefficients vanish; hence
no later round creates a kernel in \(H^2\).  The invariants persist through
successor rounds and through unions at limit rounds.  Thus the factorization
has \(Z^1=0\).
\end{proof}

All relative generators below have positive degree and augmentation zero.
Consequently their differentials lie in the augmentation ideal, and every
relative factorization is a factorization in augmented CDGAs.

\section{Pushouts and Sullivan base change}

If \(i:A\to A\otimes\Lambda Z\) is a relative Sullivan extension and
\(u:A\to C\) is any morphism, their pushout is
\[
 C\otimes_A(A\otimes\Lambda Z)
       \cong C\otimes\Lambda Z,
\]
with the differential obtained by applying \(u\) to the coefficients in
\(A\).  In particular, the base-changed map remains a relative Sullivan
extension.

\begin{lemma}[Sullivan base change]
\label{rt:lem:sullivan-base-change}
In a pushout square
\[
\begin{tikzcd}[column sep=large,row sep=large]
A \arrow[r,hook,"i"] \arrow[d,"u"']
  & A\otimes\Lambda Z \arrow[d,"\bar u"]\\
C \arrow[r,hook]
  & C\otimes_A(A\otimes\Lambda Z),
\end{tikzcd}
\]
if \(i\) is a relative Sullivan extension and \(u\) is a
quasi-isomorphism, then \(\bar u\) is a quasi-isomorphism.
\end{lemma}

\begin{proof}
First suppose that a single generator \(z\) is adjoined.  Filter the source
and target by powers of \(z\) (by exterior length when \(|z|\) is odd).  On
the associated graded objects the differential of \(z\) vanishes, and the
induced map is
\[
                   u\otimes\id_{\Lambda z}.
\]
It is a quasi-isomorphism over the field \(\Q\).  The filtration is bounded
below and, in each total degree, finite; hence the filtered comparison gives
the result for one Hirsch generator.  The same argument treats a set of
generators adjoined at one stage by filtering first by total Hirsch length.

Proceed through the Sullivan order.  Successor stages follow from the
one-stage argument and the two-out-of-three property.  At a limit stage take
the filtered union.  Filtered colimits of complexes of \(\Q\)-vector spaces
are exact, so cohomology commutes with this union.  The result follows for the
entire relative Sullivan extension.
\end{proof}

\section{The polynomial cylinder}

Except inside polynomial cylinders, all structural morphisms in this memoir
are augmented.  Homotopies and homotopy-lifting arguments, however, take place
in the category of unital, nonnegatively graded CDGAs: no augmentation of
\(B[t,dt]\) makes both endpoint evaluations augmented, because they send \(t\)
to \(0\) and \(1\), respectively.  Thus the homotopy itself need not preserve
augmentations, although both endpoint maps do.  If the homotopy is relative
to an augmented sub-CDGA \(A\), its restriction to \(A\) is the
constant-cylinder inclusion of the common endpoint map.

For a CDGA \(B\), put
\[
 B[t,dt]=B\otimes\Lambda(t,dt),
 \qquad |t|=0,\quad |dt|=1,\quad d(t)=dt,\quad d(dt)=0.
\]
The endpoint evaluations \(\varepsilon_e:B[t,dt]\to B\), \(e=0,1\), send
\(t\) to \(e\) and \(dt\) to zero.  A CDGA homotopy from \(f_0:E\to B\)
to \(f_1:E\to B\) is a morphism \(h:E\to B[t,dt]\) with
\(\varepsilon_0h=f_0\) and \(\varepsilon_1h=f_1\).

Define an operator of degree \(-1\) by
\begin{align}
 K_t(bt^j)&=0,\notag\\
 K_t(bt^jdt)&=(-1)^{|b|}b\frac{t^{j+1}}{j+1}
 \label{rt:eq:integration-operator}
\end{align}
for homogeneous \(b\in B\) and \(j\geq0\).  The sign occurs because
\(dt\) is written to the right of \(b\).

\begin{lemma}[Polynomial integration]
\label{rt:lem:polynomial-integration}
Let \(\iota_0:B\to B[t,dt]\) be the constant inclusion and set
\(K_1=\varepsilon_1K_t\).  Then
\begin{align}
 dK_t+K_td&=\id-\iota_0\varepsilon_0,
 \label{rt:eq:cylinder-contraction}\\
 dK_1+K_1d&=\varepsilon_1-\varepsilon_0.
 \label{rt:eq:endpoint-contraction}
\end{align}
Consequently, a homotopy \(h:f_0\simeq f_1\) determines a degree
\(-1\) cochain homotopy \(F=K_1h\) satisfying
\[
                      f_1-f_0=dF+Fd.
\]
\end{lemma}

\begin{proof}
It suffices to evaluate the first identity on \(bt^j\) and \(bt^jdt\).
The terms involving \(db\) cancel because of the sign in
\eqref{rt:eq:integration-operator}; the remaining derivative of
\(t^{j+1}/(j+1)\) gives \(t^jdt\).  For \(j=0\), the constant term left
after integration is exactly \(\iota_0\varepsilon_0\).  This proves
\eqref{rt:eq:cylinder-contraction}.  Applying \(\varepsilon_1\) gives
\eqref{rt:eq:endpoint-contraction}, and composing with \(h\) gives the last
formula.
\end{proof}

There is an equivalent formula on a free algebra.  If \(f_0,f_1:E\to B\)
are algebra maps and \(F:E\to B\) has degree \(-1\), extend its values on
an ordered set of free generators by
\begin{equation}
\label{rt:eq:ordered-homotopy}
 F(w_1\cdots w_r)=
 \sum_{j=1}^r(-1)^{|w_1|+\cdots+|w_{j-1}|}
 f_0(w_1)\cdots f_0(w_{j-1})F(w_j)
 f_1(w_{j+1})\cdots f_1(w_r).
\end{equation}
Equivalently,
\[
 F(ab)=F(a)f_1(b)+(-1)^{|a|}f_0(a)F(b).
\]
If \(dF+Fd=f_1-f_0\) holds on generators, the displayed formula telescopes
to the same identity on all of \(E\).  This observation is useful whenever
nonlinear differentials are written explicitly.

\section{Relative lifting with a prescribed homotopy}

We now prove the lifting statement used throughout the tower construction.
It refines the Sullivan lifting lemma
(Lemma~\ref{stab:lem:sullivan-lifting}) by prescribing a homotopy on the base
of the relative extension.  The formulation includes a nontrivial base map;
the usual retraction lemma will be an immediate special case.

\begin{lemma}[Relative lifting and homotopy extension]
\label{rt:lem:relative-lifting}
Let \(\lambda:D\to E\) be a relative Sullivan extension and let
\(\phi:M\to B\) be a quasi-isomorphism.  Suppose that maps
\[
        f:E\to B,
        \qquad r:D\to M
\]
are given together with a homotopy
\[
 h_D:D\longrightarrow B[t,dt]
\]
from \(f\lambda\) to \(\phi r\).  Then there exist a map
\(\rho:E\to M\) and a homotopy \(h_E:f\simeq\phi\rho\) such that
\[
                   \rho\lambda=r,
             \qquad h_E\lambda=h_D.
\]
\end{lemma}

\begin{proof}
Work in the underlying unital CDGA category and form the path-object pullback
\[
 P(\phi)=
 \bigl\{(m,u)\in M\times B[t,dt]\mid
             \phi(m)=\varepsilon_1(u)\bigr\}
\]
and let \(p:P(\phi)\to B\) be \(p(m,u)=\varepsilon_0(u)\).
The projection \(P(\phi)\to M\) is the pullback of the surjective
quasi-isomorphism \(\varepsilon_1:B[t,dt]\to B\).  It is therefore
surjective, and its kernel is isomorphic to the acyclic complex
\(\ker\varepsilon_1\); hence it is a quasi-isomorphism.  Moreover,
\(p\) and the composite
\[
 P(\phi)\longrightarrow M\xrightarrow{\phi}B
\]
induce the same cohomology map, because
\(\varepsilon_0\) and \(\varepsilon_1\) are chain homotopic by
Lemma~\ref{rt:lem:polynomial-integration}.  Thus \(p\) is a
quasi-isomorphism.  It is surjective: for a homogeneous \(b\in B\), the
pair \((0,(1-t)b)\) belongs to \(P(\phi)\) and maps to \(b\).

The prescribed data define a map
\[
             D\longrightarrow P(\phi),
             \qquad x\longmapsto(r(x),h_D(x)),
\]
and give a commutative square
\[
\begin{tikzcd}[column sep=large,row sep=large]
D \arrow[r,"{(r,h_D)}"] \arrow[d,hook,"\lambda"']
  & P(\phi) \arrow[d,"p",two heads,"\sim"']\\
E \arrow[r,"f"'] & B.
\end{tikzcd}
\]
Although \(P(\phi)\) need not be connected, the unital version of
Lemma~\ref{stab:lem:sullivan-lifting} applies to the preceding square.  Let
\(\ell:E\to P(\phi)\) be the resulting lift and put
\[
 \rho=\operatorname{pr}_M\ell,
 \qquad h_E=\operatorname{pr}_{B[t,dt]}\ell.
\]
Then \(\ell\lambda=(r,h_D)\), while \(p\ell=f\); the pullback identity gives
\(\varepsilon_1h_E=\phi\rho\).  Hence
\[
 \rho\lambda=r,\qquad h_E\lambda=h_D,\qquad
 \varepsilon_0h_E=f,\qquad \varepsilon_1h_E=\phi\rho.
\]
Although the path object and the lift were constructed in the unital
category, the map
\(\rho:E\to M\) is augmented: both \(E\) and \(M\) are connected, so their
augmentations are the unique unital degree-preserving maps to \(\Q\).
The endpoint maps of \(h_E\) are \(f\) and \(\phi\rho\), and hence are
augmented as asserted.
\end{proof}

\begin{corollary}[Relative retraction]
\label{rt:cor:relative-retraction}
Let \(i:M\to E\) be a relative Sullivan extension, let
\(\phi:M\to B\) be a quasi-isomorphism, and let \(f:E\to B\) satisfy
\(fi=\phi\).  Then there exist a CDGA morphism \(\sigma:E\to M\) and a
homotopy relative to \(M\) such that
\[
               \sigma i=\id_M,
          \qquad f\simeq\phi\sigma\quad\text{relative to }M.
\]
In particular, \(\sigma\) is surjective and
\[
                  E=i(M)\oplus\ker\sigma
\]
as cochain complexes.
\end{corollary}

\begin{proof}
Apply Lemma~\ref{rt:lem:relative-lifting} with \(D=M\), \(\lambda=i\),
\(r=\id_M\), and with the constant homotopy on \(\phi\).  The last
decomposition follows because \(i\sigma\) is an idempotent cochain map.
\end{proof}

This proof also displays the obstruction equation.  A retraction defined on
the subalgebra preceding a relative generator \(z\) extends over \(z\) exactly
when one can solve
\begin{equation}
\label{rt:eq:retraction-equation}
                    d\rho(z)=\rho(Dz)
                    \quad\text{in }M.
\end{equation}
The right-hand side is a cocycle.  Its cohomology class is the obstruction,
and, when it vanishes, the possible values of \(\rho(z)\) form a torsor
under \(Z^{|z|}(M)\).  Thus choices of retraction naturally branch along the
Sullivan order.

\section{Simultaneous extension}

The next result extends both maps in a given compatible retractive stage to
quasi-isomorphisms while preserving the homotopy between them.

\begin{proposition}[Simultaneous Hirsch extension]
\label{rt:prop:simultaneous-extension}
Let \(M\) be a simply connected Sullivan algebra, let
\(\phi:M\qiso B\), and let \(i:M\to E\) be a relative Sullivan
extension with \(H^1(E)=0\).  Suppose \(f:E\to B\) satisfies
\(fi=\phi\).  Choose a retraction and relative homotopy
\[
 \sigma:E\to M,
 \qquad \sigma i=\id_M,
 \qquad h:f\simeq\phi\sigma
\]
as in Corollary~\ref{rt:cor:relative-retraction}.  There exist a relative
Sullivan extension
\[
                    j:E\longrightarrow
                    \widetilde E=E\otimes\Lambda Y,
\]
extensions \(\widetilde f:\widetilde E\to B\) and
\(\widetilde\sigma:\widetilde E\to M\), and a homotopy relative to
\(M\), such that
\[
 \widetilde f\simeq\phi\widetilde\sigma,
 \qquad \widetilde f j=f,
 \qquad \widetilde\sigma j=\sigma,
 \qquad \widetilde\sigma i=\id_M,
\]
and both \(\widetilde f\) and \(\widetilde\sigma\) are
quasi-isomorphisms.  In particular,
\(\widetilde\sigma\) is a surjective quasi-isomorphism.
\end{proposition}

\begin{proof}
Orient \(h\) from \(f\) to \(\phi\sigma\), put \(F=K_1h\), and use
Lemma~\ref{rt:lem:polynomial-integration} to obtain
\begin{equation}
\label{rt:eq:chain-homotopy-F}
                    \phi\sigma-f=dF+Fd.
\end{equation}
Because \(\sigma i=\id_M\), the map \(H(\sigma)\) is surjective.  Since
\(H(f)=H(\phi)H(\sigma)\) and \(H(\phi)\) is an isomorphism,
\[
       \ker H(f)=\ker H(\sigma)=H(\ker\sigma).
\]
The last equality follows from the cochain splitting in
Corollary~\ref{rt:cor:relative-retraction}.

Choose homogeneous cocycles \(x_\lambda\in\ker\sigma\) whose classes form
a well-ordered homogeneous basis of \(\ker H(f)\), and set
\(m_\lambda=|x_\lambda|\).  The
hypothesis \(H^1(E)=0\) allows the representatives to be chosen with
\(m_\lambda\geq2\).  Adjoin generators \(z_\lambda\) of degree
\(m_\lambda-1\) and prescribe
\begin{equation}
\label{rt:eq:kernel-killing-extension}
 dz_\lambda=x_\lambda,
 \qquad \widetilde\sigma(z_\lambda)=0,
 \qquad \widetilde f(z_\lambda)=-F(x_\lambda).
\end{equation}
For a cocycle \(x_\lambda\in\ker\sigma\),
\eqref{rt:eq:chain-homotopy-F} gives
\(dF(x_\lambda)=-f(x_\lambda)\).  Hence the last prescription in
\eqref{rt:eq:kernel-killing-extension} commutes with the differential.
Extend the homotopy over each new generator by
\begin{equation}
\label{rt:eq:extended-cylinder-homotopy}
 \widetilde h(z_\lambda)=
       \widetilde f(z_\lambda)+K_t h(x_\lambda).
\end{equation}
Indeed, \eqref{rt:eq:cylinder-contraction} gives
\(d\widetilde h(z_\lambda)=h(x_\lambda)\); its value at zero is
\(\widetilde f(z_\lambda)\), while its value at one is zero, equal to
\(\phi\widetilde\sigma(z_\lambda)\).

This Hirsch extension kills the current kernel of \(H(f)\), but may create
new kernel classes.  Repeat the operation by correction rounds, with the round
number primary exactly as in
Proposition~\ref{rt:prop:relative-factorization}.  At each successor round the
representatives \(x_\lambda\) belong to the algebra constructed in the
preceding round, and each \(z_\lambda\) is ordered after every element
occurring in \(dz_\lambda=x_\lambda\).  At a limit round take the union.
Thus the total extension is Sullivan; no degree-first ordering across rounds
is used.

The invariant \(H^1(E_{[s]})=0\) is preserved.  A degree-one cycle involving
newly adjoined degree-one generators would give, after differentiation, a
linear relation among the independent kernel classes
\([x_\lambda]\in H^2(E_{[s]})\); all new coefficients therefore vanish, and
the remaining old cycle is already a boundary.  Consequently no degree-zero
generator is ever needed.  The construction is set-sized for the same reason
as in Proposition~\ref{rt:prop:relative-factorization}.

Let \(E_{[s]}\) be the algebra after \(s\) repetitions and put
\(\widetilde E=\varinjlim_sE_{[s]}\).  Every cocycle in
\(\widetilde E\) uses only finitely many generators and hence occurs at an
earlier stage.  If its image in \(B\) is a boundary, its class belongs to
\(\ker H(E_{[s]}\to B)\) at some stage and is killed at the next.  Hence
\(H(\widetilde f)\) is injective.  The map \(H(\widetilde f)\) is surjective
because \(\widetilde f\) restricts on \(M\) to the quasi-isomorphism \(\phi\).
Therefore \(\widetilde f\) is a
quasi-isomorphism.  The limiting homotopy gives
\(H(\widetilde f)=H(\phi)H(\widetilde\sigma)\); consequently
\(\widetilde\sigma\) is a quasi-isomorphism as well.  Its restriction to
\(M\) is the identity, so it is surjective.
\end{proof}

\begin{corollary}[Degree-one-free simultaneous extension]
\label{rt:cor:degree-one-free-extension}
In Proposition~\ref{rt:prop:simultaneous-extension}, assume moreover that
\(E^1=0\) and that \(H^2(f)\) is injective.  The relative extension can be
chosen with \(Y^1=0\); then \(\widetilde E^1=0\), and
\(H^2(\widetilde f)\) remains injective.
\end{corollary}

\begin{proof}
At every step, a generator of degree \(m-1\) is adjoined to kill a kernel
class of degree \(m\).  Initially there are no kernel classes in degrees one
or two.  Suppose degree-two generators have been attached to kill independent
kernel classes in degree three.  A degree-two cocycle of the form
\(e+\sum c_\lambda z_\lambda\) would satisfy
\[
                    de+\sum c_\lambda x_\lambda=0.
\]
Passing to degree-three cohomology forces every \(c_\lambda\) to vanish.
No degree-two boundary is added because no degree-one generator is present.
Thus \(H^2\) is unchanged and injectivity of the map to \(B\) persists.  This
is the induction invariant at every correction round; it is preserved under
unions because every degree-two cocycle involves only finitely many
generators.  Hence no degree-one generator is ever required.  Taking the
filtered union gives \(Y^1=0\) and the stated conclusions.
\end{proof}

\section{Strict comparison through an intermediate object}

The general lifting lemma has the following consequence, which will create
literal maps between consecutive tower stages.

\begin{lemma}[Comparison under an intermediate algebra]
\label{rt:lem:intermediate-comparison}
Let \(\lambda:E'\to \widehat E\) be a relative Sullivan extension, let
\(\phi:M\qiso B\), and let \(\widehat f:\widehat E\to B\).  Suppose a map
\(r:E'\to M\) is equipped with a homotopy from \(\widehat f\lambda\) to
\(\phi r\).  Then there exist a map \(\widehat\rho:\widehat E\to M\) and an extension
of that homotopy such that
\[
                         \widehat\rho\lambda=r.
\]
If a fixed copy of \(M\) lies in \(E'\), if \(r\) restricts to the
identity there, and if the given homotopy is constant there, then
\(\widehat\rho\) is a retraction and the extended homotopy is relative to
\(M\).
\end{lemma}

\begin{proof}
This is Lemma~\ref{rt:lem:relative-lifting} with \(D=E'\) and \(E=\widehat E\).  In
the proof of that lemma, the lift \(\ell:\widehat E\to P(\phi)\) extends the
prescribed map \(E'\to P(\phi)\) literally.  Hence both the equality on \(E'\)
and the relative condition on \(M\) are preserved.
\end{proof}

The adjective \emph{strict} is important.  If \(E'\to C'\) is factored as
\[
 E'\xrightarrow{\lambda}\widehat E
       \xrightarrow[\simeq]{\widehat\pi}\!\!\!\twoheadrightarrow C',
\]
the relative map \(\lambda\) is a literal inclusion, and the lemma permits
the new retraction to restrict literally to the old map \(r\).  Choosing an
unrelated factorization of a map defined only on \(M\) would provide no such
comparison.  The next chapter uses precisely these intermediate inclusions
to define a complete strictly coordinated tower.

\chapter{Complete Retractive Towers and Homotopical Nil-Length}
\label{rt:chap:towers}

\begin{chapterguide}[title={Chapter guide}]
This chapter turns a short CDGA model into an infinite sequence of strict
retraction problems, and recovers a short model from one complete sequence by
a filtered colimit.  The logical quantifier is essential: the criterion asks
for the existence of one coherent complete branch.  A prescribed branch may
fail even though another succeeds.  We also prove that a successful branch
may be chosen strictly coordinated, so that consecutive relative Sullivan
models are connected by literal inclusions.  The special degree-one
bookkeeping at level one is included in full.
\end{chapterguide}

The retractive-tower viewpoint originates in an observation of Cornea
\cite[Remark~(a) following Proposition~3.6, p.~108]{Cornea1994}.
There, successive homotopy-retraction requirements for augmentation-ideal
quotients are proposed as a sequence of obstructions to nil-length at
most~\(n\).  Key discussions with Daniel Tanr\'e helped shape this framework.
We make Cornea's observation precise by organizing the successive
relative Sullivan factorizations into complete---and, when needed, strictly
coordinated---retractive towers and proving the resulting characterization of
homotopical nil-length.

Let
\[
                       M=(\Lambda V,d)
\]
be a simply connected minimal Sullivan algebra of finite type.  All CDGAs and
maps below are connected, cohomologically nonnegative, and augmented over
\(\Q\).

\section{Short models and one-arrow witnesses}

The notions of an \(n\)-short augmented CDGA and of \(\nilh\) were fixed in
Definitions~\ref{def:book-short} and~\ref{div:def:nilh}.  Thus
\((B^+)^{n+1}=0\) is the shortness condition, while the definition of
\(\nilh(M)\) initially permits a zigzag of augmented quasi-isomorphisms.
For a Sullivan source, Lemma~\ref{div:lem:direct-zigzag} replaces such a
zigzag with a single arrow.  We call an augmented quasi-isomorphism from
\(M\) to an \(n\)-short CDGA an \emph{\(n\)-short witness}.  One fixed
witness will coordinate every stage of the tower.

Consequently,
\[
 \nilh(M)\leq n
 \quad\Longleftrightarrow\quad
 \text{an \(n\)-short witness \(M\qiso B\) exists}.
\]

At level \(1\), one must control cochain degree one, not merely its cohomology.

\begin{lemma}[Reduced one-short witness]
\label{rt:lem:reduced-one-short}
Suppose \(M\) is simply connected and \(\phi:M\qiso B\) is a one-short
witness.  It may be replaced by a one-short witness
\(\phi_{\mathrm{red}}:M\qiso B_{\mathrm{red}}\) satisfying
\[
                         B_{\mathrm{red}}^1=0.
\]
\end{lemma}

\begin{proof}
Since \((B^+)^2=0\), every subcomplex of \(B^+\) is a
differential ideal.  Split the cochain complex over \(\Q\) as
\[
                 B^+=H\oplus D,
\]
where \(H\) consists of chosen cocycle representatives for \(H(B^+)\),
has zero differential, and \(D\) is contractible.  The quotient
\[
                 B_{\mathrm{red}}=B/D=\Q\oplus H
\]
has square-zero augmentation ideal, zero differential, and the projection
\(B\to B_{\mathrm{red}}\) is a quasi-isomorphism.  Because
\(H^1(B)\cong H^1(M)=0\), the degree-one component of \(H\) vanishes; hence
\(B_{\mathrm{red}}^1=0\).  Compose \(\phi\) with the projection.
\end{proof}

For any augmented CDGA \(A\) and \(n\geq1\), set
\begin{equation}
\label{rt:eq:augmentation-truncation}
 \operatorname{Aug}_n(A)=A/(A^+)^{n+1},
 \qquad q_A^n:A\twoheadrightarrow\operatorname{Aug}_n(A).
\end{equation}
This is an intrinsic quotient, not a choice of model.

\section{A retractive stage and its successor}

\begin{definition}[Retractive stage of level \(n\)]
\label{rt:def:stage}
Suppose \(A_k\) is an augmented CDGA and
\(\alpha_k:M\qiso A_k\).  Put
\[
 C_k=\operatorname{Aug}_n(A_k),
 \qquad q_k=q_{A_k}^n.
\]
A \emph{stage of level \(n\)} based at \((A_k,\alpha_k)\) is a chosen
relative Sullivan factorization
\begin{equation}
\label{rt:eq:stage-factorization}
 M\xrightarrow{i_k}E_k=M\otimes\Lambda Z_k
   \xrightarrow[\simeq]{\pi_k}\!\!\!\twoheadrightarrow C_k,
 \qquad \pi_ki_k=q_k\alpha_k.
\end{equation}
It is \emph{retractive} if there is a CDGA morphism
\[
                 \rho_k:E_k\longrightarrow M,
                 \qquad \rho_ki_k=\id_M.
\]
\end{definition}

The factorization and the retraction are part of the data.  In particular,
different choices of relative generators can create different branches.
Nevertheless, every retractive stage satisfies the necessary condition
\[
 H(q_k\alpha_k)=H(\pi_k)H(i_k)
             \quad\Longrightarrow\quad
 H(q_k\alpha_k)\text{ is injective},
\]
because \(H(i_k)\) is split injective.

We now define the successor of a retractive stage.  Choose, as part of the
successor data, a factorization of its retraction
\begin{equation}
\label{rt:eq:retraction-factorization}
 E_k\xrightarrow{j_k}E'_k=E_k\otimes\Lambda Y_k
   \xrightarrow[\simeq]{r_k}\!\!\!\twoheadrightarrow M,
 \qquad r_kj_k=\rho_k,
\end{equation}
where \(j_k\) is a relative Sullivan extension.  Form the pushout
\begin{equation}
\label{rt:eq:successor-pushout}
\begin{tikzcd}[column sep=large,row sep=large]
E_k \arrow[r,hook,"j_k"]
    \arrow[d,"\pi_k"',two heads,"\sim"]
  & E'_k \arrow[d,"\beta_k","\sim"']\\
C_k \arrow[r,"\iota_k"']
  & A_{k+1}=C_k\otimes_{E_k}E'_k.
\end{tikzcd}
\end{equation}
In the towers constructed below, \(H^1(E_k)=H^1(M)=0\), so the existence of
this positive-degree factorization follows from
Proposition~\ref{rt:prop:relative-factorization}.  Every candidate branch
includes such a factorization as part of its data.
Lemma~\ref{rt:lem:sullivan-base-change} makes \(\beta_k\) a
quasi-isomorphism.  Define
\begin{equation}
\label{rt:eq:successor-alpha}
 \alpha_{k+1}=\beta_kj_ki_k:M\longrightarrow A_{k+1},
 \qquad
 t_k=\iota_kq_k:A_k\longrightarrow A_{k+1}.
\end{equation}
Since \(r_kj_ki_k=\id_M\) and \(r_k\) is a quasi-isomorphism,
\(j_ki_k\) is a quasi-isomorphism; hence so is \(\alpha_{k+1}\).  The
pushout identity gives
\begin{equation}
\label{rt:eq:direct-system-compatibility}
                      \alpha_{k+1}=t_k\alpha_k.
\end{equation}
Most importantly, \(t_k\) factors through \(C_k\), so it kills every
product of \(n+1\) elements of \(A_k^+\).

\section{Complete and strictly coordinated towers}

\begin{definition}[Complete level-\(n\) retractive tower]
\label{rt:def:complete-tower}
\index{retractive tower!complete}
A \emph{complete level-\(n\) retractive tower} of \(M\) consists of
retractive stages
\[
 (A_k,\alpha_k;E_k,i_k,\pi_k,\rho_k),
 \qquad k\geq0,
\]
with \(A_0=M\) and \(\alpha_0=\id_M\), together with factorizations
\eqref{rt:eq:retraction-factorization} and pushouts
\eqref{rt:eq:successor-pushout}, such that every successor is given by
\eqref{rt:eq:successor-alpha}.  Equivalently, such a tower is a single infinite
coherent branch through all possible stagewise choices.

The synonym \emph{complete stratified retractive tower} emphasizes that the
index \(k\) records successive strata.  It imposes no extra condition.  The
word ``level'' always refers to the fixed truncation exponent \(n\), not to
the stage number \(k\).
\end{definition}

To compare consecutive strata, we add a coordination condition.
After constructing \(A_{k+1}\), put
\(C_{k+1}=\operatorname{Aug}_n(A_{k+1})\).  There is a full map
\[
 g_k=q_{k+1}\beta_k:E'_k\longrightarrow C_{k+1}.
\]

\begin{definition}[Strictly coordinated tower]
\label{rt:def:coordinated-tower}
\index{retractive tower!strictly coordinated}
A complete retractive tower is \emph{strictly coordinated} if the next stage is
obtained by factoring the full map \(g_k\):
\begin{equation}
\label{rt:eq:strict-coordination}
 E'_k\xrightarrow{\lambda_k}E_{k+1}
   \xrightarrow[\simeq]{\pi_{k+1}}\!\!\!\twoheadrightarrow C_{k+1},
 \qquad \pi_{k+1}\lambda_k=q_{k+1}\beta_k,
\end{equation}
where \(\lambda_k\) is a relative Sullivan extension, and
\[
 i_{k+1}=\lambda_kj_ki_k,
 \qquad
 \rho_{k+1}\lambda_k=r_k.
\]
Thus \(j_k\) and \(\lambda_k\) are literal relative Sullivan inclusions,
not merely maps in the homotopy category.
\end{definition}

The equality defining \(i_{k+1}\) is compatible with the stage map:
\[
 \pi_{k+1}i_{k+1}
 =q_{k+1}\beta_kj_ki_k
 =q_{k+1}\alpha_{k+1}.
\]
Strict coordination is additional structure on a chosen branch.  Independently
selected factorizations at each stage need not admit the displayed literal
inclusions.

\section{A complete tower produces a short model}

\begin{proposition}
\label{rt:prop:tower-to-short}
If \(M\) admits a complete level-\(n\) retractive tower, then
\(\nilh(M)\leq n\).
\end{proposition}

\begin{proof}
Use the transition maps in \eqref{rt:eq:successor-alpha} to form
\[
                     A_\infty=\varinjlim_k A_k,
       \qquad \alpha_\infty:M\longrightarrow A_\infty.
\]
The maps \(\alpha_k\) are compatible by
\eqref{rt:eq:direct-system-compatibility}.  Filtered colimits of complexes
of \(\Q\)-vector spaces are exact, so
\[
                    H(A_\infty)
                  \cong\varinjlim_k H(A_k).
\]
Each \(H(\alpha_k)\) is an isomorphism, and the compatibility relation
identifies this directed system with the constant system \(H(M)\).
Therefore \(\alpha_\infty\) is a quasi-isomorphism.

Let \(x_0,\ldots,x_n\in A_\infty^+\).  Any finite collection of elements of a
filtered colimit admits representatives at a common stage; choose
\(\bar x_0,\ldots,\bar x_n\in A_k^+\).  Since \(t_k\) factors through
\[
                   C_k=A_k/(A_k^+)^{n+1},
\]
we have
\[
             t_k(\bar x_0\cdots\bar x_n)=0.
\]
Thus \(x_0\cdots x_n=0\) in the colimit, and
\[
                         (A_\infty^+)^{n+1}=0.
\]
Hence \(\alpha_\infty:M\qiso A_\infty\) is an \(n\)-short witness.
Notice that no transition map is required to be injective.
\end{proof}

\section{A short witness produces a strictly coordinated tower}

\begin{proposition}
\label{rt:prop:short-to-strong-tower}
Let \(n\geq1\).  If \(\nilh(M)\leq n\), then \(M\) admits a complete
level-\(n\) strictly coordinated retractive tower.
\end{proposition}

\begin{proof}
Choose once and for all an \(n\)-short witness
\begin{equation}
\label{rt:eq:fixed-short-witness}
                   \phi:M\qiso B,
                   \qquad (B^+)^{n+1}=0.
\end{equation}
If \(n=1\), use Lemma~\ref{rt:lem:reduced-one-short} to arrange additionally
that
\begin{equation}
\label{rt:eq:level-one-B}
                         B^1=0.
\end{equation}
We construct the tower together with quasi-isomorphisms
\begin{equation}
\label{rt:eq:coordinating-maps}
 \theta_k:A_k\qiso B,
 \qquad \theta_k\alpha_k=\phi.
\end{equation}
For \(n=1\), we maintain the stronger cochain condition
\begin{equation}
\label{rt:eq:level-one-A}
                             A_k^1=0.
\end{equation}

Start with \(A_0=M\), \(\alpha_0=\id_M\), and \(\theta_0=\phi\).
Because \(B\) is \(n\)-short, every \(\theta_k\) annihilates
\((A_k^+)^{n+1}\); hence it factors uniquely as
\begin{equation}
\label{rt:eq:theta-factorization}
             A_k\xrightarrow{q_k}C_k
                 \xrightarrow{s_k}B,
             \qquad s_kq_k=\theta_k.
\end{equation}

We first construct the initial retractive stage.  When \(n\geq2\), choose
any relative Sullivan factorization
\[
 M\xrightarrow{i_0}E_0
   \xrightarrow[\simeq]{\pi_0}\!\!\!\twoheadrightarrow C_0
\]
of \(q_0\alpha_0\).  The ideal \((A_0^+)^{n+1}\) has no component in
degrees at most two.  Therefore \(H^1(C_0)=H^1(A_0)=0\), and
\(H^1(E_0)=0\).

Suppose \(n=1\).  Condition \eqref{rt:eq:level-one-A} holds at \(k=0\)
because \(M^1=0\).  Every positive-degree element of \(A_0\) then has
degree at least two, so \((A_0^+)^2\) has no component in degrees at most
three.  Thus \(C_0^1=0\), and \(q_0\alpha_0\) induces an isomorphism on
\(H^2\).  Lemma~\ref{rt:lem:degree-one-free-factorization} supplies a
factorization with \(E_0^1=0\).

In either case, define
\[
                         f_0=s_0\pi_0:E_0\longrightarrow B.
\]
Then \(f_0i_0=\phi\).  Corollary~\ref{rt:cor:relative-retraction} gives a
retraction \(\rho_0:E_0\to M\) and a homotopy
\(f_0\simeq\phi\rho_0\) relative to \(M\).  This is the initial stage.

Assume now that a coordinated retractive stage has been constructed.  Apply
Proposition~\ref{rt:prop:simultaneous-extension} to enlarge it to
\[
 E_k\xrightarrow{j_k}E'_k,
 \qquad r_k:E'_k\qiso M,
 \qquad f'_k:E'_k\qiso B,
\]
where
\begin{equation}
\label{rt:eq:compatible-intermediate}
 r_kj_k=\rho_k,
 \qquad f'_kj_k=f_k,
 \qquad f'_k\simeq\phi r_k
       \quad\text{relative to }M.
\end{equation}
At level one, the identities
\[
 H^2(\pi_k)H^2(i_k)=H^2(q_k\alpha_k),
 \qquad
 H^2(f_k)H^2(i_k)=H^2(\phi)
\]
show that \(H^2(i_k)\) and \(H^2(f_k)\) are isomorphisms.  Indeed,
\(q_k\) is an isomorphism in degrees at most three, and \(\pi_k\) is a
quasi-isomorphism.  We may therefore use the degree-one-free
Corollary~\ref{rt:cor:degree-one-free-extension}; it gives
\begin{equation}
\label{rt:eq:level-one-Eprime}
                             (E'_k)^1=0.
\end{equation}

Form the pushout \eqref{rt:eq:successor-pushout}.  The maps \(s_k:C_k\to B\)
and \(f'_k:E'_k\to B\) agree on \(E_k\), so they induce
\[
                    \theta_{k+1}:A_{k+1}\longrightarrow B
\]
with
\begin{equation}
\label{rt:eq:theta-on-pushout}
 \theta_{k+1}\iota_k=s_k,
 \qquad \theta_{k+1}\beta_k=f'_k.
\end{equation}
Both \(\beta_k\) and \(f'_k\) are quasi-isomorphisms, so
\(\theta_{k+1}\) is a quasi-isomorphism.  Moreover,
\[
 \theta_{k+1}\alpha_{k+1}
 =f'_kj_ki_k=f_ki_k=\phi.
\]
At level one the underlying graded algebra of the pushout is
\[
 A_{k+1}\cong C_k\otimes\Lambda Y_k.
\]
Equations \eqref{rt:eq:level-one-Eprime} and
\(E'_k=E_k\otimes\Lambda Y_k\), together with \(E_k^1=C_k^1=0\), imply
\(Y_k^1=0\) and hence \(A_{k+1}^1=0\).  This completes the induction
establishing \eqref{rt:eq:level-one-A}.

It remains to choose the next stage in the strictly coordinated form.  Since \(B\)
is \(n\)-short, \(\theta_{k+1}\) factors through
\[
 q_{k+1}:A_{k+1}\longrightarrow C_{k+1};
 \qquad \theta_{k+1}=s_{k+1}q_{k+1}.
\]
Factor the entire map
\begin{equation}
\label{rt:eq:full-map-to-factor}
 q_{k+1}\beta_k:E'_k\longrightarrow C_{k+1}
\end{equation}
as in \eqref{rt:eq:strict-coordination}.  If \(n\geq2\), an arbitrary relative
Sullivan factorization works.  Its target \(E_{k+1}\) has
\(H^1(E_{k+1})=0\), because \((A_{k+1}^+)^{n+1}\) starts in degree at
least three, so \(H^1(C_{k+1})=H^1(A_{k+1})=0\).

If \(n=1\), both the source and target of
\eqref{rt:eq:full-map-to-factor} have zero degree-one component.  Moreover, by
\eqref{rt:eq:level-one-A} and the preceding degree argument, \(q_{k+1}\) is
an isomorphism in degrees at most three, while \(\beta_k\) is a quasi-isomorphism.
Thus \(H^2(q_{k+1}\beta_k)\) is an isomorphism.  Use
Lemma~\ref{rt:lem:degree-one-free-factorization} to obtain
\(E_{k+1}^1=0\).

Set \(i_{k+1}=\lambda_kj_ki_k\) and
\(f_{k+1}=s_{k+1}\pi_{k+1}\).  By
\eqref{rt:eq:theta-on-pushout},
\[
 f_{k+1}\lambda_k
 =s_{k+1}q_{k+1}\beta_k
 =\theta_{k+1}\beta_k
 =f'_k.
\]
The homotopy in \eqref{rt:eq:compatible-intermediate} is therefore precisely
the input required by Lemma~\ref{rt:lem:intermediate-comparison}.  It produces
a retraction
\(\rho_{k+1}:E_{k+1}\to M\), a homotopy
\(f_{k+1}\simeq\phi\rho_{k+1}\) relative to \(M\), and the strict
compatibility
\[
                         \rho_{k+1}\lambda_k=r_k.
\]
This is the next strictly coordinated stage.  Repeating the construction
for all \(k\) produces a complete strictly coordinated tower.
\end{proof}

\section{The characterization and its quantifiers}

\begin{theorem}[Retractive-tower characterization]
\label{rt:thm:characterization}
Let \(M=(\Lambda V,d)\) be a simply connected minimal Sullivan algebra of
finite type, and let \(n\geq1\).  The following are equivalent:
\begin{enumerate}[label=\textup{(\arabic*)},leftmargin=2.4em]
\item \(\nilh(M)\leq n\);
\item there exists a complete level-\(n\) retractive tower of \(M\);
\item there exists a complete strictly coordinated level-\(n\) retractive tower of
      \(M\).
\end{enumerate}
\end{theorem}

\begin{proof}
Proposition~\ref{rt:prop:tower-to-short} proves \((2)\Rightarrow(1)\), and
Proposition~\ref{rt:prop:short-to-strong-tower} proves
\((1)\Rightarrow(3)\).  Forgetting the comparison inclusions of a complete
strictly coordinated tower leaves a complete retractive tower; hence
\((3)\Rightarrow(2)\).
\end{proof}

The characterization is therefore genuinely global: it requires a complete
branch, whose filtered colimit supplies an \(n\)-short model.

The special normalization at \(n=1\) cannot be replaced by the weaker
statement \(H^1(A_k)=0\).  Quotienting an arbitrary CDGA with vanishing
first cohomology by \((A_k^+)^2\) can create degree-one cohomology.  The
proof instead maintains the literal condition \(A_k^1=0\), using the reduced
square-zero witness and degree-one-free relative factorizations.

\paragraph{Scope.}
The quantifier convention on page~\pageref{rt:convention:quantifiers} applies.
The new conclusion is that the colimit \(A_\infty\) is an \(n\)-short algebra
quasi-isomorphic to \(M\); it need not be a quotient of the fixed minimal
algebra.  Strict coordination only prepares the later kernel comparison, and
acyclicity of the image-ideal union is not part of
Theorem~\ref{rt:thm:characterization}.

\chapter{The Comparison Defect of a Retractive Tower}
\label{def:chap:comparison}
\begin{chapterguide}[title={Chapter guide}]
Chapter~\ref{rt:chap:towers} produced a flexible short colimit from one
complete tower.  We now ask whether strict coordination also produces the
rigid quotient required by \(\Hnil\).  The obstruction is the cohomology of
the image-ideal union \(J_\infty\), equivalently the failure of the comparison
quotient \(r_\infty\) to be a quasi-isomorphism.
\end{chapterguide}
\index{comparison defect}
\index{image ideal}
\index{retractive tower!strictly coordinated}
\index{spectral sequence!augmentation filtration}

\section{Strictly coordinated stage data}
\label{def:sec:coordinated-data}

Let \(M=(\Lambda V,d)\) be a simply connected minimal Sullivan algebra of
finite type, fix \(n\geq1\), and choose a complete strictly coordinated
level-\(n\) tower in the sense of
Definitions~\ref{rt:def:complete-tower} and
\ref{rt:def:coordinated-tower}.  The roadmap for Part~\ref{part:towers}
displays all stage and successor maps.  Here we record only the strict
identities needed for the comparison defect:
\begin{equation}
\label{def:eq:coordination-identities}
 \rho_{k+1}\lambda_k=r_k,
 \qquad r_kj_k=\rho_k.
\end{equation}

Put
\[
                 K_k:=\ker\pi_k,
 \qquad
 u_k:=(\lambda_kj_k)|_{K_k}:K_k\longrightarrow E_{k+1}.
\]
The compatibility of the tower with the next quotient gives
\(\pi_{k+1}u_k=0\), and hence \(u_k(K_k)\subseteq K_{k+1}\).  Moreover,
\begin{equation}
\label{def:eq:kernel-transition}
                       \rho_{k+1}u_k=\rho_k\quad\text{on }K_k.
\end{equation}
These identities will make all subsequent colimits actual colimits of
complexes, rather than comparisons defined only in a homotopy category.

\section{The local kernel--image butterfly}
\label{def:sec:local-butterfly}

At a fixed stage define
\begin{equation}
\label{def:eq:KSJ}
 S_k=K_k\cap\ker\rho_k,
 \qquad
 J_k=\rho_k(K_k).
\end{equation}
The complexes \(K_k\) and \(S_k\) are subcomplexes of \(E_k\), while
\(J_k\) is a differential ideal of \(M\).  Indeed, if \(a\in M\), \(x\in K_k\), and
\(y=\rho_k(x)\), then
\[
 ay=\rho_k(i_k(a)x)\in J_k,
 \qquad
 dy=\rho_k(dx)\in J_k.
\]
By definition, the restriction of \(\rho_k\) to \(K_k\) is onto \(J_k\),
and its kernel is \(S_k\).  Thus there is a short exact sequence
\begin{equation}
\label{def:eq:butterfly-sequence}
 0\longrightarrow S_k\longrightarrow K_k
   \xrightarrow{\rho_k}J_k\longrightarrow0.
\end{equation}
Since \(\pi_k\) is a surjective quasi-isomorphism, \(K_k\) is acyclic.
Consequently, the connecting morphism is an isomorphism
\begin{equation}
\label{def:eq:connecting}
 \delta_k:H^q(J_k)\xrightarrow{\cong}H^{q+1}(S_k).
\end{equation}
With the sign convention used here, if \(j\in Z^q(J_k)\) and
\(c\in K_k^q\) satisfy \(\rho_k(c)=j\), then
\[
                         \delta_k([j])=[dc].
\]
Its inverse is the intrinsic transgression
\begin{equation}
\label{def:eq:intrinsic-transgression}
 \tau^{J_k}_{\rho_k}:=\delta_k^{-1}:
 H^{q+1}(S_k)\xrightarrow{\cong}H^q(J_k).
\end{equation}
Equivalently, if \(s\in Z^{q+1}(S_k)\) and \(dc=s\) in \(K_k\), then
\[
              \tau^{J_k}_{\rho_k}([s])=[\rho_k(c)]_{J_k}.
\]
The class is independent of both choices because \(K_k\) is acyclic.

There is a second map with a different target and a different role.  The
\emph{ambient Tate obstruction} is
\begin{equation}
\label{def:eq:ambient-obstruction}
 \omega_{\rho_k}:=
 H(J_k\hookrightarrow M)\circ\tau^{J_k}_{\rho_k}:
 H^{q+1}(S_k)\longrightarrow H^q(M).
\end{equation}
The intrinsic transgression records a class inside the image ideal.  The
ambient obstruction decides whether that class has a primitive somewhere
in \(M\).  These questions must not be conflated.

\begin{definition}
\label{def:def:ambient-essential}
For \(\xi\in H^{q+1}(S_k)\), call \(\xi\) intrinsically trivial if
\(\tau^{J_k}_{\rho_k}(\xi)=0\).  It is \emph{ambient-inessential} if its
intrinsic transgression is nonzero but \(\omega_{\rho_k}(\xi)=0\), and
\emph{ambient-essential} if \(\omega_{\rho_k}(\xi)\neq0\).  The stage is
\emph{Tate-reduced} when \(H(J_k\hookrightarrow M)\) is injective.
\index{Tate-reduced}
\end{definition}

Equation~\eqref{def:eq:kernel-transition} implies
\(u_k(S_k)\subseteq S_{k+1}\).  It also implies literal inclusions in the
fixed ambient algebra
\begin{equation}
\label{def:eq:J-inclusions}
                         J_k\subseteq J_{k+1}\subseteq M.
\end{equation}

\section{Three ideals attached to a stage map}
\label{def:sec:three-ideals}

The image ideal depends on the chosen retraction.  Two related ideals depend
only on the compatible map \(\alpha_k:M\to A_k\):
\begin{equation}
\label{def:eq:three-ideals}
 I_k^{\alpha}:=\ker\alpha_k,
 \qquad
 \widehat I_k:=\ker(q_k\alpha_k)
       =\alpha_k^{-1}\bigl((A_k^+)^{n+1}\bigr),
 \qquad
 J_k=\rho_k(K_k).
\end{equation}
All three are differential ideals of \(M\).
\index{Ikalpha@$I_k^\alpha$}
\index{Ihatk@$\widehat I_k$}
\index{Jk@$J_k$}

\begin{proposition}
\label{def:prop:canonical-inclusions}
For every strictly coordinated stage,
\[
 I_k^{\alpha}\subseteq\widehat I_k\subseteq J_k,
 \qquad
 \widehat I_k\subseteq I_{k+1}^{\alpha},
 \qquad
 J_k\subseteq J_{k+1}.
\]
\end{proposition}

\begin{proof}
The first inclusion is immediate.  If \(m\in\widehat I_k\), then
\(\pi_ki_k(m)=q_k\alpha_k(m)=0\), so \(i_k(m)\in K_k\).  Applying the
retraction gives
\[
                         m=\rho_ki_k(m)\in J_k.
\]
If \(t_k:A_k\to A_{k+1}\) denotes the transition and
\(\alpha_{k+1}=t_k\alpha_k\), then the construction factors \(t_k\) through
\(q_k\).  Hence \(q_k\alpha_k(m)=0\) implies
\(\alpha_{k+1}(m)=0\), proving
\(\widehat I_k\subseteq I_{k+1}^{\alpha}\).  The last inclusion is
\eqref{def:eq:J-inclusions}.
\end{proof}

The quotients have distinct meanings.  The map \(q_k\alpha_k\) factors
through \(M/I_k^{\alpha}\), and the kernel of the induced map is
\(\widehat I_k/I_k^{\alpha}\).  Thus \(M/\widehat I_k\) is the largest
quotient of \(M\) through which \(q_k\alpha_k\) factors injectively.  By
contrast, \(M/J_k\) is the quotient selected by the retraction.  The map
\(q_k\alpha_k\) factors through \(M/J_k\) if and only if
\(J_k=\widehat I_k\).

\section{Colimits and finite death}
\label{def:sec:colimits}

Define
\begin{equation}
\label{def:eq:KSJ-colimits}
 K_\infty:=\underset{k}{\operatorname{colim}}\,K_k,
 \qquad
 S_\infty:=\underset{k}{\operatorname{colim}}\,S_k,
 \qquad
 J_\infty:=\underset{k}{\operatorname{colim}}\,J_k
           =\bigcup_kJ_k\subseteq M.
\end{equation}
Filtered colimits of complexes of rational vector spaces are exact.
Taking the colimit of~\eqref{def:eq:butterfly-sequence} yields
\[
 0\longrightarrow S_\infty\longrightarrow K_\infty
   \xrightarrow{\rho_\infty}J_\infty\longrightarrow0.
\]
Moreover,
\[
 H(K_\infty)\cong\underset{k}{\operatorname{colim}}\,H(K_k)=0.
\]
The connecting morphism and exactness of filtered colimits therefore give
natural isomorphisms
\begin{equation}
\label{def:eq:colimit-cohomology}
 H^q(J_\infty)\cong H^{q+1}(S_\infty)
 \cong\underset{k}{\operatorname{colim}}\,H^q(J_k).
\end{equation}

\begin{proposition}[Long words lie in the image ideal]
\label{def:prop:long-words}
For every \(k\),
\[
                          (M^+)^{n+1}\subseteq J_k,
\]
and consequently
\[
                          (M^+)^{n+1}\subseteq J_\infty.
\]
In particular, \(M/J_\infty\) is \(n\)-short.
\end{proposition}

\begin{proof}
Let \(m\in(M^+)^{n+1}\).  Since \(\alpha_k\) preserves augmentations,
\(\alpha_k(m)\in(A_k^+)^{n+1}\), and hence
\(q_k\alpha_k(m)=0\).  Therefore \(i_k(m)\in K_k\), and
\[
                         m=\rho_ki_k(m)\in J_k.
\]
\end{proof}

We now make the transition condition used below completely explicit.

\begin{definition}
\label{def:def:essentially-zero}
A directed system of graded vector spaces
\((V_k,v_{k\ell})\) is \emph{essentially zero} if, for every index \(k\)
and every element \(x\in V_k\), there exists an index \(\ell\geq k\),
possibly depending on \(x\), such that \(v_{k\ell}(x)=0\).  No transition
map need vanish identically, and no bound on \(\ell\) uniform in \(x\) is
required.
\index{essentially zero directed system}
\end{definition}

\begin{theorem}[Finite-death criterion]
\label{def:thm:finite-death}
\index{finite-death criterion}
The following conditions are equivalent:
\begin{enumerate}[label=\textup{(\roman*)}]
\item \(J_\infty\) is acyclic;
\item for every \(k,q\), and every \(\xi\in H^q(J_k)\), there exists
      \(\ell\geq k\) such that the image of \(\xi\) in \(H^q(J_\ell)\)
      is zero;
\item the directed system \((H^*(J_k))_k\) is essentially zero.
\end{enumerate}
\end{theorem}

\begin{proof}
By~\eqref{def:eq:colimit-cohomology}, \(H(J_\infty)\) is the filtered
colimit of the \(H(J_k)\).  An element represented at stage \(k\) is zero
in a filtered colimit precisely when it is sent to zero at some finite
later stage.  This proves the equivalence of the three conditions.
\end{proof}

Thus there is no purely infinite disappearance: if a class dies in the
union, one of its primitives already occurs at a finite stage.

\section{The short colimit and the quotient comparison}
\label{def:sec:quotient-comparison}

In the interlaced system
\[
 A_0\longrightarrow C_0\longrightarrow A_1\longrightarrow C_1
 \longrightarrow\cdots
\]
both the even and odd subsequences are cofinal.  Consequently,
\begin{equation}
\label{def:eq:A-colimit}
 A_\infty:=\underset{k}{\operatorname{colim}}\,A_k
 \cong\underset{k}{\operatorname{colim}}\,C_k.
\end{equation}
The complete-tower construction supplies a quasi-isomorphism
\begin{equation}
\label{def:eq:alpha-infinity}
                         \alpha_\infty:M\xrightarrow{\simeq}A_\infty.
\end{equation}
Since every \(C_k\) is \(n\)-short, so is \(A_\infty\).

Set
\begin{equation}
\label{def:eq:I-infinity}
 I_\infty:=\bigcup_k I_k^{\alpha}.
\end{equation}
The interlacing
\(I_k^{\alpha}\subseteq\widehat I_k\subseteq I_{k+1}^{\alpha}\)
shows that
\begin{equation}
\label{def:eq:I-infinity-identities}
 I_\infty=\bigcup_kI_k^{\alpha}
          =\bigcup_k\widehat I_k
          =\ker\alpha_\infty,
 \qquad
 (M^+)^{n+1}\subseteq I_\infty\subseteq J_\infty.
\end{equation}
Thus \(\alpha_\infty\) induces a canonical injection
\[
              \overline\alpha_\infty:M/I_\infty
                   \lhook\joinrel\longrightarrow A_\infty.
\]
It need not be surjective, and the fact that \(\alpha_\infty\) is a
quasi-isomorphism does not by itself make \(I_\infty\) acyclic.

At a finite stage define
\begin{equation}
\label{def:eq:rho-bar-k}
 \overline\rho_k:C_k\longrightarrow M/J_k,
 \qquad
 \overline\rho_k(\pi_k(e))=[\rho_k(e)]_{J_k}.
\end{equation}
This is well defined because changing \(e\) by an element of \(K_k\)
changes \(\rho_k(e)\) by an element of \(J_k\).  It is surjective, since
\[
 \overline\rho_k(q_k\alpha_k(m))=[m]_{J_k}.
\]
Equations~\eqref{def:eq:coordination-identities} make these maps
compatible with the transitions.  Passing to the colimit gives a canonical
surjection
\begin{equation}
\label{def:eq:r-infinity}
                       r_\infty:A_\infty\twoheadrightarrow M/J_\infty
\end{equation}
satisfying
\begin{equation}
\label{def:eq:r-alpha}
       r_\infty\alpha_\infty=p_{J_\infty}:M\twoheadrightarrow M/J_\infty.
\end{equation}
Equivalently, there is a strictly commutative square
\[
\begin{tikzcd}[column sep=large]
 M/I_\infty
   \arrow[r,hook,"\overline\alpha_\infty"]
   \arrow[d,two heads]
 & A_\infty \arrow[d,two heads,"r_\infty"]\\
 M/J_\infty \arrow[r,equal] & M/J_\infty .
\end{tikzcd}
\]
Define the \emph{comparison defect}
\begin{equation}
\label{def:eq:T-infinity}
                            T_\infty:=\ker r_\infty.
\end{equation}
\index{Tinfinity@$T_\infty$}

\begin{theorem}[Kernel comparison]
\label{def:thm:kernel-comparison}
The restriction of \(\alpha_\infty\) is a quasi-isomorphism
\[
               \alpha_\infty|_{J_\infty}:J_\infty
                    \xrightarrow{\simeq}T_\infty.
\]
Consequently,
\[
 r_\infty\text{ is a quasi-isomorphism}
 \quad\Longleftrightarrow\quad
 J_\infty\text{ is acyclic}.
\]
\end{theorem}

\begin{proof}
Equation~\eqref{def:eq:r-alpha} sends \(J_\infty\) to zero, so the
restriction has the displayed target.  Consider the morphism of short exact
sequences
\[
\begin{tikzcd}[column sep=small]
 0 \arrow[r]
 & J_\infty \arrow[r,hook]
   \arrow[d,"\alpha_\infty|_{J_\infty}"]
 & M \arrow[r,two heads,"p_{J_\infty}"]
   \arrow[d,"\alpha_\infty"']
 & M/J_\infty \arrow[r]
   \arrow[d,"\operatorname{id}"']
 & 0\\
 0 \arrow[r]
 & T_\infty \arrow[r,hook]
 & A_\infty \arrow[r,two heads,"r_\infty"]
 & M/J_\infty \arrow[r]
 & 0 .
\end{tikzcd}
\]
The associated long exact cohomology sequences form a commutative ladder.
The middle vertical map and the right vertical map are isomorphisms in
cohomology.  The five lemma, applied degree by degree to the ladder, shows
that the left vertical map is also an isomorphism in cohomology.  Hence it is
a quasi-isomorphism.

Since \(r_\infty\) is surjective, it is a quasi-isomorphism exactly when its
kernel \(T_\infty\) is acyclic.  The first assertion identifies this with
the acyclicity of \(J_\infty\).
\end{proof}

This theorem is the precise bridge between the flexible and rigid problems.
Every complete retractive tower supplies the short model \(A_\infty\).  For
a fixed complete strictly coordinated tower, the construction also supplies an
\(n\)-short quotient of the fixed minimal algebra by an acyclic ideal
precisely when the additional comparison \(r_\infty\) is a
quasi-isomorphism.

The inclusion \(I_\infty\subseteq J_\infty\) further separates the comparison
defect into two components.  Put
\[
 C_\infty^{\mathrm{rel}}
 :=A_\infty/\overline\alpha_\infty(M/I_\infty).
\]
The preceding commutative square and the snake lemma give a short exact
sequence
\begin{equation}
\label{def:eq:relative-defect}
 0\longrightarrow J_\infty/I_\infty\longrightarrow T_\infty
 \longrightarrow C_\infty^{\mathrm{rel}}\longrightarrow0.
\end{equation}
Thus \(T_\infty\) simultaneously records the difference between the
retraction-dependent image ideal and \(\ker\alpha_\infty\), and the failure
of \(M/I_\infty\to A_\infty\) to be onto.

\section{The finite augmentation spectral sequence}
\label{def:sec:augmentation-spectral-sequence}

We use the following page convention throughout this section.  For a
decreasing cohomological filtration \(F^pK\),
\[
 E_0^{p,q}=F^pK^{p+q}/F^{p+1}K^{p+q},
 \qquad d_r:E_r^{p,q}\longrightarrow E_r^{p+r,q-r+1}.
\]
For an increasing filtration \(G_pK\), the stage spectral sequence below is
indexed so that
\({}^{\mathrm{st}}d_r:{}^{\mathrm{st}}E_r^{p,q}\to
{}^{\mathrm{st}}E_r^{p-r,q+r+1}\).  These conventions account for the
opposite horizontal directions of the two spectral sequences.

Let \(Q_\infty=M/J_\infty\).  Both \(A_\infty\) and \(Q_\infty\) are
\(n\)-short.  Filter the comparison defect by
\begin{equation}
\label{def:eq:augmentation-filtration}
 \mathcal F^pT_\infty:=T_\infty\cap(A_\infty^+)^p,
 \qquad 1\leq p\leq n+1.
\end{equation}
Then
\[
 T_\infty=\mathcal F^1T_\infty,
 \qquad
 \mathcal F^{n+1}T_\infty=0.
\]
Every surjective augmented algebra map is strict for the augmentation-power
filtrations.  Hence
\[
 r_\infty((A_\infty^+)^p)=(Q_\infty^+)^p,
 \qquad
 \operatorname{gr}_{\mathcal F}^pT_\infty
   \cong\ker(\operatorname{gr}^p r_\infty).
\]

\begin{theorem}[Augmentation spectral sequence]
\label{def:thm:augmentation-ss}
The finite decreasing filtration~\eqref{def:eq:augmentation-filtration}
gives a strongly convergent cohomological spectral sequence
\[
 E_1^{p,q}\cong
 H^{p+q}\!\left(\ker(\operatorname{gr}^p r_\infty)\right)
 \quad\Longrightarrow\quad
 H^{p+q}(T_\infty).
\]
Its differentials have bidegree
\[
                    d_r:E_r^{p,q}\longrightarrow E_r^{p+r,q-r+1}.
\]
There are only the columns \(1\leq p\leq n\), so \(E_n=E_\infty\).
Consequently,
\[
 r_\infty\text{ is a quasi-isomorphism}
 \quad\Longleftrightarrow\quad E_n=0.
\]
\end{theorem}

\begin{proof}
The differential preserves every power of the augmentation ideal, so
\(\mathcal F^pT_\infty\) is a finite filtration by subcomplexes.  It is
therefore exhaustive, complete, and strongly convergent without any
conditional-convergence issue.  The strictness calculation identifies its
associated graded complex and hence the \(E_1\)-page.  Since a differential
\(d_r\) raises the column by \(r\), every \(d_r\) with \(r\geq n\) is zero.
Thus \(E_n=E_\infty\).  The limiting page vanishes exactly when the
associated graded of \(H(T_\infty)\) vanishes, which, for a finite
filtration, is equivalent to \(H(T_\infty)=0\).  Apply
Theorem~\ref{def:thm:kernel-comparison}.
\end{proof}

When \(n=3\), the possible differentials are
\[
 d_1:E_1^{1,*}\longrightarrow E_1^{2,*},
 \qquad
 d_1:E_1^{2,*}\longrightarrow E_1^{3,*},
 \qquad
 d_2:E_2^{1,*}\longrightarrow E_2^{3,*-1},
\]
and \(E_3=E_\infty\).  Hence \(r_\infty\) is a quasi-isomorphism exactly
when
\[
 E_2^{2,*}=0
 \quad\text{and}\quad
 d_2:E_2^{1,*}\xrightarrow{\cong}E_2^{3,*-1}.
\]
This criterion applies to the fixed strictly coordinated tower; another
short witness or compatible branch may have a different spectral sequence.

\section{Finite stages and the filtration by stage}
\label{def:sec:stage-filtrations}

At stage \(k\), put
\[
 T_k:=\ker\bigl(\overline\rho_k:C_k\twoheadrightarrow M/J_k\bigr),
 \qquad
 \mathcal F^pT_k:=T_k\cap(C_k^+)^p.
\]
Strict coordination gives compatible maps \(T_k\to T_{k+1}\).  A product of
finitely many elements in a filtered colimit is represented at a common
finite stage.  It follows that
\[
 \mathcal F^pT_\infty\cong
      \underset{k}{\operatorname{colim}}\,\mathcal F^pT_k,
 \qquad
 \operatorname{gr}^pT_\infty\cong
      \underset{k}{\operatorname{colim}}\,\operatorname{gr}^pT_k.
\]
Exactness of filtered colimits, followed by induction on the page, gives
\begin{equation}
\label{def:eq:pagewise-colimit}
 E_r^{p,q}(T_\infty)\cong
       \underset{k}{\operatorname{colim}}\,E_r^{p,q}(T_k)
 \qquad(r\geq1).
\end{equation}
Thus a class represented at a finite stage vanishes on the global
\(E_r\)-page precisely when that representative maps to zero at some later
finite stage.

There is also an increasing filtration of \(J_\infty\) by the first stage
at which an element occurs.  Set \(J_{-1}=0\) and
\(\mathcal G_pJ_\infty=J_p\).  The finite-type assumption is essential at
this point: in each cohomological degree \(m\), the increasing chain
\[
 J_0^m\subseteq J_1^m\subseteq\cdots\subseteq M^m
\]
stabilizes because \(M^m\) is finite dimensional.  The filtration is
therefore regular degree by degree and yields a strongly convergent spectral
sequence
\begin{equation}
\label{def:eq:stage-ss}
 {}^{\mathrm{st}}E_1^{p,q}\cong H^{p+q}(J_p/J_{p-1})
 \quad\Longrightarrow\quad H^{p+q}(J_\infty),
\end{equation}
with
\[
 {}^{\mathrm{st}}d_r:
 {}^{\mathrm{st}}E_r^{p,q}\longrightarrow
 {}^{\mathrm{st}}E_r^{p-r,q+r+1}.
\]
An immediate primitive is detected by a \(d_1\)-differential; a primitive
first appearing \(r\) stages later is detected, subject to the usual survival
hypotheses, by a \(d_r\)-differential.  Absent finite type or another proof of
degreewise regularity, strong convergence of this second spectral sequence
cannot be asserted.

\paragraph{Scope.}\label{def:sec:comparison-scope}
For a fixed complete strictly coordinated level-\(n\) tower,
\[
 J_\infty\text{ acyclic}
 \quad\Longleftrightarrow\quad
 (H^*(J_k))_k\text{ essentially zero}
 \quad\Longleftrightarrow\quad
 r_\infty\text{ a quasi-isomorphism}.
\]
When these conditions hold, \(M/J_\infty\) is the required \(n\)-short
quotient.  This is a statement about the fixed branch and fixed retractions;
the global quantifier convention is the one recorded at the beginning of
Part~\ref{part:towers}.  Chapter~\ref{chap:tate-selection} now determines
which classes can be killed after one such retraction has been fixed.

\chapter{Tate Selection with a Fixed Retraction}
\label{chap:tate-selection}
\begin{chapterguide}[title={Chapter guide}]
After a retraction has been chosen, a finite family of classes in \(H(J_k)\)
can be killed by persistent contractible extensions exactly when those
classes vanish in \(H(M)\).  This chapter develops the Tate-pair construction,
proves persistence and finite selection, and then passes to degreewise
reduction and acyclic envelopes.  The criterion remains relative to the
chosen retraction.
\end{chapterguide}
\index{Tate pair}
\index{Tate tower}
\index{acyclic envelope}
\index{ambient obstruction}

\section{Tate and asymptotically Tate towers}
\label{tate:sec:towers}

Let \(M\) be a simply connected minimal Sullivan algebra of finite type.
For a fixed complete strictly coordinated level-\(n\) retractive tower, retain the
notation
\[
 K_k=\ker\pi_k,\qquad
 S_k=K_k\cap\ker\rho_k,\qquad
 J_k=\rho_k(K_k)\subseteq M.
\]

\begin{definition}
\label{tate:def:towers}
A complete strictly coordinated level-\(n\) retractive tower is a \emph{Tate tower}
if every transition
\[
                         H^*(J_k)\longrightarrow H^*(J_{k+1})
\]
is zero.  It is \emph{asymptotically Tate} if the directed system
\((H^*(J_k))_k\) is essentially zero in the sense of
Definition~\ref{def:def:essentially-zero}.
\end{definition}

Naturality of the connecting isomorphisms
\(H^q(J_k)\cong H^{q+1}(S_k)\) shows that the Tate condition is equivalent
to the vanishing of every transition \(H^*(S_k)\to H^*(S_{k+1})\).  A Tate
tower is asymptotically Tate, but an asymptotically Tate tower need not have
a one-step death bound.

\begin{proposition}[Asymptotic implication]
\label{tate:prop:asymptotic-implication}
If \(M\) admits a complete level-\(n\) asymptotically Tate tower, then
\[
                         \Hnil(M)\leq n.
\]
\end{proposition}

\begin{proof}
The finite-death criterion, Theorem~\ref{def:thm:finite-death}, gives
\(H(J_\infty)=0\).  Proposition~\ref{def:prop:long-words} gives
\((M^+)^{n+1}\subseteq J_\infty\).  Hence the quotient
\[
                         M\twoheadrightarrow M/J_\infty
\]
is a quasi-isomorphism and its augmentation ideal has nilpotency at most
\(n\).  This is exactly the required bound on homology nilpotency.
\end{proof}

\begin{remark}[No converse is being asserted]
\label{tate:rem:no-converse}
Proposition~\ref{tate:prop:asymptotic-implication} is intentionally
one-directional.  Starting from an acyclic ideal \(I\) with
\((M^+)^{n+1}\subseteq I\) would require a quotient-compatible lifting
theorem that chooses all retractions and strict coordination maps, strictly
and coherently, so that every resulting image ideal is contained in \(I\).
The local Tate constructions below do not provide such a global choice.
Thus this chapter does not identify
\(\Hnil(M)\leq n\) with the existence of an
asymptotically Tate tower.
\end{remark}

\section{Stable Tate extensions and their homotopies}
\label{tate:sec:stable-extension}

Fix one retractive stage
\begin{equation}
\label{tate:eq:fixed-stage}
 M\xrightarrow{i}E\xrightarrow[\simeq]{\pi}C,
 \qquad
 \rho:E\longrightarrow M,
 \qquad
 \rho i=\id_M,
\end{equation}
where \(\pi\) is surjective.  Put
\[
                 K=\ker\pi,\qquad
                 S=K\cap\ker\rho,\qquad
                 J=\rho(K).
\]
A \emph{stable Tate extension} of this stage is obtained by tensoring \(E\)
with one or more contractible Sullivan disks
\[
 D(e):=(\Lambda(e,u),d),\qquad de=u,\qquad du=0,
\]
extending \(\pi\) by sending the new generators to zero, and extending
\(\rho\) by prescribed elements of \(M\).  The augmentation
\(D(e)\to\Q\) is a quasi-isomorphism, so the extended projection
remains a surjective quasi-isomorphism.
We call \((e,u)\) a \emph{Tate pair} and \(D(e)\) its \emph{Tate disk}.

Often the stage is coordinated by additional data
\begin{equation}
\label{tate:eq:coordinated-data}
 \phi:M\xrightarrow{\simeq}B,\qquad
 f:E\longrightarrow B,\qquad
 \mathcal H:E\longrightarrow B[t,dt],
\end{equation}
where \(B[t,dt]=B\otimes\Lambda(t,dt)\),
\(|t|=0\), \(|dt|=1\), \(d(t)=dt\), and
\[
 \varepsilon_0\mathcal H=f,\qquad
 \varepsilon_1\mathcal H=\phi\rho.
\]
The homotopy is relative to \(M\).  A stable Tate extension will be called
\emph{coordinated} when \(f\) and \(\mathcal H\) are extended with these
same endpoint conditions.

The next lemma supplies the missing explicit formula.

\begin{lemma}[Cylinder extension over a Tate disk]
\label{tate:lem:homotopy-extension}
Assume the coordinated data~\eqref{tate:eq:coordinated-data} have been
fixed.  Let \(b\in M^{m}\) and set \(j=db\).  Adjoin a Tate disk with
\[
 |e|=m,\qquad |u|=m+1,\qquad de=u.
\]
Extend the projection and the map to the witness by
\[
 \pi^+(e)=\pi^+(u)=0,\qquad
 f^+(e)=f^+(u)=0,
\]
and extend the retraction by
\[
                       \rho^+(e)=b,\qquad \rho^+(u)=j.
\]
The cylinder homotopy has a unique extension as an algebra map, determined by
\begin{align}
\label{tate:eq:homotopy-extension}
 \mathcal H^+(e)&=\phi(b)t,\\
 \label{tate:eq:homotopy-extension-u}
 \mathcal H^+(u)&=\phi(j)t+(-1)^m\phi(b)\,dt.
\end{align}
It satisfies
\[
 \varepsilon_0\mathcal H^+=f^+,\qquad
 \varepsilon_1\mathcal H^+=\phi\rho^+,
\]
and remains relative to \(M\).  If \(K_1\) denotes polynomial integration
from the initial to the terminal endpoint and \(F^+=K_1\mathcal H^+\), then
the integrated operator extends the original one by
\begin{equation}
\label{tate:eq:integrated-extension}
                         F^+(e)=0,\qquad F^+(u)=\phi(b).
\end{equation}
\end{lemma}

\begin{proof}
Coefficients from \(B\) are written to the left of \(t\) and \(dt\).  The
Leibniz rule gives
\[
 d(\phi(b)t)=\phi(db)t+(-1)^{|b|}\phi(b)\,dt
            =\phi(j)t+(-1)^m\phi(b)\,dt.
\]
Thus \(d\mathcal H^+(e)=\mathcal H^+(u)\); also
\(d\mathcal H^+(u)=0\), either by direct calculation or by applying
\(d^2\) to \(\mathcal H^+(e)\).  Since \(E_{\mathrm T}=E\otimes\Lambda(e,u)\) is
free over \(E\), these formulas and the original \(\mathcal H\) define a
CDGA morphism.

At \(t=0\), the images of both new generators vanish.  At \(t=1\), the \(dt\)-term
vanishes and the images are \(\phi(b)\) and \(\phi(j)\), exactly
\(\phi\rho^+(e)\) and \(\phi\rho^+(u)\).  Nothing was changed on \(M\), so
the extended homotopy is still relative to \(M\).

Finally, polynomial integration annihilates terms without \(dt\).  With the
coefficient placed to the left, it sends
\(\phi(b)\,dt\) to \((-1)^m\phi(b)\).
Equations~\eqref{tate:eq:homotopy-extension}--\eqref{tate:eq:homotopy-extension-u}
therefore give \(F^+(e)=0\) and \(F^+(u)=\phi(b)\).  In particular,
\[
 dF^+(e)+F^+(u)=\phi(b),\qquad
 dF^+(u)=\phi(j),
\]
which is the chain-homotopy identity on the two new generators.
\end{proof}

\section{The exact criterion for one Tate pair}
\label{tate:sec:one-pair}

Let \(s\in Z^{q+1}(S)\).  Since \(K\) is acyclic, choose
\(c\in K^q\) with \(dc=s\).  Then
\[
 \tau^J_\rho([s])=[\rho(c)]_J,\qquad
 \omega_\rho([s])=[\rho(c)]_M.
\]

\begin{theorem}[Tate-pair criterion]
\label{tate:thm:one-pair}
The following conditions are equivalent:
\begin{enumerate}[label=\textup{(\roman*)}]
\item \(\omega_\rho([s])=0\) in \(H^q(M)\);
\item there is a stable Tate extension of the stage, coordinated whenever
      data~\eqref{tate:eq:coordinated-data} are present, containing an
      element \(\widehat c\) in the double kernel of the extended stage such
      that \(d\widehat c=s\).
\end{enumerate}
More precisely:
\begin{enumerate}[label=\textup{(\alph*)}]
\item if \(\tau^J_\rho([s])=0\), no new Tate pair is needed;
\item if \(\tau^J_\rho([s])\neq0\) but
      \(\omega_\rho([s])=0\), one contractible Tate pair gives the
      required extension.
\end{enumerate}
\end{theorem}

\begin{proof}
Suppose first that \(\omega_\rho([s])=0\).  Write
\[
                         j=\rho(c)\in Z^q(J)
\]
and choose \(b\in M^{q-1}\) with \(db=j\).  Form
\[
 E_{\mathrm T}=E\otimes\Lambda(e_{q-1},u_q),
 \qquad de=u,\qquad du=0,
\]
kill \(e,u\) under the extended projection, and put
\[
                         \rho^+(e)=b,\qquad \rho^+(u)=j.
\]
Lemma~\ref{tate:lem:homotopy-extension} extends the map to the witness and
cylinder homotopy explicitly.  The element
\[
                              \widehat c=c-u
\]
lies in \(\ker\pi^+\), satisfies
\(\rho^+(\widehat c)=j-j=0\), and has
\[
                              d\widehat c=s.
\]
This proves~(ii).

If the intrinsic class is already zero, choose \(y\in J^{q-1}\) with
\(dy=j\), and then choose \(z\in K^{q-1}\) with \(\rho(z)=y\).  The
uncorrected primitive may be replaced inside the original stage by
\[
                              c'=c-dz.
\]
Indeed, \(c'\in K\), \(dc'=s\), and
\(\rho(c')=j-dy=0\).  This proves~(a), while the preceding disk proves~(b).

Conversely, suppose a stable Tate extension \(E\to E'\) contains
\(\widehat c\in K'\cap\ker\rho'\) with \(d\widehat c=s\), where \(c\) and
\(s\) are viewed in \(E'\).  Then \(c-\widehat c\) is a cocycle in
\(K'=\ker\pi'\).  The extended projection is a surjective
quasi-isomorphism, so \(K'\) is acyclic.  Choose \(v\in K'^{q-1}\) with
\[
                              dv=c-\widehat c.
\]
Applying the extended retraction gives
\[
                   \rho(c)=d\rho'(v),
\]
so \([\rho(c)]_M=0\), which is~(i).
\end{proof}

The theorem separates the intrinsic and ambient questions.  Intrinsic
triviality permits a correction already inside the original double kernel.
Ambient inessentiality permits a correction only after enlarging that
kernel by a contractible disk.

\section{Persistence under strict coordination}
\label{tate:sec:persistence}

Any such disk must be carried through all later stages.  Merely adjoining it
at one isolated stage would not define a new directed system.

\begin{lemma}[Persistent transport of Tate pairs]
\label{tate:lem:persistence}
Let a complete strictly coordinated tower be fixed, and let
\(j\in Z^q(J_k)\) satisfy \(j=db\) for some \(b\in M^{q-1}\).  Adjoin at
stage \(k\) the Tate pair
\[
 |e|=q-1,\qquad |u|=q,\qquad de=u,
\]
with retraction values \(\rho^+(e)=b\), \(\rho^+(u)=j\).  Then the pair can
be transported, by the identity on \(e,u\), through every later relative
algebra and every later intermediate algebra so that:
\begin{enumerate}[label=\textup{(\roman*)}]
\item all strict coordination identities remain literal;
\item every extended projection and every map to the fixed short witness
      kills \(e,u\);
\item every compatible cylinder homotopy extends by
      formulas~\eqref{tate:eq:homotopy-extension}--%
      \eqref{tate:eq:homotopy-extension-u};
\item the algebras \(A_\ell,C_\ell\) and the maps \(\alpha_\ell\) may be
      canonically identified with the original ones;
\item for every \(\ell\geq k\),
      \[
                     J_\ell^+=J_\ell+(b)
      \]
      as differential ideals in \(M\), and the image of \([j]\) is zero in
      \(H^q(J_\ell^+)\).
\end{enumerate}
The same assertions hold for a finite family of Tate pairs.
\end{lemma}

\begin{proof}
Tensor every later \(E_\ell\) and every later intermediate algebra
\(E'_\ell\) with the same disk \(D(e)\).  Extend each transition by the
identity on \(e,u\), each projection and each map to the witness by the
augmentation of \(D(e)\), and every retraction by the fixed assignments
\(e\mapsto b\), \(u\mapsto j\).  Since all original retractions agree on
transported elements, the strict coordination identities continue to hold on
the old algebras and on the two new generators.

Lemma~\ref{tate:lem:homotopy-extension} gives the same homotopy formula at
each stage, so the endpoint identities are preserved.  Because the new
disk is killed by every map to a short target, the corresponding pushouts
canonically remove it.  The systems of \(A_\ell,C_\ell,\alpha_\ell\) are
therefore unchanged up to the displayed canonical identifications.

It remains to compute the image ideal.  Write an element of the extended
kernel as its disk-free part plus terms containing \(e\) or \(u\).  The
disk-free part lies in \(K_\ell\), so its retraction lies in \(J_\ell\).
Every remaining term retracts into the differential ideal generated by
\(b\), because \(u\) retracts to \(db=j\in J_\ell\).  This proves
\(\rho_\ell^+(\ker\pi_\ell^+)\subseteq J_\ell+(b)\).  The reverse inclusion
follows from \(J_\ell=\rho_\ell(K_\ell)\) and
\(b=\rho_\ell^+(e)\).  Finally, \(db=j\) inside \(J_\ell+(b)\), so
\([j]\) is zero at every later extended stage.  Finite families are handled
by tensoring the corresponding disks and applying the same argument
componentwise.
\end{proof}

\section{Finite selection for a fixed retraction}
\label{tate:sec:finite-selection}

\begin{theorem}[Finite-packet criterion]
\label{tate:thm:finite-packet}
Fix a strictly coordinated tower, and let
\(V\subseteq H^*(J_k)\) be a finite-dimensional graded subspace.  The
following conditions are equivalent:
\begin{enumerate}[label=\textup{(\roman*)}]
\item a strictly coordinated continuation by stable Tate extensions,
      preserving all preceding data, sends \(V\) to zero at a finite later
      stage;
\item adjoining a finite family of Tate pairs and transporting it persistently
      as in Lemma~\ref{tate:lem:persistence} sends \(V\) to zero;
\item \(H(J_k\hookrightarrow M)\) vanishes on \(V\);
\item \(\omega_{\rho_k}\) vanishes on \(\delta_k(V)\subseteq H^{*+1}(S_k)\).
\end{enumerate}
\end{theorem}

\begin{proof}
If a class represented by \(j\in Z(J_k)\) becomes a boundary in a later
image ideal \(J_\ell\subseteq M\), it is already a boundary in \(M\).  This
proves~(i)\(\Rightarrow\)(iii).  Conditions~(iii) and~(iv) are equivalent
because
\[
 \omega_{\rho_k}\delta_k
 =H(J_k\hookrightarrow M)\tau^{J_k}_{\rho_k}\delta_k
 =H(J_k\hookrightarrow M).
\]

Assume~(iii).  Choose a homogeneous basis
\(\xi_i=[j_i]\) of \(V\), with \(j_i\in Z^{q_i}(J_k)\), and choose
\[
 b_i\in M^{q_i-1},\qquad db_i=j_i.
\]
Apply Theorem~\ref{tate:thm:one-pair} to each basis class, and transport the
resulting finite family by Lemma~\ref{tate:lem:persistence}.  Each \(j_i\)
becomes a boundary in every later enlarged image ideal.  This proves~(ii),
which is a special case of~(i).
\end{proof}

The criterion is therefore complete only after the retraction has been fixed.

\section{Degreewise Tate reduction}
\label{tate:sec:degreewise-reduction}

\begin{theorem}[Degreewise Tate reduction]
\label{tate:thm:degreewise-reduction}
Let \(M\) be simply connected and of finite type, and let
\[
 M\xrightarrow{i}E\xrightarrow[\simeq]{\pi}C,\qquad
 \rho:E\longrightarrow M
\]
be a retractive stage.  There is a degreewise locally finite stable Tate
extension \(E\to E^{\mathrm{red}}\), with the same target \(C\), for which
\[
 J^{\mathrm{red}}
 :=\rho^{\mathrm{red}}(\ker\pi^{\mathrm{red}})
\]
is Tate-reduced:
\[
             H(J^{\mathrm{red}})\longrightarrow H(M)
             \quad\text{is injective}.
\]
Only finitely many Tate pairs are required in each total degree.  If the
stage carries coordinated data, they extend by
Lemma~\ref{tate:lem:homotopy-extension}.
\end{theorem}

\begin{proof}
Proceed by increasing cohomological degree.  If every map
\(H^q(J)\to H^q(M)\) is injective, the image ideal is already Tate-reduced.
Otherwise, start with the least \(q\) for which its kernel is
nonzero.  Such a \(q\) is at least three:
\(J^0=J^1=0\), and an element of \(J^2\) that is a boundary in \(M\) has a
primitive in \(M^1=0\).  Suppose all lower degrees have been treated, and
write the current image ideal as \(J\).  Set
\[
          V_q:=\ker\bigl(H^q(J)\longrightarrow H^q(M)\bigr).
\]
The vector space \(V_q\) is finite dimensional because \(M^q\) is finite
dimensional.  Choose cycles \(j_1,\ldots,j_t\in J^q\) representing a basis
of \(V_q\), together with primitives
\[
                         b_i\in M^{q-1},\qquad db_i=j_i.
\]
Adjoin the corresponding finite family of Tate disks and put
\[
                         J^+=J+(b_1,\ldots,b_t).
\]

We verify that this operation does not disturb any already treated degree.
Since \(M^1=0\), every non-scalar product involving a new \(b_i\) has degree
at least \(q+1\).  Through degree \(q\), the quotient \(J^+/J\) is therefore
the vector space spanned by the residue classes of the \(b_i\), placed in
degree \(q-1\), with zero differential.  These residue classes are linearly
independent: a relation modulo \(J\) would, after applying \(d\), give a
linear relation among the cohomology classes \([j_i]\).  In the long exact
sequence of
\[
                  0\longrightarrow J\longrightarrow J^+
                    \longrightarrow J^+/J\longrightarrow0,
\]
the connecting morphism sends \([b_i]\) to \([j_i]\).  It is injective with
image \(V_q\).  Hence cohomology in degrees below \(q\) is unchanged and
\[
                         H^q(J^+)\cong H^q(J)/V_q.
\]
The induced map from this quotient to \(H^q(M)\) is injective.

Continuing in increasing degree produces a degreewise locally finite
filtered union.  At every finite step, the stage algebra is tensored with
finitely many contractible disks, so the projection remains a
quasi-isomorphism.  Exactness of filtered colimits preserves this property.
The preceding degree calculation shows that later steps cannot change any
already treated degree.
The limiting image ideal is therefore Tate-reduced.
\end{proof}

\begin{remark}
\label{tate:rem:reduced-not-acyclic}
Tate-reduced does not mean acyclic.  The construction removes precisely the
classes of the current image ideal that are boundaries in \(M\).  Adjoining a
primitive can create new multiplicative or contextual classes in higher
degrees.  This is why the degreewise construction must continue and why
vanishing of the initial ambient obstruction does not by itself supply an
acyclic ideal.
\end{remark}

\section{The acyclic-envelope criterion}
\label{tate:sec:acyclic-envelope}

\begin{proposition}[Stable acyclic-envelope criterion]
\label{tate:prop:acyclic-envelope}
Let \(M\) be simply connected and of finite type, and let
\(J=\rho(\ker\pi)\) be the image ideal of a fixed retractive stage.  The
following conditions are equivalent:
\begin{enumerate}[label=\textup{(\roman*)}]
\item the stage admits a degreewise locally finite stable Tate extension
      whose image ideal is acyclic;
\item there exists an acyclic differential ideal \(I\triangleleft M\) with
      \(J\subseteq I\).
\end{enumerate}
If \(I\) is generated as a differential ideal modulo \(J\) by finitely many
homogeneous elements, the extension may be chosen finite.  Any coordinated
map and homotopy extend explicitly over all the required disks.
\end{proposition}

\begin{proof}
Condition~(i) implies~(ii) by taking the image ideal of the extended
retraction.

Conversely, choose homogeneous elements \(b_\lambda\in I\) that generate
\(I\) as a differential ideal modulo \(J\).  Adjoin contractible pairs
\[
 |e_\lambda|=|b_\lambda|,\qquad
 |u_\lambda|=|b_\lambda|+1,\qquad
 de_\lambda=u_\lambda,\qquad du_\lambda=0,
\]
kill them under the projection, and define
\[
 \rho^+(e_\lambda)=b_\lambda,\qquad
 \rho^+(u_\lambda)=db_\lambda.
\]
The new projection is a surjective quasi-isomorphism and
\[
 \rho^+(\ker\pi^+)
   =J+(b_\lambda,db_\lambda)_\lambda
   =I.
\]
If coordinated data are present, extend the map to the witness by zero and
set
\[
 \mathcal H^+(e_\lambda)=\phi(b_\lambda)t,\qquad
 \mathcal H^+(u_\lambda)
  =\phi(db_\lambda)t+(-1)^{|b_\lambda|}\phi(b_\lambda)\,dt.
\]
Lemma~\ref{tate:lem:homotopy-extension} verifies that these formulas define
a CDGA morphism with the required endpoint values.  Finite type allows the
generators to be chosen so that only finitely many occur in each degree.  For
an infinite family, take the resulting degreewise filtered union.
\end{proof}

The proposition identifies the exact boundary of the fixed-retraction
method.  The phrase \emph{the image ideal may be assumed acyclic} is valid
only when an acyclic envelope has been supplied or when persistent ambient
inessentiality has been established throughout the full iterative Tate
construction.

\section{Persistent variation of the comparison map}
\label{tate:sec:r-variation}

Let a finite packet at stage \(k\) be represented by
\(j_i\in Z^{q_i}(J_k)\), and let \(db_i=j_i\).  Transport the associated
Tate pairs persistently.  Lemma~\ref{tate:lem:persistence} gives
\[
 J_\ell^+=J_\ell+(b_1,\ldots,b_m),\qquad
 J_\infty^+=J_\infty+(b_1,\ldots,b_m).
\]
The short-model system \((A_\ell,C_\ell,\alpha_\ell)\) is unchanged, and
hence so is \(A_\infty\).  Let
\[
 q_\infty^J:M/J_\infty\twoheadrightarrow M/J_\infty^+
\]
be the quotient.  The new comparison map is
\begin{equation}
\label{tate:eq:r-variation}
                         r_\infty^+=q_\infty^Jr_\infty.
\end{equation}
Thus transporting a family of Tate pairs persistently changes the comparison
by postcomposition with the quotient that adjoins the chosen ambient
primitives to the image ideal.

If
\[
 T_\infty=\ker r_\infty,\qquad
 T_\infty^+=\ker r_\infty^+,
\]
the snake lemma gives
\begin{equation}
\label{tate:eq:T-variation}
 0\longrightarrow T_\infty\longrightarrow T_\infty^+
 \longrightarrow J_\infty^+/J_\infty\longrightarrow0.
\end{equation}
In the exact sequence
\[
 0\longrightarrow J_\ell\longrightarrow J_\ell^+
 \longrightarrow J_\ell^+/J_\ell\longrightarrow0,
\]
the connecting morphism sends the residue class of \(b_i\) to \([j_i]\).
Equivalently, in the stage spectral sequence, this connecting morphism is
precisely the differential that kills the old class.

\section{Variation with the chosen retraction}
\label{tate:sec:retraction-variation}

The preceding criteria are exact only after \(\rho\) has been fixed.  The
following formula records how the ambient obstruction changes under a
homotopy of retractions.

Let \(\rho_0,\rho_1:E\to M\) be retractions of the same relative
factorization, put
\[
 S_{01}=K\cap\ker\rho_0\cap\ker\rho_1,
\]
and suppose that integrating a relative cylinder homotopy yields a
degree-\(-1\) operator \(F:E\to M\) satisfying
\begin{equation}
\label{tate:eq:F-homotopy}
                    \rho_1-\rho_0=dF+Fd.
\end{equation}

\begin{proposition}[Ambient variation]
\label{tate:prop:omega-variation}
For every \([s]\in H^{q+1}(S_{01})\),
\[
              \omega_{\rho_1}([s])-\omega_{\rho_0}([s])
                       =[F(s)]_M.
\]
\end{proposition}

\begin{proof}
Choose \(c\in K^q\) with \(dc=s\).  Equation~\eqref{tate:eq:F-homotopy}
gives
\[
 \rho_1(c)-\rho_0(c)=dF(c)+F(s).
\]
Passing to cohomology proves the formula.  If \(s\) is changed by a
boundary inside \(S_{01}\), the change in \(F(s)\) is a boundary because
both retractions vanish on \(S_{01}\).
\end{proof}

For a finite family \([s_i]\), the cohomological selection equations are
\[
                         [F(s_i)]_M=-\omega_{\rho_0}([s_i]).
\]
They are linear in the cohomological values of \(F\), but this does not make
the global choice problem linear.  The selected values must arise from a
multiplicative cylinder homotopy and from a genuine retraction satisfying
all relative Sullivan extension equations.

There is a related affine calculation after primitives have been chosen.
Suppose \(j_i\in J\), and choose primitives \(b_i^0\) satisfying
\(db_i^0=j_i\).  Every other ambient primitive
has the form
\[
                         b_i=b_i^0+z_i,\qquad dz_i=0.
\]
Let \(a_{\lambda i}\) be cocycles and let \(w_\lambda\in J\) satisfy
\[
 dw_\lambda
   =-\sum_i(-1)^{|a_{\lambda i}|}a_{\lambda i}j_i.
\]
Then
\[
 \Omega_\lambda(\mathbf b)
   =w_\lambda+\sum_i a_{\lambda i}b_i
\]
is a cycle in the enlarged ideal, and
\[
 [\Omega_\lambda(\mathbf b^0+\mathbf z)]
   =[\Omega_\lambda(\mathbf b^0)]
      +\sum_i[a_{\lambda i}z_i].
\]
Thus the dependence of any specified finite family of contextual obstructions
on the chosen primitives is governed by a linear map on the available
cohomology classes \([z_i]\).  This test is complete only for that family; a
global proof must also show that no unrecorded context creates a new
obstruction.

\paragraph{Scope and fixed data.}\label{tate:sec:scope}
The complete diagram and symbol dictionary are on
page~\pageref{rt:convention:quantifiers}.  In this chapter the short witness,
the strictly coordinated branch, and the retractions \(\rho_k\) are fixed.
For that fixed data, Theorem~\ref{tate:thm:finite-packet} characterizes the
finite packets that can be killed, Lemma~\ref{tate:lem:persistence} transports
the disks, and Theorem~\ref{def:thm:finite-death} converts finite death into
acyclicity of \(J_\infty\).  The prior global choice remains outside this
chain: a branchwise conclusion is not automatically universal, and an
acyclic quotient does not by itself produce an asymptotically Tate tower.

\begin{keyidea}{Transition to Part~\ref{part:divergence}}
The tower formalism has separated the flexible short witness from the rigid
quotient defect.  The next part constructs an obstruction that survives every
complete strictly coordinated level-three branch of one explicit model.
\end{keyidea}

\part{Divergence at Length Three}\label{part:divergence}

\chapter{The Three-Short Witness and Its Rigid Minimal Prefix}
\label{div:chap:witness}
\label{div:sec:witness}
\begin{chapterguide}[title={Chapter guide}]
Part~\ref{part:towers} isolated the strictification defect of a fixed
coordinated tower.  We now construct a length-three algebra in which that
defect survives every tower.  This chapter builds the flexible witness and
identifies the rigid low-degree prefix of its minimal Sullivan model.  Before
truncation, adjoining \(r\) with \(dr=xt\) violates the square-zero condition
by the length-four term \(d^2r=-xays\); the three-short quotient kills
precisely that term.  We verify the quotient directly, compare it degree by
degree with its minimal model, prove that \(xt\) is the unique degree-\(17\)
primitive of \(-xays\), and show that \(axt\) represents an essential class.
\end{chapterguide}
\index{three-short witness}
\index{augmentation truncation}
\index{curvature defect}

\section{The three-short witness}

Let
\[
                 F=\Lambda(s_3,x_5,y_3,a_7,t_{12},r_{16})
\]
be the free graded-commutative algebra on the displayed generators, and
let \(\mathfrak m_F=F^+\).  The generators \(s,x,y,a\) have odd degree,
whereas \(t,r\) have even degree.  Define
\[
                         T_3(F)=F/\mathfrak m_F^4.
\]

Consider the derivation of degree one on \(F\) determined by
\begin{equation}
\label{div:eq:differential}
 ds=dx=dy=0,\qquad da=xy,\qquad dt=ays,\qquad dr=xt.
\end{equation}
It is not a differential on \(F\).  The following two lemmas show that
it becomes one on \(T_3(F)\).

\begin{lemma}
\label{div:lem:stable-ideal}
The derivation in \eqref{div:eq:differential} preserves
\(\mathfrak m_F^4\).
\end{lemma}

\begin{proof}
The image under \(d\) of each generator is a sum of words of augmentation
length at least two.  Applying the derivation to a word of length at least
four replaces one factor by a word of length at least two.  Every resulting
term therefore lies in \(\mathfrak m_F^5\subseteq\mathfrak m_F^4\).
\end{proof}

\begin{lemma}
\label{div:lem:square-zero-B}
The derivation induced by \eqref{div:eq:differential} on \(T_3(F)\) has
square zero.
\end{lemma}

\begin{proof}
Because \(d\) has odd degree, \(d^2\) is a derivation.  It is therefore
enough to check the square on generators.  The only nontrivial
calculations are
\[
 d^2t=d(ays)=(xy)ys=0,
\]
because \(y\) has odd degree, and
\[
 d^2r=d(xt)=(-1)^{|x|}x\,dt=-xays.
\]
The last word has augmentation length four and is therefore zero in
\(T_3(F)\).  The square vanishes on every generator and hence on the
whole quotient.
\end{proof}

\begin{definition}
\label{div:def:B}
Let
\[
B=\bigl(T_3(F),d\bigr)
\]
denote the augmented CDGA defined by
\eqref{div:eq:differential}.
\end{definition}

\begin{proposition}
\label{div:prop:B-properties}
The CDGA \(B\) is simply connected, finite dimensional, and three-short.
More precisely,
\[
                         \dim_{\Q}B=56,\qquad (B^+)^4=0.
\]
\end{proposition}

\begin{proof}
The equality \((B^+)^4=0\) follows from the definition of \(T_3(F)\).
Also \(B^0=\Q\), and there are no elements in degrees one or two.

We give the dimension count explicitly.  Let
\[
                 O=\{s,x,y,a\},\qquad E=\{t,r\}.
\]
An odd generator can occur at most once, whereas the even generators
may be repeated.  In length two there are
\[
 \binom{4}{2}+4\cdot2+\binom{2+2-1}{2}=6+8+3=17
\]
monomials: two distinct odd generators, one odd and one even generator,
or a monomial of augmentation length two in the even generators.  In
length three,
separating according to the number \(k\) of odd factors gives
\[
\begin{array}{c|cccc}
k&0&1&2&3\\ \hline
\text{number of monomials}
 &\binom{2+3-1}{3}
 &4\binom{2+2-1}{2}
 &\binom{4}{2}\cdot2
 &\binom{4}{3}
\end{array}
\]
and hence \(4+12+12+4=32\) monomials.  The numbers in lengths
zero, one, two, and three are therefore
\[
                         1,\quad 6,\quad 17,\quad 32
\]
for a total dimension of \(56\).  Thus \(B\) is connected,
simply connected, finite dimensional, and three-short.
\end{proof}

\index{minimal Sullivan model!rigid prefix}
\index{essential cohomology class}
\index{degree calculation}

\section{The rigid low-degree part of the minimal model}
\label{div:sec:minimal-prefix}

Define the minimal Sullivan algebra
\[
P=\left(\Lambda(s_3,x_5,y_3,a_7,t_{12}),d\right),
\qquad
da=xy,\quad dt=ays,
\]
with \(ds=dx=dy=0\).  Since \(y\) is odd,
\[
                         d^2t=(xy)ys=0,
\]
so \(P\) is indeed a CDGA.  Its differential is decomposable, hence it
is minimal.  Sending each generator to the generator with the same name
defines a CDGA morphism
\[
                              \varphi_0:P\longrightarrow B.
\]
For example, \(d_P(xt)=-xays\), while \(d_B(xt)=0\), and
\(\varphi_0(xays)=0\) because this word has length four.

For calculations in \(P\) and \(B\), we fix the generator order
\[
                              s<x<y<a<t<r,
\]
with \(r\) omitted when working in \(P\).
For readability, displayed products retain the order produced by the
differential formulas.  Each displayed monomial basis is understood up to
Koszul sign; for example, \(axt=-xat\) in canonical order.  This convention
does not affect any span or cohomology class.  With this understood,
every monomial of \(P\) has the unique form
\[
 s^{\epsilon_s}x^{\epsilon_x}y^{\epsilon_y}a^{\epsilon_a}t^\alpha,
 \qquad
 \epsilon_s,\epsilon_x,\epsilon_y,\epsilon_a\in\{0,1\},
 \quad \alpha\geq0.
\]
Thus every low-degree basis calculation below reduces to the equation
\[
 3\epsilon_s+5\epsilon_x+3\epsilon_y+7\epsilon_a+12\alpha=q.
\]

\begin{lemma}[The first possible truncation relation]
\label{div:lem:first-relation}
There is no nonzero monomial of \(P\) of augmentation length at least
four in cohomological degree below \(18\).  In degree \(18\), the only
such monomial is \(sxya\), up to a nonzero scalar and a Koszul sign.
\end{lemma}

\begin{proof}
The four generators of smallest degree that can occur simultaneously
in a nonzero monomial are \(s,y,x,a\).  They are all odd and distinct,
and their degrees sum to
\[
                         3+3+5+7=18.
\]
Repeating any of them gives zero.  A monomial involving the even
generator \(t\) has even larger degree if its length is at least four.
\end{proof}

\begin{lemma}[The critical graded components]
\label{div:lem:critical-groups}
The graded components of \(P\) and \(B\) in degrees \(16\), \(17\), and
\(18\) are
\[
\begin{array}{c|c|c}
q&P^q&B^q\\ \hline
16&0&\Q\{r\}\\
17&\Q\{xt\}&\Q\{xt\}\\
18&\Q\{syt,sxya\}&\Q\{syt\}.
\end{array}
\]
Their differentials satisfy
\[
 d_Br=xt,\qquad
 d_P(xt)=-xays=-sxya,\qquad
 d_P(syt)=0.
\]
\end{lemma}

\begin{proof}
For degrees \(16\), \(17\), and \(18\), the exponent of \(t\) is at
most one.  With one \(t\), the required degrees of the odd part are
respectively \(4\), \(5\), and \(6\).  The first is impossible, the
unique degree-five odd monomial is \(x\), and the unique degree-six odd
monomial is \(sy\).  With no \(t\), there is no odd subset of degree
\(16\) or \(17\), while the unique odd subset of degree \(18\) is the
full product \(sxya\).  This proves
\[
 P^{16}=0,\qquad P^{17}=\Q\{xt\},\qquad
 P^{18}=\Q\{syt,sxya\}.
\]

In \(B\), the generator \(r\) supplies the single element in degree
\(16\).  It contributes no additional products in degrees \(17\) or
\(18\), since \(B\) has no positive-degree elements in degrees one or two.
Finally, the
length-four monomial \(sxya\) is zero in \(B\).  This gives the three
displayed bases.

The displayed differential formulas follow from the Leibniz rule.
For the last one,
\[
                         d(syt)=sy(ays)=0
\]
because the product contains repeated odd generators.
\end{proof}

\begin{proposition}[Low-degree comparison]
\label{div:prop:low-comparison}
The map
\[
                      H^q(\varphi_0):H^q(P)\longrightarrow H^q(B)
\]
is an isomorphism for every \(q\leq18\).  Moreover,
\[
                              H^{19}(P)=0.
\]
\end{proposition}

\begin{proof}
For \(q\leq14\), neither the generator \(r\) nor a truncation relation
enters the three graded components involved in computing \(H^q\), so the
relevant complexes agree.  In degree \(15\), the only additional term is
\(r\in B^{16}\).  Because \(d\) raises degree by one, only
\(d:B^{15}\to B^{16}\) could affect the cycle calculation, and no
differential formula has an \(r\)-component.  Thus \(d(B^{15})\) has
zero component along the new summand \(\Q\{r\}\).  Hence
\(H^{15}(P)\to H^{15}(B)\) is also an
isomorphism.

By Lemma~\ref{div:lem:critical-groups}, \(r\mapsto xt\) is an acyclic pair
in \(B\) across degrees \(16\) and \(17\).  In \(P\), the map
\[
                 d:P^{17}=\Q\{xt\}\longrightarrow
                    P^{18}
\]
has image \(\Q\{sxya\}\).  Therefore
\[
              H^{16}(P)=H^{16}(B)=0,\qquad
              H^{17}(P)=H^{17}(B)=0,
\]
and both \(H^{18}(P)\) and \(H^{18}(B)\) are generated by the class of
\(syt\).  The map induced by \(\varphi_0\) sends one generator to the
other.

Finally, the only monomial of \(P\) in degree \(19\) is \(at\), and
\[
                         d(at)=(xy)t-a(ays)=xyt\neq0,
\]
because \(a^2=0\).  Thus \(H^{19}(P)=0\).
\end{proof}

\begin{corollary}[Rigid minimal prefix]
\label{div:cor:rigid-prefix}
The morphism \(\varphi_0\) extends to a minimal Sullivan model
\[
             \varphi:M=P\otimes\Lambda W\xrightarrow{\simeq}B
\]
such that \(W^q=0\) for \(q\leq18\).  In particular,
\begin{equation}
\label{div:eq:degree17}
                       M^{17}=\Q\{xt\},\qquad Z^{17}(M)=0.
\end{equation}
\end{corollary}

\begin{proof}
Proposition~\ref{div:prop:low-comparison} gives the hypotheses of
Lemma~\ref{div:lem:minimal-extension} with \(N=18\): the map is an
isomorphism through degree \(18\), and
\(H^{19}(P)=0\to H^{19}(B)\) is injective.  Hence no new minimal
generator occurs in degree at most \(18\).  The statements in
\eqref{div:eq:degree17} follow from
\(P^{17}=\Q\{xt\}\) and \(d(xt)=-xays\neq0\) in the free Sullivan
algebra \(M\).
\end{proof}

\begin{remark}
\label{div:rem:finite-type}
Since \(B\) is finite dimensional and simply connected, its minimal
Sullivan model \(M\) is simply connected and of finite type; see
\cite[Proposition~12.2(ii)--(iii)]{FHT2001}.  Higher-degree generators of
\(M\) play no
role in \eqref{div:eq:degree17}, but they are essential for making
\(\varphi:M\to B\) a quasi-isomorphism in all degrees.
\end{remark}

\section{The essential class in degree 24}
\label{div:sec:essential-class}

We next identify the class that prevents ideal strictification.

\begin{lemma}
\label{div:lem:c-cycle}
The element
\[
                              c=axt\in M^{24}
\]
is a cycle.
\end{lemma}

\begin{proof}
Using \(da=xy\), \(dx=0\), and \(dt=ays\), the Leibniz rule gives
\[
\begin{aligned}
d(axt)
 &= (da)xt+(-1)^{|a|}a(dx)t+(-1)^{|a|+|x|}ax(dt)\\
 &= (xy)xt+ax(ays).
\end{aligned}
\]
The first term contains \(x^2\), and the second contains \(a^2\).
Since \(x\) and \(a\) have odd degree, both terms vanish.
\end{proof}

\begin{lemma}[The complete calculation in degrees \(23\)--\(25\)]
\label{div:lem:B2325}
The relevant graded components of \(B\) are
\[
\begin{aligned}
B^{23}&=\Q\{ar\},\\
B^{24}&=\Q\{t^2,axt,xyr,sxr\},\\
B^{25}&=0.
\end{aligned}
\]
Consequently, every element of \(B^{24}\) is a cycle.
Moreover,
\[
                             d(ar)=xyr-axt.
\]
\end{lemma}

\begin{proof}
All monomials in \(B\) have augmentation length at most three.  We sort
them by the exponents of the even generators \(t\) and \(r\).

If \(r\) occurs, then \(t\) cannot occur in degrees \(23\)--\(25\),
because \(|rt|=28\), while \(r^2\) has degree \(32\).  The odd complement
to \(r\) must have degree \(7\), \(8\), or \(9\).  Subject to the
remaining length bound of two, the corresponding lists are
\[
                   \{a\},\qquad \{sx,xy\},\qquad \varnothing.
\]
This gives \(ar\) in degree \(23\), \(sxr,xyr\) in degree \(24\),
and nothing in degree \(25\).

If \(r\) does not occur and \(t^2\) occurs, its degree is exactly
\(24\), giving the monomial \(t^2\).  If exactly one \(t\) occurs, the
odd complement must have degree \(11\), \(12\), or \(13\).  The only
odd subsets of those degrees are, respectively,
\[
                         sxy,\qquad ax,\qquad sya.
\]
The first and third would produce words of total length four and
therefore vanish in \(B\); the middle one gives \(axt\).  Finally, a
word involving only odd generators has degree at most \(18\).  These
cases prove all three displayed descriptions.

Finally,
\[
                  d(ar)=(da)r+(-1)^{|a|}a\,dr=xyr-axt.
\]
\end{proof}

\begin{proposition}
\label{div:prop:c-essential}
The class of \(c=axt\) is nonzero in \(H^{24}(M)\).
\end{proposition}

\begin{proof}
Lemma~\ref{div:lem:B2325} shows that every element of \(B^{24}\) is a cycle
and that the space of boundaries in degree \(24\) is the one-dimensional
subspace
\[
                           \Q\{xyr-axt\}.
\]
Hence
\[
 H^{24}(B)
 =\frac{\Q\{t^2,axt,xyr,sxr\}}{\Q\{xyr-axt\}}
\]
is three-dimensional.  Since \(xyr\) and \(axt\) are linearly independent,
\[
                         [axt]\neq0\in H^{24}(B).
\]
The quasi-isomorphism
\(\varphi:M\to B\) extends \(\varphi_0\) and satisfies
\(\varphi(c)=axt\).  Therefore \([c]\neq0\) in \(H^{24}(M)\).
\end{proof}

\chapter{Strict Separation and the Universal Tower Obstruction}
\label{chap:strict-separation}
\begin{chapterguide}[title={Chapter guide}]
This chapter completes the counterexample.  Acyclicity forces the unique
primitive \(xt\) into any candidate ideal; ideal closure then forces the
essential cycle \(axt\), a contradiction.  We compute both invariants exactly,
realize the example spatially, and
abstract the argument as a reusable unique-primitive dichotomy.  The final
sections apply that dichotomy uniformly to every complete strictly coordinated
level-three retractive tower.
\end{chapterguide}
\index{universal obstruction!unique primitive}
\index{homology nilpotency!strict lower bound}
\index{Sullivan realization}

\section{The strict separation}
\label{div:sec:separation}

We now prove the divergence theorem.  The obstruction begins with the
length-four element
\[
                         z=d(xt)=-xays\in\mM^4.
\]
It is a nonzero cycle in \(M^{18}\).

\begin{lemma}[The primitive is forced]
\label{div:lem:primitive-forced}
Let \(J\lhd M\) be an acyclic differential ideal satisfying
\(\mM^4\subseteq J\).  Then \(xt\in J\).
\end{lemma}

\begin{proof}
Since \(z\in\mM^4\subseteq J\) and \(dz=0\), the acyclicity of \(J\)
provides an element \(u\in J^{17}\) such that
\[
                              du=z.
\]
By \eqref{div:eq:degree17}, \(M^{17}=\Q\{xt\}\) and
\(Z^{17}(M)=0\), while \(d(xt)=z\).  Hence
\[
                         d(u-xt)=0.
\]
It follows that \(u-xt=0\), so \(u=xt\in J\).
\end{proof}

\begin{proposition}
\label{div:prop:Hnil-lower}
No acyclic differential ideal of \(M\) contains \(\mM^4\).  Consequently,
\[
                              \Hnil(M)>3.
\]
\end{proposition}

\begin{proof}
Suppose that \(J\lhd M\) were acyclic and contained \(\mM^4\).
Lemma~\ref{div:lem:primitive-forced} gives \(xt\in J\).  Since \(J\) is an
ideal and \(a\in M\), it follows that
\[
                              c=axt\in J.
\]
By Lemma~\ref{div:lem:c-cycle}, \(c\) is a cycle, and its class in \(H(J)\)
maps under \(H(J\hookrightarrow M)\) to the nonzero class of
Proposition~\ref{div:prop:c-essential}.  Hence \(H(J)\neq0\), contrary to the
assumed acyclicity.
\end{proof}

The matching upper bound is proved in
Theorem~\ref{app:thm:Hnil-exact}.  Its construction begins with
\(\mM^5\), uses the fact that every cycle of \(B\) in degree at least
\(29\) lies in \((B^+)^3\), and completes the resulting ideal degree by
degree.  Thus the obstruction above is sharp:
\[
                              \Hnil(M)=4.
\]
The detailed completion argument is kept in
Appendix~\ref{app:exact-Hnil}, where its high-degree length estimate and
inductive invariants can be checked independently of the separation
argument.

\begin{lemma}
\label{div:lem:nilh-exact}
The homotopical nil-length of \(M\) is exactly three.
\end{lemma}

\begin{proof}
The quasi-isomorphism
\[
                              \varphi:M\xrightarrow{\simeq}B
\]
and Proposition~\ref{div:prop:B-properties} give \(\nilh(M)\leq3\).

Assume that \(\nilh(M)\leq2\).  By
Lemma~\ref{div:lem:direct-zigzag}, there is an augmented quasi-isomorphism
\(f:M\to A\) with \((A^+)^3=0\).  Since
\(a,x,t\in M^+\),
\[
                         f(c)=f(a)f(x)f(t)=0.
\]
On the other hand, Proposition~\ref{div:prop:c-essential} gives
\([c]\neq0\in H(M)\), while \(H(f)\) is injective.  This is a
contradiction.  Therefore \(\nilh(M)>2\), and hence \(\nilh(M)=3\).
\end{proof}

\begin{theorem}[Divergence at length three]
\label{div:thm:algebraic-main}
There exists a simply connected minimal Sullivan algebra of finite type
\(M\) such that
\[
                         \nilh(M)=3<\Hnil(M)=4.
\]
\end{theorem}

\begin{proof}
Take \(M\) to be the minimal Sullivan model of \(B\) constructed in
Corollary~\ref{div:cor:rigid-prefix}.  It is simply connected and of finite
type by Remark~\ref{div:rem:finite-type}.  Lemma~\ref{div:lem:nilh-exact} gives
\(\nilh(M)=3\), and Theorem~\ref{app:thm:Hnil-exact} gives
\(\Hnil(M)=4\).
\end{proof}

\begin{corollary}
\label{div:cor:topological-main}
Let \(M\) be the minimal Sullivan algebra of
Theorem~\ref{div:thm:algebraic-main}.  Then its Sullivan realization \(|M|\)
is a simply connected rational space of finite rational type such that
\[
                         \Clz(|M|)=3<\Hnilz(|M|)=4.
\]
\end{corollary}

\begin{proof}
Write \(M=(\Lambda V,d)\).  Here \(H^1(M)=0\), and
\(H^*(M)\cong H^*(B)\) is finite dimensional.  Hence
\cite[Theorem~17.10]{FHT2001} shows that \(|M|\) is a simply connected
rational space and provides a canonical quasi-isomorphism
\(M\to A_{\mathrm{PL}}(|M|)\).  Thus \(M\) is the minimal Sullivan model of
\(|M|\), and
\[
 H^*(|M|;\Q)\cong H^*(M)\cong H^*(B),
 \qquad
 \pi_q(|M|)\cong (V^q)^\sharp\quad(q\geq2).
\]
Each \(V^q\), and therefore each \(\pi_q(|M|)\), is finite dimensional.
Thus \(|M|\) has finite rational type.
Cornea's
identification of nil-length with rational cone length gives
\[
                            \Clz(|M|)=\nilh(M)=3;
\]
see \cite[Theorem~1.5]{Cornea1994}.  By definition,
\(\Hnilz(|M|)=\Hnil(M)=4\).
\end{proof}

\begin{proposition}[The quadratic shadow of the example]
\label{div:prop:quadratic-shadow-infinite}
For the same minimal Sullivan algebra,
\[
 e_0(M^{(2)})=\catz(|M^{(2)}|)=\infty.
\]
Thus the quadratic-shadow theorem gives no finite upper bound for this model.
\end{proposition}

\begin{proof}
Write \(M=(\Lambda V,d)\).  The rigid-prefix result gives
\(V^{12}=\Q\{t\}\), and \(d_2t=0\) because \(dt=ays\) has word length
three.  Consider the algebra retraction
\[
 \pi:(\Lambda V,d_2)\longrightarrow(\Q[t],0),
 \qquad \pi(t)=t,
 \qquad \pi(v)=0\quad(v\ne t).
\]
To verify that \(\pi\) is a chain map, it remains only to exclude a nonzero
coefficient of \(t^2\) in \(d_2v\).  Such a coefficient necessarily forces
\(|v|=23\).  Write
\[
                         d_2v=\lambda t^2+\eta.
\]
In degree \(24\), every quadratic monomial occurring in \(\eta\) contains a
generator outside the rigid prefix and hence has image under \(\varphi\) in
\((B^+)^3\).  Indeed, such a generator maps to \((B^+)^2\), since \(B\)
has no indecomposable element above degree \(16\).  Every term of word
length at least three in \(dv\) also maps to \((B^+)^3\).  Consequently,
the word-length-two component of
\[
                         \varphi(dv)=d_B\varphi(v)
\]
is exactly \(\lambda t^2\).  On the other hand,
\[
 \varphi(v)\in B^{23}=\Q\{ar\},
 \qquad d_B\varphi(v)\in\Q\{xyr-axt\}\subseteq(B^+)^3.
\]
Comparison of the word-length-two components gives \(\lambda=0\).  Thus
\(\pi\) is a CDGA retraction of
\((\Q[t],0)\hookrightarrow M^{(2)}\).

It follows that \([t^j]\ne0\) in \(H(M^{(2)})\) for every \(j\geq1\).
For each \(n\), the nonzero class \([t^{n+1}]\) is killed by the
word-length projection \(q_n^{(2)}\), so \(e_0(M^{(2)})=\infty\).
The same polynomial family gives infinite rational cup-length and hence
\(\catz(|M^{(2)}|)=\infty\).
\end{proof}

\section{The unique-primitive dichotomy}
\label{div:sec:unique-primitive}

The preceding proof isolates a general mechanism that may be useful in
other strictification problems.

\index{primitive!unique}
\begin{proposition}[Unique-primitive dichotomy]
\label{tower:prop:primitive-or-multiple}
Let \(M\) be an augmented CDGA and let \(n\geq0\).  Suppose that there exist
homogeneous elements
\[
 p\in M^{q-1},\qquad z\in\mathfrak m_M^{n+1}\cap M^q,
 \qquad b\in M^+
\]
such that
\[
 dp=z\ne0,\qquad Z^{q-1}(M)=0,\qquad d(bp)=0,
 \qquad [bp]\ne0\in H(M).
\]
If \(J\lhd M\) is any differential ideal containing
\(\mathfrak m_M^{n+1}\), then \(H(J)\ne0\).  More precisely:
\begin{enumerate}[label=\textup{(\roman*)}]
\item if \(p\notin J\), then \([z]\ne0\) in \(H^q(J)\);
\item if \(p\in J\), then \([bp]\ne0\) in \(H^{q-1+|b|}(J)\).
\end{enumerate}
\end{proposition}

\begin{proof}
Because \(z\in\mathfrak m_M^{n+1}\subseteq J\) and \(dz=0\), the element
\(z\) defines a class in \(H^q(J)\).  Assume \(p\notin J\).  If that class
were zero, some \(u\in J^{q-1}\) would satisfy \(du=z=dp\), whence
\(u-p\in Z^{q-1}(M)=0\) and \(p=u\in J\), a contradiction.

Suppose instead that \(p\in J\).  Since \(J\) is an ideal, \(bp\in J\).
This element is a cycle, and its class maps to the nonzero class
\([bp]\in H(M)\).  Hence \([bp]\ne0\) in \(H^{q-1+|b|}(J)\).
\end{proof}

\begin{corollary}[Unique-primitive obstruction]
\label{div:crit:unique-primitive}
Suppose the hypotheses of
Proposition~\ref{tower:prop:primitive-or-multiple} hold.  Then no acyclic
differential ideal of \(M\) contains \(\mathfrak m_M^{n+1}\), and hence
\(\Hnil(M)>n\).
\end{corollary}

\begin{proof}
Every differential ideal containing \(\mathfrak m_M^{n+1}\) has nonzero
cohomology by Proposition~\ref{tower:prop:primitive-or-multiple}.
\end{proof}

\begin{remark}
The condition \(Z^{q-1}(M)=0\) in
Proposition~\ref{tower:prop:primitive-or-multiple} is deliberately stronger
than uniqueness of a primitive modulo cycles.  It makes the primitive
literally unique, rather than merely unique up to addition of an
element of \(Z^{q-1}(M)\).  The example satisfies this stronger
hypothesis exactly, by \eqref{div:eq:degree17}.
\end{remark}

For the algebra \(M\) constructed above, the proposition applies with
\[
 n=3,\qquad q=18,\qquad
 p=xt,\qquad z=-xays,\qquad b=a.
\]
The three-short witness \(B\) kills \(xt\) by adjoining the
indecomposable element \(r\) with \(dr=xt\).  This operation is possible
only after the length-four curvature \(d^2r=-xays\) has been set equal
to zero.  Within the free minimal Sullivan algebra, however, no such
indecomposable primitive is available.  An acyclic ideal is therefore
forced to contain \(xt\) itself, and multiplication by \(a\) then produces
the essential class.  This is the precise obstruction to replacing the
non-surjective short witness by a short quotient of \(M\).

\index{universal obstruction!retractive towers}

\section[Level-three tower obstruction]{Application to every complete
strictly coordinated level-three tower}

Retain the notation \(n=3\), \(q=18\), \(p=xt\), \(z=-xays\), and \(b=a\)
from the preceding section.  The required identities are recorded in
\eqref{div:eq:degree17}, Lemma~\ref{div:lem:c-cycle}, and
Proposition~\ref{div:prop:c-essential}.

\begin{corollary}[Universal tower obstruction]
\label{tower:thm:universal-obstruction}
For every complete strictly coordinated level-three retractive tower of
\(M\), let \(J_\infty\subseteq M\) be its image-ideal colimit and let
\[
 r_\infty:A_\infty\twoheadrightarrow M/J_\infty
\]
be its quotient comparison.  Then
\[
 H(J_\infty)\ne0
 \qquad\text{and}\qquad
 r_\infty\text{ is not a quasi-isomorphism}.
\]
More precisely,
\[
 xt\notin J_\infty
 \ \Longrightarrow\ \ [d(xt)]\ne0\in H^{18}(J_\infty),
\]
whereas
\[
 xt\in J_\infty
 \ \Longrightarrow\ \ [axt]\ne0\in H^{24}(J_\infty).
\]
\end{corollary}

\begin{proof}
For a complete strictly coordinated level-\(n\) tower, the image-ideal
colimit contains \((M^+)^{n+1}\); see
Proposition~\ref{def:prop:long-words}.  At level three this gives
\[
 \mathfrak m_M^4\subseteq J_\infty.
\]
These data therefore satisfy all hypotheses of
Proposition~\ref{tower:prop:primitive-or-multiple} with \(n=3\).  Its two
alternatives give the asserted class in each case and show that
\(H(J_\infty)\ne0\).  The kernel-comparison
Theorem~\ref{def:thm:kernel-comparison} now implies that \(r_\infty\) is not a
quasi-isomorphism.
\end{proof}

\begin{corollary}[Existence without strictification]
\label{tower:cor:existence-without-strictification}
The algebra \(M\) admits at least one complete strictly coordinated level-three
retractive tower, but no such tower is asymptotically Tate.
\end{corollary}

\begin{proof}
The equality \(\nilh(M)=3\), proved in
Lemma~\ref{div:lem:nilh-exact}, and the retractive-tower characterization in
Theorem~\ref{rt:thm:characterization} give the existence of at least one
complete strictly coordinated level-three tower.  For any such tower,
Corollary~\ref{tower:thm:universal-obstruction} shows that
\(J_\infty\) is not acyclic.  The finite-death criterion therefore excludes
the asymptotically Tate condition.
\end{proof}

\begin{keyidea}{Part IV summary}
The three-short CDGA \(B\) has a finite-type minimal Sullivan model with rigid
low-degree prefix.  Its unique primitive \(xt\) and essential multiple
\(axt\) force \(\Hnil(M)\geq4\), while the high-degree cubic-cycle argument
constructs an acyclic ideal containing \(\mM^5\).  Thus
\(\Hnil(M)=4\), whereas the three-short witness gives \(\nilh(M)=3\).
The same dichotomy obstructs the image ideal of every complete strictly
coordinated level-three tower, independently of all tower choices.
Part~\ref{part:stabilization}
now changes strategy.  It first proves an open-disk surjectivity theorem in
derived indecomposables and then constructs the fibre--cofibre model to which
that theorem applies.  The same hidden disk will then stabilize the separating
example after wedging with a single rational sphere.
\end{keyidea}

\part{Spherical Stabilization}\label{part:stabilization}

\chapter{The Andr\'e--Quillen and Open-Disk Toolkit}\label{chap:aq-toolkit}
\begin{chapterguide}[title={Chapter guide}]
Part~\ref{part:divergence} exhibited a rigid defect inside a fixed minimal
model; we now repair it by changing the geometry.  This chapter proves the
abstract open-disk surjectivity theorem on the mixed sector, using the Euler
contraction, homotopy transfer, and Harrison--Andr\'e--Quillen theory.
Chapter~\ref{chap:construction} constructs the required open-disk algebra from
an optimal LS root, and Chapter~\ref{chap:stabilization-proof} assembles the
three sectors to obtain stabilization.
\end{chapterguide}
\index{open-disk CDGA}
\index{Euler contraction}
\index{homotopy transfer}
\index{C-infinity algebra@$C_\infty$-algebra}
\index{Harrison complex}
\index{Andre--Quillen homology@Andr\'e--Quillen homology}
\index{derived indecomposables}
\index{marked Harrison cycle}
\index{mixed indecomposables}

\section{Operadic typing and the augmented/nonunital dictionary}
\label{stab:sec:operadic-typing}
\index{unital versus nonunital}
\index{augmentation ideal!minimal model}
\index{weak equivalence of coalgebras}

The rational-homotopy part of the memoir uses connected augmented unital
CDGAs, whereas the operad \(\Com\) in the bar--cobar theory below encodes
\emph{nonunital} commutative dg algebras.  We therefore fix the passage
between the two categories before applying any operadic construction.

Let \(\AugCDGA\) denote connected, nonnegatively graded, augmented unital
CDGAs, and let \(\ComAlgNU\) denote positively graded nonunital commutative
dg algebras.  For \(A\in\AugCDGA\), write
\[
                         \overline A=A^+=\ker(A\to\Q).
\]
For \(B\in\ComAlgNU\), its unitization is
\[
 \Q\oplus B,
 \qquad (q,b)(q',b')=(qq',qb'+q'b+bb').
\]

\begin{lemma}[Augmented/nonunital dictionary]
\label{stab:lem:aug-nonunital-dictionary}
The assignments
\[
\begin{aligned}
 \mathfrak a:\AugCDGA&\longrightarrow\ComAlgNU,
       &A&\longmapsto\overline A,\\
 \mathfrak u:\ComAlgNU&\longrightarrow\AugCDGA,
       &B&\longmapsto\Q\oplus B
\end{aligned}
\]
are inverse isomorphisms of categories.  They preserve and reflect
surjections and quasi-isomorphisms.  Transporting the standard model
structure along this isomorphism identifies cofibrations, fibrations, weak
equivalences, and cofibrant replacements in the two categories.  Moreover,
\begin{equation}\label{stab:eq:Q-nonunital}
                         Q(A)=\overline A/\overline A^{\,2}.
\end{equation}
\end{lemma}

\begin{proof}
The splitting \(A=\Q\oplus A^+\) is canonical for an augmented connected
algebra, and multiplication on \(A\) is recovered from the product on
\(A^+\) and the scalar action.  Conversely, the augmentation ideal of
\(\Q\oplus B\) is exactly \(B\).  These observations prove that the two
functors are inverse.  On underlying complexes every augmented morphism is
\(\id_\Q\oplus\overline f\); hence it is surjective, respectively a
quasi-isomorphism, exactly when \(\overline f\) is.  The model-categorical
statement is then literal transport of structure.  Formula
\eqref{stab:eq:Q-nonunital} is the definition of indecomposables after the
canonical splitting.
\end{proof}

A strictly unital, augmentation-preserving \(C_\infty\)-morphism
\(J:B\rightsquigarrow A\) means that
\[
 J_1(1)=1,\qquad J_1(B^+)\subseteq A^+,
 \qquad J_r(\ldots,1,\ldots)=0\quad(r\geq2),
\]
and that \(J_r((B^+)^{\otimes r})\subseteq A^+\).  Restriction to the
augmentation ideals is therefore a nonunital \(C_\infty\)-morphism
\(\overline J:\overline B\rightsquigarrow\overline A\); conversely,
unitization recovers \(J\).  Every operadic bar, Harrison, and cobar
construction in this chapter is applied to \(\overline A\) or
\(\overline B\).  The unit is readjoined only after applying the cobar.

The suspension symbols have the following fixed types.
\begin{center}
\small
\renewcommand{\arraystretch}{1.15}
\begin{tabularx}{0.94\textwidth}{@{}>{$}l<{$}>{\raggedright\arraybackslash}X@{}}
\toprule
\text{symbol} & \text{meaning}\\
\midrule
\bs & bar shift of cohomological degree \(-1\), so
      \(|\bs a|=|a|-1\)\\
\bs^{-1} & inverse bar shift, of degree \(+1\)\\
\ls & Lie-theoretic shift, of degree \(+1\), with
      \(\ls L_A=Z^\sharp\) for a minimal model \(\Lambda Z\to A\)\\
\tau & degree-\(+1\) shift on the square-zero suspension factor\\
\Sigma & spatial suspension\\
\bottomrule
\end{tabularx}
\end{center}
The signs attached to these shifts are collected in
Appendix~\ref{app:harrison-conventions}; their domains and codomains are fixed
here.

\section{The standard open-disk CDGA}\label{stab:sec:open-algebra}

Let \(M=(\Lambda V,d)\) be any simply connected minimal Sullivan
algebra.  Assume that \(V\) and \(U\) are positively graded and of finite
type, and let
\begin{equation}\label{stab:eq:disk-algebra}
 \Disk{U}=\Lambda(U\oplus\widehat U),
 \qquad Du=\widehat u,\quad D\widehat u=0.
\end{equation}
Put \(T=\Lambda U\) and define
\begin{equation}\label{stab:eq:open-C}
C=\Q\oplus(M\otimes\Lambda\widehat U)^+\otimes T.
\end{equation}
Equivalently, \(C^+\) is the ideal of \(M\otimes\Disk{U}\) generated by
\(M^+\) and \(\widehat U\); the pure monomials in \(T^+\) are omitted.
It is closed under multiplication.  It is also closed under the
differential: \(D(M^+)\subset M^+\), \(D\widehat U=0\), and differentiating
any \(U\)-factor introduces a factor in \(\widehat U\).  Hence \(C\) is
indeed a sub-CDGA of \(M\otimes\Disk{U}\).

\begin{lemma}[Indecomposables of the open algebra]
\label{stab:lem:QC-open}
There are natural degreewise isomorphisms
\begin{equation}\label{stab:eq:QC-open-formula}
 Q(C)\cong(V\oplus\widehat U)\otimes T,
 \qquad Q(D_C)=0.
\end{equation}
\end{lemma}

\begin{proof}
Every positive monomial of \(C\) contains a factor from \(M^+\) or
\(\widehat U\), followed by an arbitrary monomial in \(U\).  Modulo
products of two positive elements, the distinguished factor reduces to its
class in \(V=Q(M)\) or to a generator of \(\widehat U\), which proves the
first formula.  Since \(dV\subset\Lambda^{\geq2}V\) and
\(D\widehat U=0\), while \(D(U)\subset\widehat U\), the differential of
every such representative is a product of at least two positive elements
of \(C\).  Hence the induced differential on \(Q(C)\) is zero.
\end{proof}

\subsection{The suspension factor and the Euler contraction}

Define the pure open-disk algebra
\begin{equation}\label{stab:eq:C-Sigma}
 C_\Sigma=\Q\oplus(\Lambda\widehat U)^+\otimes T.
\end{equation}
For a monomial of positive total length \(N\) in
\(U\oplus\widehat U\), let
\begin{equation}\label{stab:eq:euler-homotopy}
 h=\frac{1}{N}\sum_\alpha
 u_\alpha\,\iota_{\widehat u_\alpha},
\end{equation}
where \(\iota_{\widehat u_\alpha}\) is the graded derivation contracting
one occurrence of \(\widehat u_\alpha\).  It has degree
\(-|\widehat u_\alpha|\) and is determined by
\[
 \iota_{\widehat u_\alpha}(\widehat u_\beta)=\delta_{\alpha\beta},
 \qquad \iota_{\widehat u_\alpha}(u_\beta)=0.
\]
Set \(h(1)=0\).  The graded
Euler identity gives
\begin{equation}\label{stab:eq:euler-identity}
 Dh+hD=\id-\varepsilon,
 \qquad h^2=0
\end{equation}
on \(\Disk{U}\).  In particular,
\begin{equation}\label{stab:eq:h-Domega}
 h(D\omega)=\omega,
 \qquad \omega\in T^+.
\end{equation}

Let \(\tau T^+\) denote a degree-shifted copy of \(T^+\), with the shift
chosen so that \(\tau\omega\) has the degree of \(D\omega\).  Put
\begin{equation}\label{stab:eq:H-Sigma}
 H_\Sigma=\Q\oplus \tau T^+,
 \qquad (\tau T^+)^2=0,\qquad d=0.
\end{equation}
Let \(r:\Disk{U}\to T\) be the projection killing
\(\widehat U\).

\begin{proposition}[Formality of the suspension factor]
\label{stab:prop:suspension-formality}
The formula
\begin{equation}\label{stab:eq:phi-formula}
 \varphi(1)=1,
 \qquad
 \varphi(c)=\tau\,r(hc),\quad c\in C_\Sigma^+,
\end{equation}
defines a surjective quasi-isomorphism of CDGAs
\(\varphi:C_\Sigma\twoheadrightarrow H_\Sigma\).  It satisfies
\begin{equation}\label{stab:eq:phi-Domega}
 \varphi(D\omega)=\tau\omega
 \qquad(\omega\in T^+).
\end{equation}
\end{proposition}

\begin{proof}
For \(c\in C_\Sigma^+\), one has \(r(c)=0\).  Applying \(r\) to the
Euler identity yields \(rhDc=0\), and therefore \(\varphi D=0\).  If
\(c_1,c_2\in C_\Sigma^+\), their product contains at least two
\(\widehat U\)-factors.  The homotopy \(h\) removes at most one of them, and
\(r\) then annihilates the result.  Thus
\(\varphi(c_1c_2)=0=\varphi(c_1)\varphi(c_2)\), so \(\varphi\) is a
CDGA morphism.  Equations \eqref{stab:eq:h-Domega} and
\eqref{stab:eq:phi-formula} give \eqref{stab:eq:phi-Domega}, hence surjectivity.

Finally, the short exact sequence of reduced complexes
\[
 0\longrightarrow C_\Sigma^+
 \longrightarrow\Disk{U}^+
 \longrightarrow T^+\longrightarrow0
\]
has contractible middle term.  Its connecting morphism sends
\(\omega\) to \([D\omega]\), so these classes form a basis of
\(H^+(C_\Sigma)\).  Formula \eqref{stab:eq:phi-Domega} identifies that basis
with the basis \(\tau T^+\), proving that \(\varphi\) is a
quasi-isomorphism.
\end{proof}

\subsection{A strong deformation retract}

Let
\begin{equation}\label{stab:eq:split-B}
 B_\vee=M\times_\Q C_\Sigma
   =\Q\oplus M^+\oplus C_\Sigma^+,
 \qquad M^+C_\Sigma^+=0.
\end{equation}
There is a chain decomposition
\begin{equation}\label{stab:eq:C-decomposition}
 C=M\oplus C_\Sigma^+\oplus J,
 \qquad J=M^+\otimes\Disk{U}^+.
\end{equation}
The projection
\begin{equation}\label{stab:eq:p-CB}
 p:C\twoheadrightarrow B_\vee,
 \qquad p(m+c+j)=m+c,
\end{equation}
is a CDGA morphism; \(J=\ker p\) is a differential ideal.

\begin{proposition}[Strong deformation retract]
\label{stab:prop:open-SDR}
The map \(p\) in \eqref{stab:eq:p-CB} is a surjective quasi-isomorphism.
There exist a chain inclusion \(i:B_\vee\to C\) and a degree \(-1\) homotopy
\(H:C\to C\) such that
\begin{equation}\label{stab:eq:SDR-identities}
 p i=\id_{B_\vee},\qquad
 DH+HD=\id_C-ip,\qquad
 H^2=Hi=pH=0.
\end{equation}
On the mixed ideal, \(H\) is the signed tensor homotopy
\begin{equation}\label{stab:eq:mixed-H}
 H(m\xi)=(-1)^{|m|}m\,h(\xi),
 \qquad m\in M^+,\quad\xi\in\Disk{U}^+.
\end{equation}
In particular,
\begin{equation}\label{stab:eq:H-vDomega-open}
 H(vD\omega)=(-1)^{|v|}v\omega
 \qquad(v\in V,\ \omega\in T^+).
\end{equation}
\end{proposition}

\begin{proof}
The tensor differential on \(M^+\otimes\Disk{U}^+\) is
\(d_M\otimes1+1\otimes D\).  With the Koszul sign in
\eqref{stab:eq:mixed-H}, the \(d_M\)-terms cancel in \(DH+HD\), while
\eqref{stab:eq:euler-identity} gives the identity on \(J\).  Hence \(J\) is
contractible and \(p\) is a quasi-isomorphism.  Take \(i\) to be the
linear inclusion of the two direct summands in
\eqref{stab:eq:C-decomposition}, take \(H=0\) on \(i(B_\vee)\), and use
\eqref{stab:eq:mixed-H} on \(J\).  Together with \(h^2=0\) from
\eqref{stab:eq:euler-identity}, these definitions yield all identities in
\eqref{stab:eq:SDR-identities}.
The last formula follows from \eqref{stab:eq:h-Domega}.
\end{proof}

The CDGA \(B_\vee\) models a rational wedge.  By
Proposition~\ref{stab:prop:suspension-formality}, its second factor is
quasi-isomorphic to the square-zero model of a suspension having one
reduced cohomology generator \(\tau\omega\) for every monomial
\(\omega\in T^+\).

\section{The mixed Andr\'e--Quillen edge}\label{stab:sec:aq-edge}

We now prove the central detection result: the minimal-model map reaches the
mixed indecomposable summand \(V\otimes T^+\) of the open CDGA \(C\).  The
base and disk summands are treated in Chapter~\ref{chap:stabilization-proof},
where the three pieces are assembled into full surjectivity.  All graded
duals below are degreewise duals.  Positive grading and the finite-type
assumptions ensure that every matrix and dualization is finite in each
degree.

\paragraph{Andr\'e--Quillen convention.}
For an augmented CDGA \(A\), the Andr\'e--Quillen homology used in this
chapter is
\begin{equation}\label{stab:eq:AQ-convention}
 H_*^{\mathrm{AQ}}(A/\Q;\Q)
 :=H_*\bigl(\mathbf L\Indec(\overline A)\bigr)
 =H_*\bigl(\mathbf LQ(A)\bigr),
 \qquad \overline A=A^+,
\end{equation}
where the coefficients are \(\Q\) through the augmentation.  The equality
uses the unitization equivalence of
Lemma~\ref{stab:lem:aug-nonunital-dictionary}.  Over \(\Q\), flatness is
automatic.  In the classical homological indexing, Harrison degree \(r\)
corresponds to Andr\'e--Quillen degree \(r-1\)
\cite[Section~4.2.10 and Proposition~4.2.11]{Loday1998}; the compatible
Eulerian decomposition for commutative differential graded algebras is
described in \cite[Section~5.4.8]{Loday1998}.  In our cohomological bar
convention this
shift is already incorporated in
\(Q(R_A)=\bs^{-1}\Harr(A)\), so no further suspension is inserted below.

\begin{keyidea}{Derived indecomposables and the edge map}
The ordinary indecomposable complex \(Q(A)\) is not invariant under
quasi-isomorphism.  One therefore replaces \(A\) by a cofibrant resolution
\(R_A\to A\) and forms \(Q(R_A)\); this represents the derived functor
\(\mathbf LQ(A)\), whose homology is the André--Quillen homology used here.
When \(Qd_A=0\), dualizing the map \(Q(R_A)\to Q(A)\) gives an edge
map from indecomposable coordinates of \(A\) to
\(H(Q(R_A)^\sharp)\).  Under the finite-type hypotheses, the latter
identifies with the shifted rational homotopy Lie algebra \(\ls L_A\).

Concretely, if \(f:(\Lambda Z,d)\to A\) is a minimal Sullivan model and
\(d_2:Z\to\Lambda^2Z\) is the quadratic part of its differential, then
\(L_A\) denotes the graded Lie algebra whose shifted underlying vector
space is

\[
                         \ls L_A=Z^\sharp,
\]

with bracket dual to \(d_2\) under the usual Sullivan--Quillen sign
convention.  This definition is independent of the chosen minimal model up
to graded Lie algebra isomorphism.

A \(C_\infty\)-morphism is a homotopy-coherent morphism of commutative
algebras, encoded by Taylor coefficients \(J_1,J_2,\ldots\).  Its first
coefficient is a chain map; the higher coefficients record coherent
corrections to multiplicativity.  The Harrison bar construction packages
these coefficients into a coalgebra map, and the Harrison cobar construction
turns that map into the strict cofibrant morphism used below.
\end{keyidea}
\index{derived indecomposables}
\index{Andre--Quillen homology@Andr\'e--Quillen homology}
\index{C-infinity algebra@$C_\infty$-algebra}

\subsection{Homotopy transfer and the binary mixed term}

The input is the strong deformation retract of
Proposition~\ref{stab:prop:open-SDR}; the output is a strictly unital
\(C_\infty\)-quasi-isomorphism whose binary mixed coefficient will detect
the summand \(V\otimes T^+\).

We use the following specialized form of the homotopy transfer theorem
for homotopy-commutative algebras; it follows from the rooted-tree
formulas and the unital transfer theorem of Cheng and Getzler
\cite[Theorems~6, 10, and~12]{ChengGetzler2008}.

\begin{proposition}[Transferred commutative structure]
\label{stab:prop:C-infinity-transfer}
Let \((B_\vee\mathrel{\substack{\xrightarrow{i}\\[-0.4ex]
\xleftarrow[p]{} }}C,H)\) be the strong deformation retract of
Proposition~\ref{stab:prop:open-SDR}.  Then \(B_\vee\) carries a transferred
\(C_\infty\)-structure and there is a \(C_\infty\)-quasi-isomorphism
\begin{equation}\label{stab:eq:I-C-infinity}
 I=(I_r)_{r\geq1}:B_\vee\rightsquigarrow C,
 \qquad I_1=i.
\end{equation}
The contraction preserves units and augmentations.  The morphism \(I\) does
likewise and is strictly unital.
Here \(I_r\) has cohomological degree \(1-r\).
The transferred binary operation is the original product of \(B_\vee\).  Under
the suspension and Taylor-coefficient conventions fixed in
Appendix~\ref{app:harrison-conventions},
\begin{equation}\label{stab:eq:I2-defect}
 I_2(b_1,b_2)
 =-H\bigl(i(b_1)i(b_2)-i(b_1b_2)\bigr).
\end{equation}
Consequently, for \(v\in V\), \(\omega\in T^+\), and
\(b_\omega=D\omega\in C_\Sigma^+\),
\begin{equation}\label{stab:eq:I2-mixed}
 I_2(v,b_\omega)=(-1)^{|v|+1}v\omega.
\end{equation}
\end{proposition}

\begin{proof}
The transfer formulas are finite sums over rooted trees.  Each internal
vertex is multiplication; every non-root internal edge is labelled by
\(H\); the leaves are labelled by \(i\); and the root of a transferred
operation is labelled by \(p\).  The binary operation is
\(p\mu(i\otimes i)\), which is exactly the product in the fibre product
\(B_\vee\).  Formula~\eqref{app:eq:transfer-binary-sign} gives
\eqref{stab:eq:I2-defect} with the present contraction convention.  The side
conditions for the unital transfer theorem are
\(pH=0\), \(H^2=0\), and \(H(1)=0\).  The first two are part of
\eqref{stab:eq:SDR-identities}, and the last follows from \(Hi=0\), since
\(1=i(1)\).  Thus the transferred structure and \(I\) are strictly
unital.  The maps \(p\) and \(i\) preserve augmentations.  Moreover,
\(H(C)\subset J\subset C^+\), and every tree defining \(I_r\) for
\(r\geq2\) has \(H\) at its root.  Hence
\(\varepsilon_C I_r=0\) for \(r\geq2\), while
\(\varepsilon_C I_1=\varepsilon_{B_\vee}\); therefore \(I\) is
augmentation-preserving.

Cross-products vanish in \(B_\vee\), while
\eqref{stab:eq:mixed-H} and \eqref{stab:eq:h-Domega} give
\(H(vb_\omega)=(-1)^{|v|}v\omega\).  This proves
\eqref{stab:eq:I2-mixed}.
\end{proof}

\begin{lemma}[Vanishing of higher transferred operations]
\label{stab:lem:higher-transfer-vanishing}
For the contraction of Proposition~\ref{stab:prop:open-SDR}, every
transferred operation on \(B_\vee\) of arity at least three is zero.  Hence the
transferred \(C_\infty\)-structure in
Proposition~\ref{stab:prop:C-infinity-transfer} is exactly the original
strict CDGA structure of \(B_\vee\).
\end{lemma}

\begin{proof}
A rooted tree contributing to a transferred operation of arity at least
three contains an internal edge labelled by \(H\).  Choose such an edge
closest to the leaves.  The product entering that copy of \(H\) decomposes
as an \(i(B_\vee)\)-component plus a \(J\)-component, where
\(C=i(B_\vee)\oplus J\) as complexes.  The homotopy kills \(i(B_\vee)\) and maps
\(J\) into \(J\).  Thus a surviving value above that edge lies in \(J\).
Since \(J\) is an ideal, every product farther up the tree remains in
\(J\); subsequent copies of \(H\) also preserve \(J\).  The root is
labelled by \(p\), and \(p(J)=0\), so the entire tree is zero.  The binary
transferred operation was identified with the product of \(B_\vee\) in the
proof of Proposition~\ref{stab:prop:C-infinity-transfer}; the unary
operation is its differential.  This proves the last assertion.
\end{proof}

\begin{lemma}[The \(M\)-adic marker lemma]\label{stab:lem:M-adic-marker}
Give \(M\) its multiplicative filtration
\(F^qM=(M^+)^q\).  Suppose that the inputs to a Taylor coefficient \(I_k\)
coming from the \(M\)-factor have respective filtration orders
\(q_1,\ldots,q_r\).  Then its value lies in
\begin{equation}\label{stab:eq:marker-conclusion}
 (M^+)^{q_1+\cdots+q_r}\otimes\Disk{U}.
\end{equation}
Thus a term with two \(M^+\)-markers, or one marker in \((M^+)^2\), has
zero projection to \(V\otimes T\subset Q(C)\).
\end{lemma}

\begin{proof}
Every tree defining \(I_k\) uses only multiplication and \(H\).
Multiplication adds \(M\)-adic orders, and \(H\) acts only on the disk
factor by \eqref{stab:eq:mixed-H}.  Induction from the leaves to the root
gives \eqref{stab:eq:marker-conclusion}.  In addition,
\(d_M(F^qM)\subset F^{q+1}M\) by minimality, so differential
corrections can only increase the filtration order.
\end{proof}

\subsection{The operadic Harrison coalgebra and its image model}
\label{stab:subsec:exact-Harrison-model}
\index{Harrison complex}
\index{Harrison cobar resolution}
\index{Andre--Quillen edge@Andr\'e--Quillen edge}

Let \(A\in\AugCDGA\) be simply connected.  All constructions in this
subsection are applied to the nonunital algebra \(\overline A=A^+\).  The
reduced tensor bar complex is
\[
                  \overline B(A)=T^c(\bs\overline A),
                  \qquad b=b_1+b_2,
\]
with
\begin{equation}\label{stab:eq:bar-differential-body}
 \pi_1b_1(\bs a)=-\bs(d_Aa),
 \qquad
 \pi_1b_2(\bs a\mid\bs b)=(-1)^{|a|}\bs(ab).
\end{equation}
The shuffle subcomplex is denoted by \(\mathrm{Sh}(A)\), and
\[
 \Harr^{\mathrm c}(A)
     =T^c(\bs\overline A)/\mathrm{Sh}(A),
 \qquad
 \Harr(A)=\im(e_A^{(1)}).
\]
Here \(e_A^{(1)}\) is the first Eulerian idempotent.  It commutes with the
bar differential and preserves bar weight.  In characteristic zero, Barr's
splitting gives inverse chain isomorphisms
\begin{equation}\label{stab:eq:Harrison-image-quotient-body}
\begin{gathered}
 q_A|_{\Harr(A)}:\Harr(A)\xrightarrow{\cong}\Harr^{\mathrm c}(A),
 \qquad
 s_A:\Harr^{\mathrm c}(A)\xrightarrow{\cong}\Harr(A),\\
 q_As_A=\id,
 \qquad s_Aq_A|_{\Harr(A)}=\id_{\Harr(A)}.
\end{gathered}
\end{equation}
Here the first arrow is the restriction of the quotient map.  On the full
reduced tensor bar one has \(s_Aq_A=e_A^{(1)}\).  To avoid any ambiguity, the
coalgebra structure on the image model is always the one transported through
these inverse maps; no assertion that \(e_A^{(1)}\) is a coalgebra morphism on
the full tensor bar is needed.

Let \(B_\kappa(\overline A)\) be the operadic bar construction for the
Koszul twisting morphism \(\kappa:\Comash\to\Com\).  It is a connected
conilpotent dg \(\Comash\)-coalgebra.

\begin{lemma}[Operadic bar and the Harrison quotient]\label{stab:lem:operadic-Harrison-quotient}
For every simply connected augmented CDGA \(A\), operadic d\'ecalage and
Ree's theorem give, in each weight \(r\geq1\), a natural identification
\begin{equation}\label{stab:eq:Xi-weight}
 \Comash(r)\otimes_{\Q[\mathfrak S_r]}\overline A^{\otimes r}
 \cong
 \bs^{-1}\!
 \left((\bs\overline A)^{\otimes r}/\mathrm{Sh}_r\right).
\end{equation}
Under these identifications, the internal differential corresponds to the
bar differential induced by \(d_A\), the twisting coderivation induced by
the product of \(A\) corresponds to the Harrison multiplication
differential, and the cooperad decomposition corresponds to the shifted
Harrison Lie cobracket.  Consequently the weightwise maps assemble to a
natural isomorphism of conilpotent dg \(\Comash\)-coalgebras
\begin{equation}\label{stab:eq:Xi-Harrison-operadic}
 \Xi_A:B_\kappa(\overline A)
     \xrightarrow{\cong}\bs^{-1}\Harr^{\mathrm c}(A).
\end{equation}
\end{lemma}

\begin{proof}
The Koszul duality \(\Com^!=\operatorname{Lie}\) and operadic d\'ecalage
identify the arity-\(r\) component of \(\Comash\) with the shifted dual of
\(\operatorname{Lie}(r)\); see
\cite[Sections~7.2.1--7.2.3 and Proposition~13.1.1]{LodayVallette2012}.
Ree's theorem identifies the quotient of the regular representation by the
nontrivial shuffles with \(\operatorname{Lie}(r)^\sharp\)
\cite[Theorem~1.3.6]{LodayVallette2012}.  Tensoring over
\(\Q[\mathfrak S_r]\) with \(\overline A^{\otimes r}\) yields
\eqref{stab:eq:Xi-weight}.  The operadic bar complex is
\(B_\kappa(\overline A)=\Comash\circ_\kappa\overline A\); its internal
coderivation is induced by \(d_A\), while its twisting coderivation is the
one obtained by composing the infinitesimal cooperad decomposition with the
binary product of \(A\).  Under Ree's quotient these are respectively
\(\bs^{-1}b_1\) and \(\bs^{-1}b_2\), with the signs fixed in
\eqref{stab:eq:bar-differential-body}.  This is the chain-level comparison
of operadic \(\Com\)-homology with Harrison homology
\cite[Proposition~12.1.1 and Proposition~13.1.4]{LodayVallette2012}.
Finally, the decomposition of the Koszul-dual cooperad corresponds, after
d\'ecalage, to the Harrison Lie cobracket
\cite[Corollary~13.1.5]{LodayVallette2012}.  The construction is natural in
\(A\) and respects the conilpotent weight filtrations.  This proves
\eqref{stab:eq:Xi-Harrison-operadic} as a natural isomorphism of conilpotent dg
\(\Comash\)-coalgebras.
\end{proof}
Transporting once more through \(s_A\) gives the image model used below:
\begin{equation}\label{stab:eq:Harrison-typing-diagram}
\begin{tikzcd}[column sep=large]
 B_\kappa(\overline A)
   \arrow[r,"\Xi_A","\cong"']
 & \bs^{-1}\Harr^{\mathrm c}(A)
   \arrow[r,shift left=.7ex,"\bs^{-1}s_A"]
 & \bs^{-1}\Harr(A)
   \arrow[l,shift left=.7ex,"\bs^{-1}q_A"] .
\end{tikzcd}
\end{equation}
The two arrows on the right are inverse dg \(\Comash\)-coalgebra
isomorphisms by transport.

Define first the nonunital cobar algebra and then readjoin its unit:
\begin{equation}\label{stab:eq:Harrison-resolution}
 \overline R_A
   =\Omega_\kappa\bigl(\bs^{-1}\Harr(A)\bigr),
 \qquad
 R_A=\Q\oplus\overline R_A
   =\bigl(\Lambda(\bs^{-1}\Harr(A)),d_R\bigr).
\end{equation}
Its indecomposable complex is
\begin{equation}\label{stab:eq:QR-Harrison}
 Q(R_A)=\bigl(\bs^{-1}\Harr(A),Qd_R\bigr).
\end{equation}
The counit below is the unitization of the nonunital bar--cobar counit.

\begin{proposition}[Harrison bar--cobar resolution]
\label{stab:prop:Harrison-resolution}
Let \(A\) be a simply connected augmented CDGA of finite type.  The map
\begin{equation}\label{stab:eq:bar-cobar-counit}
                 \epsilon_A^R:R_A\longrightarrow A
\end{equation}
is a quasi-isomorphism, and \(R_A\) is a quasi-free cofibrant connected
augmented CDGA.  Its indecomposable complex is
\eqref{stab:eq:QR-Harrison}.  On Harrison weight one the counit sends
\(\bs^{-1}\bs a\) to \(a\); it vanishes on the higher-weight generators
before multiplication in the cobar algebra.
\end{proposition}

\begin{proof}
By \eqref{stab:eq:Xi-Harrison-operadic}, the nonunital algebra
\(\overline R_A\) is the Koszul cobar of the operadic bar
\(B_\kappa(\overline A)\), written in the Harrison image model.  Since
\(\Com\) is Koszul, the bar--cobar counit
\(\Omega_\kappa B_\kappa(\overline A)\to\overline A\) is a
quasi-isomorphism
\cite[Corollary~11.3.5]{LodayVallette2012}.  Unitization preserves and
reflects quasi-isomorphisms by
Lemma~\ref{stab:lem:aug-nonunital-dictionary}, proving the first assertion.

The differential on the quasi-free cobar has two relevant parts.  The
internal part preserves Harrison weight and is the differential of the
generating complex.  Over \(\Q\), split that complex into homology
representatives and contractible disk pairs, placing the boundary of each
pair before its chosen preimage.  The twisting part is decomposable; on a
cogenerator it is obtained from the reduced cooperad decomposition and
therefore uses only components of strictly smaller coradical filtration.
Refining the coradical filtration by the preceding disk-pair order makes the
quasi-free algebra triangulated.  It is consequently cofibrant by
\cite[Propositions~B.6.5--B.6.6]{LodayVallette2012}; equivalently,
\cite[Theorem~12.1.6]{LodayVallette2012} recognizes this Koszul bar--cobar
construction as a cofibrant resolution computing derived indecomposables.
The same conclusion follows independently from
\cite[Theorem~4.2.4]{Fresse2009}: Koszulity supplies the required operadic
cobar weak equivalence to \(\Com\), and over \(\Q\) the operad and its
Koszul-dual cooperad are \(\Sigma_*\)-cofibrant, the cooperad is connected,
and the underlying complex of \(\overline A\) is cofibrant.

Finally, the indecomposables of the free commutative algebra underlying
\(R_A\) are its generators, and the induced linear differential is
\(Qd_R\).  The stated weight-one formula is the normalization of the
universal twisting cochain.
\end{proof}

\subsection{\texorpdfstring{Functoriality for \(C_\infty\)-morphisms}{Functoriality for C-infinity morphisms}}

Let \(J:B\rightsquigarrow A\) be strictly unital and
augmentation-preserving.  Its restriction
\(\overline J:\overline B\rightsquigarrow\overline A\) is nonunital.  There
are three distinct coalgebra maps in the construction:
\begin{enumerate}[label=\textup{(\arabic*)},leftmargin=2.5em]
\item the full tensor-bar map \(\overline B(J)\), whose corestrictions are
      the shifted Taylor coefficients;
\item the induced quotient map
      \(\Harr^{\mathrm c}(J):\Harr^{\mathrm c}(B)\to
      \Harr^{\mathrm c}(A)\), obtained because the \(C_\infty\) Taylor
      coefficients annihilate shuffle decomposables;
\item the map on the Harrison image models obtained by conjugating with the
      inverse maps in \eqref{stab:eq:Harrison-image-quotient-body}.
\end{enumerate}
Equivalently, through \(\Xi_A\) and \(\Xi_B\), the second map is the
operadic bar map \(B_\kappa(\overline J)\).  In the notation for
\(\infty\)-morphisms this map is often written \(B_\iota(\overline J)\):
for strict \(\Com\)-algebras, regarded as \(C_\infty\)-algebras with no
higher operations, the bar objects \(B_\iota\) and \(B_\kappa\) have the
same cofree \(\Comash\)-coalgebra and the same coderivation.  This equivalence
guarantees that
the quotient and transported maps preserve the complete dg
\(\Comash\)-coalgebra structure, not merely the underlying complexes.

\begin{proposition}[Functoriality for \(C_\infty\)-morphisms]
\label{stab:prop:Harrison-Cinfinity-functoriality}
Let \(A\) and \(B\) be simply connected augmented CDGAs of finite type.
An augmentation-preserving strictly unital \(C_\infty\)-morphism
\(J=(J_r):B\rightsquigarrow A\) induces the dg
\(\Comash\)-coalgebra morphism
\begin{equation}\label{stab:eq:Harrison-transported-map}
 \widetilde B(J)
  =(\bs^{-1}s_A)\circ
    (\bs^{-1}\Harr^{\mathrm c}(J))\circ
    (\bs^{-1}q_B):
 \bs^{-1}\Harr(B)\longrightarrow\bs^{-1}\Harr(A),
\end{equation}
and hence a strict augmented CDGA morphism
\begin{equation}\label{stab:eq:FJ-strict}
 F_J:=\epsilon_A^R\circ
 \bigl(\Q\oplus\Omega_\kappa\widetilde B(J)\bigr):
 R_B\longrightarrow R_A\longrightarrow A.
\end{equation}
On \(Q(R_B)=\bs^{-1}\Harr(B)\), the Harrison-weight-\(r\)
corestriction is the bar-shifted form of \(J_r\), with the sign fixed in
Appendix~\ref{app:harrison-conventions}.

If \(J_1\) is a quasi-isomorphism, then \(F_J\) is a quasi-isomorphism.
More precisely, \(B_\kappa(\overline J)\) is a \emph{weak equivalence} of
conilpotent dg \(\Comash\)-coalgebras, meaning that its image under
\(\Omega_\kappa\) is a quasi-isomorphism.  No claim is made that the
underlying bar-coalgebra map is itself a quasi-isomorphism of complexes.
\end{proposition}

\begin{proof}
An \(\infty\)-morphism is a morphism of the associated quasi-cofree
\(\Comash\)-coalgebras; its corestriction to cogenerators determines it.
This is the content of
\cite[Section~10.2.2, Proposition~10.2.1, and Theorem~10.2.3]{LodayVallette2012}.
For \(\Com\), the shuffle condition on the Taylor coefficients is
\cite[Proposition~13.1.6]{LodayVallette2012}, and the bar functor on
\(C_\infty\)-algebras is described in
\cite[Proposition~11.4.1]{LodayVallette2012}.  Conjugation by the coalgebra
isomorphisms in
\eqref{stab:eq:Harrison-typing-diagram} proves that
\eqref{stab:eq:Harrison-transported-map} is a dg
\(\Comash\)-coalgebra morphism.  Its weight-\(r\) corestriction is the
shifted Taylor coefficient because the projection to weight one commutes
with both transports.  Applying the cobar and the counit gives
\eqref{stab:eq:FJ-strict}.

If \(J_1\) is a quasi-isomorphism, then
\cite[Proposition~11.4.7]{LodayVallette2012} states that the induced
operadic bar map is a weak equivalence in exactly the sense recorded in the
statement: its cobar is a quasi-isomorphism.  This functorial passage and the
resulting weak equivalence may also be checked through
\cite[Propositions~4.2.7--4.2.8]{Fresse2009}.  The transports in
\eqref{stab:eq:Harrison-typing-diagram} are isomorphisms, so
\(\Q\oplus\Omega_\kappa\widetilde B(J)\) is a quasi-isomorphism.  The
counit \(\epsilon_A^R\) is a quasi-isomorphism by
Proposition~\ref{stab:prop:Harrison-resolution}; hence their composite
\(F_J\) is a quasi-isomorphism.
\end{proof}

\subsection{Derived indecomposables}
\index{indecomposables!derived}

Let \(\varepsilon:\Com\to\mathcal I\) be the augmentation of the operad.
On nonunital algebras it gives the adjunction
\begin{equation}\label{stab:eq:Indec-adjunction}
 \Indec=\varepsilon_!:
 \mathbf{dg\,Com\text{-}alg}
 \rightleftarrows
 \mathbf{dgMod}_{\Q}:\varepsilon^*,
 \qquad
 \varepsilon_!(A)=\mathcal I\circ_{\Com}A=A/A^2,
\end{equation}
where \(\varepsilon^*(W)\) has zero product.  Under unitization this is the
adjunction
\[
 Q:\AugCDGA\rightleftarrows\dgMod:
       \bigl(W\longmapsto\Q\oplus W,\text{ with }W^2=0\bigr).
\]
The right adjoint preserves surjective fibrations and surjective
quasi-isomorphisms because these are detected on the underlying complexes.
Thus \(Q\) is left Quillen in the transported model structure: the
extension--restriction adjunction is
\cite[Lemma~12.1.5]{LodayVallette2012}, the model structure is
\cite[Proposition~B.6.5]{LodayVallette2012}, and the resulting total left
derived functor is the derived indecomposables functor \(\mathbf LQ\).
For the Koszul operad \(\Com\), its homology is computed by the operadic,
hence Harrison, complex by \cite[Theorem~12.1.6]{LodayVallette2012}.  The fixed-operad algebra model
structure used here is also treated in \cite{Hinich1997}; no model structure
on the category of operads is used in this argument.  Mill\`es develops the
parallel cotangent-complex and cohomological form of Andr\'e--Quillen theory
for algebras over operads \cite[Theorems~1.3.6 and~1.4.2 and
Section~3.1]{Milles2011}; the present proof uses the homological left-derived
indecomposables formulation.

\begin{proposition}[Derived indecomposables]
\label{stab:prop:Harrison-derived-indecomposables}
Let \(A\) be a simply connected augmented CDGA of finite type, and let
\(f:N=(\Lambda Z,d_N)\to A\) be a minimal Sullivan model.  There are
comparison isomorphisms of graded vector spaces
\begin{equation}\label{stab:eq:derived-indec-comparison}
                         H Q(R_A)\cong H Q(N)=Z
\end{equation}
and
\begin{equation}\label{stab:eq:Harrison-edge-identification}
 H\bigl(Q(R_A)^\sharp\bigr)\cong Z^\sharp=\ls L_A.
\end{equation}
They are canonical at the level of the derived functor \(\mathbf LQ(A)\),
that is, in the homotopy category of cofibrant resolutions.  Once the
resolution maps \(R_A\to A\) and \(N\to A\) are fixed, the comparison
used below is the one constructed in
Theorem~\ref{stab:thm:Harrison-AQ-comparison}.  No strict naturality with
respect to an unfixed, nonfunctorial choice of minimal Sullivan model, and no
compatibility with a Lie bracket, is asserted.
\end{proposition}

\begin{proof}
The algebras \(R_A\) and \(N\) are cofibrant and both map by
quasi-isomorphisms to \(A\).  Therefore \(Q(R_A)\) and \(Q(N)\) represent
the same value of \(\mathbf LQ(A)\).  Their comparison is thus canonical in
the homotopy category of cofibrant resolutions, which proves the first
isomorphism in the stated sense.  Minimality gives \(Qd_N=0\), so
\(HQ(N)=Z\).  In every total degree the
Harrison complex is finite dimensional: positive grading bounds the possible
bar weight, and finite type bounds every tensor factor.  Degreewise duality
over \(\Q\) is exact, so dualizing the first isomorphism gives
\(H(Q(R_A)^\sharp)\cong Z^\sharp\), with no completed dual.  The last
identity is the underlying graded-vector-space part of Sullivan--Quillen
duality \cite[Section~21(e) and Theorem~21.6]{FHT2001}.
\end{proof}

The following lemma replaces any appeal to the assertion that a left Quillen
functor must send a chosen cylinder to a chosen cylinder.

\begin{lemma}[Indecomposables of a polynomial homotopy]
\label{stab:lem:Q-polynomial-homotopy}
Let \(R\) and \(A\) be connected augmented CDGAs, and let
\(u_0,u_1:R\to A\) be augmented morphisms.  Suppose that a polynomial
Sullivan homotopy
\[
                  \mathcal H:R\longrightarrow A[t,dt]
\]
joins them and has constant scalar part; equivalently, \(\mathcal H\) factors
through the augmented path object
\begin{equation}\label{stab:eq:augmented-path-object}
 P_\varepsilon(A)
 =A[t,dt]\times_{\Q[t,dt]}\Q.
\end{equation}
Then \(Qu_0\) and \(Qu_1\) are chain-homotopic.
\end{lemma}

\begin{proof}
Let \(K=K_1\mathcal H:R\to A\) be polynomial integration as in
Lemma~\ref{rt:lem:polynomial-integration}.  It satisfies
\begin{equation}\label{stab:eq:polynomial-chain-homotopy}
                     u_1-u_0=d_AK+Kd_R.
\end{equation}
We claim that \(K((R^+)^2)\subseteq(A^+)^2\).  For homogeneous
\(x,y\in R^+\), write
\[
 \mathcal H(x)=a_x(t)+b_x(t)dt,
 \qquad
 \mathcal H(y)=a_y(t)+b_y(t)dt.
\]
Because \(\mathcal H\) factors through \(P_\varepsilon(A)\), every
coefficient of the four polynomials \(a_x,b_x,a_y,b_y\) belongs to \(A^+\).
The \(dt\)-coefficient of the product
\(\mathcal H(xy)=\mathcal H(x)\mathcal H(y)\) is a signed sum of products of
one coefficient from the \(x\)-expression and one from the \(y\)-expression.
It therefore lies in \((A^+)^2[t]\), and its integral belongs to
\((A^+)^2\).  Linearity proves the claim.

Consequently, \(K\) induces a degree \(-1\) map
\(\overline K:Q(R)\to Q(A)\).  Projecting
\eqref{stab:eq:polynomial-chain-homotopy} to indecomposables gives
\begin{equation}\label{stab:eq:Q-homotopy-direct}
 Qu_1-Qu_0=Qd_A\,\overline K+\overline K\,Qd_R.
\end{equation}
\end{proof}

The object \(P_\varepsilon(A)\) is the standard path object in
\(\AugCDGA\).  Explicitly,
\[
 P_\varepsilon(A)=\Q\oplus A^+[t,dt].
\]
The constant inclusion \(A\to P_\varepsilon(A)\) is a quasi-isomorphism: the
usual polynomial contraction of \(\Q[t,dt]\) restricts to a contraction of
the positive summand onto \(A^+\).  The endpoint map
\[
 P_\varepsilon(A)\longrightarrow A\times_\Q A,
 \qquad u\longmapsto(\operatorname{ev}_0u,\operatorname{ev}_1u),
\]
is surjective, since a pair \((a_0,a_1)\) with the same augmentation is lifted
by \((1-t)a_0+ta_1\).  Thus two maps from a cofibrant augmented CDGA that
represent the same morphism in the homotopy category may be joined by a
polynomial homotopy of the form required in the lemma.  The proof above
records its effect on indecomposables without assuming that a chosen cylinder
is preserved by \(Q\).

\begin{theorem}[Harrison--Andr\'e--Quillen edge comparison]
\label{stab:thm:Harrison-AQ-comparison}
Let \(A\) be a simply connected augmented CDGA of finite type with
\(Q(d_A)=0\), and let \(f:N=(\Lambda Z,d_N)\to A\) be a minimal Sullivan
model.  If \(F:R\to A\) is any quasi-isomorphism from a cofibrant simply
connected augmented CDGA whose indecomposable complex is degreewise finite,
there is a comparison isomorphism, canonical once \(f\) and \(F\) are fixed
and natural with respect to morphisms of cofibrant resolutions,
\begin{equation}\label{stab:eq:edge-comparison-map}
 \Theta_F:H Q(R)\xrightarrow{\cong}Q(N)
 \quad\text{such that}\quad
 Qf\,\Theta_F=H(QF).
\end{equation}
For every \(\lambda\in Q(A)^\sharp\), the cochain
\(\lambda QF\in Q(R)^\sharp\) is a cocycle and
\begin{equation}\label{stab:eq:edge-comparison-formula}
 \Theta_F^\sharp\bigl((Qf)^\sharp\lambda\bigr)
   =[\lambda QF]\in H\bigl(Q(R)^\sharp\bigr).
\end{equation}
In particular, for \(F=\epsilon_A^R\), the edge cocycle is
\begin{equation}\label{stab:eq:edge-cocycle}
 \lambda\,Q\epsilon_A^R\in Q(R_A)^\sharp,
 \qquad \lambda\in Q(A)^\sharp,
\end{equation}
and its class corresponds to
\((Qf)^\sharp(\lambda)\in Q(N)^\sharp=\ls L_A\) under
\eqref{stab:eq:Harrison-edge-identification}.  The same statement applies to
\(F_J:R_B\to A\) from
Proposition~\ref{stab:prop:Harrison-Cinfinity-functoriality} whenever
\(J_1\) is a quasi-isomorphism.
\end{theorem}

\begin{proof}
Factor the weak equivalence \(f\) in \(\AugCDGA\) as
\begin{equation}\label{stab:eq:AQ-model-factorization}
 N\xrightarrow[\sim]{\ i\ }E
  \xrightarrow[\sim]{\ p\ }A,
 \qquad f=pi,
\end{equation}
with \(i\) an acyclic cofibration and \(p\) an acyclic fibration.  Since
\(R\) is cofibrant, lift \(F\) to \(\widetilde F:R\to E\) with
\(p\widetilde F=F\).  The lifting property of \(i\) against the augmentation
\(N\twoheadrightarrow\Q\) gives a retraction \(r:E\to N\) with
\(ri=\id_N\).  Set
\begin{equation}\label{stab:eq:AQ-resolution-comparison}
                         g=r\widetilde F:R\longrightarrow N.
\end{equation}
Both \(\widetilde F\) and \(r\) are weak equivalences, so \(g\) is a weak
equivalence between cofibrant objects.  Hence \(Qg\) is a
quasi-isomorphism by the Quillen adjunction
\eqref{stab:eq:Indec-adjunction}.

In the homotopy category, \([g]=[f]^{-1}[F]\), so \(fg\) and \(F\) are
represented by the same morphism.  Since \(R\) is cofibrant and \(A\) is
fibrant, the augmented path object \(P_\varepsilon(A)\) of
\eqref{stab:eq:augmented-path-object} supplies a polynomial Sullivan homotopy
between them.  Lemma~\ref{stab:lem:Q-polynomial-homotopy} gives a chain
homotopy between \(Q(fg)\) and \(QF\).  Define
\[
 \Theta_F=H(Qg):H Q(R)\xrightarrow{\cong}H Q(N)=Q(N).
\]
Passing to homology in that direct chain homotopy yields
\(Qf\,\Theta_F=H(QF)\).

Any two choices of \(g\) represent the same morphism
\([f]^{-1}[F]\).  They are polynomial-homotopic, and the same lemma shows
that they induce the same map on \(HQ\).  This proves canonicity.  If
\(u:R'\to R\) satisfies \(Fu=F'\), the homotopy-class identity gives
\(\Theta_{F'}=\Theta_FH(Qu)\), proving naturality.

Since \(Qd_A=0\), the composite \(\lambda QF\) is a cocycle.  All
complexes in the displayed degree are finite dimensional, so dualizing the
identity \(Qf\Theta_F=H(QF)\) degree by degree gives
\eqref{stab:eq:edge-comparison-formula}.  No completed dual occurs.  The
final identification with \(\ls L_A\) is
Proposition~\ref{stab:prop:Harrison-derived-indecomposables}.
\end{proof}

The binary Harrison coefficient used below is fixed once and for all in
formula~\eqref{app:eq:binary-unit-coefficient}.

The transferred morphism \(I:B_\vee\rightsquigarrow C\) is a strictly unital,
augmentation-preserving \(C_\infty\)-quasi-isomorphism because \(I_1=i\)
is a quasi-isomorphism.  Its bar--cobar representative
\begin{equation}\label{stab:eq:bar-cobar-specialized}
 F=F_I:R_{B_\vee}=\Omega_\kappa\bigl(\bs^{-1}\Harr(B_\vee)\bigr)\longrightarrow C
\end{equation}
is therefore a quasi-isomorphism.  By
Proposition~\ref{stab:prop:Harrison-resolution}, its source is cofibrant, and
\eqref{stab:eq:bar-cobar-specialized} is a cofibrant resolution of \(C\).
Consequently, for every minimal Sullivan model \(f:N\to C\) and every
\(\lambda\in Q(C)^\sharp\), the class represented by \(\lambda QF\)
corresponds to \((Qf)^\sharp(\lambda)\)
under \eqref{stab:eq:Harrison-edge-identification}.

Fix a homogeneous basis of \(V\).  For a basis element \(v\), let
\(x_v\in L_M\) be normalized by \(\ls x_v=v^\sharp\).  No
free-product description of \(L_C\) is needed below: independence is
proved directly by the diagonal evaluation in
Theorem~\ref{stab:thm:mixed-AQ-edge}.

\subsection{Marked Harrison cycles}

We now supplement the preceding comparison with filtered cycles carrying
one minimal-model marker and one disk marker.  Their diagonal evaluations
force surjectivity onto the mixed indecomposables.

For a bar tensor, define its total \(M\)-adic order by adding the
\(M\)-adic orders of all entries from \(M^+\).  More precisely, let
\begin{equation}\label{stab:eq:Harrison-M-filtration}
 F_M^p\Harr(M)=\Harr(M)\cap
 \left\langle \bs a_1\mid\cdots\mid \bs a_r\ \middle|\
 a_i\in(M^+)^{q_i},\ q_i\geq1,\ \sum_iq_i\geq p\right\rangle .
\end{equation}

\begin{lemma}[Filtered lift of a minimal generator]
\label{stab:lem:filtered-Harrison-lift}
For every homogeneous \(v\in V\), there is a Harrison cycle
\begin{equation}\label{stab:eq:gamma-v-shape}
 \gamma_v=\bs v+\gamma'_v,
 \qquad \gamma'_v\in F_M^2\Harr(M),
\end{equation}
whose desuspension represents the class in \(H Q(R_M)\) corresponding to
\(v\in V\) under \eqref{stab:eq:derived-indec-comparison}.
\end{lemma}

\begin{proof}
Minimality gives
\(d_M((M^+)^q)\subset(M^+)^{q+1}\), while a bar multiplication replaces
two adjacent entries of orders \(q_i,q_{i+1}\) by one entry of order at
least \(q_i+q_{i+1}\).  Hence the Harrison differential preserves
\eqref{stab:eq:Harrison-M-filtration}; its internal part strictly raises the
order on a linear entry.  In particular,
\[
 \operatorname{gr}_M^1Q(R_M)=\bs^{-1}\bs V\cong V.
\]
The indecomposable counit
\begin{equation}\label{stab:eq:Q-counit-M}
 Q\epsilon_M^R:Q(R_M)=\bs^{-1}\Harr(M)\longrightarrow Q(M)=V
\end{equation}
is the identity on that associated-graded term.  By the left-Quillen argument
in the proof of
Proposition~\ref{stab:prop:Harrison-derived-indecomposables},
\(Q\epsilon_M^R\) is a quasi-isomorphism.  It is also surjective: the counit sends the linear
Harrison generator \(\bs^{-1}\bs v\) exactly to \(v\).  Therefore its
kernel is acyclic.  To make the cycle-level choice explicit, start with
\(x\in Q(R_M)\) such that \(Q\epsilon_M^R(x)=v\).  Minimality gives
\(Qd_M(v)=0\), so \(dx\) belongs to the acyclic kernel.  Choose
\(y\in\ker Q\epsilon_M^R\) with \(dy=dx\).  Then
\(\bs^{-1}\widetilde\gamma_v=x-y\) is a cycle sent exactly to \(v\).
To identify the filtration of this cycle, note that \(Q\epsilon_M^R\) sends
a weight-one generator \(\bs^{-1}\bs a\) to the class of \(a\) in \(Q(M)\),
and it vanishes on
Harrison weights at least two.  The Eulerian idempotent preserves bar
weight and is the identity in weight one, so this description applies to
the image model \(\Harr(M)\).  Its kernel therefore consists of the
weight-at-least-two terms together with the weight-one generators
\(\bs^{-1}\bs a\) for \(a\in(M^+)^2\).  Hence
\[
 \ker Q\epsilon_M^R=\bs^{-1}F_M^2\Harr(M),
\]
and both \(\bs^{-1}\widetilde\gamma_v\) and
\(\bs^{-1}\bs v\) are sent to \(v\).  Therefore
\[
 \bs^{-1}(\widetilde\gamma_v-\bs v)
 \in\ker Q\epsilon_M^R=\bs^{-1}F_M^2\Harr(M).
\]
Thus
\(\widetilde\gamma_v=\bs v+\gamma'_v\), with
\(\gamma'_v\in F_M^2\Harr(M)\).  Taking
\(\gamma_v=\widetilde\gamma_v\) proves the claim.  Its class corresponds to
\(v\), hence pairs with \(\ls x_v=v^\sharp\).
\end{proof}

For \(\omega\in T^+\), set \(b_\omega=D\omega\in C_\Sigma^+\).  If
\(\xi\in\Harr(M)\), define the marked insertion
\begin{equation}\label{stab:eq:marked-insertion}
 \jmath_{b_\omega}(\xi)
 =e^{(1)}(\bs b_\omega\mid\xi)\in\Harr(B_\vee),
 \qquad
 \gamma_{v,\omega}=\jmath_{b_\omega}(\gamma_v).
\end{equation}

Extend \(F_M^p\) to \(\Harr(B_\vee)\) by counting only entries from the
\(M\)-factor.

\begin{lemma}[Marked Harrison cycle]\label{stab:lem:marked-Harrison-cycle}
The element \(\gamma_{v,\omega}\) is a Harrison cycle.  It contains
exactly one entry from the suspension factor, namely \(b_\omega\).  Its
initial term is \(e^{(1)}(\bs b_\omega\mid \bs v)\), and every other term
has \(M\)-adic order at least two.
\end{lemma}

\begin{proof}
Put \(L_{b_\omega}(x)=\bs b_\omega\mid x\) on the reduced tensor bar of the
\(M\)-summand.  This operator has degree
\(|\bs b_\omega|=|b_\omega|-1\).  For a homogeneous tensor
\(x=\bs a_1\mid\cdots\mid\bs a_r\), with every \(a_i\in M^+\), the
coderivation formula and \eqref{stab:eq:bar-differential-body} give
\begin{align}
 b_1L_{b_\omega}(x)
   &=(-1)^{|b_\omega|-1}L_{b_\omega}b_1(x),
   \label{stab:eq:b1-insertion-commutator}\\
 b_2L_{b_\omega}(x)
   &=(-1)^{|b_\omega|-1}L_{b_\omega}b_2(x)
     +(-1)^{|b_\omega|}
       \bs(b_\omega a_1)\mid\bs a_2\mid\cdots\mid\bs a_r.
   \label{stab:eq:b2-insertion-commutator}
\end{align}
For the first identity, the possible extra term is
\(-\bs(Db_\omega)\mid x\), which vanishes because \(Db_\omega=D^2\omega=0\).
For the second, the displayed last term is the unique multiplication crossing
the inserted bar entry; every multiplication internal to \(x\) acquires the
common Koszul factor \((-1)^{|\bs b_\omega|}\).  The crossing term vanishes
because
\[
             b_\omega a_1=0
             \quad\text{in}\quad
             B_\vee=M\times_\Q C_\Sigma,
             \qquad M^+C_\Sigma^+=0.
\]
Equivalently, with the graded commutator
\([P,Q]=PQ-(-1)^{|P||Q|}QP\), equations
\eqref{stab:eq:b1-insertion-commutator} and
\eqref{stab:eq:b2-insertion-commutator} say
\begin{equation}\label{stab:eq:marked-commutators}
 [b_1,L_{b_\omega}]=0,
 \qquad
 [b_2,L_{b_\omega}]=0.
\end{equation}
on tensors from \(M^+\).

The first Eulerian idempotent commutes with \(b=b_1+b_2\) and preserves bar
weight.  Since \(\gamma_v\) is a Harrison cycle,
\eqref{stab:eq:marked-commutators} therefore implies
\[
 b\,e^{(1)}L_{b_\omega}(\gamma_v)
 =(-1)^{|b_\omega|-1}e^{(1)}L_{b_\omega}(b\gamma_v)=0.
\]
Thus \(\gamma_{v,\omega}\) is a Harrison cycle.  The Eulerian idempotent only
permutes and recombines the existing bar entries, so every summand still has
exactly one suspension-factor entry.  Finally,
\(\gamma_v=\bs v+\gamma'_v\), with
\(\gamma'_v\in F_M^2\Harr(M)\); applying the weight-preserving operator
\(e^{(1)}L_{b_\omega}\) gives the asserted initial term and filtration bound.
\end{proof}

\begin{theorem}[Mixed Andr\'e--Quillen edge]
\label{stab:thm:mixed-AQ-edge}
Let \(f:N\to C\) be a minimal Sullivan model of the open CDGA
\eqref{stab:eq:open-C}.  Then the composite
\begin{equation}\label{stab:eq:mixed-Q-surjection}
 Q(N)\xrightarrow{Qf}Q(C)
 \longrightarrow V\otimes T^+
\end{equation}
is surjective.
\end{theorem}

\begin{proof}
Apply Proposition~\ref{stab:prop:C-infinity-transfer} to obtain
\(I:B_\vee\rightsquigarrow C\).  Let
\[
 F:\Omega_\kappa(\bs^{-1}\Harr(B_\vee))\longrightarrow C
\]
be the corresponding bar--cobar morphism.  For a basis element
\(v\omega\in V\otimes T^+\), let
\(\lambda_{v,\omega}\in(Q(C))^\sharp\) be the coordinate functional that
takes the value \(1\) on \(v\omega\) and \(0\) on every other basis element
of \(Q(C)\).

Evaluate the cocycle \(\lambda_{u,\eta}QF\) on the cycles of
Lemma~\ref{stab:lem:marked-Harrison-cycle}.  The initial term is governed by
the binary Taylor coefficient.  Equation~\eqref{stab:eq:I2-mixed} and
formula~\eqref{app:eq:binary-unit-coefficient} show that its Harrison
coefficient is a unit, and give
\begin{equation}\label{stab:eq:diagonal-evaluation}
 (\lambda_{u,\eta}QF)(\bs^{-1}\gamma_{v,\omega})
 =c_{v,\omega}\delta_{uv}\delta_{\eta\omega},
 \qquad
 c_{v,\omega}=(-1)^{|\omega|(|v|+1)+1}\in\{+1,-1\}.
\end{equation}
With the conventions of Appendix~\ref{app:harrison-conventions},
formula~\eqref{app:eq:marked-binary-unit} gives the displayed value of
\(c_{v,\omega}\).  Every remaining term arises from \(\gamma'_v\) and has
total \(M\)-adic order at least two.  The marker lemma therefore places its
image in \((M^+)^2\otimes\Disk{U}\), where every \(\lambda_{u,\eta}\)
vanishes.

The pairing in \eqref{stab:eq:diagonal-evaluation} depends only on the
cohomology classes of the cocycles and the homology classes of the cycles.
It follows that the classes \([\lambda_{v,\omega}QF]\) are linearly
independent: evaluating a finite relation
\(\sum_{u,\eta}a_{u,\eta}[\lambda_{u,\eta}QF]=0\) on
\(\bs^{-1}\gamma_{v,\omega}\) gives
\(a_{v,\omega}c_{v,\omega}=0\).  By
Theorem~\ref{stab:thm:Harrison-AQ-comparison}, these classes
correspond to \((Qf)^\sharp(\lambda_{v,\omega})\); hence the latter are
linearly independent.  Thus \((Qf)^\sharp\) is injective on
\((V\otimes T^+)^\sharp\).  In each total degree, only finitely many monomials
of \(T=\Lambda U\) occur, and the corresponding pieces of \(V\) and \(U\) are
finite dimensional.  Therefore the relevant component of \(V\otimes T^+\),
and the Harrison components paired with it, are finite dimensional.  Ordinary
degreewise duality---not a completed duality---now says exactly that
\eqref{stab:eq:mixed-Q-surjection} is surjective.
\end{proof}
\paragraph{Dictionary for the geometric construction.}
Chapter~\ref{chap:construction} realizes the abstract symbols of this
chapter as follows.
\begin{center}
\small
\begin{tabularx}{\linewidth}{@{}l>{\raggedright\arraybackslash}X@{}}
\toprule
Symbol & Meaning in Chapter~\ref{chap:construction}\\
\midrule
\(M=(\Lambda V,d)\) & Minimal Sullivan model of \(X\).\\
\(D\) & Differential on \(E_\sigma\) and on its exact relative-disk
normal form; \(\Disk{U}\) denotes the disk algebra itself.\\
\(T=\Lambda U\) & Pure fibre-coordinate algebra of the relative model.\\
\(C\) & Standardized open cofibre model corresponding to \(C_\sigma\).\\
\(V=Q(M)\) & Minimal generator space of the base.\\
\(U\) & Fibre-minimal generators in \((\Lambda U,0)\).\\
\(\widehat U\) & Disk partners \(D(U)\), replacing the relative generators \(Y\).\\
\bottomrule
\end{tabularx}
\end{center}

\begin{keyidea}{Why André--Quillen theory enters}
Ordinary cohomology detects the quasi-isomorphism between the open algebra
and the wedge model.  It does not control which indecomposable coordinates
a minimal-model map reaches.  Derived indecomposables retain exactly that
information.  Marked Harrison cycles make the decisive coordinates visible
and prevent their cancellation.
\end{keyidea}

\chapter{From an Optimal LS Root to an Open-Disk Model}\label{chap:construction}
\begin{chapterguide}[title={Chapter guide}]
We now return to the optimal \(n\)-LS root of Chapter~\ref{stab:sec:geom-prelim}.
We factor the chosen section by a relative Sullivan extension, model its
cofibre by a strict pullback, and descend the acyclic kernel to an
\(n\)-short model.  A triangular change of variables then puts the cofibre
model in the exact open-disk normal form studied abstractly in
Chapter~\ref{chap:aq-toolkit}.
\end{chapterguide}
\index{optimal LS root}
\index{relative Sullivan algebra}
\index{cofibre model}
\index{disk normal form}
\index{semifree module}

\section{The short model of the fibre--cofibre construction}
\label{stab:sec:short-model}

We retain the notation of Chapter~\ref{stab:sec:geom-prelim}; in particular,
\(n=\catz(X)\geq1\), \(M=(\Lambda V,d)\), and
\[
 M\xrightarrow{i_n}P_n\xrightarrow[\simeq]{\kappa_n}A_n
 \qquad\text{with}\qquad
 A_n=M/(M^+)^{n+1}.
\]
The strict retraction \(\rho_n:P_n\to M\) is a CDGA model of the
section \(\sigma:X\to R_n\).

\subsection{A relative model of the section}

Throughout this subsection, \(E_\sigma\) denotes the total relative Sullivan
algebra and \(F_\sigma^{\mathrm{alg}}\) its fibre algebra.

\begin{lemma}[Fibre-minimal model of the section]
\label{stab:lem:fibre-minimal-section}
The retraction \(\rho_n:P_n\to M\) has a factorization
\begin{equation}\label{stab:eq:section-relative-factor}
 P_n\lhook\joinrel\longrightarrow
 E_\sigma=(P_n\otimes\Lambda U,D)
 \xrightarrow[\simeq]{\lambda_\sigma}\!\!\!\twoheadrightarrow M,
 \qquad \lambda_\sigma|_{P_n}=\rho_n,
\end{equation}
whose fibre algebra is the base change
\begin{equation}\label{stab:eq:fibre-algebra}
 F_\sigma^{\mathrm{alg}}
 =\Q\otimes_{P_n}E_\sigma
 \cong(\Lambda U,0).
\end{equation}
Moreover, \(U\) is of finite type.
\end{lemma}

\begin{proof}
Choose a relative minimal Sullivan model of the morphism \(\rho_n\), followed
by an acyclic fibration.  The fibre theorem for relative Sullivan models
\cite[Theorem~15.3]{FHT2001} identifies its base change
\(\Q\otimes_{P_n}E_\sigma\) with a minimal Sullivan model of
\(G_\sigma=\hofib(\sigma)\).  Lemma~\ref{stab:lem:fibre-section-loop} gives
\(G_\sigma\simeq_\Q\Omega F_\rho\).  Since \(q\) has a section, \(q_*\)
is surjective on rational homotopy; the long exact sequence gives
\(\pi_1(F_\rho)\otimes\Q=0\).  Hence \(G_\sigma\) is connected and is a
rational \(H\)-space.

The minimal algebras \(M=(\Lambda V,d)\) and
\(P_n=(\Lambda(V\oplus Y),D)\) have finitely many generators in every
degree.  The long exact rational homotopy sequence for
\(F_\rho\to R_n\to X\), followed by
\(G_\sigma\simeq_\Q\Omega F_\rho\), shows that every
\(\pi_k(G_\sigma)\otimes\Q\) is finite dimensional.  Hence \(G_\sigma\)
is of finite rational type, and
\cite[Section~12(a), Example~3]{FHT2001} applies: its minimal Sullivan
model has zero
differential.  Thus the relative minimal model can be chosen
so that the base-changed differential vanishes, giving
\eqref{stab:eq:fibre-algebra}; this is a choice of the relative model, not an
assertion about an arbitrary prior factorization.

Since
\((U^k)^\sharp\cong\pi_k(G_\sigma)\otimes\Q\) for its minimal Sullivan
model, \(U\) is of finite type.
\end{proof}

\begin{corollary}[Local finiteness of the spherical factor]
\label{stab:cor:spherical-factor-locally-finite}
The wedge of rational spheres in
Lemma~\ref{stab:lem:spherical-factor} has finitely many summands in every
degree.  In particular, it belongs to \(\mathcal S_{\mathrm{lf}}\).
\end{corollary}

\begin{proof}
By Lemma~\ref{stab:lem:fibre-minimal-section}, the fibre-minimal model is
\((\Lambda U,0)\), with \(U\) of finite type.  Positive grading makes each
graded piece of \(\Lambda U\) finite dimensional; the same is therefore true
of every \(H^k(G_\sigma;\Q)\).  The number of sphere summands of a fixed
dimension in \(\Sigma G_\sigma\) is the corresponding rational Betti number,
so it is finite.
\end{proof}

Fix such a fibre-minimal factorization.  Let
\begin{equation}\label{stab:eq:pi-sigma}
 \pi_\sigma:E_\sigma\twoheadrightarrow\Lambda U,
 \qquad
 \pi_\sigma(p\otimes\omega)=\varepsilon(p)\omega,
\end{equation}
be the fibre projection.

\begin{proposition}[Cofibre pullback model]\label{stab:prop:cofibre-pullback}
The strict pullback
\begin{equation}\label{stab:eq:pullback-square}
\begin{tikzcd}[column sep=large,row sep=large]
 C_\sigma \arrow[r,hook] \arrow[d] & E_\sigma
       \arrow[d,"\pi_\sigma",two heads]\\
 \Q \arrow[r] & \Lambda U
\end{tikzcd}
\end{equation}
is a CDGA model of \(W_\sigma\).  As an augmented differential graded
vector space and as an algebra,
\begin{equation}\label{stab:eq:C-sigma-formula}
 C_\sigma
 =E_\sigma\times_{\Lambda U}\Q
 =\Q\oplus P_n^+\otimes\Lambda U.
\end{equation}
In particular, \(P_n=\Q\oplus P_n^+\otimes1\) is a CDGA subalgebra of
\(C_\sigma\).
\end{proposition}

\begin{proof}
The cofibre square
\[
\begin{tikzcd}[column sep=large,row sep=large]
 G_\sigma \arrow[r,"j"] \arrow[d] & X \arrow[d]\\
 * \arrow[r] & W_\sigma
\end{tikzcd}
\]
is a homotopy pushout.  By the Bousfield--Gugenheim polynomial-forms
equivalence \cite[Theorem~9.4]{BousfieldGugenheim1976}, this homotopy
pushout is sent contravariantly to a homotopy pullback.  Since
\(\pi_\sigma\) is surjective, it is a fibration in the standard model
structure on connected CDGAs
\cite[Chapter~4]{BousfieldGugenheim1976}; hence the strict pullback computes
that homotopy pullback.  This proves the model claim.

Every element of \(E_\sigma\) decomposes uniquely as a sum of a scalar,
a pure element of \(\Lambda^+U\), and terms carrying a coefficient in
\(P_n^+\).  The condition that its image under \(\pi_\sigma\) be scalar
removes precisely the pure \(\Lambda^+U\)-summand.  Formula
\eqref{stab:eq:C-sigma-formula} follows.  For \(p\in P_n\), the corresponding
element of the pullback is \((p,\varepsilon(p))\); therefore no positive
part of \(P_n\) is discarded.
\end{proof}

\subsection{Descent to a short model}

Extend the ideal \(K_n=\ker\kappa_n\) to \(E_\sigma\) and put
\begin{equation}\label{stab:eq:L-sigma}
 L_\sigma=K_nE_\sigma=K_n\otimes\Lambda U.
\end{equation}
Because \(K_n\subset P_n^+\), this is a differential ideal contained in
\(C_\sigma\).

\begin{lemma}[Semifree preservation of acyclicity]
\label{stab:lem:semifree-acyclicity}
Let \(K\) be an acyclic differential \(P\)-module and let
\(P\to(P\otimes\Lambda U,D)\) be a relative Sullivan algebra with
positively graded relative generators.  Then
\(K\otimes_P(P\otimes\Lambda U)\) is acyclic.
\end{lemma}

\begin{proof}
Choose the given Sullivan well-order \(U=\{u_\alpha\}_{\alpha<\lambda}\) and
put
\[
 E_\alpha=P\otimes\Lambda U_{<\alpha},
 \qquad K_\alpha=K\otimes_P E_\alpha.
\]
We prove by transfinite induction that every \(K_\alpha\) is acyclic.  The
initial term is \(K_0=K\).

At a successor stage, write
\(E_{\alpha+1}=E_\alpha\otimes\Lambda(u_\alpha)\).  The Sullivan condition
gives \(Du_\alpha\in E_\alpha\).  For \(p\geq0\), let
\(F_pK_{\alpha+1}\) be the subcomplex spanned by the monomials whose exponent
in \(u_\alpha\) is at most \(p\).  The differential preserves
\(F_pK_{\alpha+1}\), because every term involving \(Du_\alpha\) lowers
that exponent.  The finite filtration
\[
 0\subset F_0K_{\alpha+1}\subset\cdots\subset F_pK_{\alpha+1}
\]
has associated graded a finite direct sum of degree shifts of \(K_\alpha\),
with differential induced by that of \(K_\alpha\).  Equivalently, each
successive quotient is either zero or a degree shift of \(K_\alpha\).
The inductive hypothesis and the successive short exact sequences therefore
show that \(F_pK_{\alpha+1}\) is acyclic.  Finally,
\[
 K_{\alpha+1}=\varinjlim_p F_pK_{\alpha+1},
\]
and filtered colimits of complexes of \(\Q\)-vector spaces are exact.
Therefore \(K_{\alpha+1}\) is acyclic.  This argument does not require the
module \(K\) to be bounded below.

At a limit ordinal \(\beta\), one has
\(K_\beta=\varinjlim_{\alpha<\beta}K_\alpha\).  Filtered colimits of complexes
of \(\Q\)-vector spaces are exact, and therefore commute with cohomology.
Thus
\[
 H(K_\beta)=\varinjlim_{\alpha<\beta}H(K_\alpha)=0.
\]
The induction reaches
\(K_\lambda=K\otimes_P(P\otimes\Lambda U)\), proving the lemma.
\end{proof}

Applying Lemma~\ref{stab:lem:semifree-acyclicity} to \(K_n\) shows that
\(L_\sigma\) is acyclic.  Moreover,
\begin{equation}\label{stab:eq:Cpower-in-L}
 (C_\sigma^+)^{n+1}
 =(P_n^+)^{n+1}\otimes\Lambda U
 \subset K_n\otimes\Lambda U=L_\sigma.
\end{equation}
We have therefore proved the following proposition.

\begin{proposition}[The short cofibre model]\label{stab:prop:short-cofibre}
There is a surjective quasi-isomorphism
\begin{equation}\label{stab:eq:q-sigma-short}
 q_\sigma:C_\sigma\twoheadrightarrow D_\sigma=C_\sigma/L_\sigma
\end{equation}
with
\begin{equation}\label{stab:eq:D-sigma-formula}
 D_\sigma\cong
 \Q\oplus A_n^+\otimes\Lambda U,
 \qquad
 (D_\sigma^+)^{n+1}=0.
\end{equation}
For \(n\geq1\),
\begin{equation}\label{stab:eq:QD-sigma}
 Q(D_\sigma)\cong V\otimes\Lambda U,
 \qquad Q(d_{D_\sigma})=0.
\end{equation}
\end{proposition}

\begin{proof}
Only \eqref{stab:eq:QD-sigma} remains to be checked.  Since \(n\geq1\), the
quotient \(A_n=M/(M^+)^{n+1}\) does not alter indecomposables, so
\(Q(A_n)\cong V\).  If \(I=A_n^+\), then
\[
 D_\sigma^+=I\otimes\Lambda U,
 \qquad
 (D_\sigma^+)^2=I^2\otimes\Lambda U,
\]
which gives the first isomorphism.  The absolute minimality of \(P_n\)
implies \(D(P_n^+)\subset(P_n^+)^2\).  Fibre minimality gives
\(D(U)\subset P_n^+\otimes\Lambda U\).  Consequently, for
\(p\in P_n^+\) and \(\omega\in\Lambda U\),
\[
 D(p\omega)=Dp\,\omega+(-1)^{|p|}pD\omega
 \in(P_n^+)^2\otimes\Lambda U.
\]
The induced differential on the indecomposable quotient is zero.
\end{proof}

Proposition~\ref{stab:prop:short-cofibre} gives an \(n\)-short model of
\(W_\sigma\), but it does not yet compute \(\Hnilz(W_\sigma)\): the
source \(C_\sigma\) need not be minimal.  The remainder of the argument
proves that a minimal model of \(C_\sigma\) maps onto it.

\section{Exact relative disk normal form}\label{stab:sec:disk-normal-form}
\index{open-disk normal form}

The purpose of this section is to remove all nonlinear tails from the
relative disk pairs in \eqref{stab:eq:section-relative-factor}.  This is an
exact coordinate change, not an associated-graded argument.

Write, as a graded algebra under \(M\),
\begin{equation}\label{stab:eq:P-MY}
 P_n=M\otimes\Lambda Y.
\end{equation}
At this point \(Y\) is the normalized relative generator space supplied
by Theorem~\ref{stab:thm:root}, so \(\kappa_n(Y)=0\).  Before making any
further change of these generators, we record the low-degree consequence
\begin{equation}\label{stab:eq:Y12-zero}
                         Y^1=Y^2=0.
\end{equation}
Indeed, if \(Y^1\neq0\), let \(y\) be the first degree-one relative
generator in the Sullivan order.  The absolute minimality of \(P_n\)
forces \(Dy\) to be decomposable.  A decomposable element of degree two
would have to be a product of earlier degree-one generators, but
\(M^1=0\) and no such relative generator precedes \(y\).  Hence
\(Dy=0\).  Since \(P_n^0=\Q\), the element \(y\) is not a boundary,
whereas \(\kappa_n(y)=0\); this contradicts the fact that
\(\kappa_n:P_n\to A_n\) is a quasi-isomorphism.  Thus \(Y^1=0\).
Now every generator \(y\in Y^2\) is closed: its differential has degree
three, while every decomposable in \(P_n\), whose positive-degree
generators now have degree at least two, has degree at least four.  It
cannot be a boundary because \(P_n^1=0\), and again
\(\kappa_n(y)=0\) contradicts the quasi-isomorphism \(\kappa_n\).  This proves
\eqref{stab:eq:Y12-zero}, and consequently
\begin{equation}\label{stab:eq:Pn12}
                         P_n^1=0,\qquad P_n^2=M^2.
\end{equation}
The identities in \eqref{stab:eq:Pn12} are intrinsic to \(P_n\).  Because
\(Y^1=Y^2=0\), the triangular change of variables below does not affect them.

Replacing a relative generator \(y\) by
\(y-i_n\rho_n(y)\) is a triangular change under \(M\).  We may therefore
assume
\begin{equation}\label{stab:eq:rho-Y-zero}
 \rho_n(Y)=0.
\end{equation}
Likewise, replacing \(u\in U\) by \(u-i_n\lambda_\sigma(u)\) allows us
to assume \(\lambda_\sigma(U)=0\).

Suppose next that \(u\in U^1\).  All relative Sullivan generators have
positive degree, and the vanishing of the base-changed differential in
\eqref{stab:eq:fibre-algebra} gives \(\pi_\sigma(Du)=0\), thereby excluding
every term of \(Du\) that is purely in \(\Lambda^+U\).  Since
\(P_n^1=0\), no mixed term can have total degree two.  Therefore
\(Du\in P_n^2=M^2\).  Our normalization
\(\lambda_\sigma(u)=0\) gives
\[
  0=d\lambda_\sigma(u)=\lambda_\sigma(Du)=\rho_n(Du).
\]
Since \(\rho_n\) restricts to the identity on \(M\), it follows that
\(Du=0\).  The resulting nonzero class \([u]\in H^1(E_\sigma)\) would
be killed by the quasi-isomorphism \(\lambda_\sigma:E_\sigma\to M\),
contrary to \(H^1(M)=0\).  Hence \(U^1=0\), and
\begin{equation}\label{stab:eq:E-sigma-simply-connected}
                         E_\sigma^1=0.
\end{equation}
In particular, \(E_\sigma\) is simply connected, so
Proposition~\ref{stab:prop:min-contractible} applies below.

The indecomposable complex of \(E_\sigma\) is
\[
 Q(E_\sigma)=V\oplus Y\oplus U.
\]
The linear differential is zero on \(V\oplus Y\), because \(P_n\) is
absolutely minimal.  Since the fibre differential is zero, its only
possible remaining component is
\begin{equation}\label{stab:eq:D1-UY}
 D_1:U\longrightarrow Y[1].
\end{equation}
There is no component \(U\to V[1]\): applying \(\lambda_\sigma\) would
leave such a component unchanged on \(V\), whereas
\(\lambda_\sigma(U)=0\) gives
\(\lambda_\sigma D(U)=d\lambda_\sigma(U)=0\).

\begin{lemma}[Relative linear disks]\label{stab:lem:relative-linear-disks}
The map \(D_1:U\to Y[1]\) in \eqref{stab:eq:D1-UY} is an isomorphism in
every degree.
\end{lemma}

\begin{proof}
The inclusion \(M\hookrightarrow E_\sigma\) is a quasi-isomorphism:
its composite with \(\lambda_\sigma\) is \(\id_M\), and
\(\lambda_\sigma\) is a quasi-isomorphism.  Apply
Proposition~\ref{stab:prop:min-contractible} to this inclusion under \(M\).
It gives, under \(M\), a tensor decomposition
\[
 E_\sigma\cong M'\otimes\Lambda(W\oplus\widehat W),
 \qquad Dw=\widehat w,\quad D\widehat w=0,
\]
where \(M'\) is minimal and contains the displayed copy of \(M\).  The disk
factor is acyclic, so the inclusion \(M\hookrightarrow M'\) is a
quasi-isomorphism.  A quasi-isomorphism between minimal Sullivan algebras is
an isomorphism, by induction on degree and the Sullivan order; because the
map is the identity on \(M\), this identifies \(M'\) with \(M\) under
\(M\), not merely up to an abstract equivalence.  Consequently all relative
indecomposables \(Y\oplus U\) form linear disk pairs.  Their only possible
linear differential is \eqref{stab:eq:D1-UY}; hence
\(D_1:U\to Y[1]\) is an isomorphism.
\end{proof}

Choose homogeneous bases \(\{u_\alpha\}\) of \(U\) and
\(\{y_\alpha\}\) of \(Y\), compatible with the Sullivan order, such
that
\begin{equation}\label{stab:eq:Du-y-Phi}
 Du_\alpha=y_\alpha+\Phi_\alpha,
 \qquad \Phi_\alpha\in(E_\sigma^+)^2.
\end{equation}
Because the reduction of \(Du_\alpha\) to the fibre is zero, every
monomial of \(\Phi_\alpha\) contains a factor in \(M^+\) or a relative
generator in \(Y\).

\begin{theorem}[Exact disk standardization]
\label{stab:thm:exact-standardization}
Set
\begin{equation}\label{stab:eq:def-hat-u}
 \widehat u_\alpha=Du_\alpha.
\end{equation}
The substitution \(y_\alpha\mapsto\widehat u_\alpha\), fixing \(M\)
and \(U\), is a triangular graded-algebra automorphism.  It identifies
\(E_\sigma\) with
\begin{equation}\label{stab:eq:E-standard}
 E_0=M\otimes\Lambda(U\oplus\widehat U),
 \qquad
 Du_\alpha=\widehat u_\alpha,\quad
 D\widehat u_\alpha=0,\quad D|_M=d.
\end{equation}
Under this isomorphism,
\begin{equation}\label{stab:eq:C-standard}
 C_\sigma\cong
 C:=\Q\oplus(M\otimes\Lambda\widehat U)^+\otimes\Lambda U.
\end{equation}
\end{theorem}

\begin{proof}
Equation~\eqref{stab:eq:def-hat-u} and \(D^2=0\) give
\(Du_\alpha=\widehat u_\alpha\) and \(D\widehat u_\alpha=0\).
The linear term of \(\widehat u_\alpha\) is \(y_\alpha\).  Every relative
generator \(y_\beta\) occurring in the decomposable element
\(\Phi_\alpha\) has degree strictly smaller than \(|y_\alpha|\): a generator
of the same degree could not be multiplied by another positive-degree
element without exceeding \(|\Phi_\alpha|=|y_\alpha|\).  Thus the
substitution is triangular.  Its inverse is constructed recursively with
respect to degree and the Sullivan order.  This proves \eqref{stab:eq:E-standard}.

It remains to track the open ideal.  In the original coordinates the
positive part of \(C_\sigma\) is the ideal
\((M^+,Y)E_\sigma\).  Since the fibre differential is zero,
\(\widehat u_\alpha=Du_\alpha\) belongs to that ideal, and hence
\((M^+,\widehat U)\subset(M^+,Y)\).  Conversely,
\[
 y_\alpha=\widehat u_\alpha-\Phi_\alpha.
\]
Every monomial of \(\Phi_\alpha\) contains either \(M^+\) or an earlier
generator \(y_\beta\).  Induction gives
\(y_\alpha\in(M^+,\widehat U)\), whence
\[
 (M^+,Y)E_\sigma=(M^+,\widehat U)E_0.
\]
Adding the unit yields \eqref{stab:eq:C-standard}.
\end{proof}

The change of coordinates in Theorem~\ref{stab:thm:exact-standardization}
need not preserve the earlier normalization \(\kappa_n(Y)=0\).  This does not
affect the argument: Proposition~\ref{stab:prop:short-cofibre} has already provided the
surjective quasi-isomorphism \(C_\sigma\to D_\sigma\).  We use the new
coordinates only to prove that the minimal model surjects onto
\(C_\sigma\).
Thus the geometrically constructed cofibre model has exactly the open-disk
form required in Chapter~\ref{chap:stabilization-proof}.

\chapter{General and Explicit Spherical Stabilization}
\label{chap:stabilization-proof}
\begin{chapterguide}[title={Chapter guide}]
This chapter assembles the open-disk surjectivity argument and proves the
general spherical stabilization theorem.  It then specializes the construction
to the separating model of Part~\ref{part:divergence}, where puncturing the
hidden disk shows that one rational sphere suffices.  The final section
distinguishes this finite example from the degreewise finite wedge required in
general.
\end{chapterguide}
\index{surjective minimal model}
\index{spherical stabilization theorem}
\index{spherical stabilization!explicit one-sphere example}
\index{finite wedge caveat}

\section{Surjectivity for the general open-disk model}
\label{stab:sec:minimal-surjectivity}

Theorem~\ref{stab:thm:mixed-AQ-edge} treats the mixed summand.  We now establish
surjectivity on the other indecomposable summands and then pass to
algebra surjectivity.

\begin{lemma}[The base summand]\label{stab:lem:base-summand}
For a minimal model \(f:N\to C\), the summand
\(V\otimes1\subset Q(C)\) is contained in \(\im Qf\).
\end{lemma}

\begin{proof}
There are strict CDGA morphisms
\[
 M\xrightarrow{\eta}C\xrightarrow{r_M}M,
 \qquad r_M\eta=\id_M,
\]
where \(\eta(m)=m\otimes1\) and \(r_M\) kills the disk variables.
In the homotopy category, the quasi-isomorphism \(f\) is invertible, so
\(f^{-1}\eta\) is represented by a Sullivan morphism \(j:M\to N\).
Then \(fj\) is Sullivan-homotopic to \(\eta\).  Sullivan-homotopic maps
from a minimal algebra to a CDGA with \(Qd=0\) induce the same map on
indecomposables.  Indeed, if \(\mathcal H:M\to C\otimes\Lambda(t,dt)\)
is such a homotopy, let \(\pi:C^+\to Q(C)\) be the quotient and, for a
homogeneous generator \(v\), write
\[
 (\pi\otimes\id)\mathcal H(v)=a_v(t)+b_v(t)dt.
\]
Projecting the homotopy identity to \(Q(C)\otimes\Lambda(t,dt)\) and
integrating the \(dt\)-coefficient defines a degree \(-1\) map \(h\),
with the standard Koszul sign, and gives the chain-homotopy formula
\[
             Q(fj)-Q\eta=Qd_C\,h+h\,Qd_M.
\]
Here \(Qd_M=0\) by minimality and \(Qd_C=0\) by
Lemma~\ref{stab:lem:QC-open}; hence \(Q(fj)=Q\eta\).  That lemma therefore gives
\[
 Qf\,Qj(v)=v\otimes1,
\]
as required.
\end{proof}

\begin{lemma}[The disk summand]\label{stab:lem:disk-summand}
For a minimal model \(f:N\to C\),
\begin{equation}\label{stab:eq:disk-summand-captured}
 \widehat U\otimes T\subset\im Qf.
\end{equation}
\end{lemma}

\begin{proof}
For \(\widehat u\in\widehat U\), the element \(\widehat u\) is a
cocycle.  More generally,
\begin{equation}\label{stab:eq:d-hatu-omega}
 D(\widehat u\omega)
 =(-1)^{|\widehat u|}\widehat u\,D\omega.
\end{equation}
Both factors on the right are positive cocycles of \(C\):
\(D\widehat u=0\) and \(D^2\omega=0\).  Since \(Q(D_C)=0\), the
products-of-cocycles Lemma~\ref{stab:lem:products-cocycles-prelim} applies and
places the indecomposable class of every \(\widehat u\omega\) in
\(\im Qf\).
\end{proof}

\begin{theorem}[Surjective minimal model of the open algebra]
\label{stab:thm:open-surjective-minimal}
The standard open-disk CDGA \(C\) of \eqref{stab:eq:open-C} admits a
surjective quasi-isomorphism from its minimal Sullivan model:
\begin{equation}\label{stab:eq:minimal-onto-C}
 f:N_C\twoheadrightarrow C.
\end{equation}
\end{theorem}

\begin{proof}
Choose any minimal model \(f:N_C\to C\).  By
Theorem~\ref{stab:thm:mixed-AQ-edge}, the image of \(Qf\) projects onto
\(V\otimes T^+\).  Lemma~\ref{stab:lem:disk-summand} supplies preimages in
\(\widehat U\otimes T\), and Lemma~\ref{stab:lem:base-summand} supplies
preimages in \(V\otimes1\).  Subtracting preimages in these two already
treated summands gives
\[
 V\otimes T^+\subset\im Qf.
\]
Together with \eqref{stab:eq:disk-summand-captured} and the base summand, this
exhausts
\[
 Q(C)=(V\otimes1)\oplus(V\otimes T^+)
       \oplus(\widehat U\otimes T).
\]
Thus \(Qf\) is surjective.  Lemma~\ref{stab:lem:Q-surj} implies that \(f\)
is surjective as a CDGA morphism.  It was already a quasi-isomorphism by
the definition of a minimal model.
\end{proof}

\section{The general spherical stabilization theorem}
\label{stab:sec:main-proof}

We now assemble the preceding results.

\begin{theorem}\label{stab:thm:fibre-cofibre-Hnil}
Let \(X\), \(n\), and \(W_\sigma\) be as in
Theorem~\ref{stab:thm:intro-main}, and suppose \(n\geq1\).  Then
\begin{equation}\label{stab:eq:Hnil-Wsigma}
 \Hnilz(W_\sigma)=n.
\end{equation}
\end{theorem}

\begin{proof}
By Proposition~\ref{stab:prop:cofibre-pullback}, \(C_\sigma\) is a CDGA model
of \(W_\sigma\).  The exact standardization theorem identifies it with
the open algebra \(C\).  Theorem~\ref{stab:thm:open-surjective-minimal}
therefore provides a surjective quasi-isomorphism
\[
 N_\sigma=\mathcal M_{W_\sigma}
 \twoheadrightarrow C_\sigma.
\]
Compose it with the short quotient of
Proposition~\ref{stab:prop:short-cofibre}:
\begin{equation}\label{stab:eq:final-composite}
 N_\sigma\twoheadrightarrow C_\sigma
 \xrightarrow{q_\sigma}D_\sigma.
\end{equation}
Both arrows are surjective quasi-isomorphisms.  The target is \(n\)-short,
so Lemma~\ref{div:lem:surjective-characterization} gives
\(\Hnilz(W_\sigma)\leq n\).

The wedge equivalence \eqref{stab:eq:wedge-splitting} makes \(X\) a rational
homotopy retract of \(W_\sigma\).  LS category is monotone under
homotopy retracts, and Proposition~\ref{stab:prop:cat-lower} gives
\[
 n=\catz(X)\leq\catz(W_\sigma)\leq\Hnilz(W_\sigma)\leq n.
\]
Every inequality is an equality, proving \eqref{stab:eq:Hnil-Wsigma}.
\end{proof}

\begin{proof}[Proof of Theorem~\ref{stab:thm:intro-main}]
If \(n=0\), then \(X\simeq_\Q *\), its minimal Sullivan model is \(\Q\),
and \(\Hnilz(X)=0\).  With the empty wedge, the choices in the statement
give \(G_\sigma=*\) and \(W_\sigma\simeq_\Q X\), so all assertions follow.
Assume henceforth that \(n\geq1\).

The rational wedge decomposition
\(W_\sigma\simeq_\Q X\vee\Sigma G_\sigma\) is
\eqref{stab:eq:wedge-splitting}; its spherical nature follows from
Lemma~\ref{stab:lem:spherical-factor}, and degreewise finiteness from
Corollary~\ref{stab:cor:spherical-factor-locally-finite}.  The equality
\(\Hnilz(W_\sigma)=n\) is Theorem~\ref{stab:thm:fibre-cofibre-Hnil}.

For every admissible rational wedge of spheres \(S\), the collapse and
inclusion make \(X\) a homotopy retract of \(X\vee S\).  Hence
\[
 \catz(X)\leq\catz(X\vee S)\leq\Hnilz(X\vee S).
\]
Taking the infimum over \(S\in\mathcal S_{\mathrm{lf}}\) gives
\(\catz(X)\leq\sHnil(X)\).  By
Lemma~\ref{stab:lem:spherical-factor} and
Corollary~\ref{stab:cor:spherical-factor-locally-finite}, choose
an admissible wedge \(S_\sigma\in\mathcal S_{\mathrm{lf}}\) and a rational
equivalence \(S_\sigma\simeq_\Q\Sigma G_\sigma\).  Then
\(X\vee S_\sigma\simeq_\Q W_\sigma\), and the first part of the theorem
gives \(\Hnilz(X\vee S_\sigma)=n\); hence
\(\sHnil(X)\leq\catz(X)\).
Therefore
\[
 \sHnil(X)=\catz(X).
\]
\end{proof}

\begin{corollary}[Parallel spherical stabilizations]
\label{stab:cor:parallel-spherical-stabilizations}
For every simply connected rational space \(X\) of finite type,
\begin{equation}\label{stab:eq:parallel-spherical-stabilizations}
              \sClz(X)=\sHnil(X)=\catz(X).
\end{equation}
\end{corollary}

\begin{proof}
For every \(S\in\mathcal S_{\mathrm{lf}}\), the rational homotopy retraction
\(X\to X\vee S\to X\) and Proposition~\ref{stab:prop:cat-lower} give
\[
 \catz(X)\leq\catz(X\vee S)\leq\Clz(X\vee S)
       \leq\Hnilz(X\vee S).
\]
If \(\catz(X)<\infty\), taking infima and using
Theorem~\ref{stab:thm:intro-main} yields
\(\catz(X)\leq\sClz(X)\leq\sHnil(X)=\catz(X)\).  If
\(\catz(X)=\infty\), the displayed lower bounds force both stabilized
invariants to be infinite.  The equality follows in both cases.
\end{proof}

\begin{corollary}\label{stab:cor:nilh-Hnil-equality}
For the fibre--cofibre witness associated with the chosen optimal section,
\begin{equation}\label{stab:eq:all-three-equal}
 \catz(W_\sigma)=\Clz(W_\sigma)
 =\nilh(W_\sigma)=\Hnilz(W_\sigma)=n.
\end{equation}
\end{corollary}

\begin{proof}
Proposition~\ref{stab:prop:cat-lower} gives
\(\catz\leq\Clz=\nilh\leq\Hnilz\), and
Theorem~\ref{stab:thm:fibre-cofibre-Hnil} identifies the two endpoints.
\end{proof}

\section{The separating example stabilizes after one sphere}
\label{stab:sec:explicit-S17}

Let
\[
             \varphi:M\xrightarrow{\simeq}B
\]
be the fixed minimal model of Corollary~\ref{div:cor:rigid-prefix}.  Its
Sullivan realization \(|M|\) is the separating space of
Corollary~\ref{div:cor:topological-main}.  We retain the generators and
differential of Definition~\ref{div:def:B} throughout.  We show explicitly
that, after wedging with one rational sphere, its homology nilpotency attains
rational category.

\subsection{The exact category of the separating space}

\begin{proposition}\label{stab:prop:separating-category}
The separating space satisfies
\[
                         \catz(|M|)=3.
\]
\end{proposition}

\begin{proof}
The essential cycle \(c=axt\) of
Proposition~\ref{div:prop:c-essential} belongs to \((M^+)^3\).  Hence the
word-length projection
\[
                   q_2:M\longrightarrow M/(M^+)^3
\]
kills \([c]\).  If \(\catz(|M|)\leq2\), the
F\'elix--Halperin characterization recalled before
Proposition~\ref{stab:prop:cat-lower} would give a homotopy retraction of
\(q_2\).  The induced map \(H(q_2)\) would then be injective, contradicting
\([c]\neq0\).  Thus \(\catz(|M|)>2\).

On the other hand, Corollary~\ref{div:cor:topological-main} gives
\(\Clz(|M|)=3\), and
Proposition~\ref{stab:prop:cat-lower} gives
\(\catz\leq\Clz\).  Therefore \(\catz(|M|)=3\).
\end{proof}

\subsection{Puncturing the hidden disk}
\index{punctured disk}

Fix the canonical monomial basis of \(B\) associated with the generator order
used in Chapter~\ref{div:sec:witness}.

\begin{definition}\label{stab:def:punctured-witness}
Let \(B^\circ\) be the span of all canonical monomials of \(B\) except the
linear monomial \(r\).  Thus, as graded vector spaces,
\[
                         B=B^\circ\oplus\Q r.
\]
\end{definition}

For \(0\leq\ell\leq3\), let \(B_{\langle\ell\rangle}\) denote the span of
the canonical monomials of augmentation length \(\ell\).  Then
\[
 B=\bigoplus_{\ell=0}^3 B_{\langle\ell\rangle},
 \qquad
 B^\circ
 =\Q 1\oplus\Q\{s,x,y,a,t\}
   \oplus B_{\langle2\rangle}\oplus B_{\langle3\rangle}.
\]
Thus the puncture removes only the one-dimensional linear direction
\(\Q r\); it retains all nonlinear monomials containing \(r\).

\begin{proposition}\label{stab:prop:punctured-witness}
The algebra \(B^\circ\) is a three-short sub-CDGA of \(B\), of dimension
\(55\).  The inclusion \(B^\circ\hookrightarrow B\) induces an isomorphism
in cohomology in every degree other than \(17\), while
\[
                         H^{17}(B^\circ)=\Q[xt].
\]
\end{proposition}

\begin{proof}
The displayed decomposition shows that \(B^\circ\) is graded and contains the
unit.  Multiplication respects augmentation length:
\[
 B_{\langle i\rangle}B_{\langle j\rangle}
 \subseteq B_{\langle i+j\rangle},
\]
where the right-hand side is zero for \(i+j\geq4\).  A product of elements of
\(B^\circ\) could have a component along the omitted line \(\Q r\) only in
augmentation length one.  One factor would then be a scalar and the other
would already have a linear \(r\)-component, which is excluded.  Hence
\(B^\circ\) is a unital graded subalgebra of \(B\).

The differential formulas in \eqref{div:eq:differential} give
\[
 dB_{\langle1\rangle}
   \subseteq B_{\langle2\rangle}\oplus B_{\langle3\rangle},
 \qquad
 dB_{\langle2\rangle}\subseteq B_{\langle3\rangle},
 \qquad
 dB_{\langle3\rangle}=0.
\]
Indeed, differentiating a generator replaces it by a word of length at least
two, and every word of length at least four vanishes in \(B\).  This also
covers monomials containing \(r\): differentiating an \(r\)-factor replaces
it by \(xt\), while differentiating another factor leaves \(r\) accompanied
by at least two further factors.  Thus \(d(B^\circ)\subseteq B^\circ\).
Koszul signs affect only coefficients, not the length estimates.  Since
\(d^2=0\) on \(B\), the restricted differential also squares to zero.  The
augmentation of \(B\) restricts to \(B^\circ\), and
\[
                       ((B^\circ)^+)^4=0.
\]
Therefore \(B^\circ\) is an augmented three-short sub-CDGA.  Its dimension is
\(56-1=55\) by Proposition~\ref{div:prop:B-properties}.

There is a short exact sequence of complexes
\begin{equation}\label{stab:eq:punctured-B-sequence}
 0\longrightarrow B^\circ\longrightarrow B
   \longrightarrow\Q r_{16}\longrightarrow0.
\end{equation}
The quotient differential is zero because \(dr=xt\in B^\circ\).
Lemma~\ref{div:lem:critical-groups} gives
\(H^{16}(B)=H^{17}(B)=0\), and the connecting homomorphism sends the class of
\(r\) to \([xt]\).  The long exact cohomology sequence therefore gives the
asserted class in degree \(17\) and an isomorphism in every other degree.
\end{proof}

Two points about the notation are important.  First, \(B^\circ\) is neither
the quotient \(B/(r)\) nor the subalgebra generated by
\(s,x,y,a,t\).  It omits only the linear vector \(r\), while retaining every
allowed nonlinear monomial containing \(r\), including
\[
                  sr,\ yr,\ xr,\ ar,\ tr,\ r^2,\ r^3.
\]
The superscript \({}^\circ\) thus denotes a puncture in the linear
\(r\)-direction, not the removal of all \(r\)-dependence.

Second, \(B^\circ\) is not an ideal of \(B\): it contains the unit but not
\(r\).  Accordingly, \eqref{stab:eq:punctured-B-sequence} is a short exact
sequence of cochain complexes, not a quotient sequence of CDGAs.  Its last
map extracts the coefficient of the linear monomial \(r\).  The identity
\(dr=xt\) is precisely what makes that coefficient map a chain map and what
exposes the class \([xt]\) after the puncture.

Consider the Sullivan disk on the contractible pair
\[
 \mathbb D_{17}=\bigl(\Lambda(\rho_{16},\eta_{17}),d\bigr),
 \qquad d\rho=\eta,\qquad d\eta=0,
\]
and remove only its linear \(\rho\)-direction:
\begin{equation}\label{stab:eq:punctured-disk}
 \mathbb D_{17}^{\circ}
 =\Q\oplus\Q\eta\oplus
   \bigoplus_{k\geq2}
   \bigl(\Q\rho^k\oplus\Q\rho^{k-1}\eta\bigr).
\end{equation}

\begin{lemma}\label{stab:lem:punctured-disk}
The formula
\[
 \mathbb D_{17}^{\circ}\longrightarrow(\Lambda(u_{17}),0),
 \qquad \eta\longmapsto u,\qquad
 \rho^k,\rho^{k-1}\eta\longmapsto0\quad(k\geq2),
\]
defines a quasi-isomorphism.  Consequently,
\(\mathbb D_{17}^{\circ}\) is a CDGA model of \(S^{17}_{\Q}\).
\end{lemma}

\begin{proof}
The subspace in \eqref{stab:eq:punctured-disk} is closed under products and
differentials.  For every \(k\geq2\),
\[
                    d(\rho^k)=k\rho^{k-1}\eta,
                    \qquad d(\rho^{k-1}\eta)=0.
\]
These are acyclic pairs, so the positive-degree cohomology is
one-dimensional, generated by \([\eta]\).
Since \(17\) is odd, \(\Lambda(u_{17})=\Q\oplus\Q u\), and the displayed
morphism induces an isomorphism in cohomology.
\end{proof}

\subsection{A three-short model of the wedge}

Put \(E_{17}=M\otimes\mathbb D_{17}\), and let
\(p:E_{17}\twoheadrightarrow\mathbb D_{17}\) be induced by the augmentation
of \(M\).  Form the strict pullback
\begin{equation}\label{stab:eq:explicit-S17-pullback}
\begin{tikzcd}[column sep=large,row sep=large]
 O_{17} \arrow[r,hook] \arrow[d]
    & E_{17} \arrow[d,"p",two heads]\\
 \mathbb D_{17}^{\circ} \arrow[r,hook]
    & \mathbb D_{17}.
\end{tikzcd}
\end{equation}
Because \(p\) is surjective, this strict pullback computes the homotopy
pullback.  Sullivan realization sends it to the homotopy pushout
\(|M|\leftarrow *\to S^{17}_{\Q}\); hence \(O_{17}\) is a
model of \(|M|\vee S^{17}_{\Q}\).

For an entirely explicit comparison, let
\(q:\mathbb D_{17}^{\circ}\to\Lambda(u_{17})\) be the map in
Lemma~\ref{stab:lem:punctured-disk}.  The decomposition
\[
              O_{17}=\mathbb D_{17}^{\circ}
                       \oplus(M^+\otimes\mathbb D_{17})
\]
and the two augmentations define
\[
 \Psi:O_{17}\longrightarrow M\times_{\Q}\Lambda(u_{17}),
 \qquad e\longmapsto
 \bigl((\id_M\otimes\varepsilon_{\mathbb D_{17}})(e),q(p(e))\bigr).
\]
The maps on the kernels and quotients in
\[
\begin{tikzcd}[column sep=small]
0\arrow[r] & M^+\otimes\mathbb D_{17}\arrow[r]\arrow[d]
 & O_{17}\arrow[r]\arrow[d,"\Psi"]
 & \mathbb D_{17}^{\circ}\arrow[r]\arrow[d,"q"] & 0\\
0\arrow[r] & M^+\arrow[r]
 & M\times_{\Q}\Lambda(u_{17})\arrow[r]
 & \Lambda(u_{17})\arrow[r] & 0
\end{tikzcd}
\]
are quasi-isomorphisms: on the left this follows from the contractibility of
the Sullivan disk, and on the right it is
Lemma~\ref{stab:lem:punctured-disk}.  The long exact sequences show that
\(\Psi\) is a quasi-isomorphism to the standard pullback model of the wedge.

Define
\begin{equation}\label{stab:eq:Phi-explicit-S17}
 \Phi:E_{17}\longrightarrow B,
 \qquad \Phi|_M=\varphi,\qquad
 \Phi(\rho)=r,\qquad \Phi(\eta)=xt.
\end{equation}
This is a CDGA morphism because \(dr=xt\) and \(d_B(xt)=0\).  It is
surjective: \(\varphi(M)\) contains \(s,x,y,a,t\), and \(\rho\) supplies
\(r\).  The inclusion \(M\hookrightarrow E_{17}\) is a quasi-isomorphism
and its composite with \(\Phi\) is \(\varphi\), so \(\Phi\) is a
quasi-isomorphism by two-out-of-three.

\begin{proposition}\label{stab:prop:explicit-wedge-model}
The morphism \(\Phi\) restricts to a surjective quasi-isomorphism
\[
                  \Phi^\circ:O_{17}\ontoqiso B^\circ.
\]
In particular, \(B^\circ\) is a three-short CDGA model of
\(|M|\vee S^{17}_{\Q}\).
\end{proposition}

\begin{proof}
The pullback condition excludes a pure linear \(\rho\)-term.  Moreover,
\(M^{16}=0\) by Lemma~\ref{div:lem:critical-groups} and
Corollary~\ref{div:cor:rigid-prefix}, so \(\varphi(M)\) cannot contribute a
linear \(r\)-component.  Thus \(\Phi(O_{17})\subset B^\circ\).

Every monomial of \(B^\circ\) without \(r\) lifts from the subalgebra of
\(M\) on \(s,x,y,a,t\).  A monomial containing \(r\) and another positive
factor lifts from \(M^+\otimes\mathbb D_{17}\), while \(r^2\) and \(r^3\)
lift from \(\rho^2\) and \(\rho^3\), which belong to
\(\mathbb D_{17}^{\circ}\).  Hence \(\Phi^\circ\) is surjective.

Finally compare the short exact sequences
\[
\begin{tikzcd}[column sep=small]
0\arrow[r] & O_{17}\arrow[r]\arrow[d,"\Phi^\circ"]
 & E_{17}\arrow[r]\arrow[d,"\Phi"]
 & \mathbb D_{17}/\mathbb D_{17}^{\circ}
      \arrow[r]\arrow[d,"\overline\Phi","\cong"'] & 0\\
0\arrow[r] & B^\circ\arrow[r]
 & B\arrow[r]
 & B/B^\circ\arrow[r] & 0.
\end{tikzcd}
\]
The right-hand complexes are respectively \(\Q\rho_{16}\) and
\(\Q r_{16}\), both with zero differential, and
\(\overline\Phi(\rho)=r\).  The central map is a quasi-isomorphism; the long
exact cohomology sequences therefore show that \(\Phi^\circ\) is one as
well.
\end{proof}

\subsection{Surjectivity of the minimal model}

The preceding proposition supplies a short model of the wedge, but homology
nilpotency requires this short model to be a quotient of its minimal model.
The following elementary lemma provides exact, rather than merely
cohomological, lifts.

\begin{lemma}[Degreewise exact lifting]
\label{stab:lem:degreewise-exact-lifting}
Let \(g:N\to A\) be a quasi-isomorphism of connected CDGAs, let \(q\geq1\),
and suppose that \(g:N^{q-1}\to A^{q-1}\) is surjective.
\begin{enumerate}[label=\textup{(\roman*)}]
\item Every \(z\in Z^q(A)\) has a lift
      \(\widetilde z\in Z^q(N)\) with \(g(\widetilde z)=z\).
\item If \(b\in A^q\) and \(c\in Z^{q+1}(N)\) satisfy
      \(g(c)=db\), then there is \(\widetilde b\in N^q\) such that
      \(d\widetilde b=c\) and \(g(\widetilde b)=b\).
\end{enumerate}
Moreover, if \(Q(g)\) is surjective through degree \(m\), then \(g\) is
surjective through degree \(m\).
\end{lemma}

\begin{proof}
For (i), choose \(z_0\in Z^q(N)\) with
\([g(z_0)]=[z]\).  Write \(z-g(z_0)=da\), lift
\(a\in A^{q-1}\) to \(\widetilde a\in N^{q-1}\), and take
\(\widetilde z=z_0+d\widetilde a\).

For (ii), injectivity of \(H(g)\) first gives \(c=dv\) for some
\(v\in N^q\).  The element \(b-g(v)\) is a cocycle.  Part (i) gives a
cocycle \(z'\in N^q\) with \(g(z')=b-g(v)\), and
\(\widetilde b=v+z'\) has the required properties.

For the last assertion, induct on the cohomological degree.  After lifting the
indecomposable part of an element, the remainder is decomposable, and each of
its positive factors has strictly smaller degree.  This is the degreewise form
of the argument in Lemma~\ref{stab:lem:Q-surj}.
\end{proof}

\begin{proposition}\label{stab:prop:punctured-minimal-surjective}
Every minimal Sullivan model
\[
                    g:N_{17}=(\Lambda Z,d)\xrightarrow{\simeq}B^\circ
\]
is surjective.
\end{proposition}

\begin{proof}
A monomial inspection in the three-short algebra gives
\begin{equation}\label{stab:eq:QB-circ}
 Q(B^\circ)=
 \Q\{s,y,x,a,t,sr,yr,xr,ar,tr,r^2,r^3\}.
\end{equation}
Indeed, without \(r\) only the five linear generators remain
indecomposable.  Each quadratic monomial containing \(r\) is
indecomposable because the linear factor \(r\) is absent from \(B^\circ\).
Every cubic monomial containing another generator factors in \(B^\circ\),
whereas \(r^3\) does not.

The following four-column table is the complete lifting schedule.  A row with
zero target differential uses part~(i) of
Lemma~\ref{stab:lem:degreewise-exact-lifting}; every other row uses part~(ii)
and the closed element in the last column.  The two rows labelled
``intermediate'' create cycles needed by later differentials but do not add
new indecomposable directions.
\begin{center}
\footnotesize
\renewcommand{\arraystretch}{1.17}
\begin{tabularx}{\linewidth}{@{}>{\raggedright\arraybackslash}p{0.16\linewidth}
 c
 >{\raggedright\arraybackslash}p{0.21\linewidth}
 >{\raggedright\arraybackslash}X@{}}
\toprule
Target direction & Degree & Differential in \(B^\circ\)
 & Closed datum in \(N_{17}\) used for the lift\\
\midrule
\(s,y\) & 3 & \(0\) & choose cocycles
  \(\widetilde s,\widetilde y\)\\
\(x\) & 5 & \(0\) & choose a cocycle \(\widetilde x\)\\
\(a\) & 7 & \(xy\) & \(\widetilde x\widetilde y\)\\
\(t\) & 12 & \(ays\) &
  \(\widetilde a\widetilde y\widetilde s\)\\
\(xt\) (intermediate) & 17 & \(0\) & choose a cocycle \(h_{17}\)\\
\(sr,yr\) & 19 & \(-sxt,-yxt\) &
  \(-\widetilde s h_{17},-\widetilde y h_{17}\)\\
\(xr\) & 21 & \(0\) & choose a cocycle \(\xi\)\\
\(ar\) & 23 & \(xyr-axt\) &
  \(\xi\widetilde y-
    \widetilde a\widetilde x\widetilde t\)\\
\(t^2\) (intermediate) & 24 & \(0\) & choose a cocycle \(\theta\)\\
\(tr\) & 28 & \(xt^2\) & \(\widetilde x\theta\)\\
\(r^2\) & 32 & \(2xtr\) & \(2\widetilde x\tau_r\)\\
\(r^3\) & 48 & \(0\) & choose a cocycle \(R_3\)\\
\bottomrule
\end{tabularx}
\end{center}

We now verify the schedule and its induction invariant.  Start with exact
cocycle lifts
\(g(\widetilde s)=s\), \(g(\widetilde y)=y\), and
\(g(\widetilde x)=x\).  The product
\(\widetilde x\widetilde y\) is closed and maps to \(da=xy\); exact lifting
gives \(\widetilde a\in N_{17}^7\) with
\[
 g(\widetilde a)=a,
 \qquad d\widetilde a=\widetilde x\widetilde y.
\]
The product \(\widetilde a\widetilde y\widetilde s\) is closed because
\(\widetilde y^2=0\), and it maps to \(dt=ays\).  Hence there is
\(\widetilde t\in N_{17}^{12}\) satisfying
\[
 g(\widetilde t)=t,
 \qquad d\widetilde t
       =\widetilde a\widetilde y\widetilde s.
\]

The element \(xt\) is closed in \(B^\circ\), so choose
\(h_{17}\in Z^{17}(N_{17})\) with \(g(h_{17})=xt\).  Exact lifting of the
next two rows gives \(r_s,r_y\in N_{17}^{19}\) with
\begin{equation}\label{stab:eq:punctured-rs-ry-lifts}
\begin{aligned}
 g(r_s)&=sr,& dr_s&=-\widetilde s h_{17},\\
 g(r_y)&=yr,& dr_y&=-\widetilde y h_{17}.
\end{aligned}
\end{equation}
Choose \(\xi\in Z^{21}(N_{17})\) with \(g(\xi)=xr\).  The element
\[
                  c_{24}=\widetilde a\widetilde x\widetilde t
\]
is closed: its two possible differential terms contain, respectively,
\(\widetilde x^2\) and \(\widetilde a^2\).  Therefore
\(\xi\widetilde y-c_{24}\) is closed and maps to
\(xyr-axt=d(ar)\).  We obtain \(r_a\in N_{17}^{23}\) with
\begin{equation}\label{stab:eq:punctured-ra-lift}
 g(r_a)=ar,
 \qquad dr_a=\xi\widetilde y-c_{24}.
\end{equation}

The cycle \(t^2\) has differential of augmentation length four and is
therefore closed in the three-short algebra.  Choose
\(\theta\in Z^{24}(N_{17})\) with \(g(\theta)=t^2\).  The last nonclosed
rows give \(\tau_r\in N_{17}^{28}\) and \(R_2\in N_{17}^{32}\) such that
\begin{equation}\label{stab:eq:punctured-tr-r2-lifts}
\begin{aligned}
 g(\tau_r)&=tr,& d\tau_r&=\widetilde x\theta,\\
 g(R_2)&=r^2,& dR_2&=2\widetilde x\tau_r.
\end{aligned}
\end{equation}
The second right-hand side is closed because
\(d(\widetilde x\tau_r)=-\widetilde x^2\theta=0\).  Finally,
\(d(r^3)=3r^2xt=0\) by three-shortness, so exact cocycle lifting gives
\(R_3\in Z^{48}(N_{17})\) with \(g(R_3)=r^3\).

There is no circularity in these applications.  Before the row of degree
\(q\), every indecomposable direction of \eqref{stab:eq:QB-circ} in degree
less than \(q\) has already been lifted exactly.  Thus \(Q(g)\) is
surjective through degree \(q-1\), and the last assertion of
Lemma~\ref{stab:lem:degreewise-exact-lifting} makes \(g\) surjective through
that degree.  This is precisely the hypothesis needed for the next use of
part~(i) or part~(ii).  Induction through the table produces preimages of all
basis elements in \eqref{stab:eq:QB-circ}.  Hence \(Q(g)\) is surjective, and
Lemma~\ref{stab:lem:Q-surj} implies that \(g\) is surjective.
\end{proof}

\begin{theorem}[Explicit one-sphere stabilization of the separating example]
\label{stab:thm:counterexample-S17}
For the separating space \(|M|\) of Part~\ref{part:divergence},
\[
\begin{aligned}
 \catz(|M|)
 &=\catz(|M|\vee S^{17}_{\Q})
  =\Clz(|M|\vee S^{17}_{\Q})\\
 &=\nilh(|M|\vee S^{17}_{\Q})
  =\Hnilz(|M|\vee S^{17}_{\Q})=3,
\end{aligned}
\]
whereas
\[
                         \Hnilz(|M|)=4.
\]
\end{theorem}

\begin{proof}
Proposition~\ref{stab:prop:punctured-minimal-surjective} gives a surjective
quasi-isomorphism from the minimal model of
\(|M|\vee S^{17}_{\Q}\) to the three-short CDGA \(B^\circ\).
Lemma~\ref{div:lem:surjective-characterization} therefore gives
\[
                  \Hnilz(|M|\vee S^{17}_{\Q})\leq3.
\]
The inclusion and collapse exhibit \(|M|\) as a rational
homotopy retract of the wedge.  By
Propositions~\ref{stab:prop:separating-category} and
\ref{stab:prop:cat-lower},
\[
3=\catz(|M|)
 \leq\catz(|M|\vee S^{17}_{\Q})
 \leq\Hnilz(|M|\vee S^{17}_{\Q})\leq3.
\]
The comparison chain
\(\catz\leq\Clz=\nilh\leq\Hnilz\) gives the intervening equalities.
Finally, Theorem~\ref{app:thm:Hnil-exact} gives the unstabilized value.
\end{proof}

\begin{remark}\label{stab:rem:S17-not-minimal-degree}
The hidden contractible pair \(r_{16}\mapsto xt_{17}\) explains the degree:
removing its primitive exposes the class \([xt]\) that is modeled by the
added sphere.  The theorem proves that the minimum number of sphere summands
needed for this example is one, because the empty wedge leaves
\(\Hnilz(|M|)=4>3=\catz(|M|)\), while \(S^{17}_{\Q}\) suffices.  It does not claim
that \(17\) is the smallest dimension among all single spheres that might
work.
\end{remark}

\section{Scope of the general construction}
\label{stab:sec:scope}

The admissible class \(\mathcal S_{\mathrm{lf}}\) allows infinitely many
sphere summands in total, while requiring only finitely many in each
degree.  This is the natural class produced by the rational splitting of
\(\Sigma\Omega F_\rho\).  If its reduced homology is finite
dimensional, then the corresponding wedge is finite and the same proof gives
the finite-wedge version.  Without such an additional finiteness
hypothesis, the present argument does not assert that a finite total
number of spheres suffices.

The open-disk hypothesis is essential in the algebraic core.  In a
general short CDGA, the vanishing of the differential on indecomposables
does not force a minimal-model morphism to be surjective.  For example,
let
\begin{equation}\label{stab:eq:abstract-counterexample}
 A=\frac{\Lambda(b_2,c_3,a_4)}
 {(b^3,a^2,ab,ac,b^2c)},
 \qquad dc=b^2,\quad da=bc.
\end{equation}
Then \(Q(d_A)=0\).  Indeed, the six elements

\[
                         1,b,b^2,c,bc,a
\]

form a basis of \(A\), and the only nonzero differentials on this basis
are \(dc=b^2\) and \(da=bc\).  Hence

\[
                         H(A)=\Q\{1,[b]\}.
\]

The relations also make every product of three positive elements zero,
so \((A^+)^3=0\).  Its minimal model is that of \(S^2\),

\[
          \bigl(\Lambda(x_2,y_3),dy=x^2\bigr)
             \longrightarrow A,
          \qquad x\longmapsto b,\quad y\longmapsto c.
\]

The displayed morphism is a quasi-isomorphism by the cohomology
calculation.  Any quasi-isomorphism from this minimal model must
send \(x\) to \(\lambda b\), with \(\lambda\neq0\), and the chain-map
identity forces \(y\) to map to \(\lambda^2c\).  After rescaling, every
such map therefore has
\[
 x\longmapsto b,\qquad y\longmapsto c,
\]
and its image is contained in the subalgebra generated by \(b\) and \(c\),
which does not contain the indecomposable \(a\).  What rules out this defect in
the geometric construction is not shortness alone, but the exact
open-disk form \eqref{stab:eq:open-C}: the Euler contraction creates the
binary term \(I_2(v,D\omega)=(-1)^{|v|+1}v\omega\), and the marked Harrison
cycle proves that it survives in the derived indecomposables.

\begin{keyidea}{Part V summary}
The optimal LS root produces an exact open-disk model for the
fibre--cofibre construction.  Harrison--Andr\'e--Quillen comparison makes
the minimal-model map surjective on its mixed indecomposables and ultimately
yields a short quotient of the fixed model.  Together with the category lower
bound, this
proves \(\sHnil(X)=\catz(X)\); the spherical factor is degreewise finite,
although it need not be finite in total.  For the separating space itself,
puncturing the hidden disk \(r_{16}\mapsto xt_{17}\) proves the sharper
finite statement
\(\Hnilz(|M|\vee S^{17}_{\Q})
=\catz(|M|)=3<4=\Hnilz(|M|)\).  Part~\ref{part:context} compares the strict
separation with neighboring invariants and records the remaining problems,
including whether an example can genuinely require two sphere summands.
\end{keyidea}

\part{Context and Open Problems}\label{part:context}

\chapter{Context and Open Problems}\label{chap:precedents}
\begin{chapterguide}[title={Chapter guide}]
With both divergence and stabilization established, this final part locates
the minimal-model separation among neighboring presentation-dependent
phenomena and records the questions left open by the tower, product, and
stabilization theories.  In particular, it asks how large the strictification
gap can be and whether homology nilpotency satisfies a rational Ganea formula
or a full product formula.
\end{chapterguide}
\index{model dependence}
\index{sectional category}
\index{topological complexity}
\index{multiplication kernel}

\section{Earlier model-dependent separations}\label{div:sec:precedents}

\begin{remark}[Nearby presentation-dependent examples]
Three examples of Carrasquel-Vera clarify why the fixed-minimal-model
hypothesis matters, but require only a brief comparison here.  In
\cite[Example~10]{CarrasquelVera2015}, a non-Sullivan augmented CDGA has
\(\nilh=1\), whereas its presentation-dependent augmentation ideal satisfies
\(\Hnil=2\): its square is generated by a nonzero cohomology class and its
cube vanishes.  The separation disappears in the minimal Sullivan model.

In \cite[Example~18]{CarrasquelVera2015}, the source is minimal Sullivan, but
the relevant ideal is the kernel \(K\) of a morphism rather than the
augmentation ideal.  There
\[
 \operatorname{secat}=1
 <\operatorname{sc}=\operatorname{Hsc}=\operatorname{nil}K=2,
\]
and the homology nilpotency of this kernel is also two.
Finally, the non-Sullivan multiplication model in
\cite[Example~5.2]{Carrasquel2017} has the rational homotopy type of
\(S^3_{\Q}\) and gives
\[
 \operatorname{secat}=\operatorname{TC}(S^3_{\Q})=1,
 \qquad
 \operatorname{sc}=\operatorname{Hsc}=3,
\]
and the homology nilpotency of the multiplication kernel is three.
Here again \(K\) is a multiplication kernel, and the larger value disappears
in the minimal model.  These examples demonstrate dependence
on a presentation or on the chosen kernel; none compares \(\nilh\) and
\(\Hnil\) for the augmentation ideal of a fixed minimal Sullivan model.
\end{remark}

Recent work of Parent and Tanr\'e constructs, for every \(k\geq3\), a
rational space \(X_k\) with
\[
                  \catz(X_k)=k,
                  \qquad \Clz(X_k)=k+1
\]
\cite{ParentTanre2026}.  This separates rational LS category from cone
length---the first inequality in the comparison chain---and does not concern
the fixed-minimal-model separation \(\Clz<\Hnilz\) proved here.

Over a principal ideal domain \(R\) containing \(1/2\), Benzaki and Rami
extend sectional-category and topological-complexity constructions to
\((r,\rho(R))\)-mild spaces and investigate their relation to homology
nilpotency \cite{BenzakiRami2025}.  Their coefficient-general setting
concerns invariants attached to morphisms and product models; it does not
provide the minimal augmentation-ideal separation studied in this memoir.

For broader context, the geometric theory of LS category and cone length is
surveyed by Cornea, Lupton, Oprea, and Tanr\'e
\cite{CorneaLuptonOpreaTanre2003}; the rational product and
Poincar\'e-duality results surrounding these invariants are developed in
\cite{FelixHalperinLemaire1998,CorneaFelixLemaire1998}, with Hess's
module-to-rational-category theorem providing a central bridge
\cite{Hess1991}.  The join models behind module sectional category originate
in \cite{FernandezSuarezGhienneKahlVandembroucq2006} and are sharpened in the
Poincar\'e-duality and product settings in
\cite{CarrasquelVeraKahlVandembroucq2016,CarrasquelVeraParentVandembroucq2018}.
For the operadic background and the Harrison conventions used in
Chapter~\ref{chap:aq-toolkit}, compare Fresse's bar--cobar treatment
\cite{Fresse2009} and Loday's account of the Eulerian decomposition
\cite{Loday1998}.  These references locate the present comparison within the
geometric, module-theoretic, and operadic literature; they do not alter the
minimal-model distinction made above.

To the best of the author's knowledge, the construction above is the first
separation established for the augmentation ideal of the minimal Sullivan
model itself.  This claim is restricted to the minimal-model setting and does
not assert priority for arbitrary CDGA presentations or non-augmentation
kernels.

\begin{remark}[Scope and dated priority search]
\label{context:rem:priority-search}
The priority statement above is deliberately narrower than the assertion
that no model-dependent separation has appeared before.  For this revision, searches dated 1 September 2026 were run in arXiv and
zbMATH Open, supplemented by publicly indexed MathSciNet-facing records,
using combinations of the terms
\emph{homology nilpotency}, \emph{homotopical nil-length}, \emph{minimal
Sullivan model}, and \emph{augmentation ideal}.  The searches recovered the
precedents discussed in this section but no earlier example simultaneously
satisfying all three restrictions: the algebra is the fixed minimal Sullivan
model, the ideal is its augmentation ideal, and
\(\nilh(M)<\Hnil(M)\).  A bibliographic search cannot prove absolute
priority; this dated protocol records the scope of the search and explains
why the claim remains qualified by ``to the best of the author's
knowledge.''
\end{remark}
\index{priority statement}
\index{augmentation ideal!priority search}

\section{What remains open}
The example proves that the passage from an arbitrary short model to a
short quotient of the minimal Sullivan model fails at augmentation
length three.  It also shows exactly why the arguments valid at lengths
one and two cannot be extended formally.  In the cube-zero case, the top
multiplicative layer is a square-zero ideal of cycles.  In the
fourth-power-zero case, a forced primitive can still multiply with a
positive element and create an essential class.
For the separating model, Appendix~\ref{app:exact-Hnil} shows that this
obstruction is sharp: \(\Hnil(M)=4\).

The tower theory isolates a separate question.
Proposition~\ref{tate:prop:asymptotic-implication} proves that an asymptotically Tate
level-\(n\) tower implies \(\Hnil(M)\leq n\).  The present results leave open
whether an
acyclic \(n\)-short quotient can always be used to select compatible
factorizations and retractions producing such a tower.  This is a global
strict-coherence problem, not a consequence of the finite-packet Tate
criterion.

The normalized classical LS category and Ganea's strong LS category satisfy
\begin{equation}\label{context:eq:cat-strong-cat-gap}
             \operatorname{cat}(X)
             \leq \operatorname{Cat}(X)
             \leq \operatorname{cat}(X)+1
\end{equation}
for every connected space of CW type.  Strong category was introduced by
Ganea \cite{Ganea1967}; the upper estimate in
\eqref{context:eq:cat-strong-cat-gap} is the Ganea--Takens theorem, with a
homotopy-categorical proof recorded in
\cite[Section~5]{Takens1970}.  Cornea subsequently identified strong LS
category with cone length \cite{Cornea1995}.  This classical one-step
phenomenon suggests a quantitative question for the flexible-versus-rigid
comparison developed here.

For a simply connected minimal Sullivan algebra \(M\) of finite type with
\(\nilh(M)<\infty\), define its rigid strictification gap by
\begin{equation}\label{context:eq:rigid-gap}
 \delta_{\mathrm{rig}}(M)=
 \begin{cases}
  \Hnil(M)-\nilh(M),&\Hnil(M)<\infty,\\
  \infty,&\Hnil(M)=\infty.
 \end{cases}
\end{equation}
The separating algebra of this memoir has
\(\delta_{\mathrm{rig}}(M)=1\).
\index{strong LS category}
\index{Ganea--Takens inequality}
\index{strictification gap}
\index{homology nilpotency!gap from homotopical nil-length}

\begin{problem}[A universal one-step bound or unbounded gaps]
\label{context:prob:Hnil-unit-gap}
Let \(M\) range over simply connected minimal Sullivan algebras of finite
type with \(\nilh(M)<\infty\).  Is the analogue of
\eqref{context:eq:cat-strong-cat-gap}
\[
                  \nilh(M)\leq\Hnil(M)\leq\nilh(M)+1
\]
always valid?  If not, are the finite values of
\(\delta_{\mathrm{rig}}(M)\) unbounded?  More strongly, can one have
\[
                   \nilh(M)<\infty,
                   \qquad \Hnil(M)=\infty?
\]
Thus the alternatives include a universal one-step strictification theorem,
arbitrarily large finite gaps, and a finite-versus-infinite separation.
\end{problem}

The product behavior provides a second test.  For simply connected rational
spaces \(X\) and \(Y\) of finite type, F\'elix, Halperin, and Lemaire proved
\begin{equation}\label{context:eq:rational-category-product}
             \catz(X\times Y)=\catz(X)+\catz(Y)
\end{equation}
\cite[Theorem~1]{FelixHalperinLemaire1998}; its sphere-factor case is the
rational Ganea formula proved earlier by Hess \cite{Hess1991}.  If \(M\) and
\(N\) are the minimal Sullivan models of \(X\) and \(Y\), then
\(M\otimes N\) is the minimal model of \(X\times Y\).

\begin{proposition}[Product subadditivity of homology nilpotency]
\label{context:prop:Hnil-product-subadditivity}
For simply connected minimal Sullivan algebras \(M\) and \(N\) of finite
type,
\[
                 \Hnil(M\otimes N)
                 \leq \Hnil(M)+\Hnil(N)
\]
in the extended nonnegative integers.
\end{proposition}

\begin{proof}
Only the case \(p=\Hnil(M)<\infty\) and
\(q=\Hnil(N)<\infty\) requires proof.  Choose acyclic differential ideals
\(I\triangleleft M\) and \(J\triangleleft N\) such that
\[
                  (M^+)^{p+1}\subseteq I,
                  \qquad (N^+)^{q+1}\subseteq J.
\]
The tensor product of the quotient maps is a surjective quasi-isomorphism
\[
 M\otimes N\longrightarrow (M/I)\otimes(N/J),
\]
because tensoring over \(\Q\) preserves quasi-isomorphisms.  Its kernel
\[
                 K=I\otimes N+M\otimes J
\]
is therefore an acyclic differential ideal.  The augmentation ideal of the
target is generated by the images of \((M/I)^+\otimes1\) and
\(1\otimes(N/J)^+\).  Every product of \(p+q+1\) such generators contains
either \(p+1\) factors from the first ideal or \(q+1\) factors from the
second.  Hence the target is \((p+q)\)-short, so
\(((M\otimes N)^+)^{p+q+1}\subseteq K\), proving the claim.
\end{proof}

Write \(\mathcal S_n\) for the minimal Sullivan model of the rational sphere
\(S^n_{\Q}\), \(n\geq2\).  One has \(\Hnil(\mathcal S_n)=1\): this is
immediate for odd \(n\), while for even \(n\) the standard model
\((\Lambda(x_n,y_{2n-1}),dy=x^2)\) has the acyclic differential ideal
\((y,x^2)\) containing the square of its augmentation ideal.  The preceding
proposition therefore gives
\[
            \Hnil(M\otimes\mathcal S_n)\leq\Hnil(M)+1.
\]

\begin{problem}[A rational Ganea formula and product additivity]
\label{context:prob:Hnil-product-additivity}
For every simply connected minimal Sullivan algebra \(M\) of finite type and
every \(n\geq2\), does one have
\begin{equation}\label{context:eq:Hnil-Ganea-question}
              \Hnil(M\otimes\mathcal S_n)=\Hnil(M)+1?
\end{equation}
Equivalently, for the associated rational space \(X=|M|\), is
\[
              \Hnilz(X\times S^n_{\Q})=\Hnilz(X)+1?
\]
More generally, is the subadditivity bound always sharp; that is, do all
simply connected minimal Sullivan algebras \(M,N\) of finite type satisfy
\begin{equation}\label{context:eq:Hnil-product-question}
              \Hnil(M\otimes N)=\Hnil(M)+\Hnil(N)?
\end{equation}
The sphere-factor identity \eqref{context:eq:Hnil-Ganea-question} is the
first test of the full product formula
\eqref{context:eq:Hnil-product-question}.
\end{problem}
\index{Ganea formula!homology nilpotency}
\index{homology nilpotency!product subadditivity}
\index{homology nilpotency!product additivity problem}
\index{tensor product!homology nilpotency}

Section~\ref{stab:sec:scope} explains why the general stabilization theorem
does not ensure a finite total wedge.  To formulate the finite question more
sharply, for a simply connected rational
space \(Y\) of finite type with \(\catz(Y)<\infty\), set
\begin{equation}\label{stab:eq:sphere-stabilization-number}
 \sigma_{\Hnil}(Y)
 =\inf\left\{k\geq0\ \middle|\
 \begin{array}{l}
 \text{there exist }n_1,\ldots,n_k\geq2\text{ such that}\\[-2pt]
 \displaystyle
 \Hnilz\!\left(Y\vee\bigvee_{i=1}^{k}S^{n_i}_{\Q}\right)=\catz(Y)
 \end{array}\right\},
\end{equation}
where \(k=0\) denotes the empty wedge and the infimum of the empty set is
\(\infty\).  Theorem~\ref{stab:thm:counterexample-S17} and the strict
unstabilized inequality give
\[
                         \sigma_{\Hnil}(|M|)=1.
\]
\index{sphere!minimum number for stabilization}

\begin{problem}[A genuinely two-sphere stabilization]
\label{stab:prob:two-sphere-minimum}
Construct a simply connected rational space \(Y\) of finite type with
\(1\leq\catz(Y)<\infty\) such that
\[
 \Hnilz(Y\vee S^m_{\Q})>\catz(Y)
 \qquad\text{for every }m\geq2,
\]
but, for some \(p,q\geq2\), not necessarily distinct,
\[
             \Hnilz(Y\vee S^p_{\Q}\vee S^q_{\Q})=\catz(Y).
\]
Equivalently, construct \(Y\) with \(\sigma_{\Hnil}(Y)=2\).  More generally,
determine whether every prescribed positive integer can occur as
\(\sigma_{\Hnil}(Y)\).
\end{problem}

\begin{keyidea}{Concluding perspective}
The flexible-versus-rigid comparison reveals two complementary phenomena:
separation occurs inside a fixed minimal model already at length three, yet
spherical stabilization restores the rational-category value.  Retractive
towers locate the defect between these conclusions.  The remaining problems
ask how large that defect can become, how it behaves under products, and how
economically spherical directions can remove it.
\end{keyidea}

\appendix

\chapter{Harrison Suspension and Sign Conventions}
\label{app:harrison-conventions}
\begin{chapterguide}[title={Purpose of this appendix}]
The objects and their categories are typed in
Chapter~\ref{chap:aq-toolkit}.  This appendix is the corresponding sign
dictionary: it distinguishes the bar and Lie-theoretic shifts from the other
suspension symbols and records the bar differential, the first Eulerian
idempotent, the \(C_\infty\) Taylor signs, the normalization of the Koszul
twisting morphism, and the unit coefficient of the marked binary term.  The
coalgebra structure on the Harrison image model is always transported from
the shuffle quotient; it is not obtained by treating the Eulerian idempotent
as a coalgebra morphism on the full tensor bar.
\end{chapterguide}
\index{Harrison complex!sign conventions}
\index{suspension!bar}
\index{suspension!Lie-theoretic}
\index{Eulerian idempotent}
\index{Koszul twisting morphism}

\section{Gradings and the suspension symbols}
\label{app:sec:grading-suspensions}

All differential graded algebras are cohomologically graded, and their
differentials have degree \(+1\).  The bar shift \(\bs\) has degree
\(-1\):
\begin{equation}\label{app:eq:bar-shift-degree}
 |\bs a|=|a|-1,\qquad |\bs^{-1}|=+1.
\end{equation}
For homogeneous maps and tensors, we use the Koszul evaluation rule
\[
 (f_1\otimes\cdots\otimes f_r)(x_1\otimes\cdots\otimes x_r)
 =
 (-1)^{\sum_{i<j}|f_j||x_i|}
 f_1(x_1)\otimes\cdots\otimes f_r(x_r).
\]
It follows in particular that
\begin{equation}\label{app:eq:desuspension-tensor-sign}
 (\bs^{-1})^{\otimes r}
 (\bs a_1\mid\cdots\mid\bs a_r)
 =
 (-1)^{\epsilon_r(a_1,\ldots,a_r)}
 (a_1\otimes\cdots\otimes a_r),
 \quad
 \epsilon_r=\sum_{i=1}^r(r-i)(|a_i|-1).
\end{equation}

The Lie-theoretic shift \(\ls\) is a different symbol and has degree
\(+1\).  If \(N=(\Lambda Z,d)\) is a minimal Sullivan model, then
\[
 L_{N,k}=(Z^{k+1})^\sharp,\qquad
 \ls L_N=Z^\sharp.
\]
Thus, if \(x\in L_N\) is dual to \(z\in Z\), then
\(|x|=|z|-1\) and \(\ls x=z^\sharp\).  With the Koszul pairing on
\(\Lambda^2Z\) induced by \(Z^\sharp\), the bracket convention is
\begin{equation}\label{app:eq:Sullivan-Quillen-bracket-sign}
 \bigl\langle z;\ls[x,y]\bigr\rangle
 =
 (-1)^{|y|+1}
 \bigl\langle d_2z;\ls x,\ls y\bigr\rangle ,
\end{equation}
where \(d_2:Z\to\Lambda^2Z\) is the quadratic part.  This is the
Sullivan--Quillen convention of
\cite[Section~21(e) and Theorem~21.6]{FHT2001}.  The comparison used in
Chapter~\ref{chap:aq-toolkit} needs only the underlying graded-vector-space
identity \(\ls L_N=Z^\sharp\).

The symbol \(\tau\) used for the square-zero suspension factor is a
cohomological degree shift: for \(\omega\in T^+\),
\[
                  |\tau\omega|=|\omega|+1=|D\omega|.
\]
The symbol \(\Sigma\) denotes spatial suspension.  Neither \(\tau\) nor
\(\Sigma\) is identified with the Lie-theoretic shift \(\ls\).

\section{The reduced bar and Harrison models}
\label{app:sec:bar-harrison-conventions}

For an augmented CDGA \(A\), the reduced bar object is
\[
 \overline B(A)=T^c(\bs A^+).
\]
Its coderivation \(b=b_1+b_2\) is determined by
\begin{equation}\label{app:eq:bar-corestriction-signs}
 \pi_1b_1(\bs a)=-\bs(d_Aa),\qquad
 \pi_1b_2(\bs a\mid\bs b)=(-1)^{|a|}\bs(ab).
\end{equation}
The signs on longer tensors are obtained from the Koszul rule and the
coderivation identity.  In weight two the shuffle product is
\[
 (\bs a)\mathbin{\mathrm{sh}}(\bs b)
 =\bs a\mid\bs b
  +(-1)^{(|a|-1)(|b|-1)}\bs b\mid\bs a .
\]

Let \(\star\) be convolution for deconcatenation and shuffle, and put
\[
 e^{(1)}
 =\log^\star(\id)
 =\sum_{k\geq1}\frac{(-1)^{k-1}}{k}
   (\id-\eta\varepsilon)^{\star k}.
\]
On each bar weight the sum is finite.  The Harrison image model is
\(\Harr(A)=\im(e^{(1)})\).  Its weight-two idempotent is
\begin{equation}\label{app:eq:eulerian-weight-two}
 e^{(1)}_2(\bs a\mid\bs b)
 =\frac12\left(
 \bs a\mid\bs b
 -(-1)^{(|a|-1)(|b|-1)}\bs b\mid\bs a
 \right).
\end{equation}
In characteristic zero, Barr's splitting identifies this image, as a
chain complex, with the quotient by shuffle decomposables
\cite[Proposition~2.5]{Barr1968}.  For the Harrison quotient and its
classical Andr\'e--Quillen interpretation, see
\cite[Section~4.2.10 and Proposition~4.2.11]{Loday1998}; for the Eulerian
splitting, see \cite[Propositions~4.5.9--4.5.13]{Loday1998}; and for the
commutative differential graded setting, see
\cite[Section~5.4.8]{Loday1998}.
If \(q_A:\overline B(A)\twoheadrightarrow\Harr^{\mathrm c}(A)\) is the
quotient map, its inverse on the first Eulerian summand is
\[
 s_A:\Harr^{\mathrm c}(A)\longrightarrow\Harr(A),
 \qquad s_A(q_Ax)=e_A^{(1)}x,
\]
and
\[
 q_As_A=\id,
 \qquad s_Aq_A=e_A^{(1)}.
\]
Thus \(q_A|_{\Harr(A)}\) and \(s_A\) are inverse chain isomorphisms.
Throughout, the Lie coalgebra structure on the image model is transported
through this explicit pair.

\section{Cooperad and twisting normalization}
\label{app:sec:cooperad-normalization}

The symbol \(\Comash\) denotes the Koszul-dual \emph{cooperad} of
\(\Com\).  It is not the standard cocommutative cooperad
\(\mathrm{coCom}\), and it is not the Koszul-dual operad
\(\Com^!=\mathrm{Lie}\).  Operadic d\'ecalage identifies
conilpotent \(\Comash\)-coalgebras with the correspondingly shifted
conilpotent dg Lie coalgebras
\cite[Sections~7.2.1--7.2.3, Proposition~13.1.1, and
Section~13.1.8]{LodayVallette2012}.  We use the canonical
Koszul twisting morphism
\[
 \kappa:\Comash\longrightarrow\Com
\]
with arity-two coefficient \(+1\) after this d\'ecalage.  Compatibly,
the bar--cobar counit is induced by the universal twisting cochain
\[
 t_A(\bs^{-1}\bs a)=a,
 \qquad t_A=0\quad\text{on Harrison weights at least two}.
\]
Thus its weight-one part sends \(\bs^{-1}\bs a\) to \(a\), with no
additional scalar.  This normalization is used in
\[
 R_A=\Omega_\kappa(\bs^{-1}\Harr(A)).
\]

\section{\texorpdfstring{\(C_\infty\) Taylor coefficients}{C-infinity Taylor coefficients}}
\label{app:sec:Cinfinity-conventions}
\index{C-infinity algebra@$C_\infty$-algebra!Taylor coefficient}

For a strictly unital, augmentation-preserving \(C_\infty\)-morphism
\(J:B\rightsquigarrow A\), the unshifted Taylor coefficient
\[
 J_r:B^{\otimes r}\longrightarrow A
\]
has degree \(1-r\).  Its shifted corestriction is
\begin{equation}\label{app:eq:shifted-Taylor-map}
 \overline J_r
 =\bs J_r(\bs^{-1})^{\otimes r},
\end{equation}
so \eqref{app:eq:desuspension-tensor-sign} gives
\begin{equation}\label{app:eq:shifted-Taylor-evaluation}
 \overline J_r(\bs a_1\mid\cdots\mid\bs a_r)
 =
 (-1)^{\epsilon_r(a_1,\ldots,a_r)}
 \bs J_r(a_1,\ldots,a_r).
\end{equation}
The corestrictions annihilate shuffle decomposables.  For \(r=2\), this
condition and \eqref{app:eq:shifted-Taylor-evaluation} are equivalent to
the unshifted graded symmetry
\begin{equation}\label{app:eq:J2-graded-symmetry}
 J_2(a,b)=(-1)^{|a||b|}J_2(b,a).
\end{equation}

Let \(\Harr^{\mathrm c}(J)\) be the morphism induced on the quotient
Harrison Lie coalgebras.  On the image models, the transported morphism is
\begin{equation}\label{app:eq:Harrison-transported-Cinfinity-map}
 \widetilde B(J)
 =s_A\circ\Harr^{\mathrm c}(J)\circ
 \bigl(q_B|_{\Harr(B)}\bigr).
\end{equation}
The induced quotient morphism satisfies
\[
 \Harr^{\mathrm c}(J)q_B=q_A\overline B(J).
\]
Thus, for \(x\in\Harr(B)_{(r)}\),
\[
 \pi_1\widetilde B(J)(x)
 =\pi_1e_A^{(1)}\overline B(J)(x)
 =\pi_1\overline B(J)(x)
 =\overline J_r(x),
\]
where \(\pi_1e_A^{(1)}=\pi_1\) because \(e_A^{(1)}\) is the identity in
Harrison weight one.  This is the transport convention used by the Harrison
cobar functor in Chapter~\ref{chap:aq-toolkit}.

Let \(F_J:R_B\to A\) be the strict bar--cobar morphism induced by
\(\widetilde B(J)\), as in \eqref{stab:eq:FJ-strict}.  Combining
\eqref{app:eq:eulerian-weight-two},
\eqref{app:eq:shifted-Taylor-evaluation}, and the twisting normalization
gives the exact binary formula
\begin{equation}\label{app:eq:binary-unit-coefficient}
 QF_J\!\left(
 \bs^{-1}e^{(1)}_2(\bs a\mid\bs b)
 \right)
 =
 (-1)^{|a|-1}[J_2(a,b)]_{Q(A)}.
\end{equation}
Thus the coefficient is a unit in \(\mathbb Z\), not an unspecified
nonzero rational number.  Nonvanishing requires the explicit hypothesis
\([J_2(a,b)]_{Q(A)}\ne0\): a nonzero but decomposable value of \(J_2\)
vanishes in \(Q(A)\).

\section{Transfer and the marked mixed coefficient}
\label{app:sec:transfer-marked-sign}

For a strong deformation retract
\[
 \left(B\mathrel{\substack{\xrightarrow{i}\\[-0.4ex]
 \xleftarrow[p]{} }}C,H\right),
 \qquad
 DH+HD=\id_C-ip,
\]
the convention differs by a sign from that of Cheng--Getzler, who write
\(gf-\id=dH_{\mathrm{CG}}+H_{\mathrm{CG}}d\): here
\(H_{\mathrm{CG}}=-H\).  The bar differential is
\(d_{\bs C}=-\bs D\bs^{-1}\).  Hence the shifted
homotopy \(\bs H\bs^{-1}\) satisfies
\[
 d_{\bs C}(\bs H\bs^{-1})
 +(\bs H\bs^{-1})d_{\bs C}
 =\bs(ip-\id_C)\bs^{-1}.
\]
Let \(I=(I_r):B\rightsquigarrow C\) be the transferred
\(C_\infty\)-morphism.  With the rooted-tree convention of
\cite[Theorems~6, 10, and~12]{ChengGetzler2008}, desuspension therefore
gives
\begin{equation}\label{app:eq:transfer-binary-sign}
 I_2(a,b)=-H\bigl(i(a)i(b)-i(ab)\bigr).
\end{equation}

In the open-disk contraction of Chapter~\ref{chap:aq-toolkit}, let
\(b_\omega=D\omega\), where \(\omega\in T^+\).  For \(v\in V\),
\[
 H(vb_\omega)=(-1)^{|v|}v\omega,\qquad
 I_2(v,b_\omega)=(-1)^{|v|+1}v\omega.
\]
The marked Harrison cycle is ordered as
\(e^{(1)}_2(\bs b_\omega\mid\bs v)\).  Using
\eqref{app:eq:J2-graded-symmetry} and
\eqref{app:eq:binary-unit-coefficient}, its \(v\omega\)-coordinate is
\begin{equation}\label{app:eq:marked-binary-unit}
 (-1)^{|\omega|(|v|+1)+1}\in\{+1,-1\}.
\end{equation}
This is the coefficient denoted \(c_{v,\omega}\) in
\eqref{stab:eq:diagonal-evaluation}.

\chapter{The Exact Value for the Separating Model}
\label{app:exact-Hnil}
\begin{chapterguide}[title={Purpose of this appendix}]
This appendix gives the degree-by-degree completion argument that sharpens
the strict lower bound of Chapter~\ref{chap:strict-separation} to the exact
value \(\Hnil(M)=4\).  The proof is separated from the obstruction argument
so that the high-degree length estimate and the ideal construction can be
checked on their own.
\end{chapterguide}
\index{homology nilpotency!exact value for the separating model}
\index{acyclic ideal!degreewise completion}

Recall the fixed quasi-isomorphism
\[
 \varphi:M=(\Lambda(s,x,y,a,t)\otimes\Lambda W,d)
        \xrightarrow{\simeq}B,
 \qquad W^{\leq18}=0,
\]
from Corollary~\ref{div:cor:rigid-prefix}.  The target is the
three-short CDGA of Chapter~\ref{div:sec:witness}; in particular,
\((B^+)^4=0\).

\begin{lemma}[High-degree cycles have maximal length]
\label{app:lem:high-cycles-cubic}
For every \(j\geq29\),
\[
                         Z^j(B)\subseteq(B^+)^3.
\]
Consequently, \(Z^{\geq29}(B)B^+=0\).
\end{lemma}

\begin{proof}
Use the direct decomposition of \(B\) by augmentation length.  A linear
monomial has degree at most \(16\).  Among the quadratic monomials, the only
one of degree at least \(29\) is \(r^2\in B^{32}\): the preceding
quadratic monomial is \(tr\), of degree \(28\).  Moreover,
\[
                              d(r^2)=2xtr\neq0.
\]
Every cubic monomial is a cycle, since its differential has augmentation
length at least four.  Hence the quadratic component of a cycle in degree at
least \(29\) must vanish, and the cycle is cubic.  Multiplication by a
positive element then gives an element of \((B^+)^4=0\).
\end{proof}

\begin{theorem}[Exact homology nilpotency of the separating model]
\label{app:thm:Hnil-exact}
For the minimal Sullivan algebra \(M\) of
Corollary~\ref{div:cor:rigid-prefix},
\[
                              \Hnil(M)=4.
\]
Equivalently, there is an acyclic differential ideal \(J\lhd M\) containing
\(\mM^5\).
\end{theorem}

\begin{proof}
Proposition~\ref{div:prop:Hnil-lower} gives \(\Hnil(M)>3\).  We prove the
matching upper bound.

Put \(I_{29}=\mM^5\).  Since \(s,x,y,a\) have odd degree and square to
zero, a nonzero word of length five in the rigid prefix \(s,x,y,a,t\) has
degree at least
\[
                         3+3+5+7+12=30.
\]
A word involving a generator of \(W\) has larger degree because
\(W^{\leq18}=0\).  Thus \(I_{29}^{\leq29}=0\).  Also
\(\varphi(I_{29})=0\), since \(B\) is three-short.

For \(30\leq q\leq47\), we construct differential ideals
\[
                 I_{29}\subseteq I_{30}\subseteq\cdots\subseteq I_{47}
\]
such that
\begin{align}
 H^j(I_q)&=0 &&(j\leq q),
 \label{app:eq:Iq-acyclic-range}\\
 \varphi(I_q^j)&=0 &&(j\geq q+1).
 \label{app:eq:Iq-zero-image-range}
\end{align}
These assertions hold when \(q=29\).  Suppose that \(I_{q-1}\) has been
constructed.  Choose cocycles
\[
                    z_1,\ldots,z_m\in I_{q-1}^q
\]
whose classes form a basis of \(H^q(I_{q-1})\).  The second induction
hypothesis gives \(\varphi(z_i)=0\).  Since \(H(\varphi)\) is injective,
there are elements \(e_i\in M^{q-1}\) with \(de_i=z_i\).  Their images are
cycles; as \(q-1\geq29\), Lemma~\ref{app:lem:high-cycles-cubic} gives
\[
                       \varphi(e_i)\in(B^+)^3.
\]
Define
\[
                       I_q=I_{q-1}+(e_1,\ldots,e_m).
\]
This is a differential ideal because \(de_i=z_i\in I_{q-1}\).

Since \(M^1=M^2=0\), every product \(e_i\mM\) starts in degree
\((q-1)+3=q+2\).  In degrees at most \(q\), the quotient complex
\(I_q/I_{q-1}\) is therefore concentrated in degree \(q-1\), with basis
given by the classes of the \(e_i\).  These classes are independent: a
relation \(\sum_i c_i e_i\in I_{q-1}\), after differentiation, would imply
that \(\sum_i c_i[z_i]=0\) in \(H^q(I_{q-1})\), and hence every \(c_i=0\).
In the long exact cohomology sequence
associated with
\[
 0\longrightarrow I_{q-1}\longrightarrow I_q
   \longrightarrow I_q/I_{q-1}\longrightarrow0,
\]
the connecting morphism sends \([e_i]\) to \([z_i]\).  It is an
isomorphism onto \(H^q(I_{q-1})\); together with the inductive vanishing in
lower degrees, this proves \eqref{app:eq:Iq-acyclic-range}.

In degrees at least \(q+1\), each new term contains a factor \(e_i\) and a
positive-degree factor.  Its image is zero because
\[
                     \varphi(e_i)B^+\subseteq(B^+)^4=0.
\]
The old terms already have zero image, so
\eqref{app:eq:Iq-zero-image-range} follows.  This completes the induction.

At the last stage,
\[
                    H^{\leq47}(I_{47})=0,
                    \qquad \varphi(I_{47}^{48})=0.
\]
The cohomology computation of Proposition~\ref{div:prop:B-properties}
gives
\[
 H^{48}(B)=\Q\{[r^3]\},\qquad H^{>48}(B)=0.
\]
Since \(\varphi\) is a quasi-isomorphism, there is a top class
\([\Omega]\in H^{48}(M)\) satisfying
\[
 H^{48}(M)=\Q[\Omega],\qquad H^{>48}(M)=0,
 \qquad H(\varphi)([\Omega])=[r^3]\neq0.
\]
No cocycle in \(I_{47}^{48}\) can represent a nonzero multiple of
\([\Omega]\).  Lemma~\ref{calc:lem:top-completion}, applied with \(k=47\),
therefore embeds \(I_{47}\) in an acyclic differential ideal \(J\lhd M\).
Because
\[
                            \mM^5=I_{29}\subseteq I_{47}\subseteq J,
\]
Definition~\ref{div:def:Hnil} gives \(\Hnil(M)\leq4\).  The strict lower
bound proves the equality.
\end{proof}

\backmatter

\printsubjectindex

\end{document}